\documentclass[11pt,reqno]{amsart}
\usepackage[T1]{fontenc}
\usepackage[utf8]{inputenc}
\usepackage{lmodern}
\usepackage[margin=1.08in]{geometry}
\usepackage{amsmath,amssymb,amsthm,mathtools,mathrsfs}
\usepackage{booktabs,array,longtable,enumitem}
\usepackage[copyfonts]{microtype}
\usepackage{tikz-cd}
\usepackage{xcolor}
\usepackage[hidelinks]{hyperref}
\usepackage[capitalise,noabbrev]{cleveref}
\allowdisplaybreaks[2]
\numberwithin{equation}{section}

\newtheorem{theorem}{Theorem}[section]
\newtheorem{proposition}[theorem]{Proposition}
\newtheorem{lemma}[theorem]{Lemma}
\newtheorem{corollary}[theorem]{Corollary}
\theoremstyle{definition}
\newtheorem{definition}[theorem]{Definition}
\theoremstyle{remark}
\newtheorem{remark}[theorem]{Remark}
\newcommand{\F}{F}
\newcommand{\Dist}{\mathcal D'}
\newcommand{\cD}{\mathcal D}
\newcommand{\cN}{\mathcal N}

\newcommand{\Ga}{\mathbb G_a}
\newcommand{\Gm}{\mathbb G_m}
\newcommand{\CC}{\mathbb C}
\DeclareMathOperator{\Hom}{Hom}
\DeclareMathOperator{\GL}{GL}
\DeclareMathOperator{\SL}{SL}
\DeclareMathOperator{\PGL}{PGL}
\DeclareMathOperator{\GSp}{GSp}
\DeclareMathOperator{\Sp}{Sp}
\DeclareMathOperator{\GSO}{GSO}
\DeclareMathOperator{\SO}{SO}
\DeclareMathOperator{\Spin}{Spin}

\DeclareMathOperator{\GU}{GU}
\DeclareMathOperator{\Mat}{Mat}
\DeclareMathOperator{\tr}{tr}
\DeclareMathOperator{\diag}{diag}
\DeclareMathOperator{\Ad}{Ad}
\DeclareMathOperator{\Lie}{Lie}
\DeclareMathOperator{\WF}{WF}
\DeclareMathOperator{\CN}{CN}
\DeclareMathOperator{\supp}{supp}
\DeclareMathOperator{\rank}{rank}
\DeclareMathOperator{\im}{im}
\DeclareMathOperator{\vol}{vol}
\DeclareMathOperator{\NW}{NW}

\hypersetup{%
  pdftitle={Multiplicity One for Five Strongly Tempered Spherical Varieties},
  pdfauthor={Yuxuan Ge},
  pdfsubject={Multiplicity-one theorems for strongly tempered spherical varieties over non-Archimedean local fields},
  pdfkeywords={multiplicity one, strongly tempered spherical variety, p-adic group, equivariant distribution, nilpotent wavefront}
}

\usepackage{tabularx,multirow,pdflscape,graphicx,needspace}
\DeclareMathOperator{\spanop}{span}

\newcommand{\Geight}{G_8^\circ}
\newcommand{\De}{\Dist_{\mathrm e}}
\newcommand{\Do}{\Dist_{\mathrm o}}
\newcommand{\cG}{\mathcal G}
\newcommand{\cA}{\mathcal A}
\newcommand{\cP}{\mathcal P}
\newcommand{\cQ}{\mathcal Q}

\newcommand{\fu}{\mathfrak u}
\newcommand{\msevenfg}{\mathfrak{gsp}_6}

\newcommand{\calA}{\mathcal A}
\newcommand{\tA}{\widetilde{\mathcal A}}
\newcommand{\tK}{\widetilde K}
\newcommand{\teta}{\widetilde\eta}
\newcommand{\tetaK}{\widetilde\eta_K}
\newcommand{\Nil}{\mathcal N(\msevenfg^*)}
\newcommand{\Rres}{R_{\mathrm{res}}}

\newcommand{\sm}[1]{\begin{psmallmatrix}#1\end{psmallmatrix}}

\makeatletter
\renewcommand\paragraph{\@startsection{paragraph}{4}{\z@}{2ex plus .5ex minus .2ex}{.7ex}{\normalfont\normalsize\bfseries}}
\renewcommand\subparagraph{\@startsection{subparagraph}{5}{\z@}{1.5ex plus .4ex minus .2ex}{.5ex}{\normalfont\normalsize\itshape}}
\makeatother

\newcommand{\fg}{\mathfrak g}

\DeclareUnicodeCharacter{27E8}{\ensuremath{\langle}}
\DeclareUnicodeCharacter{27E9}{\ensuremath{\rangle}}
\DeclareUnicodeCharacter{2245}{\ensuremath{\simeq}}
\DeclareUnicodeCharacter{2218}{\ensuremath{\circ}}
\DeclareUnicodeCharacter{21A6}{\ensuremath{\mapsto}}
\DeclareUnicodeCharacter{03B5}{\ensuremath{\epsilon}}
\DeclareUnicodeCharacter{03BD}{\ensuremath{\nu}}
\DeclareUnicodeCharacter{03C7}{\ensuremath{\chi}}
\DeclareUnicodeCharacter{2208}{\ensuremath{\in}}
\DeclareUnicodeCharacter{00B9}{\ensuremath{{}^1}}
\DeclareUnicodeCharacter{22C9}{\ensuremath{\ltimes}}
\DeclareUnicodeCharacter{03B8}{\ensuremath{\theta}}
\DeclareUnicodeCharacter{03C4}{\ensuremath{\tau}}
\DeclareUnicodeCharacter{03B9}{\ensuremath{\iota}}
\DeclareUnicodeCharacter{2243}{\ensuremath{\simeq}}
\DeclareUnicodeCharacter{2228}{\ensuremath{\vee}}
\DeclareUnicodeCharacter{03C9}{\ensuremath{\omega}}
\DeclareUnicodeCharacter{03B6}{\ensuremath{\zeta}}
\DeclareUnicodeCharacter{03B3}{\ensuremath{\gamma}}
\DeclareUnicodeCharacter{0303}{\ensuremath{\sim}}
\DeclareUnicodeCharacter{03C8}{\ensuremath{\psi}}
\DeclareUnicodeCharacter{1D57}{\ensuremath{{}^t}}
\DeclareUnicodeCharacter{22CA}{\ensuremath{\rtimes}}
\DeclareUnicodeCharacter{2074}{\ensuremath{{}^4}}
\DeclareUnicodeCharacter{03BB}{\ensuremath{\lambda}}
\DeclareUnicodeCharacter{2194}{\ensuremath{\leftrightarrow}}
\DeclareUnicodeCharacter{03BA}{\ensuremath{\kappa}}
\DeclareUnicodeCharacter{2282}{\ensuremath{\subset}}
\DeclareUnicodeCharacter{00B1}{\ensuremath{\pm}}

\title[Multiplicity One for Five Spherical Varieties]{Multiplicity One for Five Strongly Tempered Spherical Varieties}
\author{Yuxuan Ge}
\address{School of Mathematics, University of Minnesota, Minneapolis, MN 55455, USA}
\email{ge005678@umn.edu}
\date{}
\subjclass[2020]{Primary 22E50; Secondary 22E35, 11F70, 20G41}
\keywords{multiplicity one, strongly tempered spherical varieties, $p$-adic reductive groups, equivariant distributions, nilpotent wavefront sets, Gelfand--Kazhdan method}
\begin{document}
\begin{abstract}
Let $F$ be a non-Archimedean local field of characteristic zero.
We prove multiplicity one for every irreducible smooth admissible
representation in the split forms of Wan--Zhang Models~6--10.
Combining the Gelfand--Kazhdan method with nilpotent-wavefront
restrictions on Bernstein-finite equivariant distributions, we
establish the required symmetries through rational orbit analysis,
parabolic descent, stabilizer characters, and conormal arguments. The type-$D_6$ analysis treats both the vector-similitude
group of Model~9 and the semispin form arising in the exceptional
$E_7$ model. For Model~10, the distribution theorem gives a product
bound for periods; the individual multiplicity bound follows from
Langlands classification and a vanishing theorem for standard
modules induced from nongeneric tempered data on the maximal $E_6$ Levi.
\end{abstract}
\maketitle
\tableofcontents
\section{Introduction}\label{sec:introduction}
\subsection{The multiplicity-one theorem}
Multiplicity-one problems for local periods ask whether an equivariant linear functional on an irreducible admissible representation of a reductive group over a local field is unique up to scalars. The local Gan–Gross–Prasad \cite{GGP12} conjectures provide a prominent example of the relationship between local period multiplicities and Langlands parameters. More broadly, the relative Langlands framework of Sakellaridis and Venkatesh \cite{SV17} places the harmonic analysis of spherical varieties in a systematic setting.

Within this framework, Wan and Zhang \cite{WZ} introduced ten strongly tempered spherical
varieties associated with classical and exceptional groups. For the
split forms of Models~6--10 over $p$-adic fields, we prove multiplicity one for all
irreducible smooth admissible representations. Their numbering is
recalled in the table below.

Let $F$ be a non-Archimedean local field of characteristic zero, and
fix a nontrivial smooth unitary additive character
$\psi:F\to\mathbb C^\times$. All representations are complex and smooth,
and $\pi^\vee$ denotes the smooth contragredient of $\pi$.
For the model datum $(G_i,H_i,\xi_i)$, put
\[
 m_i(\pi)=\dim\Hom_{H_i(F)}(\pi,\xi_i).
\]
In the table we list the ambient groups and the reductive parts of
their model subgroups. For Models~4--10 each subgroup also contains a
unipotent radical carrying the character specified in \cite{WZ}. We
give the precise data for Models~6--10 in their respective sections.
\begin{center}
\small
\begin{tabular}{@{}rll@{}}
\toprule
Model & $G_i$ & reductive part of $H_i$\\
\midrule
1 & $\GL_4\times\GL_2$ & $\GL_2\times\GL_2$\\
2 & $\GU_4\times\GU_2$ & $(\GU_2\times\GU_2)^0$\\
3 & $\GSp_6\times\GSp_4$ & $(\GSp_4\times\GSp_2)^0$\\
4 & $\GL_6$ & $\GL_2$\\
5 & $\GU_6$ & $\GU_2$\\
6 & $\GSp_{10}$ & $\GL_2$\\
7 & $\GSp_6\times\GL_2$ & diagonal $\GL_2$\\
8 & $\GSO_8\times\GL_2$ & diagonal $\GL_2$\\
9 & $\GSO_{12}$ & $\GL_2$\\
10 & $E_7^{\mathrm{ad}}$ & $\PGL_2$\\
\bottomrule
\end{tabular}
\end{center}
In the table, the superscript $0$ denotes pairs with the same
similitude multiplier. Models~6 and~7 are symplectic, while
Models~8 and~9 are orthogonal.
Here $\GSO_{2m}$ denotes the split connected even orthogonal similitude
group, and $E_7^{\mathrm{ad}}$ denotes the split adjoint group of type
$E_7$. The unitary models use the quadratic extensions specified in
\cite{WZ}.

Models~1 and~2 are treated together in work in preparation by the
author \cite{GeModels12}, and Model~3 in a separate work in preparation
\cite{GeModel3}. Multiplicity one for Model~4, the Ginzburg--Rallis model,
was established by Jiang, Li and Xi \cite{JLX}. Model~5 is the subject
of joint work in preparation with Mukherjee and Zhai \cite{GeModel5}.

\begin{theorem}[Multiplicity one for Models~6--10]\label{thm:main}
For $6\le i\le10$ and every irreducible smooth admissible complex
representation $\pi$ of $G_i(F)$,
\[
 \dim\Hom_{H_i(F)}(\pi,\xi_i)\le1.
\]
\end{theorem}

The model-specific statements are
Theorems~\ref{m6:main}, \ref{m7:main}, \ref{m8:main}, \ref{m9:main},
and~\ref{ep:gp1}, which also specify the central forms and model characters. We also obtain multiplicity one for the
semispin form of type $D_6$ arising in the $E_7$ descent
(Corollary~\ref{m9:semispin-gp1}).

To prove the theorem, we combine the Gelfand--Kazhdan method with
nilpotent-wavefront constraints on equivariant distributions, then
eliminate their possible supports by parabolic descent, rational orbit
geometry, stabilizer characters, and conormal arguments. For Model~10, the resulting distribution theorem
is supplemented by a Langlands-classification argument to obtain the
individual multiplicity bound.

For the classical part of the argument, the distribution theorems for
Model~7 and for the Shalika model on $\GL_4\times\GL_2$ provide the
two principal parabolic inputs to Model~6. For the orthogonal $D_6$ groups, we prove the vector-form and
semispin distribution theorems needed in Models~9 and~10.

\subsection{Distributions and support}
The generalized coefficients in the Gelfand--Kazhdan criterion satisfy
\[
 \text{generalized coefficients}
 \ \Longrightarrow\ \text{Bernstein finiteness}
 \ \Longrightarrow\ \text{nilpotent wavefront}.
\]
By \cite[Theorem A]{AGS}, Bernstein finiteness gives the
nilpotent-wavefront condition used in the subsequent geometric
arguments, including zero cotangent components in the central
directions. We denote this condition by NW in Definition~\ref{cf:native}.

In the support analysis we distinguish two kinds of vanishing. We
call vanishing intrinsic when every equivariant odd distribution on
the specified space is zero. Conditional vanishing applies to the
restriction of the original ambient NW distribution after it has
vanished on specified preceding opens. The ambient NW and preceding-support hypotheses are essential
on strata that support nonzero intrinsic odd distributions. We use the support convention of \S\ref{cf:support-convention} throughout.

For Models~6--9, let $\sigma_i$ denote the anti-involution in the
corresponding distribution theorem, and put
$H_i^{\sigma_i}=\sigma_i(H_i)$. The second factor in the
Gelfand--Kazhdan inequality is
\[
 q_i(\pi)=\dim\Hom_{H_i^{\sigma_i}(F)}
       (\pi^\vee,\chi_i^{\mathrm{op}}),\qquad
 \chi_i^{\mathrm{op}}(\sigma_i(h))=\xi_i(h)^{-1}.
\]
In the matrix realizations below, $\sigma_i$ is transpose. For each model, duality identifies $q_i(\pi)$ with $m_i(\pi)$.
For Model~10 the distribution theorem gives
\[
 m_{10}(\pi)m_{10}^{\vee}(\pi)\le1,\qquad
 m_{10}^{\vee}(\pi)=\dim\Hom_{H_{10}(F)}(\pi^\vee,\xi_{10}^{-1}).
\]
Part~V proves the individual bound, including the cases in which
one factor in the product vanishes.

\subsection{Classical geometry and central forms}
Model~6 is a split symplectic model with a self-dual flag of type
$2+2+2+2+2$. Its parabolic double-coset decomposition has ninety
rational cells (Proposition~\ref{m6:s:prop:cells-rational}), sixty of
which are excluded by stabilizer-character incompatibility. Descent handles the $2+6+2$ and $4+2+4$ parabolic opens through
Model~7 and the Shalika model, respectively. We treat the remaining
cells by the boundary, stabilizer, and conormal arguments of Part~II.

Model~8 is an orthogonal model built from a self-dual filtration and
an external $\GL_2$ framing. Here the proof proceeds through rational relative-position strata:
we determine stabilizer-character compatibility, construct the support
filtration, and use ambient conormals to remove the spectral strata.

In type $D_6$ we treat two central forms: the vector-similitude
group $\GSO_{12}$ of Model~9 and the semispin Levi form occurring in
split adjoint $E_7$. For both forms, Siegel descent to the Ginzburg--Rallis distribution
theorem clears the open part, and the $35$-cell boundary continuation
completes the argument. For the vector form, the parabolic associated
with a two-dimensional isotropic subspace has Levi
$\GL_2\times\GSO_8$, giving a geometric relation with Model~8.

\subsection{Reductions and distributional inputs}
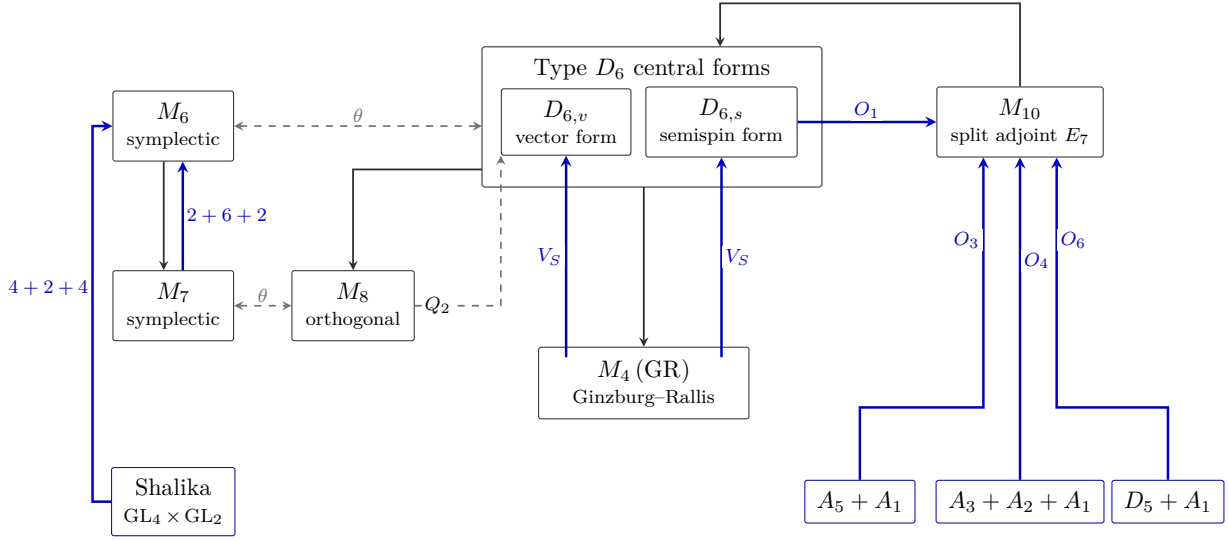
\begin{figure}[!t]
\centering
\begingroup
\resizebox{\textwidth}{!}{%
\begin{tikzpicture}[
  x=1cm,y=1cm,>=stealth,
  model/.style={draw=black!85,rounded corners=1pt,fill=white,
    align=center,inner xsep=5.8pt,inner ysep=5pt,font=\small},
  input/.style={draw=blue!60!black,rounded corners=1pt,
    align=center,inner xsep=5pt,inner ysep=4.8pt,font=\small},
  geom/.style={->,black!82,line width=.75pt},
  aux/.style={->,dashed,black!55,line width=.72pt},
  theta/.style={<->,dashed,black!55,line width=.72pt},
  dist/.style={->,blue!75!black,line width=1.05pt},
  edge/.style={fill=white,inner sep=1.2pt,font=\scriptsize},
  inputedge/.style={fill=white,inner sep=1.2pt,font=\scriptsize,
    text=blue!75!black}
]
\node[model] (m6) at (-6.50,2.85) {$M_6$\\[-1pt]{\scriptsize symplectic}};
\node[model] (m7) at (-6.50,.20) {$M_7$\\[-1pt]{\scriptsize symplectic}};
\node[model] (m8) at (-3.85,.20) {$M_8$\\[-1pt]{\scriptsize orthogonal}};
\draw[black!60,rounded corners=1.4pt,line width=.62pt]
  (-1.95,1.95) rectangle (3.06,4.02);
\node[font=\small] at (.555,3.69) {Type $D_6$ central forms};
\node[model] (d6v) at (-.71,2.91)
  {$D_{6,v}$\\[-1pt]{\scriptsize vector form}};
\node[model] (d6s) at (1.58,2.91)
  {$D_{6,s}$\\[-1pt]{\scriptsize semispin form}};
\node[model] (m10) at (5.98,2.91)
  {$M_{10}$\\[-1pt]{\scriptsize split adjoint $E_7$}};
\node[input] (shalika) at (-6.50,-2.69)
  {Shalika\\[-1pt]{\scriptsize $\GL_4\times\GL_2$}};
\node[model,minimum width=3.10cm] (gr) at (.435,-.94)
  {$M_4\,(\mathrm{GR})$\\[-1pt]{\scriptsize Ginzburg--Rallis}};
\node[input] (levia) at (3.62,-2.69) {$A_5+A_1$};
\node[input] (levib) at (5.98,-2.69) {$A_3+A_2+A_1$};
\node[input] (levic) at (8.16,-2.69) {$D_5+A_1$};
\draw[geom] ([xshift=-4pt]m6.south) -- ([xshift=-4pt]m7.north);
\draw[geom] (-1.95,2.21) -| (m8.north);
\draw[geom] (m10.north) -- (5.98,4.66) -- (1.56,4.66)
  -- (1.56,4.02);
\draw[geom] (.435,1.95) -- (gr.north);
\draw[theta] (m6.east) -- node[edge,above] {$\theta$} (-1.95,2.85);
\draw[theta] (m7.east) -- node[edge,above] {$\theta$} (m8.west);
\draw[aux] (m8.east) -| (d6v.south west);
\node[edge] at (-2.60,.20) {$Q_2$};
\draw[dist] ([xshift=4pt]m7.north) --
  node[inputedge,right] {$2+6+2$} ([xshift=4pt]m6.south);
\draw[dist] (shalika.west) -- (-7.69,-2.69) -- node[inputedge,left,pos=.57] {$4+2+4$} (-7.69,2.85)
  -- (m6.west);
\draw[dist] (-.71,-.56) -- node[inputedge,left] {$V_S$} (d6v.south);
\draw[dist] (1.58,-.56) -- node[inputedge,right] {$V_S$} (d6s.south);
\draw[dist] (d6s.east) -- node[inputedge,above] {$O_1$} (m10.west);
\draw[dist] (levia.north) -- (3.62,-1.30) -- (5.44,-1.30)
  -- node[inputedge,left,pos=.67] {$O_3$} ([xshift=-.54cm]m10.south);
\draw[dist] (levib.north) -- node[inputedge,right,pos=.68] {$O_4$} (m10.south);
\draw[dist] (levic.north) -- (8.16,-1.30) -- (6.52,-1.30)
  -- node[inputedge,right,pos=.67] {$O_6$} ([xshift=.54cm]m10.south);
\end{tikzpicture}%
}
\endgroup
\caption{Geometric reductions and distributional inputs for the five
strongly tempered spherical varieties. Black arrows indicate geometric
reductions in the Wan--Zhang framework, gray dashed arrows indicate
auxiliary relations (including theta correspondence and $Q_2$), and
blue arrows indicate distributional inputs used in the proofs.
Here $D_{6,v}$ is Model~9, and $D_{6,s}$ is the semispin form
used in the exceptional descent.
The $E_7$ inputs are labeled by their respective proper-Levi opens
$O_1,O_3,O_4,O_6$.}
\label{fig:reductions}
\end{figure}

Figure~\ref{fig:reductions} places the relevant Wan--Zhang geometric
reductions alongside the distributional inputs used in our proofs.
The two dashed theta correspondences are included for context, and
$Q_2$ marks the auxiliary geometric relation between Model~8 and
the vector form of Model~9. The common Model-4 Ginzburg--Rallis
theorem supplies the Siegel-open input for both $D_6$ central forms.
The other three proper-Levi inputs for Model~10 have types
$A_5+A_1$, $A_3+A_2+A_1$, and $D_5+A_1$; together with the
semispin $D_6$ input, they account for its four proper-Levi opens.

We establish the Ginzburg--Rallis and Shalika inputs in
Theorems~\ref{cf:gr-nw} and~\ref{cf:shalika-nw}. The Model-6 descents
appear in \S\ref{m6:model7} and \S\ref{m6:typea}; the semispin input
to the $E_7$ open $O_1$ is Theorem~\ref{m9:native-semispin}.
The boundary-support arguments complete the corresponding distribution
theorems.

\subsection{Exceptional localization and Langlands classification}
For Model~10 the finite compatibility calculation reduces $52\,416$
coarse relative positions to $98$ nonempty loci:
\[
 98=69+1+28.
\]
Of these loci, $69$ are covered by four proper-Levi opens of types
$A_5+A_1$, $A_3+A_2+A_1$, $D_6$, and $D_5+A_1$. One further locus
is the identity, and only $28$ require the residual exceptional
filtration. Most of the compatible exceptional geometry is therefore
controlled by classical proper-Levi distribution problems. The finite data determine relative positions, compatibility, and
support order, providing the geometry for the descent theorems and
conormal arguments below.

To pass from the distribution theorem to multiplicity one, we use
representation theory:
\[
 \begin{gathered}
 \text{$E_7$ distribution theorem}\ \Longrightarrow\ \text{product bound},\\
 \text{product bound and Langlands/$E_6$ argument}
 \ \Longrightarrow\ \text{individual bound}.
 \end{gathered}
\]
For the nonunitarizable Langlands quotients, proper Levi subgroups
of classical type are treated using rational duality and standard
properties of Langlands quotients. The remaining maximal $E_6$ Levi
requires a separate vanishing theorem for standard modules induced
from nongeneric tempered data. Together these arguments establish
the individual multiplicity bound.

\subsection{Organization}
After developing the common distributional framework in Part~I, we
prove the classical results in the order $7,6,8$ in Part~II and treat
both $D_6$ central forms in Part~III. Parts~IV and~V separate the $E_7$ distribution
theorem from its representation-theoretic completion. The final
section assembles Theorem~\ref{thm:main}.

We collect the rational-orbit, density, and conormal proofs in the
mathematical appendices. The model calculation PDFs contain repeated
matrix expansions and complete calculation records; we cite them by
filename and section identifier. The electronic ancillary materials
also provide the complete finite arrays and reproducible algorithms;
\path{anc/model10_gk_exhaustive/} contains the complete $E_7$ table
and Hasse data used in the geometric classification and support induction.

\subsection{AI Disclosure and Acknowledgements}
GPT-6 Astra (OpenAI) was used to assist with mathematical computations, proof development, and manuscript drafting. The author takes full responsibility for the mathematical content and the final text. The author is grateful to his advisor, Dihua Jiang, for introducing him to the study of multiplicity-one problems and for his generous guidance.

\part{Common distributional and microlocal framework}\label{part:1}
\section{Distributions, descent, and nilpotent wavefront}\label{sec:common}
We begin by fixing our conventions for distributions, actions, and
densities. Each application of Frobenius descent specifies the rational
stabilizer, quotient density, involution, and domain. In the geometric
arguments we also verify the transferred nilpotent-wavefront condition
for the descended distribution.
We write $\WF_y(T)$ for the wavefront fiber of $T$ at $y$, using
the full-cotangent convention: the zero cotangent vector is included
over every point of $\supp T$. The nilpotent-wavefront condition,
abbreviated NW, is defined in Definition~\ref{cf:native}.

\paragraph{Local notation}
In each model or explicit local chart we reset the letters
$G,H,H_0,P,Q,M,L,N,U$. We use $\xi$ for the model character and
$\eta$ for the character of the acting product group, with the stated
right-factor convention. Likewise, $\tau,\sigma,\iota$ refer to the
operations defined in that section. Cell labels are local to their
model, and we distinguish $D_6$ labels from coarse $E_7$ labels.

\subsection{Actions, test functions, and the two periods}
An $l$-space is a locally compact, totally disconnected Hausdorff space.
For an $l$-space $X$ let $\cD(X)=C_c^\infty(X)$ and
$\Dist(X)=\Hom_{\mathbb C}(\cD(X),\mathbb C)$. Thus test functions
are locally constant and compactly supported. If $A$ acts on points
and $\eta$ is its character, our convention is
\begin{equation}\label{cf:action}
 T(f\circ a)=\eta(a)^{-1}T(f),\qquad
 T(f\circ a^{-1})=\eta(a)T(f).
\end{equation}
If an involution $\sigma$ of $X$ normalizes the $A$-action and
preserves $\eta$, extend $\eta$ to $A\rtimes\langle\sigma\rangle$
by $\widetilde\eta(\sigma)=-1$. We write
$\Dist_{\widetilde\eta}(X)$ for distributions with these
equivariance equations and $\sigma$-odd parity.
For a model subgroup $H$ and anti-involution $\sigma$ of the ambient
group, the opposite character is specified by
$\chi^{\mathrm{op}}(\sigma(h))=\xi(h)^{-1}$. The point action
$g\mapsto h_1gk^{-1}$, with $k\in\sigma(H)$, has product character
$\xi(h_1)\chi^{\mathrm{op}}(k)$. Reparametrizing
$k=\sigma(h_2)^{-1}$ gives $g\mapsto h_1g\sigma(h_2)$ and
character $\xi(h_1)\xi(h_2)$. The model sections specify these
operations and the opposite subgroup explicitly.

\subsection{Support and the topology of induction}\label{cf:support-convention}
\begin{lemma}[Open restriction and closed support]\label{cf:support}
If $Z\subset X$ is closed with open complement $V$, then
\begin{equation}\label{cf:test-sequence}
0\longrightarrow\cD(V)\longrightarrow\cD(X)\longrightarrow\cD(Z)\longrightarrow0.
\end{equation}
For invariant $V,Z$, the equivariant distribution sequence is left exact:
\begin{equation}\label{cf:support-sequence}
0\longrightarrow\Dist_\eta(Z)\xrightarrow{i_*}\Dist_\eta(X)
 \longrightarrow\Dist_\eta(V).
\end{equation}
In particular, a distribution on $X$ vanishing on $V$ has a unique closed-support antecedent on $Z$.
\end{lemma}
\begin{proof}
Extend locally constant functions on $Z$ using a finite compact open refinement near their compact support. The kernel consists of functions whose compact support lies in $V$. Dualizing gives the ordinary distribution sequence, and taking equivariants preserves left exactness. Since $i_*$ is injective, the antecedent is equivariant. This is the support convention of \cite[\S2.1]{JLX}.
\end{proof}
\paragraph{Support convention}
All ordinary restrictions are to open subsets of the ambient $p$-adic
manifold. After vanishing on an open complement, the left-exact sequence
\eqref{cf:support-sequence} identifies the distribution with one supported
on the closed complement. As the invariant opens increase, we apply the NW argument to the
restriction of the original ambient distribution and intrinsic
vanishing to its closed-support antecedent. Vanishing on two opens implies vanishing
on their union by a compact-open partition of the support of a test
function. Noncharacteristic pullback, when used, is justified separately
by its wavefront hypothesis.

\subsection{Stabilizers and weighted density lines}\label{cf:stabilizer-section}
For a left Haar measure use
\[
\int_A f(xa)\,dx=\Delta_A(a)^{-1}\int_A f(x)\,dx.
\]
The orbital density condition in convention \eqref{cf:action} is
\begin{equation}\label{cf:orbital}
\eta^{-1}|_{A_x}=\Delta_A|_{A_x}\Delta_{A_x}^{-1}.
\end{equation}
We compute this condition on the $F$-rational stabilizer, including its finite components. Failure at every rational point of an algebraic locus gives vanishing of the equivariant distribution space by constructibility and the stabilizer criterion \cite[Theorems 2.3, 2.5]{JLX}. In particular, an additive stabilizer with character $t\mapsto\psi(ct)$, $c\ne0$, supplies an obstruction, since both modular characters on it are trivial. Differentiating a phase means differentiating its algebraic additive map before applying $\psi$.

If an involution preserves an orbit and has a character-one rational transporter, its action on the character-weighted density line is positive. Since the involution has square one, the positive scalar is one and the orbital line is even. For algebraic rational-point actions, constructibility passes this orbitwise conclusion to the locus by the argument of \cite[Lemma 2.7]{JLX}, using \cite[Theorem 2.3 and Remark 2.4]{JLX}.

\begin{lemma}[Constructive actions on invariant subspaces]
\label{cf:constructive-restriction}
Let an algebraic $F$-group $A$ act algebraically on an $F$-variety
$X$, and let $Y\subset X(F)$ be an $A(F)$-invariant, locally closed
subset for the analytic topology.  The action on the $l$-space $Y$
is constructive in the sense of Bernstein--Zelevinsky. Consequently,
if every orbit in $Y$ has zero distributions with the specified
character, then the corresponding distribution space on $Y$ is zero.
The same assertion holds after adjoining an algebraic finite-order
involution and using its sign character.
\end{lemma}
\begin{proof}
Constructiveness here refers to the orbit-equivalence relation. By
\cite[Definition 6.6 and \S6.15, Theorem A]{BZ76},
\[
 R_X=\{(x,a x):x\in X(F),\ a\in A(F)\}
       \subset X(F)\times X(F)
\]
is a finite union of analytically locally closed sets.
By invariance, $R_Y=R_X\cap(Y\times Y)$. Intersecting each locally
closed piece with the subspace gives a locally closed subset there,
so $R_Y$ is constructive in $Y\times Y$. Since $Y$ is locally closed
in the $l$-space $X(F)$, it is itself an $l$-space.
Apply \cite[Theorem 6.9]{BZ76} to the locally constant rank-one
$l$-sheaf with the specified character action, as in
\cite[Remark 2.4(2)]{JLX}.  An algebraic finite extension is again
an algebraic group, so the same reasoning applies to the extended
action. We can therefore apply the orbital criteria to the
analytically closed rational saturations and their locally closed
differences below, once their closedness and invariance have been
established.
\end{proof}

\begin{proposition}[Stabilizer criterion]\label{cf:stabilizer}
For an algebraic rational-point action, if the equality
\eqref{cf:orbital} fails on the stabilizer at every point,
then the equivariant distribution space is zero. In particular, a
nontrivial additive character on a unipotent stabilizer forces vanishing.
\end{proposition}
This is \cite[Theorems 2.3 and 2.5]{JLX}; both modular characters are
trivial on the unipotent obstruction. On a surviving homogeneous orbit
the weighted orbital line has the form
\begin{equation}\label{cf:orbital-integral}
 f\longmapsto\int_{A/A_x}f(ax)\eta(a)\,d\dot a,
\end{equation}
with its quotient density and the stabilizer compatibility in
\eqref{cf:orbital} understood.
\begin{proposition}[Localization]\label{cf:localization}
For an invariant map $q:X\to Y$ of the algebraic $l$-spaces used here,
with $q\circ\sigma=q$ whenever oddness under $\sigma$ is imposed,
if the specified equivariant distributions supported on every fiber
are zero, then the equivariant distributions on $X$ are zero.
\end{proposition}
Because the equivariance equations define a weakly closed space
stable under multiplication by tests from $Y$, Bernstein localization
\cite[Theorem 2.8 and Corollary 2.9]{JLX} applies.

\subsubsection{Rational transport and finite unipotent reconstruction}
\begin{lemma}[Transport from a normal form]\label{m6:s:lem:normal-transport}
Let $\tau(ax)=\theta(a)\tau(x)$, with $\chi\circ\theta=\chi$.
If $\tau(x_0)=k_0x_0$ and $\chi(k_0)=1$, then for $x=ax_0$,
\[
k_x=\theta(a)k_0a^{-1},\qquad k_xx=\tau(x),\qquad\chi(k_x)=1.
\]
Equivalently, if $nx=x_0$, then $k_x=\theta(n)^{-1}k_0n$.
\end{lemma}
\begin{proof}
The action identity gives $k_xx=\theta(a)\tau(x_0)=\tau(ax_0)$,
and $\chi(k_x)=\chi(a)\chi(k_0)\chi(a)^{-1}=1$.
Together with Lemma~\ref{m8:s:c:transporter}, this transfers an even orbital
line once its stabilizer and positive density have been computed.
\end{proof}

\begin{lemma}[Finite exponential and the algebraic phase]\label{m6:s:lem:finite-exponential}
Let $B$ be a closed algebraic matrix group over $F$. If $u\in B(F)$
is unipotent, then $\exp(t\log u)\in B(F)$ for every $t\in F$.
For nilpotent matrices $L,R$ with $Lx+xR^t=0$,
\[
\exp(tL)x\exp(tR)^t=x\qquad(t\in F).
\]
Whenever the two exponentials belong to the prescribed acting groups,
they therefore give an additive stabilizer subgroup. If its
character is $\psi\circ\Lambda$, where $\Lambda$ is an algebraic
homomorphism to $\mathbb G_a$, its value is
$\psi(t\,d\Lambda(L,R))$.
\end{lemma}
\begin{proof}
Since logarithm and exponential are finite polynomials, the equations
of $B$ evaluated on $\exp(t\log u)$ are polynomial in the parameter. They vanish
at every integer $t$, where the exponential equals $u^t$, and hence
vanish identically in characteristic zero.
The stabilizer identity follows from
$L^nx=(-1)^nx(R^t)^n$ and the finite exponential series.
The restriction of $\Lambda$ to this additive one-parameter subgroup
is an algebraic additive polynomial in characteristic zero, hence is
$t$ times its differential.
\end{proof}
In applications we first solve the reductive projection at group
level, including finite components and rational multipliers. Once
the scalar part is identified, the lemma reconstructs the specified
unipotent remainder.

\begin{lemma}[Rational homogeneous conclusion]\label{m6:s:lem:rational-homogeneous}
Suppose an algebraic $A$-variety $X$ is one geometric orbit and has a
point $x_0\in X(F)$. If $H^1(F,A_{x_0})=1$, so every stabilizer
cocycle is a coboundary, then $X(F)=A(F)x_0$. Once this rational transitivity is
established, \eqref{cf:orbital} and
Lemma~\ref{m8:s:c:transporter} determine the orbital line and its parity.
\end{lemma}
\begin{proof}
For $x=ax_0\in X(F)$ with $a\in A(\overline F)$,
$a^{-1}\sigma(a)\in A_{x_0}$ is a cocycle. Write it as
$b^{-1}\sigma(b)$ with $b\in A_{x_0}(\overline F)$.
Then $ab^{-1}\in A(F)$ carries $x_0$ to $x$, proving rational
transitivity. The density and parity assertions now follow from
the homogeneous-space statements above.
\end{proof}
In the local applications, after establishing geometric transitivity,
we use the stabilizer $\mathbb G_m\times U$, with $U$ split
unipotent. Its additive filtration and Hilbert 90 trivialize the
cocycle. In other cells, explicit rational normalizations give
rational transitivity directly.

\subsection{Frobenius descent and passive parameters}
\begin{proposition}[Frobenius descent]\label{cf:frobenius}
Suppose $A$ acts transitively on an analytic base $B$, $p:X\to B$ is equivariant, $K=A_b$, and $Y=p^{-1}(b)$. In convention \eqref{cf:action}, the descended distribution has effective character
\[
\eta|_K\Delta_A|_K\Delta_K^{-1};
\quad S(f\circ k)=\eta(k)^{-1}\Delta_A(k)^{-1}\Delta_K(k)S(f).
\]
The correspondence commutes with open restriction and is realized by integration on test functions.
\end{proposition}
This is Frobenius descent as in \cite[Theorem 2.5.7(i)]{AGdescent}. To apply it through an explicit action slice, we verify that the action map is a surjective submersion on the indicated locus \cite[Lemma 2.13]{JLX}, then identify its complete fibers and density line. These determine whether the slice correspondence is injective or an isomorphism. The formula $|\det T^*X|=|\det T^*_{\mathrm{vert}}|\otimes|\det T^*_{\mathrm{base}}|$ fixes the quotient density; every nontrivial modulus remains in the local statement.

For a free unimodular fiber $A$, equivariance extracts the line $\eta(a)\,da$. An exchange of two such fibers preserves sign only after their characters and measures have been compared. Passive parameters are handled by partial evaluation: $\cD(X\times Y)=\cD(X)\otimes\cD(Y)$, and all product tests span. This gives parameter-uniform intrinsic vanishing when the entire acting group fixes $Y$. In microlocal applications the wavefront condition must also hold in the passive variables.

\begin{lemma}[Wavefront transfer in torus Frobenius descent]
\label{cf:frobenius-wf}
Suppose $p:X\to F^\times$ is an analytic equivariant map, where
$A$ acts on the base through a surjective analytic character $q$.
Put $K=\ker q$ and $Y=p^{-1}(1)$.  Assume $q$ admits analytic local
sections and use the quotient density and character of
Proposition~\ref{cf:frobenius}.  If $S$ is the Frobenius descendant
of $T$, then at each $y\in Y$
\[
 \operatorname{WF}_y(S)\subset
 \{\omega|_{T_yY}:\omega\in\operatorname{WF}_y(T)\}.
\]
\end{lemma}
\begin{proof}
Choose a local section $s:V\to A$ with $s(1)=1$ and $q(s(b))=b$.
The analytic product chart
\[
 \Psi:V\times Y\longrightarrow p^{-1}(V),\qquad
 (b,y)\longmapsto s(b)y
\]
has inverse $x\mapsto(p(x),s(p(x))^{-1}x)$.
The test-function construction of Frobenius descent in this chart
expresses $\Psi^*T$ as $m\otimes S$, with $m$ the smooth nowhere-zero
character-weighted quotient density on $V$; any density conversion
is multiplication by a smooth nowhere-zero function and leaves
wavefront unchanged.  The prescribed modular factors are retained
in this density.

We verify the necessary direction of the product wavefront assertion
directly.  In analytic scalar coordinates near $1$, shrink $V$ so
that $m=c\,db$ with $c\ne0$.  If $(y,\xi)\in\operatorname{WF}(S)$,
$\xi\ne0$, but $((1,y),(0,\xi))$ were microlocally regular for
$m\otimes S$, choose a compact open cutoff $\alpha$ in its base
neighborhood with $\int\alpha m\ne0$.  For every cutoff $\beta$ in
the prescribed neighborhood of $y$, the product Fourier transform
at frequency $(0,\lambda\zeta)$ equals
\[
 \widehat{\alpha m}(0)\,
 \widehat{\beta S}(\lambda\zeta).
\]
Regularity of the product makes this zero for $\zeta$ near $\xi$
and all sufficiently large $|\lambda|$, contradicting
$(y,\xi)\in\operatorname{WF}(S)$.  For the full-cotangent convention,
the zero cotangent vector follows from
$\operatorname{supp}(m\otimes S)=V\times\operatorname{supp}S$.
Analytic coordinate covariance
\cite[Corollary 2.9.5]{CHLR} now gives a cotangent vector
$\omega\in\operatorname{WF}_y(T)$ with
$d\Psi_{(1,y)}^t\omega=(0,\xi)$.  Since the derivative in the $Y$
direction is the inclusion $T_yY\hookrightarrow T_yX$, its
restriction to $T_yY$ is exactly $\xi$.
\end{proof}

For the passive-parameter arguments the test-function identity is
\begin{equation}\label{cf:tensor}
 \cD(X\times Y)=\cD(X)\otimes\cD(Y).
\end{equation}
To see this, cover the compact support of a locally constant function by
finitely many product neighborhoods and refine them to compact open
rectangles on which the function is constant.
\subsubsection{Weighted products and passive parameters}
\begin{lemma}[Weighted product descent]\label{m8:s:c:product}
Let $A$ be a unimodular $l$-group, $\chi$ a continuous character, and $Y$ an $l$-space. For the left action of $A$ on $A\times Y$, every distribution satisfying \eqref{cf:action} has the unique form
\begin{equation}\label{m8:s:c:weighted-average}
 T(f)=R\left(y\longmapsto\int_A f(a,y)\chi(a)\,da\right),\qquad R\in\Dist(Y).
\end{equation}
This weak-star isomorphism may be iterated for two fiber factors. A normalizer preserving their characters and Haar measures acts on $R$ by its action on $Y$; exchange of identically normalized fibers has multiplier one.
\end{lemma}
\begin{proof}
After multiplying tests by $\chi^{-1}$, the equivariance equation becomes translation invariance. Uniqueness of the Haar line then follows from compact open coset indicators and their finite refinements. Partially evaluate on each test function of $Y$ and use \eqref{cf:tensor}; the resulting scalar coefficients define $R$. Choose $\beta_0\in\cD(A)$ with $\int\beta_0\chi=1$. Then $\phi\mapsto\beta_0\otimes\phi$ is a compactly supported right inverse of averaging. Both maps are dual to maps on tests, proving the weak-star assertion. Conversely $a=a_0a'$ in the integral gives the factor $\chi(a_0)$ for $f(a_0^{-1}a,y)$, fixing the sign. A change of fiber variables proves the last assertion.
\end{proof}

\begin{lemma}[Passive parameters]\label{m8:s:c:passive}
If the whole acting group, including its involution, fixes an $l$-space $Y$, vanishing of the specified equivariant distributions on $X$ implies their vanishing on $X\times Y$.
\end{lemma}
\begin{proof}
For $\beta\in\cD(Y)$, set $T_\beta(\phi)=T(\phi\otimes\beta)$. It has the stated equivariances on $X$ and therefore vanishes; by \eqref{cf:tensor}, $T=0$. This is \cite[Lemma 2.2]{JLX}.
\end{proof}

\subsubsection{Descent along explicitly identified fibers}
\begin{lemma}[Fiber descent]\label{m8:s:c:fiber}
Suppose an analytic slice $X:\Sigma\to Y$ gives a surjective submersion $\pi:K\times\Sigma\to Y$, $(k,s)\mapsto kX(s)$, with local sections and explicit free, locally trivial fiber actions
\(j\cdot(k,s)=(k\sigma_s(j)^{-1},b\cdot s),\) \(\sigma_s(j)X(s)=X(b\cdot s).\)
Assume the whole fiber, including ineffective scalar factors, has been identified; $K$ and the fiber groups are unimodular; the chosen density $dk\,d\sigma$ and fiber Haar density $dh$ are preserved by the action and its overlap normalizers. Assume the character is one on the stabilizer kernel and ineffective scalars, and $\eta(\sigma_s(j))=\chi(b)$ is independent of $s$.
Then the quotient density and the correspondence are
\begin{equation}\label{m8:s:c:fiber-correspondence}
 \omega=\frac{dk\,d\sigma}{dh},\qquad
 \pi^*T=\eta(k)\,dk\otimes S,\qquad
 S(f\circ b^{-1})=\chi(b)S(f).
\end{equation}
An involution lifted to $K\times\Sigma$ descends with its stated sign if it preserves both densities and the character factor. The same local Haar construction applies to a central torus quotient with explicit local sections and character one on the torus.
\end{lemma}
\begin{proof}
The tangent exact sequence $0\to T_{\rm vert}\to T(K\times\Sigma)\to\pi^*TY\to0$ gives
\begin{equation}\label{m8:s:c:density-line}
 |\det T^*(K\times\Sigma)|=|\det T_{\rm vert}^*|\otimes\pi^*|\det T^*Y|.
\end{equation}
The identified vertical torsor coordinates supply $dh$. Since translations and normalizers have modulus one, the quotient density is independent of the fiber coordinate and descends to a smooth nowhere-zero density. Relative-Haar pullback is then injective in local product charts. Lemma~\ref{m8:s:c:product} gives its unique $K$-equivariant form in \eqref{m8:s:c:fiber-correspondence}. For the residual relation $(k,s)\sim(kk_b^{-1},bs)$, inverse action on tests and change of the $K$ variable contribute $\eta(k_b)^{-1}$. Fiber invariance therefore requires $S(f\circ b^{-1})=\eta(k_b)S(f)$, with the positive character sign displayed above.

Conversely, this residual equation makes $\eta(k)dk\otimes S$ invariant under the complete fiber action. In each chart $H\times W$, Haar uniqueness identifies it with $dh\otimes T_W$. Because the overlap changes have modulus one, the $T_W$ agree and define a unique distribution on $Y$. This proves surjectivity, complementing the injectivity of the general slice map in \cite[Lemma 2.13]{JLX}. The involution assertion follows by changing the total and vertical variables. For a central torus quotient, a local section and a compactly supported fiber test of integral one give the same averaging, right inverse, and overlap argument.
\end{proof}

\subsubsection{Character-one transporters}
\begin{lemma}\label{m8:s:c:transporter}
Let an involution $\tau$ normalize an algebraic rational-point action through $\theta$, preserving its character. On a stable algebraic locus suppose every orbit carrying a nonzero even-equivariant distribution is $\tau$-stable and admits a transporter $k_xx=\tau x$ of character one, with the stated homogeneous density line. Then there are no odd equivariant distributions on that locus.
\end{lemma}
\begin{proof}
On a homogeneous orbit an equivariant distribution, when it exists, spans the character-weighted quotient-Haar line. The transporter preserves the character and contributes phase one; its action on the positive quotient density multiplies by a positive scalar. Since the involution has square one, this positive scalar is $+1$ and the orbital line is even. The extended orbits here coincide with the even orbits; constructibility for algebraic rational-point actions and \cite[Theorem 2.3 and Remark 2.4]{JLX} pass this conclusion to the locus. This is the positivity argument of \cite[Lemma 2.7]{JLX}, with constructibility as in \cite[Theorem 2.3 and Remark 2.4]{JLX}.

If $k_0x=y$ and $k_yy=\tau y$, then
\begin{equation}\label{m8:s:c:transport-back}
 \theta(k_0)^{-1}k_yk_0x=\tau x,\qquad
 \eta(\theta(k_0)^{-1}k_yk_0)=\eta(k_y).
\end{equation}
Conjugating the stabilizer also preserves its restricted character.
\end{proof}

\subsection{Elementary distribution lemmas}
\begin{lemma}[Matrix and whole-plane lemmas]\label{cf:elementary}
The following intrinsic odd distribution spaces are zero, also with passive $l$-space parameters fixed by the whole acting group:
\begin{enumerate}[label=(\roman*)]
\item simultaneous-conjugation-invariant distributions on $\GL_2(F)^2$ odd under simultaneous transpose;
\item distributions on all of $F^2$, invariant under $q(b,c)=(qb,q^{-1}c)$ for $q\in F^\times$, and odd under $(b,c)\mapsto(c,b)$;
\item distributions on $\GL_2(F)$ for the action $V\mapsto n(s)Vn(t)^t$, with character $\psi(2s+2t)$ and transpose oddness, where $n(s)=\left(\begin{smallmatrix}1&s\\0&1\end{smallmatrix}\right)$.
\end{enumerate}
\end{lemma}
We give the proofs below. For the plane, we use the symmetry on the union of the two axes from \cite[Lemma 4.6]{Prasad} and also treat the open locus $bc\ne0$ and the extension across their common origin.
\subsubsection{Matrix problems}
\begin{lemma}[Simultaneous conjugation]\label{m8:s:c:two-matrix}
Every distribution on $\GL_2(F)^2$ invariant under simultaneous conjugation and odd under simultaneous transpose is zero.
\end{lemma}
\begin{proof}
Put
$V=\mathfrak{sl}_2(F)$, $B(X,Y)=\tr(XY)$ and $Q(X)=B(X,X)/2$.
The form $B$ is nondegenerate, and $V$ has Witt index one. Writing $A=aI_2+X$ and $C=cI_2+Y$ separates the passive scalar
coordinates; invertibility is the condition
$a^2-Q(X)\ne0$, $c^2-Q(Y)\ne0$.

Conjugation and transpose on $V$ generate $O=\mathrm O(V)(F)$,
and their specified character is $\det:O\to\{\pm1\}$.
Indeed, conjugation has determinant one, transpose has determinant
minus one, and, for every anisotropic $u\in V$, reflection in $u$ is
\[
 s_u(Z)=Z-\frac{2B(Z,u)}{B(u,u)}u=-uZu^{-1}.
\]
For $\omega=\begin{psmallmatrix}0&1\\-1&0\end{psmallmatrix}$,
transpose equals $-\operatorname{Ad}(\omega)$, so $-I_V$ and then
all reflections belong to the generated group. Cartan--Dieudonn\'e
proves the assertion.

Localize by $(a,c,\Gamma)$, where
\[
 \Gamma(X,Y)=
 \begin{pmatrix}B(X,X)&B(X,Y)\\B(X,Y)&B(Y,Y)\end{pmatrix}.
\]
Since the invertible-pair domain is a union of fibers of this map,
Proposition~\ref{cf:localization} reduces the assertion to
$\Dist(E_\Gamma)^{O,\det}=0$ for each Gram fiber $E_\Gamma\subset V^2$.
If $\det\Gamma\ne0$, the span of the pair is a nondegenerate plane;
reflection in its anisotropic orthogonal line fixes the pair and has
determinant minus one. If $\Gamma=0$, the span is zero or an isotropic
line. At zero $-I_V$ fixes the pair; for an isotropic line $FN$ choose
an anisotropic $w\in N^\perp$ and use $s_w$.
In both cases every stabilizer contains an involution on which
$\det=-1$. Its modular characters are positive and equal one there,
so Proposition~\ref{cf:stabilizer} gives vanishing.

It remains to treat $\operatorname{rank}\Gamma=1$. A fixed invertible
change of the two vector coordinates, commuting with $O$, changes
$\Gamma$ to $\diag(\gamma,0)$ with $\gamma\ne0$. The fiber becomes
\[
 E=\{(X,Y): B(X,X)=\gamma,\ Y\in X^\perp,\ Q(Y)=0\}.
\]
Project to $S_\gamma=\{X:B(X,X)=\gamma\}$. This is one $O$-orbit.
For a fixed $X_0\in S_\gamma$, when $B(X,X_0)\ne-\gamma$ the
orthogonal map $-s_{X+X_0}$ sends $X_0$ to $X$. In the remaining
case $s_{X-X_0}$ does so. Near $X_0$ the first formula is a rational
continuous local section, since its denominator is
$2(\gamma+B(X,X_0))\ne0$. Translating this section covers the base.
Its stabilizer $H$ is $\mathrm O(X_0^\perp)(F)$, and the restriction
of $\det$ is the determinant on this binary plane.
Both groups are unimodular. For $O$, its identity component is
$\mathrm{PGL}_2(F)$ and its other component differs by the central
$-I_V$. The stabilizer is compact for an anisotropic binary plane; for a
hyperbolic plane it is
$F^\times\rtimes\langle\text{inversion}\rangle$, with multiplicative
Haar measure preserved by inversion.

To descend to the possibly singular fiber $Z=\{Y\in X_0^\perp:Q(Y)=0\}$,
use the local sections to identify $E$ with $O\times_H Z$. The averaging map
\[
 \operatorname{Av}\Phi(gz)=\int_H\Phi(gh,h^{-1}z)\,dh
\]
from $\cD(O\times Z)$ to $\cD(E)$ is well-defined and surjective:
local sections and a compactly supported $H$-test of integral one
lift any test, and finitely many compact open charts cover its support.
Pulling an $(O,\det)$-equivariant distribution through this map gives
$\det(g)\,dg\otimes R$ by Haar uniqueness and the test-function
tensor identity. Invariance under $(g,z)\mapsto(gh,h^{-1}z)$,
using unimodularity, gives
$R(f\circ h^{-1})=\det(h)R(f)$. Consequently vanishing on $Z$
under this $H$-character implies vanishing on $E$.

If $X_0^\perp$ is anisotropic then $Z=\{0\}$, and a reflection of
determinant minus one fixes it, giving zero. If it is hyperbolic,
choose coordinates with $Q(b,c)=bc$. Then $Z$ is the union of axes,
the determinant-one subgroup acts by $(b,c)\mapsto(tb,t^{-1}c)$,
and a determinant-minus-one element interchanges the axes.
On their punctured parts an equivariant distribution is
$\lambda_1d^\times b+\lambda_2d^\times c$, with
$\lambda_2=-\lambda_1$ by oddness. Normalize the unit annulus to
volume one. For a uniformizer $\varpi$, let
$f_0=1_{\mathcal O_F^2}|_Z$. The test
\[
 f_0(\varpi^{-1}b,\varpi c)-f_0(b,c)
\]
vanishes near zero and equals $-1_{\mathcal O_F^\times}$ on the first
axis and $1_{\varpi^{-1}\mathcal O_F^\times}$ on the second.
Torus invariance yields $0=-\lambda_1+\lambda_2$, so both coefficients
are zero. The remaining distribution is supported at the fixed origin
and vanishes by oddness. This extension calculation covers the
nonclosed orbits exchanged by transpose.

Vanishing on every Gram fiber now gives the assertion by
localization. Any further passive parameters can be included in
the base.
\end{proof}
For the one-matrix specialization, an explicit transporter is available. For scalar $K_0$, take $h=I_2$. Otherwise choose a cyclic vector $v$ among $e_1,e_2,e_1+e_2$, put $C_0=(v,K_0v)$, $t=\tr K_0$, $d=\det K_0$, and use
\begin{equation}\label{m8:s:s3:eq:8-11}
 C_0^{-1}K_0C_0=\begin{pmatrix}0&-d\\1&t\end{pmatrix},\qquad
 H_0=\begin{pmatrix}0&1\\1&t\end{pmatrix},\qquad h=C_0^{-t}H_0C_0^{-1}.
\end{equation}
Multiplication gives $hK_0h^{-1}=K_0^t$, $h=h^t$, and $\det H_0=-1$. Since the centralizer is $\GL_2$ in the scalar case and the units of a commutative two-dimensional algebra otherwise, it is unimodular. Lemma~\ref{m8:s:c:transporter} gives the one-matrix conclusion, also with passive parameters.

\begin{lemma}[The whole opposite-weight plane]\label{m8:s:c:plane}
For $q(b,c)=(qb,q^{-1}c)$, $q\in F^\times$, and $\tau(b,c)=(c,b)$,
\(\Dist(F^2)^{F^\times,\tau\text{-odd}}=0.\)
The same vanishing holds on $\GL_2(F)$ for diagonal-torus conjugation and transpose, with any passive parameters.
\end{lemma}
\begin{proof}
We recall the extension argument for the plane lemma of \cite[Lemma 4.6]{Prasad}, used in \cite[Lemma 4.3]{JLX}. On $bc\ne0$ use $(t,b)=(bc,b)$. By Lemma~\ref{m8:s:c:product}, an invariant distribution factors as a distribution in $t$ times $d^\times b$. Transpose sends $b$ to $t/b$, preserving this Haar measure, so the odd restriction is zero. Lemma~\ref{cf:support} leaves a distribution on the closed union of axes.

On the punctured axes write its restrictions as $\lambda_+d^\times b$ and $\lambda_-d^\times c$, each normalized to give $\mathcal O_F^\times$ volume one. Oddness gives $\lambda_-=-\lambda_+$. If $\varpi$ is a uniformizer and $f_0$ is the restriction of $1_{\mathcal O_F^2}$ to the axes, then $f_0\circ\varpi^{-1}-f_0$ vanishes near the origin and restricts to $-1_{\mathcal O_F^\times}$ on the first axis and $1_{\varpi^{-1}\mathcal O_F^\times}$ on the second. Invariance gives $0=-\lambda_++\lambda_-$. Both coefficients vanish. The remaining distribution is supported at the transpose-fixed origin and is therefore zero.

For matrices $V=\begin{psmallmatrix}a&b\\c&d\end{psmallmatrix}$, the diagonal entries are passive and the ratio of the torus entries runs through all $F^\times$. Lemma~\ref{m8:s:c:passive} gives vanishing on $M_2(F)$. Each closed determinant fiber injects by support into this space and therefore has the same vanishing. Localization for $\det:M_2(F)\to F$ gives vanishing on the union of nonzero fibers, $\GL_2(F)$, as in \cite[Corollary 2.10]{JLX}. Further passive factors use Lemma~\ref{m8:s:c:passive}.
\end{proof}

\begin{lemma}[The ordinary two-sided Whittaker conclusion]\label{m8:s:c:whittaker}
For $n(s)=\begin{psmallmatrix}1&s\\0&1\end{psmallmatrix}$, $\bar n(t)=n(t)^t$, the action $V\mapsto n(s)V\bar n(t)$ with character $\psi(2s+2t)$ has no transpose-odd distributions on $\GL_2(F)$, also with passive parameters.
\end{lemma}
\begin{proof}
In entries its action is
\begin{equation}\label{m8:s:v:e13-1}
 (a,b,c,d)\longmapsto(a+sc+bt+sdt,b+sd,c+dt,d).
\end{equation}
For $d\ne0$ the additive stabilizer is trivial, and, with $\Delta=\det V$, $\beta=b/d$, $\gamma=c/d$, the whole open set has product chart
\begin{align}
 V&=n(\beta)\diag(\Delta/d,d)\bar n(\gamma),
 &&\Delta,d\in F^\times,\quad\beta,\gamma\in F,\label{m8:s:v:e13-9}\\
 d^\times V&=|\Delta|^{-2}|d|\,d\Delta\,dd\,d\beta\,d\gamma.\label{m8:s:v:e13-10}
\end{align}
Since the action translates $\beta,\gamma$ while transpose exchanges them and fixes $\Delta,d$, Lemma~\ref{m8:s:c:product} gives exactly
\begin{equation}\label{m8:s:v:e13-11}
 \psi(2\beta+2\gamma)\,d\beta\,d\gamma\otimes S(\Delta,d,\text{passive parameters}),
\end{equation}
which is transpose-even.
On the closed boundary $Z=\{d=0\}$, invertibility gives $bc\ne0$. The complete stabilizer is $sc+bt=0$, hence $(s,t)=(bq,-cq)$, with character
\begin{equation}\label{m8:s:v:e13-2}
 \psi(2(b-c)q).
\end{equation}
The algebraic additive stabilizer criterion kills $Z\setminus Z_0$, where $Z_0=\{d=0,b=c\}$ is closed, stable, and pointwise fixed by transpose. Applying Lemma~\ref{cf:support} to the two pairs
\((X,Z)=(\GL_2,\{d=0\}),\) \((X,Z)=(\{d=0\},\{d=0,b=c\})\)
proves the assertion. All groups and boundary stabilizers here are additive and unimodular; the density factor in \eqref{m8:s:v:e13-10} depends only on passive variables.
\end{proof}

\subsection{Bernstein finiteness and the NW condition}
\begin{definition}\label{cf:native}
Let $M$ be a reductive $F$-group and $\mathfrak m=\Lie(M)$.
For an open $V\subset M(F)$, $\NW_M(V)$ consists of distributions
whose wavefront fibers, translated to the identity by left
translation, lie in the coadjoint nilpotent cone
$\cN(\mathfrak m^*)$. Equivalently, at $x\in V$ we identify
$T_xM(F)$ with $\mathfrak m$ by $\dot x=xA$. Under the derived
and central decomposition,
\[
\cN(\mathfrak m^*)=\cN([\mathfrak m,\mathfrak m]^*)\times\{0\}.
\]
Thus every central cotangent component is zero. We abbreviate $\NW_M(M(F))$ by $\NW(M)$.
\end{definition}
Our full-cotangent convention agrees with \cite{JLX}. A globally
Bernstein-finite distribution on $M(F)$ means that the span of its
Bernstein-center translates is finite dimensional. The connection
with wavefront is provided by the following theorem.
\begin{theorem}[Nilpotent wavefront]\label{cf:nilpotent}
A globally Bernstein-finite distribution on $M(F)$ has wavefront fibers
in $\cN(\mathfrak m^*)$, hence belongs to $\NW(M)$.
\end{theorem}
This is \cite[Theorem A]{AGS}. Open restriction retains the local wavefront fibers.

\begin{lemma}[Gauss transfer of the nilpotent cone]\label{cf:gauss}
On a Gauss open $NRN^-$ for a Levi $R$, suppose the equivariance extracts two smooth nowhere-zero character-weighted Haar factors. The descended distribution is NW on $R$ if the original distribution is NW on that open. Reconstruction from an NW descendant has the same property. The assertion also holds for the explicitly conjugated Gauss patches used below.
\end{lemma}
\begin{proof}
At the Levi slice the annihilator of the two outer tangent spaces is $\mathfrak r^*$, identified inside the ambient dual using an invariant form. Since the character-weighted factors have zero wavefront and Levi inclusion preserves Jordan decomposition, a cotangent vector in this annihilator is nilpotent in the ambient algebra exactly when it is nilpotent in the Levi algebra, including zero on its center. Coordinate covariance and group translations give the assertion over the open. This is the common calculation used in the proper-Levi and $D_6$ arguments.
\end{proof}

\subsection{Conormal translation and exclusion}
\begin{proposition}[Ambient conormal exclusion]\label{cf:conormal}
Let $Y$ be a smooth closed submanifold of an open $V\subset M(F)$. Suppose at each possible support point the ambient conormal $\CN^V_{Y,y}$ contains a nonnilpotent cotangent vector. Then $S\in\NW_M(V)$ supported on $Y$ is zero.
\end{proposition}
\begin{proof}
Around $y\in Y$ choose analytic product coordinates $(u,v)$ in which $Y=\{v=0\}$.
Restriction of locally constant compactly supported tests to $Y$ is
surjective. By local constancy and compactness, a test vanishing on
$Y$ has support disjoint from $Y$ and is annihilated by $S$. Thus
$S=R\otimes\delta_0$ in this chart, with a uniquely determined factor.
For any compactly supported cutoff $\theta$, put
$\alpha(u)=\theta(u,0)$. Then
\[
 \widehat{\theta S}(\xi,\eta)=\widehat{\alpha R}(\xi),
\]
independently of the normal frequency $\eta$.
At nonzero tangential frequency $\xi_0$, the definition of
microlocal regularity, using product frequency neighborhoods and
product cutoffs, consequently gives regularity at $(\xi_0,\eta_0)$
exactly when $R$ is regular at $\xi_0$.

To treat zero tangential frequency, suppose $u_0\in\operatorname{supp}R$. Then every cutoff $\alpha$ nonzero near
$u_0$ gives a nonzero compactly supported distribution $\alpha R$.
Injectivity of Fourier transform supplies some fixed $\kappa$
with $\widehat{\alpha R}(\kappa)\ne0$. For every neighborhood of $(0,\eta_0)$ and all
sufficiently large $|\lambda|$, the direction
$(\lambda^{-1}\kappa,\eta_0)$ lies in that neighborhood, whereas
\[
 \widehat{\theta S}\bigl(
 \lambda(\lambda^{-1}\kappa,\eta_0)\bigr)
 =\widehat{\alpha R}(\kappa)\ne0.
\]
Thus $(0,\eta_0)$ is not microlocally regular for any normal
$\eta_0$. With the full-cotangent convention, zero belongs over
support points as well. Combining these observations gives
\[
 \operatorname{WF}_{(u_0,0)}(R\otimes\delta_0)
 =\operatorname{WF}_{u_0}(R)\times(F^k)^*.
\]
Analytic coordinate covariance transports this equality back to $V$.
In particular every ambient conormal direction at a support point lies
in the wavefront fiber, and the entire fiber is invariant under
conormal translation. Its nonnilpotent conormal witness contradicts
$S\in\NW_M(V)$. Hence $S=0$. This also proves the conclusion for
residue characteristic two.
\end{proof}
The conormal is computed in the ambient reductive group, including central and similitude directions. In a cell filtration the smooth locus must be closed in the current open, and the calculation must cover its entire compatible locus.

\subsection{Rational central forms and covers}
\begin{proposition}[Connected split-torus quotients]\label{cf:torus}
Let $q:\widehat M\to M$ be one of the central quotients used below, with connected split-torus kernel. Then $q(F)$ is surjective. Generalized-function pullback, interpreted with the specified Haar densities, is injective and preserves NW, with zero wavefront in the added central directions. Descent requires trivial kernel character and the matching density line.
\end{proposition}
\begin{proof}
Hilbert 90 gives $H^1(F,\Gm^r)=0$, so the rational-point obstruction vanishes. Since the map is an analytic submersion, pullback adds a smooth Haar direction in local product charts. The agreement of the derived Lie algebras and the centrality of the new directions then give the wavefront assertion. Kernel invariance with the stated density supplies the local Haar descent and its compatibility on overlaps.
\end{proof}

\begin{proposition}[Finite central covers, sector by sector]\label{cf:finite-cover}
For a finite central isogeny in characteristic zero, the map on rational points is locally an analytic isomorphism onto its image. One must retain the rational image and every relevant coset. Transfer of an equivariant distribution requires the preimage on rational points of the acting subgroup, its character, the transported involution, and the finite-kernel equivariance. Nilpotent wavefront is preserved in these local charts. Extension by zero is used only from a clopen image sector.
\end{proposition}
In the $D_6$ application the finite kernels, rational cokernel, square classes, and central parities are calculated explicitly. The connected common cover is treated by Proposition~\ref{cf:torus}.

\subsection{Noncharacteristic fibers and central balancing}
For central balancing we use compactly supported averaging over the quotient $F^\times$ and exact Mellin inversion. After an analytic change of coordinates, the wavefront estimate reduces to compact-parameter projection, which we verify directly below. The noncharacteristic pullback and proper-pushforward theorems are \cite[Theorems 2.9.3 and 2.9.6]{CHLR}; throughout, their Section~2 applies with $\Lambda=F^\times$.

\subsubsection{Central balancing}\label{pl:o:sec:4.3}

\begin{lemma}[Central balancing]\label{pl:central-balancing} Let $X$ be a reductive $F$-group with anti-involution
$\sigma$, and let $R$ act on its two sides by
$x\mapsto r_1x\sigma(r_2)$. Let $\epsilon:R(F)\to\mathbb C^\times$
be the fixed smooth character and $\nu:R(F)\to F^\times$ a surjective homomorphism induced by an algebraic character. Suppose a surjective algebraic group character
$p:X\to\mathbb G_m$ satisfies
\begin{equation}
 p(r_1x\sigma(r_2))=(\nu(r_1)\nu(r_2))^k p(x),\qquad
 p(\sigma x)=p(x),\label{pl:o:eq:4.1}
\end{equation}
for a positive integer $k$. Assume the relevant groups are unimodular
and that there is a one-parameter subgroup in $R$ whose image in $X$ is central and whose left action changes $p$
by a nonzero integral power of its parameter. Use invariant densities
on $X$ and on the kernel of $p$, disintegrated against multiplicative
Haar measure on the base; assume the Frobenius modular ratio for
the displayed base action is one.

Assume that, for every smooth unitary character $\chi$ of $F^\times$, the
$R\times R$ odd-equivariant nilpotent-wavefront problem with character
$\epsilon(r_1)\epsilon(r_2)\chi(\nu(r_1)\nu(r_2))$ vanishes.
Then the odd-equivariant nilpotent-wavefront problem also vanishes for
\begin{equation}
 A_0=\{(r_1,r_2):\nu(r_1)\nu(r_2)=1\},\label{pl:o:eq:4.2}
\end{equation}
with character $\epsilon(r_1)\epsilon(r_2)$. The assertion remains true
with passive analytic parameters and wavefront zero in their directions.
\end{lemma}

\begin{proof}
Put $A=R(F)^2$, $q(r_1,r_2)=\nu(r_1)\nu(r_2)$ and
$\eta_0(r_1,r_2)=\epsilon(r_1)\epsilon(r_2)$, and write $a_g$ for
the point action of $g\in A$. All distributions in this proof act on
test functions, with convention~\eqref{cf:action}.
Let $S$ be a distribution in the asserted balanced odd NW space.
Multiplication by $b\circ p$, $b\in C_c^\infty(F^\times)$, preserves
these conditions. For each test function, choose such a cutoff equal to one on the
image of its support. It therefore suffices
to treat $S$ with $\operatorname{supp}S\subset p^{-1}(K)$ for compact
$K\subset F^\times$.

The quotient $A/A_0$ is $F^\times$ through $q$, including its topology.
Here are explicit local analytic sections. If $z:\Gm\to R$ is the
given one-parameter subgroup, then $\nu(z(t))=t^j$ and
$p(z(t))=t^{jk}$ for an integer $j\ne0$. The map $t\mapsto t^j$
is an analytic isomorphism on a sufficiently small neighborhood of
$1$. Its image is open; its translates, together with the surjectivity
of $q$, provide a local analytic lift $c\mapsto g_c\in A$ near every
$c\in F^\times$. In particular $q$ is open. Give this quotient the
multiplicative Haar measure $d^\times c$.

For a smooth unitary $\chi$ put $\eta_\chi(g)=\eta_0(g)\chi(q(g))$
and define
\[
 T_\chi(f)=\int_{F^\times}
       \eta_\chi(g_c)\,S(f\circ a_{g_c})\,d^\times c.
\]
This is independent of the lift: replacing $g_c$ by $g_ch$, $h\in A_0$,
multiplies its character by $\eta_0(h)$ and its distribution evaluation
by $\eta_0(h)^{-1}$. The integrand is locally constant. It vanishes
unless
\[
 c^k\in p(\operatorname{supp}f)K^{-1}.
\]
The map $c\mapsto c^k$ is proper on $F^\times$: the valuation of $c$
is bounded whenever that of $c^k$ is bounded, and the unit group is
compact. Thus the integral has compact support and defines a distribution.
Changing $g_c$ to $g'g_c$ and using Haar invariance gives
$T_\chi(f\circ a_{g'})=\eta_\chi(g')^{-1}T_\chi(f)$.
Moreover $\sigma a_{(r_1,r_2)}=a_{(r_2,r_1)}\sigma$.
The swap preserves $q$, $\eta_\chi$ and quotient Haar measure,
so $T_\chi$ is odd.

To verify the NW condition for this average, use the nonzero
wavefront convention of \cite[\S2]{CHLR}, with $\Lambda=F^\times$,
and the extended wavefront convention of Definition~\ref{cf:native}.
On a compact open target neighborhood $V\subset X(F)$,
the relevant $c$ lie in a compact set. Cover that set by finitely many
disjoint compact open pieces $C_i$, each contained in a local analytic
section domain, and set
\[
 \Phi_i:C_i\times X(F)\longrightarrow X(F),\qquad
 \Phi_i(c,x)=a_{g_c}x,\qquad
 u_i=\bigl(\eta_\chi(g_c)d^\times c\bigr)\otimes S.
\]
The contribution of $C_i$ to the displayed average is exactly
$\Phi_{i,*}u_i$; this follows from its definition on tests. The map $\Phi_i$ is
a submersion because its differential in $x$ is invertible. It is
proper on the support of $u_i$: for compact $B\subset X(F)$ its
inverse image lies in the compact set
\[
 \{(c,a_{g_c}^{-1}y):c\in C_i,\ y\in B\}
\]
and is closed there. The smooth parameter factor has empty wavefront.
We use the following precise compact-parameter projection bound:
\[
 \operatorname{WF}(\pi_*D)\subset
 \{(y,\xi):\exists c,\ ((c,y),(0,\xi))\in\operatorname{WF}(D)\}.
\]
To verify the required form directly, take a point outside the
right-hand side, cover the compact parameter support by finitely many
microlocally smooth neighborhoods for $(0,\xi)$, choose a compact open
partition $\gamma_j(c)$ subordinate to them, and shrink to common
neighborhoods in $y$ and $\xi$. For any target cutoff $\varphi$, each
Fourier transform
$\mathcal F((\gamma_j\otimes\varphi)D)(0,\lambda\xi)$
vanishes for sufficiently large $|\lambda|$, uniformly in the common
direction neighborhood. Their sum is
$\mathcal F(\varphi\pi_*D)(\lambda\xi)$.
This proves the bound, which is the projection case of
\cite[Theorem 2.9.6]{CHLR}.
Factor $\Phi_i$ as the analytic isomorphism
$(c,x)\mapsto(c,a_{g_c}x)$ followed by this projection.
By wavefront covariance and the tensor bound
\cite[Corollary 2.9.5 and Theorem 2.8.10]{CHLR}, any nonzero
cotangent vector $\lambda$ in the wavefront set of this contribution
pulls back to a source cotangent vector whose $x$-component is
$(d_xa_{g_c})^t\lambda\in\operatorname{WF}(S)$.
This component is nonzero, since $d_xa_{g_c}$ is invertible.
Left and right group translations preserve identity-translated
nilpotent cotangent vectors, by coadjoint invariance of the nilpotent cone.
Hence $\lambda$ is nilpotent. For a submersion $N_{\Phi_i}=\varnothing$. A finite sum has wavefront contained in the union of its
summands' wavefront sets, proving $T_\chi\in\NW(X)$ locally on $V$
and therefore globally. The assumed full-product theorem gives
$T_\chi=0$.

For a fixed $f$ the well-defined function
\[
 h_f(c)=\eta_0(g_c)S(f\circ a_{g_c})
\]
lies in $C_c^\infty(F^\times)$, and the displayed average says
that its integral against every smooth unitary $\chi$ is zero.
For completeness, choose a uniformizer $\varpi$ and a compact open
unit subgroup $U$ such that $h_f$ is constant on its cosets.
Only finitely many $\varpi$-valuations occur. For each character
$\alpha$ of the finite group $\mathcal O_F^\times/U$ and every
$w\in\mathbb C$ with $|w|=1$, use
$\chi(\varpi^n u)=w^n\alpha(u)$. The resulting finite Laurent
polynomial in $w$ vanishes on the unit circle, so every coefficient
vanishes. Finite Fourier inversion on $\mathcal O_F^\times/U$ now
gives $h_f=0$. Choosing the lift $g_1=1$ yields $S(f)=h_f(1)=0$.

For passive parameters $P$, evaluate against a compactly supported
test $\beta$ in those parameters. The projection to $X(F)$ is proper
on the support after this cutoff, and its differential is surjective.
The same compact-parameter projection bound shows that the partial evaluation
is NW on $X$; equivariance and oddness are unchanged. The result just proved makes every such evaluation zero; since
product tests span $C_c^\infty(X(F)\times P)$, the original
distribution is zero.
All analytic power-map derivatives used above are nonzero in $F$,
including in residue characteristic two.
\end{proof}

\subsection{The GK product bound}
\begin{theorem}[Twisted Gelfand--Kazhdan criterion]\label{cf:gk}
Let $G$ be the rational-point group of a reductive group over $F$, $H_1,H_2$ closed subgroups with characters $\chi_1,\chi_2$, and $\sigma$ a continuous anti-involution satisfying
$\sigma(H_1)=H_2$ and $\chi_2(\sigma(h))=\chi_1(h)^{-1}$ for $h\in H_1$. Put $L_hf(x)=f(h^{-1}x)$ and $R_kf(x)=f(xk)$. Suppose every Bernstein eigen-distribution in the smooth-distribution class of \cite[Theorem 2.21]{JLX} satisfying
\[
D(L_hf)=\chi_1(h)D(f),\qquad D(R_kf)=\chi_2(k)D(f)
\]
is fixed by $\sigma$. Then every irreducible admissible $\pi$ satisfies
\begin{equation}\label{cf:product}
\dim\Hom_{H_1}(\pi,\chi_1)\,
\dim\Hom_{H_2}(\pi^\vee,\chi_2)\leq1.
\end{equation}
\end{theorem}
To apply the criterion, we identify the two Hom spaces using the
explicit action convention and apply Theorem~\ref{cf:nilpotent} to
the original ambient coefficient and its antisymmetrization. The
geometric theorems, which concern all equivariant NW distributions
in their stated domains, then give the required symmetry. We describe
the smoothing construction and coefficient-kernel argument below. Models~6--9 supply the exact duality comparison for their
periods. Model~10 uses the separate representation-theoretic argument
in Part~V.

\subsubsection{Generalized coefficients and their smooth targets}\label{cf:coefficients}
Let $(\pi,V)$ be irreducible admissible, and put
$A_\pi=\Hom_{H_1}(\pi,\chi_1)$ and
$B_\pi=\Hom_{H_2}(\pi^\vee,\chi_2)$. A period functional is an
element of the algebraic dual of the smooth representation. For $\lambda_1\in A_\pi$ define
\begin{equation}\label{cf:smoothing}
 \Phi_1(f)(v)=\int_G f(x)\lambda_1(\pi(x)v)\,dx,
 \qquad f\in\cD(G),\quad v\in V.
\end{equation}
This integral is defined by local constancy and compact support. If $f$
is right invariant under a compact open subgroup $J$, then
\[
 (\pi^\vee(j)\Phi_1(f))(v)
 =\int_G f(x)\lambda_1(\pi(xj^{-1})v)\,dx
 =\int_G f(yj)\lambda_1(\pi(y)v)\,dy=\Phi_1(f)(v).
\]
It follows that $\Phi_1(f)$ belongs to the smooth contragredient $V^\vee$,
so for $\lambda_2\in B_\pi$ we may define the ordinary distribution
\begin{equation}\label{cf:coefficient}
 C_{\lambda_1,\lambda_2}(f)=\lambda_2(\Phi_1(f)).
\end{equation}
With the regular actions of Theorem~\ref{cf:gk}, changes of variables give
\begin{equation}\label{cf:coefficient-actions}
 \begin{aligned}\Phi_1(R_kf)&=\pi^\vee(k)\Phi_1(f),\\
 \Phi_1(L_hf)&=\chi_1(h)\Phi_1(f),\\
 C(L_hf)&=\chi_1(h)C(f),\quad C(R_kf)=\chi_2(k)C(f).\end{aligned}
\end{equation}
Integration is with respect to Haar measure on the unimodular ambient group.

For $e_J=\vol(J)^{-1}1_J$ and $v\in V^J$, one has
\begin{equation}\label{cf:recover-period}
 \Phi_1(e_J)(v)=\lambda_1(v).
\end{equation}
If $\lambda_1\ne0$, right-regular equivariance makes
$\Phi_1:\cD(G)\to V^\vee$ a nonzero map with $G$-stable image,
which is therefore surjective by irreducibility. The analogous formula
\[
 \Phi_2(f)(\ell)=\int_G f(x)\lambda_2(\pi^\vee(x)\ell)\,dx
\]
lands in the smooth bidual, identified with $V$ by admissibility; for a
nonzero period it is likewise surjective and satisfies
$\Phi_2(R_kf)=\pi(k)\Phi_2(f)$. In particular a pair of nonzero
periods gives a nonzero coefficient.

For completeness, the scalar form of Schur's lemma used here follows
from admissibility. A commuting endomorphism preserves every finite-dimensional $V^J$.
On a nonzero such space it has an eigenvalue over $\mathbb C$, whose
eigenspace in $V$ is a nonzero subrepresentation. Irreducibility
therefore makes the endomorphism scalar.

\subsubsection{The center and the odd coefficient}\label{cf:coefficient-center}
Write $\mathfrak Z(G)$ for the Bernstein center, to distinguish it from
the Lie algebra of the center of $G$. Naturality of the Bernstein center on the right-regular smoothing map
gives
\[
 \Phi_1(Z_zf)=z_{V^\vee}\Phi_1(f),\qquad
 C(Z_zf)=\chi_\pi(z)C(f).
\]
Here $z_{V^\vee}=\chi_\pi(z)$ is scalar. For a product group the same
argument applies to the center of either factor and to that of the
whole ambient group, acting on the irreducible product representation.

If $\sigma$ is the relevant anti-involution, its action on central
test-function operators is defined by
\[
 P_\sigma Z_z=Z_{z^\sigma}P_\sigma,
 \qquad P_\sigma f=f\circ\sigma.
\]
Since an anti-involution reverses convolution and central convolution
is commutative, the transformed operator is again essentially compact
and central. Consequently $\sigma_*C$ has eigencharacter
$z\mapsto\chi_\pi(z^\sigma)$, and
\begin{equation}\label{cf:center-span}
 \mathfrak Z(G)(C-\sigma_*C)
 \subset\mathbb C C+\mathbb C\sigma_*C.
\end{equation}
Thus the difference is Bernstein-finite, and Theorem~\ref{cf:nilpotent}
places it in $\NW(G)$, with zero cotangent component along every
central factor. This supplies the common representation-theoretic
starting point for the geometric arguments, which transfer the NW
condition through their charts.

\subsubsection{The coefficient-kernel criterion}\label{cf:kernel}
Suppose every coefficient of the two period spaces is fixed by
$\sigma$. We show that this gives the product bound. Put
$\alpha=\sigma\circ\mathrm{inv}$ and
$\theta f=f\circ\alpha$. The automorphism $\alpha$ is involutive;
its Haar modulus is positive and has square one, hence is one. Define
\[
 (f_1\star f_2)(x)=\int_G f_1(g)f_2(xg)\,dg,
 \qquad
 \mathscr B(f_1,f_2)=\langle\Phi_2(f_1),\Phi_1(f_2)\rangle.
\]
Because the support of $f_1\star f_2$ lies in the compact set
$\supp(f_2)\supp(f_1)^{-1}$ and a sufficiently small stabilizer of
$f_2$ gives local constancy, this is a test function. After compact
open refinement the defining integrals reduce to finite sums. The smoothing identities imply
\[
 \mathscr B(f_1,f_2)=C(f_1\star f_2),\qquad
 (f_1\star f_2)\circ\sigma=(\theta f_2)\star(\theta f_1).
\]
To obtain the second identity, substitute $t=\alpha(g)$ and use
$\alpha(x)=\sigma(x)^{-1}$. Coefficient symmetry therefore yields
\begin{equation}\label{cf:bilinear-symmetry}
 \mathscr B(f_1,f_2)=\mathscr B(\theta f_2,\theta f_1).
\end{equation}
If $\Phi_2(f_1)=0$, the right side vanishes for every $f_2$.
Surjectivity of $\Phi_2$ and nondegeneracy of the pairing show
$\Phi_1(\theta f_1)=0$. Interchanging the maps gives the converse, so
\begin{equation}\label{cf:kernels}
 \ker\Phi_1=\theta(\ker\Phi_2),\qquad
 \ker\Phi_2=\theta(\ker\Phi_1).
\end{equation}
The relation $\theta R_g=R_{\alpha(g)}\theta$ specifies the twist
in this comparison; the quotient maps compared next have the same
untwisted right-equivariance and target.

Fix $\lambda_2\ne0$. For any two nonzero first periods, their
surjective maps $\Phi_1,\Phi_1'$ have the same kernel by
\eqref{cf:kernels}. The induced isomorphisms from this quotient of $\cD(G)$ to
$V^\vee$ differ by a commuting endomorphism and hence by a scalar. Evaluating on $e_J$ and $V^J$ in \eqref{cf:recover-period}
recovers $\lambda_1'=c\lambda_1$. Reversing the roles gives the same
assertion for the second space. Thus either one period space is zero
or both are one dimensional.

Applying an NW odd-distribution theorem to the antisymmetrization
of every relevant Bernstein eigen-distribution verifies the universal
symmetry hypothesis of Theorem~\ref{cf:gk}. Equivalently, applying it to
every coefficient constructed above and using the kernel argument gives
the same product bound.

\subsection{Ginzburg--Rallis and Shalika NW inputs}\label{cf:gr-shalika}

\begin{theorem}[GR--NW distributional theorem]\label{cf:gr-nw}
Let $R_{\rm GR}=D_3(\GL_2)U_{2,2,2}\subset\GL_6$, where
$D_3(a)=\diag(a,a,a)$, and put
\[
 \epsilon(D_3(a)u)=\psi(\tr u_{12}+\tr u_{23}),\qquad
 \nu(D_3(a)u)=\det a.
\]
For every smooth unitary character $\chi$ of $F^\times$, a distribution
$S\in\NW(\GL_6)$, equivariant for
\[
 (r_1,r_2)g=r_1gr_2^t,\qquad
 \eta_\chi(r_1,r_2)=\epsilon(r_1)\epsilon(r_2)
                  \chi(\nu(r_1)\nu(r_2)),
\]
and odd under transpose is zero.
\end{theorem}
\begin{proof}
We follow the support order of \cite{JLX}, giving the five
distributional arguments in Appendix~\ref{app:jlx-repairs}. The
character agrees with \cite[(3.4)]{JLX}, and we use the
nilpotent-wavefront hypothesis directly. Once the intrinsic
predecessors have vanished, Section~\ref{repair:jlx4} treats the
first two open-local NW steps. Section~\ref{repair:jlx57}, including
the rational intrinsic Lemma~\ref{cf:jlx-z62}, handles the ordinary
third cell, the fourth-open NW exclusion, and the final intrinsic
complement. Applying the open--closed support sequences in this
order proves the theorem. 
\end{proof}

For the Shalika reconstruction below, we use the \emph{open-local}
form of the proof: on $G^2_{{\rm open},1}$ and $G^2_{{\rm open},2}$, Section~\ref{repair:jlx4} kills every odd GR-equivariant NW distribution on the open itself. This gives the local vanishing needed for the reconstruction.

\subsubsection{The representation-valued Shalika theorem}
Let
\[
X_{\rm Sh}=\GL_4\times \GL_2,
\quad
R_{\rm Sh}=\{(D_2(a)u,a):a\in \GL_2,
\ u=\begin{pmatrix}I_2&X\\0&I_2\end{pmatrix}\},
\]
with $\epsilon=\psi(\tr X)$, $D_2(a)=\diag(a,a)$, and $\sigma(g,y)=(g^t,y^t)$. This is the two-sided distribution problem associated with a representation-valued Shalika model.

\begin{theorem}[Representation-valued Shalika]\label{cf:shalika-nw} For every unitary $\chi$, the odd $R_{\rm Sh}\times R_{\rm Sh}$-equivariant distributions on $X_{\rm Sh}$ with character
\[
\epsilon(r_1)\epsilon(r_2)\chi(\det a_1\det a_2)
\]
and nilpotent wavefront on $X_{\rm Sh}$ are zero. In particular both central cotangent components vanish, as in Definition~\ref{cf:native}.

\end{theorem}

\textbf{Proof.} Embed $X_{\rm Sh}$ as the Levi $\GL_4\times\GL_2$ of the parabolic $Q=(4,2)$ in $\GL_6$. Its unipotent radical N consists of a 4$\times$2 block. In the factorization $N\rtimes R_{\rm Sh}$ of the GR group, the character on N is the trace of the lower 2$\times$2 block and is stabilized by $R_{\rm Sh}$. The action on N has determinant one: for the $\GL_2$ part it is
\[
\det(D_2(a))^2(\det a)^{-4}=1,
\]
and the remaining part is unipotent.

Reconstruct the distribution on the Gauss open $NX_{\rm Sh}N^t$ in $\GL_6$. By Lemma~\ref{cf:gauss}, the reconstruction is NW; it has the GR character and is transpose-odd. We can therefore apply the open-local GR support argument above.

Write
\[
g_4=\begin{pmatrix}A&B\\C&D\end{pmatrix},\qquad y\in \GL_2.
\]
On the reconstructed open, the lower-right 2$\times$2 block of the $\GL_6$ matrix is y. Its lower 4$\times$4 corner has rank $2+\rank D$, and the two relevant rectangular lower corners have rank 2. Consequently $D\ne0$ corresponds exactly to the JLX open $G^2_{{\rm open},1}$, the union of the rank matrices
\[
\begin{pmatrix}4&2\\2&2\end{pmatrix},\qquad
\begin{pmatrix}3&2\\2&2\end{pmatrix}.
\]
Applying the open-local vanishing of Section~\ref{repair:jlx4} to this ambient NW distribution shows that the original distribution on $X_{\rm Sh}$ vanishes on $D\ne0$.

It remains supported on $D=0$. Apply Frobenius descent to the external y-coordinate using the transitive $\GL_2\times\GL_2$ action on $\GL_2$. On the fiber $y=I_2$, the reductive stabilizer is $b=a^{-t}$, so the two χ factors cancel; the density character is one. We are left with an action on $\GL_4$ by left upper $U_{2,2}$, right lower $U_{2,2}$, and conjugation by $D_2(a)$. The remaining argument on this Frobenius fiber is intrinsic.

For
\[
g=\begin{pmatrix}A&B\\C&0\end{pmatrix}\in \GL_4,
\]
both B and C are invertible. The entire $U_{2,2}\times U_{2,2}$ stabilizer equation is
\begin{equation}
XC+BR=0,\qquad R=-B^{-1}XC.
\label{cf:shalika-stabilizer}
\end{equation}
Its character is
\begin{equation}
\psi\bigl(\operatorname{Tr}X+\operatorname{Tr}R\bigr)
=\psi\bigl(\operatorname{Tr}((I_2-CB^{-1})X)\bigr).
\label{cf:shalika-phase}
\end{equation}
Thus the incompatible complement vanishes unless $C=B$. On $C=B$ an upper unipotent left translation kills A, leaving
\begin{equation}
\begin{pmatrix}0&B\\B&0\end{pmatrix}.
\label{cf:shalika-compatible}
\end{equation}
Every 2$\times$2 matrix over F is similar over F to its transpose. Choose a with $aBa^{-1}=B^t$; conjugation by $D_2(a)$ takes \eqref{cf:shalika-compatible} to its transpose and has character one. The character-one property transfers back through any unipotent normalization because transpose exchanges the two acting factors and preserves their product character. The odd stabilizer therefore acts incompatibly with its positive density line. The orbital criterion, algebraic constructibility, and the ordinary open--closed sequence therefore give zero on the compatible locus as in \cite[Theorem 2.5, Lemma 2.7 and Remark 2.4]{JLX}. The closed $D=0$ contribution vanishes, completing the proof. \hfill\(\square\)

\part{Classical descendants}\label{part:2}
\section{A symplectic descendant and Model 7}\label{cl:model7}\label{m7:section}

\paragraph{Proof strategy}
Frobenius descent over the external $\GL_2$ factor reduces the
NW distribution theorem to the symplectic similitude factor.
We follow the five rational rank strata $R_{7,0},\ldots,R_{7,4}$
of \S\ref{m7:ranks} in that order. Simultaneous conjugation and
transpose of two matrices give vanishing on $R_{7,0}$, followed by
stabilizer characters and rank-one matrix arguments on $R_{7,1}$.
Together these steps supply the predecessor open for $R_{7,2}$.
Characters and transpose transport then remove $R_{7,3},R_{7,4}$.

The $R_{7,2}$ step uses nilpotent wavefront: we make its smooth
residual support closed in an ambient open and exclude it by
nonnilpotent conormals. The separated equal-eigenvalue locus is
then removed intrinsically (\S\ref{m7:residual}). Throughout this
support induction we work with the same ambient distribution.
Frobenius reconstruction yields the product-group distribution
theorem; combining the GK product bound with inverse-transpose
duality gives multiplicity one
(\S\ref{m7:endpoint}).

\subsection{The datum and the two distributional statements}\label{m7:datum}
Let $j=\left(\begin{smallmatrix}0&-1\\1&0\end{smallmatrix}\right)$,
$J_7=w_3\otimes j$, where $w_3$ is the antidiagonal permutation matrix,
and $G_7^\circ=\GSp(J_7)(\F)$. The superscript here distinguishes the
first factor, with its full similitude multiplier.
Set $G_7=G_7^\circ\times\GL_2(\F)$ and $A^\#=j^{-1}A^tj$ for
$A\in\Mat_2(\F)$. The radical of the upper $2+2+2$ parabolic is
\[
 U_7=\left\{u(X,Z)=
 \begin{pmatrix}I_2&X&Z\\0&I_2&-X^\#\\0&0&I_2\end{pmatrix}:
 \tr Z=-\det X\right\}.
\]
Write $d_3(a)=\diag(a,a,a)$ and define
\begin{equation}\label{m7:datum-formulas}
 H_7=\{(u(X,Z)d_3(a),a):u(X,Z)\in U_7,\ a\in\GL_2(\F)\},
 \qquad \xi_7(u(X,Z)d_3(a),a)=\psi(\tr X).
\end{equation}
All multipliers in the ambient first factor occur. The character on the
diagonal $\GL_2$ is trivial.

\begin{theorem}[Model 7]\label{m7:main}
For every irreducible admissible smooth complex representation $\pi$ of
$G_7$,
\[
 \dim\Hom_{H_7}(\pi,\xi_7)\leq1.
\]
Here $\F$ is a non-Archimedean local field of characteristic zero, with
arbitrary residue characteristic, and $\psi$ is a nontrivial smooth
unitary additive character.
\end{theorem}

Let $\tau_7(x,b)=(x^t,b^t)$. Put
$\mathcal A=H_7\times H_7$ and
$\eta_7(h_1,h_2)=\xi_7(h_1)\xi_7(h_2)$ for its action
$(h_1,h_2)g=h_1gh_2^t$. Write
$\widetilde{\mathcal A}=\mathcal A\rtimes\langle\tau_7\rangle$
and extend $\eta_7$ by $\widetilde\eta_7(\tau_7)=-1$.
For ordinary distributions, pullback by $(h_1,h_2)\cdot g=h_1gh_2^t$ has character
$\xi_7(h_1)^{-1}\xi_7(h_2)^{-1}$. The odd extension gives $\tau_7$
character $-1$.

\begin{theorem}[The global nilpotent-wavefront statement]\label{m7:distribution}
Every ordinary distribution $T\in\NW(G_7)$ with this equivariance
and odd parity is zero. The nilpotent cone is
$\mathcal N(\mathfrak{sp}_6^*)\times\mathcal N(\mathfrak{sl}_2^*)
\times\{0\}_{\mathfrak z(G_7)^*}$, including both zero central components.
\end{theorem}
Equivalently,
\begin{equation}\label{m7:global-target}
 \Dist_{\widetilde\eta_7}(G_7)\cap\NW(G_7)=0.
\end{equation}

After descent, the action on the first factor is
\begin{equation}\label{m7:action}
 \begin{split}
 K_7&=(U_7\times U_7)\rtimes\GL_2(\F),\\
 (u,v,a)\cdot x&=u d_3(a)x d_3(a)^{-1}v^t,\\
 \eta_7(u,v,a)&=\psi(\tr X_u+\tr X_v).
 \end{split}
\end{equation}
Here $\eta_7$ on $K_7$ denotes the restriction of the preceding
$\mathcal A$-character to the identity-fiber stabilizer.
In the detailed calculations, we abbreviate these characters by
$\eta$ and $\widetilde\eta$, and the extended acting group by
$\widetilde{\mathcal A}$.
Transpose acts by $x\mapsto x^t$ and
$(u,v,a)\mapsto(v,u,a^{-t})$.

\begin{theorem}[The nilpotent wavefront input]\label{m7:native}
An ordinary $\eta_7$-equivariant, transpose-odd distribution $S$ on
all of $G_7^\circ$ is zero if
\[
 \WF_x(S)\subset\cN(\mathfrak{gsp}_6^*)\quad(x\in G_7^\circ),
\]
where cotangent fibers use $\dot x=xA$ and the trace pairing. In
particular the central cotangent component is zero.
\end{theorem}

We will use Theorem~\ref{m7:native} in the $2+6+2$ descent for Model~6.

\subsection{Frobenius descent with the transpose-extended stabilizer}\label{m7:frobenius}
In the coordinate proofs below write $G=G_7$, $H=H_7$, $U=U_7$,
$\xi=\xi_7$, $J=J_7$, $d(a)=d_3(a)$, and $R_r=R_{7,r}$.
We now carry out the transpose-extended Frobenius descent, including
the density and wavefront transfer.

\subsubsection{Forms, the subgroup, and its character}\label{m7:s:sec:conventions}

All vector spaces, algebraic groups, and matrix calculations are over $F$. We use the same notation for an algebraic group and its group of $F$-points when this does not cause ambiguity. Set
\begin{equation}\label{m7:s:eq:forms}
j=\begin{pmatrix}0&-1\\1&0\end{pmatrix},\qquad
J=\begin{pmatrix}0&0&j\\0&j&0\\j&0&0\end{pmatrix}.
\end{equation}
Thus $J^t=-J$ and $J^2=-I_6$. Put
\[
G_7^\circ=\{x\in\GL_6(F):x^tJx=\ell(x)J\},\qquad
G=G_7^\circ\times\GL_2(F).
\]
Every block in a $3\times3$ block matrix has size two. For $A\in\Mat_2(F)$, write
\begin{equation}\label{m7:s:eq:sharp}
A^\#=j^{-1}A^tj=\operatorname{adj}(A)=\tr(A)I_2-A.
\end{equation}
This notation applies to singular matrices as well. For a $6\times6$ matrix, its symplectic adjoint is $x^\dagger=J^{-1}x^tJ$.

Let $P=LU$ be the upper-block parabolic of type $2+2+2$. Its Levi and unipotent radical are
\begin{equation}\label{m7:s:eq:Levi}
L=\{m(A,E):A,E\in\GL_2(F)\},\qquad
m(A,E)=\diag\left(A,E,\frac{\det E}{\det A}A\right),
\end{equation}
\begin{equation}\label{m7:s:eq:U}
U=\left\{u(X,Z)=\begin{pmatrix}I_2&X&Z\\0&I_2&-X^\#\\0&0&I_2\end{pmatrix}:
\tr Z=-\det X\right\}.
\end{equation}
Indeed, writing the $(2,3)$ block of an upper-unitriangular matrix as $Y$, the equation $u^tJu=J$ gives $Y=-X^\#$ and $Z+Z^\#+\det(X)I_2=0$. The group law and inverse are
\begin{align}
u(X,Z)u(X',Z')&=u(X+X',Z+Z'-X(X')^\#),\label{m7:s:eq:U-law}\\
u(X,Z)^{-1}&=u(-X,-Z-\det(X)I_2).\label{m7:s:eq:U-inverse}
\end{align}
In particular, $X\mapsto\tr X$ is additive under multiplication. Independent coordinates on $U$ are
\begin{equation}\label{m7:s:eq:U-coordinates}
X\in\Mat_2(F),\qquad Z_0=Z+\tfrac12\det(X)I_2\in\mathfrak{sl}_2(F).
\end{equation}
Thus $\dim U=7$.

For $a\in\GL_2(F)$, let $d(a)=\diag(a,a,a)$. Then
\begin{equation}\label{m7:s:eq:diag-action}
d(a)u(X,Z)d(a)^{-1}=u(aXa^{-1},aZa^{-1}).
\end{equation}
Thus \eqref{m7:datum-formulas} defines a character of the Model-7 subgroup. On this subgroup,
\[
\ell(u\,d(a))=\det a.
\]
Under the map $(A,E,c)\mapsto(m(A,E),c)$, the Levi $L\times\GL_2$ becomes $\GL_2^3$, with $H_0=\{(d(a),a)\}$ as its diagonal $\GL_2$. This is the trilinear Levi pair in the Whittaker-induction description of \cite[\S8.2]{WZ}, giving the subgroup $H_0U$.

The ambient center meets $H$ in
\begin{equation}\label{m7:s:eq:center}
Z_G\cap H=\{(zI_6,zI_2):z\in F^\times\}.
\end{equation}
Since $\xi$ is trivial there, a nonzero element of $\Hom_H(\pi,\xi)$ forces the central character of $\pi$ to be trivial on \eqref{m7:s:eq:center}. An incompatible central character therefore gives multiplicity zero.

\subsubsection{Frobenius descent to \texorpdfstring{$\GSp_6$}{GSp6}}\label{m7:s:sec:frobenius}

\paragraph{The identity fiber and its transpose-extended stabilizer}

Projection to the second factor,
\[
p:G\longrightarrow\GL_2(F),\qquad p(g,z)=z,
\]
is equivariant for the $\calA$-action
\[
(a_1,a_2)\cdot z=a_1za_2^t.
\]
This action is transitive over $F$. A section carrying $I_2$ to $z$ is
\[
z\longmapsto\bigl((d(z),z),1\bigr)\in H\times H.
\]
Above the identity, the fiber is $G_7^\circ$ and the even
stabilizer satisfies $a_2=a_1^{-t}$. It is therefore
\begin{equation}\label{m7:s:eq:K-embedding}
K=(U\times U)\rtimes\GL_2(F),\quad
(u,v,a)\longmapsto
\bigl((u\,d(a),a),(v\,d(a^{-t}),a^{-t})\bigr).
\end{equation}
The semidirect action is
\[
a\cdot(u,v)=\bigl(d(a)ud(a)^{-1},d(a^{-t})vd(a^{-t})^{-1}\bigr).
\]
The induced action and character are
\begin{align}
(u,v,a)\cdot x&=u\,d(a)x\,d(a)^{-1}v^t,\label{m7:s:eq:K-action}\\
\eta_K(u,v,a)&=\psi(\tr X_u+\tr X_v).\label{m7:s:eq:K-character}
\end{align}
Since transpose fixes the base point $I_2$, the extended stabilizer is
\begin{equation}\label{m7:s:eq:K-odd}
\tK=K\rtimes\langle\tau\rangle,\quad
\tau(x)=x^t,\quad \tau_K(u,v,a)=(v,u,a^{-t}),
\quad \tetaK(\tau)=-1.
\end{equation}
In particular, the odd element is retained in Frobenius descent.

\paragraph{The density line}

Conjugation by $d(a)$ acts on \eqref{m7:s:eq:U-coordinates} by
\[
(X,Z_0)\longmapsto(aXa^{-1},aZ_0a^{-1}).
\]
The determinant is one on both $\Mat_2$ and $\mathfrak{sl}_2$, so $H$, $\calA$, and $K$ are unimodular. The same holds for their finite extensions $\tA$ and $\tK$, because an order-two automorphism has positive Haar scaling with square one. Hence the Frobenius coefficient line is trivial:
\begin{equation}\label{m7:s:eq:Frob-density}
\Delta_{\tA}|_{\tK}\,\Delta_{\tK}^{-1}=1.
\end{equation}
At the tangent space of the homogeneous base, the even stabilizer acts by $Z\mapsto aZa^{-1}$ and the odd element by $Z\mapsto Z^t$; both have absolute determinant one. This also computes the density line directly, including the odd component.

For comparison, the modulus of the parabolic $P$ is
\begin{equation}\label{m7:s:eq:parabolic-modulus}
\delta_P(m(A,E))=\left|\frac{\det A}{\det E}\right|^5.
\end{equation}
Indeed the determinants on $X$ and $Z_0$ are respectively $(\det A/\det E)^2$ and $(\det A/\det E)^3$. Its restriction to $d(\GL_2)$ is one. The descended character is therefore exactly \eqref{m7:s:eq:K-character}.

Frobenius reciprocity for distributions, with coefficient line \eqref{m7:s:eq:Frob-density}, gives
\begin{equation}\label{m7:s:eq:Frob-isom}
\Dist_{\teta}(G)\simeq\Dist_{\tetaK}(G_7^\circ).
\end{equation}
We use \cite[Theorem~2.5.7(i)]{AGdescent}, including its compatibility with open restrictions.

\paragraph{An explicit test-function realization and nilpotent wavefront}

For $f\in C_c^\infty(G)$, define
\begin{equation}\label{m7:s:eq:Frob-average}
\mathcal I f(x)=\int_{\GL_2(F)} f(d(z)x,z)\,dz.
\end{equation}
The average belongs to $C_c^\infty(G_7^\circ)$. To prove
surjectivity, take $\varphi\in C_c^\infty(G_7^\circ)$ and
$\beta\in C_c^\infty(\GL_2(F))$ of integral one, and set
\[
f(g,z)=\beta(z)\varphi(d(z)^{-1}g).
\]
The Frobenius lift of $S$ is
\begin{equation}\label{m7:s:eq:Frob-lift}
T_S(f)=S(\mathcal I f).
\end{equation}
Since the section in \eqref{m7:s:eq:Frob-average} has character one
and the base Haar measure is bi-invariant, this lift intertwines the
descended action and character
\eqref{m7:s:eq:K-action}--\eqref{m7:s:eq:K-character} with the two-sided
Model-7 action.

To check compatibility with the involutions, substitute $z\mapsto z^t$.
The integrand for $T_S(f\circ\tau)$ then becomes $f(x^td(z),z)$. Write
\[
x^td(z)=d(z)\bigl(d(z)^{-1}x^td(z)\bigr).
\]
Conjugation by $d(z)$ has character one on $S$, so the Frobenius lift of $\tau S$ is $\tau T_S$. This verifies the odd sign in \eqref{m7:s:eq:Frob-isom}.

For the right regular action of $G_7^\circ$ on the first factor, one has directly
\begin{equation}\label{m7:s:eq:Frob-right}
\mathcal I\bigl((g,z)\mapsto f(gy,z)\bigr)(x)
=(\mathcal I f)(xy).
\end{equation}
To transfer the wavefront condition, use the global diffeomorphism
$\kappa(x,z)=(d(z)x,z)$. In these coordinates, the lift is
$\kappa^*T_S=S\otimes dz$, with a smooth, nonzero Haar factor.
At $(x,I_2)$ its differential is
\[
 d\kappa(v,w)=(v+d d(w)x,w).
\]
Consequently, for each $\lambda\in\WF_x(S)$, the product cotangent vector
$(\lambda,0)$ corresponds to a cotangent vector $(\lambda,\beta)$ of $T_S$
at $(x,I_2)$, with
$\beta(w)=-\lambda(d d(w)x)$. Its first component is exactly
$\lambda$. The first-factor part of the nilpotent cone upstairs therefore
implies
\begin{equation}\label{m7:s:eq:Frob-native}
 T_S\in\NW(G)\quad\Longrightarrow\quad S\in\NW(G_7^\circ).
\end{equation}
In particular, the first-factor central component is zero, as follows
from the global product expression
$\kappa^*T_S=S\otimes dz$.

\begin{theorem}\label{m7:s:thm:descended}
For the action and character \eqref{m7:s:eq:K-action}--\eqref{m7:s:eq:K-odd},
\begin{equation}\label{m7:s:eq:descended-target}
S\in\Dist_{\tetaK}(G_7^\circ)\cap\NW(G_7^\circ)
\quad\Longrightarrow\quad S=0.
\end{equation}
\end{theorem}
We prove this theorem in the remaining sections. Combining
\eqref{m7:s:eq:Frob-isom} with \eqref{m7:s:eq:Frob-native} then gives
\eqref{m7:global-target}, and hence \cref{m7:main}.

The invariant $\ell(g)/\det z$ on the original product becomes the full similitude multiplier $\ell(x)$ on the identity fiber.

\subsection{Five rank strata and the intrinsic families}\label{m7:ranks}
The five rational parabolic cells $R_{7,0},\ldots,R_{7,4}$ are defined
by the rank matrices in \eqref{m7:s:eq:five-ranks}. The increasing
predecessor opens in \eqref{m7:s:eq:open-sets} determine the order
of the support argument.

\begin{proposition}[Intrinsic strata]\label{m7:intrinsic}
For $i=0,1,3,4$, the space of ordinary $\eta_7$-equivariant,
transpose-odd distributions on $R_{7,i}$ is zero.
\end{proposition}

The rank-one calculations include both semisimple and nilpotent
branches and the complete stabilizer characters.

\subsubsection{Rank stratification}\label{m7:s:sec:ranks}

For $x=(x_{ij})$ in $2\times2$ blocks, define
\begin{equation}\label{m7:s:eq:rank-matrix}
\mathcal R(x)=
\begin{pmatrix}
\rank\begin{pmatrix}x_{22}&x_{23}\\x_{32}&x_{33}\end{pmatrix}&
\rank\binom{x_{23}}{x_{33}}\\[2mm]
\rank(x_{32}\ \ x_{33})&\rank x_{33}
\end{pmatrix}.
\end{equation}
We reserve the notation
\begin{equation}\label{m7:s:eq:BCD}
B=x_{23},\qquad C=x_{32},\qquad D=x_{33}.
\end{equation}

\begin{proposition}\label{m7:s:prop:ranks}
There are exactly five possible rank matrices:
\begin{equation}\label{m7:s:eq:five-ranks}
\begin{aligned}
\mathbf r_0&=\begin{pmatrix}4&2\\2&2\end{pmatrix},&
\mathbf r_1&=\begin{pmatrix}3&2\\2&1\end{pmatrix},&
\mathbf r_2&=\begin{pmatrix}3&1\\1&1\end{pmatrix},\\
\mathbf r_3&=\begin{pmatrix}2&1\\1&0\end{pmatrix},&
\mathbf r_4&=\begin{pmatrix}2&0\\0&0\end{pmatrix}.
\end{aligned}
\end{equation}
Let $R_i=\{x:\mathcal R(x)=\mathbf r_i\}$. Each $R_i$ is one $P(F)$--$P(F)^t$ double coset, and is stable under $K$ and transpose.
\end{proposition}
\begin{proof}
Write $s=\rank D$. The inverse identity
\[
x^{-1}=\ell(x)^{-1}J^{-1}x^tJ
\]
and the complementary-submatrix rank formula show that the upper-left
entry of $\mathcal R(x)$ is $2+s$ and that its off-diagonal entries
are equal. Write their common value as $r$.

Let $V_+=\langle e_1,e_2\rangle$, $V_0=\langle e_3,e_4\rangle$, and $V_-=\langle e_5,e_6\rangle$. The plane $xV_-$ is isotropic. Its intersections with $V_+$ and $V_+^\perp$ have dimensions $2-r$ and $2-s$. Its image in the symplectic plane $V_+^\perp/V_+$ is isotropic, so $0\le r-s\le1$. Together with $s\le r\le2$ this gives exactly \eqref{m7:s:eq:five-ranks}.

The parabolic $P$ stabilizes $V_+$ and acts transitively over $F$ on isotropic planes with these fixed intersection dimensions. Choose a basis of the intersection, a basis of the projection to $V_0$, and complete to a symplectic basis. The unipotent radical removes the remaining lifts into $V_+$. All these basis operations are defined over $F$. Since $P^t$ stabilizes $V_-$, this is precisely the asserted double-coset statement. The rank conditions are invariant under the two unipotent actions and simultaneous block conjugation; transpose interchanges the equal off-diagonal ranks.
\end{proof}

The following open sets organize the support argument:
\begin{equation}\label{m7:s:eq:open-sets}
\begin{aligned}
\Omega&=\{\det D\ne0\}=R_0,\\
\mathcal V&=\left\{\rank\binom{B}{D}=2\right\}=R_0\sqcup R_1,\\
\mathcal W&=\{D\ne0\}=R_0\sqcup R_1\sqcup R_2.
\end{aligned}
\end{equation}
The set $R_1$ is closed in $\mathcal V$, and $R_2$ is closed in $\mathcal W$. The complement $\{D=0\}$ is closed in $G_7^\circ$, with open part $R_3$ and closed part $R_4$.

For later calculations we use
\begin{equation}\label{m7:s:eq:w234}
\begin{split}
w_2&=\begin{pmatrix}E_{11}&0&-E_{21}\\0&I_2&0\\E_{12}&0&E_{22}\end{pmatrix},\\
w_3&=\begin{pmatrix}0&E_{21}&-E_{12}\\E_{12}&0&-E_{21}\\E_{21}&-E_{12}&0\end{pmatrix},\qquad
w_4=\begin{pmatrix}0&0&I_2\\0&I_2&0\\I_2&0&0\end{pmatrix}.
\end{split}
\end{equation}
Here the matrix units in block displays are $2\times2$ matrix units. Multiplication gives $w_i^tJw_i=J$ and $\mathcal R(w_i)=\mathbf r_i$ for $i=2,3,4$. Hence every point of $R_i$ for these indices is $U\times U$-equivalent to
\begin{equation}\label{m7:s:eq:Levi-normal}
m(A,E)w_i m(C_0,F_0)^t,
\qquad A,E,C_0,F_0\in\GL_2(F).
\end{equation}

\subsubsection{The open stratum \texorpdfstring{$R_0$}{R0}}\label{m7:s:sec:R0}

\begin{proposition}\label{m7:s:prop:R0}
One has $\Dist_{\tetaK}(R_0)=0$.
\end{proposition}
\begin{proof}
Block Gaussian elimination, starting at the invertible block $D$, gives the unique factorization
\begin{equation}\label{m7:s:eq:open-factor}
R_0=ULU^t,\qquad x=u\,m(A,E)v^t.
\end{equation}
By the symplectic equation, the diagonal factor lies in
\eqref{m7:s:eq:Levi} and the triangular factors lie in
\eqref{m7:s:eq:U}. Left translation on $u$ and $v$ then extracts
their character-weighted Haar factors, leaving distributions on
$L\simeq\GL_2(F)^2$ with action
\[
a\cdot(A,E)=(aAa^{-1},aEa^{-1}),\qquad
\tau(A,E)=(A^t,E^t).
\]
The even character is trivial, while the odd character is $-1$.
Since conjugation by $d(a)$ has Haar modulus one on $U$ and
transpose exchanges the Haar factors, the residual density
character is trivial. Applying \cref{m8:s:c:two-matrix} proves
the assertion.
\end{proof}

\subsubsection{The first boundary \texorpdfstring{$R_1$}{R1}}\label{m7:s:sec:R1}

\begin{proposition}\label{m7:s:prop:R1}
One has $\Dist_{\tetaK}(R_1)=0$. Consequently $\Dist_{\tetaK}(\mathcal V)=0$.
\end{proposition}

Following the support-reduction method of \cite[\S5]{JLX}, we treat
the semisimple and nilpotent rank-one possibilities for $D$
separately, according to their distinct residual $\GL_2(F)$ orbits.

\paragraph{The first compatibility equation}

On $R_1$ one has $\rank D=1$ and $\rank\binom BD=2$. Define
\begin{equation}\label{m7:s:eq:R1-kappa}
\kappa(x)=\tr\bigl(D^\#(C-B)\bigr),\qquad
\mathcal C_1=\{x\in R_1:\kappa(x)=0\}.
\end{equation}
The two unipotent actions change
\[
B\longmapsto B-X^\#D,\qquad
C\longmapsto C-D(X_v^\#)^t,
\]
respectively. Both unipotent contributions to $\kappa$ vanish because
$DD^\#=D^\#D=0$, while the reductive component conjugates $B,C,D$
simultaneously. Hence $\kappa$ is $K$-invariant, and
\begin{equation}\label{m7:s:eq:R1-kappa-odd}
\kappa(x^t)=-\kappa(x).
\end{equation}
In particular, $\mathcal C_1$ is closed and stable under both parities.

\paragraph{Two rank-one normalizations}

The action on $D$ is conjugation. If $\delta=\tr D\ne0$, take
\[
D_{\mathrm s}=\begin{pmatrix}0&0\\0&\delta\end{pmatrix}.
\]
If $\tr D=0$, take
\[
D_{\mathrm n}=\begin{pmatrix}0&1\\0&0\end{pmatrix}.
\]
The corresponding centralizers are the diagonal torus and
\begin{equation}\label{m7:s:eq:R1-nil-centralizer}
\{c a_h:c\in F^\times,h\in F\},\qquad
 a_h=\begin{pmatrix}1&h\\0&1\end{pmatrix},
\end{equation}
respectively. These normalizations are over $F$.

Write $\mu=\ell(x)$. In the semisimple case, the left unipotent action kills the second columns of $B,x_{13}$ and the right action kills the second rows of $C,x_{31}$. The symplectic equations then split off the first and sixth coordinates. Define the resulting embedding by
\begin{equation}\label{m7:s:eq:R1-sembedding}
\mathscr E_{\mathrm s}(\delta,\mu;M)
=\frac\mu\delta E_{11}+\delta E_{66}
+\text{the matrix $M$ in rows and columns }(2,3,4,5).
\end{equation}
Here and below $E_{ij}$ in a $6\times6$ display is a $6\times6$ matrix unit.

In the nilpotent case, kill the second columns of $B,x_{13}$ and the first rows of $C,x_{31}$. The symplectic equations give the embedding
\begin{equation}\label{m7:s:eq:R1-nembedding}
\mathscr E_{\mathrm n}(\mu;M)=
\begin{pmatrix}
0&M_{11}&M_{12}&M_{13}&M_{14}&0\\
-\mu&0&0&0&0&0\\
0&M_{21}&M_{22}&M_{23}&M_{24}&0\\
0&M_{31}&M_{32}&M_{33}&M_{34}&0\\
0&0&0&0&0&1\\
0&M_{41}&M_{42}&M_{43}&M_{44}&0
\end{pmatrix}.
\end{equation}
These operations respect \eqref{m7:s:eq:U}. On the left, first prescribe $X^\#e_2=\delta^{-1}Be_2$ in the semisimple case, or $X^\#e_1=Be_2$ in the nilpotent case. Choose the other row of $X$ to vanish, so $\det X=0$, and prescribe the appropriate column of the trace-zero matrix $Z$ to kill the corresponding column of $x_{13}+XB+ZD$. On the right, prescribe $X_v^\#e_2$ from the relevant row of $C$ and choose the second column of $Z_v$ to kill the same row of $x_{31}+CX_v^t+DZ_v^t$. A column of a trace-zero $2\times2$ matrix can be prescribed freely. These operations give rational normalization charts in the ambient symplectic space.

The core matrices satisfy
\begin{equation}\label{m7:s:eq:R1-coreforms}
M^tJ_+M=\mu J_+\quad\text{in \eqref{m7:s:eq:R1-sembedding}},\qquad
M^tJ_-M=\mu J_+\quad\text{in \eqref{m7:s:eq:R1-nembedding}},
\end{equation}
where
\[
J_+=\begin{pmatrix}0&0&0&1\\0&0&-1&0\\0&1&0&0\\-1&0&0&0\end{pmatrix},\qquad
J_-=\begin{pmatrix}0&0&0&-1\\0&0&-1&0\\0&1&0&0\\1&0&0&0\end{pmatrix}.
\]

\paragraph{The residual unipotent equations}

Put
\begin{equation}\label{m7:s:eq:R1-npm}
n_+(a,b,c)=\begin{pmatrix}1&a&b&c\\0&1&0&-b\\0&0&1&a\\0&0&0&1\end{pmatrix},\qquad
n_-(a,b,c)=\begin{pmatrix}1&a&b&c\\0&1&0&b\\0&0&1&-a\\0&0&0&1\end{pmatrix}.
\end{equation}
Their embeddings in $U$ are
\[
n_+\leftrightarrow X=\begin{pmatrix}0&0\\a&b\end{pmatrix},\quad
Z=\begin{pmatrix}0&0\\c&0\end{pmatrix};\qquad
n_-\leftrightarrow X=\begin{pmatrix}a&b\\0&0\end{pmatrix},\quad
Z=\begin{pmatrix}0&c\\0&0\end{pmatrix}.
\]
The semisimple core action is $M\mapsto n_+Mn_+^{\prime t}$ with character $\psi(b+b')$. The nilpotent core action is $M\mapsto n_-Mn_+^{\prime t}$ with character $\psi(a+b')$.

Write initially
\[
M=\begin{pmatrix}
a_0&r_1&r_2&b_0\\
s_1&e&f&p\\
s_2&g&h_0&q\\
c_0&k&r&0
\end{pmatrix}.
\]
Since the pairs $(p,q)$ and $(k,r)$ are nonzero, the residual
actions changing $b_0$ by $ap+bq$ and $c_0$ by $a'k+b'r$ allow
us to set both entries to zero. With $b_0=c_0=0$, the symplectic equations give $ps_2-qs_1=\mu$. A left operation with parameters $(a,b)=(-qz,pz)$ keeps $b_0=0$ and changes $a_0$ by $\mu z$; hence that entry can be made zero. Equation~\eqref{m7:s:eq:R1-coreforms} then forces the central $2\times2$ block to vanish. We obtain
\begin{equation}\label{m7:s:eq:R1-core}
M=\begin{pmatrix}
0&r_1&r_2&0\\
s_1&0&0&p\\
s_2&0&0&q\\
0&k&r&0
\end{pmatrix},\qquad
ps_2-qs_1=\mu,
\end{equation}
with
\begin{equation}\label{m7:s:eq:R1-relations}
rr_1-kr_2=-\mu\quad\text{in the semisimple case},\qquad
rr_1-kr_2=\mu\quad\text{in the nilpotent case}.
\end{equation}
The remaining central parameters act by
\begin{equation}\label{m7:s:eq:R1-central}
(r_1,r_2)\mapsto(r_1+zk,r_2+zr),\qquad
(s_1,s_2)\mapsto(s_1+z'p,s_2+z'q),
\end{equation}
and have character one.

Solving the residual unipotent stabilizer equations at \eqref{m7:s:eq:R1-core} gives zero central parameters and
\begin{equation}\label{m7:s:eq:R1-stabilizer}
\begin{array}{c|c|c|c}
& (a,b)&(a',b')&\text{character}\\ \hline
\text{semisimple}&(-qz,pz)&(rz,-kz)&\psi((p-k)z)\\
\text{nilpotent}&(-qz,pz)&(-rz,kz)&\psi((k-q)z).
\end{array}
\end{equation}
The determinant relations \eqref{m7:s:eq:R1-core}--\eqref{m7:s:eq:R1-relations} show that no further parameter remains. The stabilizer of the outer row--column normalization is exactly the pair of residual groups in \eqref{m7:s:eq:R1-npm}, so this is the original $U\times U$ stabilizer. It has dimension one.

At these representatives,
\begin{equation}\label{m7:s:eq:R1-kappa-core}
\kappa=\delta(k-p)\quad\text{in the semisimple case},\qquad
\kappa=q-k\quad\text{in the nilpotent case}.
\end{equation}
It follows from \cref{cf:stabilizer} that
\begin{equation}\label{m7:s:eq:R1-incompatible}
\Dist_{\eta_K}(R_1\setminus\mathcal C_1)=0.
\end{equation}
All stabilizers used here are unipotent, so their modular characters are trivial.

\paragraph{The semisimple part of the compatible locus}

Compatibility is $k=p$. The scalar
\[
p=\delta^{-1}\tr(BD^\#)
\]
is intrinsic and $K$-invariant. If $p\ne0$, use \eqref{m7:s:eq:R1-central} to set $r_1=s_1=0$. This gives
\begin{equation}\label{m7:s:eq:R1-sem-plus}
M_{\mathrm s}^+(\mu;p,q,r)=
\begin{pmatrix}
0&0&\mu/p&0\\
0&0&0&p\\
\mu/p&0&0&q\\
0&p&r&0
\end{pmatrix},\qquad \mu p\ne0,
\end{equation}
with $q,r\in F$. For $a=\diag(a_1,a_2)$, put $t=a_2/a_1$. The remaining even and odd actions are
\begin{equation}\label{m7:s:eq:R1-plane-actions}
(q,r)\mapsto(tq,t^{-1}r),\qquad(q,r)\mapsto(r,q),
\end{equation}
while $\delta,\mu,p$ are fixed.

To justify this plane reduction, fix $D_{\mathrm s}$. the quantities $\mu,p,q,r$ are $U\times U$-invariants, since the first column of $B$ and first row of $C$ are unchanged. The normalizations above show that every fiber is one $U\times U$-orbit and provide rational sections using only $\delta,\mu,p$ as denominators. They apply also at $q=r=0$. The stabilizer is
\[
\bigl(n_+(-qz,pz,0),n_+(rz,-pz,0)\bigr).
\]
The torus action has determinant one on $\Lie(U\times U)$ and fixes the parameter $z$ of this stabilizer. Transpose exchanges the two factors and changes $z$ to $-z$, with absolute modulus one. Finally, the base of rank-one idempotents $D/\delta$ is $\GL_2/T$, with both groups unimodular. The character-weighted orbital Haar line therefore has trivial
residual torus and odd density characters. This gives the action
\eqref{m7:s:eq:R1-plane-actions}, with passive parameters
$(\delta,\mu,p)$, so \cref{m8:s:c:plane} proves vanishing on the
whole plane, including its two axes.

If $p=0$, then $k=0$ and $q,r\ne0$. The central operations give
\begin{equation}\label{m7:s:eq:R1-sem-zero}
M_{\mathrm s}^0(\mu;q,r)=
\begin{pmatrix}
0&-\mu/r&0&0\\
-\mu/q&0&0&0\\
0&0&0&q\\
0&0&r&0
\end{pmatrix}.
\end{equation}
Conjugation by $d(\diag(1,r/q))$ carries its embedding \eqref{m7:s:eq:R1-sembedding} to its transpose. This transporter has character one. \Cref{cf:stabilizer} eliminates this piece as well.

\paragraph{The nilpotent part: two explicit identities}

Here $D=E_{12}$ and compatibility is $k=q$. Let
\begin{equation}\label{m7:s:eq:R1-corrected-transpose}
w=\begin{pmatrix}0&1\\1&0\end{pmatrix},\qquad
\tau_{\mathrm n}(x)=d(w)x^td(w)^{-1}.
\end{equation}
This is the odd element $(1,1,w)\tau$, so its character is $-1$. On the core it is $M\mapsto PM^tP$, where $P=\diag(1,w,1)$.

When $q\ne0$, the central operations set $r_1=s_2=0$ and give
\begin{equation}\label{m7:s:eq:R1-nil-plus}
M_{\mathrm n}^+(\mu;p,q,r)=
\begin{pmatrix}
0&0&-\mu/q&0\\
-\mu/q&0&0&p\\
0&0&0&q\\
0&q&r&0
\end{pmatrix}.
\end{equation}
Write $x_{\mathrm n}^+=\mathscr E_{\mathrm n}(\mu;M_{\mathrm n}^+)$. Then
\[
\tau_{\mathrm n}(x_{\mathrm n}^+(\mu;p,q,r))=x_{\mathrm n}^+(\mu;r,q,p).
\]
For $h\in F$, set
\begin{equation}\label{m7:s:eq:R1-nil-transporter-matrices}
X_L(h)=\begin{pmatrix}h^2q&-h^2(hq+p)\\hq&-h(hq+p)\end{pmatrix},\qquad
X_R(h)=\begin{pmatrix}h(hq-r)&hq\\0&0\end{pmatrix}.
\end{equation}
Both determinants vanish. The element
\[
k_h=\bigl(u(X_L(h),0),u(-X_R(h),0),a_h\bigr)
\]
satisfies the $6\times6$ identity
\begin{equation}\label{m7:s:eq:R1-nil-transporter}
k_h\cdot x_{\mathrm n}^+(\mu;p,q,r)
=x_{\mathrm n}^+(\mu;p+hq,q,r-hq),
\qquad
\eta_K(k_h)=\psi\bigl(h(r-p-hq)\bigr).
\end{equation}
Indeed, $\tr X_L=-hp$ and $\tr X_R=h(hq-r)$, and the two matrix products in \eqref{m7:s:eq:R1-nil-transporter} agree entrywise using \eqref{m7:s:eq:R1-nembedding}. Taking $h=(r-p)/q$ simultaneously reaches $\tau_{\mathrm n}(x)$ and makes the character one. Composition with $(1,1,w)$ gives a character-one transporter to ordinary transpose.

When $q=0$, one has $k=0$ and $p,r\ne0$. The normal form is
\begin{equation}\label{m7:s:eq:R1-nil-zero}
M_{\mathrm n}^0(\mu;p,r)=
\begin{pmatrix}0&\mu/r&0&0\\0&0&0&p\\\mu/p&0&0&0\\0&0&r&0\end{pmatrix}.
\end{equation}
Let $x_{\mathrm n}^0=\mathscr E_{\mathrm n}(\mu;M_{\mathrm n}^0)$. For $h\in F$, define
\begin{align}
X_L&=\begin{pmatrix}0&-h^2p/2\\0&-hp\end{pmatrix},&
Z_L&=\begin{pmatrix}0&h\mu/r^2\\0&0\end{pmatrix},\nonumber\\
X_R&=\begin{pmatrix}-hr&0\\h^2r/2&0\end{pmatrix},&
Z_R&=\begin{pmatrix}0&0\\h\mu/p^2&0\end{pmatrix}.\label{m7:s:eq:R1-extra-matrices}
\end{align}
They satisfy $\det X_L=\det X_R=0$ and $\tr Z_L=\tr Z_R=0$. Direct multiplication gives
\begin{equation}\label{m7:s:eq:R1-extra-stab}
n(h)=\bigl(u(X_L,Z_L),u(-X_R,-Z_R),a_h\bigr),\qquad
n(h)\cdot x_{\mathrm n}^0=x_{\mathrm n}^0,
\end{equation}
\begin{equation}\label{m7:s:eq:R1-extra-character}
n(h+h')=n(h)n(h'),\qquad
\eta_K(n(h))=\psi\bigl(h(r-p)\bigr).
\end{equation}
For $p\ne r$, this gives a nontrivial character on a unipotent
stabilizer subgroup, whose modulus is one.

To express these pieces invariantly, observe that on the nilpotent branch,
\[
q=\tr(BD).
\]
On $q=0$, the nonzero scalars $p,r$ are characterized intrinsically by
\begin{equation}\label{m7:s:eq:R1-second-compatible}
BD=pD,\qquad DC=rD.
\end{equation}
The unipotent changes vanish in these products because $D^2=0$; the reductive component conjugates the equations. Hence $p=r$ is the closed invariant condition $BD=DC$. On that locus the representative \eqref{m7:s:eq:R1-nil-zero} is fixed by $\tau_{\mathrm n}$, so $(1,1,w)$ supplies the character-one transporter to ordinary transpose.

\begin{proof}[Proof of \cref{m7:s:prop:R1}]
Equation~\eqref{m7:s:eq:R1-incompatible} reduces support to $\mathcal C_1$. Its semisimple part is open, and is eliminated by \eqref{m7:s:eq:R1-sem-plus}--\eqref{m7:s:eq:R1-sem-zero}. In its closed nilpotent part, the open set $q\ne0$ is eliminated by \eqref{m7:s:eq:R1-nil-transporter}. The remaining support lies on $q=0$, where
\eqref{m7:s:eq:R1-extra-character} forces $p=r$. The representative
is then fixed by the corrected transpose, giving vanishing on
this final closed locus.

A character-one transporter transfers through any normalization. If $y=k_0\cdot x$ and $a_y\cdot y=y^t$ with $\eta_K(a_y)=1$, then
\begin{equation}\label{m7:s:eq:transport-normalization}
a_x=\tau_K(k_0)^{-1}a_yk_0
\end{equation}
carries $x$ to $x^t$ and has character one, since $\eta_K\circ\tau_K=\eta_K$. Thus the normal forms preserve the character phase. \Cref{cf:stabilizer} and \eqref{cf:support-sequence} finish the intrinsic vanishing. The last assertion follows from \cref{m7:s:prop:R0} because $R_1$ is closed in $\mathcal V$.
\end{proof}

\subsubsection{The stratum \texorpdfstring{$R_3$}{R3}}\label{m7:s:sec:R3}

\begin{proposition}\label{m7:s:prop:R3}
One has $\Dist_{\tetaK}(R_3)=0$.
\end{proposition}

We use the stabilizer and transporter method of \cite[\S7.11]{JLX}.

\paragraph{Geometry and the unipotent stabilizer}

The rank identities of \cref{m7:s:prop:ranks} give
\begin{equation}\label{m7:s:eq:R3-geometry}
R_3=\{x\in G_7^\circ:D=0,\ \rank B=\rank C=1\}.
\end{equation}
The symplectic equations give $B^tjB=0$ and the corresponding
row equation for $C$; invertibility gives the indicated ranks.
Geometrically, the plane $xV_-$ lies in $V_+\oplus V_0$, meets
$V_+$ in a line, and projects to a line in $V_0$. The rational parabolic normalization in \cref{m7:s:prop:ranks} sends it to $\langle e_1,e_4\rangle=w_3V_-$. Thus
\[
R_3=P(F)w_3P(F)^t,
\]
with the $w_3$ in \eqref{m7:s:eq:w234}. Its columns are $e_6,e_3,e_2,-e_5,-e_4,-e_1$, which verify the required symplectic pairings directly.

For $a,b,c,z\in F$, put
\begin{equation}\label{m7:s:eq:R3-stab-w}
\begin{split}
u_0(a,b,c,z)&=u\left(\begin{pmatrix}a&0\\b&0\end{pmatrix},\begin{pmatrix}c&z\\0&-c\end{pmatrix}\right),\\
v_0(a,b,c,z)&=u\left(\begin{pmatrix}c&0\\-b&0\end{pmatrix},\begin{pmatrix}-a&z\\0&a\end{pmatrix}\right).
\end{split}
\end{equation}
To solve $u w_3v^t=w_3$ completely, we take $(u,v)=(u_0,v_0)$.
Indeed, starting with $v^t=w_3^{-1}u^{-1}w_3$ for arbitrary
$u=u(X,Z)$, membership of $v$ in $U$ forces $X_{12}=X_{22}=Z_{21}=0$; then $\tr Z=0$, and the remaining entries give \eqref{m7:s:eq:R3-stab-w}. Conversely, these matrices satisfy the equality. Multiplication is additive in $(a,b,c,z)$, since the products
$X(X')^\#$ vanish, giving the stabilizer $\Ga^4$.

At
\[
x=m(A,E)w_3m(C_0,F_0)^t,
\]
the stabilizer is obtained by conjugating $u_0,v_0$ by these two Levi factors. With $\Pi_L=E^{-1}A$ and $\Pi_R=F_0^{-1}C_0$, its character is
\begin{equation}\label{m7:s:eq:R3-stab-character}
\psi\bigl((\Pi_L)_{11}a+((\Pi_L)_{12}-(\Pi_R)_{12})b+(\Pi_R)_{11}c\bigr).
\end{equation}
Therefore the complete compatibility conditions are
\begin{equation}\label{m7:s:eq:R3-compat-P}
(\Pi_L)_{11}=(\Pi_R)_{11}=0,\qquad
(\Pi_L)_{12}=(\Pi_R)_{12}.
\end{equation}
They have the intrinsic form
\begin{equation}\label{m7:s:eq:R3-compat}
\mathcal C_3=\{x\in R_3:B=C\}.
\end{equation}
To check this equivalence, direct multiplication gives
\[
E^{-1}BF_0^{-t}=-\frac{E_{21}\Pi_R^t}{\det\Pi_R},\qquad
E^{-1}CF_0^{-t}=-\frac{\Pi_LE_{12}}{\det\Pi_L}.
\]
Equality forces the two $(1,1)$ entries to be zero, and then invertibility forces equality of the two $(1,2)$ entries. Conversely those three equations give $B=C$. On $D=0$, the unipotent actions fix $B,C$, and the reductive component conjugates them simultaneously. Thus $\mathcal C_3$ is closed and stable under both parities.

The character in \eqref{m7:s:eq:R3-stab-character} is nontrivial somewhere on the unipotent stabilizer at every point outside $\mathcal C_3$. Hence
\begin{equation}\label{m7:s:eq:R3-incompatible}
\Dist_{\eta_K}(R_3\setminus\mathcal C_3)=0
\end{equation}
by \cref{cf:stabilizer}.

\paragraph{A compatible normal form}

Every point of $\mathcal C_3$ is $K$-equivalent to
\begin{equation}\label{m7:s:eq:R3-seven-parameters}
x_3(p,q,r,s;c,\ell_0,h)=
\begin{pmatrix}
0&0&0&0&(r+\ell_0p)/c&(s+\ell_0q)/c\\
0&0&r&s&0&0\\
-h(r+\ell_0p)&-h(s+\ell_0q)&0&0&p&q\\
r+\ell_0p&s+\ell_0q&0&0&0&0\\
0&0&p&q&0&0\\
-cp&-cq&0&0&0&0
\end{pmatrix},
\end{equation}
where
\begin{equation}\label{m7:s:eq:R3-parameters}
c\ne0,\qquad \mu=ps-qr\ne0.
\end{equation}
Here $B=C=\sm{p&q\\0&0}$ and $\ell(x_3)=\mu$. The parameter $\ell_0$ in \eqref{m7:s:eq:R3-seven-parameters} is a scalar coordinate, distinct from the multiplier function $\ell$.

For explicit coverage, write
\[
\Pi_L=\begin{pmatrix}0&b\\c_0&d_0\end{pmatrix},\qquad
\Pi_R=\begin{pmatrix}0&b\\f_0&g_0\end{pmatrix},\qquad bc_0f_0\ne0,
\]
and set $H_0=A^tF_0$. Conjugation by $d(A^{-1})$ gives \eqref{m7:s:eq:R3-seven-parameters} with
\begin{equation}\label{m7:s:eq:R3-param-from-Levi}
p=\frac{(H_0)_{12}}{bc_0},\quad q=\frac{(H_0)_{22}}{bc_0},\quad
r=(H_0)_{11},\quad s=(H_0)_{21},\quad
c=f_0,\quad\ell_0=c_0g_0,\quad h=d_0/c_0.
\end{equation}
Conversely, put
\[
P_0=\begin{pmatrix}0&1\\1&h\end{pmatrix},\quad
Q_0=\begin{pmatrix}0&1\\c&\ell_0\end{pmatrix},\quad
H_0=\begin{pmatrix}r&p\\s&q\end{pmatrix}.
\]
Then
\begin{equation}\label{m7:s:eq:R3-seven-factor}
x_3=m(I_2,P_0^{-1})w_3m(H_0Q_0,H_0)^t.
\end{equation}
This proves rational coverage and the symplectic equation for every value in \eqref{m7:s:eq:R3-parameters}.

\paragraph{The semisimple common block}

Suppose $p\ne0$. Set
\[
k=q/p,\qquad \sigma=s-kr=\mu/p,\qquad A_0=r+\ell_0p,
\qquad a_0=\begin{pmatrix}1&k\\0&1\end{pmatrix}.
\]
Define
\begin{align}
X_L&=-\frac{A_0}{cp}\begin{pmatrix}1&h-k\\0&0\end{pmatrix},&
Z_L&=-\frac rp\begin{pmatrix}k&-k^2\\1&-k\end{pmatrix},\nonumber\\
X_-&=\begin{pmatrix}-ck&0\\0&0\end{pmatrix},&
Z_-&=\begin{pmatrix}0&0\\-\sigma(h-k)/p&0\end{pmatrix}.
\label{m7:s:eq:R3-ss-normalizers}
\end{align}
Their $X$ determinants and $Z$ traces are zero. The identity is
\begin{equation}\label{m7:s:eq:R3-ss-normalize}
u(X_L,Z_L)d(a_0)x_3d(a_0)^{-1}
=s_3(p,\sigma;c)u(X_-,Z_-)^t,
\end{equation}
where
\begin{equation}\label{m7:s:eq:R3-ss}
s_3(p,\sigma;c)=
\begin{pmatrix}
0&0&0&0&0&\sigma/c\\
0&0&0&\sigma&0&0\\
0&0&0&0&p&0\\
0&\sigma&0&0&0&0\\
0&0&p&0&0&0\\
-cp&0&0&0&0&0
\end{pmatrix},\qquad p\sigma c\ne0.
\end{equation}
The same formulas apply at $q=0$. Further conjugation by $d(\diag(1,c^{-1}))$ makes $c=1$. The invariant parameters are exactly
\[
p=\tr B,\qquad \mu=p\sigma.
\]
Necessity is immediate, and \eqref{m7:s:eq:R3-ss-normalize} proves sufficiency over $F$.

For every representative \eqref{m7:s:eq:R3-ss}, take
\begin{equation}\label{m7:s:eq:R3-ss-transporter}
a_{\mathrm s}=\diag\left(1,-\frac\sigma{c^2p}\right).
\end{equation}
It is invertible and satisfies
\begin{equation}\label{m7:s:eq:R3-ss-transporter-identity}
d(a_{\mathrm s})s_3(p,\sigma;c)d(a_{\mathrm s})^{-1}
=s_3(p,\sigma;c)^t.
\end{equation}
To check the identity, note that the paired entries $(2,4),(4,2)$
and $(3,5),(5,3)$ are unchanged, while the diagonal ratio in
\eqref{m7:s:eq:R3-ss-transporter} exchanges the $(1,6)$ and $(6,1)$
entries $\sigma/c$ and $-cp$. This purely reductive transporter
has character one, so \cref{cf:stabilizer} and
\eqref{m7:s:eq:transport-normalization} eliminate the semisimple
part of $\mathcal C_3$.

\paragraph{The nilpotent common block}

Now $p=0$, so $qrc\ne0$. Conjugation by $d(\diag(1,q))$ gives
\[
x_3(0,q,r,s;c,\ell_0,h)\longmapsto
x_3(0,1,\rho,s;c_1,\ell_1,h_1),
\]
where $\rho=qr$, $c_1=cq$, $\ell_1=q\ell_0$, and $h_1=h/q$. Relabel $c_1$ as $c$. With
\[
X_L=-\frac{s+\ell_1}{c}\begin{pmatrix}1&h_1\\0&0\end{pmatrix},\quad
Z_L=\begin{pmatrix}0&0\\-s&0\end{pmatrix},\quad
Z_-=\begin{pmatrix}0&-\rho h_1\\0&0\end{pmatrix},
\]
one has
\begin{equation}\label{m7:s:eq:R3-nil-normalize}
u(X_L,Z_L)x_3(0,1,\rho,s;c,\ell_1,h_1)
=n_3(\rho,c)u(0,Z_-)^t,
\end{equation}
where
\begin{equation}\label{m7:s:eq:R3-nil}
n_3(\rho,c)=
\begin{pmatrix}
0&0&0&0&\rho/c&0\\
0&0&\rho&0&0&0\\
0&0&0&0&0&1\\
\rho&0&0&0&0&0\\
0&0&0&1&0&0\\
0&-c&0&0&0&0
\end{pmatrix},\qquad \rho c\ne0,\qquad \mu=-\rho.
\end{equation}
The unipotent matrices used in \eqref{m7:s:eq:R3-nil-normalize} satisfy \eqref{m7:s:eq:U}.

On this branch, $c$ is an intrinsic invariant, characterized by
\begin{equation}\label{m7:s:eq:R3-intrinsic-c}
Bx_{31}=-cB,\qquad x_{13}B=\frac\rho c B.
\end{equation}
The unipotent changes $x_{13}\mapsto x_{13}+XB$ and $x_{31}\mapsto x_{31}+BX_v^t$ do not change these products because $B^2=0$. The reductive action conjugates the equations. Hence $(\mu,c)$ is an invariant pair, and \eqref{m7:s:eq:R3-nil-normalize} shows that it classifies the nilpotent compatible $K(F)$-orbits.

The restriction of $\eta_K$ to the $U\times U$ stabilizer is trivial here, but the centralizer of $B=E_{12}$ supplies a further direction. Let $a_z=\sm{1&z\\0&1}$. Define
\begin{align}
X_L(z)&=\begin{pmatrix}\rho z/c&-\rho z^2/(2c)\\0&0\end{pmatrix},&
Z_L(z)&=\begin{pmatrix}
\rho z^2/2&-\rho z^3(\rho+2c^2)/(6c^2)\\
\rho z&-\rho z^2/2
\end{pmatrix},\label{m7:s:eq:R3-nil-stab-L}\\
X_R(z)&=\begin{pmatrix}0&0\\cz^2/2&cz\end{pmatrix},&
Z_R(z)&=\begin{pmatrix}
-\rho z^2/2&-\rho z\\
z^3(2\rho+c^2)/6&\rho z^2/2
\end{pmatrix}.\label{m7:s:eq:R3-nil-stab-R}
\end{align}
The $X$ determinants and $Z$ traces are zero. Put
\[
\kappa_z=(u(X_L(z),Z_L(z)),u(X_R(z),Z_R(z)),a_z).
\]
Then
\begin{equation}\label{m7:s:eq:R3-nil-stab-identity}
\kappa_z\cdot n_3(\rho,c)=n_3(\rho,c),\qquad
\kappa_{z+z'}=\kappa_z\kappa_{z'}.
\end{equation}
Both identities follow algebraically by setting
\begin{align*}
N_L&=E_{12}+E_{34}+E_{56}+\frac\rho c(E_{13}-E_{46})+\rho E_{25},\\
N_R&=-E_{21}-E_{43}-E_{65}+c(E_{24}-E_{35})-\rho E_{16}.
\end{align*}
One has $N_L^4=N_R^4=0$ and
\[
\begin{aligned}
N_L n_3+n_3N_R^t&=0,\\
u(X_L(z),Z_L(z))d(a_z)&=\exp(zN_L),\\
u(X_R(z),Z_R(z))d(a_z^{-t})&=\exp(zN_R).
\end{aligned}
\]
The exponentials terminate at degree three, proving \eqref{m7:s:eq:R3-nil-stab-identity} and the asserted membership in the semidirect product.

The character is
\begin{equation}\label{m7:s:eq:R3-nil-extra-character}
\eta_K(\kappa_z)=\psi\left(z\frac{\rho+c^2}{c}\right).
\end{equation}
Every even stabilizer element has reductive component $\lambda a_z$, because it centralizes $B=E_{12}$. The subgroup just constructed lifts every $a_z$; after removing it and the ineffective scalar factor, one is in the four-dimensional unipotent stabilizer \eqref{m7:s:eq:R3-stab-w}. Thus
\[
K_{n_3(\rho,c)}\simeq(\Ga^4\rtimes\Ga)\times\Gm.
\]
This stabilizer is unimodular. Its character is trivial on the four-dimensional kernel and on scalars, so \eqref{m7:s:eq:R3-nil-extra-character} accounts for all of its possible nontriviality. A nontrivial character on the subgroup $\kappa_z$ eliminates every point with $\rho+c^2\ne0$.

The remaining compatibility condition is
\begin{equation}\label{m7:s:eq:R3-extra-compatible}
\rho=-c^2\quad\Longleftrightarrow\quad\mu=c^2
\quad\Longleftrightarrow\quad x_{13}B=Bx_{31}.
\end{equation}
The scalar $c$ is recovered regularly on every chart where an entry of $B$ is nonzero, so this is a closed invariant condition in the nilpotent branch. The scalar $c$ is part of the existing orbit data.

Let $w=\sm{0&1\\1&0}$. The identity
\begin{equation}\label{m7:s:eq:R3-corrected-transpose}
d(w)n_3(\rho,c)^td(w)^{-1}=n_3(\rho,-\rho/c)
\end{equation}
shows that transpose exchanges the two orbit parameters $c$ and $-\rho/c$. On \eqref{m7:s:eq:R3-extra-compatible}, they agree, and
\begin{equation}\label{m7:s:eq:R3-nil-transporter}
d(w)n_3(-c^2,c)d(w)^{-1}=n_3(-c^2,c)^t.
\end{equation}
This transporter has character one. Outside this locus the ordered invariant $c$ changes, but those orbits have already been eliminated by \eqref{m7:s:eq:R3-nil-extra-character}.

\begin{proof}[Proof of \cref{m7:s:prop:R3}]
Equation~\eqref{m7:s:eq:R3-incompatible} gives support in the closed locus $\mathcal C_3$. Its semisimple part is open and is eliminated by \eqref{m7:s:eq:R3-ss-transporter-identity}. In the closed nilpotent part, \eqref{m7:s:eq:R3-nil-extra-character} gives support on \eqref{m7:s:eq:R3-extra-compatible}, whose representatives have the character-one transporter \eqref{m7:s:eq:R3-nil-transporter}. Transport through the normalizations is justified by \eqref{m7:s:eq:transport-normalization}. \Cref{cf:stabilizer} and \eqref{cf:support-sequence} finish the proof. All denominators in the normalizations are explicitly nonzero. The fixed denominators $2$ and $6$ in the one-parameter subgroup are allowed in every residue characteristic because $F$ has characteristic zero.
\end{proof}

\subsubsection{The final stratum \texorpdfstring{$R_4$}{R4}}\label{m7:s:sec:R4}

\begin{proposition}\label{m7:s:prop:R4}
One has $\Dist_{\tetaK}(R_4)=0$.
\end{proposition}

We prove intrinsic vanishing by identifying the compatible locus
with a family of ordinary $\GL_2$ conjugacy classes.

\paragraph{Geometry and a unique rational factorization}

The rank pattern forces $B=C=D=0$ and $x_{22}$ invertible. Conversely, inside $G_7^\circ$ the three zero conditions already force the rank pattern. Thus
\begin{equation}\label{m7:s:eq:R4-geometry}
R_4=\{x\in G_7^\circ:x_{23}=x_{32}=x_{33}=0\}.
\end{equation}
Write
\[
x=\begin{pmatrix}H&B_0&A\\D_0&E&0\\C_0&0&0\end{pmatrix}.
\]
Invertibility makes $A,E,C_0$ invertible. Since $w_4=w_4^t=w_4^{-1}\in\Sp_6(F)$ and
\[
xw_4=\begin{pmatrix}A&B_0&H\\0&E&D_0\\0&0&C_0\end{pmatrix},
\]
one has
\begin{equation}\label{m7:s:eq:R4-parabolic}
R_4=P(F)w_4=P(F)w_4P(F)^t.
\end{equation}

Let $\lambda=\det E/\det A$. The symplectic equations give
\[
\ell(x)=\det E,\qquad C_0=\lambda A.
\]
Putting $X=B_0E^{-1}$ and $Z=H(\lambda A)^{-1}$, the remaining equations give
\[
D_0=-X^\#\lambda A,\qquad \tr Z=-\det X.
\]
Consequently the factorization
\begin{equation}\label{m7:s:eq:R4-factor}
x=u(X,Z)m(A,E)w_4
=\begin{pmatrix}
Z\lambda A&XE&A\\
-X^\#\lambda A&E&0\\
\lambda A&0&0
\end{pmatrix}
\end{equation}
is unique. Its coordinates are recovered by
\[
A=x_{13},\quad E=x_{22},\quad X=x_{12}E^{-1},\quad
Z=x_{11}(\lambda A)^{-1}.
\]
Thus $U\times\GL_2\times\GL_2\to R_4$ is an algebraic isomorphism over $F$. This proves coverage of all $F$-points and multiplier values.

\paragraph{The stabilizer and the compatibility equation}

Define an automorphism of $U$ by
\begin{equation}\label{m7:s:eq:R4-phi}
\phi(u)=w_4u^{-t}w_4,
\qquad
\phi(u(X,Z))=u((X^\#)^t,(Z^\#)^t).
\end{equation}
Here $Z^\#=-Z-\det(X)I_2$ by \eqref{m7:s:eq:U}. The stabilizer at $w_4$ is the graph
\begin{equation}\label{m7:s:eq:R4-stab-w}
(U\times U)_{w_4}=\{(u,\phi(u)):u\in U\}.
\end{equation}
Indeed $uw_4v^t=w_4$ uniquely determines $v=\phi(u)$. A block check is
\[
u(X,Z)w_4=\begin{pmatrix}Z&X&I_2\\-X^\#&I_2&0\\I_2&0&0\end{pmatrix},\qquad
\phi(u)^t=\begin{pmatrix}I_2&0&0\\X^\#&I_2&0\\Z^\#&-X&I_2\end{pmatrix};
\]
the product is $w_4$, since $Z+XX^\#+Z^\#=0$.

Write $y(A,E)=m(A,E)w_4$ and $x=u_*y(A,E)$. The complete stabilizer is
\begin{equation}\label{m7:s:eq:R4-stab-general}
(U\times U)_x=
\{(u_*m(A,E)u m(A,E)^{-1}u_*^{-1},\phi(u)):u\in U\}.
\end{equation}
Conjugation by $m(A,E)$ sends
\[
(X,Z)\longmapsto(AXE^{-1},\lambda^{-1}AZA^{-1}).
\]
With $\Pi=E^{-1}A$, its character is therefore
\begin{equation}\label{m7:s:eq:R4-character}
\psi\bigl(\tr(\Pi X)+\tr((X^\#)^t)\bigr)
=\psi\bigl(\tr((\Pi+I_2)X)\bigr).
\end{equation}
As $X$ ranges over all $2\times2$ matrices, this is identically one exactly when $A=-E$. Thus
\begin{equation}\label{m7:s:eq:R4-compatible}
\mathcal C_4=\{x\in R_4:x_{13}=-x_{22}\}.
\end{equation}
It also forces $x_{31}=-x_{22}$, since $\lambda=1$. Both unipotent actions fix $x_{13},x_{22}$ on $R_4$, and the reductive component conjugates them. Consequently $\mathcal C_4$ is closed and invariant. It is transpose-stable because both outer blocks equal $-x_{22}$.

If $\Pi+I_2\ne0$, choose a matrix unit $X_0$ with $\tr((\Pi+I_2)X_0)\ne0$. The subgroup $u(tX_0,0)$ is defined for all $t\in F$, and \eqref{m7:s:eq:R4-stab-general} gives a nontrivial stabilizer character. These stabilizers are unipotent and have trivial modulus. Hence
\begin{equation}\label{m7:s:eq:R4-incompatible}
\Dist_{\eta_K}(R_4\setminus\mathcal C_4)=0.
\end{equation}
Notice the sign in the compatibility equation: although $w_4=y(I_2,I_2)$ is symmetric, its stabilizer character is $\psi(2\tr X)$, so it is not itself character-compatible for Model~7.

\paragraph{Compatible orbits and their stabilizers}

For $E\in\GL_2(F)$, put
\begin{equation}\label{m7:s:eq:R4-yE}
y_E=y(-E,E)=\begin{pmatrix}0&0&-E\\0&E&0\\-E&0&0\end{pmatrix}.
\end{equation}
Every compatible point is uniquely $u_*y_E$. The unipotent groups fix $E=x_{22}$, and the left $U$ action is transitive on its fiber. The remaining action is
\[
d(a)y_Ed(a)^{-1}=y_{aEa^{-1}}.
\]
Therefore
\begin{equation}\label{m7:s:eq:R4-orbits}
K(F)\backslash\mathcal C_4\simeq
\{\GL_2(F)\text{-conjugacy classes in }\GL_2(F)\}.
\end{equation}
These are the scalar classes, the split classes with two distinct nonzero eigenvalues, the nontrivial Jordan blocks with nonzero eigenvalue, and the classes with irreducible quadratic characteristic polynomial. Explicit representatives are
\begin{equation}\label{m7:s:eq:R4-conjugacy-types}
\lambda I_2,\quad
\diag(\lambda_1,\lambda_2),\quad
\begin{pmatrix}\lambda&1\\0&\lambda\end{pmatrix},\quad
\begin{pmatrix}0&-\delta\\1&s\end{pmatrix},
\end{equation}
with the evident nonzero and irreducibility conditions. Completeness over $F$ follows by taking a cyclic basis $(v,Ev)$ for any nonscalar $E$ and applying Cayley--Hamilton. The characteristic polynomial $z^2-sz+\delta$ is part of the complete orbit invariant.

Let $m_E=m(-E,E)$ and $C_E=Z_{\GL_2(F)}(E)$. The even stabilizer is
\begin{equation}\label{m7:s:eq:R4-full-stab}
K_{y_E}=\{(m_Eu m_E^{-1},\phi(u),a):u\in U,\ a\in C_E\}
\simeq U\rtimes C_E.
\end{equation}
Solving the middle-block equation forces $a\in C_E$; the remaining
equations are \eqref{m7:s:eq:R4-stab-general}. The stabilizer
character is identically one, because its two traces are
$-\tr X$ and $\tr X$. For scalar $E$, $C_E=\GL_2(F)$; for nonscalar $E$, $C_E=F[E]^\times$ is commutative. Both groups are unimodular, and their conjugation action on
the independent coordinates $X\in\Mat_2$ and
$Z_0\in\mathfrak{sl}_2$ has determinant one. Hence the full
stabilizer \eqref{m7:s:eq:R4-full-stab} is unimodular, with
trivial orbital density character.

\paragraph{Explicit transpose transport}

For $E=\sm{r&s\\t&u}$, choose
\begin{equation}\label{m7:s:eq:R4-aE}
a_E=\begin{cases}
\begin{pmatrix}r-u&s\\s&0\end{pmatrix},&s\ne0,\\[3mm]
\begin{pmatrix}0&t\\t&u-r\end{pmatrix},&s=0,\ t\ne0,\\[3mm]
I_2,&s=t=0.
\end{cases}
\end{equation}
Its determinant is respectively $-s^2$, $-t^2$, or $1$. Each case satisfies $a_EE=E^ta_E$. The $6\times6$ equality is
\begin{equation}\label{m7:s:eq:R4-y-transporter}
d(a_E)y_Ed(a_E)^{-1}=y_E^t,
\end{equation}
since its three nonzero blocks are $-a_EEa_E^{-1}$, $a_EEa_E^{-1}$, and $-a_EEa_E^{-1}$. This transporter is purely reductive and has character one. For scalar and split-diagonal representatives, $y_E$ itself is symmetric and the identity transporter suffices.

For an arbitrary point $x=u_*y_E$, define
\begin{equation}\label{m7:s:eq:R4-general-transporter}
k_x=(d(a_E)u_*^{-1}d(a_E)^{-1},\ u_*,\ a_E).
\end{equation}
Both unipotent components belong to $U$. Cancellation gives
\[
k_x\cdot x=d(a_E)y_Ed(a_E)^{-1}u_*^t=y_E^tu_*^t=x^t.
\]
If $u_*=u(X_*,Z_*)$, the first component has $X$ coordinate $-a_EX_*a_E^{-1}$ and $Z$ coordinate $a_E(-Z_*-\det(X_*)I_2)a_E^{-1}$. It satisfies \eqref{m7:s:eq:U}, and the sum of the two $X$ traces is zero. Hence $\eta_K(k_x)=1$.

\begin{proof}[Proof of \cref{m7:s:prop:R4}]
By \eqref{m7:s:eq:R4-incompatible}, the distribution is supported
on $\mathcal C_4$. On each compatible orbit,
\eqref{m7:s:eq:R4-general-transporter} supplies a character-one
transpose transporter. Applying \cref{cf:stabilizer} gives
$\Dist_{\tetaK}(\mathcal C_4)=0$, and the open--closed sequence
\eqref{cf:support-sequence} proves the proposition.
\end{proof}

\subsection{The residual manifold and its ambient conormal}\label{m7:residual}
We now prove the conditional ambient exclusion on $R_{7,2}$
used in the support induction.

The tangent convention used throughout this calculation is
\begin{equation}\label{m7:s:eq:tangent-convention}
 \dot x=xA,\qquad A=x^{-1}\dot x\in\mathfrak{gsp}_6.
\end{equation}
If a matrix $Y$ represents the ambient cotangent vector
$\dot x\mapsto\tr(Y\dot x)$, its identity translate is
\begin{equation}\label{m7:s:eq:cotangent-translation}Z=Yx.\end{equation}
The other tangent trivialization gives $xY=x(Yx)x^{-1}$, with the same
nilpotence properties. The trace form includes the multiplier direction.
\subsubsection{The residual-wavefront argument on \texorpdfstring{$R_2$}{R2}}\label{m7:s:sec:R2}

\begin{proposition}\label{m7:s:prop:R2}
Let $S\in\Dist_{\tetaK}(G_7^\circ)\cap\NW(G_7^\circ)$. Then
\begin{equation}\label{m7:s:eq:R2-spectral-result}
S|_{\mathcal W}=0,\qquad \mathcal W=R_0\sqcup R_1\sqcup R_2.
\end{equation}
\end{proposition}

We first reduce the support by the method of \cite[\S7.10]{JLX},
then apply \cref{cf:conormal} to the residual family. This second
step is where the spectral hypothesis enters the stratum calculations.

\paragraph{The unipotent compatibility equations}

By \cref{m7:s:prop:ranks}, every point of $R_2$ is $U\times U$-equivalent to
\begin{equation}\label{m7:s:eq:R2-Levi}
x=m(A,E)w_2m(C_0,F_0)^t.
\end{equation}
Define
\[
\mathsf u_2(a,b,c)=u\left(\begin{pmatrix}0&0\\a&b\end{pmatrix},
\begin{pmatrix}0&0\\c&0\end{pmatrix}\right).
\]
Solving $u w_2 v^t=w_2$ gives the stabilizer
\begin{equation}\label{m7:s:eq:R2-wstab}
(U\times U)_{w_2}
=\{(\mathsf u_2(-d,c,r),\mathsf u_2(c,d,r)):c,d,r\in F\}.
\end{equation}
For example, one may solve for $v^t=w_2^{-1}u^{-1}w_2$ in the coordinates \eqref{m7:s:eq:U}; the equations forcing $v$ to lie in $U$ leave exactly the three parameters in \eqref{m7:s:eq:R2-wstab}. Substitution verifies the converse.

Put
\[
\Pi_L=E^{-1}A,\qquad \Pi_R=F_0^{-1}C_0.
\]
For a general point, conjugate the first and second coordinates of
\eqref{m7:s:eq:R2-wstab} by $m(A,E)$ and $m(C_0,F_0)$, respectively.
The resulting character is
\begin{equation}\label{m7:s:eq:R2-stab-character}
\psi\bigl(((\Pi_L)_{22}+(\Pi_R)_{12})c
+((\Pi_R)_{22}-(\Pi_L)_{12})d\bigr).
\end{equation}
Thus the complete compatibility conditions are
\begin{equation}\label{m7:s:eq:R2-compatible}
(\Pi_L)_{22}+(\Pi_R)_{12}=0,\qquad
(\Pi_R)_{22}-(\Pi_L)_{12}=0.
\end{equation}
The central parameter $r$ contributes no character.

Let $\mathcal C_2$ be the locus where the restriction of $\eta_K$ to the $U\times U$ stabilizer is trivial. It is closed in $R_2$. Indeed, the stabilizers form a constant-rank unipotent stabilizer bundle on the parabolic cell, and the restriction of the additive trace functional to its Lie algebra is a regular section of its dual bundle. Its zero set is given by \eqref{m7:s:eq:R2-compatible}.
Equivalently, the adapted-basis charts for \eqref{m7:s:eq:R2-Levi}
express compatibility by these two closed equations, with invertible
denominators and independently of the chosen representative. \Cref{cf:stabilizer} gives
\begin{equation}\label{m7:s:eq:R2-incompatible}
\Dist_{\eta_K}(R_2\setminus\mathcal C_2)=0.
\end{equation}
Both the acting unipotent group and its stabilizers have trivial modulus.

\paragraph{Normal forms on the compatible locus}

Write the matrices in \eqref{m7:s:eq:R2-compatible} as
\[
\Pi_L=\begin{pmatrix}a&b\\c&d\end{pmatrix},\qquad
\Pi_R=\begin{pmatrix}e&-d\\f&b\end{pmatrix}.
\]
Set $s=ae+cf$, $t=be+df=\det\Pi_R\ne0$, and
\[
M=\Pi_L^{-1}F_0^tE\Pi_L=\begin{pmatrix}\alpha&\beta\\\gamma&\delta\end{pmatrix}.
\]
Conjugation by $d(\Pi_L^{-1}E^{-1})$ carries \eqref{m7:s:eq:R2-Levi} to
\begin{equation}\label{m7:s:eq:R2-canonical}
x_2(s,t,M)=
\begin{pmatrix}
\alpha s+\gamma t&\beta s+\delta t&0&0&0&0\\
0&0&0&0&-\alpha s/t-\gamma&-\beta s/t-\delta\\
0&0&\alpha&\beta&0&0\\
0&0&\gamma&\delta&0&0\\
-\alpha&-\beta&0&0&0&0\\
0&0&0&0&-\alpha/t&-\beta/t
\end{pmatrix}.
\end{equation}
Conversely, every $s\in F$, $t\in F^\times$, and $M\in\GL_2(F)$ occurs. One factorization is
\begin{equation}\label{m7:s:eq:R2-canonical-factor}
x_2(s,t,M)=w_2m(M^tQ_0,M^t)^t,\qquad
Q_0=\begin{pmatrix}s&-1\\t&0\end{pmatrix}.
\end{equation}
The symplectic equation is
\[
x_2^tJx_2=(\det M)J.
\]
At this representative,
\begin{equation}\label{m7:s:eq:R2-D}
D=\begin{pmatrix}0&0\\-\alpha/t&-\beta/t\end{pmatrix}.
\end{equation}
Thus $D$ is semisimple rank one when $\beta\ne0$, and nonzero nilpotent when $\beta=0$.

Suppose first that $\beta\ne0$. Put
\begin{equation}\label{m7:s:eq:R2-sem-a}
a_M=\begin{pmatrix}
(\alpha^2-\alpha\delta+\beta\gamma)/\beta&\alpha\\
\alpha&\beta
\end{pmatrix},\qquad \det a_M=-\det M\ne0,
\end{equation}
and $\Lambda=t(\alpha-\delta)-\beta s$. Let
\begin{equation}\label{m7:s:eq:R2-sem-Z}
Z_L=\begin{pmatrix}-\Lambda/\beta&2\alpha\Lambda/\beta^2\\0&\Lambda/\beta\end{pmatrix},\qquad
Z_R=\begin{pmatrix}\Lambda/\beta&0\\-2\alpha\Lambda/\beta^2&-\Lambda/\beta\end{pmatrix}.
\end{equation}
Both are traceless. Direct block multiplication gives
\begin{equation}\label{m7:s:eq:R2-sem-transporter}
u(0,Z_L)d(a_M)x_2(s,t,M)d(a_M)^{-1}u(0,Z_R)^t=x_2(s,t,M)^t.
\end{equation}
The two middle identities are $a_MM=M^ta_M$ and $a_MD=D^ta_M$; the outer blocks agree after the corrections \eqref{m7:s:eq:R2-sem-Z}. Both unipotent characters are one. Hence the semisimple part of $\mathcal C_2$ vanishes intrinsically by \cref{cf:stabilizer}.

When $\beta=0$, write $M=\sm{\alpha&0\\\gamma&\delta}$, with $\alpha\delta\ne0$. If $\alpha=\delta$, the matrix
\begin{equation}\label{m7:s:eq:R2-equal-a}
a_0=\begin{pmatrix}\gamma+\alpha s/t&\alpha\\\alpha&0\end{pmatrix}
\end{equation}
has determinant $-\alpha^2$ and satisfies
\begin{equation}\label{m7:s:eq:R2-equal-transporter}
d(a_0)x_2\left(s,t,\begin{pmatrix}\alpha&0\\\gamma&\alpha\end{pmatrix}\right)d(a_0)^{-1}
=x_2\left(s,t,\begin{pmatrix}\alpha&0\\\gamma&\alpha\end{pmatrix}\right)^t.
\end{equation}
Again the character is one.

The remaining family is
\begin{equation}\label{m7:s:eq:R2-residual}
\Rres=K\cdot\left\{x_2\left(s,t,\begin{pmatrix}\alpha&0\\\gamma&\delta\end{pmatrix}\right):
t\alpha\delta\ne0,\ \alpha\ne\delta\right\}.
\end{equation}
For these points, transpose does not lie in the even orbit. In fact, a transporter between the cross-free representatives would have to satisfy $aD=D^ta$ and $aM=M^ta$. The first condition forces
\[
a=\begin{pmatrix}r&q\\q&0\end{pmatrix},\qquad q\ne0,
\]
and the second forces $q\alpha=q\delta$, a contradiction. The unipotent corrections do not change this middle equation: the equations keeping the $(2,3)$ and $(3,2)$ blocks zero annihilate their possible contribution to the middle block. The even stabilizer has scalar reductive projection and trivial character. The remaining support is excluded using the nilpotent wavefront of the original global distribution.

\paragraph{A closed enlargement in an ambient open}

On $R_2$, the rank conditions give
\[
\ker D\subset\ker B,\qquad \im C\subset\im D.
\]
Choose a generalized inverse $R$ of $D$, with $DRD=D$, and define
\begin{equation}\label{m7:s:eq:R2-Schur}
\mathscr M(x)=x_{22}-BRC.
\end{equation}
On a chart $D_{ij}\ne0$, one may take $R=E_{ji}/D_{ij}$. Since $B=BRD$ and $C=DRC$, the expression $BRC$ is independent of this choice. Thus \eqref{m7:s:eq:R2-Schur} is regular on the rank-one stratum.

It is invariant under $U\times U$ and transforms by conjugation under $d(a)$. For example, the left action changes $B$ to $B-X^\#D$ and $x_{22}$ to $x_{22}-X^\#C$, and these changes cancel in \eqref{m7:s:eq:R2-Schur}. Transpose takes $\mathscr M$ to $\mathscr M^t$. Removing $B,C$ in the symplectic equation gives
\begin{equation}\label{m7:s:eq:R2-Mdet}
\det\mathscr M(x)=\ell(x)\ne0.
\end{equation}

Define $\mathcal N_2\subset R_2$ by $D^2=0$ and the equations
\begin{equation}\label{m7:s:eq:R2-N2}
D\mathscr M=\alpha D,\qquad \mathscr M D=\delta D,\qquad
Dx_{31}=-\alpha D,\qquad x_{13}D=-\delta D
\end{equation}
for scalars $\alpha,\delta$. They are unique because $D\ne0$, and are regular on the nonzero-entry charts of $D$. The scalar-multiple conditions are closed, and the last two equations can be written $D(x_{31}+\mathscr M)=0$ and $(x_{13}+\mathscr M)D=0$. Therefore $\mathcal N_2$ is closed in $R_2$, hence in $\mathcal W$.

The canonical form \eqref{m7:s:eq:R2-canonical} identifies $\mathcal N_2$ with the nilpotent part of $\mathcal C_2$. Equations~\eqref{m7:s:eq:R2-N2} retain precisely its compatibility conditions. The normalization below gives the converse directly on the unequal-eigenvalue part. One has $\alpha\delta=\ell(x)$, and
\begin{equation}\label{m7:s:eq:R2-E2}
\mathcal E_2=\{x\in\mathcal N_2:\alpha=\delta\},\qquad
\Rres=\mathcal N_2\setminus\mathcal E_2.
\end{equation}
Both $\alpha$ and $\delta$ are $K$-invariant, and transpose exchanges them. The set $\mathcal E_2$ is closed in $\mathcal W$, and its intrinsic equivariant distributions vanish by \eqref{m7:s:eq:R2-equal-transporter}.

For the global distribution $S$ in \cref{m7:s:prop:R2}, \cref{m7:s:prop:R0,m7:s:prop:R1} first give support in $R_2$. Equation~\eqref{m7:s:eq:R2-incompatible} and the semisimple transporter \eqref{m7:s:eq:R2-sem-transporter} then give
\begin{equation}\label{m7:s:eq:R2-support-N2}
\supp(S|_{\mathcal W})\subset\mathcal N_2.
\end{equation}
To make the residual family closed in its ambient open, put
\begin{equation}\label{m7:s:eq:R2-Wcirc}
\mathcal W^\circ=\mathcal W\setminus\mathcal E_2.
\end{equation}
Then $\Rres$ is closed in $\mathcal W^\circ$, and
\begin{equation}\label{m7:s:eq:R2-support-res}
\supp(S|_{\mathcal W^\circ})\subset\Rres,\qquad
\WF_x(S|_{\mathcal W^\circ})\subset\Nil.
\end{equation}
The wavefront hypothesis is preserved under these open restrictions of $S$.

\paragraph{Smooth structure of the residual family}\label{m7:s:subsec:R2-smooth}

\begin{lemma}\label{m7:s:lem:R2-smooth}
The set $\Rres$ is a smooth locally closed submanifold of $G_7^\circ$ of dimension $16$.
\end{lemma}
\begin{proof}
On \eqref{m7:s:eq:R2-N2} with $\alpha\ne\delta$, choose a nonzero $\alpha$-eigenvector $v$ of $\mathscr M$ and put $w=-Dv/\alpha$. The vectors $v,w$ form a basis, determined up to a common scalar, in which
\begin{equation}\label{m7:s:eq:R2-pair-normal}
D=-\alpha E_{21},\qquad \mathscr M=\diag(\alpha,\delta).
\end{equation}
Using the eigenprojector $(\mathscr M-\delta I_2)/(\alpha-\delta)$
gives rational choices on charts. Thus the space of pairs
$(D,\mathscr M)$ is smooth of dimension
\[
2+\dim(\GL_2/\Gm)=5.
\]
The two scalar coordinates are $\alpha,\delta$.

In the basis \eqref{m7:s:eq:R2-pair-normal}, write
\[
B=\begin{pmatrix}b_1&0\\b_2&0\end{pmatrix},\qquad
C=\begin{pmatrix}0&0\\c_1&c_2\end{pmatrix}.
\]
The matrices
\begin{equation}\label{m7:s:eq:R2-cross-eliminate}
X_B=\begin{pmatrix}-b_2/\alpha&b_1/\alpha\\0&0\end{pmatrix},\qquad
X_C=\begin{pmatrix}0&0\\c_2/\alpha&-c_1/\alpha\end{pmatrix}
\end{equation}
have determinant zero and satisfy $X_B^\#D=B$, $D(X_C^\#)^t=C$. Multiplication by $u(X_B,0)$ on the left and $u(X_C,0)^t$ on the right kills $B,C$. The symplectic equation now forces
\[
x_{12}=\begin{pmatrix}0&0\\l_1&l_2\end{pmatrix}.
\]
Left multiplication with $X=-x_{12}\mathscr M^{-1}$ and $Z=0$ removes $x_{12}$ without changing the preceding zero blocks, and the remaining cross equation gives $x_{21}=0$. The six coordinates $b_1,b_2,c_1,c_2,l_1,l_2$ are recovered uniquely by reversing these operations.

The resulting cross-free matrix has the form
\begin{equation}\label{m7:s:eq:R2-cross-free}
x=\begin{pmatrix}
A&0&P_0\\0&\diag(\alpha,\delta)&0\\Q_0&0&-\alpha E_{21}
\end{pmatrix},\quad
P_0=\begin{pmatrix}p_1&0\\p_3&-\delta\end{pmatrix},\quad
Q_0=\begin{pmatrix}-\alpha&0\\q_3&q_4\end{pmatrix},
\end{equation}
where the outer symplectic equations give
\begin{equation}\label{m7:s:eq:R2-cross-free-A}
A=\begin{pmatrix}
-p_3-q_3p_1/\alpha&\delta-q_4p_1/\alpha\\
a_{21}&-(\delta q_3+p_3q_4)/\alpha
\end{pmatrix}.
\end{equation}
The five quantities $p_1,p_3,q_3,q_4,a_{21}$ are arbitrary. Conversely, \eqref{m7:s:eq:R2-cross-free}--\eqref{m7:s:eq:R2-cross-free-A} satisfy $x^tJx=\alpha\delta J$.

To place these last five directions in the even orbit, set
\begin{equation}\label{m7:s:eq:R2-final-normalize}
T_R=\begin{pmatrix}q_3/\alpha&0\\q_4/\alpha&-q_3/\alpha\end{pmatrix},\qquad
T_L=\begin{pmatrix}
-p_3/\alpha&p_1/\alpha\\
(a_{21}+p_3q_3/\alpha)/\alpha&p_3/\alpha
\end{pmatrix}.
\end{equation}
Both are traceless, and
\[
u(0,T_L)x\,u(0,T_R)^t
=x_2\left(0,1,\begin{pmatrix}\alpha&0\\0&\delta\end{pmatrix}\right).
\]
Since all these operations have rational inverses on the specified
charts, they give local coordinates of dimensions $5+6+5=16$.
This proves embedded smoothness and the claimed dimension, as well
as the unequal-eigenvalue identification in \eqref{m7:s:eq:R2-E2}.
\end{proof}

\paragraph{Orbit tangent, parameter directions, and their intersection}

Let $S_{\mathrm{res}}$ be the five-parameter family of representatives
inside the braces in \eqref{m7:s:eq:R2-residual}.
For the tangent calculation, write
\[
x=x_2\left(s,t,\begin{pmatrix}\alpha&0\\\gamma&\delta\end{pmatrix}\right),\qquad
\Delta=\alpha-\delta,\qquad \rho=\alpha s+\gamma t.
\]
Thus
\begin{equation}\label{m7:s:eq:R2-res-matrix}
x=\begin{pmatrix}
\rho&\delta t&0&0&0&0\\
0&0&0&0&-\rho/t&-\delta\\
0&0&\alpha&0&0&0\\
0&0&\gamma&\delta&0&0\\
-\alpha&0&0&0&0&0\\
0&0&0&0&-\alpha/t&0
\end{pmatrix},\qquad t\alpha\delta\Delta\ne0.
\end{equation}
Define
\[
\nu(X,Z)=\begin{pmatrix}0&X&Z\\0&0&-X^\#\\0&0&0\end{pmatrix},\qquad \tr Z=0,
\]
so $\nu(X,Z)\in\fu$. The differential of the action on the five-parameter family is
\begin{equation}\label{m7:s:eq:R2-differential}
\begin{split}
&\nu(X_L,Z_L)x+x\nu(X_R,Z_R)^t+[d(b),x]\\
&\hspace{8mm}+\dot s\,x_s+\dot t\,x_t+\dot\alpha\,x_\alpha
+\dot\gamma\,x_\gamma+\dot\delta\,x_\delta.
\end{split}
\end{equation}
Its complete kernel has seven free scalars $r_1,r_2,r_3,\lambda,\sigma,\theta,h$:
\begin{equation}\label{m7:s:eq:R2-kernel-parameters}
\dot s=\sigma,\quad \dot t=t\theta,\quad \dot\alpha=\dot\delta=0,
\quad\dot\gamma=\Delta h-\gamma\theta,
\end{equation}
\begin{equation}\label{m7:s:eq:R2-kernel-bX}
b=\begin{pmatrix}\lambda-\theta&0\\-h&\lambda\end{pmatrix},\quad
X_R=\begin{pmatrix}r_1&r_2\\0&0\end{pmatrix},\quad
X_L=\begin{pmatrix}0&0\\(\gamma r_1+\delta r_2)/\alpha&-r_1\end{pmatrix},
\end{equation}
\begin{equation}\label{m7:s:eq:R2-kernel-Z}
Z_R=\begin{pmatrix}th&r_3\\0&-th\end{pmatrix},\qquad
Z_L=\begin{pmatrix}
\sigma+th&0\\-\delta r_3/\alpha-2h\rho/\alpha&-\sigma-th
\end{pmatrix}.
\end{equation}
For completeness, the $(3,3)$ block of \eqref{m7:s:eq:R2-differential} forces $b_{12}=0$ and $b_{22}-b_{11}=\dot t/t$. The middle block gives $\dot\alpha=\dot\delta=0$ and the displayed formula for $\dot\gamma$. The $(2,3)$ and $(3,2)$ blocks force the first row of $X_L$ and second row of $X_R$ to be zero. The $(1,2)$ block gives the remaining relation between them. The outer blocks then give \eqref{m7:s:eq:R2-kernel-Z}. Conversely, substitution into \eqref{m7:s:eq:R2-differential} gives zero. These eliminations are valid throughout the parameter domain $t\alpha\Delta\ne0$.

Since $\dim K=18$, the differential rank is $18+5-7=16$, agreeing with \cref{m7:s:lem:R2-smooth}. Setting the parameter velocities to zero leaves a four-dimensional kernel, so
\begin{equation}\label{m7:s:eq:R2-tangent-summary}
\begin{split}
T_x(Kx)&=\{\nu(X_L,Z_L)x+x\nu(X_R,Z_R)^t+[d(b),x]\},\qquad \dim T_x(Kx)=14,\\
T_xS_{\mathrm{res}}\cap T_x(Kx)&=\operatorname{span}\{x_s,x_t,x_\gamma\},\\
T_x\Rres&=T_x(Kx)\oplus\operatorname{span}\{x_\alpha,x_\delta\}.
\end{split}
\end{equation}
Indeed, the parameter projections of the kernel span $\partial_s$, $t\partial_t-\gamma\partial_\gamma$, and $\Delta\partial_\gamma$. Thus three of the five parameter directions are orbit directions, leaving $\alpha,\delta$ as transverse coordinates.

\paragraph{The conormal calculation}

Under \eqref{m7:s:eq:tangent-convention}, the ambient conormal is represented by
\begin{equation}\label{m7:s:eq:R2-full-conormal}
Z=\begin{pmatrix}
z_0&0&0&0&0&0\\
z_1&0&0&0&0&0\\
z_2&0&-z_0&0&0&0\\
z_3&0&z_4&z_0&0&0\\
z_5&0&0&0&0&0\\
z_6&-z_5&z_3&-z_2&z_1&-z_0
\end{pmatrix},
\qquad
(\delta-\alpha)z_4=2\alpha z_1-\frac{2\alpha z_5}{t}+2\gamma z_0.
\end{equation}
The cotangent vector on $T_xG_7^\circ$ is $\dot x\mapsto\tr(Zx^{-1}\dot x)$, with ambient representative $Y=Zx^{-1}$.

To verify the annihilation equations in \eqref{m7:s:eq:R2-tangent-summary},
first note that $Z^tJ+JZ=0$, so $Z\in\mathfrak{sp}_6$. The right-unipotent pairing vanishes because $Z$ and $\nu(X_R,Z_R)^t$ are block lower triangular and the latter has zero diagonal blocks. The matrix $xZx^{-1}$ is block upper triangular, so its pairing with $\nu(X_L,Z_L)$ is also zero. For the reductive tangent directions, the sum of the three diagonal $2\times2$ blocks of $xZx^{-1}-Z$ is
\[
\begin{pmatrix}0&0\\
-2z_1+2z_5/t-2\gamma z_0/\alpha+(\delta-\alpha)z_4/\alpha&0
\end{pmatrix},
\]
which is zero by \eqref{m7:s:eq:R2-full-conormal}.

Each $x^{-1}x_p$ for $p=s,t,\alpha,\gamma,\delta$ is lower triangular. Its diagonal is
\begin{equation}\label{m7:s:eq:R2-param-diagonals}
\begin{array}{c|c}
p&\diag(x^{-1}x_p)\\ \hline
s&(0,0,0,0,0,0)\\
t&(0,1/t,0,0,-1/t,0)\\
\alpha&(1/\alpha,0,1/\alpha,0,1/\alpha,0)\\
\gamma&(0,0,0,0,0,0)\\
\delta&(0,1/\delta,0,1/\delta,0,1/\delta).
\end{array}
\end{equation}
In every case, the trace pairing with
\eqref{m7:s:eq:R2-full-conormal} vanishes. Its six free parameters
give a six-dimensional annihilator; since
$\dim G_7^\circ-\dim\Rres=22-16=6$, this is the full conormal.

Take $z_0=\Delta$, $z_4=-2\gamma$, and all the other $z_i$ equal to zero. We obtain
\begin{equation}\label{m7:s:eq:R2-witness}
Z_*=\Delta\diag(1,0,-1,1,0,-1)-2\gamma E_{43},\qquad
Y_*=Z_*x^{-1}.
\end{equation}
Thus $Y_*x=Z_*$ is the identity translate required by \eqref{m7:s:eq:cotangent-translation}. One may also check its annihilation by putting
\[
A_0=\begin{pmatrix}\Delta&0\\0&0\end{pmatrix},\quad
A_1=\begin{pmatrix}0&0\\0&-\Delta\end{pmatrix},\quad
B_0=\begin{pmatrix}-\Delta&0\\-2\gamma&\Delta\end{pmatrix}.
\]
Then
\begin{equation}\label{m7:s:eq:R2-witness-conjugate}
Z_*=\diag(A_0,B_0,A_1),\qquad
xZ_*x^{-1}=\begin{pmatrix}
A_1&0&\dfrac{\Delta\rho}{\alpha}\diag(-1,1)\\
0&B_0&0\\0&0&A_0
\end{pmatrix}.
\end{equation}
This has exactly the left/right triangular pattern used for the conormal witnesses in \cite[Propositions~6.1 and~4.8]{JLX}.

The non-nilpotence certificate is
\begin{equation}\label{m7:s:eq:R2-witness-nonnil}
\det(zI_6-Z_*)=z^2(z-\Delta)^2(z+\Delta)^2,
\qquad \tr(Z_*^2)=4\Delta^2\ne0.
\end{equation}
At an arbitrary point $x'=LxR$ of the same $K$-orbit, the transported ambient cotangent vector is $Y_*'=R^{-1}Y_*L^{-1}$, and
\[
Y_*'x'=R^{-1}(Y_*x)R.
\]
Hence non-nilpotence persists at every point of $\Rres$.

\begin{proof}[Proof of \cref{m7:s:prop:R2}]
By \eqref{m7:s:eq:R2-support-N2}--\eqref{m7:s:eq:R2-support-res},
we have both the residual support condition and nilpotent wavefront.
Moreover, \cref{m7:s:lem:R2-smooth} makes this support a smooth closed
submanifold of $\mathcal W^\circ$, and
\eqref{m7:s:eq:R2-witness-nonnil} supplies a nonnilpotent conormal
direction at every point. Applying \cref{cf:conormal} gives
$S|_{\mathcal W^\circ}=0$.

Returning to $\mathcal W$, the support is now contained in the closed equal-eigenvalue locus $\mathcal E_2$. Equation~\eqref{m7:s:eq:R2-equal-transporter}, \cref{cf:stabilizer}, and \eqref{cf:support-sequence} eliminate this remaining distribution. Thus $S|_{\mathcal W}=0$.
\end{proof}

\begin{remark}\label{m7:s:rem:R2-scope}
The conormal exclusion in \cref{m7:s:prop:R2} applies to distributions on $\mathcal W^\circ$ supported in $\Rres$ with nilpotent wavefront fibers. The final closed-support step on $\mathcal E_2$ uses its intrinsic odd vanishing.
\end{remark}

\subsection{The period conclusion}\label{m7:endpoint}
The following duality identifies the two periods in the GK product
bound and completes Theorem~\ref{m7:main}.

\subsubsection{Transpose and contragredients}

Let $\tau$ be componentwise transpose on $G$ and ordinary transpose on $G_7^\circ$. It preserves $G_7^\circ$: from $x^tJx=\ell(x)J$ and $J^{-1}=-J$, one gets $xJx^t=\ell(x)J$. Define the automorphism
\begin{equation}\label{m7:s:eq:iota}
\iota(g,z)=(g^{-t},z^{-t})=\tau((g,z)^{-1}).
\end{equation}
Transpose takes $U$ to $U^t$. The two periods use $H$ and $\iota(H)$, respectively.

\begin{proposition}\label{m7:s:prop:RV}
For every irreducible admissible representation $\pi$ of $G$, one has
\[
\pi\circ\iota\simeq\pi^\vee,\qquad
\pi^\vee\circ\iota\simeq\pi.
\]
\end{proposition}
\begin{proof}
For an alternating space $V$ and an anti-symplectic involution $h$,
we use Roche--Vinroot's dualizing automorphism
\begin{equation}\label{m7:s:eq:RV-auto}
g\longmapsto\ell(g)^{-1}hgh^{-1}.
\end{equation}
It sends every irreducible smooth representation of the similitude group to its contragredient~\cite[Theorem~B]{RVfactor}. Their pointwise conjugacy statement and character argument also apply to the product of the two formed spaces; local integrability is supplied by \cite[Theorem~1, \S3.3]{Clozel}.

Take $h_2=\diag(1,-1)$ and $h_6=d(h_2)$. They are anti-symplectic involutions for $j$ and $J$. On our two groups,
\[
g^{-t}=\ell(g)^{-1}JgJ^{-1},\qquad
z^{-t}=(\det z)^{-1}jzj^{-1}.
\]
Thus \eqref{m7:s:eq:iota} differs from the product of \eqref{m7:s:eq:RV-auto} by conjugation by
\[
(Jh_6^{-1},jh_2^{-1})\in G_7^\circ\times\GL_2(F).
\]
Inner conjugation preserves representation isomorphism classes,
giving the first isomorphism. Composing with $\iota^2=1$ gives
the second.
\end{proof}

Let $\calA=H\times H$ act on $G$ by
\begin{equation}\label{m7:s:eq:global-action}
(h_1,h_2)\cdot g=h_1gh_2^t,
\qquad \eta(h_1,h_2)=\xi(h_1)\xi(h_2).
\end{equation}
Transpose exchanges the two factors, and preserves $\eta$. Put
\begin{equation}\label{m7:s:eq:global-extension}
\tA=\calA\rtimes\langle\tau\rangle,
\qquad \teta|_{\calA}=\eta,\qquad \teta(\tau)=-1.
\end{equation}

\subsubsection{Completion of the multiplicity-one theorem}\label{m7:s:sec:completion}

\begin{proof}[Proof of \cref{m7:s:thm:descended}]
Let $S\in\Dist_{\tetaK}(G_7^\circ)\cap\NW(G_7^\circ)$. \Cref{m7:s:prop:R2}, whose proof uses \cref{m7:s:prop:R0,m7:s:prop:R1}, gives
\[
S|_{\{x_{33}\ne0\}}=0.
\]
Thus $S$ is supported on the closed set
$\{x_{33}=0\}=R_3\sqcup R_4$, with $R_3$ open and $R_4$ closed
there. Applying the intrinsic vanishing of \cref{m7:s:prop:R3}
leaves support on $R_4$, where \cref{m7:s:prop:R4} gives $S=0$.

This proves \eqref{m7:s:eq:descended-target}; only the $R_2$ step uses NW.
\end{proof}

\begin{proof}[Proof of \cref{m7:main}]
Let $T\in\Dist_{\teta}(G)\cap\NW(G)$. Under the Frobenius isomorphism for the transpose-extended stabilizer \eqref{m7:s:eq:Frob-isom}, it corresponds to $S\in\Dist_{\tetaK}(G_7^\circ)$. The density character is trivial by \eqref{m7:s:eq:Frob-density}, and the wavefront transfer \eqref{m7:s:eq:Frob-native} gives $S\in\NW(G_7^\circ)$. \Cref{m7:s:thm:descended} gives $S=0$, hence $T=0$.

Thus \eqref{m7:global-target} holds. Generalized coefficients satisfy NW
by \eqref{cf:center-span} and Theorem~\ref{cf:nilpotent}. Applying
\Cref{cf:gk} with \eqref{cf:coefficient-actions}--\eqref{cf:product}
and the duality of \cref{m7:s:prop:RV} gives
\[
\dim_\CC\Hom_H(\pi,\xi)\le1
\]
for every irreducible admissible smooth complex representation of $G$.
The common-center condition \eqref{m7:s:eq:center} merely identifies
the cases with zero period.
\end{proof}

\section{Two flag coarsenings and their boundaries: Model 6}\label{cl:model6}\label{m6:section}

\paragraph{Proof strategy}
We descend along the $2+6+2$ and $4+2+4$ coarsenings of the
$2+2+2+2+2$ symplectic flag. Their Gauss opens reduce to
Model~7 and representation-valued Shalika, respectively; their
distribution theorems clear these opens. We then enlarge this open successively:
its closed complement contains the support, and a boundary argument
clears the next piece. The rank opens of the ninety rational cells
in \S\ref{m6:sieve} organize this induction. Unipotent and full
stabilizer characters remove incompatible cells; the remaining
boundary steps form the complete twenty-three-step filtration in
Table~\ref{m6:s:tab:assembly-filtration}.

For cells $7,8$, an opposite-weight plane and a stabilizer jump give
intrinsic odd vanishing. We then treat the pairs $(9,14)$
and $(12,22)$ by ambient conormal exclusion; intervening character
steps supply the predecessor opens needed for the second pair.
Since these pairs admit intrinsic odd distributions, the conormal
arguments apply to the original global NW distribution after the
earlier opens have vanished. The orbital-family criterion and the
remaining stabilizer-character arguments then remove the middle and
late families, while a relative support filtration handles the
entire outer interval.
Applying the closed-support sequences in this order proves the
distribution theorem (\S\ref{m6:filtration}). Generalized
coefficients then give the GK product bound, and similitude-group
duality identifies its two period spaces.

\subsection{The similitude datum}\label{m6:datum}
Let $w_5$ be the five-by-five antidiagonal permutation matrix,
$j=\left(\begin{smallmatrix}0&-1\\1&0\end{smallmatrix}\right)$, and
$d_5(a)=\diag(a,a,a,a,a)$. Put $J_6=w_5\otimes j$ and
$G_6=\GSp(J_6)(F)$. The five graded blocks have size two.
Let $P_6=L_6U_6$ stabilize the standard
isotropic flag of dimensions $2\subset4$. We take
$H_6=U_6\rtimes d_5(\GL_2(F))$, with
$\xi_6(ud_5(a))=\psi(\operatorname{tr}u_{12}+\operatorname{tr}u_{23})$.
We derive the Levi embedding, character stabilizer, and density
characters below.

\begin{theorem}[Model 6]\label{m6:main}
For every irreducible admissible smooth complex representation $\pi$
of $G_6$,
\[
 \dim\Hom_{H_6}(\pi,\xi_6)\leq1.
\]
This holds over every non-Archimedean local field $F$ of
characteristic zero, including residue characteristic two, for the
fixed nontrivial smooth additive character $\psi$.
\end{theorem}

\begin{theorem}[The NW distribution statement]\label{m6:distribution}
Every ordinary $\eta_6$-equivariant, transpose-odd distribution
$T\in\NW(G_6)$ is zero. Explicitly, the hypothesis is
\[
 \WF(T)\subset\cN(\mathfrak{sp}_{10}^*)
                   \times\{0\}_{\mathrm{central\ cotangent}}.
\]
\end{theorem}

In this section write $G=G_6$, $H=H_6$, $U=U_6$, $P=P_6$, $J=J_6$,
and $\xi=\xi_6$. The action is $(h_1,h_2)g=h_1gh_2^t$ and
$\eta_6(h_1,h_2)=\xi_6(h_1)\xi_6(h_2)$. For the rational parabolic double cosets we use
$\mathcal P_i=P(F)W_iP^t(F)$, $1\leq i\leq90$, numbered by
decreasing lexicographic order of their sixteen corner-rank entries
(Proposition~\ref{m6:s:prop:cells-rational}). The signed monomial
representatives $W_i$, including their column conventions, are given
in Section \texttt{M6-CELLS} of \path{anc/model6_calculations.pdf}
and in \path{anc/model6/cell_compatibility.json}.
The codes \texttt{P001},\ldots,\texttt{P090} denote these same
cells. In particular, all cell numbers used in the Gauss-open arguments
below refer to this fixed ordering.

\Needspace{10\baselineskip}
\subsubsection{The Wan--Zhang Model 6}\label{m6:s:sec:setup}

\paragraph{The alternating form and the similitude group}
Let $w_5$ denote the five-by-five reversal matrix and put
\begin{equation}\label{m6:s:eq:setup-form}
j=\begin{pmatrix}0&-1\\1&0\end{pmatrix},\qquad
J=w_5\otimes j=
\begin{pmatrix}
0&0&0&0&j\\0&0&0&j&0\\0&0&j&0&0\\0&j&0&0&0\\j&0&0&0&0
\end{pmatrix}.
\end{equation}
Every block in this display has size two. This is the recursive Chapter~8 realization of \cite[\S8.1]{WZ}; $J^t=-J$ and $J^2=-I_{10}$. The group is
\[
G=\{g\in\operatorname{GL}_{10}(F):g^tJg=\ell(g)J\}.
\]
Its multiplier is the homomorphism $\ell:G\to F^\times$, and its dimension is $56$.

For $A\in M_2(F)$ define
\[
A^\sharp=A^\#=j^{-1}A^tj=\operatorname{adj}(A)=\operatorname{tr}(A)I_2-A.
\]
This definition applies to singular matrices as well. If $a$ is invertible, the notation $a^*$ of \cite{WZ} is $ja^{-t}j^{-1}=a/\det a$. Thus the sharp operation used in Lie coordinates and the star operation in the original Levi definition have different meanings.

\paragraph{The symplectic flag and its two coarsenings}\label{m6:flag-covector}
Let $0\subset V_1\subset V_2\subset V_2^\perp\subset V_1^\perp\subset V$
be the self-dual filtration, with $\dim V_1=2$, $\dim V_2=4$.
Write $E_1=V_1$, $E_2=V_2/V_1$, and $W=V_2^\perp/V_2$.
If $L$ is the one-dimensional target of the alternating form,
the graded pieces are
\[
 E_1,\quad E_2,\quad W,\quad E_2^\vee\otimes L,
                         \quad E_1^\vee\otimes L.
\]
Here $W$ is a two-dimensional alternating space with values in $L$.
The radical abelianization is
\[
 U/[U,U]=\Hom(E_2,E_1)\oplus\Hom(W,E_2).
\]
For $X:E_2\to E_1$ and $Y:W\to E_2$, the fixed identifications
$E_1\simeq E_2\simeq W$ define the character
$\psi(\operatorname{tr}X+\operatorname{tr}Y)$ by the two trace pairings.
We prove this abelianization by the commutator calculation below;
the nondegenerate trace pairing then identifies the Levi stabilizer
with the diagonal $\GL_2$.

Retaining $V_1$ and forgetting $V_2$ gives the $2+6+2$ parabolic:
the first arrow belongs to its radical and the second survives on
$V_1^\perp/V_1$. Retaining $V_2$ and forgetting $V_1$ gives the
$4+2+4$ parabolic: the first arrow is internal to $\GL(V_2)$ and
the second belongs to its radical. This gives the two Levi extractions and splits the character in the
full similitude group. To treat their boundaries, we first classify
the relative positions of the underlying self-dual flags by parabolic
double cells. The coefficient families retain the characters on the
graded pieces, whose trace functionals must cancel on common unipotent
directions. We then analyze the full stabilizers and boundary
wavefronts on these families.

\paragraph{The parabolic and its unipotent radical}
Let $P=LU$ stabilize the coordinate isotropic flag of dimensions $2\subset4$. Its matrices are upper triangular in the five two-dimensional blocks. The diagonal symplectic equations give
\begin{equation}\label{m6:s:eq:setup-levi}
m(A,B,C)=\operatorname{diag}
\left(A,B,C,\frac{\det C}{\det B}B,\frac{\det C}{\det A}A\right).
\end{equation}
Indeed the middle equation is $C^tjC=\ell j$, so $\ell=\det C$; the paired outer equations then uniquely determine the fourth and fifth blocks. Thus $L\simeq\operatorname{GL}_2^3$, with $A,B,C$ arbitrary invertible matrices. The opposite parabolic is $P^t=LU^t$.

We specify all unipotent entries using finite exponential coordinates. Put
\begin{equation}\label{m6:s:eq:setup-N}
\mathcal N(X,Y,Z,W,T_0,S_0)=
\begin{pmatrix}
0&X&Z&W&S_0\\
0&0&Y&T_0&-W^\#\\
0&0&0&-Y^\#&-Z^\#\\
0&0&0&0&-X^\#\\
0&0&0&0&0
\end{pmatrix},
\end{equation}
where $X,Y,Z,W\in M_2(F)$ and $T_0,S_0\in\mathfrak{sl}_2(F)$. The equation $\mathcal N^tJ+J\mathcal N=0$ pairs the $(a,b)$ block with $(6-b,6-a)$ and imposes the negative sharp in \eqref{m6:s:eq:setup-N}. On the antidiagonal it reads $T_0+T_0^\#=S_0+S_0^\#=0$, which is precisely the two trace-zero conditions. This proves both the sufficiency and completeness of the displayed linear coordinates.

Every strictly upper block triangular matrix here has fifth power zero. Consequently
\begin{equation}\label{m6:s:eq:setup-exp}
U=\{\exp\mathcal N:\mathcal N\text{ as in \eqref{m6:s:eq:setup-N}}\},
\qquad
\exp\mathcal N=I+\mathcal N+\frac{\mathcal N^2}{2}
 +\frac{\mathcal N^3}{6}+\frac{\mathcal N^4}{24}.
\end{equation}
Finite logarithm gives the inverse, so these are polynomial inverse
maps over every characteristic-zero field considered here. All
denominators are valid in $F$, for every residue characteristic,
and $\dim U=22$.

The adjacent blocks of $u=\exp\mathcal N$ are $u_{12}=X$, $u_{23}=Y$, $u_{34}=-Y^\#$, and $u_{45}=-X^\#$. Multiplication of upper unitriangular matrices adds these blocks, so $u\mapsto(X,Y)$ is a homomorphism onto $M_2(F)\oplus M_2(F)$. Its kernel is the commutator subgroup. To see the latter equality, grade \eqref{m6:s:eq:setup-N} by block distance. Bracketing an $X$ direction with a $Y$ direction gives the $Z$ block $XY$, spanning $M_2$. Brackets of $Y$ directions give the trace-zero $T_0$ block; the matrix units span this three-dimensional space. Further brackets with the $X$ directions generate the $W$ block and then the trace-zero $S_0$ block. These are exactly the degrees at least two. The exponential correspondence therefore gives
\[
U/[U,U]\simeq M_2(F)\oplus M_2(F).
\]

\paragraph{The subgroup, character, and densities}
Write $d_5(a)=\operatorname{diag}(a,a,a,a,a)$. The subgroup and character are
\begin{equation}\label{m6:s:eq:setup-model}
H_0=d_5(\operatorname{GL}_2(F)),\qquad H=U\rtimes H_0,\qquad
\xi(u d_5(a))=\psi(\operatorname{tr}X+\operatorname{tr}Y).
\end{equation}
These are exactly the data of \cite[\S8.1]{WZ}. Conjugation by $m(A,B,C)$ sends $X$ to $AXB^{-1}$ and $Y$ to $BYC^{-1}$. Preservation of the character on every scalar multiple of every $X,Y$ gives
\[
\operatorname{tr}((B^{-1}A-I)X)=0,\qquad
\operatorname{tr}((C^{-1}B-I)Y)=0.
\]
Nondegeneracy of the trace pairing on $M_2(F)$ forces $A=B=C$,
identifying the Levi stabilizer of the character with $H_0$.

Set $a=\det A$, $b=\det B$, and $c=\det C$ for this density calculation. The determinants of conjugation on the six independent coordinate spaces in \eqref{m6:s:eq:setup-N} are
\[
\begin{array}{c|cccccc}
\text{coordinate}&X&Y&Z&W&T_0&S_0\\ \hline
\det\operatorname{Ad}(m)&(a/b)^2&(b/c)^2&(a/c)^2&
a^2b^2c^{-4}&(b/c)^3&(a/c)^3.
\end{array}
\]
For example, the determinant of $M_2\ni X\mapsto AXB^{-1}$ is $(a/b)^2$. The $W$ action is $W\mapsto c^{-1}AWB^\#$, of determinant $a^2b^2c^{-4}$. On $T_0$ the action is $(b/c)\operatorname{Ad}(B)$, whose determinant is $(b/c)^3$ because conjugation has determinant one on $\mathfrak{sl}_2$; the $S_0$ calculation has $A$ in place of $B$. Multiplication gives the parabolic modulus
\begin{equation}\label{m6:s:eq:setup-modulus}
\delta_P(m(A,B,C))=|a|^9|b|^5|c|^{-14}.
\end{equation}
Its restriction to $H_0$ is one. The group $H_0$ is unimodular, $U$ is unipotent, and $H_0$ has determinant one on every independent $U$ coordinate. Therefore
\begin{equation}\label{m6:s:eq:setup-unimodular}
\Delta_H=\Delta_{\mathcal A}=1.
\end{equation}
Here $\Delta_K$ is characterized, for a left Haar measure, by $\int f(xk)\,dx=\Delta_K(k)^{-1}\int f(x)\,dx$. The determinant $\delta_P$ on the unipotent radical is distinguished from this group modular convention.

Every $zI_{10}$ lies in $H$, and $\xi(zI_{10})=1$. If $\omega_\pi$ is the central character of an irreducible admissible representation, then a nonzero period forces $\omega_\pi=1$: the two equalities $\lambda(\pi(z)v)=\omega_\pi(z)\lambda(v)$ and $\lambda(\pi(z)v)=\lambda(v)$ give the assertion. 

\paragraph{The two-sided action and transpose}
The identity $g^tJg=\ell(g)J$ also gives $gJg^t=\ell(g)J$, so transpose preserves $G$. It is a continuous anti-involution. Put
\begin{equation}\label{m6:s:eq:setup-action}
\mathcal A=H\times H,\qquad (h_1,h_2)\cdot g=h_1gh_2^t,\qquad
\eta(h_1,h_2)=\xi(h_1)\xi(h_2).
\end{equation}
Transpose satisfies $\tau((h_1,h_2)\cdot g)=(h_2,h_1)\cdot\tau(g)$. It therefore acts on the character extension $\widetilde{\mathcal A}$ with $\widetilde\eta(\tau)=-1$. Both $D'_\eta$ and $D'_{\widetilde\eta}$ below use the convention in \eqref{cf:action}.

For the eventual representation-theoretic reduction, let $\alpha(g)=g^{-t}$. The second subgroup is $\alpha(H)=H^t$ and its character is specified by $\chi_2(\alpha(h))=\xi(h)$, or equivalently $\chi_2(h^t)=\xi(h)^{-1}$. Section~\ref{m6:filtration} constructs the smoothing maps, derives both test-function characters, and uses similitude-group duality to identify the two period spaces. This fixes the sign in \eqref{m6:s:eq:setup-action} before any local calculation is used.

\subsection{Retaining the first isotropic subspace: Model 7}\label{m6:model7}
Retaining the first step of the flag gives the $2+6+2$ parabolic.
Its middle symplectic quotient carries the second trace term of the character,
so the descendant is Model 7 in the Levi subgroup $\GSp_6\times\GL_2$.
The extraction chart identifies the descended datum with Model~7.
Its character, transpose, density, and wavefront calculations
preserve the hypotheses needed to apply Theorem~\ref{m7:native}.
\subsubsection{The initial determinant open and the Model 7 geometry}\label{m6:s:sec:model7}

\begin{proposition}\label{m6:s:prop:model7-open}
Let $T\in D'_{\widetilde\eta}(G)\cap\NW(G)$. Then
\[
T|_\Omega=0,\qquad
\Omega=\{g:\det g_{55}\ne0\}=\mathcal P[1,2,3,4,5].
\]
\end{proposition}

We first construct the descent and transfer the wavefront condition,
including its central cotangent components. This allows us to apply
the distributional statement of Theorem~\ref{m7:native}, proved in
Section~\ref{m7:section}.

\paragraph{The outer parabolic and the intersection model}
Let $Q=MN$ be the upper parabolic with block sizes $2,6,2$, and put
$J_{\mathrm{mid}}=w_3\otimes j$. Its Levi embedding is
\begin{equation}\label{m6:s:eq:model7-embedding}
\iota(B,A)=\operatorname{diag}
\left(A,B,\frac{\ell(B)}{\det A}A\right),\qquad
B\in\operatorname{GSp}(J_{\mathrm{mid}})(F),\ A\in\operatorname{GL}_2(F).
\end{equation}
Thus $M$ is the direct product $\operatorname{GSp}_6\times\operatorname{GL}_2$,
with $A$ and $B$ independent. Its radical $N$ has Lie algebra
\eqref{m6:s:eq:setup-N} with $Y=T_0=0$, and hence dimension $15$.
Let $U_7$ be the unipotent group in the central three two-dimensional
blocks. The two semidirect products are
\[
U=N\rtimes U_7,\qquad H=N\rtimes H_M,
\]
where
\begin{equation}\label{m6:s:eq:model7-intersection}
H_M=H\cap M=
\{\iota(u_7d_3(a),a):u_7\in U_7,\ a\in\operatorname{GL}_2(F)\}.
\end{equation}
Under \eqref{m6:s:eq:model7-embedding} this is exactly the seventh model in
\cite[\S8.2]{WZ}. The characters are
\begin{equation}\label{m6:s:eq:model7-characters}
\xi_N(n)=\psi(\operatorname{tr}X_n),\qquad
\xi_M(\iota(u_7d_3(a),a))=\psi(\operatorname{tr}Y_{u_7}).
\end{equation}
The block $Y$ occupies position $(2,3)$ in the five-block
decomposition and is the first adjacent block in the central
three-block realization. Hence the trace character retains its sign
and normalization.

The determinants of conjugation on the twelve-dimensional outer
rectangular space and the three-dimensional remaining radical space
are, respectively,
\[
\left(\frac{\det A}{\ell(B)}\right)^6,\qquad
\left(\frac{\det A}{\ell(B)}\right)^3.
\]
For the first, the action on $\operatorname{Hom}(F^6,F^2)$ has determinant
$(\det A)^6(\det B)^{-2}$, and $\det B=\ell(B)^3$.
For the second, it is the scalar $(\det A)/\ell(B)$ times conjugation
on $\mathfrak{sl}_2$. Therefore
\begin{equation}\label{m6:s:eq:model7-modulus}
\delta_Q(\iota(B,A))=
\left|\frac{\det A}{\ell(B)}\right|^9,\qquad
\delta_Q|_{H_M}=1.
\end{equation}
The group $H_M$ preserves $\xi_N$: on the top $2$-by-$6$ radical
block, right multiplication by an upper block $u_7^{-1}$ leaves
its first two-column block unchanged, and the diagonal $a$ acts
there by conjugation.

\paragraph{The rational big-cell coordinates}
Multiplication gives an $F$-rational isomorphism
\begin{equation}\label{m6:s:eq:model7-big-cell}
N\times M\times N\longrightarrow\Omega,\qquad
(n_1,m,n_2)\longmapsto n_1mn_2^t.
\end{equation}
To prove this, we use block elimination. In $2+6+2$
blocks write
\[
g=\begin{pmatrix}g_{aa}&g_{ab}&g_{ac}\\
g_{ba}&g_{bb}&g_{bc}\\g_{ca}&g_{cb}&D\end{pmatrix},
\qquad \det D\ne0.
\]
An element of $N$ has the form
\[
n(C,E)=\begin{pmatrix}
I_2&jC^tJ_{\mathrm{mid}}&E\\0&I_6&C\\0&0&I_2
\end{pmatrix},\qquad
E^tj+jE+C^tJ_{\mathrm{mid}}C=0.
\]
Reading the last column gives $C_1=g_{bc}D^{-1}$ and
$E_1=g_{ac}D^{-1}$. Its isotropy equation in $g^tJg=\ell(g)J$
is precisely the displayed equation for $n(C_1,E_1)$.
Likewise, the last row gives
$C_2^t=D^{-1}g_{cb}$ and $E_2^t=D^{-1}g_{ca}$.
The identity $gJg^t=\ell(g)J$ supplies the corresponding equation
for $n(C_2,E_2)$.
Removing these two factors leaves a matrix with bottom diagonal
block $D$ and zero other blocks in its last row and column.
Its symplectic equations force the remaining off-diagonal
$2$-by-$6$ and $6$-by-$2$ blocks to vanish, and its diagonal
blocks have the form \eqref{m6:s:eq:model7-embedding}.
The resulting inverse is unique and rational, since every
inversion uses $\det D\ne0$.

In these coordinates, the equivariant extraction on test functions is
\begin{equation}\label{m6:s:eq:model7-extraction}
\mathcal J f(m)=\int_{N\times N}
f(n_1mn_2^t)\xi_N(n_1)\xi_N(n_2)\,dn_1\,dn_2.
\end{equation}
Equivariance in each $N$ variable forces a character-weighted Haar
distribution, leaving an ordinary distribution on $M$.
By \eqref{m6:s:eq:model7-modulus}, conjugation by $H_M$ preserves
both Haar factors, as well as their characters. The residual action
and character are
\[
(h_1,h_2)m=h_1mh_2^t,\qquad
\eta_M(h_1,h_2)=\xi_M(h_1)\xi_M(h_2).
\]
Transpose exchanges the two Haar factors and sends $m$ to $m^t$.
Its density Jacobian has absolute value one. Thus extraction gives
the untwisted correspondence
\begin{equation}\label{m6:s:eq:model7-correspondence}
D'_{\widetilde\eta}(\Omega)
\simeq D'_{\widetilde\eta_M}(M).
\end{equation}
To see that test-function averaging is surjective, work in \eqref{m6:s:eq:model7-big-cell}: choose compactly supported
functions in the two $N$ variables with weighted integral one.
Their product with a test function on $M$ maps to that function.

\paragraph{The Frobenius fiber and its density}
Alternatively, descend over the bottom-right block map
$p(g)=g_{55}$, which transforms as $D\mapsto aDb^t$ under
the two reductive coordinates of $\mathcal A$.
Over $\Omega$ this base is $\operatorname{GL}_2(F)$, with rational
section $z\mapsto(d_5(z),1)$ and base point $I_2$.
The stabilizer at that point is
\[
K_\Omega=(U\times U)\rtimes\operatorname{GL}_2(F),\qquad
(u,v,a)\longmapsto(u d_5(a),v d_5(a^{-t})).
\]
The fiber is $Y_I=\{g:g_{55}=I_2\}$, of dimension $52$.
It has action
\[
(u,v,a)y=u d_5(a)y d_5(a)^{-1}v^t
\]
and character $\xi(u)\xi(v)$. The group $K_\Omega$ is
unimodular: each independent unipotent coordinate transforms
by conjugation, of determinant one. On the tangent of the base,
$Z\mapsto aZa^{-1}$ also has determinant one. Hence the relative
Frobenius character is trivial; transpose has absolute determinant
one on that tangent and adds no density sign.

Using the global coordinate map $(z,y)\mapsto d_5(z)y$, we identify
$\operatorname{GL}_2\times Y_I$ with $\Omega$ and average by
\[
\mathcal I f(y)=\int_{\operatorname{GL}_2(F)}f(d_5(z)y)\,dz.
\]
Proposition~\ref{cf:frobenius} gives $T(f)=S(\mathcal I f)$ with
exactly the fiber character just stated.
In the big-cell coordinates that fiber is
\begin{equation}\label{m6:s:eq:model7-fiber}
Y_I\simeq N\times\operatorname{GSp}_6(F)\times N,\qquad
(n_1,B,n_2)\mapsto n_1j_0(B)n_2^t,\quad
j_0(B)=\operatorname{diag}(\ell(B)I_2,B,I_2).
\end{equation}
These three blocks satisfy the similitude equation with multiplier $\ell(B)$.
After the two $N$ factors have been extracted, the action is
\begin{equation}\label{m6:s:eq:model7-K}
\begin{split}
K_7&=(U_7\times U_7)\rtimes\operatorname{GL}_2(F),\\
(u,v,a)B&=u d_3(a)B d_3(a)^{-1}v^t,\\
\eta_7(u,v,a)&=\psi(\operatorname{tr}Y_u+\operatorname{tr}Y_v).
\end{split}
\end{equation}
Transpose acts by $B\mapsto B^t$ and
$(u,v,a)\mapsto(v,u,a^{-t})$. This is the descended Model $7$
action in Section~\ref{m7:section}. The base coordinate in
\eqref{m6:s:eq:model7-fiber} is the last outer block, while the
$\operatorname{GL}_2$ coordinate in \eqref{m6:s:eq:model7-embedding}
is the first. Exchanging the two outer spaces by the symplectic
block permutation relates the two realizations and fixes the
intersection model.

\paragraph{Wavefront transfer, including the similitude direction}
We retain the assumed NW distribution $T$ on the original $G$.
Open restriction preserves its nilpotent wavefront fibers.
Since the factors $\xi_N(n_i)\,dn_i$ in
\eqref{m6:s:eq:model7-extraction} are smooth and nowhere zero,
their wavefront components vanish in product coordinates. At $(1,m,1)$ the translated radical
tangents are $\mathfrak n$ and $\mathfrak n^t$, since the Levi
normalizes them. The trace-orthogonal complement is
\[
(\mathfrak n+\mathfrak n^t)^\perp=\mathfrak m.
\]
A cotangent vector of the middle distribution therefore lifts to a
block-diagonal cotangent vector $Z\in\mathfrak m$ of the original
distribution. Nilpotence of $Z$ forces each of its diagonal
blocks to be nilpotent.

To include the center explicitly, write $\nu(B_0)$ for the
infinitesimal multiplier on $\mathfrak{gsp}_6$. A general Levi
Lie element has the form
\begin{equation}\label{m6:s:eq:model7-Levi-tangent}
Z=\operatorname{diag}
\bigl(A_0,B_0,A_0+(\nu(B_0)-\operatorname{tr}A_0)I_2\bigr).
\end{equation}
If $Z$ is nilpotent, $A_0,B_0$ are nilpotent,
$\operatorname{tr}A_0=0$, and $\nu(B_0)=0$.
For a general Levi tangent $Z'$ determined by $(a_0,b_0)$,
the trace pairing is then
\[
\operatorname{tr}(ZZ')=
2\operatorname{tr}(A_0a_0)+\operatorname{tr}(B_0b_0).
\]
The cotangent vector on the $\mathfrak{gsp}_6$ factor is therefore
represented by the nilpotent $B_0$, with zero central component.

The Frobenius coordinate change is an analytic diffeomorphism
with a smooth base Haar factor. Thus the same statement survives
passage to \eqref{m6:s:eq:model7-fiber}. Differentiating its first
block $\ell(B)I_2$ contributes a scalar outer tangent, whose
pairing with a nilpotent outer cotangent vector is zero. Consequently
the extracted distribution $S_7$ satisfies
\begin{equation}\label{m6:s:eq:model7-WF-transfer}
WF(S_7)\subset\mathcal N(\mathfrak{gsp}_6^*).
\end{equation}
This is the NW condition needed for the Model-7 distribution theorem.

\subparagraph{Vanishing after extraction.}
The action computed above and the wavefront transfer
\eqref{m6:s:eq:model7-WF-transfer} allow us to apply
Theorem~\ref{m7:native} to the extracted distribution. Injectivity
of the two averaging correspondences then proves
Proposition~\ref{m6:s:prop:model7-open}.

\subsection{Retaining the four-plane: the Shalika descendant}\label{m6:typea}
For the $4+2+4$ parabolic, the first trace term of the character lies in the
internal $2+2$ unipotent of the $\GL_4$ factor, while the second
belongs to the extracted radical. The Levi subgroup is
$\GL_4\times\GL_2$, with the second factor acting diagonally in
the period. Applying Theorem~\ref{cf:shalika-nw} gives vanishing on
the entire Gauss open $\Omega'=\mathcal P_1\sqcup\mathcal P_6\sqcup\mathcal P_{34}$.
This auxiliary Shalika descendant corresponds to the $4+2+4$
coarsening. Independent rank-one and rank-zero calculations are
given in the appendices.
\subsubsection{The \texorpdfstring{\(4+2+4\)}{4+2+4} extraction}\label{m6:s:p006:s2}\label{m6:s:sec:extraction-424}

\begin{proposition}[The Shalika Gauss open]\label{m6:s:prop:extraction-424}\label{m6:rank-dictionary}
Character-weighted averaging in the Gauss chart
\eqref{m6:s:eq:p006-2-2} identifies the two-sided equivariant distributions
on $\Omega'$ with the Shalika distribution problem on
$\GL_4\times\GL_2$. It preserves transpose oddness, has trivial
residual density character, and transfers the NW condition.
Consequently, for every $T\in D'_{\widetilde\eta}(G)\cap\NW(G)$,
\[
 T|_{\Omega'}=0,\qquad
 \Omega'=\mathcal P_1\sqcup\mathcal P_6\sqcup\mathcal P_{34}.
\]
The external-factor descent and the rank dictionary are given in
\S\ref{supp:m6-shalika-coordinates}.
\end{proposition}
\begin{proof}

\paragraph{The maximal Levi and its character factors}\label{m6:s:p006:s2-1}

Put

\[
j=\begin{pmatrix}0&-1\\1&0\end{pmatrix},\quad
J_4=\begin{pmatrix}0&j\\j&0\end{pmatrix},\quad
J_{10}=w_5\otimes j.
\]

Let \(Q'=M'N'\) be the parabolic with block sizes \(4,2,4\), with upper unipotent radical \(N'\). Its opposite big cell is
\[
\Omega'=\left\{g\in G:\det g_{\{7,8,9,10\},\{7,8,9,10\}}\ne0\right\}.
\]
Use the Levi embedding

\begin{equation}
\iota(A,C)=
\operatorname{diag}\bigl(A,C,\det(C)J_4^{-1}A^{-t}J_4\bigr)
\quad (A,C\in GL_4\times GL_2).
\label{m6:s:eq:p006-2-1}
\end{equation}

The big-cell coordinate map is

\begin{equation}
N'\times(GL_4\times GL_2)\times N'\longrightarrow\Omega',
\quad (n_1,(A,C),n_2)\longmapsto n_1\iota(A,C)n_2^t.
\label{m6:s:eq:p006-2-2}
\end{equation}

Block Gaussian decomposition and the rank classification in Section~\ref{m6:s:sec:cells} prove that this is an \(F\)-rational analytic isomorphism and that
\begin{equation}
\Omega'=\mathcal P_1\sqcup\mathcal P_6\sqcup\mathcal P_{34}.
\label{m6:s:eq:p034-1-8}
\end{equation}
Here \(\dim N'=18\). The external factor \(C\) ranges over all of \(GL_2(F)\): its determinant is the similitude multiplier, ranging over $F^\times$.

The residual modulus is trivial, while each radical measure carries the restriction of the fixed character. In the coordinates of Section~\ref{m6:s:sec:setup}, \(N'\) is obtained by setting the adjacent \(X\)-coordinate equal to zero, and

\begin{equation}
\xi_{N'}(n)=\psi(\operatorname{tr}Y_n).
\label{m6:s:eq:p006-2-3}
\end{equation}

The averaging map is therefore

\begin{equation}
\mathcal J'f(A,C)=\int_{N'\times N'}
 f(n_1\iota(A,C)n_2^t)\xi_{N'}(n_1)\xi_{N'}(n_2)
 \,dn_1\,dn_2.
\label{m6:s:eq:p006-2-4}
\end{equation}

The intersection group is

\begin{equation}
H_S=\{(n(X)d(a),a):X\in M_2(F),\ a\in GL_2(F)\},
\quad d(a)=\operatorname{diag}(a,a),
\quad n(X)=\begin{pmatrix}I&X\\0&I\end{pmatrix},
\label{m6:s:eq:p006-2-5}
\end{equation}

with character \(\psi(\operatorname{tr}X)\). The parabolic modulus is

\begin{equation}
\delta_{Q'}(A,C)=|\det A|^7|\det C|^{-14},
\qquad \delta_{Q'}|_{H_S}=1.
\label{m6:s:eq:p006-2-6}
\end{equation}

The modulus follows by writing an arbitrary radical Lie element in \(4,2,4\) blocks as
\[
\begin{pmatrix}0&R&S\\0&0&-j^{-1}R^tJ_4\\0&0&0\end{pmatrix},
\qquad R\in M_{4,2},\quad S^tJ_4+J_4S=0.
\]
There are eight coordinates in \(R\) and ten in \(S\). Conjugation on the first quotient is \(R\mapsto ARC^{-1}\), with determinant \((\det A)^2(\det C)^{-4}\); on the second, symmetric-square quotient its determinant is \((\det A)^5(\det C)^{-10}\). Their product gives \eqref{m6:s:eq:p006-2-6}. If \(A=u(X)d(a)\) and \(C=a\), then \(\det A=(\det a)^2\), so the restriction equals one. For \(R=\binom ZY\), the upper element \(u(X)\) replaces \(Z\) by \(Z+XY\) and fixes \(Y\), while the diagonal element conjugates \(Y\) by \(a\). Thus the radical character is preserved.

It follows that \eqref{m6:s:eq:p006-2-4} leaves exactly the two-sided \(H_S\)-character. Transpose exchanges the two identical radical Haar factors.

The averaging in \eqref{m6:s:eq:p006-2-4} is onto: in product coordinates choose compactly supported functions in each radical variable with weighted integral one. Conversely, translation equivariance in each radical forces the character-weighted Haar functional. Hence the averaging identifies the distribution spaces, compatibly with open restrictions and closed-support injections.

\paragraph{Vanishing on the Gauss open}
To prove vanishing, work on the Levi subgroup $\GL_4\times\GL_2$
before external Frobenius descent. The subgroup in
\eqref{m6:s:eq:p006-2-5} is exactly $R_{\rm Sh}$ of
Theorem~\ref{cf:shalika-nw}: switching $u(X)D_2(a)$ to
$D_2(a)u(a^{-1}Xa)$ preserves the trace character.
The rational Gauss chart, modulus-one graph action, and transpose
exchange proved above give an odd Shalika descendant with $\chi=1$.
Lemma~\ref{cf:gauss} gives its NW condition on the entire Levi.

Explicitly, the Levi Lie algebra embedding is
\[
 (P,Q)\longmapsto
 \diag(P,Q,\tr(Q)I_4-J_4^{-1}P^tJ_4).
\]
Nilpotence forces $P,Q$ nilpotent and $\tr P=\tr Q=0$. Restricting
the ambient trace pairing then represents the Levi cotangent vector by
$(2P,Q)$, so both derived components are nilpotent and both central
components are zero. Smooth weighted Haar factors contribute no
nonzero wavefront direction. Theorem~\ref{cf:shalika-nw} kills the
descendant, and injectivity of \eqref{m6:s:eq:p006-2-4} gives
$T|_{\Omega'}=0$.
\end{proof}

The rank-one stratum supports intrinsic odd distributions. Its
independent conormal analysis and the rank-zero classification are
given in \S\ref{supp:m6-sharpness}; the proof here uses
the Shalika theorem on the entire open.

\subsection{Rational relative positions and the character sieve}\label{m6:sieve}
The corner-rank classification gives ninety rational cells. Of these,
unipotent characters exclude sixty, and full $H\times H$ stabilizers
exclude eight further families. We give both tests below, together
with the complete finite tables.
\subsubsection{Rational parabolic cells and compatibility}\label{m6:s:sec:cells}

The main argument gives the rational rank classification and derives the
complete stabilizer-character test. Appendix~\S\ref{m6:s:app:compatibility}
specifies the complete finite generator and determinant obstructions;
\texttt{anc/model6/} contains all ninety systems, thirty witnesses and
an exact reconstruction verifier.
The open assembly is proved in \S\ref{m6:filtration}.

\paragraph{Corner ranks and the six outer positions}
Number the two-dimensional coordinate blocks from one to five. For $1\leq a,b\leq4$ set
\begin{equation}\label{m6:s:eq:cells-ranks}
R_{ab}(g)=\operatorname{rank}
g_{\{2a+1,\ldots,10\},\,\{2b+1,\ldots,10\}}.
\end{equation}
Multiplication on the left by an upper block triangular matrix preserves each lower row quotient, and multiplication on the right by a lower block triangular matrix preserves each right column quotient. Thus these sixteen ranks are invariant under $P\times P^t$, and hence under $\mathcal A$.

The symplectic complementary-rank identities are
\begin{equation}\label{m6:s:eq:cells-complement}
R_{ab}=10-2a-2b+R_{5-a,5-b}.
\end{equation}
One way to obtain them is to identify the annihilator of each coordinate flag space with its opposite complementary flag space using $J$. The rank of the induced map on the two quotients equals their dimensions minus the dimension of the corresponding annihilator intersection. Applying this statement to the complementary pair interchanges $(a,b)$ with $(5-a,5-b)$ and gives \eqref{m6:s:eq:cells-complement}. Equivalently, the complementary minor identity for an invertible matrix applied to $g^{-1}=\ell(g)^{-1}J^{-1}g^tJ$ gives the same rank equality. In particular the first two rows of $R$ determine the matrix. Transpose sends $R$ to $R^t$.

For the outer $2+6+2$ parabolic it is useful to write
\[
s=\operatorname{rank}g_{55},\qquad
r=\operatorname{rank}\begin{pmatrix}g_{25}\\g_{35}\\g_{45}\\g_{55}\end{pmatrix}.
\]
The corresponding corner matrix is
\[
R_Q(g)=\begin{pmatrix}6+s&r\\r&s\end{pmatrix}.
\]
There are six possibilities:
\[
(r,s)=(2,2),(2,1),(1,1),(2,0),(1,0),(0,0).
\]
Here is a rational description. Let $E$ be the coordinate isotropic two-plane, let $E^-$ be its paired opposite, and put $W=gE^-$. Then
\[
\dim(W\cap E)=2-r,\qquad
\dim(W\cap E^\perp)=2-s.
\]
Choose a Witt basis with $e_1,e_2$ spanning $E$ and
$e_3,e_4,e_5,f_3,f_4,f_5$ projecting to a Witt basis of $E^\perp/E$. Representatives for $W$ in the six cases are, respectively,
\[
\langle f_1,f_2\rangle,\quad
\langle f_1,e_3\rangle,\quad
\langle e_1,f_2\rangle,\quad
\langle e_3,e_4\rangle,\quad
\langle e_1,e_3\rangle,\quad
\langle e_1,e_2\rangle.
\]
The quotient contribution has dimension $r-s$ and is isotropic in $E^\perp/E$. Extending bases of the intersections to adapted Witt bases proves transitivity over $F$ in each outer double cell. The paired quotients carry alternating forms. The case $(2,0)$ is present because this quotient has dimension six.

\paragraph{The intersection table}
Augment the ranks by $R_{0b}=10-2b$, $R_{a0}=10-2a$, and $R_{5b}=R_{a5}=0$. The intersection table $C=(c_{ij})_{1\leq i,j\leq5}$ is obtained by consecutive differences:
\begin{equation}\label{m6:s:eq:cells-table-difference}
c_{ij}=R_{i-1,j-1}-R_{i,j-1}-R_{i-1,j}+R_{ij}.
\end{equation}
Its entries have the following complete description:
\begin{equation}\label{m6:s:eq:cells-tables}
c_{ij}\geq0,\qquad
\sum_jc_{ij}=\sum_ic_{ij}=2,\qquad
c_{ij}=c_{6-i,6-j},\qquad c_{33}\in\{0,2\}.
\end{equation}
For a signed monomial representative, $c_{ij}$ counts the nonzero entries in block $(i,j)$. Conversely,
\[
R_{ab}=\sum_{i>a,\ j>b}c_{ij}.
\]

\begin{proposition}\label{m6:s:prop:cells-rational}
The tables \eqref{m6:s:eq:cells-tables} classify the rational double cells
$P(F)\backslash G/P^t(F)$. Each table has a signed symplectic monomial representative over $F$. There are exactly ninety such cells.
\end{proposition}
\begin{proof}
Consider the two self-dual flags determined by a point of $G$. Their common associated graded space is indexed by the intersections of the two flags; its dimensions are the consecutive differences \eqref{m6:s:eq:cells-table-difference}. Positivity and the row and column sums follow from the fact that each flag quotient has dimension two. Orthogonal complementation pairs the $(i,j)$ quotient with the $(6-i,6-j)$ quotient, so their dimensions agree. The only fixed quotient is $(3,3)$, on which the form is nondegenerate alternating. Its dimension is therefore even.

Conversely, for a table satisfying \eqref{m6:s:eq:cells-tables}, take $c_{ij}$ basis vectors in each paired quotient, and choose their dual vectors in the opposite quotient over $F$. On the central quotient choose a symplectic basis; its dimension is zero or two. Lift these vectors successively through the two filtrations. At each step subtract linear combinations of previously chosen dual vectors to remove unwanted pairings. Each correction solves a linear equation with an already prescribed nonzero pairing, so is defined over $F$. The resulting basis is simultaneously adapted to both flags and has the fixed alternating Gram matrix. Two pairs of flags with the same table are therefore related by rational symplectic transformations preserving the respective flags. In this basis the relative position is represented by a monomial matrix whose entries are signs. This proves existence and completeness over $F$. A multiplier of an arbitrary point is absorbed into one of the Levi factors.

It remains to count the tables. Let their first two rows be $a,b$, each a composition of two into five nonnegative parts. Central symmetry fixes the last two rows as $b^\mathrm{rev},a^\mathrm{rev}$. The middle row is forced to be
\begin{equation}\label{m6:s:eq:cells-middle-row}
c_j=2-a_j-b_j-a_{6-j}-b_{6-j}.
\end{equation}
Its nonnegativity requires $a_3+b_3\in\{0,1\}$.
If this sum is zero, both totals in the paired columns $\{1,5\}$ and $\{2,4\}$ equal two. If the first row contributes $k$ to the first pair, it contributes $2-k$ to the second. Splitting contributions between the two columns in each pair gives
\[
\sum_{k=0}^2\bigl((k+1)(3-k)\bigr)^2=9+16+9=34.
\]
If $a_3+b_3=1$, choose the row containing the central entry and the paired-column total which is one rather than two. There are four symmetric choices. In one such choice, let the central entry lie in row $a$ and the first paired-column total be one. The remaining entry of $a$ lies in that first pair or the second pair. In the first case its two placements, together with the three allocations of row $b$ in the other pair, give $2\cdot3=6$. In the second case its two placements, the two placements of the first-pair entry of $b$, and the two placements of the other entry of $b$ give $2\cdot2\cdot2=8$. Thus this case contributes $4(6+8)=56$. The total is $34+56=90$.
\end{proof}

The dual-basis construction gives explicit rational representatives and rational transformations for any flags with the same table. These classify the $P(F)$--$P(F)^t$ double cells; the $H(F)\times H(F)$ orbit families inside them are studied below.

\paragraph{Representatives, transpose, and rank opens}
Order the matrices $R$ lexicographically decreasing in row-major order, and denote the resulting cells by $\mathcal P_i$, $1\leq i\leq90$. Their identifiers in the tables are P001 through P090. Write $R_i$ for the rank matrix of $\mathcal P_i$ and $W_i$ for its signed representative. The public table specified in \S~\ref{m6:s:app:ledger} lists the first two rank rows, every signed physical column of $W_i$, and the transpose partner. The last two rank rows follow from \eqref{m6:s:eq:cells-complement}, so the table contains the complete invariant.

The signed column convention is $W_ie_k=\epsilon_k e_{a_k}$ when the $k$th entry is $\epsilon_ka_k$. Each $W_i$ has multiplier one. This prescription, the table, and \eqref{m6:s:eq:setup-form} specify each representative as an ten-by-ten matrix. They make the signs in the compatibility equations reproducible. The transpose permutation has forty fixed cells and twenty-five two-element pairs. It is read by transposing $R_i$ and finding the equal rank matrix in the list.

For an index set $I$, put $\mathcal P[I]=\bigcup_{i\in I}\mathcal P_i$ with its subspace topology. Define
\begin{equation}\label{m6:s:eq:cells-opens}
\mathcal O_j=\mathcal P[1,\ldots,j],\qquad
\mathcal S_j=\mathcal O_j\cup\tau(\mathcal O_j).
\end{equation}
To see that every $\mathcal O_j$ is open in $G$, note that a
lexicographic initial segment is upward closed under componentwise
comparison of the finitely many possible rank matrices. It is
therefore the union, over its matrices $R_i$, of the opens
\begin{equation}\label{m6:s:eq:cells-rank-open}
V_i=\bigcap_{a,b}\{g:R_{ab}(g)\geq (R_i)_{ab}\}.
\end{equation}
Each rank lower bound is open, as it is the union of the nonvanishing loci of the relevant minors. Set
\begin{equation}\label{m6:s:eq:cells-boundary}
B_i=\bigcup_{a,b}\{g\in V_i:R_{ab}(g)>(R_i)_{ab}\}.
\end{equation}
This too is open in $G$, and $V_i\setminus B_i=\mathcal P_i$. In particular a cell is closed in its rank open. All $\mathcal S_j$ are transpose-stable.

The first five cells make up the determinant open
$\Omega=\{\det g_{55}\ne0\}=\mathcal O_5$. The terminal rank matrix is
\[
R_{90}=\begin{pmatrix}
6&4&2&0\\4&2&0&0\\2&0&0&0\\0&0&0&0
\end{pmatrix}.
\]
The symbol $W_{90}$ always refers to its signed representative. The different block reversal $w_0$ used in the terminal calculation is related to it by a displayed Levi change in Section~\ref{m6:s:sec:p090}.

\paragraph{The rational coefficient atlas and stabilizers}\label{m6:s:sec:coefficient-atlas}
For a signed representative $w=W_i$, put
\[
I_w=P\cap wP^tw^{-1},\qquad q^\natural=(w^{-1}qw)^t.
\]
The rational classification gives an open quotient map
$P(F)^2\to\mathcal P_i$, $(p,r)\mapsto pwr^t$, with equivalence
\begin{equation}\label{m6:s:eq:atlas-quotient}
(p,r)\sim(pq,r(q^\natural)^{-1}),\qquad q\in I_w(F).
\end{equation}
Since the stabilizer is smooth in characteristic zero, the orbit map
has local analytic sections. For the coefficient coordinates, define
\[
\mathsf l(\mathsf P,\mathsf Q)=m(\mathsf Q\mathsf P,\mathsf Q,I_2).
\]
The identity
\begin{equation}\label{m6:s:eq:atlas-factor}
m(A,B,C)=d_5(C)\mathsf l(B^{-1}A,C^{-1}B)
\end{equation}
gives rational coordinates $P\simeq H\times\mathrm{GL}_2^2$.
This factorization uses the matrix $C$ itself, retaining every
multiplier. Hence four invertible coefficient matrices and two
$H$ factors cover every rational point of every cell.

At $x=lwr^t$, the even stabilizer consists exactly of
\begin{equation}\label{m6:s:eq:atlas-stabilizer}
\left(lql^{-1},\ r(q^\natural)^{-1}r^{-1}\right),\qquad
q\in I_w(F),\quad
lql^{-1},\ r(q^\natural)^{-1}r^{-1}\in H(F).
\end{equation}
Indeed the stabilizer equation forces the second component from the
first, and the two membership conditions are also sufficient.
We use this right inverse-transpose convention in the projected
gauge and root calculations. Replacing $q$ by $q^{-1}$ reverses
the orientation of a gauge equation.

If a closed coefficient locus has saturated preimage under the above
quotient, its image is closed inside the cell: the complementary
saturated open maps to an open under the quotient. Smoothness follows from the local submersion charts. In each
application below, we establish saturation and the coefficient
equations throughout the invertible domain, retaining the exceptional
pivots, linked diagonal blocks, and opposite-root factors.

\paragraph{The unipotent compatibility criterion}
Every rational point of $\mathcal P_i$ can be written
\begin{equation}\label{m6:s:eq:cells-param}
g=u_0lW_ir^tv_0^t,\qquad u_0,v_0\in U,\quad l,r\in L.
\end{equation}
To avoid conflict with the parabolic $P$, we use upright coefficient matrices
$\mathsf P,\mathsf Q,\mathsf R,\mathsf S$. For $l=m(A_l,B_l,C_l)$ and
$r=m(A_r,B_r,C_r)$ they are
\begin{equation}\label{m6:s:eq:cells-coefficients}
\mathsf P=B_l^{-1}A_l,\quad
\mathsf Q=C_l^{-1}B_l,\quad
\mathsf R=B_r^{-1}A_r,\quad
\mathsf S=C_r^{-1}B_r.
\end{equation}
These are four independent invertible matrices. Conversely, the rational choices
$l=m(\mathsf Q\mathsf P,\mathsf Q,I)$ and
$r=m(\mathsf S\mathsf R,\mathsf S,I)$ realize any such quadruple.

For a signed representative $W_i$ put
\[
\mathfrak u_i=\mathfrak u\cap\operatorname{Ad}(W_i)\mathfrak u^t,\qquad
Z'=-(W_i^{-1}ZW_i)^t\quad(Z\in\mathfrak u_i).
\]
Then $ZW_i+W_i(Z')^t=0$, and
\[
\left(l\exp(tZ)l^{-1},\,r\exp(tZ')r^{-1}\right)
\]
is an $U\times U$ stabilizer element of $lW_ir^t$ for every $t\in F$.
Its character is $\psi(tF_i(Z))$, where
\begin{equation}\label{m6:s:eq:cells-phase}
F_i(Z)=
\operatorname{tr}(\mathsf P Z_{12}+\mathsf Q Z_{23})
+\operatorname{tr}(\mathsf R Z'_{12}+\mathsf S Z'_{23}).
\end{equation}
The signs follow from the inverse transpose in $Z'$ and the plus sign in the original character.

\begin{proposition}\label{m6:s:prop:cells-compatibility}
The character on the unipotent stabilizer at \eqref{m6:s:eq:cells-param} is trivial if and only if
\begin{equation}\label{m6:s:eq:cells-compatible-system}
F_i|_{\mathfrak u_i}=0,\qquad
\det\mathsf P\det\mathsf Q\det\mathsf R\det\mathsf S\ne0.
\end{equation}
\end{proposition}
\begin{proof}
Conjugation by $u_0,v_0$ does not change a character of $U$, so it suffices to use $lW_ir^t$. The first projection identifies its unipotent stabilizer with
$l(U\cap W_iU^tW_i^{-1})l^{-1}$. This unipotent algebraic group is connected in characteristic zero, and the finite exponential and logarithm identify it with its Lie algebra. Formula \eqref{m6:s:eq:cells-phase} is the character on every one-parameter subgroup; it therefore determines the character on the whole group. Conversely, a nonzero value $F_i(Z)$ gives a nontrivial character on the displayed rational additive subgroup.
\end{proof}

For independent reproduction of the finite systems, take the twenty-two basis elements obtained by placing $E_{11},E_{12},E_{21},E_{22}$ in each of $X,Y,Z,W$, and $E_{12},E_{21},E_{11}-E_{22}$ in each of $T_0,S_0$ in \eqref{m6:s:eq:setup-N}. Their dependent entries are fixed by the negative sharp operation. They are root vectors for a diagonal torus. Signed monomial conjugation permutes those root spaces. Retain precisely the basis elements for which $W_i^{-1}ZW_i$ is strictly lower block triangular. Substituting each retained element in \eqref{m6:s:eq:cells-phase} gives the complete linear system in the sixteen scalar entries of the four coefficient matrices. Appendix~\S\ref{m6:s:app:compatibility} specifies the complete certificate, whose row-reduced equations for every cell are in \texttt{anc/model6/cell\_compatibility.json}. This basis-selection rule and the signed table specify all entries of each system.

\paragraph{Sixty whole-cell vanishings}
Before imposing invertibility, the equations in \eqref{m6:s:eq:cells-compatible-system} define a linear subspace. Its invertible locus is nonempty precisely when none of the four determinant polynomials vanishes identically on that subspace. The forward implication is immediate. For the converse their product is a nonzero polynomial, and an infinite field has a rational point off its zero locus. The certificate of Appendix~\S\ref{m6:s:app:compatibility} gives and checks an explicit rational invertible quadruple for each nonempty system.

The counts by outer position are
\[
\begin{array}{c|r|r|r}
(r,s)&\text{cells}&\text{forced incompatible}&\text{invertible solutions}\\\hline
(2,2)&5&0&5\\
(2,1)&23&11&12\\
(1,1)&5&2&3\\
(2,0)&29&23&6\\
(1,0)&23&20&3\\
(0,0)&5&4&1\\\hline
\text{total}&90&60&30.
\end{array}
\]
For example, $\mathcal P_{10}$ gives $q_{11}=q_{21}=r_{12}=r_{22}=0$. Thus the first column of $\mathsf Q$ and the second column of $\mathsf R$ vanish, contrary to invertibility. Each incompatible row of the complete electronic table supplies a determinant contradiction of this kind.

\begin{proposition}\label{m6:s:prop:incompatible}
Let
\[
\begin{split}
I_{\mathrm{inc}}=\{&
10,11,13,15,19,20,21,23,25,26,28,29,30,35,36,37,38,39,\\
&43,44,45,46,48,49,50,51,52,54,55,56,57,58,59,60,61,62,\\
&64,65,66,67,68,70,71,72,73,74,75,76,77,78,79,80,81,82,83,\\
&85,86,87,88,89\}.
\end{split}
\]
For each $i\in I_{\mathrm{inc}}$, $D'_\eta(\mathcal P_i)=0$, without imposing a transpose sign.
\end{proposition}
\begin{proof}
For each listed cell, the reconstructed systems and determinant
identities of Appendix~\S\ref{m6:s:app:compatibility}, recorded in
\path{anc/model6/compatibility_tables.md}, force at least one
coefficient determinant to vanish. Hence at every rational point its unipotent stabilizer character is nontrivial. Proposition~\ref{cf:stabilizer}, applied to the unipotent action, gives the conclusion. The displayed basis rule proves that these equations exhaust the unipotent stabilizer; the determinant contradictions are exact polynomial identities over $F$.
\end{proof}

An invertible solution identifies a nonempty unipotently compatible locus. Its full $H\times H$ stabilizers impose further conditions, and its orbit families have continuous parameters. Indeed $\dim\mathcal A=52$ and the subgroup $(zI_{10},z^{-1}I_{10})$ acts trivially, so every orbit has dimension at most $51$, whereas $\dim G=56$. On the open Levi slice the residual action already contains simultaneous conjugation on two invertible two-by-two matrices. The remaining sections therefore study families of rational orbits inside the thirty surviving parabolic cells.

\subsubsection{Eight stabilizer-character exclusions}
\label{m6:s:sec:full-stabilizer-exclusions}
\label{m6:s:sec:intrinsicpackets}

We apply the stabilizer criterion to explicit additive subgroups.
The equations cited in the table are proved in
\S\ref{m6:s:app:stabilizers}, which gives the complete rational
coefficient families, both group memberships, nonzero character
coefficients, and relative modular characters. These calculations
supply the inputs for the vanishing argument below.

The initial $U\times U$ sieve leaves nonempty compatible loci on the
eight cells below. Their $H\times H$ stabilizers contain the additional
additive subgroups recorded in Table~\ref{m6:s:tab:full-stabilizer-exclusions}.
The local coefficient names in different rows are independent.

\begin{theorem}\label{m6:s:thm:full-stabilizer-exclusions}\label{m6:full-stabilizer}
\label{m6:s:prop:intrinsicpackets}
For $i\in\{18,24,32,41,42,47,63,84\}$,
\[
D'_\eta(\mathcal P_i)=0.
\]
This is an assertion about all ordinary $\eta$-equivariant distributions
on each entire rational cell.
\end{theorem}

For the six primary subgroup calculations, use the row's coefficient
frame $x=lwr^t$ and put
\begin{equation}\label{m6:s:eq:full-exclusion-subgroup}
\Gamma_x(N;t)=
\left(l(I+tN)l^{-1},\;
r\bigl(I-t(w^{-1}Nw)^t\bigr)r^{-1}\right),\qquad t\in F.
\end{equation}
Since every listed $N$ is square-zero, the two factors satisfy the
additive parameter law, with matrix identity
\[
(I+tN)w\bigl(I-t(w^{-1}Nw)^t\bigr)^t=w.
\]
The cited calculations prove both $H$ memberships for all allowed
coefficients, placing one factor in $H_0$ and the other in
$U$. These additional $H\times H$ stabilizer directions give a
second character test after the double-unipotent sieve.

\begin{longtable}{@{}>{\raggedright\arraybackslash}p{.055\linewidth}>{\raggedright\arraybackslash}p{.255\linewidth}>{\raggedright\arraybackslash}p{.365\linewidth}>{\raggedright\arraybackslash}p{.240\linewidth}@{}}
\caption{All-parameter additive subgroup certificates. In the last
column $n$ means that the stabilizer is
$\mathbb G_m\times U_n$ and its relative modular and orbital density
characters are both one.}\label{m6:s:tab:full-stabilizer-exclusions}\\
\hline
Cell & Complete coefficient family; mandatory units & Subgroup for every $t\in F$ & Phase $\psi(c(x)t)$; orbital density \\
\hline
\endfirsthead
\hline
Cell & Complete coefficient family; mandatory units & Subgroup for every $t\in F$ & Phase $\psi(c(x)t)$; orbital density \\
\hline
\endhead
18 & \eqref{m6:s:eq:packets-p18-1}; $d,e,f,h,k,ad-bc\ne0$
& $\Gamma_x(R_{18};t)$, $R_{18}=S_f-(d/k)Y_{21}-(hd/(fk))W_{12}$; its right factor is $d_5(I+tE_{21})$, by \eqref{m6:s:eq:packets-p18-4}.
& $c=-de/k$; $n=8$, \eqref{m6:s:eq:packets-f-9}--\eqref{m6:s:eq:packets-f-10}. \\[3pt]
24 & Rational epsilon image \eqref{m6:s:eq:packets-t-2} of the complete 18-family; $r_{22},s_{12},q_{22}\ne0$
& Transpose exchanges the factors of the preceding subgroup; the partner frame is exactly \eqref{m6:s:eq:packets-t-3}.
& $c=-r_{22}s_{12}/q_{22}$, \eqref{m6:s:eq:packets-t-4}; $n=8$, with the same density. \\[3pt]
32 & \eqref{m6:s:eq:sector-32-family}; $b,c,e,f,g,i\ne0$
& $\Gamma_x(E_0;t)$, $E_0=\frac{cf}{bg}A_{21}+\frac fgB_{12}+C_{21}$; its left factor is $d_5(I+tE_{21})$, by \eqref{m6:s:eq:sector-32-subgroup}.
& $c=fi/g$; $n=13$, Section~\ref{m6:s:sec:sector-p032}. \\[3pt]
41 & \eqref{m6:s:eq:late-a41-1}; $a,c,d,f,h,i,k\ne0$
& $\Gamma_x(D;t)$, $D=\frac{ad}{cf}A_{21}+\frac dfB_{21}+C_{21}$; its left factor is $d_5(I+tE_{21})$, by \eqref{m6:s:eq:late-a41-11}.
& $c=-adi/(cf)$; $n=9$, \eqref{m6:s:eq:late-a41-9} and \S\ref{m6:s:sec:exclusion-density-41}. \\[3pt]
42 & \eqref{m6:s:eq:packets-p42-1}; $a,c,d,f,h,j\ne0$
& $\Gamma_x(L_{42};t)$, $L_{42}=L_C+\frac dfL_B+\frac{ad}{cf}L_A$; its left factor is $d_5(I+tE_{21})$, by \eqref{m6:s:eq:packets-p42-5}.
& $c=-had/(cf)$; $n=10$, \eqref{m6:s:eq:packets-f-9}--\eqref{m6:s:eq:packets-f-10}. \\[3pt]
47 & Rational epsilon image \eqref{m6:s:eq:packets-t-2} of the complete 42-family; $q_{11},r_{11},s_{11},r_{22},s_{22}\ne0$ by lower triangular invertibility
& Transpose exchanges the factors of the 42-subgroup; use the partner frame \eqref{m6:s:eq:packets-t-3}, including its signed fifth block.
& $c=-q_{11}r_{11}s_{11}/(r_{22}s_{22})$, \eqref{m6:s:eq:packets-t-5}; $n=10$, with the same density. \\[3pt]
63 & \eqref{m6:s:eq:late-63-a}; $a,d,e,f,h,j\ne0$
& $\Gamma_x(L;t)$, $L=\frac daA_{12}+B_{12}+\frac feC_{21}$; its left factor is $d_5(I+(f/e)tE_{21})$, by \eqref{m6:s:eq:late-63-h}.
& $c=jd/a$; $n=14$, \eqref{m6:s:eq:late-63-f} and \S\ref{m6:s:sec:exclusion-density-63}. \\[3pt]
84 & \eqref{m6:s:eq:late-84-a}; $a,c,d,e,g\ne0$
& $\Gamma_x(L;t)$, $L=\frac{cd}{ae}A_{12}+\frac deB_{12}+C_{21}$; its left factor is $d_5(I+tE_{21})$, by \eqref{m6:s:eq:late-84-k}.
& $c=cdg/(ae)$; $n=19$, \eqref{m6:s:eq:late-84-i} and \S\ref{m6:s:sec:exclusion-density-84}. \\
\hline
\end{longtable}

\begin{proof}[Proof of Theorem~\ref{m6:s:thm:full-stabilizer-exclusions}]
The family equations in each row cover the entire rational
$U\times U$-compatible locus, with the original similitude factors
retained by the Levi factorization. The group calculations
solve both projected membership equations before applying finite
logarithm and exponential. Their scalar parameter is a matrix
entry; the resulting central-times-unipotent descriptions include all
components and all allowed parameter specializations.

Each displayed coefficient $c(x)$ is nonzero by the mandatory units.
Hence $t\mapsto c(x)t$ is onto $F$, and $\psi(c(x)t)$ is a nontrivial
character of the entire additive subgroup, regardless of its conductor
or of the valuation of $c(x)$. The two transpose-derived rows use the
rational epsilon substitution and preserve $\eta$.

On each compatible locus, the relative modular character is one
and the subgroup character is nontrivial. They therefore differ in
\eqref{cf:orbital}. At every initially incompatible point,
the established sieve already supplies an unipotent stabilizer
subgroup with nontrivial character. The relative modular character
also restricts to one on that subgroup: conjugation is unipotent on the
relevant Lie algebras and their stable quotients. Thus
Proposition~\ref{cf:stabilizer}, with its rational-point algebraic orbit
assembly, applies at every point of each cell.
\end{proof}

The detailed certificates in \S\ref{m6:s:app:stabilizers}
retain the complete group equations, signed representatives, rational
coverage, and separate extraction determinants. In particular the
additional phases $\psi(fv)$ for cell 63 and
$\psi(-ade\,v/(cg))$ for cell 84 remain in their stabilizer-character formulas.
This theorem supplies local cell vanishings; the open support
order remains the 23-step order of Section~\ref{m6:s:sec:assembly}.

\subsection{Boundary arguments and intrinsic orbital families}\label{m6:packets}
We prove the boundary steps in their support order: the intrinsic
incidence arguments on $7,8$, the two conditional spectral pairs,
and the remaining orbital families.
\subsubsection{\texorpdfstring{The framed incidence cells $\mathcal P_7$ and $\mathcal P_8$}{The framed incidence cells P7 and P8}}\label{m6:s:sec:p007-p008}

Both cells lie outside the determinant opens for the $2+6+2$ and
$4+2+4$ coarsenings. Their common framed $Sp_4$ incidence slice
has two distinct residual mechanisms: an entire opposite-weight
plane for $\mathcal P_7$, and a stabilizer-character condition on
the boundary of $\mathcal P_8$. The rational normalizations,
complete stabilizers and density calculations are given in
\S\ref{supp:m6-incidence-calculations}.

\begin{proposition}\label{m6:s:prop:p007-p008}
The intrinsic odd equivariant distribution spaces on \(\mathcal P_7\) and \(\mathcal P_8\) vanish. Their union \(\mathcal O_8\setminus\mathcal O_6\), with its induced topology, has the same property. Consequently, any original global \(\eta\)-equivariant, transpose-odd distribution \(T\) satisfying \(T|_{\mathcal O_6}=0\) also satisfies \(T|_{\mathcal O_8}=0\).
\end{proposition}

\paragraph{The residual correspondence}
Let $\mathcal A=H\times H$, $\mathcal U=U\times U$ and
$G_0=H_0\times H_0\simeq\GL_2\times\GL_2$.
For $i=7,8$, the rational parabolic factorization identifies
$\mathcal P_i$ with $(P\times P)/\mathcal J_i$, where
$\mathcal J_i$ is the graph of $P\cap W_iP^tW_i^{-1}$.
Its image $\mathcal R_i$ in $L\times L$ gives the
residual space $\mathscr Y_i=(L\times L)/\mathcal R_i$.
Integration along the unipotent fibers gives a $G_0$-equivariant
correspondence valued in the unipotent orbital-distribution line.
It retains that line's relative modular character and commutes
with open restriction and closed-support injection.

Write the normalized Levi pair as
$l=m(A,B,I_2)$, $r=m(D,E,I_2)$, and put
$\mathsf P=B^{-1}A$, $Q=B$, $R=E^{-1}D$, $S=E$.
The kernel-character equations are
\begin{equation}\label{m6:s:eq:p78-3-15}
\begin{array}{ll}
\mathcal P_7:&\mathsf P_{12}=R_{12},\quad Q_{12}=S_{12},\\
\mathcal P_8:&\mathsf P_{12}=R_{12},\quad Q_{12}=-S_{22},\quad Q_{22}=S_{12}.
\end{array}
\end{equation}
Off this closed locus a one-parameter stabilizer has nontrivial
character. All normalizations and quotient maps used here have
rational local sections, including at the stabilizer jumps;
\S\ref{supp:m6-incidence-calculations} also records their phases
and positive density factors.

\paragraph{The opposite-weight plane on $\mathcal P_7$}
On the compatible locus set
$\wp=(B^{-1}A)_{12}=(E^{-1}D)_{12}$ and $q=B_{12}=E_{12}$.
The complete stabilizer-character condition is
\begin{equation}
{\text{full-compatible }\mathcal P_7\text{ locus}
=\{\mathsf P_{12}=R_{12},\ Q_{12}=S_{12},\ (\mathsf P_{12},Q_{12})\ne(0,0)\}.}
\label{m6:s:eq:p78-5-13}
\end{equation}
For $q\ne0$, normalize $B=E=f$, where
$f=\left(\begin{smallmatrix}0&1\\1&0\end{smallmatrix}\right)$,
and write $p=A_{22}=D_{22}$.
On $p\ne0$, the remaining rational normal forms have passive
parameters $(p,k)\in(F^\times)^2$ and two affine parameters $(b,e)$.
The effective residual action is
\begin{equation}
h:(p,k,b,e)\longmapsto(p,k,hb,h^{-1}e),\qquad h\in F^\times.
\label{m6:s:eq:p78-6-6}
\end{equation}
Transpose swaps $b,e$. Frobenius reciprocity identifies
the intrinsic odd distribution problem with
\begin{equation}\label{m6:s:eq:p78-6-15}
D'((F^\times)^2\times F^2)^{F^\times,\ \tau=-1},\quad
h(b,e)=(hb,h^{-1}e),\quad\tau(b,e)=(e,b).
\end{equation}
The orbital line on this plane is untwisted, and the associated
bundle remains smooth at $(b,e)=(0,0)$. In particular, both
punctured axes have the included origin in their closure inside
the same cell. Lemma~\ref{m8:s:c:plane} applies to the entire
plane, with passive $(p,k)$, and excludes an invariant extension
of the odd difference of the two axis measures. On the other
compatible branches, $q\ne0,p=0$ and $q=0,\wp\ne0$, the rational
normal forms have character-one transpose transporters.

\begin{proposition}\label{m6:s:prop:p007}One has \(D'_{\widetilde\eta}(\mathcal P_7)=0\).\end{proposition}
\par\noindent\textit{Proof.}

Remove the unipotent-incompatible open first, using
Proposition~\ref{cf:stabilizer}. On its closed compatible antecedent,
remove $q\ne0$ in two steps: the open $p\ne0$ branch vanishes by
\eqref{m6:s:eq:p78-6-15} and Lemma~\ref{m8:s:c:plane}; the closed
$p=0$ branch vanishes by \eqref{m6:s:eq:p78-6-18},
\eqref{m6:s:eq:p78-5-5} and Lemma~\ref{m8:s:c:transporter}.
On the remaining $q=0$ antecedent, the open $\wp\ne0$ branch
vanishes by \eqref{m6:s:eq:p78-6-22}; the final closed $q=\wp=0$
layer has nontrivial stabilizer character \eqref{m6:s:eq:p78-5-12}.
The left-exact support sequence \eqref{cf:support-sequence}
therefore gives, including both stabilizer jumps,
\begin{equation}\label{m6:s:eq:p78-7-2}
D'_{\widetilde\eta}(\mathcal P_7)=0.
\end{equation}
\hfill\(\square\)\par

\paragraph{The stabilizer jump on $\mathcal P_8$}
Every unipotent-compatible point admits a zero-phase rational
normalization $B=I_2$, $E=f$, $A_{12}=D_{22}=p$.
Using the whole invertible matrix $B$ gives this normalization
even if a component of $Be_2$ vanishes. The remaining action is
\begin{equation}
{
A\mapsto C_{h,z}^{-1}A\operatorname{diag}(r,1),
\qquad
D\mapsto C_{h,z}^tD\operatorname{diag}(r^{-1},h^{-1}),
\quad C_{h,z}=\begin{pmatrix}h&0\\z&1\end{pmatrix},
}
\label{m6:s:eq:p78-8-4}
\end{equation}
where $h,r\in F^\times$ and $z\in F$. Its fiber density is
$d_8=|h|^{-3}$ and its dual relative density is $|h|^3$;
point stabilizers have $h=1$, and transpose restoration has $h=-1$.
The $p\ne0$ branch has the complete rational parameters
$(k,m)\in F^\times\times F$. The fixed character-one translate
$\theta(g)=d_5(f)g^td_5(f)$ fixes each normal representative
$X_8(k,m)$ of \eqref{m6:s:eq:p78-9-2}. Its complete
stabilizer is unimodular, and its
positive orbital measure is $\theta$-even.

On $p=0$, the complete ordered parameters are
$(a,b)\in(F^\times)^2$. An additional unipotent stabilizer
has character
\begin{equation}
{
\eta\bigl(\exp(\sum_{j=1}^4t_jK_j+zV)\bigr)
=\psi\!\left(z\left(\frac1a-\frac1b\right)\right).
}
\label{m6:s:eq:p78-9-15}
\end{equation}
where the $K_j$ span the character-trivial graph kernel.
Thus its compatibility condition is exactly $a=b$.
The involution $\theta$ swaps these ordered parameters;
on $a=b$ its transporter has character one and preserves
the positive orbital line. The full rational normalizations and
these measure comparisons hold over dyadic fields as well:
the positive Jacobian $|2|$ is retained.

\begin{proposition}\label{m6:s:prop:p008}
Every intrinsic odd equivariant distribution on \(\mathcal P_8\) is zero. The complete rational orbit parameters are \((k,m)\in F^\times\times F\) on the nonzero branch and the ordered \((a,b)\in(F^\times)^2\) on its boundary; stabilizer-character compatibility on that boundary is exactly \(a=b\).
\end{proposition}
\begin{proof}
Remove the unipotent-incompatible open. On its closed compatible
antecedent, the open $p\ne0$ branch vanishes by localization and
the character-one involution just described. On the remaining
closed $p=0$ layer, the open $a\ne b$ part vanishes by the
nontrivial stabilizer character; its closed complement $a=b$
vanishes by Lemma~\ref{m8:s:c:transporter}.
The support sequence gives
\begin{equation}\label{m6:s:eq:p78-9-21}
D'_{\widetilde\eta}(\mathcal P_8)=0.
\end{equation}
The rational normalizers transport this conclusion to every point;
their own characters may be nontrivial, but the transported
transpose transporter has character one by
Lemma~\ref{m6:s:lem:normal-transport}.
\end{proof}

\begin{proof}[Proof of Proposition~\ref{m6:s:prop:p007-p008}]
The rank-open filtration is
$\mathcal O_7=\mathcal O_6\sqcup\mathcal P_7$ and
$\mathcal O_8=\mathcal O_7\sqcup\mathcal P_8$.
In $\mathcal O_8\setminus\mathcal O_6$, the first cell is open
and the second is its closed complement. The two intrinsic
vanishing results and the support sequence therefore give vanishing
on that union. It is closed in $\mathcal O_8$ with open complement
$\mathcal O_6$, proving the asserted implication for the same
original global distribution.
\end{proof}

\subsubsection{The spectral pair \texorpdfstring{$\mathcal P_9,\mathcal P_{14}$}{P9 and P14}}
\label{m6:s:sec:p009-p014}

This pair has intrinsic odd orbital distributions. Its full
stabilizer-compatible locus is not locally closed, so the original
global NW distribution is excluded on two successive smooth
support manifolds. The rational orbit classification and complete
group calculations are in \S\ref{supp:m6-p009-calculations}.

\paragraph{Compatible families and their two branches}
The upright letters $\mathsf P,\mathsf Q,\mathsf R,\mathsf S$
denote two-by-two coefficient matrices. Put
$\ell_0(\mathsf P,\mathsf Q)=m(\mathsf Q\mathsf P,\mathsf Q,I_2)$,
and use representatives
$x=\ell_0(\mathsf P,\mathsf Q)W_9\ell_0(\mathsf R,\mathsf S)^t$.
Normalizing the central Levi blocks by $H_0$ retains every
similitude. The unipotent-compatible coefficient family is
\begin{equation}\label{m6:s:eq:p009-family}
 \mathsf Q=\begin{pmatrix}q&0\\u&v\end{pmatrix},\qquad
 \mathsf R=\begin{pmatrix}r&0\\w&q\end{pmatrix},\qquad
 qvr\det\mathsf P\det\mathsf S\ne0,
\end{equation}
Let $p=p_{12}$, $s=s_{12}$ and $F_0=rp_{11}+sp_{22}$.
The full stabilizers are central $\mathbb G_m$ times connected
unipotent groups, with relative density one. Their remaining
character is described by the following table, where $z,t$ are the
extra additive stabilizer parameters:
\begin{equation}\label{m6:s:eq:p009-character-table}
\begin{array}{c|c|c}
\text{stratum}&\dim\mathcal A_x&\text{additional phase}\\ \hline
p\ne0,\ s\ne0&4&0\\
p\ne0,\ s=0&5&-(s_{11}/s_{22})p\,t\\
p=0,\ s\ne0&5&-qF_0z/(vp_{22})\\
p=s=0&6&-qrp_{11}z/(vp_{22})
\end{array}
\end{equation}
The nonzero coefficients on the second and fourth rows are forced
by invertibility. Therefore the stabilizer-compatible locus is
\begin{equation}\label{m6:s:eq:p009-full-compatible}
 \{ps\ne0\}\ \sqcup\ \{p=0,\ F_0=0\}.
\end{equation}
In its second branch, $s\ne0$ follows from invertibility.
This union contains the dense open $ps\ne0$ but is not open
at a point of $p=F_0=0$; it is consequently not locally closed.

\begin{proposition}\label{m6:s:prop:p009-rational-orbits}
On the character-zero locus, the rational orbits consist precisely
of the following two disjoint families:
\[
\begin{array}{ll}
x_{\rho,\sigma}=
\ell_0\left(\begin{pmatrix}0&1\\\sigma&0\end{pmatrix},I\right)
W_9\ell_0\left(I,\begin{pmatrix}0&1\\\rho&0\end{pmatrix}\right)^t,
&\rho,\sigma\in F^\times,\\[4pt]
y_\kappa=
\ell_0(\operatorname{diag}(-1,1),I)
W_9\ell_0\left(I,\begin{pmatrix}0&1\\\kappa&0\end{pmatrix}\right)^t,
&\kappa\in F^\times.
\end{array}
\]
The first family lies on \(ps\ne0\), with complete invariants
\begin{equation}\label{m6:s:eq:p009-generic-invariants}
 \rho=-\frac{qvp\det\mathsf S}{s},\qquad
 \sigma=-\frac{r\det\mathsf P\,s}{p^2}.
\end{equation}
The second lies on \(p=F_0=0\), with complete invariant
\begin{equation}\label{m6:s:eq:p009-exceptional-invariant}
 \kappa=-qv p_{22}\det\mathsf S.
\end{equation}
The respective stabilizers are
\(\mathbb G_m\times\mathbb G_a^3\) and
\(\mathbb G_m\times\operatorname{Heis}_3\times\mathbb G_a\).
\end{proposition}
Lifting the rational reductions to $H(F)^2$ proves completeness of
these invariants over $F$. To see that each displayed orbit is
closed in the cell, use the saturated equations
\eqref{m6:s:eq:p009-closed-generic} and
\eqref{m6:s:eq:p009-closed-exceptional}, which force the required
nonzero coefficients.
Because the homogeneous quotient topology agrees with the subspace
topology, the orbital integral $\mathscr T_x$ of
\eqref{m6:s:eq:p009-orbital-distribution} is finite on compact
test supports and has the required equivariance.
Transpose swaps the two cells and transports their characters
and positive measures.

\begin{proposition}\label{m6:s:prop:p009-intrinsic}
There are nonzero ordinary transpose-odd \(\eta\)-equivariant
distributions on \(\mathcal P_9\sqcup\mathcal P_{14}\).
\end{proposition}
\begin{proof}
Take the difference $\mathscr T_{x_{1,1}}-\tau_*\mathscr T_{x_{1,1}}$.
It has the required equivariance and parity, with disjoint closed
orbits as supports. A compact-open test function near a
phase-one point in the first cell gives a positive first integral
and zero second integral; hence the difference is nonzero.
The same holds for $y_1$.
\end{proof}

\paragraph{The two support manifolds and their ambient conormals}

Let $\mathcal C_9$ be the $H(F)^2$-saturation of the coefficient
family \eqref{m6:s:eq:p009-family}. It is the closed, smooth
unipotent-compatible envelope of dimension fifty; its rational
quotient charts are established in \S\ref{supp:m6-p009-calculations}.

The invariant equations \(p=0\) and \(s=0\) descend to closed
subsets \(C_p,C_s\) of \(\mathcal C_9\). Define \(\mathcal C_9^\circ=\mathcal C_9\setminus(C_p\cup C_s),\qquad \mathcal E_9=\{p=0,\ rp_{11}+sp_{22}=0\}\subset\mathcal C_9.\)
The first is smooth of dimension fifty. The second is closed
and smooth of dimension forty-eight: in the coefficient space,
\(p=0\) is an independent equation, and on that divisor the
derivative of the second equation with respect to \(s\) is
\(p_{22}\ne0\). The same smooth quotient subtracts fifteen
fiber dimensions. It consists exactly of the \(y_\kappa\)-orbits.

We use right translation and the invariant trace pairing on $\mathfrak{gsp}_{10}$,
\begin{equation}\label{m6:s:eq:p009-trace-trivialization}
 \lambda_Z(\dot x)=\operatorname{tr}(Z\dot x x^{-1}),\qquad
 Z\in\mathfrak{gsp}_{10}.
\end{equation}
The ambient space is the fifty-six-dimensional
\(\mathfrak{gsp}_{10}=\mathfrak{sp}_{10}\oplus FI_{10}\).

\begin{proposition}\label{m6:s:prop:p009-conormals}
The ambient conormals of \(\mathcal C_9^\circ\) and
\(\mathcal E_9\) have dimensions six and eight, respectively.
Their complete bases are \eqref{m6:s:eq:p009-generic-conormal} and
\eqref{m6:s:eq:p009-exceptional-conormal} in
\S\ref{m6:s:cert:p009-conormals}. Both contain the uniform witness
\[
 D_9=\operatorname{diag}(0,0,1,0,-1,1,0,-1,0,0),\qquad
 \operatorname{tr}(D_9^2)=4\ne0.
\]
At a saturated point \(h_1xh_2^t\) the witness is
\(\operatorname{Ad}_{h_1}D_9\).
\end{proposition}

\begin{proof}
To compute the tangent spaces, add the full parameter derivatives
to the two action images:
\[
\{Ax+xB^t\}+F\partial_\rho x+F\partial_\sigma x,\qquad
\{Ay+yB^t\}+F\partial_\kappa y.
\]
Their dimensions are $52-4+2=50$ and $52-5+1=48$.
\S\ref{m6:s:cert:p009-conormals} pairs all six/eight basis vectors
against both actions and every parameter derivative, including
the scalar tangent $I_{10}$ ($\operatorname{tr}Z=0$).
Independence and these codimensions prove completeness.
The witness has polynomial $t^6(t-1)^2(t+1)^2$ and is
nonnilpotent. Conjugation transports it to every rational point
and multiplier.
\end{proof}

Transpose transports the full conormal fiber.
If \(v=\dot x x^{-1}\), then its new right-trivialized tangent is
\((x^{-1}vx)^t\). The trace-dual formula is consequently
\begin{equation}\label{m6:s:eq:p009-cotangent-transpose}
 Z\longmapsto(x^{-1}Zx)^t,
\end{equation}
because the trace of its product with \((x^{-1}vx)^t\) is
\(\operatorname{tr}(Zv)\). Applying
\eqref{m6:s:eq:p009-cotangent-transpose} to every basis vector in
\eqref{m6:s:eq:p009-generic-conormal} and \eqref{m6:s:eq:p009-exceptional-conormal}
gives the entire partner conormal in \(\mathcal P_{14}\);
characteristic polynomials and nilpotence are preserved.

\paragraph{The two-stage spectral exclusion}

For \(i=9,14\), let \(V_i\) be the rank-open envelope defined by
\(R_{ab}(g)\ge(R_i)_{ab}\) for all \(a,b\), and put
\(B_i=V_i\setminus\mathcal P_i\). Each rank-at-least condition is
an open union of nonvanishing minor loci. The rank classification
gives \(B_i=\bigcup_{a,b}\{g\in V_i:R_{ab}(g)>(R_i)_{ab}\};\)
thus \(B_i\) is open and \(\mathcal P_i\) is closed in
\(V_i\). Put \(V_{9,14}=V_9\cup V_{14},\qquad B_{9,14}=B_9\cup B_{14}.\)
Incomparability of the two rank matrices gives
\(V_{9,14}\setminus B_{9,14}=\mathcal P_9\sqcup\mathcal P_{14}\).

\begin{proposition}\label{m6:s:prop:p009-spectral}
Let \(T\in D'_\eta(G)\cap\NW(G)\). Then \(T|_{B_i}=0\ \Longrightarrow\ T|_{V_i}=0\quad(i=9,14).\)
Consequently \(T|_{B_{9,14}}=0\) implies \(T|_{V_{9,14}}=0\).

\end{proposition}

\begin{proof}
Assume $T|_{B_9}=0$. The open--closed support sequence gives a
unique equivariant cell antecedent whose closed pushforward is
$T|_{V_9}$. Its open unipotent-incompatible part vanishes by
\cite[Theorem 2.5 and Remark 2.6]{JLX}; another closed-support
injection gives $\operatorname{supp}(T|_{V_9})\subset\mathcal C_9$.

The subset $D_\partial=C_p\cup C_s$ is closed in $V_9$.
On the ambient open $M=V_9\setminus D_\partial$, the restriction
of the same original $T$ is supported on the closed smooth
$\mathcal C_9^\circ$. Its six-dimensional ambient conormal has
a nonnilpotent vector at every point. Proposition~\ref{cf:conormal}
with $V=M,Y=\mathcal C_9^\circ$ gives $T|_M=0$.

Take the resulting closed-support antecedent on $D_\partial$.
Its open part $C_s\setminus C_p$ is killed by row two of
\eqref{m6:s:eq:p009-character-table}, with density one.
Left exactness puts support on $C_p$, whose open complement
of $\mathcal E_9$ is killed by rows three and four. Thus
$\operatorname{supp}(T|_{V_9})\subset\mathcal E_9$.
The singular union $D_\partial$ serves as an ordinary support
space in these character steps.

Now $\mathcal E_9$ is closed and smooth in $V_9$ and its
eight-dimensional ambient conormal has the nonnilpotent witness.
Proposition~\ref{cf:conormal} with $V=V_9,Y=\mathcal E_9$
again applies to the original $T$, giving $T|_{V_9}=0$.
The complete cotangent transport
\eqref{m6:s:eq:p009-cotangent-transpose} gives the same argument
on $V_{14}$. The paired zero restrictions glue.
\end{proof}

The two conormal exclusions apply to open restrictions of the original
\(T\) on \(G\), with its nilpotent wavefront hypothesis. The
intervening character steps use ordinary closed-support antecedents.
The complete predecessor lists for support assembly are
\begin{equation}\label{m6:s:eq:p009-boundary-list}
B_9=B_{14}=\mathcal P_1\sqcup\mathcal P_2\sqcup\mathcal P_6
\sqcup\mathcal P_7\subset\mathcal O_8.
\end{equation}
Thus the already zero $\mathcal O_8$ supplies both antecedents,
and the proposition enlarges the zero open to $\mathcal S_9$.

Every division above is by a stated nonzero field element or a
nonzero rational integer. The arguments therefore include residue
characteristic two.

\subsubsection{The spectral pair \texorpdfstring{$\mathcal P_{12},\mathcal P_{22}$}{P12 and P22}}
\label{m6:s:sec:p012-p022}

Here the entire unipotent-compatible envelope is smooth and has
a uniform nonnilpotent ambient conormal, including its
stabilizer-character-incompatible divisor. A single
whole-envelope wavefront step therefore applies to the original
global distribution. The rational gauges, group stabilizers and
expanded tangents are given in \S\ref{supp:m6-p012-calculations}.

\paragraph{The coefficient envelope and residual action}
Use the coefficient section $\ell_0$ of
\S\ref{m6:s:sec:p009-p014}. The entire unipotent-compatible
family is
\begin{equation}\label{m6:s:eq:p012-family}
 \mathsf P=\begin{pmatrix}p&q\\r&t\end{pmatrix},\quad
 \mathsf Q=\begin{pmatrix}0&b\\c&d\end{pmatrix},\quad
 \mathsf R=\begin{pmatrix}e&0\\f&c\end{pmatrix},\quad
 \mathsf S=\begin{pmatrix}g&b\\h&i\end{pmatrix},
 \qquad bce(pt-qr)(gi-bh)\ne0.
\end{equation}
Put $\Delta=gi-bh$ and denote this twelve-dimensional
coefficient space by $\mathcal Q_{12}$. Set
$l=\ell_0(\mathsf P,\mathsf Q)$,
$r_0=\ell_0(\mathsf R,\mathsf S)$ and $x=lW_{12}r_0^t$,
and let $\mathcal C_{12}$ be the $H(F)^2$-saturation of
$\mathcal Q_{12}$. The coefficient space is smooth and
closed relative to $\GL_2^4$; all entries not forced nonzero
by the displayed determinants may vanish, including $t$.
A rational gauge with a rational lift to the intersection
normalizes
$\mathsf Q=\mathsf S=J_0=\left(\begin{smallmatrix}0&1\\-1&0\end{smallmatrix}\right)$
and $\mathsf R=A_0=\operatorname{diag}(1,-1)$.
On the remaining matrix the entire effective residual action is
\begin{equation}\label{m6:s:eq:p012-residual-action}
 \mathsf P\longmapsto
 \begin{pmatrix}1&-h_0/\rho\\0&1/\rho\end{pmatrix}
 \mathsf P
 \begin{pmatrix}1&0\\t_0&1\end{pmatrix},
 \qquad \rho\in F^\times,\quad h_0,t_0\in F;
\end{equation}
Its additional scalar action kernel is $F^\times$.

\begin{proposition}\label{m6:s:prop:p012-orbits}
The compatible coefficient envelope has exactly two branches of
rational orbits. On \(t\ne0\), their normal forms are
\[
 (\mathsf P,\mathsf Q,\mathsf R,\mathsf S)
 =(\operatorname{diag}(\kappa,1),J_0,A_0,J_0),\qquad
 \kappa=-\frac{c\Delta e\det\mathsf P}{t}\in F^\times.
\]
On \(t=0\), their normal forms are
\[
 \left(\begin{pmatrix}0&\upsilon\\1&0\end{pmatrix},J_0,A_0,J_0\right),
 \qquad \upsilon=-c\Delta p_{12}\in F^\times.
\]
Each invariant is necessary and sufficient over \(F\) on its
branch, and all its nonzero values occur.
\end{proposition}
On $t\ne0$, the stabilizer is a central $\mathbb G_m$ times a
split connected six-dimensional unipotent group, with trivial
character. At $t=0$, it gains one unipotent direction with
character $\psi(-ebr v/(cq))$, where $v\in F$ is its additive
parameter; its coefficient is nonzero
because $b,c,e,q,r\ne0$. Both stabilizers are unimodular and
have relative density one. Thus the full character-compatible
locus is exactly $t\ne0$.

Let $x_\kappa$ denote the first normal representative in
Proposition~\ref{m6:s:prop:p012-orbits}. For fixed $\kappa\ne0$,
its coefficient fiber is
$\kappa t+c\Delta e\det\mathsf P=0$.
Since the second term never vanishes, the saturated equation excludes
$t=0$. Each fiber is therefore one closed rational
orbit with its homogeneous quotient topology, so the
orbital integral $\mathscr T_\kappa$ of
\eqref{m6:s:eq:p012-orbital-distribution} is well defined on
compact test supports. Transpose sends it to the partner cell
and preserves the parameter $\kappa$, the character and the
positive quotient measure.

\begin{proposition}\label{m6:s:prop:p012-intrinsic}
The ordinary distribution space
\(D'_{\widetilde\eta}(\mathcal P_{12}\sqcup\mathcal P_{22})\)
is nonzero. More precisely,
\(\mathscr T_\kappa-\tau_*\mathscr T_\kappa\ne0\) for every
\(\kappa\in F^\times\).
\end{proposition}
\begin{proof}
Transpose transports the phase and positive measure together,
exchanges the $H$ factors and preserves $\eta$; its transported
logarithmic Jacobian is one and scalar inversion preserves
multiplicative Haar measure. The two supports are disjoint closed
orbits in the pair. A compact-open test near $x_\kappa$, where
the phase is one, gives a positive first integral and zero
transposed integral. Their closed pushforwards therefore give
a nonzero odd distribution.
\end{proof}

To see how the compatible union meets the incompatible divisor,
fix the normalized triple. Then
$\mathsf P(s)=\left(\begin{smallmatrix}0&q_0\\1&s\end{smallmatrix}\right)$,
$q_0\ne0$, has determinant $-q_0$ and for $s\ne0$ invariant
$\kappa=-q_0/s$, but tends to the incompatible divisor as $s\to0$.
The compatible branch is therefore not closed in the ambient space.
For the wavefront argument, we retain the entire smooth envelope,
including this divisor.

\paragraph{The smooth spectral support}

For the root basis $N_1,\ldots,N_{22}$ of $\mathfrak u$, put
$M_U(x)=[\operatorname{vec}(N_jx),\operatorname{vec}(xN_j^t)]_{j=1}^{22}$
and let $L_U$ be the row of the two adjacent-trace phases on
each side. Both parabolic factors normalize $U$, so the kernel
has dimension six throughout the cell and $\operatorname{rank}M_U=38$.
Finite exponential identifies group-character triviality with
$L_U|_{\ker M_U}=0$. Hence
\begin{equation}\label{m6:s:eq:p012-intrinsic-envelope-equations}
\mathcal C_{12}=\left\{x\in\mathcal P_{12}:
\operatorname{rank}\begin{pmatrix}M_U(x)\\L_U\end{pmatrix}\le38\right\}.
\end{equation}
The $39$-minor equations intrinsically define the closed envelope,
independently of a coefficient factorization.

Put $I_{12}=P\cap W_{12}P^tW_{12}^{-1}$. The intersection calculation
in \S\ref{supp:m6-p012-calculations} gives $\dim I_{12}=18$.
Since $\dim P=34$, the cell has dimension $50$.
The saturated inverse image is smooth of dimension $52+12=64$.
Subtracting the eighteen-dimensional quotient fiber in local product
charts shows that $\mathcal C_{12}$ is a smooth embedded manifold of
dimension $46$ everywhere, including $t=0$ and after scalar
extension. The descended divisor $\mathcal D_{12}:t=0$ is gauge
invariant and smooth of dimension $52+11-18=45$, since $qr\ne0$.
For a basis $B_1,\ldots,B_{26}$ of $\Lie(H)$, write
$M_H(x)=[\operatorname{vec}(B_jx),\operatorname{vec}(xB_j^t)]_{j=1}^{26}$.
Equivalently $\mathcal D_{12}=\{x\in\mathcal C_{12}:\operatorname{rank}M_H(x)\le44\}$,
so its closedness is expressed by $45$-minors.
The complement $\mathcal C_{12}^\circ$ is smooth but not closed,
as the degeneration above shows. These statements transfer by
transpose to $\mathcal C_{22},\mathcal D_{22}$.

\paragraph{The tangent space of the envelope and ten-dimensional ambient conormal}

Use the right-trivialized trace pairing on $\mathfrak{gsp}_{10}$
\eqref{m6:s:eq:p009-trace-trivialization}, conjugated by \(l^{-1}\)
at \(x=lW_{12}r_0^t\).

\begin{proposition}\label{m6:s:prop:p012-conormal}
The entire smooth envelope \(\mathcal C_{12}\), including \(t=0\),
has ambient conormal of dimension \(56-46=10\).
At \(x\) it is
\(\operatorname{Ad}_l\operatorname{span}_F(C_1,\ldots,C_{10})\),
where the complete basis is \eqref{m6:s:eq:p012-conormal-basis} in
\S\ref{m6:s:cert:p012-conormals}; at \(h_1xh_2^t\), apply
\(\operatorname{Ad}_{h_1l}\) to that basis. A uniform witness is
\[
\begin{split}
D_{12}=C_9={}&\operatorname{diag}(0,0,1,1,0,0,-1,-1,0,0)\\
&-(f/e+g/b)(E_{25}+E_{69})-(2i/b)E_{47},
\end{split}
\]
with
\[
\det(sI-D_{12})=s^6(s-1)^2(s+1)^2,\qquad
\operatorname{tr}(D_{12}^2)=4\ne0.
\]
\end{proposition}

\begin{proof}
The tangent formula \eqref{m6:s:eq:p012-full-tangent-span}
includes all forty-four unipotent and eight reductive-action
directions, and all twelve coefficient tangents
\((p,q,r,t,b,c,d,e,f,g,h,i)\). In particular, the derivatives
in \(b,c\) act simultaneously on every block containing them.
The scalar similitude tangent is retained; pairing with it gives
\(\operatorname{tr}Z=0\).
\S\ref{m6:s:cert:p012-conormals} gives the logarithmic derivatives,
the complete pairing lists, eight independent pivot equations,
and all ten independent conormal matrices and dual tangent equations.
The variable denominators of the resulting basis are only \(b,e\),
which are mandatory nonzero. Its rank is uniform at \(t=0\)
and when \(f=g=i=0\). The displayed witness is symplectic and
upper triangular in physical coordinates, giving its characteristic
polynomial and trace-square. Conjugation transports it to every
rational saturated point, with arbitrary multiplier.
\end{proof}

On \(t=0\) the orbit tangent loses one dimension, but the twelve
coefficient derivatives supply the extra transverse direction
needed for the same forty-six-dimensional support tangent.
Transpose carries the whole ten-dimensional conormal to the
partner by \eqref{m6:s:eq:p009-cotangent-transpose}, preserving
\eqref{m6:s:eq:p012-witness-polynomial}.

\paragraph{Conditional spectral exclusion and its antecedents}

Define the rank opens, using \eqref{m6:s:eq:p012-rank-matrices}, \(V_i=\bigcap_{a,b}\{x:R_{ab}(x)\ge(R_i)_{ab}\},\qquad B_i=\bigcup_{a,b}\{x\in V_i:R_{ab}(x)>(R_i)_{ab}\}, \quad i=12,22.\)
The rank classification gives \(V_i\setminus B_i=\mathcal P_i\).
Both sets are open in \(G\), and the cell is closed
in \(V_i\). Put \(V_{12,22}=V_{12}\cup V_{22},\qquad B_{12,22}=B_{12}\cup B_{22}.\)
Incomparability gives
\(V_{12,22}\setminus B_{12,22}
=\mathcal P_{12}\sqcup\mathcal P_{22}\).
Each cell is itself closed in the union open: for example its
complement there is \(B_{12}\cup V_{22}\) for \(\mathcal P_{12}\).
The corresponding \(\mathcal C_i\) and \(\mathcal D_i\) are
therefore closed in their rank opens and in the pair open.

\begin{proposition}\label{m6:s:prop:p012-spectral}
Let \(T\in D'_\eta(G)\cap\NW(G)\).
Then \(T|_{B_i}=0\ \Longrightarrow\ T|_{V_i}=0\qquad(i=12,22).\)
In particular \(T|_{B_{12,22}}=0\) implies
\(T|_{V_{12,22}}=0\).

\end{proposition}

\begin{proof}
Assume $T|_{B_{12}}=0$. The support sequence
\cite[(2.2)--(2.3)]{JLX} gives the unique equivariant cell
antecedent whose closed pushforward is $T|_{V_{12}}$.
Its open unipotent-incompatible part vanishes by
\cite[Theorem 2.5 and Remark 2.6]{JLX}: the action is algebraic
and constructible on $F$-points
\cite[Theorem 2.3 and Remark 2.4]{JLX}, and both unipotent
modular functions are one. Another closed-support injection yields
$\operatorname{supp}(T|_{V_{12}})\subset\mathcal C_{12}$,
including the incompatible divisor.

We can now apply the conormal criterion: the entire envelope is
closed and smooth, and Proposition~\ref{m6:s:prop:p012-conormal}
supplies a nonnilpotent ambient conormal at every point, including
the divisor.
Proposition~\ref{cf:conormal} with $V=V_{12}$,
$Y=\mathcal C_{12}$ and the same original global $T$
therefore gives $T|_{V_{12}}=0$.
Transpose transfers the support geometry and full cotangent fiber
to $V_{22}$, so the identical argument applies to the same $T$
there. Under the paired boundary hypothesis the two zero
restrictions glue.
\end{proof}

For the global support induction the predecessor
lists are
\begin{equation}\label{m6:s:eq:p012-boundary-lists}
\begin{aligned}
B_{12}&=\bigsqcup_{j\in\{1,2,3,4,6,7,8,9,10,11\}}\mathcal P_j,\\
B_{22}&=\bigsqcup_{j\in\{1,2,3,4,6,7,8,14,15,21\}}\mathcal P_j.
\end{aligned}
\end{equation}
Both are contained in the preceding
\(\mathcal S_{11}=\bigsqcup_{j\in\{1,\ldots,11,14,15,21\}}\mathcal P_j\).
The passage from \(\mathcal S_9\) to \(\mathcal S_{11}\)
uses the incompatible cells \(10,11\) and their reflected
partners \(21,15\). Once those
previous opens have vanished, Proposition~\ref{m6:s:prop:p012-spectral}
enlarges the zero open to \(\mathcal S_{12}\).

Throughout the conditional argument the NW condition belongs
to the original \(T\) on \(G\). The trace form and finite unipotent
exponential and logarithm are valid over every characteristic-zero
field considered here, including residue characteristic two.

\subsubsection{The intrinsic orbital families}\label{m6:orbital-families}

We begin with a common criterion for intrinsic odd vanishing.
Table~\ref{m6:tab:orbit-families} records the rational and
group-theoretic inputs for its applications: saturated compatible
loci, complete invariants and fibers, stabilizers, relative densities,
and character-one transpose transporters. The first eight rows are
proved at the appendix locations cited there; we treat the ninth
row in \S\ref{m6:s:sec:p090}.

Although these families have different rational quotients and
stabilizers, their orbital lines are all transpose-even. This is the
common argument developed below. The spectral packets, which admit
intrinsic odd orbital distributions, are treated by the separate
ambient conormal arguments above.

\begin{proposition}[Rational orbital-family criterion]\label{m6:orbit-family}
Let $X$ be an invariant, transpose-stable locally closed $F$-space in
the present action. Suppose that its unipotently compatible locus
$C$ is closed in $X$, that the incompatible complement has ordinary
equivariant vanishing, and that there is a continuous invariant map
$q:C\to B$ to an $l$-space on which transpose acts trivially.
Assume, for every $b\in B$, that:
\begin{enumerate}
\item $q^{-1}(b)$ is one rational $\mathcal A(F)$-orbit, with the
homogeneous-space topology supplied by a surjective orbit submersion;
\item its stabilizer has trivial $\eta$ and trivial relative
modular character;
\item a rational transporter from a point to its transpose
has character one.
\end{enumerate}
Then $D'_{\widetilde\eta}(X)=0$. Every nonempty fiber has a nonzero
even orbital distribution whose closed pushforward to $X$ is defined.
If $B$ consists of one point, $\dim D'_\eta(X)=1$.
\end{proposition}
\begin{proof}
By the homogeneous distribution criterion, each fiber has one
equivariant line, represented by character-weighted positive quotient
Haar measure. A character-one transporter makes transpose act on this
line by a positive scalar; since transpose squares to the identity,
that scalar is $1$. Every fiber therefore has zero odd space, and
localization
(Proposition~\ref{cf:localization}) gives zero odd space on $C$.
Together with vanishing on $X\setminus C$, the left-exact support
sequence gives the conclusion on $X$. Since $q$ is continuous and
points of $B$ are closed, each fiber is closed in the compatible
locus and hence in $X$. Its locally finite Haar line therefore
has a nonzero closed pushforward.
When there is only one fiber, the same support sequence identifies
the whole even space with that line.
\end{proof}

\begin{proposition}[Complete intrinsic family data]\label{m6:orbit-verification}
The compatible loci in Table~\ref{m6:tab:orbit-families} satisfy the
criterion above with the displayed complete invariant bases. Each
stabilizer in the table is the $F$-rational stabilizer group. Every row has
trivial algebraic component group, trivial character on the stabilizer,
and trivial relative orbital density; intermediate coefficient
densities are retained separately in its calculation.
\end{proposition}

Here connectedness refers to the algebraic stabilizer group. The criterion
uses the rational orbit with its homogeneous topology, including
singular centralizers. It is an intrinsic ordinary-distribution statement.
In the global proof it is applied only to the unique closed-support
antecedent of the original distribution after the indicated predecessor
open has vanished, as in \S\ref{m6:filtration}.

\begingroup\small
\begin{longtable}{@{}p{.065\linewidth}p{.19\linewidth}p{.31\linewidth}>{\raggedright\arraybackslash}p{.32\linewidth}@{}}
\caption{One orbital argument and its distinct rational inputs. Here
$U_n$ denotes the specified split unipotent group of dimension $n$,
with the group law specified in the corresponding calculation. The notation $T^2$ denotes a split torus
isomorphic to $\mathbb G_m^2$.}\label{m6:tab:orbit-families}\\
\toprule Cell&Complete base&Stabilizer&Geometry and verification\\\midrule
\endfirsthead
\toprule Cell&Complete base&Stabilizer&Geometry and verification\\\midrule
\endhead
16&$(\nu,\kappa)\in F\times F^\times$
&$1\to U_5\to K\to Z_{\GL_2}(N_\nu)\to1$
&Rank-one conjugation, including $\nu=0$; p.~\pageref{m6:s:sec:middle-p016}.\\
17&$\kappa\in F^\times$&$\mathbb G_m\times U_6$
&Regular $\mathrm{PGL}_2$ conclusion; p.~\pageref{m6:s:sec:middle-p017}.\\
27&$t\in F^\times$&$T^2\ltimes U_{10}$
&Rank-one idempotents; p.~\pageref{m6:s:sec:middle-p027}.\\
31&$\rho\in F^\times$&$\mathbb G_m\times U_{11}$
&Rational parameter orbits; p.~\pageref{m6:s:sec:sector-p031}.\\
33&Point&$\mathbb G_m\times U_{15}$
&Residual group $J_0$; p.~\pageref{m6:s:sec:sector-p033}.\\
40&$\kappa\in F^\times$&$T^2\ltimes U_6$
&Rank-one, nonzero-trace presentation; p.~\pageref{m6:s:sec:p040}.\\
53&Point&$\mathbb G_m\times U_{12}$
&One rational homogeneous orbit; p.~\pageref{m6:s:sec:p053}.\\
69&Point&$\mathbb G_m\times U_{15}$
&One rational homogeneous orbit; p.~\pageref{m6:s:sec:p069}.\\
90&Point&Graph isomorphic to $H$
&The closed translate $Hc$; Section~\ref{m6:s:sec:p090}.\\
\bottomrule
\end{longtable}\endgroup

\begin{proof}
The coefficient calculations cited in the last column verify the
four pieces of input: the entire compatible zero set, its complete
rational invariant, the group equations with character and
modular character, and the transpose transporter. The representative
quotient calculations below explain why the first, second, third,
and sixth rows have different geometries.

To prove closedness, note that the $U\times U$ stabilizer kernels
have constant rank on each parabolic cell: this normal subgroup is
transported by $P\times P$. Restricting the additive character gives
a regular section of the dual kernel bundle, whose zero set is
precisely the compatible locus. In the complete coefficient
atlas this is exactly the displayed saturated zero set. Smooth
surjective quotient charts give its stated topology and dimensions.
The explicit invariants descend through the overlap action,
and their rational normalizations identify each fiber with one
orbit. The divisors allowed in each coefficient family remain in
those normalizations.

For the singleton rows $53$ and $69$, we prove transitivity using
geometric irreducibility of the compatible varieties. The stabilizers
have constant dimension at every geometric point, so all geometric
orbits have the dimension of the variety and are open. Irreducibility
then gives one orbit. At the displayed rational point the stabilizer is
$\mathbb G_m\times U_n$. A split additive central series and
Hilbert 90 give $H^1(F,\mathbb G_m\times U_n)=1$, which kills the
transporter cocycle and proves rational transitivity. 

The remaining singleton rows have explicit rational transporters.
Every density computation is on the entire stabilizer, including
its reductive or scalar factor. The signed transpose formulas in
the cited calculations have character one. These facts verify all
hypotheses of Proposition~\ref{m6:orbit-family} independently of
any wavefront condition.
\end{proof}

The finite calculations retain the jump from the torus to the
nonreductive centralizer in row $16$, the additive isotropy in row $33$,
and the two independent extra stabilizer directions in row $69$.
The rank-zero Shalika family has additional singular fibers and is
treated separately in \S\ref{supp:m6-sharpness}. The incidence pair has an opposite-weight
plane and a stabilizer jump; its two arguments also remain separate.

\subsubsection{The remaining cells in the \texorpdfstring{$(2,1)$}{(2,1)} sector}
\label{m6:s:sec:middle21}

We apply Proposition~\ref{m6:orbit-family} through three quotient
descriptions. The complete coefficient gauges, stabilizer groups,
and intermediate density characters are computed on
p.~\pageref{m6:s:sec:middle-p016},
p.~\pageref{m6:s:sec:middle-p017}, and p.~\pageref{m6:s:sec:middle-p027},
respectively, including the nilpotent, nonsquare, and free-zero branches.

Their residual geometries are distinct. For
\(\mathcal P_{16}\) the residual matrix has rank one and can be nilpotent;
\(\mathcal P_{17}\) leads to the regular action on \(\mathrm{PGL}_2\);
\(\mathcal P_{27}\) leads to conjugation of rank-one idempotents. The
arguments establish intrinsic odd vanishing.

We use upright letters \(\mathsf P,\mathsf Q,\mathsf R,\mathsf S\)
for coefficient matrices, and put
\[
 \mathsf l(\mathsf P,\mathsf Q)=m(\mathsf Q\mathsf P,\mathsf Q,I_2),\qquad
 x_i(\mathsf P,\mathsf Q,\mathsf R,\mathsf S)
 =\mathsf l(\mathsf P,\mathsf Q)W_i\mathsf l(\mathsf R,\mathsf S)^t.
\]
In each individual calculation, \(w\) denotes the indicated \(W_i\).
For the root lists below, define
\[
\mathcal N(X,Y,Z,W,T,S)=
\begin{pmatrix}
0&X&Z&W&S\\0&0&Y&T&-W^\sharp\\0&0&0&-Y^\sharp&-Z^\sharp\\
0&0&0&0&-X^\sharp\\0&0&0&0&0
\end{pmatrix},\qquad T,S\in\mathfrak{sl}_2.
\]
Thus \(X_{21}\) is the matrix with \(X=E_{21}\) and all other inputs
zero; \(T_h,T_e,T_f\) mean respectively \(T=\operatorname{diag}(1,-1),E_{12},E_{21}\),
and the notation for \(S\) is identical. Levi root matrices are the inputs in
\[
\mathcal L(A,B,C)=\operatorname{diag}\bigl(A,B,C,
B+(\operatorname{tr}C-\operatorname{tr}B)I_2,
A+(\operatorname{tr}C-\operatorname{tr}A)I_2\bigr).
\]
In particular \(A_{21},B_{21},C_{21}\) denote the three indicated inputs;
we also write \(L_A=A_{21},L_B=B_{21},L_C=C_{21}\).

The quotient and stabilizer conventions are
\eqref{m6:s:eq:atlas-quotient}--\eqref{m6:s:eq:atlas-stabilizer}.
For each family below, the compatible equations are closed in the
invertible coefficient space and their preimage is saturated;
Section~\ref{m6:s:sec:coefficient-atlas} therefore gives the stated cell
closedness and the local smooth charts. Finite logarithms are used only
after the complete scalar and reductive projections have been solved,
as in Lemma~\ref{m6:s:lem:finite-exponential}.

\begin{proposition}\label{m6:s:prop:middle21}
For \(i=16,17,27\), one has
\(D'_{\widetilde\eta}(\mathcal P_i)=0\).
Each of these cells nevertheless has nonzero \(\eta\)-equivariant
ordinary distributions.
\end{proposition}
\paragraph{Complete quotient data.}
The rational calculations in the cited appendix sections verify the
following three inputs to Proposition~\ref{m6:orbit-family}.
They retain the nontrivial intermediate densities before taking the
homogeneous quotient.
\begin{center}
\begin{tabular}{@{}p{.08\linewidth}p{.29\linewidth}>{\raggedright\arraybackslash}p{.25\linewidth}p{.298\linewidth}@{}}
\toprule Cell & Complete invariant and quotient & Full stabilizer & Density and transpose \\\midrule
16 & $\nu=\operatorname{tr}N\in F$, $\kappa\in F^\times$;
$\operatorname{rank}N=1$, conjugation & Extension of $Z_{GL_2}(N)$ by
$U_5$; dimension $7$ & Raw $D(f\circ k)=|q|^6D(f)$, $d'=qd$;
$|d|^6D$ is invariant. Sections $d=1$ make the residual line trivial;
$N\mapsto N^t$. \\
17 & $\kappa=\det\mathsf P\det\mathsf R\,gl\in F^\times$;
regular $GL_2$ action on $PGL_2$ & $\mathbb G_m\times U_6$;
dimension $7$ & Full orbital density one;
$[b]\mapsto[b^{-t}d_0]$, $d_0=\operatorname{diag}(1,-1)$, preserves Haar. \\
27 & $t=eDD_R\in F^\times$, $D=ad-bc$, $D_R=fh-ik$;
conjugation on rank-one idempotents & $T^2\ltimes U_{10}$;
dimension $12$ & Six weights $1$, two $\lambda/\mu$, two $\mu/\lambda$;
product one. $E\mapsto E^t$, with symmetric original representative $x_t$. \\
\bottomrule
\end{tabular}
\end{center}
For row 16 the definitions are
$v=\mathsf Pe_2$, $\phi=e_1^t\mathsf Q$,
$d=\det\mathsf P\det\mathsf R$, $N=v\phi/d$,
$\kappa=\det\mathsf Q\det\mathsf S$.
To identify the rational quotient, we use the free gauge action
and its exact fibers. The section
$N_\nu=\left(\begin{smallmatrix}\nu&1\\0&0\end{smallmatrix}\right)$
and symmetric transporter
$b_\nu=\left(\begin{smallmatrix}\nu&1\\1&0\end{smallmatrix}\right)$,
$\det b_\nu=-1$, work also at $\nu=0$. The centralizer changes from
split torus to $\mathbb G_m\times\mathbb G_a$ there, with unchanged
dimension and no new component.
The scalar coefficients in rows 17 and 27 refer to the normal forms
in \S\ref{m6:s:sec:middle-p017} and \S\ref{m6:s:sec:middle-p027},
respectively.
For row 17 the rational normalization divides only by the forced units
$b,e,g,l,\det\mathsf P,\det\mathsf R$, so every nonsquare $\kappa$
is included. For row 27 it divides only by $d,e,f,D$ and includes the
divisor $di=fm$ and all square classes of $t$. Put
$T_\Delta=\{(d,d^{-1}):d\text{ is invertible diagonal}\}$ and
$\mathcal E=\{E\in\Mat_2(F):E^2=E,\ \operatorname{rank}E=1\}$.
Its exact quotient map is
\[
 [(a,b)]\mapsto(ab^t,aE_{11}a^{-1}):
 (GL_2\times GL_2)/T_\Delta\simeq GL_2\times\mathcal E.
\]
Nonzero columns of $E$ and $I-E$ give rational sections; Frobenius at
$z=I$ has trivial relative density and leaves $E\mapsto E^t$.
The appendix proves full rational coverage, full group stabilizers,
character-one transporters, local quotient topology and closedness of
every invariant fiber in its cell for all three rows.

\begin{proof}[Proof of Proposition~\ref{m6:s:prop:middle21}]
These data verify the middle rows of
Proposition~\ref{m6:orbit-verification}. The common orbital argument
therefore proves intrinsic odd vanishing on each whole cell.
Every compatible fiber is closed in the cell, with full character and
relative density trivial; its nonzero weighted Haar line has a closed
pushforward. This proves even nonvanishing as well.
\end{proof}

\subsubsection{The \texorpdfstring{$(1,1)$}{(1,1)} sector and its degeneration geometry}
\label{m6:s:sec:sector11}

We identify the orbit invariants and transpose action, then assemble
the three-cell support argument. The complete rational group and
density calculations appear on p.~\pageref{m6:s:sec:sector-p031}
and p.~\pageref{m6:s:sec:sector-p033}; the extra stabilizer character
is computed in \S\ref{m6:s:cert:p032-full}. The incidence cone
and its submersion, including the compatibility jump, are described
in \S\ref{m6:s:cert:three-cell-cone}.

We use the coefficient lift and root matrices defined in
Section~\ref{m6:s:sec:middle21}. The signed columns of the three representatives are
\[
\begin{split}
W_{31}&=(1,-9,3,-7,5,6,4,8,2,10),\\
W_{32}&=(1,-9,8,5,4,-7,-6,-3,2,10),\\
W_{33}&=(1,-9,8,-7,5,6,4,-3,2,10).
\end{split}
\]
For all three cells the projection gauge, with \((A,B,C)\) the first
three blocks of \(k\in I_w\) and \((D,E,F)\) those of \(k^\natural\), is
\begin{equation}\label{m6:s:eq:sector11-gauge}
(\mathsf P,\mathsf Q,\mathsf R,\mathsf S)
\longmapsto(B^{-1}\mathsf P A,C^{-1}\mathsf Q B,
E\mathsf R D^{-1},F\mathsf S E^{-1}).
\end{equation}
The stabilizer membership conditions in the corresponding frame are
\begin{equation}\label{m6:s:eq:sector11-membership}
B\mathsf P=\mathsf P A,\quad C\mathsf Q=\mathsf Q B,\quad
E\mathsf R=\mathsf R D,\quad F\mathsf S=\mathsf S E.
\end{equation}
In these six projection equations the last letter \(F\) denotes a matrix; occurrences such as \(F^\times\) continue to denote the field.
In the individual projection lists we use \(F\) throughout.

\begin{proposition}\label{m6:s:prop:sector11}
One has
\[
D'_{\widetilde\eta}(\mathcal P_{31})=0,\qquad
D'_\eta(\mathcal P_{32})=0,\qquad
D'_{\widetilde\eta}(\mathcal P_{33})=0.
\]
Moreover, in its topology as a locally closed subset of \(G\),
\[
D'_{\widetilde\eta}(\mathcal P_{31}\sqcup\mathcal P_{32}\sqcup\mathcal P_{33})=0.
\]
\end{proposition}

\paragraph{The parameter family on cell $31$}
The complete compatible coefficients and their rational gauge give the
invariant $\rho=p_{11}r_{11}\in F^\times$ and the normal form
\[
 \mathsf P=\mathsf Q=I,\qquad
 \mathsf R=\operatorname{diag}(\rho,-1),\qquad
 \mathsf S=\begin{pmatrix}0&1\\1&0\end{pmatrix}.
\]
Every fiber is a closed rational orbit with stabilizer
$\mathbb G_m\times U_{11}$. Its algebraic additive character is
zero, and its relative density is one. The normalization covers the nilpotent,
nonsplit and free-zero coefficient specializations for every $\rho\ne0$. Transpose
preserves $\rho$. Write $\sigma(h_1,h_2)=(h_2,h_1)$ for factor
exchange on $\mathcal A=H\times H$. For $h_i=u_id_5(a_i)$, define
\[
\lambda(h_1,h_2)=\sum_{i=1}^2
\bigl(\operatorname{tr}(u_i)_{12}+\operatorname{tr}(u_i)_{23}\bigr),
\qquad \eta=\psi\circ\lambda.
\]
If $a$ carries a point to its transpose, then $\sigma(a)a$ stabilizes
the point. Since $\lambda$ vanishes on the stabilizer,
$2\lambda(a)=0$, whence $\lambda(a)=0$ in the characteristic-zero
field. Thus the transporter has character one even in residue
characteristic two. p.~\pageref{m6:s:sec:sector-p031} prints
the complete equations; Proposition~\ref{m6:orbit-family} supplies
the odd conclusion.

\paragraph{The cell \texorpdfstring{$\mathcal P_{32}$}{P32}}\label{m6:s:sec:sector-p032}

The P032 row of Table~\ref{m6:s:tab:full-stabilizer-exclusions} applies to
the entire family \eqref{m6:s:eq:sector-32-family}. Its additive subgroup
\eqref{m6:s:eq:sector-32-subgroup} has phase \(\psi(tfi/g)\), with
\(f,g,i\ne0\) forced by invertibility, and relative modulus one.
\S\ref{m6:s:cert:p032-full} retains the complete intersection,
normalization, scalar-unipotent stabilizer, fixing identities and signed
transpose transport. In particular all \(\lambda=cef\ne0\), including
nonsquare \(-\lambda\), remain covered. The common exclusion theorem
gives \(D'_\eta(\mathcal P_{32})=0\).

\paragraph{The single compatible orbit on cell $33$}
Here $\mathsf R e_2=-\mathsf P e_2$ and $\mathsf S=-j\mathsf Q$.
Since $\mathsf P e_2\ne0$, the gauge first gives
$(I,I,\operatorname{diag}(\kappa,-1),-j)$ with $\kappa\ne0$,
then the legal diagonal choice $A=B=C=\operatorname{diag}(1,\kappa^{-1})$
gives $\kappa=1$.
Thus the whole compatible locus is one closed rational orbit.
Its stabilizer is $\mathbb G_m\times U_{15}$, with trivial
character, relative density, and algebraic component group.

The intermediate quotient retains the nonreductive group
\[
 J_0=\{(c(I+tE_{21}),c^{-1}(I-tE_{12})):
                     c\in F^\times, t\in F\}.
\]
The full group $J_0=\mathbb G_m\times\mathbb G_a$, including its additive factor,
acts on this intermediate quotient. For the
canonical point $x_0$, with $M=\left(\begin{smallmatrix}0&1\\1&0\end{smallmatrix}\right)$,
the identity $x_0^t=d_5(M)x_0d_5(M)^t$ supplies a character-one
transporter. Proposition~\ref{m6:orbit-family} now applies.
p.~\pageref{m6:s:sec:sector-p033} contains the rational
normalization, the group solution and the density calculation.

\paragraph{The incidence cone and the three-cell support theorem}
Put \(Y=\left(\begin{smallmatrix}a&b\\c&d\end{smallmatrix}\right)\),
\(C=Y^tj\), \(L=-\operatorname{adj}(Y)^t\). The local incidence cone is \(\mathcal Z=\{(Y,D):\det Y=0,\ D=D^t, \ \operatorname{rank}\binom C D\le1, \ \operatorname{rank}(L\ D)\le1\}.\)
Its strata \(\operatorname{rank}D=1\), \(D=0,\operatorname{rank}Y=1\),
and \(Y=D=0\) are respectively P031, P032 and P033, of slice dimensions
four, three and zero. The cone is irreducible of dimension four and has
seven-dimensional tangent space at its singular origin. A local
analytic chart is its product with a smooth 42-dimensional factor;
the seven transverse directions include all the displayed slice variables.

\S\ref{m6:s:cert:three-cell-cone} contains the polynomial matrix,
exhaustive rank equations \eqref{m6:s:eq:degeneration-vii-1}, rational charts
\eqref{m6:s:eq:degeneration-vii-2} (including arbitrary nonsquare rank-one
scalars), the submersion ranks, and every degeneration curve.
In particular \eqref{m6:s:eq:degeneration-vii-3} is a fully compatible P031
curve to the compatible P033 point. Equation~\eqref{m6:s:eq:degeneration-vii-4}
gives, for every \(\lambda\in F^\times\), a fully compatible P031
curve with \(\rho=-\lambda t\) to the unipotently compatible P032 orbit.
At the boundary, new character-bearing stabilizer directions
appear, so stabilizer-character compatibility fails to be closed
across the jump. We obtain the distributional conclusion by the
open support argument below.

\subparagraph{Distributional assembly in the topology}

Use only the established legal lexicographic opens \(\mathcal O_i\). The family is the locally closed l-space \(X=\mathcal O_{33}\setminus\mathcal O_{30} =\mathcal P_{31}\sqcup\mathcal P_{32}\sqcup\mathcal P_{33}.\)

Here \(\mathcal P_{31}=X\cap\mathcal O_{31}\) is open in \(X\), and its complement \(Z=\mathcal P_{32}\sqcup\mathcal P_{33}\) is closed in \(X\). In \(Z\), \(\mathcal P_{32}\) is open and \(\mathcal P_{33}\) is closed. All these spaces are stable under \(\widetilde{\mathcal A}\). Therefore

\[
0\longrightarrow D'_{\widetilde\eta}(Z)\longrightarrow D'_{\widetilde\eta}(X)
\longrightarrow D'_{\widetilde\eta}(\mathcal P_{31}),
\]
\[
0\longrightarrow D'_{\widetilde\eta}(\mathcal P_{33})\longrightarrow D'_{\widetilde\eta}(Z)
\longrightarrow D'_{\widetilde\eta}(\mathcal P_{32})
\]

are the permitted left-exact sequences. The last term of the first vanishes by \(\mathcal P_{31}\); the outer terms of the second vanish by the stronger \(\mathcal P_{32}\) theorem and the \(\mathcal P_{33}\) odd theorem. These exact sequences prove \(D'_{\widetilde\eta}(X)=0\).
These cell arguments and the two support sequences prove Proposition~\ref{m6:s:prop:sector11}.

\subsubsection{The outer cells and the eight-cell interval}
\label{m6:s:sec:outerlate}

We prove two orbital statements and then assemble the eight-cell
support filtration. The $F$-rational stabilizer groups, invariants,
and densities are computed on p.~\pageref{m6:s:sec:p040} and
p.~\pageref{m6:s:sec:p053}, while \S\ref{m6:s:cert:p041-full}
gives the intervening character obstruction. Explicit degeneration
curves and their minors appear in \S\ref{m6:s:cert:outer-degenerations}.

We retain the coefficient lift \(\mathsf l\), root matrices, and two-sided
action from Section~\ref{m6:s:sec:middle21}. Here \(U_n\) denotes the split unipotent group obtained by finite
exponential from the displayed nilpotent Lie algebra, with its
resulting group law.

\begin{proposition}\label{m6:s:prop:outerlate}
The cells \(\mathcal P_{40},\mathcal P_{41},\mathcal P_{53}\) satisfy
\[
D'_{\widetilde\eta}(\mathcal P_{40})=0,\qquad
D'_\eta(\mathcal P_{41})=0,\qquad
\dim D'_\eta(\mathcal P_{53})=1,\quad
D'_{\widetilde\eta}(\mathcal P_{53})=0.
\]
Moreover \(D'_\eta(\mathcal P_{40})\ne0\).
\end{proposition}

\paragraph{The orbital argument and its finite input}
Proposition~\ref{m6:orbit-family} supplies the common distributional
argument. The two surviving compatible geometries are a family of
rank-one semisimple matrices on cell $40$ and a single rational
homogeneous orbit on cell $53$. The intervening cell $41$ is removed
by an extra character-bearing subgroup. Stabilizer, rationality and density calculations are printed in
\S\ref{m6:supp:outer-orbits}; the interval support argument
below is retained in the main text.

\paragraph{Cell $40$: a nonzero-trace rank-one family}
The compatible coefficients have the form
\[
 \mathsf P=\begin{pmatrix}a&0\\b&c\end{pmatrix},\quad
 \mathsf Q=\begin{pmatrix}d&e\\f&g\end{pmatrix},\quad
 \mathsf R=\begin{pmatrix}a&0\\h&e\end{pmatrix},\quad
 \mathsf S=\begin{pmatrix}i&c\\j&k\end{pmatrix},\qquad
 ace\det\mathsf Q\det\mathsf S\ne0.
\]
Their complete invariant is $\kappa=a\det\mathsf Q\det\mathsf S$.
Rational operations divide only by the parameters forced to be nonzero and lead to
$\mathsf P=\mathsf R=I$,
$\mathsf Q=\left(\begin{smallmatrix}0&1\\1&0\end{smallmatrix}\right)$,
$\mathsf S=\left(\begin{smallmatrix}0&1\\\kappa&0\end{smallmatrix}\right)$.
The associated ten-dimensional representative $x_\kappa$ is symmetric.
The stabilizer is $T^2\ltimes U_6$, of dimension eight. Its
unipotent weights are
\[
 1,\ \alpha/\beta,\ 1,\ 1,\ \beta/\alpha,\ 1;
\]
the product is one, and the two stabilizer character phases cancel
on the entire group. Thus each $\kappa$-fiber has one even orbital
line and no odd line, by Proposition~\ref{m6:orbit-family}.

There is a precise rank-one interpretation. For a rank-one matrix
$M$ with $t=\operatorname{tr}M\ne0$, put $E=M/t$, $F_0=I-E$.
Then $E^2=E$, and
\[
 \sigma(M)=\begin{pmatrix}
 0&t^{-1}F_0&E&0&0\\t^{-1}F_0&E&0&0&0\\
 E&0&0&0&F_0\\0&0&0&F_0&M\\0&0&F_0&M&0
 \end{pmatrix}
\]
is an equivariant family of representatives with
$\sigma(M)^t=\sigma(M^t)$ and $t=\kappa^{-1}$.
Its saturation map is a submersion: the 44 orbit directions and
$d\kappa=-t^{-2}dt$ supply all 45 directions. This identifies the
residual geometry and retains the six-dimensional unipotent
isotropy of the original orbit.

The residual parameter domain excludes trace-zero $M$ because
$\kappa\ne0$. A raw block of an unnormalized compatible point can
nevertheless be nilpotent; the degeneration curves below exhibit
this. The divisors $ad+be=0$ and
$ai+ch=0$, individually and together, are included in the rational
normalization. p.~\pageref{m6:s:sec:p040} supplies the equations
and density calculation at those divisors as well.

\paragraph{The cell \texorpdfstring{$\mathcal P_{41}$}{P41}}\label{m6:s:sec:p041}

The P041 row of Table~\ref{m6:s:tab:full-stabilizer-exclusions} records the
additive stabilizer and its nonzero phase throughout the complete
family \eqref{m6:s:eq:late-a41-1}. \S\ref{m6:s:cert:p041-full} gives its
19-dimensional intersection, group equations, fixing identity
\eqref{m6:s:eq:late-a41-12}, and all preliminary extraction determinants.
Since the relative modular character is one,
Theorem~\ref{m6:s:thm:full-stabilizer-exclusions} gives
\(D'_\eta(\mathcal P_{41})=0\).

\paragraph{Cell $53$: a single rational homogeneous orbit}
The compatible family is
\[
 \mathsf P=\begin{pmatrix}0&b\\c&d\end{pmatrix},\quad
 \mathsf Q=\begin{pmatrix}f&-j\\h&i\end{pmatrix},\quad
 \mathsf R=\begin{pmatrix}0&b\\-c&e\end{pmatrix},\quad
 \mathsf S=\begin{pmatrix}k&i\\l&j\end{pmatrix},\qquad
 bc(fi+hj)(kj-il)\ne0.
\]
The stabilizer is $\mathbb G_m\times U_{12}$ on this entire
family, with algebraic additive character zero and relative density
one. In particular neither $i$ nor $j$ is individually required
to be nonzero. The saturated compatible variety is smooth and geometrically
irreducible, closed in the cell, and of dimension $39$. The
uniform stabilizer dimension $13$ makes it one geometric orbit.
The rational base point has this split scalar-unipotent stabilizer,
whose $H^1$ vanishes by multiplicative Hilbert 90 and a split additive
central series. Hence it is one rational orbit as well.

At the coefficient tuple
$(\mathsf P,\mathsf Q,\mathsf R,\mathsf S)=(M,I,-j,M)$,
where $M=\left(\begin{smallmatrix}0&1\\1&0\end{smallmatrix}\right)$,
the resulting point satisfies
$x_{53,0}^t=d_5(M)x_{53,0}d_5(M)$.
Both factors have character one. Proposition~\ref{m6:orbit-family}
therefore gives one even line and zero odd space. All these assertions,
including the group solution for $i=0$ or $j=0$, are proved
in p.~\pageref{m6:s:sec:p053} of the appendices.

\paragraph{Parabolic and compatible degenerations}\label{m6:s:sec:outerdegeneration}
The parabolic degenerations P040 to P041 and P041 to P053 hold for
every nonzero curve parameter. \S\ref{m6:s:cert:outer-degenerations}
records the polynomial curves, torus multipliers and full-rank minors.
It also gives the stronger compatible curve \(b(t)=(I+t(E_{52}+E_{96}-E_{74}))x_{53},\qquad b(0)=x_{53}.\)
For every \(t\ne0\), the character-one rational transporter
\eqref{m6:s:eq:late-p40-10} sends \(b(t)\) to \(x_{1/t}\). Hence
\(\kappa(b(t))=t^{-1}\) and
\(C_{53}\subset\overline{C_{40}}^{\,G(F)}\), also algebraically.

The normalized residual matrix is \(tE_{11}\); a conjugate tends to
the nonzero nilpotent \(E_{12}\), but its conjugator and
\eqref{m6:s:eq:late-p40-9}--\eqref{m6:s:eq:late-p40-10} are singular at zero.
Thus neither limit defines a trace-zero value of the normalized section.
The raw block \(b(t)_{45}=b(t)_{54}=-E_{21}\) remains distinct from that
normalized matrix. 

\paragraph{The eight-cell interval}\label{m6:s:sec:outerinterval}

\subparagraph{interval membership}

Retain the established ranks \(R_{ab}(g)=\operatorname{rank}g_{\{2a+1,\ldots,10\},\{2b+1,\ldots,10\}}, \qquad 1\le a,b\le4.\)

The rank matrices at the two ends of the interval are

\[
R_{40}=\begin{pmatrix}6&5&4&2\\5&5&4&2\\4&4&3&1\\2&2&1&0\end{pmatrix},
\quad
R_{53}=\begin{pmatrix}6&5&3&1\\5&4&2&1\\3&2&2&1\\1&1&1&0\end{pmatrix}.
\]

Each row of the following table is obtained from the rank classification, with the complementary rows recovered by its established identity \(R_{ab}=10-2a-2b+R_{5-a,5-b}\). The finite classification in Section~\ref{m6:s:sec:cells} makes the two boundary inequalities a comparison of the displayed rank entries. The rows passing both comparisons are exactly:

\begin{center}
\begin{tabular}{@{}p{0.18\linewidth} p{0.61\linewidth} p{0.13\linewidth}@{}}
\hline
Cell & First two rows of \(\mathsf{R}\) & Transpose \\
\(\mathcal P_{40}\) & \((6,5,4,2;\ 5,5,4,2)\) & \(\mathcal P_{40}\) \\
\(\mathcal P_{41}\) & \((6,5,4,2;\ 5,4,3,2)\) & \(\mathcal P_{41}\) \\
\(\mathcal P_{42}\) & \((6,5,4,2;\ 5,4,3,1)\) & \(\mathcal P_{47}\) \\
\(\mathcal P_{47}\) & \((6,5,3,2;\ 5,4,3,2)\) & \(\mathcal P_{42}\) \\
\(\mathcal P_{48}\) & \((6,5,3,2;\ 5,4,2,1)\) & \(\mathcal P_{48}\) \\
\(\mathcal P_{51}\) & \((6,5,3,1;\ 5,5,3,1)\) & \(\mathcal P_{51}\) \\
\(\mathcal P_{52}\) & \((6,5,3,1;\ 5,4,3,1)\) & \(\mathcal P_{52}\) \\
\(\mathcal P_{53}\) & \((6,5,3,1;\ 5,4,2,1)\) & \(\mathcal P_{53}\) \\
\hline
\end{tabular}
\end{center}

Thus

\begin{equation}
\label{m6:s:eq:late-i1}
\mathcal L=\{g:R_{53}\le \mathsf{R}(g)\le R_{40}\}
=\mathcal P_{40}\sqcup \mathcal P_{41}\sqcup \mathcal P_{42}\sqcup \mathcal P_{47}
\sqcup \mathcal P_{48}\sqcup \mathcal P_{51}\sqcup \mathcal P_{52}\sqcup \mathcal P_{53}.
\end{equation}

The support filtration includes all eight cells of this interval, including the character-incompatible cells \(\mathcal P_{48},\mathcal P_{51},\mathcal P_{52}\).

Define

\begin{equation}
\label{m6:s:eq:late-i2}
V_i=\bigcap_{a,b}\{g:R_{ab}(g)\ge(R_i)_{ab}\},\qquad
\mathcal B=V_{53}\cap\bigcup_{a,b}\{g:R_{ab}(g)>(R_{40})_{ab}\}.
\end{equation}

Rank at least \(r\) is a union of nonvanishing loci of \(r\)-minors. The sets \(V_i\) and \(\mathcal B\) are therefore algebraic and analytic opens in the original \(G\). The symmetry of the two boundary matrices makes \(\mathcal B\) and \(V_{53}\) transpose-stable. Formula \eqref{m6:s:eq:late-i1} gives \(\mathcal L=V_{53}\setminus\mathcal B\), hence \(\mathcal L\) is closed in \(V_{53}\).

\subparagraph{An explicit transpose-stable support filtration}

Use the layers

\begin{equation}
\label{m6:s:eq:late-i3}
\begin{gathered}
E_1=\mathcal P_{40},\quad E_2=\mathcal P_{41},\quad
E_3=\mathcal P_{42}\sqcup\mathcal P_{47},\\
E_4=\mathcal P_{48},\quad E_5=\mathcal P_{51},\quad
E_6=\mathcal P_{52},\quad E_7=\mathcal P_{53}.
\end{gathered}
\end{equation}

Let \(\mathsf{S}_k\) be the identifiers in the first \(k\) layers, and put

\begin{equation}
\label{m6:s:eq:late-i4}
O_0=\varnothing,\qquad
O_k=\mathcal L\cap\bigcup_{i\in \mathsf{S}_k}V_i,\qquad
F_k=\mathcal L\setminus O_k.
\end{equation}

The table comparisons show that \(O_k\) is exactly the union of those selected layers. Therefore \(O_k\) is open in \(\mathcal L\), \(F_k\) is closed there, and
\(E_k=F_{k-1}\setminus F_k\) is open in the current closed remainder. Each \(O_k,F_k,E_k\) is transpose-stable.

The cell \(\mathcal P_{51}\) is incomparable with \(\mathcal P_{41}\), \(\mathcal P_{42}\), \(\mathcal P_{47}\) and \(\mathcal P_{48}\). Each open layer in \eqref{m6:s:eq:late-i4} consists of some of the currently maximal cells. The maximal-cell sets before each removal are

\begin{center}
\begin{tabular}{@{}p{0.07\linewidth} p{0.36\linewidth} p{0.23\linewidth} p{0.25\linewidth}@{}}
\hline
Step & Maximal cells in the current remainder & Layer removed & Input \\
1 & \(\mathcal P_{40}\) & \(\mathcal P_{40}\) & Intrinsic odd vanishing \\
2 & \(\mathcal P_{41}\), \(\mathcal P_{51}\) & \(\mathcal P_{41}\) & Even vanishing \\
3 & \(\mathcal P_{42}\), \(\mathcal P_{47}\), \(\mathcal P_{51}\) & \(\mathcal P_{42}\), \(\mathcal P_{47}\) & Even vanishing \\
4 & \(\mathcal P_{48}\), \(\mathcal P_{51}\) & \(\mathcal P_{48}\) & Compatibility sieve \\
5 & \(\mathcal P_{51}\) & \(\mathcal P_{51}\) & Compatibility sieve \\
6 & \(\mathcal P_{52}\) & \(\mathcal P_{52}\) & Compatibility sieve \\
7 & \(\mathcal P_{53}\) & \(\mathcal P_{53}\) & Intrinsic odd vanishing \\
\hline
\end{tabular}
\end{center}

In \(E_3\) the two cells are individually open and closed, and their \(\mathcal A\)-equivariant spaces vanish. Hence the entire paired layer has no even equivariant distribution, and in particular no odd one.

For \(\mathsf{S}_0\in D'_{\widetilde\eta}(\mathcal L)\), restrict to the open \(E_1\). Its intrinsic theorem gives zero. Sequence \eqref{m6:s:eq:late-c3} then gives a unique odd equivariant antecedent \(\mathsf{S}_1\) on \(F_1\). Restrict \(\mathsf{S}_1\) to the open \(E_2\), kill it, and take the next closed antecedent. Repeat through \eqref{m6:s:eq:late-i3}. The last remainder is \(\mathcal P_{53}\) and its odd theorem kills the last antecedent. Consequently

\begin{equation}
\label{m6:s:eq:late-i5}
{D'_{\widetilde\eta}(\mathcal L)=0.}
\end{equation}

At each step we use an open restriction or closed-support injection.
The uniqueness of the antecedent preserves oddness, since all
inclusions and opens are transpose-stable.

The even equivariant space on \(\mathcal L\) is nonzero. Indeed, \(\mathcal P_{53}\) is closed in \(\mathcal L\), since \(\mathcal L\subset V_{53}\) and \(\mathcal P_{53}\) is closed in \(V_{53}\). Its nonzero even orbital distribution has a closed pushforward to \(\mathcal L\).

\subparagraph{Conditional implication for the same original global distribution}

For a global odd equivariant \(T\), the sequence is

\begin{equation}
\label{m6:s:eq:late-i6}
0\longrightarrow D'_{\widetilde\eta}(\mathcal L)
\longrightarrow D'_{\widetilde\eta}(V_{53})
\longrightarrow D'_{\widetilde\eta}(\mathcal B).
\end{equation}

If \(T|_{\mathcal B}=0\), the restriction of this same \(T\) to \(V_{53}\) has a unique closed-support antecedent on \(\mathcal L\), killed by \eqref{m6:s:eq:late-i5}. Hence

\begin{equation}
\label{m6:s:eq:late-i7}
{T|_{\mathcal B}=0\ \Longrightarrow\ T|_{V_{53}}=0.}
\end{equation}

These calculations prove Proposition~\ref{m6:s:prop:outerlate} and the odd-vanishing theorem for the interval \(\mathcal L\).

\subsubsection{The late boundary cells}
\label{m6:s:sec:lateboundary}
\label{m6:s:sec:late}

We use the full group calculations of
\S\ref{m6:s:cert:p063-full}, \S\ref{m6:s:cert:p069-full}, and
\S\ref{m6:s:cert:p084-full}, together with the rational orbit and
transpose calculation for cell $69$ in \S\ref{m6:s:sec:p069}.

Two arguments handle the remaining three cells. On
\(\mathcal P_{63}\) and \(\mathcal P_{84}\), an extra additive subgroup
in the stabilizer has nontrivial character at every unipotently
compatible point. On \(\mathcal P_{69}\), the entire compatible locus is
one closed rational orbit, and its orbital line is even under transpose.

\begin{proposition}\label{m6:s:prop:lateboundary}
One has
\[
D'_\eta(\mathcal P_{63})=D'_\eta(\mathcal P_{84})=0,
\qquad \dim D'_\eta(\mathcal P_{69})=1,
\qquad D'_{\widetilde\eta}(\mathcal P_{69})=0.
\]
\end{proposition}

\paragraph{The cell \texorpdfstring{$\mathcal P_{63}$}{P63}}\label{m6:s:sec:p063}

The P063 row of Table~\ref{m6:s:tab:full-stabilizer-exclusions} applies on
the whole compatible family \eqref{m6:s:eq:late-63-a}. Its additive
subgroup \eqref{m6:s:eq:late-63-h} has nonzero phase on every allowed point.
\S\ref{m6:s:cert:p063-full} gives the group coordinates,
intersection, all stabilizer branches and modular/density check.
Thus \(D'_\eta(\mathcal P_{63})=0\).

\paragraph{Cell $69$: the second scalar-unipotent orbit}
The complete compatible coefficient family is
\[
 \mathsf P=\begin{pmatrix}0&p\\-b&t\end{pmatrix},\quad
 \mathsf Q=\begin{pmatrix}q&0\\d&c\end{pmatrix},\quad
 \mathsf R=\begin{pmatrix}0&q\\b&e\end{pmatrix},\quad
 \mathsf S=\begin{pmatrix}p&0\\f&g\end{pmatrix},\qquad pbqcg\ne0.
\]
All four free entries $t,d,e,f$ may vanish. The stabilizer is
$\mathbb G_m\times U_{15}$ throughout the family;
its two directions beyond the thirteen-dimensional kernel are
\[
 C_{21}+(q/c)B_{21}-(pq/(bc))A_{12},\qquad
 S_f+(p/g)T_h-(pq/(bg))Z_{12}.
\]
All denominators are nonzero on the coefficient family. The
stabilizer has zero algebraic additive character and relative
orbital density one.
The compatible variety is smooth, geometrically irreducible, closed
in the cell and of dimension $36$. Every geometric orbit has
dimension $52-16=36$; the split scalar-unipotent stabilizer has
trivial $H^1$, so there is exactly one rational orbit.

For the canonical coefficient point $(\mathsf P,\mathsf Q,
\mathsf R,\mathsf S)=(-j,I,M,I)$ with
$M=\left(\begin{smallmatrix}0&1\\1&0\end{smallmatrix}\right)$,
the identity is
\[
 x_{69,0}^t=d_5(D_0)x_{69,0}d_5(D_0)^t,
 \qquad D_0=\operatorname{diag}(1,-1).
\]
Both factors have character one, including their multiplier $-1$.
Proposition~\ref{m6:orbit-family} yields precisely one even orbital
line and zero odd space. The complete group and density calculation
appears on p.~\pageref{m6:s:sec:p069} of the appendices.

\paragraph{The cell \texorpdfstring{$\mathcal P_{84}$}{P84}}\label{m6:s:sec:p084}

The P084 row of Table~\ref{m6:s:tab:full-stabilizer-exclusions} uses the
subgroup \eqref{m6:s:eq:late-84-k} throughout the family
\eqref{m6:s:eq:late-84-a}. \S\ref{m6:s:cert:p084-full} contains the
intersection, uniform group stabilizer, nonzero character coefficient,
and complete density and transpose checks. They establish
\(D'_\eta(\mathcal P_{84})=0\) on every allowed rational branch.

\paragraph{Rank-open consequences}
For each \(i\in\{40,41,53,63,69,84\}\), define \(V_i\) as in \eqref{m6:s:eq:late-i2} and

\begin{equation}
\label{m6:s:eq:late-gsi1}
B_i=V_i\cap\bigcup_{a,b}\{g:R_{ab}(g)>(R_i)_{ab}\}.
\end{equation}

Since all six target rank matrices are symmetric, \(V_i,B_i\) are \(\mathcal A\)-invariant and transpose-stable ambient opens. The established complete rational invariant gives \(V_i\setminus B_i=\mathcal P_i\). Hence the target cell is closed in \(V_i\).

Assume \(T|_{B_i}=0\), for one fixed original global distribution. Sequence \eqref{m6:s:eq:late-c3} gives a unique equivariant distribution \(S_i\) on the closed cell such that \(i_*S_i=T|_{V_i}\). For the odd cases it is odd. The intrinsic theorem then gives

\begin{equation}
\label{m6:s:eq:late-gsi2}
{T|_{B_i}=0\ \Longrightarrow\ T|_{V_i}=0,\qquad
i=40,41,53,63,69,84.}
\end{equation}
The complete boundaries, obtained by entrywise comparison of the
rank matrices, are listed in Table~\ref{m6:s:tab:late-boundaries}.
The three local calculations prove Proposition~\ref{m6:s:prop:lateboundary}. Their rank-open consequences will be used in the global induction of Section~\ref{m6:s:sec:assembly}.

\subsubsection{The terminal graph orbit}\label{m6:s:sec:p090}

We finish with the graph orbit, its stabilizer and Haar line, and
the resulting terminal support step. To identify the orbit, we use
the signed coefficient dictionary and complete unipotent character
equations in \S\ref{m6:s:cert:p090-compatibility}, valid for arbitrary
invertible coefficients and all square classes.

\begin{proposition}\label{m6:s:prop:p090}
The terminal cell satisfies $D'_{\widetilde\eta}(\mathcal P_{90})=0$.
Its compatible locus is a single closed rational orbit $Hc$, and its
even-character equivariant distribution space is one-dimensional and
transpose-even.
\end{proposition}

\begin{proof}
Put
\[
 w_0=w_5\otimes I_2,\qquad
 b=\operatorname{diag}(I_2,-I_2,I_2,-I_2,I_2),\qquad c=bw_0.
\]
These matrices are symmetric involutions of multiplier one;
$b$ commutes with $w_0$ and both commute with $d_5(\GL_2)$.
The signed representative is related to the block reversal by
\begin{equation}\label{m6:s:eq:p090-3-1}
 W_{90}=d_*w_0,\qquad d_*=m(j,j,I_2)=\operatorname{diag}(j,j,I_2,j,j).
\end{equation}
Since $w_0P^tw_0^{-1}=P$, the cell is the closed translate
$\mathcal P_{90}=Pw_0$. In particular its upper/opposite
parabolic intersection is $P$. Its unique rational chart is
\begin{equation}\label{m6:s:eq:p090-3-4}
 U\times\GL_2^3\longrightarrow\mathcal P_{90},\qquad
 (u,A,B,C)\longmapsto u\,m(A,B,C)w_0.
\end{equation}
To obtain the inverse, read the three diagonal blocks of $gw_0$ and
remove their Levi factor. The chart therefore covers all rational
points and multiplier values. It also shows that the cell has
dimension $34$; transpose preserves it.

The 22-dimensional unipotent stabilizer and its character are
computed in \eqref{m6:s:eq:p090-4-2}--\eqref{m6:s:eq:p090-4-3}.
The two adjacent matrix blocks vary independently, so nondegeneracy
of their trace pairings gives exactly $B=-A$, $C=A$.
Scaling an elementary matrix tests the character itself, irrespective
of the conductor of $\psi$. In the signed chart the equations are
$B=-A$, $C=Aj$; the conversion is
\eqref{m6:s:eq:p090-3-1}. The coefficient dictionary, with arbitrary
invertible coefficients and all square classes, is printed in
\S\ref{m6:s:cert:p090-compatibility}.
It follows that the compatible set is
\begin{equation}\label{m6:s:eq:p090-5-1}
 \mathcal C=\{u\,m(A,-A,A)w_0\}=Hc.
\end{equation}
It is closed in $G$, smooth of dimension $26$, and has the global
rational section $hc\mapsto(h,1)$ for the orbit map.

Write $\lambda(ud_5(a))=\operatorname{tr}u_{12}+
\operatorname{tr}u_{23}$, so $\xi=\psi\circ\lambda$.
The algebraic anti-involution $\tau_H(k)=ck^tc$ preserves $H$ and
$\lambda$: the signs of $b$ cancel those from block reversal in
the two adjacent coordinates. Thus $\beta(k)=\tau_H(k)^{-1}$ is
an involutive automorphism with $\lambda(\beta(k))=-\lambda(k)$.
Solving the original stabilizer equation at every $hc$ gives
\begin{equation}\label{m6:s:eq:p090-6-1}
 K_{hc}=\{(h\beta(k)h^{-1},k):k\in H\}\simeq H.
\end{equation}

It retains the $\GL_2$ quotient, the 22-dimensional radical,
every scalar pair and every determinant value. The stabilizer is
geometrically connected of dimension $26$ throughout the family. The algebraic
character on this graph is zero because
\[
 \lambda(h\beta(k)h^{-1})+\lambda(k)=0.
\]
Each $M_2$ and $\mathfrak{sl}_2$ coordinate of $U$ transforms
under $d_5(a)$ by conjugation, of determinant one. Consequently
$H$, $\mathcal A$ and the graph $K_{hc}$ are unimodular;
the relative orbital density is exactly one.

At every point the transpose transporter is explicit:
\[
 (h^{-1},h)\cdot hc=ch^t=(hc)^t,\qquad
 \eta(h^{-1},h)=1.
\]
These facts verify the last row of
Proposition~\ref{m6:orbit-verification}. Its even line is concretely
\begin{equation}\label{m6:s:eq:p090-9-1}
 \mu_\xi(f)=\int_{H(F)}f(hc)\xi(h)\,dh.
\end{equation}
By the closed embedding and the homeomorphism $H\simeq Hc$, this
is a locally finite ordinary distribution. To see that it is nonzero,
choose a test supported on a compact open subgroup where $\xi=1$. Changing variables
$k=a_1h\tau_H(a_2)$ yields
$\mu_\xi(f\circ(a_1,a_2))=\eta(a_1,a_2)^{-1}\mu_\xi(f)$,
which fixes the convention. Under $H\simeq Hc$, transpose is the
single global map $\tau_H$. It preserves $\xi$ and sends positive
Haar measure to a positive scalar multiple of itself; involutivity
forces that scalar to be one. Thus the line is even.

On the complement of $\mathcal C$, the unipotent stabilizer
character is nontrivial, so ordinary equivariant distributions vanish.
Applying
Proposition~\ref{m6:orbit-family} and the closed-support sequence
proves exactly the stated odd vanishing and one-dimensional even
space. The nonzero even line also pushes forward to $G$ because
$Hc$ is closed there.
\end{proof}

The open complement is $\mathcal O_{89}=G\setminus\mathcal P_{90}$.
Hence the same left-exact sequence gives the terminal implication
\begin{equation}\label{m6:s:eq:p090-11-1}
 T|_{\mathcal O_{89}}=0\quad\Longrightarrow\quad T=0
 \qquad(T\in D'_{\widetilde\eta}(G)).
\end{equation}
Section~\ref{m6:s:sec:assembly} establishes this predecessor condition.

\subsection{The support order and the period theorem}\label{m6:filtration}
We assemble the local results along the twenty-three invariant opens
of Table~\ref{m6:s:tab:assembly-filtration}, using the original global
NW distribution at each conditional step. Duality then completes
the period theorem.
\subsubsection{Support induction}
\label{m6:short-assembly}
We prove Theorem~\ref{m6:distribution}. Fix the original equivariant odd
distribution $T\in\NW(G)$. Propositions~\ref{m6:s:prop:model7-open}
and~\ref{m6:rank-dictionary} give
\[
 T|_{\Omega\cup\Omega'}=0.
\]
In particular $T|_{\mathcal O_6}=0$, since
$\mathcal O_6\subset\Omega\cup\Omega'$. The intrinsic theorem for
cells $7,8$ then gives $T|_{\mathcal O_8}=0$.

For the first conormal argument, the boundaries are
\[
 B_9=B_{14}=\mathcal P[1,2,6,7]\subset\mathcal O_8.
\]
Proposition~\ref{m6:s:prop:p009-spectral} therefore gives vanishing
on $V_9\cup V_{14}$ and hence on $\mathcal S_9$.
Cells $10,11,15,21$ are character-incompatible, so the next open
$\mathcal S_{11}$ is also cleared. Now
\[
\begin{split}
 B_{12}&=\mathcal P[1,2,3,4,6,7,8,9,10,11],\\
 B_{22}&=\mathcal P[1,2,3,4,6,7,8,14,15,21]
\end{split}
\]
are contained in $\mathcal S_{11}$. Applying
Proposition~\ref{m6:s:prop:p012-spectral} to $T$ gives vanishing
on $\mathcal S_{12}$. These are the only remaining conormal
exclusions after the two Gauss opens.

The rest of the induction uses intrinsic vanishing. The increasing
invariant opens and their closed differences are listed in
Table~\ref{m6:s:tab:assembly-filtration}; their rank
descriptions and all boundary containments are proved in
Proposition~\ref{m6:s:prop:assembly-filtration}.
The middle-family theorem first treats cells $16,17,27$;
cell $34$ has already been removed with $\Omega'$.
The three-cell theorem treats $31,32,33$ together. On the outer
boundary the relevant closed difference is the entire interval
\[
 \mathcal P[40,41,42,47,48,51,52,53],
\]
to which the interval theorem of \S\ref{m6:s:sec:outerlate} applies.
The remaining nonempty compatible families are $63,69,84,90$;
their predecessors consist of previously removed families and
character-incompatible cells. The cited support table verifies each
containment, including the extra cell $77$ before the terminal step.
Applying the left-exact sequence of Lemma~\ref{cf:support} at each
step gives vanishing on $G$, as required.

\subsubsection{The opposite period and similitude-group duality}
Let $(\pi,V)$ be an irreducible admissible smooth complex representation of
$G$. Define
\begin{equation}\label{m6:s:eq:gk-periods}
\alpha(g)=g^{-t},\qquad
H_1=H,\quad \chi_1=\xi,\qquad
H_2=\alpha(H)=H^t,\quad \chi_2(\alpha(h))=\xi(h).
\end{equation}
Since $h^t=\alpha(h^{-1})$, this last equality is equivalently
\begin{equation}\label{m6:s:eq:gk-inverse}
\chi_2(h^t)=\xi(h)^{-1}.
\end{equation}
Put
\[
A_\pi=\operatorname{Hom}_H(\pi,\xi),\qquad
B_\pi=\operatorname{Hom}_{H^t}(\pi^\vee,\chi_2).
\]

We use the similitude-group dualizing theorem of Roche and Vinroot
\cite[Theorem~B and \S7]{RVfactor}. For a nondegenerate alternating space
over a non-Archimedean local field, the automorphism
$\iota_h(g)=\ell(g)^{-1}hgh^{-1}$, with $h$ an anti-symplectic involution,
takes every irreducible smooth representation to its contragredient.
Its character-theoretic proof applies to the characteristic-zero fields
here, including dyadic fields; its analytic character input is \cite[Theorem~1, \S3.3]{Clozel}. In our coordinates take
\[
h=d_5(\operatorname{diag}(1,-1)),\qquad h^2=1,\qquad \ell(h)=-1.
\]
The symplectic identity gives the comparison
\begin{equation}\label{m6:s:eq:gk-duality-comparison}
\alpha(g)=\ell(g)^{-1}JgJ^{-1}
=\operatorname{Ad}(Jh^{-1})\bigl(\iota_h(g)\bigr).
\end{equation}
The element $Jh^{-1}$ belongs to the ambient similitude group.
The difference between the two automorphisms is therefore inner in the
same group as $\pi$. Consequently
\begin{equation}\label{m6:s:eq:gk-duality}
\pi\circ\alpha\simeq\pi^\vee,\qquad
\pi^\vee\circ\alpha\simeq\pi.
\end{equation}
The central-character formula
$\pi^\vee\simeq\pi\otimes(\omega_\pi\circ\ell^{-1})$ in
\cite[\S7.3(c)]{RVfactor} is consistent with
\eqref{m6:s:eq:gk-duality-comparison}.

Pullback along $\alpha:H\to H^t$, followed by
\eqref{m6:s:eq:gk-duality}, gives
\begin{equation}\label{m6:s:eq:gk-period-identification}
B_\pi\simeq\operatorname{Hom}_H(\pi^\vee\circ\alpha,\xi)
\simeq\operatorname{Hom}_H(\pi,\xi)=A_\pi.
\end{equation}

\subsubsection{The individual multiplicity}
If $A_\pi=0$, Theorem~\ref{m6:main} follows immediately.
Otherwise \eqref{m6:s:eq:gk-period-identification} makes $B_\pi$ nonzero
and identifies its dimension with that of $A_\pi$.
The product bound therefore implies that both spaces are
one-dimensional. Equivalently,
\[
\left(\dim_{\mathbb C}\operatorname{Hom}_H(\pi,\xi)\right)^2\leq1.
\]
This proves Theorem~\ref{m6:main} for every irreducible admissible
smooth complex representation of $\operatorname{GSp}_{10}(F)$.
All algebraic coordinate changes and
finite exponentials used in the local proof are over $F$ in
characteristic zero, and the distributional and similitude-group
duality inputs have the stated dyadic scope.

\section{A framed orthogonal flag: Model 8}\label{cl:model8}\label{m8:section}

\paragraph{Proof strategy}
Frobenius descent over the external $\GL_2$ factor sends an
ambient NW distribution to an NW distribution on the identity
fiber. We prove vanishing on this fiber; injective reconstruction
then gives ambient vanishing. Using the eighteen rational strata $S_1,\ldots,S_{18}$ defined
in \S\ref{m8:intrinsic-families}, we first clear the rank
open with strata $S_1,S_2$, exclude the residual support on $S_3$
by an ambient conormal argument, and add the Siegel cell $S_8$.
Character vanishing on the exchanged pair $S_4,S_5$ and intrinsic
odd vanishing on $S_6,S_7$ then clear the initial open containing
the first eight strata.

We then add the remaining strata in the rank-open order
\[
 9,13,10,11,12,15,14,16,17,18.
\]
Most steps use stabilizer characters, matrix conjugation or
character-one transpose transport. At $S_{11}$ a second conormal exclusion applies once the
distribution vanishes on its predecessor open, which includes
$S_{10}$. Both $S_3$ and $S_{11}$ admit intrinsic odd
distributions; only the supports induced by the same ambient NW
distribution are excluded. For the intermediate steps exchanged by transpose we use ordinary
equivariant vanishing, while $S_{18}$ supplies the final odd step.
This support induction, assembled in \S\ref{m8:filtration}, proves
the distribution theorem. Frobenius reconstruction and connected
inverse-transpose duality then turn the GK product bound into
multiplicity one.

\Needspace{7\baselineskip}
\subsection{The split connected datum}\label{m8:datum}
\subsubsection{Preliminaries}\label{m8:s:sec:preliminaries}
\paragraph{The split connected orthogonal similitude group}
Every block in a four-by-four block matrix below has size two. Set
\begin{equation}\label{m8:s:p:sharp}
 J=\begin{pmatrix}0&-1\\1&0\end{pmatrix},\qquad
 A^\#=JA^tJ^{-1}=\tr(A)I_2-A.
\end{equation}
For two-by-two matrices,
\(AA^\#=A^\#A=(\det A)I_2,\) \((AB)^\#=B^\#A^\#,\) \(A^\#B+B^\#A=\tr(A^\#B)I_2.\)
The fixed form is
\begin{equation}\label{m8:s:p:Q}
 Q=\begin{pmatrix}
 0&0&0&J\\0&0&J&0\\0&-J&0&0\\-J&0&0&0
 \end{pmatrix},\qquad Q^t=Q,\qquad Q^2=I_8.
\end{equation}
Write
\(G_8^\circ=\{x\in\GL_8(F):x^tQx=\mu(x)Q,\) \(\det x=\mu(x)^4\}.\)
This determinant equality selects the connected similitude group,
of dimension $29$, whose Lie algebra decomposes as
\(\fg_8=F I_8\oplus\mathfrak{so}(Q),\) \(\mathfrak{so}(Q)=\{A:A^tQ+QA=0\}.\)
Since the trace pairing is nondegenerate in characteristic zero,
we use it to identify the Lie algebra with its dual. We specify
the translation convention when computing tangent and conormal spaces.

Let $P$ be the block upper triangular parabolic stabilizing the standard isotropic flag of dimensions $2$ and $4$. Its Levi subgroup is
\begin{equation}\label{m8:s:p:levi}
 m(A,E,t)=\diag\left(A,E,\frac{t}{\det E}E,
                       \frac{t}{\det A}A\right),
 \quad A,E\in\GL_2(F),\ t\in F^\times.
\end{equation}
We write this subgroup as $M$, and put $D(a)=m(a,a,\det a)=\diag(a,a,a,a)$. In particular $\mu(m(A,E,t))=t$ and $\mu(D(a))=\det a$.

The unipotent radical $U$ of $P$ consists of the matrices
\begin{equation}\label{m8:s:p:U}
 u(X,y,Z,z)=\begin{pmatrix}
 I_2&X&Z&zI_2+\tfrac12(XZ^\#-ZX^\#)\\
 0&I_2&yI_2&Z^\#-yX^\#\\
 0&0&I_2&-X^\#\\
 0&0&0&I_2
 \end{pmatrix},
 \quad X,Z\in\Mat_2(F),\quad y,z\in F.
\end{equation}
These ten parameters give affine coordinates over $F$, including
when the residue characteristic is two, since $1/2\in F$. The group law gives
\(X_{uu'}=X_u+X_{u'},\) \(y_{uu'}=y_u+y_{u'},\)
so
\(\xi(u(X,y,Z,z))=\psi(\tr X+2y)\)
is a character. We give the multiplication and inverse formulas below. Conjugation
by $D(a)$ replaces $X,Z$ by $aXa^{-1},aZa^{-1}$ and fixes $y,z$,
so it preserves $\xi$.

\paragraph{The Model 8 subgroup and the common center}
Set $B=\GL_2(F)$ and $\cG=G_8^\circ\times B$. The subgroup and character are
\(H=\{(uD(a),a):u\in U,\ a\in B\},\) \(\xi((uD(a),a))=\xi(u).\)
By the preceding conjugation formula, this extension is a group
character. It gives the split realization in \cite[\S8.4]{WZ},
which we use throughout; its scalar middle coordinate is $Y=yI_2$,
so $\tr Y=2y$.

The center relevant to periods is
\begin{equation}\label{m8:s:p:common-center}
 C_H=Z(\cG)\cap H=\{(tI_8,tI_2):t\in F^\times\},
 \qquad \xi|_{C_H}=1.
\end{equation}
Indeed, centrality in the second factor forces $a=tI_2$, and centrality in the first factor then forces $u=1$. If $\pi$ has a nonzero period, its central character is trivial on $C_H$. Writing its restrictions to the two scalar centers as $\omega_8$ and $\omega_2$, this condition is
\begin{equation}\label{m8:s:p:center-condition}
 \omega_8(t)\omega_2(t)=1\qquad(t\in F^\times).
\end{equation}

\paragraph{Transpose, opposite periods, and wavefront sets}
Put
\(\sigma(x,b)=(x^t,b^t),\) \(\alpha(x,b)=(x^{-t},b^{-t})=\sigma((x,b)^{-1}).\)
Transpose is an anti-involution; inverse transpose is an involutive automorphism. The equations defining $G_8^\circ$ give
\(xQx^t=\mu(x)Q,\) \(x^{-t}=\mu(x)^{-1}QxQ,\)
so both operations preserve the connected group. Because both factors of $\cG$ are unimodular, inversion preserves
Haar measure. The involutive automorphism $\alpha$ also preserves
it: its positive Haar modulus has square one. Thus both $\alpha$
and $\sigma$ preserve the chosen Haar measure.

The opposite subgroup is $H^\sigma=\sigma(H)=\alpha(H)$. It is not $H$ in these coordinates. Define
\begin{equation}\label{m8:s:p:opposite-character}
 \chi_2(k)=\xi(\alpha(k))=\xi(\sigma(k))^{-1}\quad(k\in H^\sigma).
\end{equation}
Then
\begin{equation}\label{m8:s:p:character-inverses}
 \chi_2(\sigma(h))=\xi(h)^{-1},\qquad
 \chi_2(\sigma(h)^{-1})=\xi(h),\qquad
 \chi_2(\alpha(h))=\xi(h).
\end{equation}
For an irreducible admissible $\pi$ the two period spaces are
\begin{equation}\label{m8:s:p:periods}
 \cP(\pi)=\Hom_H(\pi,\xi),\qquad
 \cQ(\pi)=\Hom_{H^\sigma}(\pi^\vee,\chi_2).
\end{equation}
Their identification is proved in \S\ref{m8:s:sec:duality}; the coefficient argument is given in \S\ref{cf:kernel}.

We use the full-cotangent wavefront convention of Definition~\ref{cf:native}, retaining the zero cotangent vector over support points. In an analytic chart, absence of a cotangent vector means that a cutoff nonzero at the point has Fourier transform vanishing asymptotically in a neighborhood of that direction; the precise definition and coordinate covariance are those of \cite[\S2.2]{JLX}.

For an ambient NW distribution on $G_8^\circ$, restriction to an
open subset retains its wavefront fibers. By conormal translation, as in Proposition~\ref{cf:conormal}, a distribution $R$ supported on a smooth closed $Y$ in an ambient open satisfies
\begin{equation}\label{m8:s:p:conormal-wf}
 N^*_{Y,y}\subset\WF_y(R)\qquad(y\in\supp R).
\end{equation}
Thus a nonnilpotent ambient conormal at every possible support point excludes that support for the original global distribution. The computations for $S_3$ and $S_{11}$ include the central and similitude directions and verify closedness in their stated opens after the specified predecessors vanish.

\subsubsection{The framed orthogonal filtration}
\label{m8:intrinsic-datum}
To describe the subgroup intrinsically, we first determine the
Levi stabilizer of the trace character and then use the external
framing to select its graph. The commutator calculation uses
\eqref{m8:s:t:eq1.3}, derived from \eqref{m8:s:p:U} in Appendix,
p.~\pageref{m8:s:t:sub1.2}, together with inversion and the Levi action.

In this intrinsic description $A,B$ denote planes and $L$ denotes
the similitude line. Let $V$ be the split quadratic space
with form $Q$, regarded as a
pairing with values in a one-dimensional similitude space $L$.
The parabolic $P$ stabilizes
\[
 0\subset E_2\subset E_4\subset E_6=E_2^\perp\subset V,
 \qquad E_4=E_4^\perp.
\]
Writing $A=E_2$ and $B=E_4/E_2$, its associated graded space is
\[
 A\oplus B\oplus(B^\vee\otimes L)
       \oplus(A^\vee\otimes L).
\]
Each summand has dimension two, and the quadratic form pairs the
first with the fourth and the second with the third.

\begin{proposition}[The intrinsic character datum]
\label{m8:covector}
There is a Levi-equivariant identification
\[
 U/[U,U]\simeq\operatorname{Hom}(B,A)
                  \oplus(\bigwedge^2 B\otimes L^{-1}).
\]
The character $\xi$ is the additive character associated with a
linear functional specified by an isomorphism $\alpha:A\longrightarrow B$
and a nonzero functional on $\bigwedge^2B\otimes L^{-1}$.
In the chosen coordinates it is $X,y\mapsto\operatorname{tr}X+2y$.
Its stabilizer in the Levi $M$ is exactly
\[
 M_\xi=\{m(a,a,\det a):a\in\GL_2(F)\}=D(\GL_2(F)).
\]
In the product Levi $M\times\GL_2$, the character stabilizer is
$M_\xi\times\GL_2$. The reductive part of $H$ is the graph of
$a\mapsto D(a)$, using the specified identification of the external
two-dimensional space with $B$.
\end{proposition}
\begin{proof}
The degree-one upper blocks are $X:B\to A$ and the skew map
$B^\vee\otimes L\to B$. The latter is
$\bigwedge^2B\otimes L^{-1}$; in dimension two it is the scalar
$y$-coordinate. For the remaining blocks we use the coordinates $Z,z$. By
\eqref{m8:s:t:eq1.3}, the coordinates $(X,y)$ add under multiplication. Commutators between an $X$-direction and the $y$-direction
span all $Z$-directions modulo the $z$-line. Commutators between
$X$ and $Z$ give
\[
 (X,Z)\longmapsto\operatorname{tr}(XZ^\#)
\]
on the $z$-line, a nondegenerate pairing. Hence $u\mapsto(X,y)$ has kernel $[U,U]$, both on $F$-points and
as an algebraic group, proving the displayed abelianization.

Through the trace pairing, the dual of $\operatorname{Hom}(B,A)$
is $\operatorname{Hom}(A,B)$, so the character is defined by the
asserted pair of linear functionals. Under $m(A_0,B_0,t)$ the two coordinates transform as
\[
 X\longmapsto A_0XB_0^{-1},\qquad
 y\longmapsto\frac{\det B_0}{t}y.
\]
Preservation of their independent linear functionals is equivalent
to $B_0^{-1}A_0=I$ and $t=\det B_0$. This is the stated stabilizer.
The external factor acts trivially on $U$, so the graph condition is
additional data supplied by its identification with $B$.
\end{proof}

\paragraph{Relative positions and the meaning of the rank parameter}
Let $F_b^-$ denote the standard opposite subspace spanned by the last
$b$ coordinates. For every $x\in G_8^\circ$,
\begin{equation}\label{m8:flag-ranks}
 r_{a,b}(x)=b-\dim(xF_b^-\cap E_{8-a}).
\end{equation}
Thus the rank matrix records the relative position of two self-dual
flags. Each flag carries the trace character and alternating form, and the external factor retains their common two-dimensional frame. After
Frobenius descent this is exactly why the residual reductive action
is simultaneous $D(a)$-conjugation. The ranks forget the characters and the frame; a rank stratum therefore generally contains a family of $K$-orbits.

In particular the parameter in the rank table is
\[
 c=r_{2,2}=2-\dim(xF_2^-\cap E_6)\in\{0,1,2\}.
\]
This discrete incidence rank records orthogonality: $c=0$ means
$xF_2^-\subset E_6$, or equivalently that $E_2$ and $xF_2^-$
are orthogonal. On this locus the intersection of the two unipotent
radicals changes, leading to the compatibility equations analyzed below. We distinguish this rank from the scalar entries of residual matrices
and the parameters of one-parameter degenerations.

The terminal cells have representatives normalizing $M$: they
permute the two-dimensional graded pieces as whole pieces. One
Levi parameter then describes their unipotent double quotient.
In the middle cells the joint graded pieces are one-dimensional;
mixed roots in $P\cap wP^tw^{-1}$ couple the two Levi parameters.
This filtration explains the different compatibility equations.
We solve them using the rational stabilizer and density calculations
below, and treat the two spectral cases by ambient-conormal exclusion.

In the coordinate arguments below, letters such as $A,B$ again
denote matrices when explicitly declared, and $\nu\in F^\times$
denotes a scalar multiplier in the Siegel chart.
Write $G_8=G_8^\circ\times\GL_2(F)=\mathcal G$, $H_8=H$,
$U_8=U$ and $\xi_8=\xi$ for the datum just defined. We also write
$S_r=S_{8,r}$ within this section.
\begin{theorem}[Model 8]\label{m8:main}
For every irreducible admissible smooth complex representation $\pi$
of $G_8$,
\[
 \dim\Hom_{H_8}(\pi,\xi_8)\leq1.
\]
Here $\F$ is a non-Archimedean local field of characteristic zero
with arbitrary residue characteristic, and $\psi$ is a nontrivial
continuous additive character. The individual central characters
are arbitrary; incompatibility on the common center makes the
period space zero.
\end{theorem}

We work with the connected split similitude point group and the
graph subgroup and common center in \eqref{m8:s:p:common-center}.

\subsection{The NW distribution theorem}\label{m8:distribution-class}
The ambient group is $G_8=G_8^\circ\times\GL_2(F)$, with
$\sigma_8(x,b)=(x^t,b^t)$. Frobenius descent over the external factor
has identity fiber $G_8^\circ$ and stabilizer
\[
 K_8=(U_8\times U_8)\rtimes\GL_2(F),\qquad
 (u,v,a)x=uD(a)xD(a)^{-1}v^t,
 \quad\eta_8(u,v,a)=\xi_8(u)\xi_8(v).
\]
Its relative density is trivial and its second reductive coordinate
is $a^{-t}$. Extend $\eta_8$ by value $-1$ on transpose.
The two nilpotent cones are
\begin{equation}\label{m8:native-cones}
 \begin{aligned}
 \cN(\Lie(G_8^\circ)^*)&=\cN(\mathfrak{so}_8^*)\times\{0\},\\
 \cN(\Lie(G_8)^*)&=\cN(\mathfrak{so}_8^*)
                    \times\cN(\mathfrak{sl}_2^*)\times\{0\}^2.
 \end{aligned}
\end{equation}
\begin{theorem}[The Model-8 distribution theorem]\label{m8:distribution}
Every ordinary, transpose-odd, $\eta_8$-equivariant distribution
$T\in\NW(G_8^\circ)$ is zero. Consequently every ordinary
$\sigma_8$-odd, two-sided Model-8 equivariant distribution
$S\in\NW(G_8)$ is zero.
\end{theorem}
We first prove the fiber assertion by open induction. For an ambient
$S\in\NW(G_8)$, the Frobenius correspondence below gives
$T\in\NW(G_8^\circ)$; once $T=0$, injectivity gives $S=0$. The two
spectral steps, at $S_3$ and $S_{11}$, concern the distributions
induced by the same ambient NW distribution after clearing their
respective predecessor opens.
\subsubsection{Frobenius descent}\label{m8:s:sec:frobenius}
\paragraph{The ambient two-sided action and its identity fiber}
Let $\mathcal L=H\times H$ act on $\cG$ by
\((h_1,h_2)g=h_1g\sigma(h_2),\) \(\eta_{\mathcal L}(h_1,h_2)=\xi(h_1)\xi(h_2).\)
Because $\sigma$ reverses multiplication, this is a left action. With
\((L_hf)(g)=f(h^{-1}g),\) \((R_kf)(g)=f(gk),\)
the inverse point action on tests is
\(f\longmapsto L_{h_1}R_{\sigma(h_2)^{-1}}f.\)
Thus the two translation identities
\begin{equation}\label{m8:s:f:translation-character}
 S(L_hf)=\xi(h)S(f),\qquad
 S(R_kf)=\chi_2(k)S(f)
 \quad(h\in H,\ k\in H^\sigma)
\end{equation}
give
\(S(L_{h_1}R_{\sigma(h_2)^{-1}}f) =\xi(h_1)\chi_2(\sigma(h_2)^{-1})S(f) =\xi(h_1)\xi(h_2)S(f).\)
The inverse character in \eqref{m8:s:p:opposite-character} therefore produces two plus signs in the descended character.

Under the projection $p:\cG\to B$, $p(x,b)=b$, the base action is
\((a,c)b=abc^t,\) which is transitive. At $b=I_2$, the stabilizer condition is $c=a^{-t}$. Writing $h_1=(uD(a),a)$ and $h_2=(vD(c),c)$, the action on that fiber is
\(K=(U\times U)\rtimes B,\) \((u,v,a)x=uD(a)xD(a)^{-1}v^t.\)
Here $B$ conjugates the first unipotent factor by $D(a)$ and the
second by $D(a^{-t})$, and the character is
\(\eta_K(u,v,a)=\xi(u)\xi(v) =\psi(\tr X_u+2y_u+\tr X_v+2y_v).\)
Transpose $\tau(x)=x^t$ normalizes this action through
\(\theta(u,v,a)=(v,u,a^{-t}),\) \(\eta_K\circ\theta=\eta_K.\)
Put $\widetilde K=K\rtimes\langle\tau\rangle$, with
\(\widetilde\eta_K|_K=\eta_K,\) \(\widetilde\eta_K(\tau)=-1.\)
For a $K$-stable $l$-space write $\De(X)=\Dist_{\eta_K}(X)$. If it is also transpose-stable, put $\Do(X)=\Dist_{\widetilde\eta_K}(X)$. The subscript $\mathrm e$ denotes $K$-equivariance. Individual strata exchanged by transpose are used with $\De$, and their stable unions may be used with $\Do$.

\paragraph{Modular characters and the base density}
The additive coordinate density $d^4X\,dy\,d^4Z\,dz$ on $U$ is Haar: left and right translations in the ordered unipotent coordinates have triangular Jacobian with diagonal entries one. In particular $U$ is unimodular. Conjugation by $D(a)$ acts on each of the two matrix spaces as $M\mapsto aMa^{-1}$, whose determinant is
\begin{equation}\label{m8:s:f:conjugation-det}
 (\det a)^2(\det a^{-1})^2=1.
\end{equation}
The scalar coordinates are fixed. Hence the diagonal subgroup preserves Haar measure on $U$.

On $B$, the density $db=|\det b|^{-2}d^4b$ is both left and right invariant: left or right multiplication by $a$ has determinant $(\det a)^2$ on $\Mat_2(F)$, canceled by the determinant factor. Therefore
\begin{equation}\label{m8:s:f:unimodular}
 \Delta_H=\Delta_{H^\sigma}=\Delta_{H\times H}=\Delta_K=1.
\end{equation}
Both conjugations defining $K$, by $D(a)$ and by $D(a^{-t})$,
have determinant one. Transpose preserves Haar measure as well,
since its automorphism has finite order and positive modulus.

For comparison the parabolic modulus is
\begin{equation}\label{m8:s:f:general-parabolic-modulus}
 \delta_P(m(A,E,t))=|\det A|^5|\det E|\,|t|^{-6}.
\end{equation}
The coordinate calculation on its Lie algebra is
\begin{equation}\label{m8:s:v:e9-1}
\begin{array}{c|c|c}
\text{coordinate}&\text{transformation}&\text{determinant}\\ \hline
X&AXE^{-1}&(\det A)^2(\det E)^{-2}\\
y&(\det E/t)y&\det E/t\\
Z&(\det E/t)AZE^{-1}&(\det A)^2(\det E)^2t^{-4}\\
z&(\det A/t)z&\det A/t.
\end{array}
\end{equation}
The product of their absolute values is \eqref{m8:s:f:general-parabolic-modulus}. Its restriction to $m(a,a,\det a)$ is one.

On the base, the stabilizer of $I_2$ acts by $b\mapsto aba^{-1}$,
whose tangent action on $\Mat_2(F)$ has determinant one by
\eqref{m8:s:f:conjugation-det}; the unipotent groups act trivially. Transpose fixes $E_{11},E_{22}$ and exchanges $E_{12},E_{21}$, so its determinant on the tangent space is $-1$ and its absolute density multiplier is one. Thus the base density line is trivial, including on the transpose
component, in agreement with the relative modulus in
\eqref{m8:s:f:unimodular}.

\paragraph{Global coordinates and the correspondence}
Consider the global algebraic isomorphism
\begin{equation}\label{m8:s:f:kappa}
 \kappa:G_8^\circ\times B\longrightarrow\cG,
 \qquad \kappa(y,b)=(yD(b),b),
 \qquad\kappa^{-1}(x,b)=(xD(b)^{-1},b).
\end{equation}
Under the second ambient diagonal element $(D(c),c)$, points transform as
\((x,b)\longmapsto(xD(c^t),bc^t),\) \((y,b)\longmapsto(y,bc^t)\)
in the new coordinates. Its character is one. Consequently, if $S$ is ambient equivariant, $\kappa^*S$ is invariant under every right translation in its $B$ coordinate.

Lemma~\ref{m8:s:c:product} with $A=B$, character one, gives a unique $T\in\Dist(G_8^\circ)$ with
\begin{equation}\label{m8:s:f:haar-factor}
 \kappa^*S=T\otimes db.
\end{equation}
The correspondence is defined by transport under the global isomorphism
\eqref{m8:s:f:kappa}.

Define
\begin{equation}\label{m8:s:f:averaging}
 \cA f(y)=\int_B f(yD(b),b)\,db.
\end{equation}
Since $f\circ\kappa$ is a finite sum of product test functions, $\cA f$ is a test function; its support lies in the projection of the compact support of $f\circ\kappa$. If $\beta_0\in\cD(B)$ has integral one, put
\(f_\phi(x,b)=\phi(xD(b)^{-1})\beta_0(b).\)
This is a test function and $\cA f_\phi=\phi$. Hence $\cA$ is surjective. Equation \eqref{m8:s:f:haar-factor} is exactly
\begin{equation}\label{m8:s:f:reconstruction}
 S(f)=T(\cA f),\qquad T(\phi)=S(f_\phi).
\end{equation}
The second formula is independent of $\beta_0$. Since $\cA$ is
surjective, it determines $T$ uniquely and makes
$T\mapsto T\circ\cA$ injective. The Haar-line argument accounts
for every ambient distribution with right diagonal invariance; in
particular, a zero descendant gives a zero ambient distribution.

We now establish all remaining equivariances in both directions. Left multiplication by $u\in U$ sends $(y,b)$ to $(uy,b)$, so
\begin{equation}\label{m8:s:f:left-U}
 T(y\mapsto\phi(u^{-1}y))=\xi(u)T(\phi).
\end{equation}
Left multiplication by $(D(a),a)$ has coordinate action
\((y,b)\longmapsto(D(a)yD(a)^{-1},ab).\)
Left invariance of $db$ and triviality of the diagonal character give
\begin{equation}\label{m8:s:f:diagonal}
 T(\phi\circ\Ad(D(a)^{-1}))=T(\phi).
\end{equation}

For $v\in U$ in the second acting factor the point action is $(x,b)\mapsto(xv^t,b)$. In coordinates it becomes
\((y,b)\longmapsto(yw_b^t,b),\) \(w_b=D(b^{-t})vD(b^t).\)
Indeed $w_b^t=D(b)v^tD(b)^{-1}$. The conjugation formula for $U$ gives
\begin{equation}\label{m8:s:f:right-U-character}
 \tr X_{w_b}=\tr X_v,\qquad y_{w_b}=y_v,
 \qquad\xi(w_b)=\xi(v).
\end{equation}
Testing ambient equivariance on $\phi(y)\beta(b)$ gives
\begin{equation}\label{m8:s:f:right-U-tests}
 \int_B\beta(b)T(y\mapsto\phi(yw_b^{-t}))\,db
 =\xi(v)T(\phi)\int_B\beta(b)\,db.
\end{equation}
Because $\phi$ has a compact open right stabilizer and
$b\mapsto w_b^{-t}$ is continuous, the translated test function
is locally constant in $b$. So is the scalar function
$b\mapsto T(y\mapsto\phi(yw_b^{-t}))$. Testing
\eqref{m8:s:f:right-U-tests} on sufficiently small compact open
neighborhoods therefore gives equality at every $b$. Evaluating
at $b=I_2$ yields
\begin{equation}\label{m8:s:f:right-U}
 T(y\mapsto\phi(yv^{-t}))=\xi(v)T(\phi).
\end{equation}

Conversely, suppose $T$ satisfies \eqref{m8:s:f:left-U}, \eqref{m8:s:f:diagonal}, and \eqref{m8:s:f:right-U}, and define $S=T\circ\cA$. The left $U$ and left diagonal calculations give the corresponding ambient equivariances. For the right $U$ action, apply \eqref{m8:s:f:right-U} to $w_b$ for each $b$ and use \eqref{m8:s:f:right-U-character}; integration gives \eqref{m8:s:f:right-U-tests}. Right diagonal invariance follows from $(y,b)\mapsto(y,bc^t)$ and right invariance of $db$. These four subgroups generate $H\times H$. We have proved the correspondence
\begin{equation}\label{m8:s:f:correspondence}
 \operatorname{Fr}_8=\operatorname{Fr}:\Dist_{\eta_{\mathcal L}}(\cG)
       \xrightarrow{\ \sim\ }\De(G_8^\circ),
 \qquad S\longmapsto T,
 \qquad\operatorname{Fr}^{-1}(T)=T\circ\cA.
\end{equation}

For fixed $b$, the coordinate change $y\mapsto yD(b)$ is right
translation on the unimodular group $G_8^\circ$, so it preserves
the product Haar density:
\begin{equation}\label{m8:s:f:product-density}
 \int_B\int_{G_8^\circ}f(yD(b),b)\,dy\,db
 =\int_B\int_{G_8^\circ}f(x,b)\,dx\,db.
\end{equation}
Thus \eqref{m8:s:f:correspondence} uses the product Haar densities.

\paragraph{Transpose compatibility}
The coordinate formula for transpose is
\begin{equation}\label{m8:s:f:transpose-coordinate}
 \kappa^{-1}\sigma\kappa(y,b)
 =\bigl(D(b^t)y^tD(b^t)^{-1},b^t\bigr).
\end{equation}
If $F=f\circ\kappa$, transport of $\sigma S$ therefore gives
\begin{equation}\label{m8:s:f:transpose-integral}
 (\kappa^*\sigma S)(F)
 =\int_B T\bigl(y\mapsto
 F(D(b^t)y^tD(b^t)^{-1},b^t)\bigr)\,db.
\end{equation}
The distribution $\tau T$ is invariant under every $\Ad(D(a))$, since
\(\tau\circ\Ad(D(a))=\Ad(D(a^{-t}))\circ\tau.\)
For each fixed $b$ we may thus remove the conjugation in \eqref{m8:s:f:transpose-integral} when evaluating $\tau T$. Transpose on $B$ preserves Haar measure, giving
\(\kappa^*(\sigma S)=\tau T\otimes db.\)
The families involved are compactly supported test families, so these evaluations are finite sums in product coordinates. Equivalently, the derivative of the coordinate change is block triangular: its first diagonal multiplier is that of an inner automorphism and its second is that of transpose on $B$, both of absolute value one. Hence
\begin{equation}\label{m8:s:f:transpose-intertwining}
 \operatorname{Fr}(\sigma S)=\tau\operatorname{Fr}(S).
\end{equation}
The even correspondence \eqref{m8:s:f:correspondence} restricts to an isomorphism of the odd spaces. Here the odd sign is the value $\widetilde\eta_K(\tau)=-1$.

\paragraph{Nilpotent wavefront under the external Frobenius map}\label{m8:fr-native}
This correspondence also has the one-way NW implication needed
here. By \eqref{m8:s:f:haar-factor} and the smooth nowhere-zero Haar
factor, local product cutoffs give
\begin{equation}\label{m8:s:f:native-product}
 \WF(T\otimes db)=
 \{((y,b),(\lambda,0)):(y,\lambda)\in\WF(T),\ b\in B\}.
\end{equation}
At $(y,I_2)$ the differential of the global coordinate is
\begin{equation}\label{m8:s:f:native-differential}
 d\kappa_{(y,I_2)}(v,w)=(v+y\,dD(w),w).
\end{equation}
Thus a cotangent vector $(\lambda,0)$ on the product-coordinate side
corresponds to an ambient cotangent vector $(\lambda,\beta)$, with
$\beta(w)=-\lambda(y\,dD(w))$. Its first component is exactly
$\lambda$. Coordinate covariance \cite[Lemma~2.15]{JLX} and NW
wavefront on the ambient product force this component to lie
in $\cN(\mathfrak{so}_8^*)\times\{0\}$. Hence
\begin{equation}\label{m8:s:f:native-transfer}
 S\in\NW(G_8^\circ\times\GL_2)
 \quad\Longrightarrow\quad
 \operatorname{Fr}_8(S)\in\NW(G_8^\circ).
\end{equation}
At $(y,I_2)$ the translated tangent convention agrees on the first
component, and the specified translations relate the other base
points to this one. Thus the global Haar factorization gives the
NW implication with the rational stabilizers, characters, transpose
corrections and densities established above.

\subparagraph{Two-factor coordinates.}\label{m8:s:v:s1-3}

Associate to \(k=(u,v,h)\) the pair
\begin{equation}
\bigl(g_L(k),g_R(k)\bigr)=\bigl(uD(h),vD(h^{-t})\bigr).
\label{m8:s:v:e1-3}
\end{equation}
Then \(kx=g_L(k)xg_R(k)^t\), multiplication is componentwise in these pairs, and \(\theta\) exchanges the two factors. We use this description for the right inverse, the residual subgroup on \(S_7\), and the variable-compensator formulas.

\subsection{Compatible relative positions and intrinsic descendants}
\label{m8:intrinsic-families}
We organize the eighteen self-dual relative positions by the rank
$c$ of the terminal two-by-two block, as in the table below. After treating the strata with
$c=2$ or $1$, we divide those with $c=0$ into the Siegel-open
cell $S_8$, the terminal strata $S_9,S_{13},S_{14},S_{17},S_{18}$,
and the middle strata $S_{10},S_{11},S_{12},S_{15},S_{16}$. On $S_8$ the affine-fiber character forces the
two residual matrices to agree, leaving conjugation of one matrix.
Whole graded-piece permutations form the terminal family of
Theorem~\ref{m8:s:t:main}; one-dimensional joint graded pieces form
the middle family of Theorem~\ref{m8:middle-family}. Within each family, the rational intersections and character
compatibility admit a uniform description; the residual actions
distinguish the two families.

\subsubsection{Algebra for the normal families}
\begin{lemma}[Levi--Weyl--Levi equivalence]\label{m8:s:c:normal-family}
For $P=LU$ and $w\in G_8^\circ$, put $J_w=P\cap wP^tw^{-1}$, with Levi projections $\pi_L,\pi_L^-$. Two points $lwr,l'wr'$ are $K$-equivalent exactly when some $h\in\GL_2(F)$ and $j\in J_w(F)$ satisfy
\begin{equation}\label{m8:s:c:normal-equations}
 l'=D(h)l\pi_L(j)^{-1},\qquad r'=\pi_L^-(w^{-1}jw)rD(h)^{-1}.
\end{equation}
\end{lemma}
\begin{proof}
An even transporter gives $j=l'^{-1}uD(h)l\in J_w$; taking both projections gives the equations. Conversely put
\begin{equation}\label{m8:s:c:normal-transporter}
 u=l'jl^{-1}D(h)^{-1},\qquad v^t=D(h)r^{-1}(w^{-1}jw)^{-1}r'.
\end{equation}
Both Levi projections are the identity; hence $u,v\in U$, and
substitution gives $uD(h)lwrD(h)^{-1}v^t=l'wr'$. This proves necessity and sufficiency with the required right inverse.
\end{proof}

We pass between Lie and group intersections using finite exponential
and logarithm. These are mutually inverse polynomials on the
characteristic-zero unipotent matrix groups below and commute with
conjugation and transpose; at each application we retain the
nonlinear exponential corrections. We also use the common affine calculation: for nonzero column $b$ and row $c$, the solutions of $pc+br=0$ are exactly
\begin{equation}\label{m8:s:c:affine-line}
 (p,r)=(\lambda b,-\lambda c),\qquad\lambda\in F.
\end{equation}
The semisimple and nilpotent applications below retain their different trace characters.
\begin{lemma}[Rational simultaneous Witt bases]\label{m8:simultaneous-witt}
Let $E_\bullet,F_\bullet$ be self-dual four-step flags in the split
quadratic space $(F^8,Q)$, with successive dimensions $2,4,6,8$.
In this lemma the indices refer to flag steps, so
$\dim E_i=\dim F_i=2i$ for $0\le i\le4$, with $E_0=F_0=0$.
If every nonzero joint graded multiplicity is one, the flags admit
a simultaneous adapted Witt basis over $F$. Two such flag pairs
with the same multiplicities are related by an $F$-rational isometry;
if their first flags coincide, that isometry has determinant one.
The construction gives rational local sections on nonzero-minor charts.
\end{lemma}
\begin{proof}
Choose a basis $e_1,\ldots,e_4$ of $E_2$ adapted to $E_1$ and
$E_2\cap F_j$, with joint grades $(a_i,b_i)$, $a_i\in\{1,2\}$.
The coordinate functional of $e_i$ vanishes on
$E_{a_i-1}+(E_2\cap F_{b_i-1})$. For
$W_i=E_{5-a_i}\cap F_{5-b_i}$, self-duality gives
\[
 W_i^\perp=E_{a_i-1}+F_{b_i-1},\qquad
 E_2\cap W_i^\perp=E_{a_i-1}+(E_2\cap F_{b_i-1}).
\]
By nondegeneracy we can solve for $f_i\in W_i$ with
$Q(e_j,f_i)=\delta_{ji}$ using a nonzero pairing minor. Replacing
\[
 f_i\longmapsto f_i'=f_i-\tfrac12\sum_jQ(f_i,f_j)e_j
\]
gives $Q(f_i',f_j')=0$ without changing the dual pairings.
The correction stays in $E_{5-a_i}$. A nonzero $Q(f_i,f_j)$
requires $(5-b_i)+(5-b_j)>4$, hence $b_j\le5-b_i$;
thus it stays in $F_{5-b_i}$ as well. The nonzero dual pairing
prevents a drop in the opposite joint grade. We thus obtain a simultaneous Witt basis, and matching two such
bases gives the asserted isometry. A flag-preserving isometry stabilizes $E_2$
and has diagonal blocks $A,A^{-t}$ in a Witt decomposition, so
its determinant is one. All choices use linear systems, nonzero
minors and $1/2$, giving rational local sections over $F$.
\end{proof}

\subsubsection{The inverse-corner rank identity}\label{m8:s:c:inverse-corner}
For an invertible matrix $g=\begin{psmallmatrix}A&B\\C&D\end{psmallmatrix}$ with four square blocks of the same size, \(\operatorname{rank}(g^{-1})_{11}=\operatorname{rank}(D)\): the kernel of \((g^{-1})_{11}\) identifies with \(\ker D\) by
\(v\in\ker D\longmapsto Bv,\) \(g(0,v)=(Bv,0).\)
Injectivity follows from invertibility of \(g\), and the converse follows by applying \(g\) to a vector whose first component is zero.

\subsubsection{Rank stratification}\label{m8:s:sec:ranks}
\paragraph{Rank invariants}
For $p,q\in\{6,4,2\}$, let $r_{p,q}(x)$ be the rank of the lower-right $p$-by-$q$ submatrix of $x$. The rank matrix, with its rows and columns ordered by $(6,4,2)$, is
\begin{equation}\label{m8:s:r:matrix}
 R(x)=\bigl(r_{p,q}(x)\bigr)_{p,q\in(6,4,2)}
 =\begin{pmatrix}4+c&2+e&b\\2+d&a&d\\b&e&c\end{pmatrix},
\end{equation}
where
\(c=r_{2,2},\) \(a=r_{4,4},\) \(b=r_{6,2}=r_{2,6},\) \(d=r_{4,2},\) \(e=r_{2,4}.\)
Left multiplication by a block upper triangular matrix performs invertible row operations on each terminal block of rows; right multiplication by a block lower triangular matrix performs invertible column operations on each terminal block of columns. Therefore these ranks are invariant under $P\times P^t$, and in particular under $K$. Transpose exchanges $d,e$ and fixes $c,a,b$.

For completeness, the complementary-rank identity for an invertible $n$-by-$n$ matrix $x$ is
\(\rank(x^{-1})_{I,J}=|I|+|J|-n+\rank x_{J^c,I^c}.\)
It follows by computing the kernel of the map from the coordinate subspace indexed by $J$ to the coordinate quotient indexed by $I$: that kernel is identified, under $x$, with the intersection of the corresponding coordinate subspace with the image of its complement. For $x\in G_8^\circ$, the identity $x^{-1}=\mu(x)^{-1}Qx^tQ$ turns complementary corners into terminal corners. In particular $r_{6,6}=4+c$, $r_{6,4}=2+e$, and $r_{4,6}=2+d$, giving \eqref{m8:s:r:matrix}. The rational cell descriptions in the following sections use the ranks themselves and include rational coverage.

\paragraph{The eighteen strata}
Let $S_i$ denote the locally closed locus with rank matrix $R_i$ in Table~\ref{m8:s:r:table}. For each of the eighteen possible rank patterns, the table records
the dimension and transpose permutation. By \eqref{m8:s:r:matrix},
the displayed five-tuple determines the whole rank matrix.

\begin{table}[htbp]
\centering
\caption{The complete rank stratification of $G_8^\circ$.}\label{m8:s:r:table}
\renewcommand{\arraystretch}{1.12}
\begin{tabular}{c c r r c}
\toprule
Stratum & $(c,a,b,d,e)$ & Dimension & Codimension & Transpose\\
\midrule
$S_1$ &(2,4,2,2,2)&29&0&$S_1$\\
$S_2$ &(2,2,2,2,2)&28&1&$S_2$\\
$S_3$ &(1,4,2,2,2)&28&1&$S_3$\\
$S_4$ &(1,2,2,1,2)&27&2&$S_5$\\
$S_5$ &(1,2,2,2,1)&27&2&$S_4$\\
$S_6$ &(1,2,2,1,1)&26&3&$S_6$\\
$S_7$ &(1,2,1,1,1)&25&4&$S_7$\\
$S_8$ &(0,4,2,2,2)&25&4&$S_8$\\
$S_9$ &(0,2,2,0,2)&24&5&$S_{13}$\\
$S_{10}$&(0,2,2,1,1)&25&4&$S_{10}$\\
$S_{11}$&(0,2,1,1,1)&24&5&$S_{11}$\\
$S_{12}$&(0,2,1,0,1)&23&6&$S_{15}$\\
$S_{13}$&(0,2,2,2,0)&24&5&$S_9$\\
$S_{14}$&(0,0,2,0,0)&23&6&$S_{14}$\\
$S_{15}$&(0,2,1,1,0)&23&6&$S_{12}$\\
$S_{16}$&(0,0,1,0,0)&22&7&$S_{16}$\\
$S_{17}$&(0,2,0,0,0)&20&9&$S_{17}$\\
$S_{18}$&(0,0,0,0,0)&19&10&$S_{18}$\\
\bottomrule
\end{tabular}
\end{table}

\paragraph{Rational completeness and the connected relative positions}
\label{m8:s:r:completeness}
We prove that the eighteen rank patterns in
Table~\ref{m8:s:r:table} exhaust $G_8^\circ(F)$.  Put $\bar\imath=9-i$ and
$\beta(i)=\lceil i/2\rceil$.  The diagonal torus of $G_8^\circ$ is
\[
 T=\{\diag(t_1,t_2,t_3,t_4,t/t_4,t/t_3,t/t_2,t/t_1)\}.
\]
The upper triangular subgroup $B_0$ is a split Borel contained in $P$.
Its Weyl group consists of the permutations $w$ of $\{1,\ldots,8\}$
such that $w(\bar\imath)=\overline{w(i)}$ and
$\#\{i\leq4:w(i)>4\}$ is even.  To verify the connectedness condition explicitly, write $q_i=Q_{i,\bar\imath}$; these are
$(-1,1,-1,1,1,-1,1,-1)$.  For any pair-preserving permutation,
define a monomial matrix $n_w$ by
\[
 n_we_i=e_{w(i)}\quad(1\leq i\leq4),\qquad
 n_we_{\bar\imath}=(q_i/q_{w(i)})e_{\overline{w(i)}}.
\]
It satisfies $n_w^tQn_w=Q$, and its determinant is
$(-1)^{\#\{i\leq4:w(i)>4\}}$: permutation of the four paired
coordinates has even sign, each interchange within a pair has odd
sign, and $\prod_{i=1}^4q_i/q_{w(i)}=1$.
More generally, a monomial similitude with this permutation and multiplier
$\mu$ has determinant $(-1)^{\#\{i\leq4:w(i)>4\}}\mu^4$, by
multiplying its four paired coefficient identities.  Thus exactly the
even pair-interchange permutations have representatives in the connected
group.  This is $W(D_4)$, of order $4!2^3$.
Conversely a normalizer of $T$ permutes its eight distinct coordinate
weight lines and respects their paired form, giving exactly these Weyl
elements.

Let $n_0e_i=e_{\bar\imath}$.  It lies in $G_8^\circ(F)$ and
$n_0Pn_0^{-1}=P^t$.  The split Bruhat decomposition, including its
$F$-rational factorization, gives
\[
 G_8^\circ(F)=\bigcup_{w\in W(D_4)} B_0(F)n_wB_0(F),\qquad
 G_8^\circ(F)=\bigcup_{w\in W(D_4)}P(F)n_wP^t(F).
\]
For the first equality one may use the multiplication isomorphism of
\cite[Theorem~21.80(a),(b)]{MilneAG}; taking $F$-points of that
isomorphism proves the asserted rational factorization.  The second equality follows by multiplying the first on the
right by $n_0$, using $B_0\subset P$ and $n_0B_0n_0^{-1}\subset P^t$.

For a Weyl permutation define the block-incidence matrix
\[
 N(w)_{ij}=\#\{k:\beta(k)=j,\ \beta(w(k))=i\}
 \qquad(1\leq i,j\leq4).
\]
It has nonnegative integral entries, every row and column sums to two,
and $N_{ij}=N_{5-i,5-j}$.  A terminal corner rank of $n_w$ is the
sum of the corresponding entries of $N(w)$.  Consequently the five
ranks in \eqref{m8:s:r:matrix} determine $N(w)$ uniquely:
\begin{equation*}\tag{M8-R1}\label{m8:s:r:incidence}
 N(w)=
 \begin{pmatrix}
 c&e-c&b-e&2-b\\
 d-c&a-d-e+c&2-a-b+d+e&b-d\\
 b-d&2-a-b+d+e&a-d-e+c&d-c\\
 2-b&b-e&e-c&c
 \end{pmatrix}.
\end{equation*}
Conversely, nonnegativity in this matrix is exactly
\begin{equation*}\tag{M8-R2}\label{m8:s:r:inequalities}
 0\leq c\leq d,e\leq b\leq2,
 \qquad d+e-c\leq a\leq2-b+d+e.
\end{equation*}
The number of pair interchanges is $4-a$, since $a$ counts those of
the four coordinates in the first half whose images remain in that
half.  Hence connectedness is exactly the additional condition
$a\in2\mathbb Z$.

Every integral matrix in \eqref{m8:s:r:incidence} satisfying
\eqref{m8:s:r:inequalities} and this parity condition is attained.
Indeed, regard coordinates $1,2$ and $3,4$ as two source groups of
paired coordinates, and similarly regard the four target pairs as
two target groups.  Place $N_{ij}$ positive edges and $N_{5-i,j}$
negative edges from source group $j$ to target group $i$, for
$i,j\in\{1,2\}$.  Every source and target group has degree two by
the column sums and centrosymmetry.  Give its two edges the two
distinct coordinate labels in that group.  A positive edge sends its
source label to its target label; a negative edge sends it to the
paired target label.  Extend by $w(\bar\imath)=\overline{w(i)}$.
This constructs a permutation with incidence matrix $N$, with an
even number $4-a$ of pair interchanges, and the above $n_w$ is an
$F$-rational representative in $G_8^\circ$.

For clarity, solving the inequalities gives the following complete
list; entries in the last column are $(d,e;a)$.
\[
\begin{array}{c|c|l}
 b&c& (d,e;a)\\\hline
 0&0&(0,0;0),\ (0,0;2)\\
 1&0&(0,0;0),\ (0,1;2),\ (1,0;2),\ (1,1;2)\\
 1&1&(1,1;2)\\
 2&0&(0,0;0),\ (0,2;2),\ (1,1;2),\ (2,0;2),\ (2,2;4)\\
 2&1&(1,1;2),\ (1,2;2),\ (2,1;2),\ (2,2;4)\\
 2&2&(2,2;2),\ (2,2;4)
\end{array}
\]
The list is precisely the eighteen tuples in
Table~\ref{m8:s:r:table}, proving exhaustion. We can also read
off the rational parabolic double cosets from the incidence matrix.
The Levi Weyl group $W_M$ is generated by $(12)(78)$ and
$(34)(56)$ and has order four.  Two pair-preserving permutations
with the same incidence matrix are related by left and right
elements of $W_M$: match their edges separately for each source
group, target group and sign in the construction above.  Matching
the source labels and the target labels gives permutations within
each two-element group; their extensions to paired labels are
exactly elements of $W_M$; different monomial lifts with the same
permutation differ by $T(F)\subset M(F)$.  Conversely these left and right
permutations preserve incidence.  Thus $N$ identifies
$W_M\backslash W(D_4)/W_M$, and, by the rational coverage and rank
invariance already proved,
\[
 S_i(F)=P(F)n_iP^t(F)\qquad(1\leq i\leq18),
\]
for any of the above representatives having the tuple of $S_i$.
Thus each listed rank pattern is a single rational parabolic double coset.  Transpose replaces $N$ by $N^t$, which exchanges
$d,e$ and gives exactly the transpose column of the table.

To compute the dimensions directly, put
\[
 C(w)=\#\{(i,j):i\ne j,\ i\ne\bar\jmath,
 \ \beta(i)\leq\beta(j),\
 \beta(w^{-1}(i))\geq\beta(w^{-1}(j))\}.
\]
The off-diagonal matrix positions $(i,j)$ and
$(\bar\jmath,\bar\imath)$ describe the same one-dimensional root
space in $\mathfrak{so}(Q)$.  Therefore
$\dim(P\cap n_wP^tn_w^{-1})=5+C(w)/2$, where five is the torus
dimension, and
\[
 \dim S_i=2\dim P-5-C(w)/2=33-C(w)/2.
\]
We obtain the dimensions in the table over every characteristic-zero
$F$, including residue characteristic two.

The transpose pairs are precisely $S_4\leftrightarrow S_5$, $S_9\leftrightarrow S_{13}$, and $S_{12}\leftrightarrow S_{15}$. Every other stratum is stable. An exact-rank locus is locally closed, since rank upper bounds are closed determinantal conditions and rank lower bounds are open.

\paragraph{Support strategy}
For a rank matrix $R_i$ write
\(\Omega_i=\{x\in G_8^\circ:R(x)\geq R_i\},\)
where the inequality is entrywise. As a finite intersection of rank lower-bound conditions, $\Omega_i$
is open; imposing the additional upper bounds $R(x)\leq R_i$ makes
$S_i$ closed in it. The notation $\Omega_2$ and $\Omega_4$ used for the two initial maximal-rank regions below will always be given by their explicit defining equations; in particular the initial $\Omega_4=\{r_{4,4}=4\}$ denotes the corner-rank open.

The first two useful regions are
\(\Omega_2=\{\det x_{44}\ne0\}=S_1\sqcup S_2,\) \(\Omega_4=\{r_{4,4}=4\}=S_1\sqcup S_3\sqcup S_8.\)
Here $x_{44}$ is the terminal two-by-two block, and $x_{44}\ne0$
means that at least one of its entries is nonzero.

We prove intrinsic odd vanishing on $S_1,S_2,S_6,S_7,S_8,S_{10},S_{18}$ and intrinsic even-group vanishing on $S_4,S_5,S_9,S_{12},S_{13},S_{14},S_{15},S_{16},S_{17}$. The even vanishings already follow from $U\times U$ in the indicated cases. On $S_3$ and $S_{11}$ the odd spaces are nonzero, as is
$\De(S_{18})$. We therefore use the ambient wavefront input from
\S\ref{m8:s:sec:preliminaries} to eliminate the support induced by
the original global distribution at the two spectral steps. This
proves the theorem for ambient NW distributions, which we later
apply to generalized coefficients. We assemble the open filtration
and its support sequences in \S\ref{m8:s:sec:assembly}.

\subsubsection{The open stratum}\label{m8:s:sec:s1}
\begin{proposition}\label{m8:s:s1:vanishing}
The open stratum has intrinsic odd vanishing: $\Do(S_1)=0$.
\end{proposition}
\begin{proof}
Write $a^*=a/\det a$ and $m(a,b,t)=\diag(a,b,tb^*,ta^*)$ as in \eqref{m8:s:p:levi}. The open rank profile in Table~\ref{m8:s:r:table} gives the $F$-product chart
\begin{equation}\label{m8:s:s1:factorization}
 S_1=ULU^t,\qquad U\times L\times U\longrightarrow S_1,
 \quad(p,m,q)\longmapsto pmq^t,
 \quad L=\GL_2^2\times F^\times.
\end{equation}
For uniqueness, $pmq^t=p'm'(q')^t$ gives $(p')^{-1}p=m'(q')^tq^{-t}m^{-1}$. The first side is upper unitriangular, the second lower triangular with diagonal $m'm^{-1}$; comparison gives $m'=m$, then $p'=p,q'=q$. Thus the two-unipotent fibers are free, with dimensions $10+9+10=29$.

The actions are $(p,m,q)\mapsto(up,m,q)$ and $(p,m,q)\mapsto(p,m,vq)$, the latter because $pmq^tv^t=pm(vq)^t$. With the same Haar measure on both copies of $U$, set
\begin{equation}\label{m8:s:s1:average}
 \mathcal A_1f(m)=\int_{U^2} f(pmq^t)\xi(p)\xi(q)\,dp\,dq.
\end{equation}
In the product chart, this is a compactly supported locally constant
function. Choosing a test $\beta_0$ with $\int_U\beta_0\xi=1$ gives
the right inverse $\phi\mapsto\beta_0\otimes\phi\otimes\beta_0$. Lemma~\ref{m8:s:c:product} therefore gives the exact weak-star correspondence
\(T(f)=R(\mathcal A_1f),\) \(R(\phi)=T(\beta_0\otimes\phi\otimes\beta_0).\)
For $h\in\GL_2$, the chart action is
\((p,m,q)\longmapsto (D(h)pD(h)^{-1},D(h)mD(h)^{-1},D(h)^{-t}qD(h)^t).\)
On each matrix coordinate of $U$, conjugation has determinant $(\det h)^2(\det h^{-1})^2=1$; it fixes the scalar coordinates and preserves $\xi$. Thus
\begin{align}
 D(h)m(a,b,t)D(h)^{-1}&=m(hah^{-1},hbh^{-1},t),\label{m8:s:s1:diagonal-levi}\\
 \mathcal A_1(f\circ(1,1,h)^{-1})(m)&=\mathcal A_1f(D(h)^{-1}mD(h)).\label{m8:s:s1:diagonal-intertwining}
\end{align}
Transpose sends $pm(a,b,t)q^t$ to $qm(a^t,b^t,t)p^t$, using $(a^*)^t=(a^t)^*$. Exchange of the two identical weighted Haar factors gives
\(\mathcal A_1(f\circ\tau)(m)=\mathcal A_1f(m^t).\)
We are left with a distribution invariant under simultaneous
conjugation of $(a,b)$ and odd under their simultaneous transpose,
with $t$ passive. Applying Lemmas~\ref{m8:s:c:two-matrix} and
\ref{m8:s:c:passive} gives $R=0$, hence $T=0$. We use this
vanishing in the open support argument of
\S\ref{m8:s:a:entrance-section}.
\end{proof}

\subsubsection{The invertible-final-block boundary}\label{m8:s:sec:s2}
We use the big-cell factorization to descend the character and
reduce odd vanishing to the two-matrix problem. The rational chart
and distributional argument follow below; Appendix~\S\ref{m8:s:app:B}
gives the expanded stabilizer matrices for this action.

\begin{proposition}\label{m8:s:s2:vanishing}
The stratum $S_2$ has intrinsic odd vanishing: $\Do(S_2)=0$.
\end{proposition}

\paragraph{The final-plane chart and its rational inverse}
Put $\mathcal Q_4=\begin{psmallmatrix}0&J\\-J&0\end{psmallmatrix}$ and
$G_4=\{g:g^t\mathcal Q_4g=\mu(g)\mathcal Q_4,\ \det g=\mu(g)^2\}$.
In sizes $2+4+2$, the fixed form is $Q=\begin{psmallmatrix}0&0&J\\0&\mathcal Q_4&0\\-J&0&0\end{psmallmatrix}$. Let $Q_2$ stabilize the final two-plane and $Q_2^+$ the initial two-plane; they are lower and upper triangular, respectively. Their product Levi is
\(M_2=\{m(a,g)=\diag(a,g,\mu(g)a^*):a\in\GL_2(F),\ g\in G_4\} \simeq\GL_2(F)\times G_4.\)
Set $T=zI_2+\tfrac12(XZ^\#-ZX^\#)$ and
\begin{equation}\label{m8:s:s2:e1-4}
 n(X,Z,z)=u(X,0,Z,z)=
 \begin{pmatrix}I_2&X&Z&T\\0&I_2&0&Z^\#\\0&0&I_2&-X^\#\\0&0&0&I_2\end{pmatrix}.
\end{equation}
Then $N=R_u(Q_2^+)$ is this nine-dimensional group and $R_u(Q_2)=N^t$. Indeed, for $B=(X\ Z)$ the grouped upper unipotent matrix has form equations
\begin{equation}\label{m8:s:s2:e1-6}
 W=-\mathcal Q_4B^tJ=\binom{Z^\#}{-X^\#},\qquad
 T-T^\#=BW,
\end{equation}
so $T=\tfrac12\tr(T)I_2+\tfrac12BW$. Its law and inverse are
\begin{equation}\label{m8:s:s2:e1-7}
\begin{aligned}
 n(X,Z,z)n(X',Z',z')&=n\left(X+X',Z+Z',z+z'+\tfrac12\tr(XZ'^\#-ZX'^\#)\right),\\
 n(X,Z,z)^{-1}&=n(-X,-Z,-z).
\end{aligned}
\end{equation}
Translation has triangular Jacobian one, so $dX\,dZ\,dz$ is left and right Haar.
With $v(y)=u(0,y,0,0)$, multiplication gives
\begin{align}
 u(X,y,Z,z)&=n(X,Z-yX,z)v(y),\qquad U=N\rtimes\{v(y)\},\label{m8:s:s2:e1-9}\\
 v(y)n(X,Z,z)v(y)^{-1}&=n(X,Z-yX,z),\label{m8:s:s2:e1-10}\\
 \zeta_N(n(X,Z,z))&=\psi(\tr X),\qquad \xi(v(y))=\psi(2y).\label{m8:s:s2:e1-11}
\end{align}
These retain the fixed coefficient $2$ in the Model 8 character \cite[\S8.4]{WZ}.

\begin{proposition}\label{m8:s:s2:big-cell}
Multiplication $N\times M_2\times N\to\Omega_2=\{\det x_{44}\ne0\}$, $(n_L,m,n_R)\mapsto n_Lmn_R^t$, is an algebraic isomorphism over $F$.
\end{proposition}
\begin{proof}
Use grouped indices $1,M,4$ and put $d=x_{44}$. The inverse radical coordinates are
\begin{align}
 W_L&=x_{M4}d^{-1},&T_L&=x_{14}d^{-1},&
 W_R&=(d^{-1}x_{4M})^t,&T_R&=(d^{-1}x_{41})^t,\label{m8:s:s2:e2-2}\\
 X_\bullet&=-W_{\bullet,2}^\#,&Z_\bullet&=W_{\bullet,1}^\#,&
 z_\bullet&=\tfrac12\tr T_\bullet,&n_\bullet&=n(X_\bullet,Z_\bullet,z_\bullet).\label{m8:s:s2:e2-3}
\end{align}
The final-column isotropy equation is $T_L^tJ+W_L^t\mathcal Q_4W_L-JT_L=0$. With $B_L=-JW_L^t\mathcal Q_4$ this is $T_L-T_L^\#=B_LW_L$, verifying \eqref{m8:s:s2:e1-6}. The identity $xQx^t=\mu(x)Q$, applied to $x^t$, verifies the right coordinates by the same final-column equation under transpose.

Now $n_L^{-1}xn_R^{-t}=\begin{psmallmatrix}A&B_0&0\\C_0&G&0\\0&0&d\end{psmallmatrix}$. Its form equations give $A^tJd=\mu(x)J$ and $B_0^tJd=0$, hence $A=\mu(x)d^*$ and $B_0=0$. Invertibility gives $G\in\GL_4$, and $C_0^t\mathcal Q_4G=0$ gives $C_0=0$. Finally $G^t\mathcal Q_4G=\mu(x)\mathcal Q_4$ and $\det G=\mu(x)^2$, since $\det A\det d=\mu(x)^2$ and $\det x=\mu(x)^4$. Thus $G\in G_4$, with
\(g=x_{MM}-x_{M4}d^{-1}x_{4M},\) \(a=\mu(x)d^*.\)
Conversely $n_Lm(a,g)n_R^t$ has final block $\mu(g)a^*$. Its final rows and columns recover the two radicals, then the Levi; all inverse formulas are regular where $\det d\ne0$. The radical fibers are free.
\end{proof}
The lower-right four-block rank of $m(a,g)$ is $2+\rank g_{22}$. Consequently
\(S_2\simeq N\times\mathcal Y_2\times N,\) \(\mathcal Y_2=\{(a,g)\in\GL_2\times G_4:g_{22}=0\}.\)
Invertibility then forces $g_{12},g_{21}$ invertible and gives the remaining ranks in Table~\ref{m8:s:r:table}.

\paragraph{The concrete distribution and density data}
With equal Haar measures on the two copies of $N$, Lemma~\ref{m8:s:c:product} applies to
\begin{equation}\label{m8:s:s2:e3-2}
 \mathcal A f(m)=\int_{N^2}f(n_Lmn_R^t)\psi(\tr X_L+\tr X_R)\,dn_L\,dn_R.
\end{equation}
The chart gives compact support, and two fiber tests with weighted integral one give a right inverse. Substituting $n_L=pn_L'$, $n_R=qn_R'$ gives the two positive characters. The surviving group is $\mathcal R=F^2\times\GL_2$, with
\begin{equation}\label{m8:s:s2:e3-6}(s,t,h)(a,g)=\left(hah^{-1},n_4(s)d_4(h)g d_4(h)^{-1}n_4(t)^t\right),\qquad
 \chi_{\mathcal R}(s,t,h)=\psi(2s+2t),
\end{equation}
where $n_4(s)=\begin{psmallmatrix}I_2&sI_2\\0&I_2\end{psmallmatrix}$ and $d_4(h)=\diag(h,h)$. The lift $(v(s),v(t),h)$ conjugates the left radical by $v(s)D(h)$ and the right by $v(t)D(h^{-t})$. The shear \eqref{m8:s:s2:e1-10} and $(X,Z,z)\mapsto(hXh^{-1},hZh^{-1},z)$ both preserve Haar and $\zeta_N$.

The action of $M_2$ on $N/Z(N)$ is $B\mapsto aBg^{-1}$, with determinant $(\det a)^4(\det g)^{-2}=(\det a)^4\mu(g)^{-4}$. On the center its scaling factor is $\det a/\mu(g)$. Hence
\begin{equation}\label{m8:s:s2:e3-10}\delta_+(m)=\left|\frac{\det a}{\mu(g)}\right|^5,\qquad
 \delta_{Q_2}=\delta_+^{-1},\qquad
 dx=\delta_+(m)^{-1}\,dn_L\,dm\,dn_R.
\end{equation}
With $\int f(xb)dx=\Delta(b)^{-1}\int f(x)dx$, the group modular functions of $Q_2^+,Q_2$ are $\delta_+^{-1},\delta_+$ respectively. The ambient density follows by left-trivializing the differential: its first radical factor is $\Ad(m^{-1})$, the others triangular with determinant one. The acting groups $N,N^2,M_2,\mathcal R,K$ are unimodular, and the radical fiber stabilizer is trivial. Both $\det a$ and $\mu(g)$ are residual and transpose invariants, so the nonconstant ambient density changes only invariant normalization. Intrinsic stratum distributions use the product-chart density.

Transpose sends $(n_L,m,n_R)$ to $(n_R,m^t,n_L)$ and $(a,g)$ to $(a^t,g^t)$, since $(a^*)^t=(a^t)^*$. It exchanges $s,t$ and sends $h$ to $h^{-t}$, preserving the character and the two identical fiber Haar lines. Thus
\begin{equation}\label{m8:s:s2:e3-15}
 \Do(S_2)\simeq\Dist_{\widetilde\chi_{\mathcal R}}(\mathcal Y_2).
\end{equation}

\paragraph{The rational tensor quotient}
Since $\mathcal Q_4=(-J)\otimes J$, the map
\(p(A,B)=A\otimes B,\) \(\mu(p(A,B))=\det A\det B,\) \(\ker p=\{(zI_2,z^{-1}I_2):z\in F^\times\}\)
maps $\GL_2^2$ to $G_4$: the form scales by $\det A\det B$ and the determinant is its square. The kernel follows by comparing tensor blocks. This is a central torus quotient.

To prove surjectivity on $F$-points, proceed as follows. Nonzero isotropic tensors in $F^2\otimes F^2$ are decomposable. Their maximal linear spaces are the two rulings $\ell\otimes F^2$ and $F^2\otimes\ell$. A map preserving both rulings has $g(e_i\otimes v)=p_i\otimes B_iv$, with $p_1,p_2$ independent and $B_i$ invertible. Decomposability of the image of $(e_1+e_2)\otimes v$ forces $B_1v,B_2v$ collinear for every $v$, hence $B_2=cB_1$ over $F$. Thus $g=A\otimes B$. The tensor flip has multiplier one and determinant $-1$, so a ruling-exchanging map has determinant $-\mu^2$ and is excluded by connectedness.

On a chart with block $g_{ij}$ invertible and $(g_{ij})_{uv}\ne0$, the section is
\(B=g_{ij},\) \(a_{ij}=1,\) \(a_{kl}=(g_{kl})_{uv}/(g_{ij})_{uv}.\)
Since $g_{kl}=a_{kl}B$, these charts cover every rational point.
Thus $(\GL_2^2)/F^\times\simeq G_4$, with the stated analytic
local sections.

Writing $u_2(s)=\begin{psmallmatrix}1&s\\0&1\end{psmallmatrix}$ gives $n_4(s)=p(u_2(s),I_2)$ and $d_4(h)=p(I_2,h)$. Upstairs the actions are
\begin{align}
 (s,t,h)(a,A,B)&=(hah^{-1},u_2(s)Au_2(t)^t,hBh^{-1}),\label{m8:s:s2:e4-7}\\
 \tau(a,A,B)&=(a^t,A^t,B^t),\qquad
 z(a,A,B)=(a,zA,z^{-1}B),\quad\chi_Z=1.\label{m8:s:s2:e4-8}\end{align}
For $q(a,A,B)=(a,p(A,B))$ the fiber averaging is
\begin{equation}\label{m8:s:s2:e4-11}q_!\phi(a,g)=\int_{F^\times}\phi(a,zA,z^{-1}B)\,d^\times z,
 \qquad q^\natural S(\phi)=S(q_!\phi).
\end{equation}
Using the local sections and a compactly supported torus test of
integral one, we see that $q_!$ is surjective and $q^\natural$
injective. Lemma~\ref{m8:s:c:fiber} identifies the image of
$q^\natural$ with all invariant distributions upstairs. The actions commute with the kernel and transpose fixes $z$; the multiplicative fiber measure is preserved and has character and modulus one. Thus this is an untwisted equivariant correspondence.

\paragraph{Compatibility and the two-matrix problem}
The preimage of $g_{22}=0$ has $A=\begin{psmallmatrix}\alpha&\beta\\\gamma&0\end{psmallmatrix}$, $\beta\gamma\ne0$. The global kernel normalization $z=\beta^{-1}$ gives unique coordinates
\begin{equation}\label{m8:s:s2:e5-2}r=\alpha/\beta,\quad\kappa=\gamma/\beta,\quad C=\beta B,
 \qquad g(r,\kappa,C)=\begin{pmatrix}rC&C\\\kappa C&0\end{pmatrix},\quad
 \mu(g)=-\kappa\det C.
\end{equation}
Indeed $C=g_{12}$, $g_{11}C^{-1}=rI_2$, $g_{21}C^{-1}=\kappa I_2$. Every point has the unique form
\begin{equation}\label{m8:s:s2:e5-5}
 x=n_L\diag\bigl(a,g(r,\kappa,C),-\kappa\det(C)a^*\bigr)n_R^t,
 \quad a,C\in\GL_2, r\in F, \kappa\in F^\times.
\end{equation}
Its dimension is $9+(4+1+1+4)+9=28$. Multiplication and transpose give
\begin{align}
 (s,t,h)(a,r,\kappa,C)&=(hah^{-1},r+\kappa s+t,\kappa,hCh^{-1}),\label{m8:s:s2:e5-6}\\
 \tau(a,r,\kappa,C)&=(a^t,r/\kappa,\kappa^{-1},\kappa C^t).\label{m8:s:s2:e5-7}
\end{align}
The latter follows by transposing the two-by-two tensor factor and renormalizing its upper off-diagonal entry; it preserves $\mu$ and squares to one. The complete additive stabilizer and its character are
\begin{equation}\label{m8:s:s2:e5-8}E_x=\{(s,-\kappa s):s\in F\},\qquad
 \chi_{\mathcal R}(s,-\kappa s,1)=\psi(2(1-\kappa)s).
\end{equation}
The residual stabilizer is $F\times Z_{\GL_2}(a,C)$, with the same line. It is unimodular: the common centralizer is $\GL_2$ when both matrices are scalar, and otherwise is contained in a commutative nonscalar two-by-two centralizer. Uniqueness of the chart determines its complete ordinary lift
(Appendix~\S\ref{m8:s:app:B}).

For Haar $da,dC$, additive $dr$, and multiplicative $d^\times\kappa$, transpose has Jacobian $|\kappa|^{-1}$: only $r\mapsto r/\kappa$ contributes, while inversion and scalar multiplication preserve the multiplicative Haar factors. Thus
\(|\kappa|^{-1/2}\,da\,dr\,d^\times\kappa\,dC\)
is invariant under both residual action and transpose. The central coordinate is recovered as $A=\beta\begin{psmallmatrix}r&1\\\kappa&0\end{psmallmatrix}$, $B=\beta^{-1}C$; its density $d^\times\beta$ is preserved by $\beta\mapsto z\beta$ and by transpose $\beta\mapsto\kappa\beta$.

On $\kappa\ne1$, the character in \eqref{m8:s:s2:e5-8} is nontrivial on the complete unipotent stabilizer, whose relative modulus is one. Applying \eqref{cf:orbital} gives vanishing there. We then use
Lemma~\ref{cf:support} with $X=\mathcal Y_2$ and the closed stable
subset $\mathcal Y_{2,1}=\{\kappa=1\}$ to confine the support to
that subset. Intrinsically this is
\(\mathcal C_2=\{x\in S_2:x_{23}-x_{24}d^{-1}x_{43} =x_{32}-x_{34}d^{-1}x_{42}\},\) \(d=x_{44},\)
the equality of the off-diagonal Schur blocks. It is closed, $K$-stable and transpose-stable.

On this support, $(s,t,h)$ acts by $(a,r,C)\mapsto(hah^{-1},r+s+t,hCh^{-1})$ with character $\psi(2s+2t)$; transpose is $(a,r,C)\mapsto(a^t,r,C^t)$. The ineffective line $(s,-s)$ has character and modulus one, so the quotient additive fiber has exactly $\psi(2r)dr$ by Lemma~\ref{m8:s:c:product}. The exact reconstruction is
\begin{equation}\label{m8:s:s2:e6-8}\mathcal I(S_0)(f)=S_0\left((a,C)\longmapsto\int_F f(a,r,1,C)\psi(2r)\,dr\right),
\end{equation}
where $S_0$ is invariant under simultaneous conjugation of $(a,C)$ and odd under $(a,C)\mapsto(a^t,C^t)$. By the weighted product lemma, this correspondence holds in both
directions and reconstructs the test functions. Applying
Lemma~\ref{m8:s:c:two-matrix} gives $S_0=0$, which proves
Proposition~\ref{m8:s:s2:vanishing} through \eqref{m8:s:s2:e3-15}.
For the composed reconstruction see Appendix,
\S~\ref{m8:s:app:A}; we use the open support consequence in
\S\ref{m8:s:a:entrance-section}.

\subsubsection{\texorpdfstring{The paired strata \(S_4\) and \(S_5\)}{The paired strata S4 and S5}}\label{m8:s:sec:s4-s5}
The main proof gives rational coverage and the character obstruction
on both strata. The complete block conjugates cited from
Appendix~\S\ref{m8:s:app:B} are essential verification of the
intersection equations and of the stabilizer direction used here.

\begin{theorem}\label{m8:s:s45:vanishing}
The strata \(S_4\) and \(S_5\) satisfy
\(D'_{\eta_K}(S_4)=D'_{\eta_K}(S_5)=0,\) \(D'_{\widetilde\eta_K}(S_4\sqcup S_5)=0.\)
The two individual vanishings hold already under \(U\times U\).
\end{theorem}
We use $L(A,B,t)=m(A,B,t)$ from \eqref{m8:s:p:levi} and the fixed $U,K,\eta_K,\theta$. The proof uses the complete two-unipotent stabilizer after rational normalization.
\paragraph{Rational rank geometry}\label{m8:s:s45:sub1}
\subparagraph{The rectangular rank conditions}\label{m8:s:s45:sub1.2}

Let \(r_{m,n}(x)\) denote the lower-right \(m\times n\) rank, and retain \(x_{ij}\) for the two-by-two blocks. On \(S_4\), put
\[
d=x_{44},\qquad
H=\begin{pmatrix}x_{33}&x_{34}\\x_{43}&x_{44}\end{pmatrix},
\qquad S=(x_{43}\ \ x_{44}).
\]
Then
\begin{equation}\label{m8:s:s45:eq1.7}
\operatorname{rank}d=1,\qquad \operatorname{rank}H=2,
\qquad\operatorname{rank}S=2.
\end{equation}
Conversely, inside the fixed \(G_8^\circ\), \eqref{m8:s:s45:eq1.7} implies exactly the \(S_4\) rank matrix. Indeed, the first two block rows of \(H\) are uniquely expressible as
\begin{equation}\label{m8:s:s45:eq1.8}
(x_{33}\ \ x_{34})=T(x_{43}\ \ x_{44})
\end{equation}
for a matrix \(T\in M_2(F)\). Hence
\(\operatorname{rank}\binom{x_{34}}{x_{44}}=1\). The complementary-rank identities give the other entries
\[
R(x)=\begin{pmatrix}5&4&2\\3&2&1\\2&2&1\end{pmatrix}.
\]

Thus \(x_{34}\) has rank at most one and annihilates \(\ker x_{44}\), while the two bottom block rows in \eqref{m8:s:s45:eq1.8} have the specified incidence. The rank of \(x_{43}\) can be one or two. We therefore use a nonzero two-column minor of the rectangular matrix \(S\).

For \(S_5\), transpose gives the analogous statement with the lower-right block column of \(H\) of rank two and a unique factorization of the left block column through it.

The terminal two-by-two block and the lower-right four-by-four block are singular, so the normalization below uses only the rectangular minors in \eqref{m8:s:s45:eq1.7}.

The pair itself is the locally closed, \(K\)- and transpose-stable set described by
\begin{equation}\label{m8:s:s45:eq1.9}
\begin{gathered}
r_{2,2}=1,\quad r_{4,4}=2,\quad r_{6,2}=r_{2,6}=2,\\
(r_{4,2},r_{2,4})\in\{(1,2),(2,1)\}.
\end{gathered}
\end{equation}
Its two components are clopen; their ambient opens occur in \S\ref{m8:s:a:entrance-section}.

\paragraph{A rational slice from the rank conditions}\label{m8:s:s45:sub2}
\begin{proposition}\label{m8:s:s45:slice}
Every point of \(S_4(F)\) has the unique rational coordinate factorization \eqref{m8:s:s45:eq2.12}, with parameters \eqref{m8:s:s45:eq2.10}. The inverse is regular on the nonzero-minor charts and these inverses agree on overlaps.
\end{proposition}
\begin{proof}\leavevmode

\subparagraph{Eliminating the dependent rows}\label{m8:s:s45:sub2.1}

Set
\(\ell(X)=u(X,0,0,0),\) \(c(z)=u(0,0,0,z).\)
Choose any two columns \(I\subset\{1,2,3,4\}\) for which the two-by-two matrix \(S_I\) is invertible. Formula \eqref{m8:s:s45:eq1.8} is explicitly
\begin{equation}\label{m8:s:s45:eq2.1}
T=(x_{33}\ \ x_{34})_I\,S_I^{-1}.
\end{equation}
Different choices of \(I\) give the same matrix. Left multiplication by \(\ell(T^\#)\) makes blocks \((3,3)\) and \((3,4)\) zero.

Let \(M\) denote the first block row in the last four columns after this operation. Orthogonality of those last four columns says
\begin{equation}\label{m8:s:s45:eq2.2}
M^tJS-S^tJM=0.
\end{equation}
Because \(S:F^4\to F^2\) is surjective, \eqref{m8:s:s45:eq2.2} first implies \(\ker S\subset\ker M\). Thus \(M=Z_0S\), and then
\(Z_0^tJ=JZ_0.\)
For a two-by-two matrix this means \(Z_0^\#=Z_0\), hence \(Z_0=zI_2\). Therefore
\(M=zS,\) \(zI_2=M_I S_I^{-1}.\)
Left multiplication by \(c(-z)\) now makes blocks \((1,3)\) and \((1,4)\) zero as well.

\subparagraph{Orthogonal block elimination}\label{m8:s:s45:sub2.2}
Set
\begin{equation}\label{m8:s:s45:eq2.4}
 \mathcal J=\begin{pmatrix}0&J\\J&0\end{pmatrix},\qquad
 \mathcal L=\begin{pmatrix}0&J\\-J&0\end{pmatrix},
 \qquad\mathcal J^{-1}=-\mathcal J,\quad\mathcal L^{-1}=\mathcal L.
\end{equation}
In the temporary output row order $(1,3\mid2,4)$, the matrix
just obtained has the form
\begin{equation}\label{m8:s:s45:eq2.5}
 \begin{pmatrix}g&0\\C_0&b\end{pmatrix},\qquad
 g,b\in\GL_4(F),\quad
 g^t\mathcal Lb=\nu\mathcal J,\quad b=\nu\mathcal Lg^{-t}\mathcal J.
\end{equation}
Invertibility gives the two invertible diagonal blocks.
With
\begin{equation}\label{m8:s:s45:eq2.7}
 V=-b^{-1}C_0,\quad\mathcal JV=V^t\mathcal J,\quad
 \begin{pmatrix}I_4&0\\V&I_4\end{pmatrix}=n_S(V^t)^t,
\end{equation}
the last identity identifies the correction with the transpose
of the six-dimensional radical in $U$, with
$n_S(W)=(\begin{smallmatrix}I_4&W\\0&I_4\end{smallmatrix})$; here
$V^t=(\begin{smallmatrix}Z&z_1I_2\\y_1I_2&Z^\#\end{smallmatrix})$.
Right multiplication eliminates $C_0$. Thus the full rational slice is
\begin{equation}\label{m8:s:s45:eq2.9}
 x(g,\nu)=\begin{pmatrix}
 g_{11}&g_{12}&0&0\\0&0&b_{11}&b_{12}\\
 g_{21}&g_{22}&0&0\\0&0&b_{21}&b_{22}
 \end{pmatrix},\qquad b=\nu\mathcal Lg^{-t}\mathcal J,
\end{equation}
with exactly the conditions
\begin{equation}\label{m8:s:s45:eq2.10}
 g\in\GL_4(F),\quad\operatorname{rank}g_{22}=1,\quad\nu\ne0.
\end{equation}
Indeed $b_{22}=-\nu J(g^{-t})_{11}J$ and the inverse-corner
identity of \S\ref{m8:s:c:inverse-corner} gives the rank condition
in both directions. Conversely the slice satisfies the form equation,
all required ranks and $\det x=\det g\det b=\nu^4$;
the row permutation has sign one. The normalization therefore
remains in the specified connected group.
The complete block substitutions are in
\texttt{anc/model8\_calculations.pdf}, calculation M8-C09,
under ``The triangular form and its elimination''.

\subparagraph{Uniqueness of the coordinate factorization}\label{m8:s:s45:sub2.3}

The construction yields
\begin{equation}\label{m8:s:s45:eq2.12}
x=\ell(T^\#)^{-1}c(z)\,x(g,\nu)\,n_S(V^t)^{-t}.
\end{equation}
The factors \(T,z,g,\nu,V\) are uniquely recovered by \eqref{m8:s:s45:eq2.1}--\eqref{m8:s:s45:eq2.7}. They are regular on the minor charts used above and agree on overlaps. Conversely, arbitrary values of
\(T\in M_2,\) \(z\in F,\) \((g,\nu)\text{ as in \eqref{m8:s:s45:eq2.10}},\) \(V^t\text{ as in \eqref{m8:s:s45:eq2.7}}\)
reconstruct a point of \(S_4\). This gives a global product chart whose inverse is expressed locally by the displayed minors. Its dimension is
\(4+1+15+1+6=27.\)
Thus it gives rational local sections over $F$.

\end{proof}
\paragraph{Normal forms and rational equivalence}\label{m8:s:s45:sub3}
\begin{proposition}\label{m8:s:s45:normal-forms}
The family \eqref{m8:s:s45:eq3.2} meets every \(U(F)\times U(F)\)-orbit in \(S_4(F)\). The unused orbit-equivalence catalogue is retained in Appendix~\S\ref{m8:s:app:B}.
\end{proposition}
\begin{proof}\leavevmode

\subparagraph{The representative and the normal family}\label{m8:s:s45:sub3.1}

Let \(E_{ij}\) be the two-by-two matrix units and set
\begin{equation}\label{m8:s:s45:eq3.1}
w=
\begin{pmatrix}
E_{11}&E_{21}&0&0\\
0&0&E_{21}&-E_{11}\\
E_{22}&E_{12}&0&0\\
0&0&-E_{12}&E_{22}
\end{pmatrix}.
\end{equation}
Its nonzero entries, in ordinary eight-by-eight notation, are
\(w_{11}=w_{23}=w_{45}=w_{54}=w_{62}=w_{88}=1,\) \(w_{37}=w_{76}=-1.\)
The displayed entries give \(w^tQw=Q\), \(\det w=1\), and the \(S_4\) rank matrix.

A covering family of final normal forms is
\begin{equation}\label{m8:s:s45:eq3.2}
x_{A,B,C,E;t,s}=L(A,B,t)\,w\,L(C,E,s),
\end{equation}
where
\(A,B,C,E\in \operatorname{GL}_2(F),\) \(t,s\in F^\times.\)
These are the parameter conditions, and the multiplier is $ts$. 

\subparagraph{Rational coverage by unipotent operations}\label{m8:s:s45:sub3.2}

We now give further reductions of the coordinate \(g\) in \eqref{m8:s:s45:eq2.9}. Left multiplication of \(g\) by \(n(Z)=\begin{psmallmatrix}I_2&Z\\0&I_2\end{psmallmatrix}\) is realized on the original eight-by-eight slice by \(u(0,0,Z,0)\). Right multiplication by \(\overline n(R)=\begin{psmallmatrix}I_2&0\\R&I_2\end{psmallmatrix}\) is realized by \(\ell(R^t)^t\). Thus these are even-group operations.

Factor the rank-one terminal block as
\(g_{22}=uv,\)
with a nonzero column \(u\) and row \(v\), both over \(F\). Choose complementary vectors to obtain
\[
L_0=(u\ \ u')\in \operatorname{GL}_2(F),\qquad
R_0=\binom{v'}v\in \operatorname{GL}_2(F).
\]
Then \(g_{22}=L_0E_{12}R_0\). Writing
\(g_1=\operatorname{diag}(I_2,L_0^{-1})\,g\, \operatorname{diag}(I_2,R_0^{-1})\)
is a change of factor coordinates.

If \(g_1=\begin{psmallmatrix}a_1&b_1\\c_1&E_{12}\end{psmallmatrix}\), the explicit unipotents
\(P_1=-(b_1e_2)e_1^t,\) \(R_1=-e_2(e_1^tc_1)\)
remove the second column of \(b_1\) and the first row of \(c_1\). The resulting upper-left block is
\(a_0=a_1-(b_1e_2)(e_1^tc_1),\)
and the matrix is
\[
\begin{pmatrix}
a_0&b_*e_1^t\\
e_2c_*&E_{12}
\end{pmatrix},
\qquad b_*=b_1e_1\ne0,\quad c_*=e_2^tc_1\ne0.
\]
Choose
\[
T_0=(u_0\ \ b_*),\qquad S_0=\binom{v_0}{c_*}
\]
in \(\operatorname{GL}_2(F)\). Write
\[
T_0^{-1}a_0S_0^{-1}=\begin{pmatrix}\alpha&\beta\\\gamma&\delta\end{pmatrix}.
\]
Invertibility of \(g\) forces \(\alpha\ne0\). The remaining unipotents, with the normalized column and row
\(\widehat p=(-\beta,-\delta)^t,\) \(\widehat r=(-\gamma,0),\)
remove the last three entries, leaving \(\alpha E_{11}\). In the unnormalized coordinates use \(p=T_0\widehat p\), \(r=\widehat rS_0\), and \(P=pe_2^t\), \(R=e_1r\). Absorb \(\alpha\) into the first column of \(T_0\).

We have obtained, up to the unipotents described above,
\begin{equation}\label{m8:s:s45:eq3.5}
g=\operatorname{diag}(A_0,L_0)
\begin{pmatrix}E_{11}&E_{21}\\E_{22}&E_{12}\end{pmatrix}
\operatorname{diag}(C_0,R_0).
\end{equation}
The block-diagonal factor-coordinate changes normalize the corresponding upper and lower unipotent groups, so their conjugates are still precisely the allowed \(n(Z)\) and \(\overline n(R)\). Their conjugates therefore give the required operations in \(U\times U\).

Taking \(t=1\), \(B=L_0^*\), \(A=A_0\), \(C=C_0\), \(E=R_0\), and \(s=\nu\) turns \eqref{m8:s:s45:eq3.5} into \eqref{m8:s:s45:eq3.2}; the other four-by-four block agrees by the unique form equation \eqref{m8:s:s45:eq2.5}. In particular, \(t\) can be fixed to one for coverage. Allowing arbitrary \(t\) is useful for symmetric formulas below.

All choices use a nonzero minor and an \(F\)-linear complement. They give rational sections on finite minor charts. This proves coverage of every \(F\)-point under \(U\times U\subset K\).

\end{proof}
\paragraph{The complete unipotent stabilizer and character}\label{m8:s:s45:sub4}
\begin{proposition}\label{m8:s:s45:unipotent-stabilizer}
The \(U\times U\)-stabilizer of a normal point \eqref{m8:s:s45:eq3.2} is \eqref{m8:s:s45:eq4.2}. It is an additive group of dimension two, with character \eqref{m8:s:s45:eq4.5}, which is nontrivial at every normal point.
\end{proposition}
\begin{proof}

Fix a point \(x=x_{A,B,C,E;t,s}\) as in \eqref{m8:s:s45:eq3.2}, and define
\begin{equation}\label{m8:s:s45:eq4.1}
\kappa=\frac{\det B}{t},\qquad
F_0=C^tE_{21}E^{-t},\qquad T_0=CE^{-1}.
\end{equation}
Then \(\kappa\ne0\).

\subparagraph{Solution of the stabilizer equations}\label{m8:s:s45:sub4.1}

The stabilizer inside \(U\times U\) is
\begin{equation}\label{m8:s:s45:eq4.2}
\begin{aligned}
u_{y,r}&=u\left(0,\kappa y,\kappa rAE_{22}B^{-1},0\right),\\
v_{y,r}&=u\left(-rF_0,0,-\frac{\det E}{s}yF_0,0\right),
\qquad y,r\in F.
\end{aligned}
\end{equation}
Every pair in \eqref{m8:s:s45:eq4.2} satisfies \(u_{y,r}xv_{y,r}^t=x\), and gives every stabilizer pair. The parameters identify it with $\mathbb G_a^2$.

To solve the stabilizer completely, write
\[
n(X,y,Z,z)=
\begin{pmatrix}
0&X&Z&zI_2\\0&0&yI_2&Z^\#\\0&0&0&-X^\#\\0&0&0&0
\end{pmatrix}\in\mathfrak u,
\]
imposing \(w^{-1}n(X,y,Z,z)w\in\mathfrak u^t\) gives exactly
\begin{equation}\label{m8:s:s45:eq4.3}
X=0,\quad z=0,\quad Z=rE_{22},\quad y\text{ arbitrary}.
\end{equation}
Moreover,
\begin{equation}\label{m8:s:s45:eq4.4}
\bigl(w^{-1}n(0,y,rE_{22},0)w\bigr)^t
=n(rE_{21},0,yE_{21},0).
\end{equation}
The matrices in \eqref{m8:s:s45:eq4.3} and \eqref{m8:s:s45:eq4.4} square to zero. The common finite exponential/logarithm argument in \S\ref{sec:common} therefore identifies the entire group intersection. Conjugating by \(l=L(A,B,t)\), and by \(r_0^t=L(C^t,E^t,s)\) on the right, with the required inverse on the right, yields precisely \eqref{m8:s:s45:eq4.2}.

\subparagraph{The stabilizer character}\label{m8:s:s45:sub4.2}

Since
\(\operatorname{tr}(C^tE_{21}E^{-t})=(CE^{-1})_{21},\)
formula \eqref{m8:s:s45:eq4.2} gives
\begin{equation}\label{m8:s:s45:eq4.5}
\eta_K(u_{y,r},v_{y,r},1)
=\psi\left(2\kappa y-(T_0)_{21}r\right).
\end{equation}
The inverse in the right stabilizer equation gives the minus sign
in the second term. Since the right \(y\)-coordinate and left
\(X\)-coordinate vanish, and the character is trivial on the
right \(Z\)-coordinate, these are all its contributions.

The \(r\)-line is character-trivial exactly when \((T_0)_{21}=0\). But the \(y\)-line is always nontrivial. Triviality on the stabilizer would require
\begin{equation}\label{m8:s:s45:eq4.6}
2\kappa=0,\qquad (T_0)_{21}=0,
\end{equation}
and the first equation is impossible. Accordingly, the compatible locus is empty, including where the \(r\)-line is character-trivial.

\end{proof}
\paragraph{Density and vanishing}\label{m8:s:s45:sub6}
The algebraic group $U\times U$ and its complete additive stabilizer \eqref{m8:s:s45:eq4.2} are unimodular, so their relative modular character is one; $dy\,dr$ is its Haar measure. Rational coverage gives this stabilizer at every $F$-point up to conjugation. Since $2\kappa\ne0$, \eqref{m8:s:s45:eq4.5} and the criterion \eqref{cf:orbital} imply
\begin{equation}\label{m8:s:s45:eq6.5}\Dist_{\eta_K}(S_4)=0,
\end{equation}
already under $U\times U$. The conjugation identity $(U\times U)_x=k_0^{-1}(U\times U)_{k_0x}k_0$ preserves the restricted character, so the vanishing transports through every initial normalization. The change $y\mapsto\kappa y$ multiplies Haar by $|\kappa|$
and preserves the nontriviality of $\psi(2\kappa y)$.

To compute the modulus of the $K$-stabilizer, use the diagonal
projections in Appendix~\S\ref{m8:s:app:B}: a split torus,
scalar matrices times a one-parameter unipotent, or scalars.
Each is unimodular. The two-dimensional kernel is central: the intersection unipotents with $X=z=0$ commute, and
\begin{equation}\label{m8:s:s45:eq6.2}
 \frac{\det B_0}{\alpha\delta}=1,\qquad
 \begin{pmatrix}\alpha&0\\\gamma&\delta\end{pmatrix}E_{22}
 \begin{pmatrix}\alpha&\phi\\0&\delta\end{pmatrix}^{-1}=E_{22},
 \quad B_0=\begin{pmatrix}\alpha&\phi\\0&\delta\end{pmatrix}.
\end{equation}
Conjugation fixes both kernel coordinates, so $\Delta_{K_x}=1$
and $\Delta_K|_{K_x}\Delta_{K_x}^{-1}=1$, in agreement with the
two-unipotent calculation. The diagonal stabilizer preserves the
same incompatibility.

\paragraph{\texorpdfstring{Transposition and \(S_5\)}{Transposition and S5}}\label{m8:s:s45:sub7}
\begin{proposition}\label{m8:s:s45:transpose}
Transposition carries the rational normal forms and stabilizers on \(S_4\) to those on \(S_5\). It follows that \(D'_{\eta_K}(S_5)=0\) and \(D'_{\widetilde\eta_K}(S_4\sqcup S_5)=0\).
\end{proposition}
\begin{proof}\leavevmode

\subparagraph{\texorpdfstring{Normal forms and stabilizers on \(S_5\)}{Normal forms and stabilizers on S5}}\label{m8:s:s45:sub7.1}

Transpose gives the covering family
\(x_5=x_{A,B,C,E;t,s}^t =L(C^t,E^t,s)\,w^t\,L(A^t,B^t,t).\)
This covers every \(F\)-point of \(S_5\): transpose it to \(S_4\), apply the rational even normalization already proved, and transpose the equation back. The residual equivalence is the transpose of \eqref{m8:s:s45:eq3.8}.

The complete stabilizer is
\(K_{x_5}=\theta(K_{x_4}),\) \(\theta(u,v,h)=(v,u,h^{-t}).\)
In particular, its two-unipotent stabilizer is \((v_{y,r},u_{y,r})\), with the same exponent \eqref{m8:s:s45:eq4.5}. All stabilizer moduli remain one. Hence
\(D'_{\eta_K}(S_5)=0.\)

\subparagraph{The two transpose components}\label{m8:s:s45:sub7.2}

The even \(K\)-action preserves all entries of the rank matrix. Since the rank matrices of \(S_4\) and \(S_5\) are different,
\(x^t\notin Kx\) \(\text{for every }x\in S_4\cup S_5.\)
The vanishing here is therefore supplied by the nontrivial unipotent stabilizer character.

More precisely, because the two components are clopen in their union, odd extension by transpose identifies
\(D'_{\eta_K}(S_4) \simeq D'_{\widetilde\eta_K}(S_4\sqcup S_5),\) \(T\longmapsto T-\tau T,\)
where the two summands are extended by zero from the two clopen components. This space is zero by \eqref{m8:s:s45:eq6.5}.

\end{proof}

\subsubsection[The two remaining strata with c = 1]{The two remaining strata with \texorpdfstring{\(c=1\)}{c = 1}}\label{m8:s:sec:s6-s7}
The two remaining rank-one incidence strata have different
character-compatible descendants: diagonal-torus conjugation on
$S_6$, and two-sided Whittaker translation on $S_7$. Their rational
coverage, complete residual relations and fiber densities are
proved in \S\ref{m8:supp:s6s7}; the intersection and transpose
calculations are in \S\ref{m8:supp:s6s7-details}.

\begin{proposition}\label{m8:s:v:vanishing}
For the strata \(S_6\) and \(S_7\), every \((K,\eta_K)\)-equivariant distribution invariant under transpose up to the sign \(-1\) vanishes:
\(D'_{\widetilde\eta_K}(S_6)=D'_{\widetilde\eta_K}(S_7)=0.\)
\end{proposition}

Write $L(A,B,t)=m(A,B,t)$. The representatives are
\begin{equation}\label{m8:s:v:e2-7}
\begin{aligned}
w_6&=E_{11}+E_{26}+E_{37}+E_{44}-E_{55}+E_{62}+E_{73}-E_{88},\\
w_7&=E_{11}+E_{27}+E_{33}+E_{45}+E_{54}+E_{66}+E_{72}+E_{88}.
\end{aligned}
\end{equation}
They are symmetric involutions in the prescribed connected
similitude group, with multipliers $-1$ and $1$.
Lemma~\ref{m8:s:v:rational-coverage} gives
$S_i(F)=P(F)w_iP(F)^t$ with rational pivot sections.
For a normal form $lw_ir$, put
\[
 l=L(A,B,t),\quad r=L(C,E,s),\quad
 M=B^{-1}A,\quad N=CE^{-1},\quad
 \kappa=\det B/t,\quad\ell=\det E/s.
\]
All four matrices are invertible and $\kappa\ell\ne0$.

\paragraph{Compatible support}
\begin{lemma}\label{m8:s:v:compatible}
The $U\times U$ stabilizers are a three-dimensional Heisenberg
group on $S_6$ and a four-dimensional additive group on $S_7$.
In the coordinates of \eqref{m8:s:v:e5-4}, their characters are
\begin{equation}\label{m8:s:v:e5-11}
 \eta_K(u,v,1)
 =\psi\bigl((M_{22}-2\ell)p+(2\kappa-N_{22})y\bigr)
 \quad (S_6),
\end{equation}
\begin{equation}\label{m8:s:v:e5-15}
 \eta_K(u,v,1)
 =\psi\bigl(M_{12}p+(M_{22}+N_{22})q-N_{21}z\bigr)
 \quad (S_7).
\end{equation}
Their zero loci define smooth closed subvarieties
$\mathcal C_6\subset S_6$ and $\mathcal C_7\subset S_7$, of
dimensions $24$ and $22$. On the normal family they are given by
\begin{equation}\label{m8:s:v:e5-12}
 M_{22}=2\ell,\qquad N_{22}=2\kappa,
\end{equation}
\begin{equation}\label{m8:s:v:e5-16}
 M_{12}=0,\qquad N_{21}=0,\qquad N_{22}=-M_{22},
\end{equation}
respectively. Every $(K,\eta_K)$-equivariant distribution on
$S_i$ is supported on $\mathcal C_i$.
\end{lemma}
\begin{proof}
The complete intersection equations \eqref{m8:s:v:e5-3}--\eqref{m8:s:v:e5-4}
and their finite exponentials are proved in
Appendix~\S\ref{m8:supp:s6s7-details}. Taking both $X$ traces
and scalar $2y$ terms gives the displayed characters.
Their coefficients vanish independently, including in residue
characteristic two, and rational normal-family sections give
smoothness and the stated dimensions. On the complement the
unipotent stabilizer is unimodular and its character is nontrivial.
The orbital criterion \eqref{cf:orbital} gives even vanishing there;
Lemma~\ref{cf:support} then gives the closed-support reduction.
\end{proof}

\begin{lemma}\label{m8:s:v:full-stabilizers}
On $\mathcal C_6$ the character on the complete $K$-stabilizer is
trivial. On $\mathcal C_7$, write
$EB=\left(\begin{smallmatrix}a&b\\c&d\end{smallmatrix}\right)$.
An additional unipotent direction occurs precisely when $d=0$;
its character is trivial precisely when $\kappa b=\ell c$.
All surviving stabilizers have relative modular character one.
\end{lemma}
\begin{proof}
The general stabilizer projections and characters are calculated
in Appendix~\S\ref{m8:s:app:B} and on
pp.~\pageref{m8:s:v:s6}--\pageref{m8:s:v:s6-3}.
For the canonical stabilizers, including the $d=0$ extension,
\eqref{m8:s:v:e12-1}--\eqref{m8:s:v:e12-2} and the calculation
after \eqref{m8:s:v:e13-5} give modular function one.
\end{proof}

\paragraph{The two canonical descendants}
\begin{proposition}[Distinct descendants of rank-one incidence]
\label{m8:rank-one-descendants}
Put $\Sigma=F^\times\times\GL_2(F)$ and define
\begin{align}
 X_6(\rho,V)&=L(I,I,2)w_6L(\operatorname{diag}(\rho,1)V,V,2\det V),
                 \label{m8:s:v:e7-3}\\
 X_7(\rho,V)&=w_7L(\operatorname{diag}(\rho,-1)V,V,\det V).
                 \label{m8:s:v:e7-8}
\end{align}
They cover $\mathcal C_6,\mathcal C_7$ under $K(F)$. Their complete
residual relations, with $\rho$ fixed, are respectively
$V\mapsto hVh^{-1}$ for diagonal $h\in T_2$, with character one,
and $V\mapsto n(a)V\bar n(b)$, with character $\psi(2a+2b)$.
Here $n(a)=\left(\begin{smallmatrix}1&a\\0&1\end{smallmatrix}\right)$
and $\bar n(b)=\left(\begin{smallmatrix}1&0\\b&1\end{smallmatrix}\right)$.
For $\pi_i(k,s)=kX_i(s)$, the fibers are free torsors for
\[
 H_6=H_3\times T_2,\qquad
 H_7=(\mathbb G_a^4\rtimes_\beta\mathbb G_a^2)\times\mathbb G_m.
\]
Here $H_6,H_7$ denote fiber groups, $H_3$ is the Heisenberg group
in \eqref{m8:s:v:e5-7}, and $\beta$ is \eqref{m8:s:v:e8-4}.
They have dimensions five and seven and are unimodular. Dividing
$dk\,d^\times\rho\,d^\times V$ by their Haar densities gives a
smooth nowhere-zero $K$-invariant density $\omega_i$ on
$\mathcal C_i$ with trivial relative density character.
\end{proposition}
\begin{proof}
Proposition~\ref{m8:s:v:slices} proves the rational normalizations,
complete residual relations and exact characters. Its divisions
are by entries forced nonzero by compatibility or invertibility.
Proposition~\ref{m8:s:v:fiber-torsors} constructs the free fiber
actions and local trivializations, retaining the kernel and
ineffective scalar also at residual stabilizer jumps.
Since the residual torus centralizes $H_3$ and $\beta$ has
determinant one, both fiber groups are unimodular.
Proposition~\ref{m8:s:v:density} gives
$\pi_i^*\omega_i=(dk\,d^\times\rho\,d^\times V)/dh_i$.
The $K$ and fiber-group actions preserve the quotient density.
\end{proof}

\begin{lemma}\label{m8:s:v:transpose-density}
The compensator $a(V)=(1,1,V^{-t})$ has character one and satisfies
\[
 a(V)X_i(\rho,V)^t=X_i(\rho,V^t).
\]
Transpose lifts to
$K\times\Sigma$ by
\begin{equation}\label{m8:s:v:e10-1}
 \Phi_i(k,\rho,V)=\bigl(\theta(k)a(V)^{-1},\rho,V^t\bigr).
\end{equation}
It preserves the total and fiber densities, hence also $\omega_i$.
\end{lemma}
\begin{proof}
The lift follows from \eqref{m8:s:v:e7-14}. Its differential is
triangular in the $K$ and $\Sigma$ variables: the diagonal blocks
preserve $dk$ and $d^\times\rho\,d^\times V$.
The kernel maps \eqref{m8:s:v:e10-3}--\eqref{m8:s:v:e10-4},
residual parameter exchange and scalar inversion are computed on
pp.~\pageref{m8:s:v:s10}--\pageref{m8:s:v:s10-3}. Their absolute Jacobians are
one; the factors $1/2$ on $S_6$ cancel in the canonical coordinates.
The tangent density identity \eqref{m8:s:c:density-line} therefore
proves invariance of the quotient density over every field in scope.
\end{proof}

\begin{proposition}\label{m8:s:v:distribution-correspondence}
The slice maps give the exact ordinary-distribution correspondences
\begin{equation}\label{m8:s:v:e11-4}
 D'_{\widetilde\eta_K}(\mathcal C_6)
 \simeq D'(F^\times\times\GL_2)^{T_2,\ \text{transpose-odd}},
\end{equation}
\begin{equation}\label{m8:s:v:e11-5}
 D'_{\widetilde\eta_K}(\mathcal C_7)
 \simeq D'(F^\times\times\GL_2)^{F^2,\ \psi(2a+2b),\ \text{transpose-odd}}.
\end{equation}
The actions are diagonal conjugation and two-sided Whittaker
translation, respectively. In both cases $\rho$ is passive and
transpose acts by $(\rho,V)\mapsto(\rho,V^t)$ and preserves
the density and character.
\end{proposition}
\begin{proof}
Apply Lemma~\ref{m8:s:c:fiber} to the preceding free torsors and
quotient densities. The kernel and ineffective scalar characters
are one; the residual transporter characters are $1$ and
$\psi(2a+2b)$. The determinant-one overlap maps give descent on
all rational charts; the complete descent is given on
pp.~\pageref{m8:s:v:s11}--\pageref{m8:s:v:s11-4}.
Under \eqref{m8:s:v:e10-1},
$\eta_K(\theta(k)a(V)^{-1})=\eta_K(k)$.
Together with the preceding lemma, this identifies transpose on the
parameter distribution with matrix transpose and proves both isomorphisms.
\end{proof}

\paragraph{Intrinsic vanishing}
\begin{proposition}\label{m8:s:v:s6-vanishing}
One has $D'_{\widetilde\eta_K}(S_6)=0$.
\end{proposition}
\begin{proof}
Modulo its scalar center, $T_2$ acts on
$V=\left(\begin{smallmatrix}a&b\\c&d\end{smallmatrix}\right)$
by $(b,c)\mapsto(qb,q^{-1}c)$, $q\in F^\times$.
Transpose exchanges $b,c$ and fixes $a,d,\rho$.
The two punctured axes are exchanged even orbits whose closures
contain their diagonal limit inside $\GL_2$. Thus the argument
uses the entire opposite-weight plane, including its origin.
Lemma~\ref{m8:s:c:plane}, with determinant-fiber localization to
$\GL_2$, makes the right side of \eqref{m8:s:v:e11-4} zero.
The closed-support reduction on $(S_6,\mathcal C_6)$ proves
the assertion. The complete stabilizer and transporter formulas
are \eqref{m8:s:v:e12-1}--\eqref{m8:s:v:e12-2} and the
transpose transporters on p.~\pageref{m8:s:v:s12-2}.
\end{proof}

\begin{proposition}\label{m8:s:v:s7-vanishing}
One has $D'_{\widetilde\eta_K}(S_7)=0$.
\end{proposition}
\begin{proof}
Use the two-sided Whittaker action with character $\psi(2a+2b)$
in \eqref{m8:s:v:e11-5}. On its closed boundary $V_{22}=0$,
the additive stabilizer is compatible precisely when
$V_{12}=V_{21}$. By \eqref{m8:s:v:e7-7}, this is exactly
$\kappa(EB)_{12}=\ell(EB)_{21}$, the condition of
Lemma~\ref{m8:s:v:full-stabilizers}.
The incompatible boundary has even vanishing. On the compatible
boundary the matrix is symmetric; on $V_{22}\ne0$, transpose
has a character-one residual transporter.
Lemma~\ref{m8:s:c:whittaker}, applied to this closed boundary
and its open complement, gives zero in \eqref{m8:s:v:e11-5}.
Combining the density and full stabilizer calculations on
pp.~\pageref{m8:s:v:s13}--\pageref{m8:s:v:s13-3} with the
closed-support reduction on $(S_7,\mathcal C_7)$ finishes the proof
and hence proves Proposition~\ref{m8:s:v:vanishing}.
\end{proof}

\subsubsection*{The common Siegel chart}\label{m8:common-siegel}
We use the following chart both for intrinsic vanishing on $S_8$
and for the later spectral argument on $S_3$. Its residual variables
$g,A,B,C,D,P,R,h$ and scalar $\nu$ are local to this chart: $P$ is
an additive matrix parameter here, while $Q_4$ denotes the parabolic.
We first establish the two descents, with their characters and
density factors, and then apply them to the two strata.

\paragraph{Siegel geometry}\label{m8:s:s3:geometry}

\subparagraph{Forms, orientation, and the Levi subgroup}\label{m8:s:s3:forms}

Set
\[
\mathcal J=\begin{pmatrix}0&J\\J&0\end{pmatrix},
\qquad \mathcal J^t=-\mathcal J,\qquad \mathcal J^2=-I_4.
\]
With the first two and last two original blocks grouped together,
\begin{equation}
Q=\begin{pmatrix}0&\mathcal J\\-\mathcal J&0\end{pmatrix}.
\label{m8:s:s3:eq:1-1}
\end{equation}
Choose the upper Siegel parabolic
\(Q_4=\operatorname{Stab}_{G_8^\circ}(E_+),\) \(E_+=F^4\oplus0.\)
Its opposite \(Q_4^t\) stabilizes \(E_-=0\oplus F^4\). This upper orientation contains the required radical in $U$ and remains in the fixed maximal-isotropic component.

For \(g\in \mathrm{GL}_4(F)\), define
\(g^\dagger=\mathcal J^{-1}g^{-t}\mathcal J.\)
The common Levi is
\(M_4=\{m(g,\nu)=\operatorname{diag}(g,\nu g^\dagger): g\in \mathrm{GL}_4(F),\ \nu\in F^\times\}.\)
The map \((g,\nu)\mapsto m(g,\nu)\) is a group isomorphism
\(M_4\simeq \mathrm{GL}_4(F)\times F^\times,\) \(\mu(m(g,\nu))=\nu.\)
Indeed,
\begin{equation}
g^t\mathcal J(\nu g^\dagger)=\nu\mathcal J,
\qquad \det m(g,\nu)=\det(g)\nu^4\det(g)^{-1}=\nu^4.
\label{m8:s:s3:eq:1-3}
\end{equation}

\subparagraph{The six-dimensional radical}\label{m8:s:s3:radical}

Let
\(\mathcal W=\{W\in \Mat_4(F):W^t\mathcal J=\mathcal J W\}.\)
Every such matrix has the unique form
\begin{equation}
W(Z,y,z)=\begin{pmatrix}Z&zI_2\\yI_2&Z^\#\end{pmatrix},
\qquad Z\in M_2(F),\quad y,z\in F.
\label{m8:s:s3:eq:1-4}
\end{equation}
Equivalently, \(W\mathcal J^{-1}\) is alternating. Hence
\begin{equation}
N_4=\left\{n_S(W)=\begin{pmatrix}I_4&W\\0&I_4\end{pmatrix}:W\in\mathcal W\right\},
\qquad \dim N_4=6.
\label{m8:s:s3:eq:1-5}
\end{equation}
This is an abelian vector group; its Haar measure will be \(dZ\,dy\,dz\). Direct multiplication shows that its form-preservation equation is precisely \(W^t\mathcal J=\mathcal J W\).

\subparagraph{The big cell and its explicit inverse}\label{m8:s:s3:bigcell}

Write an element of \(G_8^\circ\) in four-by-four blocks as
\[
x=\begin{pmatrix}E&F_0\\G_0&H\end{pmatrix}.
\]
Then
\(\Omega_4=\{\det H\ne0\}=N_4M_4N_4^t,\)
and multiplication is an algebraic isomorphism
\[
\mathcal W\times \mathrm{GL}_4\times\mathbb G_m\times\mathcal W
\longrightarrow\Omega_4.
\]
Its inverse is
\(W_L=F_0H^{-1},\) \(W_R=(H^{-1}G_0)^t,\) \(g=E-F_0H^{-1}G_0,\) \(\nu=\mu(x).\)
The factorization is
\begin{equation}
x=n_S(W_L)m(g,\nu)n_S(W_R)^t.
\label{m8:s:s3:eq:1-8}
\end{equation}

To verify the inverse, begin with the isotropy equation for the
final four columns,
\(F_0^t\mathcal JH-H^t\mathcal JF_0=0,\)
so \(W_L^t\mathcal J=\mathcal JW_L\). Since \(Q^2=I_8\), transpose preserves the same similitude group. Applying the same calculation to \(x^t\) proves that \(W_R\in\mathcal W\). Direct multiplication gives
\(n_S(-W_L)x n_S(-W_R)^t =\operatorname{diag}(g,H).\)
Invertibility of \(x,H\) implies invertibility of \(g\), and the remaining form equation gives
\(H=\nu\mathcal J^{-1}g^{-t}\mathcal J=\nu g^\dagger.\)
Conversely, \eqref{m8:s:s3:eq:1-3} proves that every product \eqref{m8:s:s3:eq:1-8} belongs to the connected \(\mathrm{GSO}_8\). Its lower-right four-by-four block is invertible, and the inverse
formulas are unique and regular for $\det H\ne0$. Thus the
radical fibers are free.

\subparagraph{The two rank strata in the slice}\label{m8:s:s3:strata}

Write
\begin{equation}
g=\begin{pmatrix}A&B\\C&D\end{pmatrix}\in \mathrm{GL}_4(F),
\label{m8:s:s3:eq:1-9}
\end{equation}
with two-by-two blocks. The terminal block of \(x\) is the terminal block of \(H=\nu g^\dagger\). Thus
\(x_{44}=\nu J^{-1}(g^{-t})_{11}J,\) \(\operatorname{rank}(x_{44})=\operatorname{rank}(D).\)
The rank equality follows from the kernel identification in \S\ref{m8:s:c:inverse-corner}.

Consequently, with
\[
\mathcal X_1=\left\{\begin{pmatrix}A&B\\C&D\end{pmatrix}\in \mathrm{GL}_4:
\operatorname{rank}D=1\right\},\qquad
\mathcal X_0=\left\{\begin{pmatrix}A&B\\C&0\end{pmatrix}\in \mathrm{GL}_4\right\},
\]
we have the exact identifications
\begin{equation}
S_3\simeq N_4\times(\mathcal X_1\times F^\times)\times N_4,
\qquad
S_8\simeq N_4\times(\mathcal X_0\times F^\times)\times N_4.
\label{m8:s:s3:eq:1-11}
\end{equation}

\paragraph{Character decomposition and distribution descent}\label{m8:s:s3:descent}

\subparagraph{\texorpdfstring{Decomposition of $U$}{Decomposition of U}}\label{m8:s:s3:decomposition}

Put
\[
n(P)=\begin{pmatrix}I_2&P\\0&I_2\end{pmatrix},
\qquad \ell(P)=m(n(P),1).
\]
Explicitly,
\[
\ell(P)=
\begin{pmatrix}
I_2&P&0&0\\
0&I_2&0&0\\
0&0&I_2&-P^\#\\
0&0&0&I_2
\end{pmatrix}.
\]
From the fixed coordinates of \(U\), multiplication gives
\begin{equation}
u(X,y,Z,z)=
n_S\left(Z,y,z+\tfrac12\operatorname{tr}(XZ^\#)\right)\ell(X).
\label{m8:s:s3:eq:2-1}
\end{equation}
Thus \(U=N_4\rtimes\ell(M_2)\). For the scalar shift in \eqref{m8:s:s3:eq:2-1} we use
\(XZ^\#+ZX^\#=\operatorname{tr}(XZ^\#)I_2.\)
Restricting the character gives
\(\xi(n_S(Z,y,z))=\psi(2y),\) \(\xi(\ell(X))=\psi(\operatorname{tr}X).\)
Thus the integrated characters are $\psi(2y_L+2y_R)$, and the two trace characters survive.

The normalizer formulas are
\begin{equation}
\ell(P)n_S(Z,y,z)\ell(P)^{-1}
=n_S\bigl(Z+yP,\ y,\ z+\operatorname{tr}(PZ^\#)+y\det P\bigr),
\label{m8:s:s3:eq:2-3}
\end{equation}
and
\begin{equation}
D(h)n_S(Z,y,z)D(h)^{-1}
=n_S(hZh^{-1},y,z).
\label{m8:s:s3:eq:2-4}
\end{equation}
Both preserve \(\psi(2y)\) and have Jacobian one on the six coordinates.

\subparagraph{Averaging and reconstruction}\label{m8:s:s3:averaging}

Use convention \eqref{cf:action}.

In the global chart \eqref{m8:s:s3:eq:1-8}, define
\begin{equation}
\mathcal A f(g,\nu)=
\int_{\mathcal W^2}
f\bigl(n_S(W_L)m(g,\nu)n_S(W_R)^t\bigr)
\psi(2y_L+2y_R)\,dW_L\,dW_R.
\label{m8:s:s3:eq:2-6}
\end{equation}
The same formula applies intrinsically to either subset in \eqref{m8:s:s3:eq:1-11}. A change of variables proves
\(\mathcal A(f\circ(n_S(V_L),n_S(V_R))^{-1}) =\psi(2v_{L,y}+2v_{R,y})\mathcal A f.\)
Both signs are positive. The free fibers have trivial stabilizer. Lemma~\ref{m8:s:c:product}, with a compactly supported fiber test of weighted integral one, gives the correspondence
\begin{equation}\label{m8:s:s3:eq:2-7}
 S\longmapsto[f\longmapsto S(\mathcal A f)].
\end{equation}

\subparagraph{The residual group and character}\label{m8:s:s3:residual}

It is convenient to transpose the right surviving parameter once and for all. Set
\[
\overline n(R)=\begin{pmatrix}I_2&0\\R&I_2\end{pmatrix}=n(R^t)^t,
\qquad d(h)=\operatorname{diag}(h,h).
\]
The residual group is
\(\mathscr L=(M_2(F)\oplus M_2(F))\rtimes \mathrm{GL}_2(F),\)
where \(h\) acts on both additive summands by conjugation. Its action is
\begin{equation}
(P,R,h)\cdot(g,\nu)
=\bigl(n(P)d(h)g d(h)^{-1}\overline n(R),\nu\bigr),
\qquad
\chi(P,R,h)=\psi(\operatorname{tr}P+\operatorname{tr}R).
\label{m8:s:s3:eq:2-9}
\end{equation}
The corresponding element of the group \(K\) is \((\ell(P),\ell(R^t),h)\). The group law in these coordinates is
\((P,R,h)(P',R',h') =(P+hP'h^{-1},\ R+hR'h^{-1},\ hh').\)
The right parameter has thus been reparametrized by transpose, $R=B_{\mathrm{right}}^t$.

Transpose gives
\(\tau(g,\nu)=(g^t,\nu),\) \(\theta(P,R,h)=(R^t,P^t,h^{-t}).\)
Here \((g^\dagger)^t=(g^t)^\dagger\). Extend \(\chi\) by \(\widetilde\chi(\tau)=-1\).

\subparagraph{Modular characters and density lines}\label{m8:s:s3:densities}

Conjugation by \(m(g,\nu)\) acts on the alternating matrix \(W\mathcal J^{-1}\) as
\[W\mathcal J^{-1}\longmapsto \nu^{-1}g(W\mathcal J^{-1})g^t.\]
Therefore the parabolic modulus is
\(\delta_{Q_4}(m(g,\nu)) =|\det g|^3|\nu|^{-6}.\)
The opposite parabolic has inverse modulus. With the group modular convention
\(\int f(xb)\,dx=\Delta(b)^{-1}\int f(x)\,dx\), the group modular functions of the two parabolics are the inverses of their respective parabolic moduli.

For the groups used in the descent,
\(\Delta_{N_4}=\Delta_{N_4^2}=\Delta_{M_4} =\Delta_{\mathrm{GL}_2}=\Delta_{\mathscr L}=\Delta_K=1.\)
The additive groups are unimodular; conjugation on \(M_2\) has determinant one; and \eqref{m8:s:s3:eq:2-3}--\eqref{m8:s:s3:eq:2-4} account for the action on \(N_4\). Since both fiber stabilizers are trivial, all relative modular
characters in \eqref{m8:s:s3:eq:2-7} are one.

The ambient big-cell Haar density has a nonconstant factor:
\begin{equation}
dx=|\det g|^{-3}|\nu|^6\,dW_L\,dg\,d^\times\nu\,dW_R,
\label{m8:s:s3:eq:2-13}
\end{equation}
up to the one fixed overall normalization. At \((1,m,1)\), left trivialization of the differential applies \(\operatorname{Ad}(m^{-1})\) to the first radical, giving \(\delta_{Q_4}(m)^{-1}\). This gives the stated inverse-modulus factor.

Both $\det g$ and $\nu$ are residual and transpose invariants, so the nonconstant factor \eqref{m8:s:s3:eq:2-13} changes only invariant normalization. Intrinsic distributions use the product chart \eqref{m8:s:s3:eq:1-11}; its equivariance is untwisted.

Transpose exchanges \(W_L,W_R\), preserves their identical Haar measures, and exchanges the two character factors in \eqref{m8:s:s3:eq:2-6}. Consequently,
\(D'_{\widetilde\eta_K}(S_3) \simeq D'_{\widetilde\chi}(\mathcal X_1\times F^\times),\) \(D'_{\widetilde\eta_K}(S_8) \simeq D'_{\widetilde\chi}(\mathcal X_0\times F^\times).\)
The entire extended group fixes $\nu\in F^\times$. Scalar $h$ is ineffective with character one and remains in the acting group.

\paragraph{The residual block action and stabilizers}\label{m8:s:s3:block-action}

For \(g\) as in \eqref{m8:s:s3:eq:1-9}, the additive subgroup \(E=M_2\oplus M_2\) acts by
\begin{equation}
(A,B,C,D)\longmapsto
(A+PC+BR+PDR,\ B+PD,\ C+DR,\ D).
\label{m8:s:s3:eq:3-1}
\end{equation}
The diagonal \(\mathrm{GL}_2\) conjugates all four blocks simultaneously. Thus the additive stabilizer is
\(E_g=\{(P,R):PD=0,\ DR=0,\ PC+BR=0\}.\)
The first two equations imply $PDR=0$.

For a residual stabilizer, first require \(hDh^{-1}=D\). Writing \(A^h=hAh^{-1}\), and similarly for the other blocks, the remaining equations are
\begin{equation}
B^h+PD=B,\quad C^h+DR=C,\quad
A^h+PC^h+BR=A.
\label{m8:s:s3:eq:3-3}
\end{equation}

\subsubsection{\texorpdfstring{Intrinsic odd vanishing on $S_8$}{Intrinsic odd vanishing on S8}}\label{m8:s:s3:S8}
\begin{proposition}[Intrinsic odd vanishing on the zero-block stratum]\label{m8:s:s3:prop:S8}
One has $D'_{\widetilde\eta_K}(S_8)=0$ for ordinary distributions.
\end{proposition}
\par\smallskip\noindent\textit{Proof.}\par

\paragraph{Block equations and normal form}\label{m8:s:s3:S8-normal}

Here \(D=0\), so
\begin{equation}
g=\begin{pmatrix}A&B\\C&0\end{pmatrix}\in \mathrm{GL}_4(F)
\quad\Longleftrightarrow\quad B,C\in \mathrm{GL}_2(F).
\label{m8:s:s3:eq:8-1}
\end{equation}
Indeed, the last two columns and the last two rows force invertibility of \(B,C\), and conversely the matrix is invertible for every \(A\) when they are invertible.

Taking \(P=-AC^{-1}\), \(R=0\) gives
\begin{equation}
g_0(B,C)=\begin{pmatrix}0&B\\C&0\end{pmatrix}.
\label{m8:s:s3:eq:8-2}
\end{equation}
This normalizing element has character
\(\psi(-\operatorname{tr}(AC^{-1}))\).

The additive stabilizer at any such point is
\(P\in M_2(F),\) \(R=-B^{-1}PC.\)
Its character is
\begin{equation}
\psi\left(\operatorname{tr}\bigl((I_2-CB^{-1})P\bigr)\right).
\label{m8:s:s3:eq:8-4}
\end{equation}
The nondegeneracy of the trace pairing therefore forces
\begin{equation}
C=B.
\label{m8:s:s3:eq:8-5}
\end{equation}
Since this condition is closed and stable, the unipotent stabilizer
criterion eliminates its complement.

\paragraph{Distribution descent and density}\label{m8:s:s3:S8-density}

On the compatible set write \(B=C=K_0\in \mathrm{GL}_2(F)\). The action is
\(A\longmapsto A+PK_0+K_0R,\) \(K_0\longmapsto hK_0h^{-1},\) \(A\longmapsto hAh^{-1}\)
for the additive and diagonal parts, respectively. The additive action is transitive on \(A\), and
\[
\operatorname{tr}\bigl(K_0^{-1}(PK_0+K_0R)\bigr)
=\operatorname{tr}P+\operatorname{tr}R.
\]
Thus the exact affine-fiber distribution is
\begin{equation}
|\det K_0|^{-2}\,
\psi\bigl(\operatorname{tr}(K_0^{-1}A)\bigr)\,dA.
\label{m8:s:s3:eq:8-7}
\end{equation}
To check the density, use \((P,A)\) as coordinates on the two additive parameters. Then
\(R=K_0^{-1}A-K_0^{-1}PK_0,\) \(dP\,dR=|\det K_0|^{-2}dP\,dA.\)
The stabilizer is parametrized by \(dP\). Formula \eqref{m8:s:s3:eq:8-7} is the resulting quotient density; it is invariant under simultaneous conjugation and under \((A,K_0)\mapsto(A^t,K_0^t)\).

We are therefore left with distributions on
\(\GL_2(F)\times F^\times\), with coordinates \((K_0,\nu)\),
invariant under conjugation of \(K_0\) and odd under
\(K_0\mapsto K_0^t\), with \(\nu\) passive. 

At the normal form \(g_0(K_0,K_0)\), the stabilizer is
\(h\in Z_{\mathrm{GL}_2}(K_0),\) \(P\in M_2(F),\) \(R=-K_0^{-1}PK_0.\)
Its character is trivial. It is unimodular: the centralizer is either \(\mathrm{GL}_2\) or the units in a commutative two-dimensional algebra, and its conjugation action on \(M_2\) has determinant one.

\paragraph{Transpose transporters}\label{m8:s:s3:S8-transpose}

Use the cyclic-vector matrix $h$ in \eqref{m8:s:s3:eq:8-11}, satisfying $hK_0h^{-1}=K_0^t$ (and $h=I_2$ for scalar $K_0$). Thus
\begin{equation}
(0,0,h)g_0(K_0,K_0)=g_0(K_0,K_0)^t,
\qquad\chi(0,0,h)=1.
\label{m8:s:s3:eq:8-12}
\end{equation}
For a non-normalized point use \eqref{m8:s:c:transport-back}.

Lemma~\ref{m8:s:c:transporter} and the passive-parameter Lemma~\ref{m8:s:c:passive} prove vanishing of this residual odd space. 

Combining \eqref{m8:s:s3:eq:8-4}--\eqref{m8:s:s3:eq:8-12} with the first descent proves
\(D'_{\widetilde\eta_K}(S_8)=0.\)

\par\hfill$\square$\par\smallskip

\subsubsection{Terminal relative positions and one-matrix descent}
\label{m8:s:sec:terminal}
Put $W=\left(\begin{smallmatrix}0&1\\1&0\end{smallmatrix}\right)$,
$\Sigma=\operatorname{diag}(1,-1)$, and
$\ell(A,B,t)=m(A,B,t)$; write $L=M$ for this Levi group.
Because the five representatives normalize the Levi, their
compatibility equations admit the uniform treatment below. Four
strata are incompatible with the unipotent character; the fifth
descends to conjugation on $\GL_2$. We prove the rational
factorizations and complete stabilizers in \S\ref{m8:supp:terminal}.

\begin{theorem}[Terminal compatibility and descent]\label{m8:s:t:main}
For $i=9,13,14,17,18$, put
$\mathcal N_i=U\cap w_iU^tw_i^{-1}$ and
$\mathcal V_i=U\cap w_iUw_i^{-1}$.
Every point of $S_i(F)$ has a unique factorization
$v_0\ell(A,B,t)w_iv^t$ with
$v_0\in\mathcal V_i(F)$, $v\in U(F)$, and $A,B\in\GL_2(F)$,
$t\in F^\times$. For $R=B^{-1}A$ and $\kappa=\det B/t$, the
complete unipotent stabilizer and its character are determined by
the following table:
\begin{equation}\label{m8:s:t:eq5.1}
\begin{array}{c|c|l}
 i&\mathcal N_i&\Lambda_i(X,y,Z,z)\\\hline
9&u(0,y,Z,0)&2\kappa y-\operatorname{tr}(WZ)\\
13&u(X,0,0,z)&\operatorname{tr}(RX)-2z\\
14&u(0,y,Z,z)&2\kappa y-2z\\
17&u(X,0,Z,z)&\operatorname{tr}(RX)+\operatorname{tr}(WZ)\\
18&U&\operatorname{tr}((R+I)X)+2(\kappa-1)y.
\end{array}
\end{equation}
Its embedding is
$(\ell n\ell^{-1},(w_i^{-1}n^{-1}w_i)^t)$, and its character
is $\psi(\Lambda_i(n))$.
Consequently
\[
 D'_{\eta_K}(S_9)=D'_{\eta_K}(S_{13})
 =D'_{\eta_K}(S_{14})=D'_{\eta_K}(S_{17})=0.
\]
These vanishings already hold under $U\times U$. The only compatible
terminal locus is
\begin{equation}\label{m8:s:t:eq5.2}
 \mathcal C_{18}=\{\ell(-B,B,\det B)w_{18}v^t:B\in\GL_2(F),v\in U\}.
\end{equation}
It is closed, smooth and fourteen-dimensional. There is an exact
ordinary-distribution descent
\begin{equation}\label{m8:s:t:eq9.8}
 D'_{\eta_K}(\mathcal C_{18})
       \simeq D'(\GL_2(F))^{\mathrm{conj}},
\end{equation}
under which transpose is matrix transpose. Hence
\[
 D'_{\widetilde\eta_K}(S_{18})=0,
 \qquad D'_{\eta_K}(S_{18})\ne0.
\]
The locally closed family $\mathcal B=S_9\sqcup S_{13}\sqcup S_{14}$
has even vanishing, and the closed family
$\mathcal D=S_{17}\sqcup S_{18}$ has odd vanishing.
\end{theorem}

\begin{proof}
\emph{Factorization and compatibility.}
Proposition~\ref{m8:s:t:factorization} gives the global rational
chart by terminal-rank pivots and orthogonal block elimination.
The forbidden-block equations in \eqref{m8:s:t:eq4.2} determine
$\mathcal N_i$ and prove sufficiency as group equations. Their
dimensions are $5,5,6,9,10$; $\mathcal N_{17}$ is Heisenberg
with pairing \eqref{m8:s:t:eq4.8}. The expanded conjugates are
in \texttt{anc/model8\_calculations.pdf}, M8-C03,
``The full stabilizer equations''.
Levi conjugation \eqref{m8:s:t:eq1.5} and the full inverse
\eqref{m8:s:t:eq1.4}, with its $yX$ correction, give the right
factor $(w_i^{-1}n^{-1}w_i)^t$. Adding its $X,y$ character to
the left one gives \eqref{m8:s:t:eq5.1}.

For $i=9,13,14,17$, respectively, the one-parameter choices
\[
 u(0,r,0,0),\quad u(0,0,0,r),\quad
 u(0,0,0,r),\quad u(0,0,rE_{12},0)
\]
have character exponents $2\kappa r,-2r,-2r,r$. Each is nonzero
as an $F$-linear functional. On $S_{18}$, independence of $X$ and
$y$ gives exactly $R=-I$, $\kappa=1$. Every nonzero $F$-linear functional is surjective onto $F$,
so these are the complete compatibility conditions. In particular the argument applies in residue
characteristic two.

Since the acting unipotent group and its stabilizers are unimodular,
the rational orbit criterion gives ordinary even vanishing
on the four incompatible strata and on
$S_{18}\setminus\mathcal C_{18}$. This uses the criterion with
$\eta^{-1}$ on the stabilizer; a nontrivial character and its inverse
are both distinct from the trivial relative modulus.

\emph{Full stabilizers.}
Proposition~\ref{m8:s:t:full-stabilizers} gives
$K_x=\mathcal N_i\rtimes Z_{\GL_2}(AW^{\epsilon_A},BW^{\epsilon_B})$,
with $(\epsilon_A,\epsilon_B)=(0,1),(1,0),(1,1),(1,0),(1,1)$.
The residual relation is simultaneous rational conjugation with
$t$ fixed. The centralizer list includes nonsplit quadratic tori
and repeated-root unipotents; \eqref{m8:s:t:eq6.4} has determinant
one on each present matrix coordinate and fixes the scalar
coordinates. Thus $\Delta_K=\Delta_{K_x}=1$ at every rational point.
The Levi modulus \eqref{m8:s:t:eq6.9} also equals one on
$A=-B$, $t=\det B$.

\emph{The compatible descendant and its density.}
Write $M=BW$ and set
\begin{equation}\label{m8:s:t:eq8.2}
 x(M)=\begin{pmatrix}
 0&0&0&-M\\0&0&M&0\\0&M&0&0\\-M&0&0&0
 \end{pmatrix}.
\end{equation}
It has similitude $-\det M$ and determinant $(\det M)^4$.
Proposition~\ref{m8:s:t:chart} proves the global product chart
\begin{equation}\label{m8:s:t:eq8.3}
 \GL_2\times U\longrightarrow\mathcal C_{18},
 \qquad(M,v)\longmapsto x(M)v^t,
\end{equation}
with inverse $M=x_{23}$, $v^t=x(M)^{-1}x$. For $C=MJ$, its
fiber anti-automorphism is
\begin{equation}\label{m8:s:t:eq8.7}
 \theta_M(u(X,y,Z,z))
 =u(C^{-1}XC,y,-C^{-1}(Z-yX)C,z).
\end{equation}
The chart actions \eqref{m8:s:t:eq8.10}--\eqref{m8:s:t:eq8.13}
satisfy $\xi\circ\theta_M=\xi$, $(\theta_M)_*dv=dv$ and
$\theta_{M^t}\theta_M=1$. Thus transpose and both unipotent actions
preserve the fiber density.

Define the averaging map
\begin{equation}\label{m8:s:t:eq9.1}
 \mathcal Af(M)=\int_U f(x(M)v^t)\xi(v)\,dv.
\end{equation}
It is a surjection onto $D(\GL_2)$: take a compactly supported
$\beta_0$ with $\int\beta_0\xi=1$ and lift
$\phi(M)$ to $\phi(M)\beta_0(v)$ in the product chart.
Every test function there is a finite sum of such products.
Right-$U$ equivariance forces each fiber functional to be a
multiple of $\xi(v)\,dv$, by Haar uniqueness after twisting.
Consequently every such distribution is uniquely
\begin{equation}\label{m8:s:t:eq9.3}
 T(f)=R(\mathcal Af),\qquad R\in D'(\GL_2).
\end{equation}
Conversely every $R$ defines a distribution by this formula.
Right Haar invariance and $\xi\circ\theta_M=\xi$ give the left
character $\xi(a)$ as well. The preceding density calculation gives
\[
 \mathcal A(f\circ h^{-1})(M)=\mathcal Af(h^{-1}Mh),\qquad
 \mathcal A(f\circ\tau)(M)=\mathcal Af(M^t).
\]
Surjectivity now proves \eqref{m8:s:t:eq9.8}, and, in both directions,
\begin{equation}\label{m8:s:t:eq9.9}
 D'_{\widetilde\eta_K}(\mathcal C_{18})
 \simeq D'(\GL_2)^{\mathrm{conj},\,\tau\text{-odd}}.
\end{equation}

\emph{Transpose on every rational characteristic-polynomial fiber.}
For a nonscalar $M=\left(\begin{smallmatrix}a&b\\c_0&d\end{smallmatrix}\right)$,
choose a cyclic vector $e_1$ when $c_0\ne0$, $e_2$ when
$c_0=0$, $b\ne0$, and $e_1+e_2$ otherwise. The cyclic determinants
are $c_0,-b,d-a$ in these cases. The bases $[v,Mv]$ and
$[w,M^tw]$ identify both matrices with
$\left(\begin{smallmatrix}0&-\det M\\1&\operatorname{tr}M\end{smallmatrix}\right)$.
Their ratio gives a rational transporter $hMh^{-1}=M^t$.
This construction also covers an irreducible quadratic polynomial
and a nonsemisimple repeated root. The diagonal element
$k_h=(1,1,h)$ transports $x(M)$ to its transpose with character one.
If $x=k_0x(M)$, its transporter is
\begin{equation}\label{m8:s:t:eq10.4}
 k_x=\iota(k_0)k_hk_0^{-1},\qquad \eta_K(k_x)=1.
\end{equation}

Localize distributions under
$M\mapsto(\operatorname{tr}M,\det M)$. With nonzero discriminant
each fiber is one rational conjugacy orbit. With zero discriminant
it is the closed scalar orbit and the nonscalar Jordan orbit,
separated by the ordinary support sequence. Conjugacy orbits have
unimodular stabilizers by the centralizer list above. Transpose
preserves each orbit, acts on its invariant measure line by a
positive scalar, and has square one; that scalar is one.
There is therefore no odd invariant distribution on either orbit.
Localization proves the right side of \eqref{m8:s:t:eq9.9} is zero.
Together with even vanishing of the open incompatible complement
and the closed-support sequence, this proves odd vanishing on
$S_{18}$.

The even space is nonzero: the distribution
\begin{equation}\label{m8:s:t:eq11.3}
 \mathscr T_0(f)=\int_U f(x(I_2)v^t)\xi(v)\,dv
\end{equation}
is well defined on the closed product-chart fiber $M=I_2$,
is even-equivariant by the action formulas, and is fixed by transpose.
A product test function whose $U$-factor is a compact open set
on which $\xi=1$ proves it is nonzero.

\emph{The two support families.}
In the open $\{r_{6,2}=2\}$ the two closed loci
$x_{34}=x_{44}=0$ and $x_{43}=x_{44}=0$ are respectively
$S_9\sqcup S_{14}$ and $S_{13}\sqcup S_{14}$.
Their intersection is $S_{14}$, so $\mathcal B$ is locally closed,
with $S_9,S_{13}$ clopen in its complement of $S_{14}$.
Also $\mathcal D=\{r_{6,2}=0\}$ is closed and $S_{18}$ is closed
in it. The two left-exact support sequences, with even equivariance
on $\mathcal B$ and odd equivariance on $\mathcal D$, prove the
family assertions. The support arguments use these two families separately; their union
fails to be locally compact by
\eqref{m8:s:t:eq12.7}--\eqref{m8:s:t:eq12.8}.
\end{proof}

\paragraph{Two open support extensions}\label{m8:s:t:sub13}
\begin{lemma}\label{m8:s:t:open-extension}
If a $(K,\eta_K)$-equivariant distribution $T$ vanishes on
\begin{equation}\label{m8:s:t:eq13.1}
\mathcal O_{\rm old}=S_1\sqcup S_2\sqcup S_3\sqcup S_4\sqcup S_5\sqcup S_8,
\end{equation}
then it vanishes on
$\mathcal O_{\rm new}=\{r_{2,4}=2\}\cup\{r_{4,2}=2\}$.
\end{lemma}
\begin{proof}
By the rank table, the two maximal-rank opens are
\begin{equation}\label{m8:s:t:eq13.2}
\begin{aligned}
\mathcal O_L=\{r_{2,4}=2\}
 &=S_1\sqcup S_2\sqcup S_3\sqcup S_4\sqcup S_8\sqcup S_9,\\
\mathcal O_R=\{r_{4,2}=2\}
 &=S_1\sqcup S_2\sqcup S_3\sqcup S_5\sqcup S_8\sqcup S_{13}.
\end{aligned}
\end{equation}
In $\mathcal O_L$, the equations $x_{34}=x_{44}=0$ cut out
$S_9$: the open rank condition makes $x_{43}$ invertible.
Its complement lies in $\mathcal O_{\rm old}$. Transpose gives
the same assertion for $S_{13}\subset\mathcal O_R$.
Apply the left-exact sequences
\begin{equation}\label{m8:s:t:eq13.3}
0\longrightarrow D'_{\eta_K}(S_i)
\longrightarrow D'_{\eta_K}(\mathcal O_\bullet)
\longrightarrow D'_{\eta_K}(\mathcal O_\bullet\setminus S_i),
\qquad (i,\bullet)=(9,L),(13,R).
\end{equation}
By Theorem~\ref{m8:s:t:main} the first space is zero, while $T$
has zero image in the last. Hence
\begin{equation}\label{m8:s:t:eq13.5}
T|_{\mathcal O_{\rm new}}=0,\qquad
\mathcal O_{\rm new}=\mathcal O_{\rm old}\sqcup S_9\sqcup S_{13}.
\end{equation}
\end{proof}
For the ambient NW distribution, the spectral argument supplies
the premise on $S_3$.

For reference, the rank table also describes the larger natural opens
\(\Omega_b=\{r_{6,2}=2\}\) and \(\Omega_a=\{r_{4,4}\ge2\}\) by
\begin{equation}\label{m8:s:t:eq13.6}
\Omega_b\setminus(S_{14}\cup\mathcal O_{\mathrm{new}})
 =S_6\sqcup S_{10},
\end{equation}
\begin{equation}\label{m8:s:t:eq13.7}
\Omega_a\setminus(S_{17}\cup\mathcal O_{\mathrm{new}})
 =S_6\sqcup S_7\sqcup S_{10}\sqcup S_{11}\sqcup S_{12}\sqcup S_{15}.
\end{equation}
The vanishing on \(S_{14}\) or \(S_{17}\) enters a support argument only
after its corresponding open complement has been cleared. The complete
global filtration will make these premises explicit.

\subsubsection{Middle relative positions and the two residual characters}
\label{m8:s:sec:middle}
The five middle cells have one-dimensional joint graded pieces.
Evaluating the intersection characters leaves two residual matrix
problems: one gives odd vanishing on $S_{10}$, while the other
produces the closed $S_{11}$ support family excluded in
\S\ref{m8:spectral}.
Rational flag coverage, intersections and nonlinear lifts are proved
in \S\ref{m8:supp:middle}.

Write $n_+(s)=\left(\begin{smallmatrix}1&s\\0&1\end{smallmatrix}\right)$
and $n_-(v)=\left(\begin{smallmatrix}1&0\\v&1\end{smallmatrix}\right)$.
The representatives $w_i\in\mathrm{SO}(Q)(F)$ are specified in
\S\ref{m8:s:m:representatives}.

\begin{theorem}[Middle-cell compatibility and residual equivalence]
\label{m8:middle-family}
For each of the five indicated labels, $S_i(F)=P(F)w_iP(F)^t$.
Let $I_i=P\cap w_iP^tw_i^{-1}$ and write a normal form as
$mw_in$, where $m=m(A,B,\mu)$ and $n=m(C,D,\nu)$.
The complete unipotent stabilizer is the exponential of
\begin{equation}\label{m8:middle-kernel}
 \{(\operatorname{Ad}(m)N,
     -\operatorname{Ad}(n^t)(w_i^{-1}Nw_i)^t):
       N\in\mathfrak u\cap\operatorname{Ad}(w_i)\mathfrak u^t\}.
\end{equation}
Its dimensions are $4,5,6,6,7$. On this group the character is defined by the
linear functional
\[
 \operatorname{tr}(RX)+2\kappa y+
 \operatorname{tr}(SX')+2\lambda y',
\]
where
\[
 R=B^{-1}A,\quad S=(CD^{-1})^t,
 \quad\kappa=\det B/\mu,\quad\lambda=\det D/\nu.
\]
The complete five-row evaluation is
\eqref{m8:s:m:eq:characters}. Its zero locus is empty for
$i=12,15,16$. For $i=10,11$ it is respectively
\[
 R_{12}=S_{12}=0,\ R_{11}=S_{11};\qquad
 R_{11}=S_{11}=0,\ R_{12}=S_{12},\ R_{21}=S_{21}.
\]
Every compatible point has a rational canonical representative
$F_\epsilon(E)$, respectively $G(E)$, defined below, with
$E\in\GL_2(F)$ and $\epsilon\in F^\times$ on $S_{10}$.
The complete residual equivalence is
\begin{equation}\label{m8:s:m:eq:residual}
 E'=n_+(s)En_-(v),
\end{equation}
with $\epsilon$ fixed on $S_{10}$ and exact characters
\begin{equation}\label{m8:s:m:eq:residual-character}
 \chi_{10,\epsilon}(s,v)=\psi(2(s+v)/\epsilon),
 \qquad \chi_{11}(s,v)=1.
\end{equation}
All statements concern the $F$-point groups, including their
components and all singular residual parameter loci.
\end{theorem}
\begin{proof}
The rank-nullity formula \eqref{m8:flag-ranks} gives the joint
multiplicities in \S\ref{m8:s:m:flags}. Lemma~\ref{m8:simultaneous-witt}
and Proposition~\ref{m8:s:m:prop:cells} give $F$-point coverage
and rational pivot sections, retaining the original similitude.

The intersections $I_i$ are split connected solvable groups. Their
common torus has five parameters; their unipotent kernels are
obtained by the finite polynomial logarithm from preservation of
both flags. In the two Levi projections their images are
$(l t,t_i r)$, with $t_i=w_i^{-1}tw_i$ and
\[
\begin{array}{c|c|c}
i&L_i&R_i\\\hline
10&N_-\times N_-&N_+\times N_+\\
11&N_+\times N_-&N_-\times N_+\\
12&N_+\times N_+&N_-\times N_+\\
15&N_+\times N_-&N_-\times N_-\\
16&N_+\times N_+&N_-\times N_-
\end{array}
\]
The root conjugations in Appendix, p.~\pageref{m8:s:m:roots} realize all these
projections on $F$-points. Their common kernel is
$U\cap w_iU^tw_i^{-1}$, giving \eqref{m8:middle-kernel} by
logarithm and exponential. The latter retains the correction
$Z\mapsto Z+yX/2$ and the scalar term
$z\mapsto z-y\det X/6$. Thus the character on the stabilizer,
which depends only on $X,y,X',y'$, is precisely the stated
linear functional.

We evaluate the five characters and determine their zero loci
below. On the two compatible loci, the Levi projection equations
give every parameter change. Using only entries forced nonzero,
the normalizations in Appendix, p.~\pageref{m8:s:m:normalization}
then reach the canonical families. Before the
similitude gauge their exact relation is
\[
 E'=n_+(s)E\begin{pmatrix}1&0\\-h&t\end{pmatrix},
 \quad \rho'=\rho/t,
\]
with character $\psi(2s/\epsilon-2h/\rho)$ on $S_{10}$ and
one on $S_{11}$. Fixing $\rho=\epsilon$, respectively $\rho=1$,
forces $t=1$. Setting $v=-h$ proves
\eqref{m8:s:m:eq:residual}--\eqref{m8:s:m:eq:residual-character}.
The nonlinear group lifts in Appendix, p.~\pageref{m8:s:m:lifts} prove
sufficiency with these exact signs. Necessity follows from the
intersection projections just used. The residual action on $S_{11}$ is exactly the two-unipotent
action in \eqref{m8:s:m:eq:residual} for every invertible $E$.
\end{proof}

\paragraph{Canonical representatives}
Set $W=\left(\begin{smallmatrix}0&1\\1&0\end{smallmatrix}\right)$,
$H=\operatorname{diag}(1,-1)$ (only in this paragraph, a matrix
rather than the model subgroup), $m_0=m(W,I,1)$ and
$V=m_0w_{11}m_0$. Then
\[
 F_\epsilon(E)=w_{10}m(\operatorname{diag}(1,\epsilon)E,E,
                        \epsilon\det E),\qquad
 G(E)=m_0w_{11}m(WE,E,\det E)=VD(E).
\]
Both $V$ and $V_\epsilon=F_\epsilon(I)$ are symmetric.
The similitudes are $\epsilon\det E$ and $\det E$.
For later use the $S_{10}$ residual lift is
\[
 k_{10}(s,v)=
 (u(0,s/\epsilon,0,s),u(0,v/\epsilon,0,v),n_-(-v)).
\]
It has character $\psi(2(s+v)/\epsilon)$. The $S_{11}$ lift
uses $Z$-coordinates, has character one, and is given with all
nonlinear terms in Appendix, p.~\pageref{m8:s:m:lifts}.

\paragraph{Stabilizer characters and three even-vanishing strata}\label{m8:s:m:characters}
\begin{proposition}[Even vanishing]\label{m8:s:m:prop:even}
The strata $S_{12},S_{15},S_{16}$ satisfy
\[
D'_{\xi\boxtimes\xi}(S_{12})=D'_{\xi\boxtimes\xi}(S_{15})=D'_{\xi\boxtimes\xi}(S_{16})=0.
\]

\end{proposition}
\par\smallskip\noindent\textit{Proof.}\par

\subparagraph{Character formulas}\label{m8:s:m:character-table}

Put
\[
R=B^{-1}A,\quad S=(CD^{-1})^t,\quad
\kappa=\frac{\det B}{\mu},\quad
\lambda=\frac{\det D}{\nu}.
\]
The left conjugation sends $X$ to \(AXB^{-1}\) and $y$ to \(\kappa y\). The right conjugation through \(n^t\) sends $X'$ to \(C^tX'D^{-t}\) and $y'$ to \(\lambda y'\). Taking traces therefore gives the exponent
\[
{\operatorname{tr}(RX)+2\kappa y+
       \operatorname{tr}(SX')+2\lambda y'.}
\]
Exponentiation preserves $X,y$, so these formulas give the group stabilizer character.

Substitution of the stabilizer tables gives
\begin{equation}
\begin{array}{c|l}
S_{10}&(R_{11}-S_{11})p+R_{12}r-S_{12}a\\
S_{11}&R_{11}p-S_{11}z+(R_{21}-S_{21})q+(R_{12}-S_{12})r\\
S_{12}&R_{21}q-S_{11}z-S_{21}a-S_{12}d+(2\kappa-S_{22})y\\
S_{15}&R_{11}p+R_{21}q+R_{12}r+(R_{22}+2\lambda)h-S_{21}a\\
S_{16}&(R_{21}-S_{21})q+2\kappa y+2\lambda c.
\end{array}
\label{m8:s:m:eq:characters}
\end{equation}
 Here \(S_{ij}\) denotes a matrix entry in \eqref{m8:s:m:eq:characters}.

\subparagraph{Compatible loci}\label{m8:s:m:compatibility}

Because the free stabilizer coordinates range over all $F$, a nonzero linear coefficient gives a nontrivial character for the fixed nontrivial additive $\psi$. Thus
\[
\mathcal C_{10}:\quad R_{12}=S_{12}=0,\quad R_{11}=S_{11},
\]
\[
\mathcal C_{11}:\quad R_{11}=S_{11}=0,\quad
R_{12}=S_{12},\quad R_{21}=S_{21}.
\]
On $S_{10}$, invertibility makes both matrices lower triangular with the same nonzero upper-left entry. On $S_{11}$, their common off-diagonal entries are both nonzero.

On $S_{12}$, compatibility would require \(S_{11}=S_{12}=S_{21}=0\), making $S$ singular. On $S_{15}$, it would require \(R_{11}=R_{12}=R_{21}=0\), making $R$ singular. On $S_{16}$, the free $y$-direction has coefficient \(2\kappa\ne0\). Hence none of these three strata has a compatible $F$-point.

\subparagraph{\texorpdfstring{Vanishing on $S_{12},S_{15},S_{16}$}{Vanishing on S12,S15,S16}}\label{m8:s:m:even-vanishing}

Apply the rational-point criterion \eqref{cf:orbital} to
$U\times U$ on each exact-rank stratum. These are algebraic locally
closed $F$-varieties, and the complete stabilizers just computed
are unipotent. Both acting group and stabilizers are therefore
unimodular. At every point the character, and hence its inverse
in the criterion, is nontrivial by \S\ref{m8:s:m:compatibility}.
The criterion applies to ordinary distributions on the entire
rational-point locus.
It follows that
\[
{
D'_{\xi\boxtimes\xi}(S_{12})=
D'_{\xi\boxtimes\xi}(S_{15})=
D'_{\xi\boxtimes\xi}(S_{16})=0.
}
\]
Restricting $K$-equivariance to $U\times U$ proves the claimed even-group vanishing for all three strata.

\par\hfill$\square$\par\smallskip

\paragraph{Stabilizers and quotient densities}\label{m8:s:m:K-stab}

Let
\[
E=\begin{pmatrix}p&q\\r&d\end{pmatrix}.
\]
The solved residual equations show that every stabilizer element has diagonal component
\[
a=c\,n_-(-v),\qquad c\in F^\times,
\]
and residual parameters satisfying
\[
n_+(s)En_-(v)=E.
\]
Expanding the matrix equation gives exactly
\begin{equation}
sd=0,\qquad vd=0,\qquad sr+vq+svd=0.
\label{m8:s:m:eq:residual-stabilizer}
\end{equation}
The fiber over any such diagonal component consists of one explicit lift from Appendix, p.~\pageref{m8:s:m:lifts} times the entire $U\times U$ kernel in Appendix, p.~\pageref{m8:s:m:stab}. This proves completeness of the following descriptions.

\subparagraph{\texorpdfstring{The locus $d\ne0$}{The locus d != 0}}\label{m8:s:m:K-open}

Equation \eqref{m8:s:m:eq:residual-stabilizer} gives $s = v = 0$. Therefore
\[
K_{F_\epsilon(E)}=\mathbb G_m(F)\times U_{4,E}(F),
\qquad
K_{G(E)}=\mathbb G_m(F)\times U_{5,E}(F),
\]
where \(U_{4,E}\) and \(U_{5,E}\) are exactly the conjugated kernels in \eqref{m8:s:m:eq:stabilizer}. The character is trivial on each whole stabilizer.

The scalar factor is the subgroup \((1,1,cI_2)\). The left-factor
flag equations in Appendix, p.~\pageref{m8:s:m:left-equations}
determine the complete stabilizer for every $E$.

\subparagraph{\texorpdfstring{The locus $d=0$}{The locus d=0}}\label{m8:s:m:K-boundary}

Now \(qr\ne0\), and the only residual condition is
\[
sr+vq=0.
\]
It has the parameterization
\[
(s,v)=(qz,-rz),\qquad z\in F.
\]
For $S_{10}$ its character is
\[
{\psi\left(\frac{2(q-r)z}{\epsilon}\right).}
\]
Therefore the $K$-compatible part of this boundary is exactly $q = r$.

For $S_{11}$ this direction has character one for every nonzero $q,r$.

The respective stabilizers have the form
\[
K_{F_\epsilon(E)}\simeq\mathbb G_m\times\mathcal U_{5,E},
\qquad
K_{G(E)}\simeq\mathbb G_m\times\mathcal U_{6,E},
\]
where the second factors are connected unipotent extensions of
the kernel by the extra one-dimensional unipotent direction.

The central scalar factor splits off as follows. Extracting $c$ from the forced diagonal matrix gives a group homomorphism onto \(\mathbb G_m\) with central section \(c\mapsto(1,1,cI)\). Its kernel projects into the unipotent group \(N_-\); it is a unipotent extension of the known kernel. The explicit lifts parameterize it over the affine line, so it is connected and contains every component. The scalar parameter ranges over all $F^\times$.

For reference, the dimensions and orbit dimensions are

\begingroup
\small
\setlength{\tabcolsep}{3pt}
\renewcommand{\arraystretch}{1.18}
\begin{longtable}{p{0.29\linewidth}p{0.16\linewidth}p{0.18\linewidth}p{0.28\linewidth}}
\toprule
Family & Condition & \(\dim K_x\) & \(\dim Kx\), with \(\dim K=24\) \\
\midrule
\endfirsthead
Family & Condition & \(\dim K_x\) & \(\dim Kx\), with \(\dim K=24\) \\
\midrule
\endhead
$S_{10}$ canonical & $d \ne 0$ & 5 & 19 \\
$S_{10}$ canonical & $d = 0$ & 6 & 18 \\
$S_{11}$ canonical & $d \ne 0$ & 6 & 18 \\
$S_{11}$ canonical & $d = 0$ & 7 & 17 \\
\bottomrule
\end{longtable}
\endgroup

\subparagraph{Modular functions}\label{m8:s:m:modulus}

Conjugation by $D$(a) acts on each of the $X$ and $Z$ spaces by ordinary conjugation of two-by-two matrices and fixes $y,z$. Its determinant on each matrix space is
\[
(\det a)^2(\det a^{-1})^2=1.
\]
The second $U$ factor has the same property with \(a^{-t}\). Hence
\[
{\Delta_K=1.}
\]
The stabilizers just described are products of a central torus with unipotent groups; therefore
\[
{\Delta_{K_x}=1,\qquad
\Delta_K|_{K_x}\Delta_{K_x}^{-1}=1.}
\]
This proves that the relative density is trivial on the two
compatible canonical families. On $S_{12},S_{15},S_{16}$ the
complete $U\times U$ stabilizers give the vanishing in
Section~\ref{m8:s:m:characters}.

\paragraph{\texorpdfstring{Intrinsic odd vanishing on $S_{10}$}{Intrinsic odd vanishing on S10}}\label{m8:s:m:S10}
\begin{proposition}[Odd vanishing on the first middle stratum]\label{m8:s:m:prop:odd10}
The stratum $S_{10}$ satisfies $D'_{\widetilde\eta_K}(S_{10})=0$.
\end{proposition}
\par\smallskip\noindent\textit{Proof.}\par

\subparagraph{Transpose transporters}\label{m8:s:m:transpose10}

Symmetry of \(V_\epsilon\) gives
\[
{D(E^{-t})F_\epsilon(E)^tD(E^t)=F_\epsilon(E^t).}
\]
For $d \ne 0$, take
\[
s=(r-q)/d,\qquad v=(q-r)/d=-s.
\]
Then
\[
n_+(s)En_-(v)=E^t.
\]
The lift \(k_{10}(s,v)\) has character one, and
\[
h_E=(1,1,E^t)\,k_{10}(s,v)
\]
satisfies
\[
h_EF_\epsilon(E)=F_\epsilon(E)^t,\qquad \eta_K(h_E)=1.
\]
If $d = 0$ and $q \ne r$, the extra stabilizer character in Section~\ref{m8:s:m:K-boundary} is nontrivial. If $d = 0$ and $q = r$, then $E$ is symmetric and the pure diagonal transporter \((1,1,E^t)\) has character one.

To transport this conclusion to an arbitrary point, suppose \(k_0x=x_c\) and \(hx_c=\tau x_c\) with \(\eta(h)=1\), then
\[
\theta(k_0)^{-1}hk_0x=\tau x,
\]
and its character is still one because \(\eta\circ\theta=\eta\). This transports the character-one transpose relation through the
$S_{10}$ normalizations in Appendix,
p.~\pageref{m8:s:m:normalization10}.

\subparagraph{The distribution argument}\label{m8:s:m:odd10}

The $U\times U$ stabilizer rank is constant on $S_{10}$. Vanishing of the character functional on that stabilizer defines a closed algebraic subset \(\mathcal C_{10}\), with the three equations in Section~\ref{m8:s:m:compatibility} in the smooth Levi presentation. Its complement is killed by the unipotent criterion of \cite[Remark~2.6]{JLX}.

On \(\mathcal C_{10}\), consider the algebraic extended group
\[
\widetilde K=K\rtimes\langle\tau\rangle.
\]
At every point, either an even stabilizer element has a nontrivial character with relative modulus one, or a character-one even transporter $h$ exists. In the latter case
\[
h^{-1}\tau\in\widetilde K_x,
\qquad \widetilde\eta_K(h^{-1}\tau)=-1.
\]
Because modular characters have positive real values, the inverse
extended character, also -1 at this element, differs from the
relative modular character. This is the stabilizer inequality
required by \cite[Theorem~2.5]{JLX} for \(\widetilde K\).

The extended action is algebraic, including its two components, so the constructibility hypothesis is satisfied. We can apply the theorem on the whole compatible locus and then use
\[
0\longrightarrow D'_{\widetilde\eta_K}(\mathcal C_{10})
\longrightarrow D'_{\widetilde\eta_K}(S_{10})
\longrightarrow D'_{\widetilde\eta_K}(S_{10}\setminus\mathcal C_{10}).
\]
Both outside terms vanish. Therefore
\[
{D'_{\widetilde\eta_K}(S_{10})=0.}
\]

Thus algebraic constructibility and the extended-stabilizer test
turn the pointwise transporter calculation into distributional vanishing.

\par\hfill$\square$\par\smallskip

\paragraph{\texorpdfstring{The compatible locus in $S_{11}$}{The compatible locus in S11}}\label{m8:s:m:compatible11}
\begin{proposition}[The compatible support locus]\label{m8:s:m:prop:Y}
The compatible locus $\mathcal C_{11}\subset S_{11}$ is closed and smooth of dimension $20$. Its closed transpose-stable subvariety $\mathcal Y$ defined in \eqref{m8:s:m:eq:Y} is smooth of dimension $19$, and $D'_{\widetilde\eta_K}(S_{11}\setminus\mathcal Y)=0$.
\end{proposition}
\par\smallskip\noindent\textit{Proof.}\par

\subparagraph{A smooth locus of dimension twenty}\label{m8:s:m:smooth20}

Let \(\mathsf A_x\) be the matrix of
\[
\mathfrak u\oplus\mathfrak u\longrightarrow T_xG_8^\circ,
\qquad (N,N')\longmapsto Nx+xN'^t,
\]
using its 64 matrix-entry coordinates. On $S_{11}$ its rank is constantly 15, because its kernel is the complete five-dimensional stabilizer in Appendix, p.~\pageref{m8:s:m:stab}. Let $\ell$ be the row representing
\[
\operatorname{tr}X+2y+\operatorname{tr}X'+2y'.
\]
Then
\[
\mathcal C_{11}
 =\left\{x\in S_{11}:
   \operatorname{rank}\begin{pmatrix}\mathsf A_x\\\ell\end{pmatrix}=15\right\}.
\]
Equivalently, all its $16\times16$ minors vanish. This is a closed algebraic condition on the stratum.

The parabolic presentation \(U\times M\times M\times U\to S_{11}\) is a smooth surjective presentation in characteristic zero: it is a Levi-coordinate form of the transitive $P \times P^t$ action, whose stabilizer is the smooth group $I_{11}$. Its $F$-point surjectivity was proved in Appendix, p.~\pageref{m8:s:m:cells}.

In this presentation, the compatible equations are exactly
\[
R_{11}=0,\quad S_{11}=0,\quad R_{12}-S_{12}=0,\quad R_{21}-S_{21}=0.
\]
They are four independent equations on the two independently variable invertible matrices $R,S$. The remaining $B,D$ and similitude coordinates are free. Smoothness and codimension therefore descend through the smooth presentation. Thus
\[
{\mathcal C_{11}\text{ is smooth of dimension }24-4=20.}
\]

\subparagraph{A smooth closed locus of dimension nineteen}\label{m8:s:m:smooth19}

When \(x_{44}=0\), the blocks \(x_{34}\) and \(x_{43}\) are unchanged by $U\times U$. This follows directly from the block multiplication: the only possible extra terms contain \(x_{44}\). Diagonal $\mathrm{GL}_2$ conjugates both blocks simultaneously.

At a canonical point $G$($E$),
\[
x_{34}=x_{43}=\begin{pmatrix}0&0\\r&d\end{pmatrix}.
\]
Throughout \(\mathcal C_{11}\), these blocks are consequently
equal and have rank one. We define
\begin{equation}
{
\mathcal Y=\mathcal C_{11}\cap\{\operatorname{tr}x_{34}=0\}
 =K\cdot\left\{G\!\left(\begin{pmatrix}p&q\\r&0\end{pmatrix}\right):qr\ne0\right\}.
}
\label{m8:s:m:eq:Y}
\end{equation}
It is closed in $S_{11}$. At a canonical point the extra equation is $d = 0$; varying the $E_{22}$ entry in the invertible canonical family has derivative 1 on \(\operatorname{tr}x_{34}\). Hence this is a regular equation on the smooth compatible locus. Thus
\[
{\mathcal Y\text{ is smooth of dimension }19.}
\]
Both \(\mathcal C_{11}\) and \(\mathcal Y\) are transpose-stable. For the second statement, transpose exchanges \(x_{34}\) and \(x_{43}^t\), and those blocks agree on \(\mathcal C_{11}\).

\subparagraph{Odd vanishing outside the closed locus}\label{m8:s:m:complement11}

Since $V$ is symmetric,
\[
{D(E^{-t})G(E)^tD(E^t)=G(E^t).}
\]
For $d \ne 0$, the parameters \(s=(r-q)/d,v=-s\) and the $S_{11}$ lift in Appendix, p.~\pageref{m8:s:m:lifts} give a character-one transporter to $G$($E^t$). Following it by diagonal conjugation through $E^t$ gives one to $G$($E$)$^t$. The extended-stabilizer criterion from Section~\ref{m8:s:m:odd10} therefore kills odd distributions on \(\mathcal C_{11}\setminus\mathcal Y\). The incompatible complement \(S_{11}\setminus\mathcal C_{11}\) is killed by $U\times U$.

The closed/open sequences are
\[
0\longrightarrow D'_{\widetilde\eta_K}(\mathcal C_{11}\setminus\mathcal Y)
\longrightarrow D'_{\widetilde\eta_K}(S_{11}\setminus\mathcal Y)
\longrightarrow D'_{\widetilde\eta_K}(S_{11}\setminus\mathcal C_{11}),
\]
followed by
\begin{equation}
0\longrightarrow D'_{\widetilde\eta_K}(\mathcal Y)
\longrightarrow D'_{\widetilde\eta_K}(S_{11})
\longrightarrow D'_{\widetilde\eta_K}(S_{11}\setminus\mathcal Y).
\label{m8:s:m:eq:Y-sequence}
\end{equation}
Thus the rightmost term in \eqref{m8:s:m:eq:Y-sequence} is zero. This uses the ordinary left-exact equivariant support sequence. \cite[Section~2.1]{JLX} explicitly states the unequivariant exact sequence and the left-exact equivariant sequence. 

The microlocal support argument uses the entire closed smooth
$\mathcal Y$, including its $q=r$ locus.

\par\hfill$\square$\par\smallskip

\paragraph{Rank boundaries and subsequent support steps}
The union $S_{12}\sqcup S_{15}\sqcup S_{16}$ is locally closed,
with the paired strata exchanged by transpose. Its even vanishing
follows from the character theorem and the closed/open sequences
in Appendix, p.~\pageref{m8:s:m:familyB}. The degenerations and predecessor
envelopes are verified in Appendix, p.~\pageref{m8:s:m:topology} of the
mathematical appendix. On $S_{11}$ the remaining closed smooth space $\mathcal Y$
supports nonzero odd orbital distributions; it is removed only
by the conditional nilpotent-wavefront theorem below, on its stated open
and after its predecessors have been cleared.

\subsection{The two spectral support steps}\label{m8:spectral}
Although $S_3$ and $S_{11}$ support nonzero intrinsic odd
distributions, their intrinsic reductions confine the support to
smooth closed families with nonnilpotent ambient conormals. We
exclude these families for the same ambient NW distribution after
clearing the specified predecessor opens. For $S_3$ we calculate
the conormal in the original ambient group; for $S_{11}$ the
calculation covers the whole nineteen-dimensional family, including
its diagonal parameter locus.
\subsubsection{\texorpdfstring{Spectral support on the Siegel stratum $S_3$}{Spectral support on the Siegel stratum S3}}\label{m8:s:sec:s3-s8}
We now analyze the rank-one locus of the common Siegel chart
\eqref{m8:s:s3:eq:1-8}--\eqref{m8:s:s3:eq:3-3}, whose rational
coordinates, residual action and density were established on
p.~\pageref{m8:common-siegel}. In contrast to the intrinsic odd
vanishing on $S_8$ proved in Proposition~\ref{m8:s:s3:prop:S8},
$S_3$ supports intrinsic odd distributions. We exclude the support
induced by the ambient NW distribution through the conormal
calculation below.

Split the rank-one subset into
\(\mathcal X_{\rm ss}=\{\operatorname{rank}D=1,\ \operatorname{tr}D\ne0\},\) \(\mathcal X_{\rm nil}=\{\operatorname{rank}D=1,\ \operatorname{tr}D=0\}.\)
These are the open and closed parts of $\mathcal X_1$, respectively.

\paragraph{\texorpdfstring{The semisimple terminal block in $S_3$}{The semisimple terminal block in S3}}\label{m8:s:s3:semisimple}
\begin{proposition}[The semisimple locus]\label{m8:s:s3:prop:semisimple}
The semisimple terminal-block locus satisfies $D'_{\widetilde\chi}(\mathcal X_{\rm ss}\times F^\times)=0$.
\end{proposition}
\par\smallskip\noindent\textit{Proof.}\par

\subparagraph{Normal form and first descent}\label{m8:s:s3:semisimple-normal}

Put \(d=\operatorname{tr}D\ne0\). By simultaneous conjugation, take
\(D=\operatorname{diag}(0,d).\)
To obtain this normalization over $F$, choose an ordered basis
from nonzero columns of the projectors \(I-D/d\) and \(D/d\).
These choices give algebraic sections on the corresponding open
charts. The stabilizer is the diagonal torus
\(T\subset \mathrm{GL}_2\), and both groups are unimodular;
hence Frobenius descent over \(\mathrm{GL}_2/T\) is untwisted,
with \(d\) retained as a free parameter.

Write \(B=(b\ B_2)\), \(C=\binom{c}{C_2}\), with \(b\) a column and \(c\) a row. Set
\(P=-d^{-1}B_2e_2^t,\) \(R=-d^{-1}e_2C_2.\)
Then \eqref{m8:s:s3:eq:3-1} sends \(g\) to
\begin{equation}
\begin{pmatrix}
A_0&b&0\\
c&0&0\\
0&0&d
\end{pmatrix},
\qquad A_0=A-d^{-1}B_2C_2,
\label{m8:s:s3:eq:4-1}
\end{equation}
where the displayed sizes are \(2+1+1\). Thus
\[
Y=\begin{pmatrix}A_0&b\\c&0\end{pmatrix}\in \mathrm{GL}_3(F).
\]
The removed four-dimensional subgroup has coordinates
\(P=p_0e_2^t,R=e_2r_0\). Its forward coordinates are
\(B_2=dp_0,\) \(C_2=dr_0,\) \(A=A_0+dp_0r_0,\)
and character \(\psi((p_0)_2+(r_0)_2)\). Its action is free, with Haar \(dp_0\,dr_0=|d|^{-4}dB_2\,dC_2\). The factor depends only on the fixed parameter \(d\).

The remaining unipotents are
\(P=pe_1^t,\) \(R=e_1r,\)
acting on \(Y\) by
\begin{equation}
A_0\longmapsto A_0+pc+br,
\qquad \chi=\psi(p_1+r_1).
\label{m8:s:s3:eq:4-2}
\end{equation}
If \(h=\operatorname{diag}(t_1,t_2)\), write \(q=t_1/t_2\). Its action on \(Y\) is conjugation by \(\operatorname{diag}(q,1,q)\). Transpose is ordinary transpose of \(Y\).  

The weights \((q,1)\) and \((q^{-1},1)\) of the removed
coordinate pairs cancel in their combined density. On the four
surviving additive coordinates the determinant is one, and transpose
exchanges the two sides with absolute Jacobian one. Thus
\eqref{m8:s:s3:eq:4-2} has invariant density.

\subparagraph{Unipotent stabilizer and compatibility}\label{m8:s:s3:semisimple-compatible}

Invertibility gives $b,c\ne0$. Applying \eqref{m8:s:c:affine-line}, the complete stabilizer is
\begin{equation}
p=\lambda b,\qquad r=-\lambda c,\qquad\lambda\in F.
\label{m8:s:s3:eq:4-3}
\end{equation}
Its character is
\[
\psi\bigl(\lambda(b_1-c_1)\bigr).
\]
The unipotent stabilizer criterion therefore forces
\begin{equation}
b_1=c_1.
\label{m8:s:s3:eq:4-4}
\end{equation}
The rational-point criterion \eqref{cf:orbital} applies on the entire incompatible open set.

Intrinsically before choosing the eigenbasis, the two scalars are
\[
b_1=\operatorname{tr}\bigl(B(I-D/d)\bigr),\qquad
c_1=\operatorname{tr}\bigl(C(I-D/d)\bigr).
\]
They are unchanged by the additive actions and simultaneous conjugation. Thus \eqref{m8:s:s3:eq:4-4} defines a closed stable subset of this locus.

\subparagraph{The nonzero compatible scalar and the whole opposite-weight plane}\label{m8:s:s3:plane}

Suppose \(b_1=c_1=\beta\ne0\), and write \(b=(\beta,z)^t\), \(c=(\beta,w)\). For \(A_0=(a_{ij})\), choose
\(p_1=0,\) \(r_1=-a_{11}/\beta,\) \(r_2=-a_{12}/\beta,\) \(p_2=-a_{21}/\beta+za_{11}/\beta^2.\)
This gives the normal form
\begin{equation}
Y=\begin{pmatrix}0&0&\beta\\0&\alpha&z\\\beta&w&0\end{pmatrix},
\quad
\alpha=a_{22}-\frac{wa_{21}+za_{12}}\beta+
\frac{zw a_{11}}{\beta^2}.
\label{m8:s:s3:eq:4-5}
\end{equation}
Here \(\det Y=-\alpha\beta^2\), so \(\alpha\ne0\). The character of this normalizing element is \(\psi(-a_{11}/\beta)\).

The three-dimensional affine orbit in \(A_0\) contributes the character function
\(\psi(a_{11}/\beta)\)
and, using \(d\lambda\) on \eqref{m8:s:s3:eq:4-3}, quotient Haar density
\begin{equation}
|\beta|^{-2}\,da_{11}\,da_{12}\,da_{21}.
\label{m8:s:s3:eq:4-6}
\end{equation}
For example, the complement \(p_1=0\) has Jacobian \(|\beta|^3\) onto these three entries, while adjoining the stabilizer coordinate has Jacobian \(|\beta|\). Both the phase and \eqref{m8:s:s3:eq:4-6} are invariant under the induced torus and transpose transformations.

The exact remaining problem has parameters \((\alpha,\beta,d,\nu)\in(F^\times)^4\), on which the group acts trivially, and a plane \((z,w)\in F^2\), with
\begin{equation}
q(z,w)=(q^{-1}z,qw),\qquad \sigma(z,w)=(w,z),\qquad
\chi(q)=1,\quad\widetilde\chi(\sigma)=-1.
\label{m8:s:s3:eq:4-7}
\end{equation}
Here \(F^\times\rtimes\langle\sigma\rangle\), with
\(\sigma q\sigma=q^{-1}\), acts on the whole plane.

Apply Lemma~\ref{m8:s:c:plane} with its torus parameter replaced by $q^{-1}$ and passive $(\alpha,\beta,d,\nu)$. Its domain is the $(z,w)$-plane, including both punctured axes and their common origin. On $zw\ne0$ the transporter is $q=z/w$; the axes are exchanged distinct orbits, covered by the lemma's extension argument.

The residual stabilizer at \eqref{m8:s:s3:eq:4-5} consists of the line \eqref{m8:s:s3:eq:4-3}, scalar \(h\), and a possible ratio \(q\) satisfying \(qz=z\), \(qw=w\). Thus \(q=1\) unless \(z=w=0\). These stabilizers are unimodular; the extra torus when both vanish acts trivially on the stabilizer-line coordinate.

\subparagraph{The zero compatible scalar}\label{m8:s:s3:zero-scalar}

Now \(b_1=c_1=0\). Invertibility forces \(z=b_2\ne0\), \(w=c_2\ne0\), and \(\alpha=a_{11}\ne0\). Choose
\(r_1=-a_{21}/z,\) \(p_1=-a_{12}/w,\) \(p_2=-a_{22}/w,\) \(r_2=0.\)
The normal form is
\begin{equation}
Y=\begin{pmatrix}\alpha&0&0\\0&0&z\\0&w&0\end{pmatrix}.
\label{m8:s:s3:eq:4-8}
\end{equation}
The removed directions contribute
\(\psi(a_{12}/w+a_{21}/z),\) \(|zw|^{-1}\,da_{12}\,da_{21}\,da_{22}.\)
Both transpose and the torus preserve this phase and density.
The residual problem is therefore \eqref{m8:s:s3:eq:4-7} on
\((F^\times)^2\), with \((\alpha,d,\nu)\) passive. At every
point the character-one transporter \(q=z/w\) gives vanishing,
as in the nonzero-product part of the opposite-weight-plane lemma.

The incompatible complement, the open compatible locus \(\beta\ne0\), and its closed remainder \(\beta=0\) therefore all vanish in the required order. We have proved
\begin{equation}
D'_{\widetilde\chi}(\mathcal X_{\rm ss}\times F^\times)=0.
\label{m8:s:s3:eq:4-10}
\end{equation}

\par\hfill$\square$\par\smallskip

\paragraph{\texorpdfstring{The nilpotent terminal block in $S_3$}{The nilpotent terminal block in S3}}\label{m8:s:s3:nilpotent}

\subparagraph{Rational normalization}\label{m8:s:s3:nilpotent-normal}

All nonzero nilpotent two-by-two matrices are conjugate over \(F\) to
\[
N=E_{12}=\begin{pmatrix}0&1\\0&0\end{pmatrix}.
\]
For explicit local sections, choose \(v\) with \(Dv\ne0\).
Since \(D(Dv)=0\), the basis \((Dv,v)\) represents \(D\) by
\(N\); on each corresponding open chart we may choose \(v\)
to be a standard basis vector.

The centralizer is
\begin{equation}
C_N=\{\zeta(I_2+kN):\zeta\in F^\times,\ k\in F\}.
\label{m8:s:s3:eq:5-1}
\end{equation}
Both \(\mathrm{GL}_2\) and \(C_N\) are unimodular, as are
\(E\rtimes \mathrm{GL}_2\) and \(E\rtimes C_N\). Frobenius
descent from the nilpotent conjugacy orbit to \(D=N\) therefore
has relative modular character one.

Write \(B=(b\ B_2)\), \(C=\binom{C_1}{c}\). Set
\(P=-B_2e_1^t,\) \(R=-e_2C_1.\)
The exact normalized matrix is
\begin{equation}
g(A_0,b,c)=
\begin{pmatrix}
(A_0)_{11}&(A_0)_{12}&b_1&0\\
(A_0)_{21}&(A_0)_{22}&b_2&0\\
0&0&0&1\\
c_1&c_2&0&0
\end{pmatrix},
\quad A_0=A-B_2C_1.
\label{m8:s:s3:eq:5-2}
\end{equation}
The four removed coordinates are \(P=p_0e_1^t,R=e_2r_0\), with
\begin{equation}
B_2=p_0,\quad C_1=r_0,\quad A=A_0+p_0r_0,
\quad\chi=\psi((p_0)_1+(r_0)_2).
\label{m8:s:s3:eq:5-3}
\end{equation}
Their Jacobian is one. The element used to reach \eqref{m8:s:s3:eq:5-2} has character \(\psi(-B_{12}-C_{12})\).

Put
\begin{equation}
w=\begin{pmatrix}0&1\\1&0\end{pmatrix},\qquad
\sigma(g)=d(w)g^t d(w).
\label{m8:s:s3:eq:5-4}
\end{equation}
This involution preserves \(D=N\). Since \(d(w)\) is part of the residual diagonal action and its character is one, oddness under transpose is equivalent, in this fixed-\(D\) problem, to oddness under \(\sigma\).

On \eqref{m8:s:s3:eq:5-2},
\begin{equation}
\sigma:(A_0,b,c)\longmapsto(wA_0^tw,\ wc^t,\ b^tw).
\label{m8:s:s3:eq:5-5}
\end{equation}
It exchanges the removed coordinate pairs, with the fixed internal swaps by \(w\), and preserves their Haar densities and character. On \eqref{m8:s:s3:eq:5-1}, it induces \((\zeta,k)\mapsto(\zeta^{-1},-k)\), preserving \(d^\times\zeta\,dk\). Thus the nilpotent-orbit descent has trivial relative density.

\subparagraph{The additive stabilizer}\label{m8:s:s3:additive-stabilizer}

The surviving additive coordinates are now
\(P=pe_2^t,\) \(R=e_1r.\)
They act by
\(A_0\longmapsto A_0+pc+br,\) \(\chi=\psi(p_2+r_1).\)
Invertibility gives $b,c\ne0$; \eqref{m8:s:c:affine-line} gives the additive stabilizer
\begin{equation}
p=\lambda b,\qquad r=-\lambda c.
\label{m8:s:s3:eq:5-7}
\end{equation}
Its character is
\begin{equation}
\psi\bigl(\lambda(b_2-c_1)\bigr).
\label{m8:s:s3:eq:5-8}
\end{equation}
Thus every equivariant distribution is supported on
\(b_2=c_1=\rho.\)
Before normalization, the two quantities are the invariant scalars
\(\rho_L=\operatorname{tr}(BD),\) \(\rho_R=\operatorname{tr}(CD).\)
Indeed, \(D^2=0\) makes these unchanged by \eqref{m8:s:s3:eq:3-1}, and they are invariant under simultaneous conjugation. Transpose exchanges them. Hence their equality is a closed stable support condition.

The remaining residual stabilizer direction $I+kN$ will supply a second compatibility equation.

\subparagraph{\texorpdfstring{The locus $\rho\ne0$}{The locus rho != 0}}\label{m8:s:s3:nonzero-rho}

Write \(b=(u,\rho)^t\), \(c=(\rho,v)\), with \(\rho\ne0\). For \(A_0=(a_{ij})\), choose
\(r_1=0,\) \(p_1=-a_{11}/\rho,\) \(p_2=-a_{21}/\rho,\) \(r_2=-a_{22}/\rho+va_{21}/\rho^2.\)
This yields
\begin{equation}
g_\rho(\alpha,u,v)=
\begin{pmatrix}
0&\alpha&u&0\\
0&0&\rho&0\\
0&0&0&1\\
\rho&v&0&0
\end{pmatrix},
\quad
\alpha=a_{12}-\frac{va_{11}+ua_{22}}\rho+
\frac{uv a_{21}}{\rho^2}.
\label{m8:s:s3:eq:5-11}
\end{equation}
Here \(\det g_\rho=-\alpha\rho^2\), so \(\alpha\ne0\). The normalization character is \(\psi(-a_{21}/\rho)\); the corresponding affine-orbit phase and density are
\(\psi(a_{21}/\rho),\) \(|\rho|^{-2}\,da_{11}\,da_{21}\,da_{22}.\)
Both are preserved by \eqref{m8:s:s3:eq:5-5}, which swaps \(u,v\), swaps \(a_{11},a_{22}\), and fixes \(a_{21}\).

Let
\[
k=\frac{v-u}{\rho},\qquad
h=I_2+kN,
\]
\begin{equation}
P=k\begin{pmatrix}v&-kv\\\rho&-k\rho\end{pmatrix},
\qquad
R=-k\begin{pmatrix}0&0\\\rho&u\end{pmatrix}.
\label{m8:s:s3:eq:5-13}
\end{equation}
Direct substitution in \eqref{m8:s:s3:eq:3-1}--\eqref{m8:s:s3:eq:3-3} gives
\((P,R,h)g_\rho(\alpha,u,v) =\sigma(g_\rho(\alpha,u,v))=g_\rho(\alpha,v,u).\)
Moreover,
\begin{equation}
\operatorname{tr}P+\operatorname{tr}R
=k(v-u)-k^2\rho=0.
\label{m8:s:s3:eq:5-15}
\end{equation}
Composing with \((0,0,w)\) gives a transpose transporter of character one.

The ordinary stabilizer at \eqref{m8:s:s3:eq:5-11} consists of \eqref{m8:s:s3:eq:5-7} and scalar \(h\). Indeed, a centralizer element \(I+kN\) changes \(b_1\) to \(b_1+k\rho\), so fixing the point forces \(k=0\). The stabilizer is unimodular and has trivial character.

Lemma~\ref{m8:s:c:transporter} now gives vanishing on the whole compatible open locus. Its formula \eqref{m8:s:c:transport-back} transports these explicit character-one elements to the original points, regardless of the normalization character.

\subparagraph{\texorpdfstring{The locus $\rho=0$}{The locus rho=0}}\label{m8:s:s3:zero-rho}

Suppose now \(\rho=0\). Invertibility forces
\(b=(b,0)^t,\) \(c=(0,c),\) \(b,c\in F^\times.\)
Choose
\(r_1=-a_{11}/b,\) \(p_1=-a_{12}/c,\) \(p_2=-a_{22}/c,\) \(r_2=0.\)
This gives
\begin{equation}
g(a,b,c)=
\begin{pmatrix}
0&0&b&0\\
a&0&0&0\\
0&0&0&1\\
0&c&0&0
\end{pmatrix},
\qquad a,b,c\in F^\times.
\label{m8:s:s3:eq:5-17}
\end{equation}
Its determinant is \(-abc\). The removed affine directions contribute
\begin{equation}
\psi(a_{11}/b+a_{22}/c),\qquad
|bc|^{-1}\,da_{11}\,da_{12}\,da_{22}.
\label{m8:s:s3:eq:5-18}
\end{equation}
Transpose in the fixed-\(D\) slice acts by
\begin{equation}
\sigma(g(a,b,c))=g(a,c,b).
\label{m8:s:s3:eq:5-19}
\end{equation}
It preserves the phase and density in \eqref{m8:s:s3:eq:5-18}.

For $H=I+kN$, conjugation followed by the four-coordinate normalization \eqref{m8:s:s3:eq:5-3} sends
\begin{equation}
b\mapsto b'=Hb,\qquad c\mapsto c'=cH^{-1},
\qquad
A_0\mapsto HA_0H^{-1}+k^2b'c'.
\label{m8:s:s3:eq:5-20}
\end{equation}
The correcting element has
\(P_0=kb'e_1^t,\) \(R_0=-ke_2c',\)
and trace character exponent
\(k(b_1-c_2)+k^2(b_2+c_1).\)
The Jacobian of \eqref{m8:s:s3:eq:5-20} is one. The stabilizer calculation includes this character cocycle and the subsequent affine normalization.

\subparagraph{The stabilizer and second compatibility equation}\label{m8:s:s3:full-stabilizer}

Solving \eqref{m8:s:s3:eq:3-3} at \eqref{m8:s:s3:eq:5-17} gives all stabilizer elements:
\begin{equation}
\begin{aligned}
h&=\zeta\begin{pmatrix}1&k\\0&1\end{pmatrix},\qquad \zeta\in F^\times,\quad k,v\in F,\\
P&=\begin{pmatrix}
bk&-\dfrac bc v+\left(\dfrac ac-b\right)k^2\\
0&\dfrac ac k
\end{pmatrix},\\
R&=\begin{pmatrix}
-\dfrac ab k&v\\
0&-ck
\end{pmatrix}.
\end{aligned}
\label{m8:s:s3:eq:5-22}
\end{equation}
The additive stabilizer \eqref{m8:s:s3:eq:5-7} is the subfamily \(k=0\). The trace exponent on the stabilizer is
\begin{equation}
\operatorname{tr}P+\operatorname{tr}R
=k\,\frac{(b-c)(a+bc)}{bc}.
\label{m8:s:s3:eq:5-23}
\end{equation}
As a product of scalar \(F^\times\) with a two-dimensional
unipotent group, the stabilizer is unimodular. Its effective
centralizer direction \(k\) gives the additional compatibility test.

Consequently,
\begin{equation}
b=c\quad\text{or}\quad a=-bc
\label{m8:s:s3:eq:5-24}
\end{equation}
is necessary for a compatible orbital distribution. If \(b=c\), \eqref{m8:s:s3:eq:5-19} fixes the normal form, and \((0,0,w)\) is a character-one transporter to its ordinary transpose. The transporter criterion gives odd vanishing on these orbits.

But on
\begin{equation}
a=-bc,\qquad b\ne c,
\label{m8:s:s3:eq:5-25}
\end{equation}
the character on the stabilizer is identically one and no transpose transporter exists. To see the latter directly, every transporter to \eqref{m8:s:s3:eq:5-19} must have \(h\in C_N\). On the \(\rho=0\) normal forms, the two quantities \(b,c\) are individually invariant under the residual action preserving \(D=N\). Equation \eqref{m8:s:s3:eq:5-20} shows this explicitly. Thus \((b,c)\) cannot become \((c,b)\) unless \(b=c\).

\paragraph{\texorpdfstring{Closed orbits and intrinsic odd distributions on $S_3$}{Closed orbits and intrinsic odd distributions on S3}}\label{m8:s:s3:obstruction}
\begin{proposition}[The intrinsic obstruction]\label{m8:s:s3:prop:intrinsic}
Every odd equivariant distribution on $\mathcal X_1\times F^\times$ is supported on the closed smooth locus $\mathcal N\times F^\times$ defined below. The ordered parameters $(b,c)$ classify its residual orbits. If $b\ne c$, the corresponding closed $K$-orbits in $S_3$ give nonzero odd orbital distributions. In particular, $D'_{\widetilde\eta_K}(S_3)\ne0$.
\end{proposition}
\par\smallskip\noindent\textit{Proof.}\par

\subparagraph{The intrinsic closed support locus}\label{m8:s:s3:closed-locus}

Let
\(\mathcal X^*=\{g\in \mathrm{GL}_4(F):D\ne0\}.\)
On the locus
\begin{equation}
\det D=\operatorname{tr}D=0,\qquad
\operatorname{tr}(BD)=\operatorname{tr}(CD)=0,
\label{m8:s:s3:eq:6-1}
\end{equation}
there are unique scalars \(b,c\) such that
\begin{equation}
BD=bD,\qquad DC=cD.
\label{m8:s:s3:eq:6-2}
\end{equation}
These are regular functions locally wherever a chosen entry of \(D\) is nonzero. Define
\begin{equation}
a=\operatorname{tr}(DA)-
\bigl(\operatorname{tr}B-b\bigr)
\bigl(\operatorname{tr}C-c\bigr).
\label{m8:s:s3:eq:6-3}
\end{equation}
The normalizations \eqref{m8:s:s3:eq:5-2} and
\eqref{m8:s:s3:eq:5-17} show that \(a,b,c\ne0\) and identify
them with the three normal-form parameters.

For example, at \(D=N\), \eqref{m8:s:s3:eq:6-1} says
\[
B=\begin{pmatrix}b&\beta_0\\0&\delta_0\end{pmatrix},\qquad
C=\begin{pmatrix}\alpha_0&\gamma_0\\0&c\end{pmatrix},
\]
and \eqref{m8:s:s3:eq:6-3} is
\(a=A_{21}-\delta_0\alpha_0.\)
This is unchanged by the additive transformations \eqref{m8:s:s3:eq:3-1}. Simultaneous conjugation preserves \eqref{m8:s:s3:eq:6-2}--\eqref{m8:s:s3:eq:6-3}, so these are invariant functions.

Define
\(\mathcal N=\{g\in\mathcal X^*:\text{\eqref{m8:s:s3:eq:6-1} holds and }a+bc=0\}.\)
This is a closed subset of \(\mathcal X^*\), including its
\(b=c\) points. To see that it is smooth of dimension eleven,
use local Jordan-basis sections together with the four free
coordinates in \eqref{m8:s:s3:eq:5-3}, the three affine coordinates
removed in \eqref{m8:s:s3:eq:5-18} and the two parameters
\((b,c)\). They give local product charts of dimension
\begin{equation}
2+4+3+2=11.
\label{m8:s:s3:eq:6-5}
\end{equation}
The two-dimensional first term is the nilpotent conjugacy orbit. All denominators in the subsequent normalizations are nonzero on this locus.

Transpose acts on \(\mathcal N\) by \((b,c)\mapsto(c,b)\),
fixing both \(a\) and the scalar \(\nu\).

\subparagraph{The support filtration}\label{m8:s:s3:support}
Put
\[
\begin{aligned}
Z_0&=\mathcal X_1\times F^\times,& Z_1&=\mathcal X_{\rm nil}\times F^\times,\\
Z_2&=Z_1\cap\{\rho_L=\rho_R\},&Z_3&=Z_2\cap\{\rho_L=\rho_R=0\},\\
Z_4&=\mathcal N\times F^\times.
\end{aligned}
\]
Each next set is closed in the preceding one: the first by $\tr D=0$, the next two by the invariant scalar equalities, and the last by the regular equation $a+bc=0$. Apply Lemma~\ref{cf:support} to these four specified pairs:
\begin{center}
\begin{tabular}{c p{.30\linewidth} p{.49\linewidth}}
\toprule
$(X,Z)$&Open complement&Vanishing input\\
\midrule
$(Z_0,Z_1)$&$\mathcal X_{\rm ss}\times F^\times$&The two semisimple cases, \eqref{m8:s:s3:eq:4-10}.\\
$(Z_1,Z_2)$&$\rho_L\ne\rho_R$&Additive character \eqref{m8:s:s3:eq:5-8}.\\
$(Z_2,Z_3)$&$\rho\ne0$&Transporter \eqref{m8:s:s3:eq:5-13}--\eqref{m8:s:s3:eq:5-15}.\\
$(Z_3,Z_4)$&$a+bc\ne0$&For $b\ne c$, the stabilizer trace exponent \eqref{m8:s:s3:eq:5-23}; for $b=c$, the transporter from \eqref{m8:s:s3:eq:5-19}.\\
\bottomrule
\end{tabular}
\end{center}
Constructibility in the last row gives vanishing on the entire open complement. Thus every odd equivariant distribution on $Z_0$ is supported on $Z_4$, including its $b=c$ points.

\subparagraph{Nonzero odd orbital distributions}\label{m8:s:s3:odd-orbital}

For fixed \(b,c\in F^\times\), \(b\ne c\), and \(\nu\in F^\times\), set
\(x_{b,c,\nu}=m(g(-bc,b,c),\nu).\)
The product $m(g(-bc,b,c),\nu)$ specifies its full matrix;
\eqref{m8:s:s3:eq:1-11} proves it lies in $S_3$.
Its expanded $8\times8$ entries are preserved in
\texttt{anc/model8\_calculations.pdf}, calculation M8-C07,
under ``Nonzero odd orbital distributions''.

By the explicit normalizations, every point of \(\mathcal N\)
with fixed ordered \((b,c)\) lies in the residual orbit of
\(g(-bc,b,c)\). Each fiber is therefore a single orbit, closed
in \(\mathcal X^*\). After restoring the two free radicals and fixing \(\nu\), the corresponding \(K\)-orbit
\(\mathcal O_{b,c,\nu}=Kx_{b,c,\nu}\)
is closed in \(\Omega_4\setminus S_8\), and hence in \(S_3\).

Its stabilizer is \eqref{m8:s:s3:eq:5-22}, embedded as
\((\ell(P),\ell(R^t),h)\in K,\)
with \(a=-bc\). By uniqueness in \eqref{m8:s:s3:eq:1-8}, this is the full
stabilizer. It is unimodular and has trivial character. For example, in \eqref{m8:s:s3:eq:5-22} its two traces separately vanish when \(a=-bc\).

Consequently, \(K/K_{x_{b,c,\nu}}\) has an invariant quotient measure, and
\begin{equation}
T_{b,c,\nu}(f)=
\int_{K/K_{x_{b,c,\nu}}}
f(kx_{b,c,\nu})\eta_K(k)\,d\dot k
\label{m8:s:s3:eq:6-9}
\end{equation}
is a nonzero equivariant distribution on that closed orbit, hence on \(S_3\). The integral is finite for a compactly supported test function: the orbit is closed in \(S_3\), and the quotient parametrizes it homeomorphically.

Transpose exchanges \(\mathcal O_{b,c,\nu}\) and \(\mathcal O_{c,b,\nu}\), which are disjoint for \(b\ne c\). Therefore
\(T_{b,c,\nu}-\tau T_{b,c,\nu}\ne0\)
is an intrinsic \((\widetilde K,\widetilde\eta_K)\)-equivariant distribution on \(S_3\). Thus $D'_{\widetilde\eta_K}(S_3)\ne0$.

The orbit closure in \(\Omega_4\) meets \(S_8\): the family
\begin{equation}
g_\varepsilon=
\begin{pmatrix}
0&\operatorname{diag}(b,c)\\
\operatorname{diag}(b,c)&\varepsilon E_{12}
\end{pmatrix},\qquad\varepsilon\ne0,
\label{m8:s:s3:eq:6-11}
\end{equation}
belongs to the same residual orbit and tends as \(\varepsilon\to0\) to the rank-zero locus.

\par\hfill$\square$\par\smallskip

\paragraph{\texorpdfstring{The conormal and spectral exclusion on $S_3$}{The conormal and spectral exclusion on S3}}\label{m8:s:s3:conormal}
\begin{proposition}[The spectral exclusion]\label{m8:s:s3:prop:spectral}
The closed submanifold $\widehat{\mathcal N}$ of the open $\mathcal M$ has dimension $24$ and codimension $5$. Its ambient conormal is given by \eqref{m8:s:s3:eq:7-7}, and every point has a nonnilpotent conormal direction. Consequently an original odd equivariant distribution on $G_8^\circ$ with nilpotent wavefront fibers on $\mathcal M$ vanishes on $\mathcal M$.
\end{proposition}
\par\smallskip\noindent\textit{Proof.}\par

\subparagraph{The residual tangent space}\label{m8:s:s3:tangent}

Write \(g=g(-bc,b,c)\) as in \eqref{m8:s:s3:eq:5-17} with $a=-bc$, and set
\[
L(P)=\begin{pmatrix}0&P\\0&0\end{pmatrix},\qquad
R(Q_0)=\begin{pmatrix}0&0\\Q_0&0\end{pmatrix},\qquad
D_0(H_0)=\operatorname{diag}(H_0,H_0).
\]
The complete tangent space is
\begin{equation}
\begin{aligned}
T_g\mathcal N={}&
\{L(P)g+gR(Q_0)+[D_0(H_0),g]:
P,Q_0,H_0\in M_2(F)\}\\
&+F\,\partial_b g+F\,\partial_c g.
\end{aligned}
\label{m8:s:s3:eq:7-1}
\end{equation}
Here
\(\partial_b g=E_{13}-cE_{21},\) \(\partial_c g=E_{42}-bE_{21}.\)
The orbit term has dimension \(12-3=9\), by the stabilizer \eqref{m8:s:s3:eq:5-22}, and the two invariant parameters contribute two additional directions. Thus \eqref{m8:s:s3:eq:7-1} has dimension eleven, agreeing with \eqref{m8:s:s3:eq:6-5}.

\subparagraph{The five-dimensional conormal}\label{m8:s:s3:residual-conormal}

Represent a cotangent vector by \(\dot g\mapsto\operatorname{tr}(Z\dot g)\). Solving the annihilation equations from \eqref{m8:s:s3:eq:7-1} gives
\begin{equation}
Z=
\begin{pmatrix}
0&-\dfrac{\delta}{bc}&0&\dfrac{b\alpha-\beta+\epsilon}{c}\\
0&0&0&-\dfrac{\delta}{c}\\
-\dfrac{\delta}{b}&\alpha&\beta&\gamma\\
0&0&\delta&\epsilon
\end{pmatrix},
\quad\alpha,\beta,\gamma,\delta,\epsilon\in F.
\label{m8:s:s3:eq:7-2}
\end{equation}
This five-dimensional space is the full conormal: direct multiplication
annihilates every term in \eqref{m8:s:s3:eq:7-1}.

Take \(\delta=1\) and the other four parameters zero. Then
\begin{equation}
Z_*=
\begin{pmatrix}
0&-1/(bc)&0&0\\
0&0&0&-1/c\\
-1/b&0&0&0\\
0&0&1&0
\end{pmatrix},
\qquad
Z_*g=H_*:=\operatorname{diag}(1,-1,-1,1).
\label{m8:s:s3:eq:7-3}
\end{equation}
In particular, the residual conormal already has a nonnilpotent representative after left trivialization. Every matrix in \eqref{m8:s:s3:eq:7-2} satisfies \(\operatorname{tr}(Zg)=0\).

\subparagraph{The conormal in the orthogonal group}\label{m8:s:s3:orthogonal-conormal}

Set
\(\mathcal M=\Omega_4\setminus S_8.\)
This is open in \(G_8^\circ\): equivalently, the lower-right four-by-four block is invertible and \(x_{44}\ne0\). By \eqref{m8:s:s3:eq:1-11}, it corresponds to \(g\in\mathcal X^*\).

The closed support enlargement in the group is
\begin{equation}
\widehat{\mathcal N}
=N_4\cdot\{m(g,\nu):g\in\mathcal N,\ \nu\in F^\times\}\cdot N_4^t
\subset\mathcal M.
\label{m8:s:s3:eq:7-4}
\end{equation}
The global chart proves that it is a closed smooth submanifold of the open ambient set \(\mathcal M\), of dimension
\(6+11+1+6=24.\)
Its codimension in \(G_8^\circ\) is five.

At \(m=m(g,\nu)\), the tangent space is exactly
\begin{equation}
T_m\widehat{\mathcal N}
=\mathfrak n_4m+m\mathfrak n_4^t+
\{dm_{(g,\nu)}(\dot g,\dot\nu):\dot g\in T_g\mathcal N,\ \dot\nu\in F\}.
\label{m8:s:s3:eq:7-5}
\end{equation}
For a Levi variation, left trivialization gives
\begin{equation}
m^{-1}\dot m=
\operatorname{diag}\left(Y,\lambda I_4-\mathcal J^{-1}Y^t\mathcal J\right),
\quad Y=g^{-1}\dot g,\quad\lambda=\dot\nu/\nu.
\label{m8:s:s3:eq:7-6}
\end{equation}
For \(Z\) in \eqref{m8:s:s3:eq:7-2}, put \(H=Zg\), which is traceless. Under the trace pairing on \(\mathfrak{gso}_8\), the corresponding ambient conormal element is
\begin{equation}
\mathcal H_Z=\frac12\operatorname{diag}
\left(H,-\mathcal J^{-1}H^t\mathcal J\right).
\label{m8:s:s3:eq:7-7}
\end{equation}
It pairs with \eqref{m8:s:s3:eq:7-6} as
\[
\operatorname{tr}(\mathcal H_Zm^{-1}\dot m)
=\operatorname{tr}(HY)-\tfrac\lambda2\operatorname{tr}H
=\operatorname{tr}(Z\dot g).
\]
It annihilates the two radical tangent directions because it is block diagonal. Thus \eqref{m8:s:s3:eq:7-7} is the five-dimensional ambient conormal of \eqref{m8:s:s3:eq:7-4} after trivialization. The similitude tangent direction forces the central cotangent component to vanish.

For the witness \eqref{m8:s:s3:eq:7-3}, this is the explicit matrix
\begin{equation}
\mathcal H_*=
\tfrac12\operatorname{diag}(1,-1,-1,1,-1,1,1,-1).
\label{m8:s:s3:eq:7-8}
\end{equation}
It lies in \(\mathfrak{so}(Q)\), satisfies \(\mathcal H_*^tQ+Q\mathcal H_*=0\), and has eigenvalues \(\pm1/2\). It is therefore nonnilpotent. This computation is valid in residue characteristic two because it takes place in the characteristic-zero field \(F\).

Since the normal forms meet every \(K\)-orbit, we can translate
the tangent and conormal spaces to any point of
\(\widehat{\mathcal N}\). The corresponding coadjoint
identifications preserve nilpotence and nonnilpotence, so the
argument applies throughout this locus.

\subparagraph{The conditional nilpotent-wavefront implication}\label{m8:s:s3:spectral}

Let \(T\) be a \(\widetilde K\)-equivariant distribution on \(G_8^\circ\), with the sign character, whose wavefront fibers on \(\mathcal M\) are nilpotent. Open restriction preserves these ambient wavefront fibers. 

The vanishing theorem for \(S_1\) implies that \(T|_{\mathcal M}\) is supported on the closed subset \(S_3\). Its corresponding intrinsic distribution on \(S_3\) has, by Section~\ref{m8:s:s3:support}, support in \(\widehat{\mathcal N}\). Hence
\(\operatorname{supp}(T|_{\mathcal M})\subset\widehat{\mathcal N}.\)
The common wavefront implication \eqref{m8:s:p:conormal-wf} and the ambient conormal in \eqref{m8:s:s3:eq:7-7}--\eqref{m8:s:s3:eq:7-8} contradict the nilpotent-wavefront condition at any point of nonzero support. 

Therefore
\begin{equation}
T|_{\mathcal M}=0.
\label{m8:s:s3:eq:7-9}
\end{equation}
This proves the conditional vanishing at the $S_3$ step for the
same ambient NW distribution on the open
\(\mathcal M\subset G_8^\circ\).

\par\hfill$\square$\par\smallskip

\subsubsection{\texorpdfstring{Closed orbits and spectral exclusion on $S_{11}$}{Closed orbits and spectral exclusion on S11}}\label{m8:s:sec:spectral-s11}
To isolate the spectral step, we first construct intrinsic odd
orbital distributions on the smooth closed locus $\mathcal Y$ of
Section~\ref{m8:s:m:smooth19}. We then compute its ambient
conormal and use it to prove the conditional NW support implication.
The orbit topology and uniform conormal ranks are established in
\S\ref{m8:s:m:conormal}; the rational flag normalizations cover all
$F$-points by Appendix, p.~\pageref{m8:s:m:witt}.

\paragraph{Closed orbits and nonzero odd orbital distributions}\label{m8:s:m:orbits}
\begin{proposition}[Odd orbital distributions]\label{m8:s:m:prop:orbits}
For every $(q,r)\in(F^\times)^2$, the orbit $\mathcal O_{q,r}$ defined below is closed in $S_{11}$ and has dimension $17$. The ordered pair $(q,r)$ determines the orbit, and transpose exchanges $\mathcal O_{q,r}$ with $\mathcal O_{r,q}$. Each orbit carries a nonzero $\eta_K$-equivariant ordinary orbital distribution. Consequently $D'_{\widetilde\eta_K}(S_{11})\ne0$.
\end{proposition}
\par\smallskip\noindent\textit{Proof.}\par

\subparagraph{The ordered invariants}\label{m8:s:m:ordered}
On $d=0$, the complete residual relation gives
\[
 \begin{pmatrix}p&q\\r&0\end{pmatrix}\mapsto
 \begin{pmatrix}p+sr+vq&q\\r&0\end{pmatrix}.
\]
Thus $q,r$ are individual invariants and $p$ is removed by
$s=-p/r$, $v=0$. Put
$E_{q,r}=(\begin{smallmatrix}0&q\\r&0\end{smallmatrix})$,
$x_{q,r}=G(E_{q,r})$, and $\mathcal O_{q,r}=Kx_{q,r}$.
By the product definition, $\mu(x_{q,r})=-qr$, and transpose
exchanges the ordered invariants.

\subparagraph{Closedness and quotient topology}\label{m8:s:m:closed}
The surjective map
\[
 \Phi:K\times\{(p,q,r):qr\ne0\}\to\mathcal Y,
 \qquad(k,p,q,r)\mapsto kG\!\left(\begin{pmatrix}p&q\\r&0\end{pmatrix}\right)
\]
has differential rank 19 by the uniform tangent calculation below.
Since $\mathcal Y$ is smooth of dimension 19, it is a surjective
analytic submersion, hence open and with local analytic sections.
By complete residual equivalence, $(k,p,q,r)\mapsto(q,r)$ is
constant on fibers and therefore descends to
\begin{equation}\label{m8:s:m:eq:ordered-map}
 \mathcal Y\longrightarrow(F^\times)^2,\qquad x\longmapsto(q,r).
\end{equation}
This map is continuous and analytic by the local sections, and is
a submersion by varying the canonical parameters. Its fibers are
exactly $\mathcal O_{q,r}$ and are smooth of dimension 17.
They are closed in $\mathcal Y$, hence in $S_{11}$.
The action differential has rank 17 by the same tangent calculation;
thus the continuous bijection $K/K_x\to\mathcal O_{q,r}$ is
locally open, and is a homeomorphism with the inherited topology.

\subparagraph{Full stabilizer and orbital distribution}\label{m8:s:m:orbit-stab}
The five-dimensional kernel is the whole $S_{11}$ row of
\eqref{m8:s:m:eq:stabilizer}, exponentiated and conjugated by
the canonical Levi factors. Its character is one, since $R=S=W$,
$\kappa=\lambda=1$. The remaining unipotent direction is the
full lift $k_{11}(qz,-rz;E_{q,r})$ of
\S\ref{m8:s:m:lifts}; the ineffective scalar is
$(1,1,cI_2)$ for every $c\in F^\times$. By the complete
residual-stabilizer calculation,
\[
 K_x\simeq\mathbb G_m(F)\times\mathcal U_{6,q,r}(F),\qquad
 \eta_K|_{K_x}=1,\qquad\Delta_K=\Delta_{K_x}=1.
\]
The expanded representative and nonlinear stabilizer lift are in
\texttt{anc/model8\_calculations.pdf}, calculation M8-C06,
under ``The ordered invariants'' and ``The stabilizer''.

\subparagraph{Ordinary orbital distributions}\label{m8:s:m:orbital-distributions}
The modular equality supplies an invariant quotient measure. Define
\begin{equation}\label{m8:s:m:eq:orbital}
 \mathcal T_{q,r}(f)=\int_{K/K_x}f(kx_{q,r})\eta_K(k)\,d\dot k.
\end{equation}
Since the character is trivial on $K_x$, it descends to the
quotient. For a compactly supported test function, its support
meets the closed orbit in a compact set, whose inverse image in
$K/K_x$ is compact by the homeomorphism above. Thus the integral
defines an ordinary distribution. A small compact open orbit neighborhood
on which the character is constant supplies a nonzero test
function; it extends from the closed orbit to $S_{11}$.
Changing variables $k=gk'$ gives
$\mathcal T_{q,r}(f\circ g^{-1})=\eta_K(g)\mathcal T_{q,r}(f)$,
and its support is the whole orbit. For $q\ne r$,
$\mathcal T_{q,r}-\tau\mathcal T_{q,r}$ is nonzero and odd,
since its two supports are disjoint closed orbits. For example,
$q=1,r=2$ works over every field in scope.

\subparagraph{The residual distribution model}\label{m8:s:m:small-model}
On $d=0$ the left $N_+$ translates $p/r$, while the right
$N_-$ also preserves the ordered $q,r$. Haar uniqueness therefore
gives precisely $dp\otimes S(q,r)$, with arbitrary
$S\in D'((F^\times)^2)$ (the nonzero $|r|$ factor is absorbed
in $S$). Transpose fixes $p$ and exchanges $q,r$;
$dp\otimes(\delta_{(1,2)}-\delta_{(2,1)})$ is the corresponding
odd residual example. This proves intrinsic nonvanishing on $S_{11}$.
\par\hfill$\square$\par\smallskip

\paragraph{\texorpdfstring{The tangent space and conormal of $\mathcal Y$}{The tangent space and conormal of Y}}\label{m8:s:m:conormal}
\begin{proposition}[Complete conormal calculation]\label{m8:s:m:prop:conormal}
In the trivialization $\dot x=xA_8$ and the trace-pairing identification of the Lie-algebra dual, the conormal of $\mathcal Y$ at $G(\begin{psmallmatrix}p&q\\r&0\end{psmallmatrix})$, $qr\ne0$, is the ten-dimensional space \eqref{m8:s:m:eq:conormal}. The tangent space has dimension $19$ for every such $p,q,r$, including $q=r$. Every point of $\mathcal Y$ has a nonnilpotent conormal direction.
\end{proposition}
\par\smallskip\noindent\textit{Proof.}\par

We compute the annihilator in the orthogonal Lie algebra, with the trace pairing; this retains the similitude direction until its annihilation is established.

\subparagraph{Tangent directions}\label{m8:s:m:tangent}

At
\[
E=\begin{pmatrix}p&q\\r&0\end{pmatrix},\qquad qr\ne0,
\qquad x=G(E)=VD(E),
\]
the differential of the action-and-parameter map is
\begin{equation}
{
Nx+xN'^t+D(A)x-xD(A)+VD(\dot E),
\quad N,N'\in\mathfrak u,\ A\in\mathfrak{gl}_2,
\quad \dot E=\begin{pmatrix}\dot p&\dot q\\\dot r&0\end{pmatrix}.
}
\label{m8:s:m:eq:tangent}
\end{equation}
This expression combines the 10 left $U$ directions, the 10 right
$U$ directions, all 4 diagonal-$\mathrm{GL}_2$ directions and all
3 residual directions.

Use the explicit trivialization
\[
\dot x=xA_8,\qquad A_8=x^{-1}\dot x.
\]
Then the generated tangent subspace is
\begin{equation}
{
\operatorname{Ad}(x^{-1})\mathfrak u+\mathfrak u^t
+\{x^{-1}D(A)x-D(A):A\in\mathfrak{gl}_2\}
+\{D(E^{-1}\dot E):\dot E_{22}=0\}.
}
\label{m8:s:m:eq:trivialized-tangent}
\end{equation}
We show below that its rank is 19. Since $\mathcal Y$ is smooth of
dimension 19, \eqref{m8:s:m:eq:tangent} therefore gives its full
tangent space.

\subparagraph{The orthogonal Lie algebra and annihilator equations}\label{m8:s:m:annihilator}

The trace pairing on \(\mathfrak{gso}_8\) is nondegenerate. One independent basis is
\[
Q(E_{ij}-E_{ji}),\quad 1\le i<j\le8,
\qquad I_8.
\]
Of these vectors, the first 28 satisfy \(B^tQ+QB=0\), while
$I_8$ supplies the similitude direction. In this basis the Gram
matrix of the trace pairing has rank 29.

The annihilator equations for a general $B$ in this space can be written without choosing coordinates:
\begin{equation}
\begin{aligned}
&\operatorname{tr}(BN'^t)=0 &&\text{for every }N'\in\mathfrak u,\\
&\operatorname{tr}(xBx^{-1}N)=0 &&\text{for every }N\in\mathfrak u,\\
&\operatorname{tr}\bigl((xBx^{-1}-B)D(A)\bigr)=0
 &&\text{for every }A\in\mathfrak{gl}_2,\\
&\operatorname{tr}\bigl(BD(E^{-1}\dot E)\bigr)=0
 &&\text{for }\dot E_{22}=0.
\end{aligned}
\label{m8:s:m:eq:annihilator}
\end{equation}
These are exactly the four summands of \eqref{m8:s:m:eq:trivialized-tangent}, paired with $B$. Notice that \(\dot E=E\) is allowed at $d = 0$; it contributes the scalar tangent \(I_8\), so every conormal $B$ has trace zero. Hence the conormal lies in the orthogonal Lie-algebra dual.

\subparagraph{The ten-dimensional conormal}\label{m8:s:m:conormal-matrix}

The complete solution at every $p,q,r$ with $qr\ne0$ is
\begin{equation}\label{m8:s:m:eq:conormal}
 N^*_{\mathcal Y,x}
 =\operatorname{span}_F\{C_1,C_2,C_3,C_4,C_5,C_6,C_9,C_{10},H_1,H_2\},
\end{equation}
where
\[
\begin{aligned}
C_1&=E_{31}-E_{86},&C_2&=E_{41}+E_{85},\\
C_3&=E_{42}-E_{75},&C_4&=E_{61}-E_{83},\\
C_5&=E_{51}+E_{62}+E_{73}+E_{84},&C_6&=E_{71}+E_{82},\\
C_9&=E_{21}+E_{87}-E_{51}-E_{84},&
C_{10}&=E_{43}+E_{65}-E_{51}-E_{84},\\
H_1&=\operatorname{diag}(-1,0,1,0,0,-1,0,1),&
H_2&=\operatorname{diag}(0,-1,0,1,-1,0,1,0).
\end{aligned}
\]
These ten independent vectors lie in $\mathfrak{so}(Q)$ and
satisfy the four annihilator equations at every canonical boundary
point. To prove completeness we use the uniform minor below. The expanded $8\times8$
parameter matrix and entrywise substitution are in
\texttt{anc/model8\_calculations.pdf}, calculation M8-C06,
under ``The ten-dimensional conormal''.

\subparagraph{A uniform rank calculation}\label{m8:s:m:rank-certificate}

Let the Lie-algebra basis order be:

\begin{enumerate}
\item the ten $U$-coordinate roots \(X_{11},X_{12},X_{21},X_{22},Y,Z_{11},Z_{12},Z_{21},Z_{22},Z_0\);
\item their ten transposes;
\item the five Levi torus directions $A_{11},A_{22},B_{11},B_{22}$ and the similitude direction;
\item the four Levi roots \(A_+,A_-,B_+,B_-\).
\end{enumerate}
Let the tangent generator order be the ten left-$U$ directions, ten right-$U$ directions, the four diagonal $\mathrm{GL}_2$ directions $E_{11},E_{12},E_{21},E_{22}$, and the three residual directions $p,q,r$. Define the $27\times29$ pairing matrix by
\[
\mathsf P_{jk}=\operatorname{tr}(B_kx^{-1}T_j).
\]
With zero-based indices, select rows
\[
(0,1,2,3,4,5,6,7,8,9,12,14,15,17,18,20,21,25,26)
\]
and columns
\[
(0,1,2,3,4,5,6,7,8,9,12,14,15,17,20,21,24,25,27).
\]
The determinant of this fixed $19\times19$ minor, computed without specializing $p$, is
\begin{equation}
{-524288\frac{r^5}{q^7}=-2^{19}\frac{r^5}{q^7}.}
\label{m8:s:m:eq:minor}
\end{equation}
Since this determinant is nonzero for every \(qr\ne0\), the
pairing rank is at least 19. The ten independent matrices in
\eqref{m8:s:m:eq:conormal} annihilate all 27 tangent generators,
giving the upper bound \(29-10=19\). Thus the rank is exactly
19 and the displayed conormal is complete.

Specializing $r = q$ leaves rank 19. The determinant remains nonzero in residue characteristic $2$, since $F$ has characteristic zero.

As a separate dimension check, the first 24 action generators have rank 17 on this boundary. Their kernel has dimension 7, agreeing with the $S_{11}$ stabilizer in Section~\ref{m8:s:m:K-stab}. Adding the three residual directions increases rank by two, because the $p$-direction is already an orbit direction. Thus
\[
\dim T_x\mathcal Y=17+2=19.
\]

\subparagraph{\texorpdfstring{Removal of $p$ and conormal transport}{Removal of p and conormal transport}}\label{m8:s:m:transport}

The parameter $p$ can also be removed by the following identity:
\[
u(0,0,-(p/r)H,0)\,G\!\left(\begin{pmatrix}p&q\\r&0\end{pmatrix}\right)
 =G\!\left(\begin{pmatrix}0&q\\r&0\end{pmatrix}\right).
\]
It is a pure left translation. Under the trivialization \(x^{-1}\dot x\), its cotangent matrix is unchanged.

More generally, for a fixed $K$-action written \(x'=LxR\),
\[
x'^{-1}\dot x'=R^{-1}(x^{-1}\dot x)R.
\]
The corresponding conormal matrix transforms by
\begin{equation}
{B\longmapsto R^{-1}BR.}
\label{m8:s:m:eq:transport}
\end{equation}
This transports the complete conormal to every point of $\mathcal Y$ and preserves nilpotence and nonnilpotence.

\subparagraph{A nonnilpotent conormal vector}\label{m8:s:m:witness}

The displayed \(H_1\) lies in the conormal and satisfies \(H_1^tQ+QH_1=0\). It is diagonal with nonzero eigenvalues $\pm1$; in particular
\[
\operatorname{tr}(H_1^2)=4\ne0.
\]
It is not nilpotent. Every nonzero linear combination of $H_1$ and $H_2$ is also a nonzero diagonal matrix and therefore nonnilpotent. By \eqref{m8:s:m:eq:transport}, each point of $\mathcal Y$ has a nonnilpotent conormal direction in the orthogonal Lie-algebra dual.

\par\hfill$\square$\par\smallskip

\paragraph{Spectral exclusion after clearing the predecessor open}\label{m8:s:m:spectral}
\begin{proposition}[The conditional NW implication]\label{m8:s:m:prop:spectral}
Let $T\in\Do(G_8^\circ)\cap\NW(G_8^\circ)$ be the same ambient distribution throughout. With the opens defined below,
\[
T|_{\Omega_{\rm above}}=0\quad\Longrightarrow\quad T|_{\Omega_{11}}=0.
\]
The NW hypothesis belongs to the ambient distribution on $G_8^\circ$, including its zero central cotangent component.
\end{proposition}
\par\smallskip\noindent\textit{Proof.}\par

\subparagraph{The predecessor opens}\label{m8:s:m:opens}

For a selected rank matrix define
\[
\Omega_i=\{x\in G_8^\circ:R(x)\ge R_i\text{ entrywise}\}.
\]
Every condition rank $\ge n$ is Zariski open, so these are analytic open submanifolds of $G_8^\circ$. The rank table gives
\[
\Omega_{11}=
\mathcal O_{\rm cleared}\sqcup S_6\sqcup S_7\sqcup S_{10}\sqcup S_{11}.
\]
Define
\[
\Omega_{\rm above}=
\Omega_{11}\cap\bigl(\{\operatorname{rank}x_{44}\ge1\}
 \cup\{r_{4,4}=4\}\bigr)
 =\mathcal O_{\rm cleared}\sqcup S_6\sqcup S_7,
\]
\[
\Omega_{\rm pre}=\Omega_{11}\setminus S_{11}
 =\Omega_{\rm above}\sqcup S_{10}.
\]
$S_{11}$ is closed in \(\Omega_{11}\); $S_{10}$ is closed in \(\Omega_{\rm pre}\). All three opens are $K$- and transpose-stable.

Retain the original NW odd equivariant distribution $T$. The predecessor hypothesis is
\begin{equation}
{T|_{\Omega_{\rm above}}=0.}
\label{m8:s:m:eq:premise}
\end{equation}
In Section~\ref{m8:s:sec:assembly} we establish this premise from
the open-stratum results and the equality
$\Omega_{\rm above}=\mathcal O_7$. 

\subparagraph{Trapping the support}\label{m8:s:m:trapping}

Restricting $T$ to the open \(\Omega_{\rm pre}\), we apply
\[
0\longrightarrow D'_{\widetilde\eta_K}(S_{10})
\longrightarrow D'_{\widetilde\eta_K}(\Omega_{\rm pre})
\longrightarrow D'_{\widetilde\eta_K}(\Omega_{\rm above}).
\]
By Section~\ref{m8:s:m:S10} the left term vanishes, while
\eqref{m8:s:m:eq:premise} makes the restriction of $T$ in the
right term zero. Hence \(T|_{\Omega_{\rm pre}}=0\).

Next use
\[
0\longrightarrow D'_{\widetilde\eta_K}(S_{11})
\longrightarrow D'_{\widetilde\eta_K}(\Omega_{11})
\longrightarrow D'_{\widetilde\eta_K}(\Omega_{\rm pre}).
\]
It follows that \(T|_{\Omega_{11}}\) is the closed-embedding pushforward of an ordinary equivariant distribution on $S_{11}$. Apply \eqref{m8:s:m:eq:Y-sequence} to that induced distribution. We obtain
\begin{equation}
{\operatorname{supp}(T|_{\Omega_{11}})\subset\mathcal Y.}
\label{m8:s:m:eq:trapping}
\end{equation}
$\mathcal Y$ is closed in $S_{11}$ and hence in \(\Omega_{11}\), and it is a smooth 19-dimensional submanifold. It includes the locus $q=r$.

We now apply the conormal argument to the original ambient NW
distribution $T$.

\stepcounter{table}
\subparagraph{The common wavefront input}\label{m8:s:m:wavefront-inputs}
Nilpotent wavefront is local on the open $\Omega_{11}$.
\subparagraph{The zero cotangent vector convention}\label{m8:s:m:zero}
We use the full-cotangent convention of Definition~\ref{cf:native},
including zero over the support, exactly as in
Proposition~\ref{cf:conormal}.
\subparagraph{Exclusion at every support point}\label{m8:s:m:contradiction}
By \eqref{m8:s:m:eq:trapping}, the restriction is supported on the
smooth closed nineteen-dimensional $\mathcal Y$. The ambient
conormal contains a transport of $H_1$ at every point by
\S\ref{m8:s:m:witness}; its trace-square is nonzero. This includes
$q=r$. Proposition~\ref{cf:conormal} therefore gives
$T|_{\Omega_{11}}=0$.

We have proved exactly
\[
{
T|_{\mathcal O_{\rm cleared}\sqcup S_6\sqcup S_7}=0
\quad\Longrightarrow\quad
T|_{\mathcal O_{\rm cleared}\sqcup S_6\sqcup S_7\sqcup S_{10}\sqcup S_{11}}=0.
}
\]
Both displayed unions are the opens above. The premise is discharged in Section~\ref{m8:s:sec:assembly}. Thus the proposition applies there to the same original global
NW distribution after clearing its predecessor open.

\par\hfill$\square$\par\smallskip

\subsection{The open induction}\label{m8:filtration}
We verify the rank opens and predecessor containments for the
support induction described at the start of this section.
\subsubsection{Global support assembly}\label{m8:s:sec:assembly}
\paragraph{The NW fiber theorem}
The precise assertion on the identity fiber is
\begin{equation}\label{m8:s:a:native-class}
 \Do(G_8^\circ)\cap\NW(G_8^\circ)=0.
\end{equation}
\begin{theorem}\label{m8:s:a:global}
Every NW, odd, $\eta_K$-equivariant ordinary distribution on
$G_8^\circ$ is zero.
\end{theorem}
Fix such a $T$; throughout the proof, all open restrictions and
induced closed-support distributions come from this same ambient
distribution. Once we prove its vanishing, the fiber assertion of
Theorem~\ref{m8:distribution} follows. Wavefront transfer by
\eqref{m8:s:f:native-transfer} and injective Frobenius reconstruction
then give the assertion for $G_8^\circ\times\GL_2(F)$.

\paragraph{The initial open: one support assembly}\label{m8:s:a:entrance-section}
Keep the same original NW distribution $T$. The rank table gives
\begin{align}
 \Omega_2&=\{\det x_{44}\ne0\}=S_1\sqcup S_2,\label{m8:s:s1:omega2}\\
 \Omega_4&=\{r_{4,4}=4\}=S_1\sqcup S_3\sqcup S_8,\qquad
 \mathcal M=\Omega_4\setminus S_8=S_1\sqcup S_3.\label{m8:s:s1:omega4}
\end{align}
These are open subsets of $G_8^\circ$; for $\mathcal M$, the additional condition is $x_{44}\ne0$. In each case the displayed terminal stratum is closed in the indicated open. Intrinsic vanishing on $S_1,S_2$ and Lemma~\ref{cf:support} give $\Do(\Omega_2)=0$. On $\mathcal M$, the $S_1$ result gives a closed-support element on $S_3$; the ambient conormal calculation in Proposition~\ref{m8:s:s3:prop:spectral} kills this particular element from $T$. Then the closed $S_8$ step gives $T|_{\Omega_4}=0$.

For the transpose pair set
\begin{equation}\label{m8:s:s45:eq8.1}\begin{aligned}
 \mathcal O_4&=\{x_{44}\ne0,r_{2,4}=2\}=S_1\sqcup S_2\sqcup S_3\sqcup S_4,\\
 \mathcal O_5&=\{x_{44}\ne0,r_{4,2}=2\}=S_1\sqcup S_2\sqcup S_3\sqcup S_5,\\
 \mathcal O_4\cap\mathcal O_5&=S_1\sqcup S_2\sqcup S_3.
\end{aligned}
\end{equation}
These $K$-stable rank opens are exchanged by transpose, and their
closed differences are $S_4$ and $S_5$. Since the already cleared
$\Omega_2$ and $\mathcal M$ cover the common open complement,
we can apply the even support sequence on each $\mathcal O_i$
with Theorem~\ref{m8:s:s45:vanishing}. We obtain
\begin{equation}\label{m8:s:a:old}\begin{aligned}
 \mathcal O_{\rm old}&=\Omega_4\cup\mathcal O_4\cup\mathcal O_5\\
 &=\{r_{4,4}=4\}\cup\{r_{2,2}\ge1,\ r_{2,4}=2\text{ or }r_{4,2}=2\}\\
 &=S_1\sqcup S_2\sqcup S_3\sqcup S_4\sqcup S_5\sqcup S_8,
 \qquad T|_{\mathcal O_{\rm old}}=0.
\end{aligned}
\end{equation}
For the final two steps put
\begin{align}
 \mathcal U_6&=\{r_{2,2}\ge1,r_{2,6}=2\},&
 \mathcal O_6&=\mathcal U_6\cup\Omega_4=\mathcal O_{\rm old}\sqcup S_6,\label{m8:s:a:O6}\\
 \mathcal U_7&=\{r_{2,2}\ge1\},&
 \mathcal O_7&=\mathcal U_7\cup\Omega_4=\mathcal O_6\sqcup S_7=\bigcup_{i=1}^8S_i.\label{m8:s:a:O7}
\end{align}
Table~\ref{m8:s:r:table} gives the label sets $\mathcal U_6=\bigcup_{i=1}^6S_i$, $\mathcal U_7=\bigcup_{i=1}^7S_i$, and $\Omega_4=S_1\cup S_3\cup S_8$. The conditions require a nonzero entry of $x_{44}$, a maximal two-by-six minor, or a four-by-four determinant, respectively, so they are open. Lower-right ranks are $K$-invariant, and $r_{2,6}=r_{6,2}$ gives transpose stability. Since $\mathcal O_{\rm old}$ and $\mathcal O_6$ are open, $S_6$ is closed in $\mathcal O_6$ and $S_7$ in $\mathcal O_7$. Also $S_6\sqcup S_7$ is closed in $\mathcal O_7$, with $S_6$ open in that closed subset.

The following table records the applications of
Lemma~\ref{cf:support} in the initial step. In each left-exact
sequence, the right term is the restriction of the same $T$.
\begin{center}
\begin{tabular}{llll}
\toprule
Ambient $X$&Closed $Z$&Open $X\setminus Z$&Input on $Z$\\
\midrule
$\Omega_2$&$S_2$&$S_1$&intrinsic odd vanishing\\
$\mathcal M$&$S_3$&$S_1$&spectral exclusion for $T$\\
$\Omega_4$&$S_8$&$\mathcal M$&intrinsic odd vanishing\\
$\mathcal O_4$&$S_4$&$\mathcal O_4\cap\mathcal O_5$&intrinsic even vanishing\\
$\mathcal O_5$&$S_5$&$\mathcal O_4\cap\mathcal O_5$&intrinsic even vanishing\\
$\mathcal O_6$&$S_6$&$\mathcal O_{\rm old}$&intrinsic odd vanishing\\
$\mathcal O_7$&$S_7$&$\mathcal O_6$&intrinsic odd vanishing\\
\bottomrule
\end{tabular}
\end{center}
The last two intrinsic inputs are Proposition~\ref{m8:s:v:vanishing}. Consequently
\begin{equation}\label{m8:s:a:entrance}
 T|_{\mathcal O_6}=T|_{\mathcal O_7}=0,\qquad
 \mathcal O_7=\{x:x_{44}\ne0\}\cup\{x:r_{4,4}=4\}.
\end{equation}
Here $x_{44}\ne0$ means that at least one matrix entry is nonzero.
The $S_3$ step uses the ambient NW hypothesis. All restrictions above are to opens, and all stratum distributions enter by closed pushforward.

\paragraph{Named rank opens and the predecessor equality}
For $i=10,11,12,15,16$ use the rank envelopes $\Omega_i=\{R(x)\ge R_i\}$. Define further
\begin{equation}\label{m8:s:a:named-opens}
 \mathcal O_L=\{e=2\},\quad\mathcal O_R=\{d=2\},\quad
 \Omega_b=\{b=2\},\quad\Omega_a=\{a\ge2\},\quad
 \Omega_{\mathrm{pre}}=\Omega_{11}\setminus S_{11}.
\end{equation}
Since $S_{11}$ is closed in $\Omega_{11}$, the last set is open as well. In the following table $[1,m]$ means the set of labels $1,\ldots,m$.
\begin{center}
\small
\begin{tabular}{lll}
\toprule
Open & Exact labels & Transpose\\
\midrule
$\mathcal O_7$ & $[1,8]$ & stable\\
$\mathcal O_L$ & $1,2,3,4,8,9$ & $\mathcal O_R$\\
$\mathcal O_R$ & $1,2,3,5,8,13$ & $\mathcal O_L$\\
$\Omega_{10}$ & $1,2,3,4,5,6,8,10$ & stable\\
$\Omega_{\mathrm{pre}}$ & $[1,8],10$ & stable\\
$\Omega_{11}$ & $[1,8],10,11$ & stable\\
$\Omega_{12}$ & $[1,12]$ & $\Omega_{15}$\\
$\Omega_{15}$ & $[1,8],10,11,13,15$ & $\Omega_{12}$\\
$\Omega_b$ & $1,2,3,4,5,6,8,9,10,13,14$ & stable\\
$\Omega_{16}$ & $[1,16]$ & stable\\
$\Omega_a$ & $[1,13],15,17$ & stable\\
\bottomrule
\end{tabular}
\end{center}
Substituting the rows of Table~\ref{m8:s:r:table} into the defining
conditions gives these identities. Since the table exhausts the
rank map of $G_8^\circ$, they hold as equalities of subsets. In particular
\begin{equation}\label{m8:s:a:pre-explicit}
 \Omega_{\mathrm{pre}}=\mathcal O_7\cup\Omega_{10}
                     =\mathcal O_7\sqcup S_{10}.
\end{equation}

The predecessor open used for $S_{11}$ in \S\ref{m8:s:sec:spectral-s11} is
\[
 \Omega_{\mathrm{above}}
 =\Omega_{11}\cap\bigl(\{\rank x_{44}\ge1\}\cup\{r_{4,4}=4\}\bigr).
\]
For a matrix $x$, $x_{44}\ne0$ is equivalent to $\rank x_{44}\ge1$. Moreover every rank matrix with label $1,\ldots,8$ dominates $R_{11}$, so $\mathcal O_7\subset\Omega_{11}$. It follows that
\begin{equation}\label{m8:s:a:predecessor-equality}
 \Omega_{\mathrm{above}}=\Omega_{11}\cap\mathcal O_7
                       =\mathcal O_7.
\end{equation}
In particular \eqref{m8:s:a:entrance} supplies the predecessor premise
for the original $T$, including the strata $S_4,S_5,S_6,S_7$.

\paragraph{The cumulative filtration}
Define
\begin{equation}\label{m8:s:a:filtration-definition}
\begin{aligned}
 X_0&=\mathcal O_7,&X_1&=X_0\cup\mathcal O_L,&X_2&=X_1\cup\mathcal O_R,\\
 X_3&=X_2\cup\Omega_{\mathrm{pre}},&X_4&=X_3\cup\Omega_{11},&X_5&=X_4\cup\Omega_{12},\\
 X_6&=X_5\cup\Omega_{15},&X_7&=X_6\cup\Omega_b,&X_8&=X_7\cup\Omega_{16},\\
 X_9&=X_8\cup\Omega_a,&X_{10}&=G_8^\circ.
\end{aligned}
\end{equation}
Each $X_j$ is open in $G_8^\circ$; its strata and the group used
in the support sequence are listed below.
\begin{center}
\small
\begin{tabular}{clccl}
\toprule
$j$ & Labels in $X_j$ & $X_j\setminus X_{j-1}$ & Group & $X_j$ transpose-stable\\
\midrule
0 &$[1,8]$&---&$\widetilde K$&yes\\
1 &$[1,9]$&$S_9$&$K$&no\\
2 &$[1,9],13$&$S_{13}$&$K$&yes\\
3 &$[1,10],13$&$S_{10}$&$\widetilde K$&yes\\
4 &$[1,11],13$&$S_{11}$&$\widetilde K$&yes\\
5 &$[1,13]$&$S_{12}$&$K$&no\\
6 &$[1,13],15$&$S_{15}$&$K$&yes\\
7 &$[1,15]$&$S_{14}$&$K$&yes\\
8 &$[1,16]$&$S_{16}$&$K$&yes\\
9 &$[1,17]$&$S_{17}$&$K$&yes\\
10&$[1,18]$&$S_{18}$&$\widetilde K$&yes\\
\bottomrule
\end{tabular}
\end{center}
Each new stratum is closed in the current open, since its complement is the preceding open. Nonemptiness of every stratum makes the inclusions strict. The order is
\begin{equation}\label{m8:s:a:order}
 X_0\subset X_1\subset\cdots\subset X_{10}=G_8^\circ,
 \qquad 9,13,10,11,12,15,14,16,17,18.
\end{equation}
For the second support sequence we use $K$, because neither the
preceding open $X_1$ nor the new stratum $S_{13}$ is transpose-stable,
even though $X_2$ itself is stable. The same choice is required
at $X_6$ for its preceding open $X_5$ and new stratum $S_{15}$.

\paragraph{A common argument for the intrinsic support steps}
For $j=1,2,3,5,6,7,8,9$, let $S_{i_j}=X_j\setminus X_{j-1}$,
and let $E_j$ and the equivariance choice $\mathcal D_j$ be given
in the following table. Here $\mathcal D_j=\De$ for $K$ and
$\mathcal D_j=\Do$ for $\widetilde K$.
\begin{center}\small
\begin{tabular}{cclll}
\toprule
$j$&$i_j$&Group&$E_j$&Labels in $E_j\cap X_{j-1}$\\
\midrule
1&9&$K$&$\mathcal O_L$&$1,2,3,4,8$\\
2&13&$K$&$\mathcal O_R$&$1,2,3,5,8$\\
3&10&$\widetilde K$&$\Omega_{\mathrm{pre}}$&$[1,8]$\\
5&12&$K$&$\Omega_{12}$&$[1,11]$\\
6&15&$K$&$\Omega_{15}$&$[1,8],10,11,13$\\
7&14&$K$&$\Omega_b$&$1,2,3,4,5,6,8,9,10,13$\\
8&16&$K$&$\Omega_{16}$&$[1,15]$\\
9&17&$K$&$\Omega_a$&$[1,13],15$\\
\bottomrule
\end{tabular}\end{center}
The named-open table and \eqref{m8:s:a:filtration-definition}
show that $X_j=X_{j-1}\cup E_j$, that $S_{i_j}$ is closed in both
$X_j$ and $E_j$, and that $E_j\setminus S_{i_j}=E_j\cap X_{j-1}$.
Every space used with $\Do$ is transpose-stable. The two instances
of Lemma~\ref{cf:support} are therefore
\begin{equation}\label{m8:s:a:sequence1}
 0\longrightarrow\mathcal D_j(S_{i_j})
 \longrightarrow\mathcal D_j(X_j)
 \longrightarrow\mathcal D_j(X_{j-1})
\end{equation}
and
\begin{equation}\label{m8:s:a:local1}
 0\longrightarrow\mathcal D_j(S_{i_j})
 \longrightarrow\mathcal D_j(E_j)
 \longrightarrow\mathcal D_j(E_j\cap X_{j-1}).
\end{equation}
In each listed row the left term vanishes: the terminal theorem
handles $9,13,14,17$, the middle even theorem handles $12,15,16$,
and the middle odd theorem handles $10$. Hence $T|_{X_{j-1}}=0$ implies
$T|_{X_j}=0$. Applying this implication first for $j=1,2,3$ gives
$T|_{X_3}=0$ from \eqref{m8:s:a:entrance}. The $j=3$ local
complement is precisely $\Omega_{\mathrm{above}}=\mathcal O_7$,
so this also clears the predecessor required for $S_{11}$.
The $K$ choice for the transpose-exchanged steps is essential,
even when the current cumulative open happens to be stable.

\paragraph{The nonzero kernel at \texorpdfstring{$S_{11}$}{S11}}
Now $X_4=X_3\sqcup S_{11}$. The two support sequences are
\begin{align}
0&\longrightarrow\Do(S_{11})\xrightarrow{i_{4,*}}\Do(X_4)
  \longrightarrow\Do(X_3),\label{m8:s:a:sequence4}\\
0&\longrightarrow\Do(S_{11})\xrightarrow{i_{11,*}}\Do(\Omega_{11})
  \longrightarrow\Do(\Omega_{\mathrm{pre}}).\label{m8:s:a:local4}
\end{align}
By the closed-orbit construction in
\S\ref{m8:s:sec:spectral-s11}, the left space is nonzero. Since
$T|_{X_3}=0$, left exactness nevertheless gives a unique element
\begin{equation}\label{m8:s:a:R11}
 R_{11}\in\Do(S_{11}),\qquad T|_{X_4}=i_{4,*}R_{11}.
\end{equation}
The element $R_{11}$ is defined by this closed-support pushforward.

The original global $T$ vanishes on $\Omega_{\mathrm{above}}=\mathcal O_7$ by \eqref{m8:s:a:entrance}. The conditional spectral theorem of \S\ref{m8:s:sec:spectral-s11} therefore applies and gives
\begin{equation}\label{m8:s:a:S11-spectral}
 T|_{\Omega_{11}}=0.
\end{equation}
Here the nilpotent wavefront of $T$ persists on the open subset.
The conditional theorem confines the induced support to the smooth
nineteen-dimensional locus and excludes it using its nonnilpotent conormal.

The open $\Omega_{11}\subset X_4$ contains all of $S_{11}$. Restricting \eqref{m8:s:a:R11} to this open gives
\[
 0=T|_{\Omega_{11}}=i_{11,*}R_{11}.
\]
Indeed, extension by zero of a test function from $\Omega_{11}$ to $X_4$ has the same restriction to $S_{11}$. The map $i_{11,*}$ is injective by Lemma~\ref{cf:support}: for any $\varphi\in\cD(S_{11})$, choose an extension $f\in\cD(\Omega_{11})$ and obtain $R_{11}(\varphi)=i_{11,*}R_{11}(f)=0$. Thus $R_{11}=0$, and
\begin{equation}\label{m8:s:a:X4-zero}
 T|_{X_4}=0.
\end{equation}
Thus, after clearing the predecessor, we eliminate the particular
closed-support element induced by the original NW distribution.

\paragraph{Continuation after the spectral step}
Equation~\eqref{m8:s:a:X4-zero} starts the same induction for
$j=5,6,7,8,9$ in the table above. It gives $T|_{X_9}=0$.
The local complements are all contained in the preceding cleared
open; in particular $S_6$ and $S_{10}$ precede $S_{14}$, while
$S_{12},S_{14},S_{15}$ precede $S_{16}$. These five kernels vanish
already under $U\times U$ as ordinary equivariant distribution spaces.

\paragraph{The final closed stratum}
The rank table gives
\begin{equation}\label{m8:s:a:S18-closed}
 S_{18}=\{r_{6,2}=0,\ r_{4,4}=0\}.
\end{equation}
These are closed polynomial conditions. Also $\Omega_{16}=\{r_{6,2}\ge1\}$, and the possible values of $r_{4,4}$ are $0,2,4$, so
\begin{equation}\label{m8:s:a:final-complement}
 \begin{aligned}
 G_8^\circ\setminus S_{18}
 &=\{r_{6,2}\ge1\}\cup\{r_{4,4}\ge1\}\\
 &=\Omega_{16}\cup\{r_{4,4}\ge2\}
 =\Omega_{16}\cup\Omega_a=X_9.
 \end{aligned}
\end{equation}
Thus the final step is the open/closed decomposition $G_8^\circ=X_9\sqcup S_{18}$. It has the odd support sequence
\begin{equation}\label{m8:s:a:sequence10}
 0\longrightarrow\Do(S_{18})\longrightarrow\Do(G_8^\circ)
       \longrightarrow\Do(X_9).
\end{equation}
We have already shown $T|_{X_9}=0$, and intrinsic odd vanishing
on $S_{18}$ makes the left term zero. Hence $T=0$, proving
Theorem~\ref{m8:s:a:global}.

Every local complement in the support-step table and
\eqref{m8:s:a:local4} is the intersection of its envelope with the
preceding cumulative open. Restrictions are to open subsets of
$G_8^\circ$, and support elements are defined by injective closed
pushforward. The two spectral steps use the same ambient NW
distribution throughout.

\subsection{Connected inverse-transpose duality}\label{m8:duality}
To turn the GK product bound into multiplicity one, we prove
inverse-transpose duality on the connected group. A rational
cyclic-vector construction gives a determinant-one conjugator,
to which we apply the character argument below.
\subsubsection{Inverse-transpose duality for connected \texorpdfstring{$G_8^\circ$}{G8}}\label{m8:s:sec:duality}
To pass from the product bound supplied by the distribution theorem
to an individual bound for every irreducible admissible representation,
we identify the opposite period space through
$\pi^\alpha\simeq\pi^\vee$.
For the connected orthogonal similitude group, we verify the required
rational conjugacy within its identity component and then apply
Harish--Chandra character theory.

\paragraph{The connected component and a full-measure open}
The algebraic identities
\begin{equation}\label{m8:s:d:inverse-transpose}
 xQx^t=\mu Q,\qquad x^{-t}=\mu^{-1}QxQ
\end{equation}
follow from $x^tQx=\mu Q$ and $Q^2=I_8$. Inverse transpose has multiplier $\mu^{-1}$ and determinant $\mu^{-4}$, so it preserves the connected condition defining $G_8^\circ$. In ordinary coordinates $Q$ is antidiagonal with four paired signed reversals. Each paired two-dimensional coordinate subspace contributes determinant $-1$, whence
\begin{equation}\label{m8:s:d:Q-component}
 \det Q=(-1)^4=1,\qquad Q^tQQ=Q,\qquad Q\in\SO(Q)(F).
\end{equation}

Consider the subset
\begin{equation}\label{m8:s:d:very-regular}
 G_8^{\circ,\mathrm{vr}}=
 \{x:\operatorname{disc}(\operatorname{charpol}_x)\ne0,
       \ \det(x^2-\mu(x)I_8)\ne0\}.
\end{equation}
It is a nonempty Zariski-open regular semisimple subset. For example, if $\varpi$ is a uniformizer, the diagonal matrix with entries
\begin{equation}\label{m8:s:d:diagonal-example}
 \varpi,\varpi^2,\varpi^3,\varpi^4,
 \varpi^{-4},\varpi^{-3},\varpi^{-2},\varpi^{-1}
\end{equation}
preserves the antidiagonal pairing, has determinant and multiplier one, has eight distinct eigenvalues, and has no eigenvalue with square one.

The complement of \eqref{m8:s:d:very-regular} has Haar measure zero. Each excluded condition is a proper algebraic subset of the smooth geometrically connected group. In analytic coordinate charts its defining nonzero regular function has zero set of measure zero. This follows from the zero-set property for nonzero analytic functions over a characteristic-zero local field: in one variable a nonzero analytic function cannot vanish on a ball; induction on the number of variables and Fubini apply after expanding in one variable and separating the locus on which all its coefficient functions vanish. No analytic open subset of a smooth irreducible algebraic group is contained in a proper algebraic subset. A countable cover by analytic charts gives the asserted measure-zero complement. Thus the open used for the character argument is full measure in $G_8^\circ(F)$.

\paragraph{The etale algebra and its adjoint involution}
Fix $x\in G_8^{\circ,\mathrm{vr}}$ and write $V_0=F^8$. Its characteristic polynomial has eight distinct roots. Since every eigenvalue is a root of the minimal polynomial, the
minimal and characteristic polynomials agree. Consequently
\begin{equation}\label{m8:s:d:etale}
 A=F[x]
\end{equation}
is an etale $F$-algebra of dimension eight.

The $A$-module $V_0$ is free of rank one over $A$, with a cyclic vector defined over $F$. To prove this rationally, write $A=\prod_i E_i$ with $E_i/F$ finite separable fields, and decompose $V_0=\bigoplus_i V_i$ by the idempotents. Faithfulness of the action implies $V_i\ne0$ for every $i$. On the other hand,
\[
 \dim_F V_0=8=\dim_F A=\sum_i[E_i:F].
\]
Every $V_i$ must therefore have dimension one over $E_i$. Choosing a nonzero vector in each $V_i$ gives a cyclic vector $v_0\in V_0$ and the $F$-linear isomorphism
\begin{equation}\label{m8:s:d:cyclic}
 A\xrightarrow{\ \sim\ }V_0,\qquad a\longmapsto av_0.
\end{equation}
In particular, this construction is defined over $F$ for nonsplit
etale algebras as well.

For the symmetric pairing $\langle v,w\rangle=v^tQw$, the adjoint of $x$ is
\begin{equation}\label{m8:s:d:adjoint}
 x^*=Q^{-1}x^tQ=\mu x^{-1}\in A.
\end{equation}
Since $x$ is invertible, $x^{-1}$ is a polynomial in $x$. Thus adjoint preserves $A$; on this commutative algebra its anti-involution is an $F$-algebra involution. Define
\begin{equation}\label{m8:s:d:h}
 h(av_0)=a^*v_0.
\end{equation}
It is an $F$-linear involution. Put $\ell(a)=\langle v_0,av_0\rangle$. Symmetry of the form gives $\ell(a^*)=\ell(a)$, and
\[
 \langle av_0,bv_0\rangle=\ell(a^*b).
\]
Hence
\begin{equation}\label{m8:s:d:isometry}
 \begin{aligned}
 \langle h(av_0),h(bv_0)\rangle
 &=\ell(ab^*)=\ell((ab^*)^*)\\
 &=\ell(a^*b)=\langle av_0,bv_0\rangle.
 \end{aligned}
\end{equation}
Thus $h$ is a rational isometry. Acting on an arbitrary vector $av_0$,
\[
 hxh^{-1}(av_0)=(xa^*)^*v_0=x^*av_0,
\]
so
\begin{equation}\label{m8:s:d:adjoint-conjugacy}
 hxh^{-1}=\mu x^{-1}.
\end{equation}

\paragraph{The determinant-one calculation}
After extending scalars to an algebraic closure, $A$ is the product of eight primitive-idempotent lines. On these lines the involution $*$ acts by
\[
 \lambda\longleftrightarrow\mu/\lambda.
\]
A fixed idempotent would give an eigenvalue with $\lambda^2=\mu$, excluded in \eqref{m8:s:d:very-regular}. Thus there are exactly four pairs. The basis in \eqref{m8:s:d:cyclic} given by $e_i v_0$ is permuted by $h$ exactly as $e_i\mapsto e_i^*$, without any scalar factors. Therefore
\begin{equation}\label{m8:s:d:h-component}
 \det h=(-1)^4=1,\qquad h\in\SO(Q)(F).
\end{equation}
Although this determinant was computed after extension of scalars, the map $h$ itself was constructed over $F$. Combining \eqref{m8:s:d:inverse-transpose}, \eqref{m8:s:d:Q-component}, and \eqref{m8:s:d:adjoint-conjugacy} gives the rational conjugacy
\begin{equation}\label{m8:s:d:connected-conjugacy}
 (hQ)x^{-t}(hQ)^{-1}
 =\mu^{-1}hxh^{-1}=x^{-1},\qquad hQ\in\SO(Q)(F).
\end{equation}
The conjugator belongs to the determinant-one component of the group.

\paragraph{Harish--Chandra characters}
\begin{proposition}\label{m8:s:d:G8-duality}
If $\rho$ is an irreducible admissible smooth complex representation of $G_8^\circ$, and $\alpha_8(x)=x^{-t}$, then $\rho^{\alpha_8}\simeq\rho^\vee$.
\end{proposition}
\begin{proof}
By \cite[Theorem~1, \S3.3]{Clozel}, in characteristic zero the character distribution of an irreducible admissible representation is represented by a locally integrable function and is locally constant on the regular semisimple set; \S3.5 explicitly includes dyadic fields. Equality of character distributions determines the representation by \cite[Corollary~2.3.2]{Casselman}. The relevant character identities are
\begin{equation}\label{m8:s:d:character-identities}
 \Theta_{\rho^\vee}(g)=\Theta_\rho(g^{-1}),\qquad
 \Theta_{\rho^\beta}(g)=\Theta_\rho(\beta(g))
\end{equation}
for a Haar-preserving automorphism $\beta$.

To explain the distribution identities, admissibility makes $\rho(f)$ finite rank for every test function $f$. The character distribution is its trace. Dualizing this finite-rank operator gives the inversion identity; change of variables under $\beta$ gives the automorphism identity. Local integrability across the singular set is the analytic content needed to pass from almost-everywhere equality to equality of distributions.

For $x\in G_8^{\circ,\mathrm{vr}}$, \eqref{m8:s:d:connected-conjugacy} and conjugation invariance give
\[
 \Theta_{\rho^{\alpha_8}}(x)
 =\Theta_\rho(x^{-t})
 =\Theta_\rho(x^{-1})
 =\Theta_{\rho^\vee}(x).
\]
The two locally integrable functions agree on a full-measure subset, so their distribution characters agree everywhere. Character uniqueness proves the proposition. The algebraic construction and character theorem both apply in characteristic zero with arbitrary residue characteristic.
\end{proof}

\paragraph{The second factor and the product group}
Every nonscalar $b\in\GL_2(F)$ has a rational cyclic basis. In that basis it is
\[
 C=\begin{pmatrix}0&-d\\1&t\end{pmatrix},\qquad
 d=\det b\ne0,\quad t=\tr b.
\]
The symmetric invertible matrix
\[
 M=\begin{pmatrix}0&1\\1&t\end{pmatrix},\qquad\det M=-1,
\]
satisfies $MC=C^tM$. If $b=PCP^{-1}$, the matrix $g=P^{-t}MP^{-1}$ therefore satisfies $gbg^{-1}=b^t$. Scalar matrices equal their transposes. Applying the result to $b^{-1}$ gives rational conjugacy of $b^{-t}$ and $b^{-1}$, including repeated-root and nonsplit cases. The same character argument proves $\rho_2^{\alpha_2}\simeq\rho_2^\vee$ for $\alpha_2(b)=b^{-t}$.

More directly, on the full-measure regular semisimple subset $G_8^{\circ,\mathrm{vr}}\times B^{\mathrm{rs}}$, combine the two rational conjugators to carry $\alpha(x,b)$ to $(x,b)^{-1}$ within $\cG$. Apply the character argument to the irreducible admissible product representation itself. This gives
\begin{equation}\label{m8:s:d:product-duality}
 \pi^\alpha\simeq\pi^\vee,\qquad
 (\pi^\vee)^\alpha\simeq\pi.
\end{equation}
The second identity follows from involutivity and smooth biduality.

\paragraph{Identification of the period spaces}
Pullback under the group isomorphism $\alpha:H\to H^\sigma$, together with the character identity in \eqref{m8:s:p:character-inverses}, gives
\begin{equation}\label{m8:s:d:period-isomorphism}
 \begin{aligned}
 \cQ(\pi)
 &=\Hom_{H^\sigma}(\pi^\vee,\chi_2)\\
 &\simeq\Hom_H((\pi^\vee)^\alpha,\chi_2\circ\alpha)\\
 &=\Hom_H((\pi^\vee)^\alpha,\xi)
 \simeq\Hom_H(\pi,\xi)=\cP(\pi).
 \end{aligned}
\end{equation}
For an explicit map, choose an intertwiner
\[
 A_\pi:V\longrightarrow V^\vee,\qquad
 A_\pi\pi(g)=\pi^\vee(\alpha(g))A_\pi.
\]
Then $\lambda_1\mapsto\lambda_1\circ A_\pi^{-1}$ is a linear bijection from $\cP(\pi)$ to $\cQ(\pi)$. Indeed, for $k\in H^\sigma$,
\[
 \begin{aligned}
 \lambda_1(A_\pi^{-1}\pi^\vee(k)w)
 &=\lambda_1(\pi(\alpha(k))A_\pi^{-1}w)\\
 &=\chi_2(k)\lambda_1(A_\pi^{-1}w).
 \end{aligned}
\]
The map is a linear bijection for any choice of the intertwiner $A_\pi$.

The common-center condition is exactly \eqref{m8:s:p:center-condition}. As a further description, $b^{-t}=(\det b)^{-1}JbJ^{-1}$ and \eqref{m8:s:d:inverse-transpose} express $\alpha$ as an inner conjugation followed by scalar multiplication. Thus the duality is compatible with
\[
 \pi^\vee\simeq\pi\otimes
 \bigl[(x,b)\mapsto\omega_8(\mu(x))^{-1}\omega_2(\det b)^{-1}\bigr].
\]
On $H$, both scalar arguments equal $\det a$, so this twist is trivial whenever \eqref{m8:s:p:center-condition} holds. The individual central characters remain arbitrary, subject to the
common-center condition for a nonzero period.

\subsubsection{Application to generalized coefficients}
We now apply the smoothing construction in \S\ref{cf:coefficients}
to the two period spaces in \eqref{m8:s:p:periods}. For
$S=C_{\pi,\lambda_1,\lambda_2}$ put
$\Delta S=S-\sigma_*S$. The opposite character in
\eqref{m8:s:p:periods} gives exactly the plus-plus equivariance
\eqref{m8:s:f:translation-character}.

By the common-center calculation \eqref{cf:center-span}, $\Delta S$
is finite under the center of the ambient product group, with
possibly different eigencharacters. Applying
Theorem~\ref{cf:nilpotent} before descent gives
$\Delta S\in\NW(G_8)$. Thus
\eqref{m8:s:f:transpose-intertwining} and
\eqref{m8:s:f:native-transfer} give
\begin{equation}\label{m8:s:c:descendant}
 T=\operatorname{Fr}_8(\Delta S)
   \in\Do(G_8^\circ)\cap\NW(G_8^\circ).
\end{equation}
Equivalently,
\begin{equation}\label{m8:spectral-class}
 \mathscr V_8^-=
 \spanop_{\CC}\{\operatorname{Fr}_8
 (C_{\pi,\lambda_1,\lambda_2}
       -\sigma_*C_{\pi,\lambda_1,\lambda_2})\}=0,
\end{equation}
where the span consists of finite sums. This is an application of
the NW distribution theorem.

For every pair of periods, Theorem~\ref{m8:distribution} gives
$T=0$; reconstruction by \eqref{m8:s:f:reconstruction} then yields
$\Delta S=0$. We obtain the product bound from the coefficient-kernel
argument in \S\ref{cf:kernel}. By the period isomorphism
\eqref{m8:s:d:period-isomorphism}, this product is the square of
the Model-8 multiplicity. This proves Theorem~\ref{m8:main}.

\part{The two central forms of \texorpdfstring{$D_6$}{D6} and Model 9}\label{part:3}
\section{The \texorpdfstring{$D_6$}{D6} theorem and connected Model 9}
\label{m9:section}

\paragraph{Proof strategy}
We treat both rational central forms of type $D_6$: the vector-similitude
group $\GSO_{12}$ of Model~9 and the semispin proper-Levi form needed
for Model~10. For each form, the distribution proof has two stages.
We begin on the Siegel Gauss open, where Frobenius descent extracts
the two free unipotent factors. The resulting Levi distribution satisfies the balanced
Ginzburg--Rallis theorem of Lemma~\ref{m9:balanced-gr}, with the
central parameters retained and their cotangent components zero.
Applying this theorem clears the initial open (\S\ref{m9:siegel});
the rational covers and subgroup preimages retain the different
central forms and all their rational sectors.

We then process the complement in increasing cell codimension.
Within each cell, unipotent stabilizer characters confine support
to its full compatible locus. The $35$ nonempty compatible boundary
loci have nonnilpotent ambient conormals at every point
(Lemma~\ref{m9:conormal-family}). After their predecessors have
vanished, the NW condition of the original ambient distribution
excludes these supports (\S\ref{m9:boundary-section}). Continuing in this way through the boundary proves the
distribution theorem and hence the GK product bound; rational duality
then proves multiplicity one on both forms
(\S\ref{m9:gp1-section} and \S\ref{m9:semispin-gp1-section}).

\subsection{The two central forms and the literal datum}
\label{m9:datum}

Let $J=\left(\begin{smallmatrix}0&-1\\1&0\end{smallmatrix}\right)$,
and let $Q$ be the $12$ by $12$ matrix whose six antidiagonal $2$ by $2$
blocks are $J,J,J,-J,-J,-J$.  Thus $Q^t=Q$, $Q^2=1$, and $\det Q=1$.
The literal group, subgroup, and character are
\begin{align}
 G_{9,v}&=\{g:g^tQg=\mu(g)Q,\ \det g=\mu(g)^6\},\label{m9:gso}\\
 H_{9,v}&=D_6(\GL_2)\ltimes U_9,
 &D_6(a)&=\operatorname{diag}(a,a,a,a,a,a),\label{m9:h-vector}\\
 \xi_{9,\kappa}(u)&=\psi(\operatorname{tr}X+
           \operatorname{tr}Y+\kappa z),
 &\kappa&\in\F^\times.\label{m9:phase}
\end{align}
Here $U_9$ is the unipotent radical of the parabolic whose Levi is
\begin{equation}\label{m9:o:eq:2.3}
 \left\{\operatorname{diag}\left(g_1,g_2,g_3,
 \frac{\mu g_3}{\det g_3},\frac{\mu g_2}{\det g_2},
 \frac{\mu g_1}{\det g_1}\right)\right\}.
\end{equation}
Here $X,Y,Z$ are the adjacent upper blocks, with $Z=zI_2$, and we
extend the character trivially to the reductive subgroup. In the literal Wan--Zhang datum
$\kappa=2$; see \cite[\S8.3]{WZ}.  For the inherited $D_6$ pinning in
Model~10, the coefficient is $\kappa=1$. Proving the statement for every nonzero
$\kappa$ treats these two coefficients over all fields under consideration,
including fields of residue characteristic two.

Use simple roots
\begin{equation}\label{m9:o:eq:2.6}
 \beta_i=e_i-e_{i+1}\ (1\leq i\leq5),\qquad
 \beta_6=e_5+e_6,\qquad I=\{\beta_1,\beta_3,\beta_5\}.
\end{equation}
The exceptional roots $\alpha_i$ use the Bourbaki numbering fixed
in \S\ref{m10:datum}; their coroots are denoted by $H_i$ in the
exceptional lattice calculations below.
With this choice, $U_9$ has the $27$ positive roots outside $I$, and
its reductive stabilizer contains the diagonal $\SL_2$ in the three
commuting root groups in $I$. The correspondence with the exceptional pinning is
\begin{equation}
 (\beta_1,\ldots,\beta_6)
 =(\alpha_7,\alpha_6,\alpha_5,\alpha_4,\alpha_2,\alpha_3).
 \label{m9:root-bridge}
\end{equation}
The transpose-preserving change of basis
$\operatorname{diag}(1,1,1,1,1,1,1,-1,1,-1,1,-1)$ carries $Q$
to the antidiagonal identity form.  In the ordered basis
$(1,2,3,4,5,6,-6,-5,-4,-3,-2,-1)$ it carries $zI_2$ to
$zE_{e_5+e_6}$.

Put $S=\operatorname{Spin}_{12}$ and let
$i:\SL_2\to S$ be the diagonal embedding in the three commuting
root-$\SL_2$ subgroups indexed by $I$.
Let $c_v$ denote the nonidentity element of the kernel of the vector
representation and put $c_s=i(-I_2)$. These central
elements are distinct: their vector images are $I_{12}$ and $-I_{12}$.
The two central forms are
\begin{align}
 G_{9,v}&=(S\times\mathbb G_m)/
       \langle(c_v,1),(c_s,-1)\rangle,\label{m9:vector-form}\\
 G_{9,s}=M_1&=(S\times\mathbb G_m)/
       \langle(c_s,1),(c_v,-1)\rangle.\label{m9:semispin-form}
\end{align}
Write $D_s=S/\langle c_s\rangle$ and $\bar c_v$ for the image of $c_v$.
Denote the vector map $S\to\SO_{12}$ by $\rho_v$ and the quotient
map $S\to D_s$ by $\rho_s$.
Then
\begin{equation}\label{m9:o:eq:2.10}
 G_{9,s}=(D_s\times\mathbb G_m)/\langle(\bar c_v,-1)\rangle,
\end{equation}
and its subgroup is
\begin{equation}
 H_{9,s}=\operatorname{PGL}_2\ltimes U_9.
 \label{m9:h-semispin}
\end{equation}
With $z_1=\omega_1^\vee=(2,2,3,4,3,2,1)$ in the exceptional coroot
basis, $z_1(-1)=\bar c_v$ and
\begin{equation}\label{m9:o:eq:2.12}
 \chi_1(z_1(t))=t^2,
\end{equation}
where $\chi_1=\omega_1$ is the primitive character of $M_1$.

\paragraph{The central lattice calculation.}\label{m9:o:sec:2.3}
In the exceptional coroot basis, the Cartan matrix $A$ has determinant
$2$ and
\[
 v=\tfrac12(H_2+H_5+H_7),\qquad
 Av=(0,1,0,-1,1,-1,1).
\]
Moreover
\[
 z_1-(H_2+H_3)-(H_2+H_5+H_7)=(2,0,2,4,2,2,0).
\]
This even difference gives $z_1(-1)=\bar c_v$ in $D_s$.
In standard $D_6$ coordinates one may take
\[
 c_s=(e_1-e_2+e_3-e_4+e_5-e_6)(-1),\qquad
 c_v=(2e_5)(-1).
\]
Their pairings with the two half-spin weights are $(0,1)$ and $(1,1)$,
respectively; the vector weights pair oddly with $c_s$ and evenly with
$c_v$. These identities determine the two kernels in
\eqref{m9:vector-form}--\eqref{m9:semispin-form}.

The natural finite cover
\begin{equation}
 \widehat G_{9,s}=D_s\times\mathbb G_m,
 \qquad q_s(d,t)=dz_1(t)
 \label{m9:natural-map}
\end{equation}
is a twofold cover of $G_{9,s}$.

\subsection{A common cover, rational sectors, and the pinned operation}
\label{m9:common-cover-geometry}

The appropriate common group is the fiber product over the pinned
adjoint group:
\begin{equation}
 \mathcal C_9=G_{9,v}\times_{\operatorname{PSO}_{12}}G_{9,s}
 =\frac{S\times\mathbb G_{m,t}\times\mathbb G_{m,u}}
 {\langle(c_v,1,-1),(c_s,-1,1)\rangle}.
 \label{m9:common-cover}
\end{equation}
Its projections
\begin{equation}\label{m9:o:eq:2.17}
 p_v:\mathcal C_9\to G_{9,v},\qquad p_s:\mathcal C_9\to G_{9,s}
\end{equation}
are induced by
\begin{equation}\label{m9:o:eq:2.16}
 (s,t,u)\mapsto(t\rho_v(s),\rho_s(s)z_1(u)).
\end{equation}
The kernel of $p_v:\mathcal C_9\to G_{9,v}$ is the $u$-center,
while that of $p_s:\mathcal C_9\to G_{9,s}$ is the $t$-center.
Since both kernels are connected split copies of $\mathbb G_m$,
Hilbert~90 makes both maps surjective on $\F$-points.
The derived group of $\mathcal C_9$ is $S$, its center is a split
two-dimensional torus, and
\begin{equation}
 \mathfrak c_9=\mathfrak{so}_{12}\oplus\F\oplus\F,
 \qquad \cN(\mathfrak c_9^*)
 =\cN(\mathfrak{so}_{12}^*)\times\{0\}\times\{0\}.
 \label{m9:two-central-covectors}
\end{equation}

If $h_s:\operatorname{PGL}_2\to G_{9,s}$ is the indicated embedding,
set
\begin{equation}\label{m9:o:eq:2.19}
 R_{9,0}=\{(D_6(a),h_s([a])):a\in\GL_2\},
 \qquad R_9=R_{9,0}\ltimes U_9.
\end{equation}
Then
\begin{equation}\label{m9:o:eq:2.20}
 p_s^{-1}(H_{9,s})=R_9,\qquad p_s(R_9(\F))=H_{9,s}(\F).
\end{equation}
The second map has scalar split-torus kernel, whereas $p_v$ identifies
$R_9$ with $H_{9,v}$. In addition,
\begin{equation}\label{m9:o:eq:2.21}
 p_v^{-1}(H_{9,v})=R_9Z_u.
\end{equation}
Extend \eqref{m9:phase} trivially across $R_{9,0}$ to obtain $\xi_R$.
These statements retain all determinant classes of $a\in\GL_2(\F)$.

Let $V_\F=\F^\times/(\F^\times)^2$.  The rational cokernels of
$S\times\mathbb G_m\to G_{9,v}$ and
$S\times\mathbb G_m\to G_{9,s}$ are both $V_\F^2$;
those of $S\times\mathbb G_m\to\widehat G_{9,s}$ and $q_s$ are
$V_\F$.  The rational cokernel of
$S\times\mathbb G_m^2\to\mathcal C_9$ is $V_\F^2$.
In contrast, each of $p_v,p_s$ and each scalar-center quotient to
$\operatorname{PSO}_{12}$ is surjective on rational points.
For the finite rows these assertions follow by restricting to split
maximal tori and applying Hilbert~90 there.

More explicitly,
\begin{equation}
 J_s=q_s(\widehat G_{9,s}(\F))
 =\{m:\chi_1(m)\in(\F^\times)^2\},\qquad
 G_{9,s}(\F)=\coprod_{a\in V_\F}J_s s_1(a),
 \label{m9:sector-cosets}
\end{equation}
where representatives of $V_\F$ are chosen and
$s_1(a)=\alpha_1^\vee(a)$.  The full inverse image of the acting subgroup is
\begin{equation}
 \widehat H_{9,s}(\F)=
 \{(h\bar c_v^{\epsilon},(-1)^{\epsilon}):
       h\in H_{9,s}(\F),\ \epsilon\in\{0,1\}\}.
 \label{m9:full-subgroup-preimage}
\end{equation}
It maps onto all of $H_{9,s}(\F)$, including every determinant
square class in $\operatorname{PGL}_2(\F)$.

The NW Chevalley anti-involution $\iota_9$ fixes the torus and
exchanges the pinned positive and negative root generators.  It is
ordinary transpose on the literal matrix group.  Its doubled point
action is
\begin{equation}
 (h,k)x=hx\iota_9(k),\qquad
 \eta(h,k)=\xi(h)\xi(k).
 \label{m9:plus-plus}
\end{equation}
The operation exchanges the two abstract acting factors, sending $H$
to its opposite.  On \eqref{m9:common-cover} it fixes $t,u$ and applies
$\iota_S$ to $s$.  It fixes each finite central kernel pointwise.
On the $a$-sector coordinate $m=q_s(d,t)s_1(a)$ it becomes
\begin{equation}
 \widehat\iota_{9,a}(d,t)=
  (\Ad(s_1(a))\iota_{D_s}(d),t).
 \label{m9:sector-operation}
\end{equation}
The conjugation factor in this formula is an $\F$-automorphism of $D_s$.

Deck invariance under $(\bar c_v,-1)$ permits two central parities:
\begin{equation}
 L_{\bar c_v}T_\varepsilon=\varepsilon T_\varepsilon,
 \qquad L_{t\mapsto-t}T_\varepsilon=\varepsilon T_\varepsilon,
 \qquad\varepsilon=\pm1.
 \label{m9:central-parity}
\end{equation}
Since the two finite projectors preserve nilpotent wavefront and
all equivariance conditions, we retain every
$(a,\varepsilon)\in V_\F\times\{\pm1\}$.

The explicit rational representatives and subgroup lifts are recorded in
\S\ref{m9:o:sec:2.6}. The finite-cover sectors will be used in
Corollary~\ref{m9:natural-cover}.

\subsection{Balanced Ginzburg--Rallis vanishing}
\label{m9:balanced-section}

Let
\[
 R_{\mathrm{GR}}=D_3(\GL_2)U_{2,2,2}\subset\GL_6,
 \quad \epsilon(D_3(a)u)=\psi(\operatorname{tr}u_{12}
                        +\operatorname{tr}u_{23}),
 \quad \nu(D_3(a)u)=\det a.
\]
Here $D_3(a)=\diag(a,a,a)$, and $U_{2,2,2}$ is the unipotent
radical of the upper block parabolic of $\GL_6$ with three blocks
of size two.
We use Theorem~\ref{cf:gr-nw} as the NW input for
$R_{\mathrm{GR}}\times R_{\mathrm{GR}}$. For the Siegel descent we need
the following balanced form, including control of the passive directions.

\begin{lemma}[Balanced Ginzburg--Rallis theorem]
\label{m9:balanced-gr}
Put
$A_0=\{(r_1,r_2)\in R_{\mathrm{GR}}^2:
\nu(r_1)\nu(r_2)=1\}$.
Every transpose-odd, $(A_0,\epsilon\boxtimes\epsilon)$-equivariant
distribution on $\GL_6(\F)$ with wavefront in
$\cN(\mathfrak{gl}_6^*)$ is zero.  The conclusion also holds with
passive analytic parameters, provided all wavefront components in
those parameter directions are zero.
\end{lemma}

\begin{proof}
Apply Lemma~\ref{pl:central-balancing} with $X=\GL_6$, $p=\det$,
and $k=3$. The character theorem is Theorem~\ref{cf:gr-nw}.
Conjugation by $a$ on each matrix block of $U_{2,2,2}$ has determinant
one, so $R_{\rm GR}$ and the balanced acting groups are unimodular.
As torsors under $\SL_6$, the determinant fibers carry invariant
densities that disintegrate against multiplicative Haar measure. Scalar
translation by $tI_6$ changes determinant by $t^6$, giving the required
central coordinate. Its derivative is nonzero over $F$, including in
residue characteristics two and three. All hypotheses of the balancing
lemma, including the zero wavefront in fixed parameter directions,
are therefore satisfied.
\end{proof}
To prove the passive-parameter assertion, pair against an arbitrary
compactly supported test in those parameters. Compact-parameter projection
\cite[Theorem~2.9.6]{CHLR} bounds the resulting wavefront by cotangent vectors
with zero parameter components, hence by
$\cN(\mathfrak{gl}_6^*)$. The balancing argument applies to every such
pairing, and product tests span the test-function space.

\subsection{The compatible boundary and its conormal}
\label{m9:boundary-section}

Both central forms admit the following $35$-cell continuation.

Let $Q_S=R_SN_S$ be the Siegel parabolic obtained by deleting
$\beta_6$, and set $V_S=N_SR_SN_S^-$.  Put $P_0=L_0U_9$ and
$P_0^-=\iota_9(P_0)$.  The $D_6$ rational cells are
\begin{equation}\label{m9:o:eq:5.5}
 X_v=P_0(\F)n_vP_0^-(\F),\qquad
 v\in W_I\backslash W(D_6)/W_I,
\end{equation}
with minimal representatives and pinned Tits lifts. Split Bruhat
decomposition gives these representatives and root-group charts over
$F$ with the full rational torus of each central form.
The $D_6$ Weyl
group consists of the $23{,}040$ even signed permutations.  Minimality
means that $v(\beta_i)$ and $v^{-1}(\beta_i)$ are positive for
$i=1,3,5$.  The finite localization is
\begin{center}
\begin{tabular}{lrrr}
\hline
Region & Compatible & Incompatible & Total\\
\hline
Siegel Gauss open &13&8&21\\
Siegel boundary &35&436&471\\
Whole group &48&444&492\\
\hline
\end{tabular}
\end{center}
Membership in $V_S$ means that the representative is unsigned, or
equivalently belongs to the Siegel $A_5$ Weyl subgroup.

For the ordered phase roots use
\begin{equation}\label{m9:o:eq:5.3}
 \begin{aligned}
 (\gamma_0,\ldots,\gamma_8)=(&e_1-e_3,e_1-e_4,e_2-e_3,e_2-e_4,\\
                           &e_3-e_5,e_3-e_6,e_4-e_5,e_4-e_6,e_5+e_6).
 \end{aligned}
\end{equation}
Let $\lambda_\kappa\in\mathfrak u_9^*$ be the phase functional
$\operatorname{tr}X+\operatorname{tr}Y+\kappa z$ of
\eqref{m9:phase}. For $l\in L_0$, put
$p_l=\lambda_\kappa\circ\Ad(l)$ and
$p_j=p_l(E_{\gamma_j})$ for $0\le j\le8$.
The phase image is
\begin{equation}
 \mathcal O_9=\left\{p:
 p_8\det\begin{pmatrix}p_0&p_1\\p_2&p_3\end{pmatrix}
       \det\begin{pmatrix}p_4&p_5\\p_6&p_7\end{pmatrix}\ne0\right\}.
 \label{m9:phase-open}
\end{equation}
To see that the phase map is a surjective submersion on each central
form, first consider the vector form. A rational section is
$g_3=I$, $g_2=Q'^t$, $g_1=Q'^tP^t$, $\mu=\kappa/p_8$,
where $P,Q'$ are the two displayed phase matrices.  The semispin
statement follows through the surjective connected-center maps of the
two Levis to their common adjoint Levi, including their differentials.

At $x=ln_v\iota_9(r)$ put $p=p_l$, $q=p_r$.  The entire unipotent
stabilizer is parametrized by
\begin{equation}\label{m9:o:eq:5.6}
 Y\longmapsto
 (l\exp(-\Ad(n_v)Y^t)l^{-1},\ r\exp(Y)r^{-1}),
 \quad Y\in\mathfrak u_9,\quad\Ad(n_v)Y^t\in\mathfrak u_9.
\end{equation}
Its character is
\begin{equation}\label{m9:o:eq:5.7}
 \psi(q(Y)-p(\Ad(n_v)Y^t)).
\end{equation}
Thus the complete
compatible locus $C_v$ is defined by
\begin{equation}
 c(v,\gamma)p(E_{-v\gamma})=q(E_\gamma),
 \quad \Ad(n_v)E_\gamma^t=c(v,\gamma)E_{-v\gamma},
 \quad c(v,\gamma)\in\{\pm1\}.
 \label{m9:compatibility}
\end{equation}
Unlisted phase coordinates are zero.  Each nontrivial equation is a
forced zero or a signed one-to-one match.  These equations are
independent throughout the open \eqref{m9:phase-open} on both
sides.

\paragraph{Rational compatibility.}\label{m9:o:sec:5.4}
They also give a complete rational nonemptiness test: neither
$p_8$ nor $q_8$ may be forced to zero, and each required determinant
must have an unblocked monomial. Indeed, since signed matches identify each coordinate at most once,
each unblocked determinant is a nonzero polynomial in independent
variables. Their product therefore has a nonvanishing
point over the infinite field $F$.
The chart
\begin{equation}\label{m9:o:eq:5.9}
 U_9\times L_0\times L_0\times U_9\longrightarrow X_v,
 \qquad(u,l,r,z)\longmapsto ul n_v\iota_9(r)\iota_9(z)
\end{equation}
is a surjective submersion. Under the phase submersion, compatibility becomes a system of
independent linear equations. Intrinsically, it is the zero locus of
the stabilizer-character differential on the constant-rank
stabilizer bundle. Thus $C_v$ is smooth and closed in
$X_v$, with codimension equal to the number of nontrivial equations.
Off $C_v$,
unipotent stabilizer-character incompatibility gives ordinary
equivariant vanishing.

Use root matrices
$E_{e_i-e_j}=E_{i,j}-E_{-j,-i}$,
$E_{e_i+e_j}=E_{i,-j}-E_{j,-i}$,
$E_{-\gamma}=E_\gamma^t$.
The Cartan matrices and invariant form are
\begin{equation}\label{m9:o:eq:5.2}
 H_i=E_{i,i}-E_{-i,-i},\qquad
 B(A,B)=\tfrac12\operatorname{Tr}(AB).
\end{equation}
For $F_Y=p(\Ad(n_v)Y^t)-q(Y)$ and Levi tangent variables $A,B$,
\begin{equation}\label{m9:o:eq:5.10}
 dF_Y(A,B)=p([A,\Ad(n_v)Y^t])-q([B,Y]).
\end{equation}
Writing $D_v$ for all these differentials, the tangent space of the compatible locus, after right
translation and then $\Ad(l^{-1})$, is
\begin{equation}
 \mathcal T_v(p,q)=\mathfrak u_9+\Ad(n_v)\mathfrak u_9^-
 +\{A+\Ad(n_v)B^t:(A,B)\in\ker D_v\}.
 \label{m9:full-tangent}
\end{equation}
Here $A,B$ include all Levi root, Cartan, and central directions.

\begin{lemma}[Uniform boundary conormal]
\label{m9:conormal-family}
For each of the $35$ compatible cells outside $V_S$, there is an
explicit vector $C$ annihilating \eqref{m9:full-tangent} and a phase
coordinate $p_j$ such that
\begin{equation}
 B(C,C)=2p_j^2\ne0
 \label{m9:nonzero-norm}
\end{equation}
at every point of the complete compatible locus.
\end{lemma}

\begin{proof}
The formulas in Table M9-C35 of \path{anc/model9_calculations.pdf}
specify a selected equation $F_Y$ and satisfy
\begin{align}
 B(C,\mathfrak u_9)=B(C,\Ad(n_v)\mathfrak u_9^-)&=0,
 \label{m9:o:eq:5.12}\\
 B(C,A+\Ad(n_v)B^t)&=dF_Y(A,B).
 \label{m9:o:eq:5.13}
\end{align}
for every Levi pair $(A,B)$, modulo the equations
\eqref{m9:compatibility}.  For example, cell $1$, with word $6$,
has $q_8=-p_8$ and $C=-p_8(H_5+H_6)$.
Cell $10$, with word $4364$, has
\[
 C=p_4(H_3-H_5)-p_5E_{e_6-e_5}+q_5E_{e_3+e_4}.
\]
Here the root terms account for the non-Cartan Levi derivatives, and
the open determinant conditions together with $p_6=0$ force $p_4\ne0$.
The coordinate $p_8$ is nonzero by the phase-open condition. Rows
selecting $p_0$ or $p_3$ have $p_2=0$; rows selecting $p_4$ or $p_7$
have $p_6=0$. The corresponding triangular phase matrix is invertible,
so the selected diagonal entry is nonzero. In cell~$19$, $q_6=0$
forces $q_4\ne0$, and $q_4=-p_6$ gives $p_6\ne0$.
These cases exhaust the table. To verify the identities, one checks the
$54$ outer root directions and every Cartan, Levi-root and central
tangent direction, modulo the full compatibility equations;
Section M9-V gives the matrix calculation.
At a general point of the submersion chart the conormal is
$\Ad(ul)C$, using the rational adjoint action.
\end{proof}

\begin{theorem}[$D_6$ boundary continuation]
\label{m9:boundary-continuation}
Let $G$ be $G_{9,v}$ or $G_{9,s}$ and let $\kappa\ne0$.
Suppose a distribution $T$ on $G(\F)$ has nilpotent wavefront
and is equivariant under $U_9\times U_9$ for
$(u,v)x=ux\iota_9(v)$ with character
$\xi_{9,\kappa}(u)\xi_{9,\kappa}(v)$.
If $T|_{V_S}=0$, then $T=0$.
\end{theorem}

\begin{proof}
Process the boundary cells in increasing codimension; for a minimal
representative,
\begin{equation}\label{m9:o:eq:5.15}
 \operatorname{codim}_G X_v=\ell(v).
\end{equation}
Indeed $\mathfrak p_0\cap\Ad(n_v)\mathfrak p_0^-$ has the full
rank plus $6+\ell(v)$ root directions, while
$\dim P_0=\operatorname{rank}(G)+33$.
Algebraic orbit closure gives the open set
\begin{equation}\label{m9:o:eq:5.16}
 W_v=V_S\ \cup\!
 \bigcup_{\operatorname{codim}X_w<\operatorname{codim}X_v}X_w
 \ \cup X_v.
\end{equation}
Once all predecessor cells have vanished, the open/closed support
sequence confines the ambient restriction
$T|_{W_v}$ first to $X_v$ and then, by stabilizer incompatibility, to
$C_v$.  This clears the empty compatible loci. For the other $35$ cells, Lemma~\ref{m9:conormal-family} supplies a
conormal outside the nilpotent cone at every support point.  Conormal
translation invariance of the wavefront contradicts the inherited
NW condition; equivalently the quadratic norm of a translate by
$tC$ has nonzero leading coefficient $2p_j^2$.
Thus $T|_{W_v}=0$, and the finite induction proves the assertion.
\end{proof}

The boundary theorem uses the unipotent equivariance and nilpotent
wavefront of the ambient distribution. At each step the conormal
argument applies to its restriction to an open, after all predecessor
cells have vanished.

The complete compatible-cell and conormal arrays are specified in
\S\ref{supp:d6-expanded}, with stable identifiers M9-K48 and M9-C35.

\subsection{Siegel descent on the vector and semispin forms}
\label{m9:siegel}

The Siegel radical $N_S$ has $15$ roots $e_i+e_j$, and
$U_9=U_{2,2,2}\ltimes N_S$.  Both outer factors in
$N_SR_SN_S^-$ act freely.  Extracting their character-weighted Haar
lines yields an injective Frobenius correspondence with an NW
Levi distribution.  The relative density character is trivial: on the vector form,
\[
 \mathfrak n_S\simeq\bigwedge^2\F^6\otimes\mu^{-1},\qquad
 (\det D_3(a))^5(\det a)^{-15}=1,
\]
and on the semispin form the residual reductive group is
$\operatorname{PGL}_2$, which has no algebraic characters.
Unipotent factors have determinant one.  Transpose exchanges the
two Haar lines with scalar one.  The inherited Levi cotangent vectors are
nilpotent with zero components in every central direction, by
Lemma~\ref{cf:gauss}.

For $G_{9,v}$ the Siegel Levi is $\GL_6\times\mathbb G_m$, with
coordinates $(g,\mu)$ and transpose $(g^t,\mu)$.  The residual subgroup
is $(D_3(a)u,\det a)$.  Its doubled action changes $\mu$ by
$\det a\det b$.  The stabilizer of the fiber $\mu=1$ is exactly
the balanced GR group, with density ratio one.
Lemma~\ref{cf:frobenius-wf} transfers wavefront to its Frobenius
descendant: restriction from the ambient product cone
$\cN(\mathfrak{gl}_6^*)\times\{0\}$ to the literal fiber
$\GL_6\times\{1\}$ leaves a nilpotent $\GL_6$ cotangent vector.
Lemma~\ref{m9:balanced-gr} clears the Siegel open.

The semispin Siegel Levi has the connected Siegel-Levi form
\begin{equation}
 R_s=\frac{\mathbb G_{m,A_1}\times\mathbb G_{m,A_2}\times\GL_6}
 {\{(z^{-2},z^{-2},zI_6):z\in\mathbb G_m\}}.
 \label{m9:siegel-cover}
\end{equation}
It is surjective on rational points.  For
$g=\operatorname{diag}(B_1,\ldots,B_6)$, the precise root-coordinate map is
\begin{equation}
 \begin{array}{lll}
 \alpha_1=A_1/A_2,&\alpha_2=B_5/B_6,&\alpha_3=A_2B_5B_6,\\
 \alpha_4=B_4/B_5,&\alpha_5=B_3/B_4,&\alpha_6=B_2/B_3,
 \quad\alpha_7=B_1/B_2.
 \end{array}
 \label{m9:siegel-torus-map}
\end{equation}
The primitive kernel vector is $(-2,-2,1,1,1,1,1,1)$ and the
seven nonzero Smith invariants are one.  The subgroup preimage is
\begin{equation}
 \widetilde H_R=
 \{((\det a)^{-1},(\det a)^{-1},D_3(a)u):
                   a\in\GL_2,\ u\in U_{2,2,2}\}.
 \label{m9:siegel-full-subgroup}
\end{equation}
Its image contains all $F$-points of the $\operatorname{PGL}_2$ torus:
$a=\operatorname{diag}(t,1)$ has root values
\[
(1,t,1,t^{-1},t,t^{-1},t).
\]
Positive and negative unipotent generators map diagonally to the three
leaf root groups. Hence the subgroup map is the full quotient
$\GL_2(F)\to\operatorname{PGL}_2(F)$; its scalar kernel, obtained by
$a=zI_2$, is exactly the cover kernel. The dimensions and these root
maps prove that \eqref{m9:siegel-full-subgroup} is the full preimage.
The lifted anti-involution is $(A_1,A_2,g)\mapsto(A_1,A_2,g^t)$,
and the pulled-back wavefront lies in
$\{0\}^2\times\cN(\mathfrak{gl}_6^*)$.

The doubled action scales both $A_i$ by $(\det a\det b)^{-1}$.
Normalize $A_2=1$; then $r=A_1/A_2$ is a passive parameter and the
stabilizer is balanced GR.
Lemma~\ref{cf:frobenius-wf} supplies an ambient cotangent vector in the wavefront set
for every descendant cotangent vector.  Restricting
$(0,0,\xi_g)$ to the fiber $(r,1,g)$ gives $(0,\xi_g)$, with
$\xi_g$ nilpotent and zero cotangent component in $r$.
The passive-parameter form of Lemma~\ref{m9:balanced-gr} therefore
kills the descendant.
Moreover,
\begin{equation}
 \chi_1=A_1^2A_2\det g=r^2\det g\quad(A_2=1),
 \qquad \bar c_v:r\longmapsto-r.
 \label{m9:retained-sectors}
\end{equation}
Thus the odd central parity remains present as odd dependence on $r$.
The rational sector $[\chi_1]=a$ is precisely $[\det g]=a$ in these
coordinates.  Each such determinant-square-class locus is clopen and
stable under the balanced action and transpose, so extension by zero
applies separately in every sector.

\begin{theorem}[Distribution theorem on the semispin central form]
\label{m9:native-semispin}
\label{m9:native}
Let $G_{9,s}$, $H_{9,s}$, $\xi_{9,\kappa}$, and $\iota_9$ be as
above, with $\kappa\in\F^\times$.  An $\iota_9$-odd distribution
on $G_{9,s}(\F)$, equivariant for \eqref{m9:plus-plus}, is zero if
\[
 \WF(T)\subset G_{9,s}(\F)\times
       \bigl(\cN(\mathfrak{hspin}_{12}^*)\times\{0\}\bigr)
\]
under cotangent trivialization.
\end{theorem}

\begin{proof}
After semispin Siegel descent clears $V_S$, apply
Theorem~\ref{m9:boundary-continuation} to its complement.
\end{proof}

\begin{theorem}[Distribution theorem on the literal connected form]
\label{m9:native-vector}
For every $\kappa\ne0$, every transpose-odd distribution on
$G_{9,v}(\F)=\GSO_{12}(\F)$ equivariant under
$H_{9,v}\times H_{9,v}$ with character
$\xi_{9,\kappa}\boxtimes\xi_{9,\kappa}$ for
\eqref{m9:plus-plus}, and with nilpotent wavefront, is zero.
\end{theorem}

\begin{proof}
After vector-form Siegel descent clears $V_S$, apply
Theorem~\ref{m9:boundary-continuation} to its complement.
\end{proof}

\begin{remark}
For the vector form, the parabolic stabilizing a two-dimensional
isotropic subspace has Levi factor $\GL_2\times\GSO_8$. This
explains the geometric relation to Model~8 in
Figure~\ref{fig:reductions}.
\end{remark}

\subsection{Natural-cover and common-cover formulations}

\begin{corollary}[The exact natural-cover distribution theorem]
\label{m9:natural-cover}
For every $a\in\F^\times/(\F^\times)^2$, let $T_a$ be a
distribution on $\widehat G_{9,s}(\F)$, invariant under
$\ker q_s(\F)=\langle(\bar c_v,-1)\rangle$, odd under
\eqref{m9:sector-operation}, and equivariant under the subgroup
\eqref{m9:full-subgroup-preimage} on both sides with the pulled-back
character and point action
$(h,k)x=hx\widehat\iota_{9,a}(k)$.
If its wavefront lies in
$\cN(\mathfrak{hspin}_{12}^*)\times\{0\}$ after trivialization,
then $T_a=0$.  This includes both parities in
\eqref{m9:central-parity}.
\end{corollary}

\begin{proof}
Finite analytic inverse charts and deck averaging descend $T_a$ to
the clopen coset $J_ss_1(a)$.  The maps intertwine the transported
action and anti-involution and identify the cotangent cones.
Extend by zero to $G_{9,s}(\F)$ and apply
Theorem~\ref{m9:native-semispin}.  Pullback is injective.
\end{proof}

\begin{corollary}[Common-cover distribution theorem]
\label{m9:common-cover-theorem}
Every odd $(R_9\times R_9,\xi_R\boxtimes\xi_R)$-equivariant
distribution on $\mathcal C_9(\F)$ for the common pinned operation
and doubled point action is zero if its wavefront lies in
$\cN(\mathfrak{so}_{12}^*)\times\{0\}^2$ after trivialization.
The same conclusion holds on $\mathcal C_9/Z$ for any connected
split central subtorus $Z$, with the image subgroup and descended
character and operation.
\end{corollary}

\begin{proof}
Since the kernel of $p_s$ lies in $R_9$ and has trivial character,
central Frobenius descent gives a distribution to which
Theorem~\ref{m9:native-semispin} applies, with the additional zero central
cotangent component prescribed by \eqref{m9:two-central-covectors}.
For a further connected split central quotient, Hilbert~90 makes
pullback injective on rational points and adds zero kernel cotangent vectors.
\end{proof}

The finite natural cover is treated sector by sector in
Corollary~\ref{m9:natural-cover}. Pullback of a literal distribution
along $p_v$ belongs to the common theorem's class and is $Z_u$-invariant.

\subsection{Connected orthogonal duality and literal GP1}
\label{m9:gp1-section}

Set $G=G_{9,v}$, $H=H_{9,v}$, $\xi=\xi_{9,2}$, and
$\alpha(g)=g^{-t}$.  The two coefficient spaces in the GK argument are
\begin{equation}
 \mathcal P(\pi)=\Hom_H(\pi,\xi),\qquad
 \mathcal Q(\pi)=\Hom_{H^{\iota_9}}(\pi^\vee,\chi_2),
 \quad\chi_2(\iota_9(h))=\xi(h)^{-1}.
 \label{m9:two-hom-spaces}
\end{equation}
Their dimensions are $m_9(\pi)$ and $q_9(\pi)$ in the notation
of the introduction. The opposite subgroup and its character are
retained until the explicit duality comparison below.
For functionals $\lambda_1,\lambda_2$ in these spaces, define
\[
 \Phi_1(f)(v)=\int_G f(g)\lambda_1(\pi(g)v)\,dg,
 \qquad S(f)=\lambda_2(\Phi_1(f)).
\]
Compact-open smoothing gives $\Phi_1(f)\in V^\vee$.  For
$L_hf(g)=f(h^{-1}g)$ and $R_kf(g)=f(gk)$, the inverse test action of
$g\mapsto h_1g\iota_9(h_2)$ is
$L_{h_1}R_{\iota_9(h_2)^{-1}}$; its character is plus-plus because
$\chi_2(\iota_9(h_2)^{-1})=\xi(h_2)$.
As a Bernstein eigendistribution, the generalized coefficient $S$
has nilpotent wavefront by \cite[Theorem~A]{AGS}. Applying
Theorem~\ref{m9:native-vector} to its odd part and then the
coefficient-kernel GK argument \cite[Theorem~2.21]{JLX} gives
\begin{equation}
 \dim\mathcal P(\pi)\,\dim\mathcal Q(\pi)\leq1.
 \label{m9:product-bound}
\end{equation}

\begin{proposition}[Duality in connected $\GSO_{12}$]
\label{m9:connected-duality}
For every irreducible admissible smooth complex representation of
$G_{9,v}(\F)$ one has $\pi^\alpha\simeq\pi^\vee$.
\end{proposition}

\begin{proof}
Let $G^{\mathrm{vr}}$ be the nonempty Zariski-open set of matrices
having twelve distinct eigenvalues and satisfying
$\det(g^2-\mu(g)I)\ne0$.  Its complement has Haar measure zero.
For $g$ in this set, $A_g=\F[g]$ is an etale algebra of degree twelve,
and the defining vector space is free of rank one over $A_g$.
Choose an $\F$-rational cyclic vector $v_0$.  The adjoint operation
with respect to $Q$ preserves $A_g$ and gives
$g^*=\mu(g)g^{-1}$.  Define the $\F$-linear map
$\kappa(av_0)=a^*v_0$.  With
$\ell(a)=\langle v_0,av_0\rangle$, symmetry yields
$\ell(a^*)=\ell(a)$, and hence
\[
 \langle av_0,bv_0\rangle=\ell(a^*b)
 =\ell(ab^*)=\langle\kappa av_0,\kappa bv_0\rangle.
\]
Thus $\kappa$ is an isometry and
$\kappa g\kappa^{-1}=\mu(g)g^{-1}$.

Over an algebraic closure the involution on $A_g$ pairs the twelve
primitive idempotents by $\lambda\leftrightarrow\mu/\lambda$,
with no fixed idempotents.  On their cyclic-vector basis $\kappa$
permutes exactly six pairs, so $\det\kappa=(-1)^6=1$.
Also $Q\in\SO(Q)(\F)$.  Consequently $h=\kappa Q$ belongs to the
connected group $\SO(Q)(\F)$, and
\begin{equation}
 h g^{-t}h^{-1}=g^{-1}.
 \label{m9:so-conjugator}
\end{equation}
To pass from rational conjugacy to representations, use local
integrability of Harish--Chandra characters
\cite[Theorem~1, \S3.3]{Clozel}, including dyadic fields by its
\S3.5, and character determination
\cite[Corollary~2.3.2]{Casselman}. The automorphism and inversion
identities recalled in \cite[\S7.1--7.2]{RVfactor}, together with
\eqref{m9:so-conjugator}, give
$\Theta_{\pi^\alpha}(g)=\Theta_{\pi^\vee}(g)$ on the full-measure
set $G^{\mathrm{vr}}$.  Their character distributions agree, which
proves the proposition.
\end{proof}

\begin{theorem}[Literal Model 9, all admissible representations]
\label{m9:gp1}
\label{m9:main}
Let $\F$ be a non-Archimedean local field of characteristic zero,
with arbitrary residue characteristic.  For literal split
Wan--Zhang Model~9 on connected $\GSO_{12}(\F)$, every irreducible
admissible smooth complex representation satisfies
\[
 \dim\Hom_{H_{9,v}(\F)}(\pi,\xi_{9,2})\leq1.
\]
The same holds for the outer-related second realization after
transporting the entire datum by its $\F$-defined outer automorphism.
\end{theorem}

\begin{proof}
Since $\chi_2(\alpha(h))=\xi(h)$, pulling back the second Hom space
along $\alpha$ gives
\[
 \mathcal Q(\pi)\simeq
 \Hom_H((\pi^\vee)^\alpha,\xi)\simeq\mathcal P(\pi)
\]
by Proposition~\ref{m9:connected-duality}.  The product bound
\eqref{m9:product-bound} is therefore a square bound.
\end{proof}

\subsection{Semispin multiplicity one}
\label{m9:semispin-gp1-section}
The semispin distribution theorem and connected-center duality give the following multiplicity-one result.

\begin{corollary}[Semispin multiplicity one]
\label{m9:semispin-gp1}
For the group $G_{9,s}$ and subgroup $H_{9,s}=\PGL_2\ltimes U$
above, every irreducible admissible smooth complex representation of
$G_{9,s}(F)$ satisfies
\[
 \dim\Hom_{H_{9,s}(F)}(\pi,\xi_{9,1})\leq1.
\]
The same argument applies to each nonzero coefficient $\kappa$ in
Theorem~\ref{m9:native-semispin}.
\end{corollary}

\begin{proof}
Write $G_s=G_{9,s}$, $H_s=H_{9,s}$, $\xi_s=\xi_{9,\kappa}$ and
$\alpha_s=\iota_9\circ\mathrm{inv}$. On the defining cover this
involutive $F$-automorphism is
$(s,u)\mapsto(\iota_{\Spin}(s^{-1}),u^{-1})$; it preserves both
$(c_s,1)$ and $(c_v,-1)$. Its center $Z_s$ is the split torus
$\Gm$: the coordinate of a central class $[(c_v^a c_s^b,u)]$ is
$(-1)^a u$. Thus
\[
 1\longrightarrow Z_s\longrightarrow G_s
 \xrightarrow{q}\operatorname{PSO}_{12}\longrightarrow1
\]
has connected split kernel, and is onto on $F$-points by Hilbert~90.

The connected rational conjugacy in \eqref{m9:so-conjugator}
descends through $\GSO_{12}\to\operatorname{PSO}_{12}$, also onto on
$F$-points. Its scalar-stable full-measure open gives rational inverse
conjugacy downstairs for the same pinned torus-inverting operation.
Apply the rational lifting argument of \S\ref{ep:o:sec:3} to $q$.
Explicitly, on the open where both centralizers are tori, put
$S=Z_{G_s}(g)$ and $\bar S=Z_{\operatorname{PSO}_{12}}(q(g))$.
The exact sequence $1\to Z_s\to S\to\bar S\to1$ and
$H^1(F,Z_s)=0$ give
\[
 H^1(F,S)\hookrightarrow H^1(F,\bar S).
\]
Torus inversion supplies a geometric conjugator between
$\alpha_s(g)$ and $g^{-1}$. By rational conjugacy downstairs, its cocycle in $S$ becomes a
coboundary; the injection then makes it a coboundary already in $S$.
Thus the conjugator can be chosen in $G_s(F)$. Character determination, as in
Proposition~\ref{m9:connected-duality}, proves
\begin{equation}
 \pi^{\alpha_s}\simeq\pi^\vee.
 \label{m9:semispin-duality}
\end{equation}

The two GK factors are
\[
 \mathcal P_s=\Hom_{H_s}(\pi,\xi_s),\qquad
 \mathcal Q_s=\Hom_{H_s^{\iota_9}}
 (\pi^\vee,\chi_s^{\mathrm{op}}),\qquad
 \chi_s^{\mathrm{op}}(\iota_9(h))=\xi_s(h)^{-1}.
\]
The coefficient calculation preceding \eqref{m9:product-bound}
works on this unimodular group with precisely these characters:
the inverse test action $L_{h_1}R_{\iota_9(h_2)^{-1}}$ has character
$\xi_s(h_1)\xi_s(h_2)$. The odd part of the Bernstein
eigencoefficient is Bernstein-finite, hence NW with all central
cotangent components zero. Theorem~\ref{m9:native-semispin} and
Theorem~\ref{cf:gk} give $\dim\mathcal P_s\dim\mathcal Q_s\leq1$.
Finally $\alpha_s(H_s)=H_s^{\iota_9}$ and
$\chi_s^{\mathrm{op}}(\alpha_s(h))=\xi_s(h)$, so
\[
 \mathcal Q_s\simeq
 \Hom_{H_s}((\pi^\vee)^{\alpha_s},\xi_s)
 \simeq\mathcal P_s
\]
by \eqref{m9:semispin-duality}. The product bound is therefore a square bound.
The automorphism inverts the entire connected center, as does
contragredience.
\end{proof}

\part{Exceptional \texorpdfstring{$E_7$}{E7} geometry and the distribution theorem}\label{part:4}
\section{Exceptional localization and NW proper-Levi descendants}
\label{m10:distribution-section}

\paragraph{Proof strategy}
We prove the distribution theorem by first clearing the proper-Levi
opens and then continuing through the residual Bruhat cells. Testing unipotent stabilizer characters over all Levi
parameters leaves $98$ compatible loci among $52\,416$ coarse
cells: $69$ proper-Levi, one identity and $28$ residual loci
(Proposition~\ref{m10:finite-localization}).

We first clear the opens of types $A_5+A_1$, $A_3+A_2+A_1$ and
$D_6$ by descent to balanced Ginzburg--Rallis, extended Shalika and
semispin distribution theorems. Their union supplies the predecessor
for the relative $D_5+A_1$ argument. With all four opens cleared,
we continue in descending Bruhat length. On the identity locus, a character-one transporter gives vanishing.
Among the residual cells, $24$ have intrinsic odd vanishing by
characters, transporters and opposite-weight-plane descent. The other
four form two exchanged pairs, whose ambient conormals exclude the
support of the original NW distribution once their predecessors vanish (Proposition~\ref{m10:residual-conditional}).

This gives the GK product bound. The individual bound follows in
\S\ref{ep:section}: Hermitian duality handles unitarizable
representations, and Langlands classification reduces the others to
rational duality on type-$A/D$ Levis and generic tempered $E_6$ data,
or to zero period for positive standard modules from nongeneric
$E_6$ data. The latter uses a Mackey filtration, normalized $D_5$
descent, tempered exponents and the compact-boundary operator.

\subsection{The datum and its sign operation}
\label{m10:datum}
Let $G_{10}$ be split adjoint $E_7$ over $\F$, with Bourbaki edges
$1-3-4-5-6-7$ and $2-4$. Put $I=\{2,5,7\}$ and
$P_{10}=L_{10}U_{10}$. The derived Levi is
$\SL_2^3/\mu_{2,\mathrm{diag}}$, and the diagonal leaf subgroup is
$H_{10,0}=\PGL_2$. Thus $H_{10}=H_{10,0}\ltimes U_{10}$, where
$\dim U_{10}=60$. In the physical Wan--Zhang coordinates,
\[
 \mathfrak u_{10}/[\mathfrak u_{10},\mathfrak u_{10}]
 =\F v_1\oplus\F e\oplus\Mat_2(x)\oplus\Mat_2(z),\qquad
 \xi_{10}=\psi(v_1+e+\tr x+\tr z).
\]
The six character-bearing roots are
\[
 \alpha_1,\ \alpha_3,\ \alpha_6+\alpha_7,\
 \alpha_5+\alpha_6,\ \alpha_4+\alpha_5,\ \alpha_2+\alpha_4.
\]
In particular the simple-root support is $\{1,3\}$; the coefficients on
$E_4,E_6$ are zero. The physical-to-Chevalley coefficient signs are fixed in
\eqref{re:o:eq:3.7}.

Write $\omega_{10}$ for the relative opposition element; its Weyl
action is $w_{\mathrm{long}}w_I=-w_I$. Let $r\in H_{10,0}(F)$
be the image of $\left(\begin{smallmatrix}0&1\\-1&0\end{smallmatrix}\right)$
in the diagonal leaf $\PGL_2$, so that it represents $w_I$, and put
$n=r\omega_{10}$. Then $r^2=n^2=1$ in the adjoint group. The pinned
Chevalley anti-involution and the character-preserving operation are
\[
 \iota_{10}(g)=ng^{-1}n^{-1},\qquad
 \tau_{10}(g)=b_-rg^{-1}r^{-1}b_-^{-1}.
\]
Here $b_-\in Z(L_{10})(\F)$ has square one, acts by $-1$ on all four
phase summands, and has simple-root values $-1$ on $\{1,3,4,6\}$ and
$1$ on $I$. Consequently $\tau_{10}^2=1$,
$\tau_{10}(H_{10})=H_{10}$ and $\xi_{10}\circ\tau_{10}=\xi_{10}$.
The operation $\iota_{10}$ fixes the split torus and interchanges
pinned positive and negative root generators. On Levi slices it
exchanges the two acting factors in the plus-plus convention.

On $G_{10}$ use
\[
 (h_1,h_2)g=h_1gh_2^{-1},\qquad
 \eta_{10}(h_1,h_2)=\xi_{10}(h_1)\xi_{10}(h_2)^{-1},
\]
and extend by the sign $-1$ on $\tau_{10}$. All distribution spaces below
use this convention and its explicitly transported descended form.

In the detailed coordinate calculations write $G=G_{10}$,
$H=H_{10}$, $H_0=H_{10,0}$, $P=P_{10}$, $L=L_{10}$, $U=U_{10}$,
$\xi=\xi_{10}$ and $\eta=\eta_{10}$. There $w_0=\omega_{10}$ is
relative opposition; the longest Weyl element remains
$w_{\mathrm{long}}$. Likewise $\iota$ and $\tau$ denote the two
operations just defined.

\subsubsection{The character coordinates and the anti-involution}\label{gk:o:sec:A}

\paragraph{Pinning and character support}\label{gk:o:sec:A.1}
Write $j$ for the embedding of \cite[Chapter~7, pp.~62--64]{WZ}, reserving
$\iota$ for the anti-involution. Put
\[
 J=\begin{pmatrix}0&1\\-1&0\end{pmatrix},\qquad
 \Gamma=\diag(J,J,J),\qquad D_3=E_{55}+E_{66}.
\]
The tensor basis and its dual follow the physical coordinates of this
embedding; their root-coordinate identification is given in Section
M10-P of \path{anc/model10_calculations.pdf}.
The positive generators compatible with the physical coordinates are
\[
 E_1=E_{12}\in\mathfrak{sl}_3,\qquad E_3=w_2\otimes D_3,
 \qquad(E_7,E_6,E_5,E_4,E_2)=(E_{12},E_{23},E_{34},E_{45},E_{56})
\]
in the $\mathfrak{sl}_6$ factor. For the matrix roots $F_i=E_i^t$;
for the tensor root the normalization is $F_3=-w_2^\vee\otimes D_3^\vee$.
Indeed, \cite[\S4.2.1]{Pollack} gives
\[
 [E_3,F_3]=h_3=
 \left(\tfrac13\diag(-1,2,-1),\
 \tfrac13\diag(-1,-1,-1,-1,2,2)\right).
\]
Its weights on $E_3,E_1,E_4$ are $2,-1,-1$, and the other simple
weights are zero. This determines the sign of the tensor generator.

Root addition gives the ten-dimensional
abelianization in \S\ref{m10:datum}. With
$\xi(\exp X)=\psi(\lambda_\xi(X))$, the phase is
$\lambda_\xi=v_1+e+\tr x+\tr z$; hence its root support is the six roots
listed there and $\lambda_\xi(E_4)=\lambda_\xi(E_6)=0$.

\paragraph{The Chevalley operation}\label{gk:o:sec:A.2}
The element $w_0$ in \cite[Chapter~7]{WZ} has square one, centralizes
$H_0$, normalizes $L$ and takes $U$ to $U^-$. Its Weyl action is
$w_{\mathrm{long}}w_I=-w_I$, with $w_I=s_2s_5s_7$. For
$r=j(I_3,\Gamma,1)$ and $n=rw_0$ one has $r^2=n^2=1$ in the
adjoint group. The differential of $\iota(g)=ng^{-1}n^{-1}$ is
transpose on the two matrix factors and
\[
 v\otimes X\longmapsto -v^\vee\otimes X^t
\]
on the tensor factor, with the corresponding reverse map on its dual.
Here $T$ is the split maximal torus determined by the pinning.
It therefore fixes $T$ and exchanges the pinned $E_i,F_i$.

\paragraph{Character inversion and its correction}\label{gk:o:sec:A.3}
The operation $\tau_0=\Ad(w_0)\iota=\Ad(r)\mathrm{inv}$ acts by
\[
 (v_1,x,z,e)\longmapsto(-v_1,-JxJ^{-1},-JzJ^{-1},-e),
 \qquad\xi\circ\tau_0=\xi^{-1}.
\]
\label{gk:o:sec:A.4}
Set
\[
 b_-=a_{\mathrm{WZ}}(-1)=
 j\left(\diag(-1,1,1),\diag(-I_2,I_2,-I_2),-1\right).
\]
This element is central in $L$, has square one, commutes with $r$, and
acts by $-1$ on each phase summand. Thus
$\tau=\Ad(b_-r)\mathrm{inv}$ has square one and acts by
\[
 (v_1,x,z,e)\longmapsto(v_1,JxJ^{-1},JzJ^{-1},e).
\]
It preserves $H$ and $\xi$. On $H_0$ it is projective transpose, since
$Jg^{-1}J^{-1}=g^t/\det g$. These identities on the abelianization
prove character invariance on all of $U$.

\paragraph{Central scaling}\label{supp:e7-scaling}
For general $t$, the map $a_{\mathrm{WZ}}(t)$ scales $e$ by $t^{-1}$.
The cocharacter
\[
 a_{\mathrm{corr}}(t)=j\left(\diag(t^3,t^2,t^{-8}),
                  \diag(tI_2,I_2,t^{-1}I_2),t\right)
\]
scales all four phase summands by $t$: substitution verifies the
determinant constraint and each weight. Its value at $-1$ is $b_-$.

\paragraph{The Levi quotient, including its central form}\label{gk:o:sec:D.1}

Wan--Zhang's parametrization is
\[
\ell=
j\left(
\operatorname{diag}(a,b,c),
\operatorname{diag}(g_1,g_2,g_3),
\lambda
\right),
\qquad
\lambda^3abc\prod_i\det g_i=1.
\]
The kernel of \(j\) is a split two-dimensional torus. Restricting to the three \(\SL_2\) factors gives \(\mu_{2,\mathrm{diag}}\), proving
\[
[L,L]\simeq \SL_2^3/\mu_{2,\mathrm{diag}}.
\]

Its action on the abelianization is
\[
(v,e,x,z)\longmapsto
\left(
\frac abv,\;
b\lambda\det(g_3)e,\;
g_1xg_2^{-1},\;
g_2zg_3^{-1}
\right).
\]
Using coefficient matrices paired by \(\operatorname{tr}(A^tx)\), define
\[
p(\ell)=\lambda_\xi\circ\operatorname{Ad}(\ell)
=
\left(
\frac ab,\;
b\lambda\det g_3,\;
(g_2^{-1}g_1)^t,\;
(g_3^{-1}g_2)^t
\right).
\]
This identifies
\[
H_0(F)\backslash L(F)
\simeq
\mathcal O:=(F^\times)^2\times \GL_2(F)^2.
\]

An explicit rational section is obtained as follows: for
\((s,t,A,B)\in\mathcal O\), take
\[
\lambda=1,\quad g_3=I,\quad g_2=B^t,\quad g_1=B^tA^t,
\]
\[
b=t,\quad a=st,\quad
c=(st^2\det A(\det B)^2)^{-1}.
\]
The stabilizer of \((1,1,I,I)\) is exactly Wan--Zhang's \(H_0\).

\subsection{The finite sieve and the four maximal patches}
\label{m10:localization}
The coarse relative positions are indexed by
$W_I\backslash W(E_7)/W_I$; compatibility is tested on the full
Levi-parameter family.
\begin{proposition}[Finite localization]
\label{m10:finite-localization}
The relative Bruhat decomposition has $52\,416$ coarse cells. Exactly
$98$ have nonempty full character-compatibility loci, with decomposition
\[
 98=69\ \text{proper-Levi cells}+1\ \text{identity cell}
       +28\ \text{residual }E_7\text{ cells}.
\]
There are only $23$ compatible chosen Tits representatives: $75$ additional
cells become compatible after varying their Levi parameters. The empty
loci are eliminated by unipotent stabilizer-character incompatibility in
the Bruhat support filtration.
\end{proposition}

The complete classification and Hasse diagram are provided in the
electronic ancillary materials, in
\path{anc/model10_gk_exhaustive/model10_gk_bundle.zip.xz}.
Section~\ref{supp:e7-enumeration} describes the enumeration and closure convention.

To test compatibility over the full coefficient family, we use linear
matching on the phase space
$\mathcal P=(\F^\times)^2\times\GL_2(\F)^2$. For a chosen lift $n_w$,
the equations are
\[
 c(w,\beta)p(e_{w\beta})=q(e_\beta),\qquad
 \beta\in\Phi(U_{10})\cap w^{-1}\Phi(U_{10}),
\]
followed by nonvanishing of the scalar coordinates and four matrix
determinants. A failed equation gives an unipotent stabilizer with
nontrivial character and trivial modular character. Bruhat closure
and its finite closed unions supply the topology for the support induction.

For $j\in\{1,3,4,6\}$, set $K_j=\Delta\setminus\{j\}$,
$M_{10,j}=M_{K_j}$ and
\[
 O_j=Q_{K_j}\omega_{10}Q_{K_j}
    =N_jM_{10,j}\omega_{10}N_j.
\]
In the detailed patch calculations, write $M_j=M_{10,j}$,
$H_j=H_{10}\cap M_j$, $U_j=U_{10}\cap M_j$, and $b=b_-$.
The extra adjoint cocharacter is
$v=\tfrac12(H_2+H_5+H_7)$, as in \S\ref{m9:o:sec:2.3}.
We cover the compatible proper-Levi loci by four opens and treat
the residual loci in the later Bruhat induction. The table records
both geometric coverage and a disjoint assignment for counting;
the support order remains the Bruhat order.
\[
\begin{array}{c|c|r|r|l}
 \text{patch}&\text{derived type}&\text{in patch}&\text{assigned}&\text{argument}\\\hline
 O_3&A_5+A_1&26&26&\text{balanced GR and Whittaker boundary}\\
 O_4&A_3+A_2+A_1&12&6&\text{balanced extended Shalika}\\
 O_1&D_6&48&32&\text{NW semispin theorem}\\
 O_6&D_5+A_1&18&5&\text{relative B5 theorem}
\end{array}
\]
Every compatible proper cell satisfies $ww_{\mathrm{long}}w_I\in W_{K_j}$
for at least one patch. In terms of the relative opposition this is
$w\omega_{10}\in W_{K_j}$. The four opens cover
\[
 \Omega_{\mathrm{proper}}
 =\bigcup_{I\subseteq K\subsetneq\Delta}Q_K\omega_{10}Q_K.
\]

The lift signs and the Bruhat closure convention are fixed in
\S\ref{gk:o:sec:C.1}--\ref{gk:o:sec:C.3}.
\paragraph{Reducing \(C_w\) modulo \(H_0\times H_0\)}\label{gk:o:sec:D.2}

Here $C_w$ is the locus of Levi pairs satisfying stabilizer-character
compatibility in this relative position. Its defining root equations
are derived below.
For
\[
x=\ell_1n_w\ell_2,
\]
substitute \(X=\operatorname{Ad}(\ell_2^{-1})Y\), and set
\[
p=p(\ell_1),\qquad q=p(\ell_2^{-1}).
\]
Then
\[
{
\frac{C_w}{H_0\times H_0}
=
\left\{
(p,q)\in\mathcal O^2:
p(\operatorname{Ad}(n_w)Y)=q(Y)
\quad
\forall Y\in
\mathfrak u\cap\operatorname{Ad}(n_w^{-1})\mathfrak u
\right\}.
}
\]

Write
\[
\operatorname{Ad}(n_w)e_{\beta}
=c(w,\beta)e_{w\beta}.
\]
The complete equations are
\[
{
c(w,\beta)p(e_{w\beta})=q(e_{\beta}),
\qquad
\beta\in\Phi(U)\cap w^{-1}\Phi(U).
}
\]

Because both functionals vanish on \([\mathfrak u,\mathfrak u]\), each equation is either zero, a forced-zero coordinate, or a signed match between one left and one right coordinate. Each coefficient occurs in at most one matching equation.

Thus the initial compatibility problem over the phase space is a linear coordinate-matching problem in \(F^{20}\), followed by the open conditions
\[
s_Lt_Ls_Rt_R
\det A_L\det B_L\det A_R\det B_R\ne0.
\]

\paragraph{Nonemptiness test}\label{gk:o:sec:D.3}

The locus is empty exactly when either a required scalar is forced to zero, or some \(2\times2\) matrix has both determinant monomials blocked by forced-zero entries.

For sufficiency, observe that the matching equations identify only
distinct left/right coordinates. An unblocked determinant is therefore
a nonzero polynomial in the remaining independent variables, and the
product of these determinants has a nonvanishing point over the
infinite field \(F\).

The result is:

\begingroup\small\setlength{\tabcolsep}{3pt}
\begin{longtable}{@{}>{\raggedright\arraybackslash}p{0.470\linewidth}>{\raggedright\arraybackslash}p{0.470\linewidth}@{}}
\toprule
Test & Number of cells \\
\midrule\endhead
Chosen Tits representative passes & \textbf{23} \\
\(C_w\) is nonempty & \textbf{98} \\
Representative fails but \(C_w\ne\varnothing\) & \textbf{75} \\
\(C_w\) is empty & \textbf{52,318} \\
\bottomrule\end{longtable}\endgroup

The exact cell IDs, reduced words, integral witnesses and exclusion
certificates are provided in the electronic ancillary materials at
\path{anc/model10_gk_exhaustive/model10_gk_bundle.zip.xz}.

\paragraph{An explicit representative failure with nonempty \(C_w\)}\label{gk:o:sec:D.4}

Cell 2526 has length \(14\) and reduced word
\[
w=s_6s_4s_5s_3s_4s_1s_3s_2s_4s_5s_6s_1s_3s_4.
\]
Its equations are
\[
(A_L)_{12}=(B_L)_{12}=(A_R)_{12}=(B_R)_{12}=0,
\]
\[
s_R=(B_L)_{22},\quad t_R=t_L,\quad
(A_R)_{11}=-(A_L)_{11},
\]
\[
(B_R)_{11}=-(B_L)_{11},\qquad
(B_R)_{22}=-s_L.
\]

The base pair $(1,1,I,I)$ fails, whereas the integral witness
in Section~2, cell $2526$, of \path{anc/model10_calculations.pdf} satisfies all
these equations and has determinants $21,33,-3,121$. The quotient
locus has dimension $11$. This illustrates why compatibility must
be tested on the full coefficient family.

\paragraph{The explicit incompatibility stabilizer}\label{gk:o:sec:D.5}

If an equation fails at \(x=\ell_1n_w\ell_2\), choose a witnessing root \(\beta\) and put
\[
X=\operatorname{Ad}(\ell_2^{-1})e_{\beta}.
\]
Then
\[
\left(
\exp(t\operatorname{Ad}(x)X),\,
\exp(tX)
\right)
\]
fixes \(x\), and its character is
\[
\psi\!\left(
t\bigl[\lambda_\xi(\operatorname{Ad}(x)X)-\lambda_\xi(X)\bigr]
\right).
\]
For a suitable \(t\) the character is nontrivial. Since the
stabilizer is unipotent and has trivial modular character,
\cite[Theorem 2.5 and Remark 2.6]{JLX} give zero equivariant distributions on each of the 52,318 empty-\(C_w\) strata. 

These intrinsic vanishings enter the ambient argument through the
Bruhat support filtration and its distribution exact sequences.

\subsection{Frobenius descent on the proper-Levi opens}
\label{m10:free-descent}
Let $T\in\NW(G_{10})$ be the original odd equivariant ambient distribution.
Using the product coordinates on $O_j$, extract two character-weighted
Haar lines to obtain a descendant $S_j$ on $M_{10,j}(\F)$. With
$H_{10,j}=H_{10}\cap M_{10,j}$, put $\xi_{10,j}=\xi_{10}|_{H_{10,j}}$.
The initial slice involution is $\Ad(b_-)\iota_{10}$. Changing the right
variable by $k=\tau_{10}(h_2^{-1})$ and writing $m=xb_-$ gives exactly
\[
 (h,k)x=hx\iota_{10}(k),\qquad
 \eta_j(h,k)=\xi_{10,j}(h)\xi_{10,j}(k),\qquad
 \iota_{10}S_j=-S_j.
\]
Both Haar lines have trivial $H_{10,j}$ density character: the determinant
is an algebraic character, trivial on the unipotent radical and on
$\PGL_2$. Since the involution has square one and a positive density scalar,
that scalar is one. Free descent gives the wavefront bound
\[
 \WF(S_j)\subset\cN(\mathfrak m_{10,j}^*),
 \quad\text{with zero component on }\mathfrak z(\mathfrak m_{10,j})^*.
\]
Thus each descendant has the equivariance, oddness and nilpotent
wavefront needed below. Injectivity of free descent then transfers
its vanishing to the original patch.

\paragraph{Free unipotent descent and the right action}\label{pl:o:sec:1.2}

For each of the four $K=K_j$, the specified factorization is an analytic isomorphism
\begin{equation}
N_K\times M_K\times N_K\longrightarrow O_K,
\qquad (u,m,v)\longmapsto u m w_0v.
\label{pl:o:eq:1.4}
\end{equation}
Since $H=(H\cap M_K)\ltimes N_K$, extracting the two
character-weighted Haar factors gives a distribution $S_K$ on $M_K$
by Frobenius descent through \eqref{pl:o:eq:1.4}.

Put $H_K=H\cap M_K$ and $\theta=\Ad(w_0)$. In the $m$-coordinate the original $h_1gh_2^{-1}$ action becomes
\[
m\longmapsto h_1m\theta(h_2)^{-1}.
\]
The induced anti-involution is
\begin{equation}
\sigma_K=\operatorname{Ad}(b)\iota|_{M_K},
\qquad \tau(mw_0)=\sigma_K(m)w_0.
\label{pl:o:eq:1.5}
\end{equation}
To convert to the plus-plus convention, let $\tau_H$ be the
restriction of $\tau$ to $H_K$. Since $\sigma_K\tau_H=\theta$, setting
\[
k=\tau_H(h_2^{-1})
\]
gives $\theta(h_2)^{-1}=\sigma_K(k)$ and $\xi(k)=\xi(h_2)^{-1}$. Consequently the descended action and character are
\begin{equation}
(h,k)\cdot m=h m\sigma_K(k),\qquad
\eta_K(h,k)=\xi_K(h)\xi_K(k).
\label{pl:o:eq:1.6}
\end{equation}
The odd involution acts by $\sigma_K$ and has sign $-1$.

\paragraph{Right translation by $b$}\label{pl:o:sec:1.3}

Set $m=xb$. Since $\iota(b)=b$ and $b^2=1$,
\begin{equation}
\sigma_K(xb)=\iota(x)b,\qquad
h(xb)\sigma_K(k)=h x\iota(k)b.
\label{pl:o:eq:1.7}
\end{equation}
Thus after right translation by b the problem is exactly
\begin{equation}
(h,k)\cdot x=h x\iota(k),\qquad
\eta_K(h,k)=\xi_K(h)\xi_K(k),\qquad \iota\text{ has sign }-1.
\label{pl:o:eq:1.8}
\end{equation}
In the matrix covers constructed below, $\iota$ lifts to ordinary transpose in every $\GL$ factor.

\paragraph{Density}\label{pl:o:sec:1.4}

For $h\in H_K$, the determinant of $\Ad(h)$ on $\Lie(N_K)$ is an algebraic character of $H_K$. It is trivial on its unipotent radical and on $H_0\simeq\PGL_2$, which has no nontrivial algebraic characters. Hence
\begin{equation}
|\det(\operatorname{Ad}(h)|\mathfrak n_K)|_F=1.
\label{pl:o:eq:1.9}
\end{equation}
The same calculation applies to the second N-factor, so both weighted
Haar lines have trivial relative density character. Since the
involution exchanges them and has square one, its positive density
scalar is one.

\subsubsection{The Levi groups, including their centers}\label{pl:o:sec:2}

Let $D_j=[M_j,M_j]$. The derived cocharacter lattice is $Q_{K_j}^\vee+\mathbb Z v$. Therefore
\begin{equation}
\begin{array}{c|c|c|c}
j&D_j&\dim M_j&\dim U_j\\ \hline
1&\operatorname{HSpin}_{12}&67&27\\
3&(\SL_6/\mu_2)\times \SL_2&39&13\\
4&((\SL_4\times \SL_2)/\mu_{2,\Delta})\times \SL_3&27&7\\
6&(\operatorname{Spin}_{10}\times \SL_2)/\mu_{2,\Delta}&49&18.
\end{array}
\label{pl:o:eq:2.1}
\end{equation}
In each row $H_j=\PGL_2\ltimes U_j$, of dimension $3+\dim U_j$. The diagonal $\mu_2$ in the last row is generated by $(-1_{\Spin_{10}},-I_2)$. In the $j=4$ row it is generated by $(-I_4,-I_2)$. In the $j=3$ row it is the scalar $-I_6$ in the $\SL_6$ factor.

\paragraph{Rank-seven root data}\label{pl:o:sec:2.1}

The root datum of $M_j$ is
\begin{equation}
X^*(T)=Q(E_7),\quad X_*(T)=P^\vee(E_7),\quad
\Phi(M_j)=\Phi(E_7)\cap\mathbb Z K_j,
\label{pl:o:eq:2.2}
\end{equation}
with the corresponding inherited coroots and pinning. Since the group is split, \eqref{pl:o:eq:2.2} determines its F-isomorphism class and its rational-point group.

The center $Z(M_j)$ is a connected split $\mathbb G_m$, since
$Q(E_7)/\mathbb ZK_j$ is freely generated by the omitted simple root. Its primitive cocharacter is $z_j=\omega_j^\vee$. The primitive character of $M_j$ with kernel $D_j$ is $\chi_j=\omega_j$; for $j=1,3,4,6$ these fundamental weights lie in $Q(E_7)$. Directly from the inverse Cartan matrix,
\begin{equation}
\chi_j(z_j(t))=t^{d_j},\qquad
(d_1,d_3,d_4,d_6)=(2,6,12,4).
\label{pl:o:eq:2.3}
\end{equation}
Thus $Z(D_j)\simeq\mu_{d_j}$ and the central product describing $M_j$ is
\begin{equation}
M_j\simeq(D_j\times\mathbb G_m)/
\{(z_j(\zeta)^{-1},\zeta):\zeta^{d_j}=1\}.
\label{pl:o:eq:2.4}
\end{equation}
Here $z_j(\zeta)$, for $\zeta\in\mu_{d_j}$, lies in $D_j$. Formula \eqref{pl:o:eq:2.4}, together with \eqref{pl:o:eq:2.1}, specifies the finite central gluing.

\paragraph{Rational cosets of the central product}\label{pl:o:sec:2.2}

Let $q_j(d,t)=dz_j(t)$. Then
\begin{equation}
q_j(D_j(F)\times F^\times)
=\{m\in M_j(F):\chi_j(m)\in(F^\times)^{d_j}\}.
\label{pl:o:eq:2.5}
\end{equation}
For the reverse inclusion, choose $t$ with $t^{d_j}=\chi_j(m)$; then $mz_j(t)^{-1}\in D_j(F)$. Moreover $s_j(a)=\alpha_j^\vee(a)$ has $\chi_j(s_j(a))=a$. Therefore
\begin{equation}
M_j(F)/q_j(D_j(F)\times F^\times)
\simeq F^\times/(F^\times)^{d_j},
\label{pl:o:eq:2.6}
\end{equation}
with representatives $s_j(a)$. These are finite groups for the fields under consideration, including residue characteristic 2 or 3.

This also identifies the cohomological obstruction explicitly. Under the kernel identification in \eqref{pl:o:eq:2.4}, the connecting map is
\begin{equation}
\delta_j:M_j(F)\longrightarrow H^1(F,\mu_{d_j})
\simeq F^\times/(F^\times)^{d_j},\qquad
\delta_j(m)=[\chi_j(m)].
\label{pl:o:eq:2.6a}
\end{equation}
Indeed a lift over an algebraic closure is obtained by choosing a $d_j$-th root $t$ of $\chi_j(m)$, and its Galois cocycle is $(z_j(\gamma(t)/t)^{-1},\gamma(t)/t)$. Formula \eqref{pl:o:eq:2.5} identifies the kernel, and $s_j(a)$ proves surjectivity of this particular connecting map.

Because $H_j$ is contained in $D_j$, each coset in
\eqref{pl:o:eq:2.6} is stable under the two-sided $H_j$ action.
It is also preserved by the anti-involution, which fixes $\chi_j$ and
$s_j(a)$. In the coordinate $m=q_j(x)s_j(a)$, its lift is
\begin{equation}
\widehat\iota_{j,a}(d,t)
=\bigl(\operatorname{Ad}(s_j(a))\iota_{D_j}(d),\ t\bigr).
\label{pl:o:eq:2.7}
\end{equation}
This involution fixes the kernel in \eqref{pl:o:eq:2.4} pointwise.

For $j=3,4$ we use the connected-kernel matrix covers of
\S\ref{pl:o:sec:7} and \S\ref{pl:o:sec:8}. For $j=1,6$ the
distribution arguments apply separately in every rational sector.

\paragraph{Finite central covers}\label{pl:o:sec:2.3}

Suppose $q:\widetilde G\to M$ is a finite central isogeny whose $F$-image $J$ is open of finite index, $H(F)\subset J$, and $q^{-1}(H)(F)\to H(F)$ is surjective. Suppose the anti-involution preserves each coset $Js$, and its transported anti-involution on $J$ lifts to $\widetilde G$. Then distributions on $Js$ identify, by translation and étale pullback, with distributions on $\widetilde G(F)$ invariant under $\ker(q)(F)$. Under this identification:

\begin{itemize}
\item the acting subgroup is the \textbf{full} preimage $q^{-1}(H)(F)$, and its character is the pullback of the original character;
\item the involution is the transported, lifted involution;
\item pullback is injective on the coset and preserves the nilpotent-wavefront condition;
\item the finite deck density has trivial character and involution scalar one.

\end{itemize}
Consequently vanishing for all these lifted problems implies vanishing on $M(F)$. Conversely, a kernel-invariant lifted distribution descends locally and then globally to its coset, so the equivalence is exact with these hypotheses.

Proof: a finite isogeny in characteristic zero is an analytic local isomorphism; its rational fibers over $J$ are precisely the $\ker(q)(F)$-orbits. Use local inverse charts and finite averaging. Equivariance and the involution identities follow from the intertwining maps. On Lie algebras the differential is an isomorphism preserving Jordan decomposition, so the wavefront condition is unchanged. Finally sum over the finitely many clopen cosets.

For \eqref{pl:o:eq:2.4}, $q^{-1}(H_j)(F)$ consists exactly of
\[
(hz_j(\zeta)^{-1},\zeta),\quad h\in H_j(F),\quad\zeta\in\mu_{d_j}(F),
\]
so the subgroup-surjectivity hypothesis is satisfied. 

\subsubsection{Nilpotent wavefront on the Levi descendants}\label{pl:o:sec:3}

We use $\NW(M)$ in the sense of Definition~\ref{cf:native}: the derived component is nilpotent and the central cotangent component is zero.

\paragraph{Descent lemma}\label{pl:o:sec:3.1}

\textbf{Lemma} Under \eqref{pl:o:eq:1.4}, if the ambient $T$ belongs to $\NW(G)$, its descendant belongs to $\NW(M_K)$. Conversely, reconstruction from a distribution in $\NW(M_K)$, with the two smooth character-weighted N-factors, has nilpotent wavefront on the open patch $O_K$.

\textbf{Proof.} Open restriction preserves the assumed NW condition.
To apply Lemma~\ref{cf:gauss}, we check the cotangent identification
in the present coordinates. A character of a p-adic unipotent group is locally constant, so each outer character-weighted Haar factor is smooth and nonvanishing. In the product coordinates its wavefront contribution is zero.

At $g=mw_0$, left translation identifies the first $N$-tangent with $\mathfrak n_K^-$, the last $N$-tangent with $\mathfrak n_K$, and the $M$-tangent with $\Ad(w_0^{-1})\mathfrak m_K$. With an invariant nondegenerate form,
\[
(\mathfrak n_K\oplus\mathfrak n_K^-)^\perp=\mathfrak m_K.
\]
Thus a descendant cotangent vector is embedded in $\mathfrak g^*$ by the Levi inclusion, followed by $\Ad(w_0)$. An element of $\mathfrak m_K$ is nilpotent in $\mathfrak g$ if and only if it is nilpotent in $\mathfrak m_K$: Jordan decomposition in a Levi agrees with Jordan decomposition in the ambient reductive group. This proves the assertion at the slice. Left and right unipotent translations transport it throughout the patch and preserve nilpotence. The same calculation proves the converse reconstruction statement. \hfill\(\square\)

For an ordinary Gauss patch $NMN^-$, the proof is identical without $w_0$. This version will be used inside $\GL_6$ in \S\ref{pl:o:sec:5}.

\paragraph{Central-torus covers}\label{pl:o:sec:3.2}

For the connected central-torus quotients used below, $q$ is a
submersion and $q(F)$ is surjective. Generalized-function pullback along $q$ is therefore defined and
injective. Its cotangent map identifies a nilpotent cotangent vector downstairs with a nilpotent cotangent vector upstairs, with zero component along the extra central kernel directions. This follows directly from the derived-Lie-algebra isomorphism and the central nature of the kernel. The kernel Haar density and its conjugation character are trivial.

The lower-group theorems apply to the NW condition just proved.

\subsection{The two type-A arguments}
\label{m10:type-a}
The matrix covers in \S\ref{pl:o:sec:7} and
\S\ref{pl:o:sec:8} have connected split central-torus kernels and
are surjective on $F$-points. Their subgroup preimages, characters and
transpose operations are determined there from the inherited pinning.

\begin{proposition}[The type-A patch theorem]
\label{m10:type-a-patches}
Odd equivariant distributions satisfying NW for these two
data vanish. In particular the original $T$ satisfies
$T|_{O_3}=T|_{O_4}=0$.
\end{proposition}

We reduce these patches by $\GL_2$ and $\GL_3$ Gauss descent to
balanced GR and balanced extended Shalika, respectively. Since the
passive torus directions have zero wavefront,
Lemma~\ref{pl:central-balancing} supplies the balanced distribution
theorems. The small Whittaker boundary lemma then clears the
complements of the Gauss opens. We give the
exact covers and density calculations below.

\subsubsection{A representation-valued Shalika distribution lemma}\label{pl:o:sec:5}
The ambient group and graph subgroup are $X_{\rm Sh}=\GL_4\times\GL_2$
and $R_{\rm Sh}$ of Theorem~\ref{cf:shalika-nw}, with phase
$\epsilon=\psi(\tr X)$ and product transpose.

\paragraph{Shalika--NW}\label{pl:o:sec:5.1}
\begin{theorem}[Representation-valued Shalika]\label{pl:shalika-nw}
For every smooth unitary $\chi$, the odd
$R_{\rm Sh}\times R_{\rm Sh}$-equivariant distributions on $X_{\rm Sh}$,
with character $\epsilon(r_1)\epsilon(r_2)\chi(\det a_1\det a_2)$
and nilpotent wavefront, are zero.
\end{theorem}
\begin{proof}
Apply Theorem~\ref{cf:shalika-nw}, whose proof uses the $4+2$
Gauss reconstruction and the open-local GR support argument. To
compare with the balanced statement below, recall its intrinsic
rank-zero conclusion:
\begin{equation}
 XC+BR=0,\qquad R=-B^{-1}XC,
 \label{pl:o:eq:5.1}
\end{equation}
\begin{equation}
 \psi(\tr X+\tr R)=\psi(\tr((I_2-CB^{-1})X)),
 \label{pl:o:eq:5.2}
\end{equation}
which forces $C=B$ and leaves, after unipotent normalization,
\begin{equation}
 \begin{pmatrix}0&B\\B&0\end{pmatrix}.
 \label{pl:o:eq:5.3}
\end{equation}
Rational similarity of $B$ to $B^t$ supplies the character-one
transporter. The ambient-group, density, and NW hypotheses
are those proved in the common theorem.
\end{proof}

\paragraph{Extended and balanced Shalika}\label{pl:o:sec:5.2}

Enlarge the acting subgroup to
\begin{equation}
R_{\rm eSh}=\{(cD_2(a)u,a):c\in F^\times,\ a\in \GL_2,\ u\in U_{2,2}\},
\qquad \nu(c,a)=c\det a.
\label{pl:o:eq:5.4}
\end{equation}
The phase remains $\psi(\tr X)$. An odd NW distribution equivariant under $R_{\rm eSh}$ on both sides with the specified character vanishes by Theorem~\ref{pl:shalika-nw}: restrict its equivariance to $c=1$, where the auxiliary character is $\chi(\det a)$.

Set
\[
p(g,y)=\det g\,(\det y)^2.
\]
Then
\begin{equation}
p(r_1(g,y)\sigma(r_2))=
\bigl(c_1c_2\det a_1\det a_2\bigr)^4p(g,y).
\label{pl:o:eq:5.5}
\end{equation}
The central subgroup $(tI_4,I_2)$ lies in $R_{\rm eSh}$, has
$\nu=t$, and changes $p$ by $t^4$. It supplies a section of the
quotient character used in Lemma~\ref{pl:central-balancing}.
The group is unimodular: $c$ is central and $a$ acts by conjugation
on $\Mat_2$ with determinant one. Quotient averaging and Mellin
inversion in that lemma give the balanced statement.

Lemma~\ref{pl:central-balancing}, now with $k=4$, proves
\begin{equation}
{\text{balanced extended-Shalika odd NW distributions vanish}.}
\label{pl:o:eq:5.6}
\end{equation}
This is the second lower-model theorem needed below.

\subsubsection{Removing the $\GL_2$ and $\GL_3$ Whittaker boundaries}\label{pl:o:sec:6}

\begin{lemma}[Small Whittaker boundaries]\label{pl:whittaker-boundary} Let $n=2$ or $3$. Let $N$ be the upper maximal unipotent subgroup of $\GL_n$, with character $\psi$ composed with the sum of the simple upper entries. Let Y be a reductive group or a passive analytic parameter space with a specified closed cotangent cone. Suppose $S$ on $\GL_n\times Y$ is equivariant for N-left and Nᵗ-right with this character on both sides and its $\GL_n$ wavefront component is nilpotent. If S vanishes on
\[
(N T N^t)\times Y,
\]
then S=0.

\end{lemma}

\textbf{Proof.} Use the $\GL_2$ or $\GL_3$ opposite-Bruhat cells $BwB^-$, with permutation representatives w. At $tw$, a stabilizer root satisfies
\[
R=-w^{-1}t^{-1}Ltw.
\]
For an inversion root $E_{ik}$, character compatibility compares whether $E_{ik}$ is simple with whether its image is a negative simple root. A simple/composite mismatch makes the cell incompatible; if both are simple, compatibility requires $t_i=t_k$.

The complete list is:

\begingroup\small\setlength{\tabcolsep}{3pt}
\begin{longtable}{@{}>{\raggedright\arraybackslash}p{0.235\linewidth}>{\raggedright\arraybackslash}p{0.235\linewidth}>{\raggedright\arraybackslash}p{0.235\linewidth}>{\raggedright\arraybackslash}p{0.235\linewidth}@{}}
\toprule
n & Representative & Compatibility condition & Consequence after preceding opens vanish \\
\midrule\endhead
2 & identity & none & Gauss open \\
2 & $s_1$ & $t_1=t_2$ & Nonzero diagonal semisimple conormal \\
3 & identity & none & Gauss open \\
3 & $s_1$ & $t_1=t_2$ & Nonzero diagonal semisimple conormal \\
3 & $s_2$ & $t_2=t_3$ & Nonzero diagonal semisimple conormal \\
3 & $s_1s_2$ & simple/composite mismatch & Ordinary character vanishing \\
3 & $s_2s_1$ & simple/composite mismatch & Ordinary character vanishing \\
3 & longest element & $t_1=t_2=t_3$ & Nonzero traceless diagonal semisimple conormal \\
\bottomrule\end{longtable}\endgroup

For the conormal calculation, consider $g=tw$, all left-translated $N\times N^t$ orbit tangent vectors are root-space vectors, and the torus tangents on the compatible submanifold satisfy the displayed equalities. For a relation $t_i=t_k$, the diagonal element $w^{-1}(E_{ii}-E_{kk})w$ annihilates both sets of tangents. This element is nonzero semisimple. For the longest $\GL_3$ cell, any nonzero traceless diagonal element works. Give this cotangent vector zero Y-component. Unipotent transport preserves its semisimplicity.

The compatible submanifold is smooth: in the ordinary Bruhat coordinates it is a product of the fixed unipotent coordinate spaces and the smooth subtorus defined by the stated equalities. Process the cells by increasing length. The union of all shorter opposite-Bruhat cells with the cell currently considered is an open subset of $\GL_n$. Once the shorter pieces vanish, the restriction to this open is supported in the current cell; ordinary character vanishing first confines that support to the compatible submanifold, which is closed in the open under consideration. The conormal theorem \cite[Theorem 2.19]{JLX} then contradicts the nilpotent-wavefront hypothesis. This proves the compatible-cell vanishing. The incompatible cells use only the ordinary stabilizer criterion.

At each step the conormal theorem applies to the original distribution restricted to the current ambient open. \hfill\(\square\)

\paragraph{The Gauss descendant and passive torus ratios}\label{pl:o:sec:6.2}

Since the $N\times N^t$ action on $NTN^t$ is free, its
Frobenius descendant lies on $T\times Y$. A torus cotangent vector, embedded in $\mathfrak{gl}_n^*$ at the diagonal representative, is diagonal; nilpotence forces it to be zero. Thus the descended wavefront is contained in
\begin{equation}
\{0\}\text{ on }T\quad\times\quad\text{the allowed cone on }Y.
\label{pl:o:eq:6.1}
\end{equation}
Normalizing a common scalar torus coordinate below uses Frobenius
descent and Lemma~\ref{cf:frobenius-wf}. The restriction of the
displayed cotangent vector to the fiber remains zero in the torus-ratio
directions and retains its allowed $Y$-component. Thus the remaining
ratios are passive. The Gauss density depends only on ratios of
diagonal entries, so its relative character in this normalization is trivial.

\subsubsection{The A5+A1 patch: quotient and closure}\label{pl:o:sec:7}

\paragraph{The rational group and an F-surjective cover}\label{pl:o:sec:7.1}

Let $A\in\GL_2$ and $g\in\GL_6$. There is a pinned central quotient
\begin{equation}
{
M_3\simeq
(\GL_2\times \GL_6)/
\{(z^{-2}I_2,zI_6):z\in\mathbb G_m\}.
}
\label{pl:o:eq:7.1}
\end{equation}
The kernel is a connected split torus.

To specify the F-form directly on the torus, take $A=\diag(A_1,A_2)$ and $g=\diag(B_1,\ldots,B_6)$. The simple E7 root values of the image are
\begin{equation}
\begin{aligned}
\alpha_1&=A_1/A_2,&\alpha_2&=B_5/B_6,&\alpha_3&=A_2B_5B_6,\\
\alpha_4&=B_4/B_5,&\alpha_5&=B_3/B_4,&\alpha_6&=B_2/B_3,&\alpha_7&=B_1/B_2.
\end{aligned}
\label{pl:o:eq:7.2}
\end{equation}
The $\SL_{2}$ root maps to $\alpha_{1}$, and the $\SL_{6}$ simple-root chain maps to ($\alpha_{7}$,$\alpha_{6}$,$\alpha_{5}$,$\alpha_{4}$,$\alpha_{2}$). These torus and root-group maps are compatible with the inherited pinning and determine the homomorphism. Their only central kernel is that in \eqref{pl:o:eq:7.1}. Any seven nonzero simple-root values in the adjoint torus have a preimage in \eqref{pl:o:eq:7.2}: first set $B_6=1$ and solve successively for $B_5,\ldots,B_1$, then for $A_2$ and $A_1$. Root-group lifting and the split Bruhat decomposition then give surjectivity on $M_{3}(F)$. Equivalently, it follows from Hilbert 90 for the split central kernel.

The induced cocharacter matrix is integral and surjective; all seven nonzero Smith invariants are one. The centers of $\GL_{2}$ and $\GL_{6}$ map to $z_3$ and $2z_3$, respectively. The primitive character $\chi_3$ pulls back to
\begin{equation}
(\det A)^3(\det g)^2.
\label{pl:o:eq:7.3}
\end{equation}
These formulas also verify the derived group ($\SL_{6}$/$\mu_{2}$)$\times$$\SL_{2}$ and the central degree $d_3=6$.

One can obtain \eqref{pl:o:eq:7.1} directly from the WZ j-map: take its $\GL_{3}$ entry block diagonal of type (2,1), use the split central kernel to set $\lambda=1$, eliminate the last $\GL_{3}$ scalar using the determinant equation, and retain the residual action $(A,g)\mapsto(z^{-2}A,zg)$. This gives \eqref{pl:o:eq:7.2}.

\paragraph{Inverse image of H3 and the character}\label{pl:o:sec:7.2}

With $n(v)=\begin{psmallmatrix}1&v\\0&1\end{psmallmatrix}$, the inverse image is
\begin{equation}
{
\widetilde H_3=
\{((\det a)^{-1}n(v),\ D_3(a)u):
 a\in \GL_2,
 v\in F,
 u\in U_{2,2,2}\}.
}
\label{pl:o:eq:7.4}
\end{equation}
Its character is
\begin{equation}
\widetilde\xi_3=\psi(v+\operatorname{Tr}u_{12}+\operatorname{Tr}u_{23}),
\label{pl:o:eq:7.5}
\end{equation}
trivial on a. The kernel in \eqref{pl:o:eq:7.1} is already present by taking $a=zI_2$.

To check the $\PGL_{2}$ embedding, take $a=\diag(t,1)$. In \eqref{pl:o:eq:7.4} its image has simple-root values
\[
(1,t,1,t^{-1},t,t^{-1},t),
\]
whose cocharacter is exactly the defining diagonal cocharacter $v$ in \S\ref{m9:o:sec:2.3}. The unipotent $\SL_{2}$ generators map to the diagonal generators for $\alpha_{2}$,$\alpha_{5}$,$\alpha_{7}$. These root maps identify the full rational subgroup $H_0(F)$.

Transpose on both GL factors lifts $\iota$ and fixes the central
kernel pointwise, preserving its Haar density. Conjugation by a on the trace blocks has determinant one, and the scalar first-factor action centralizes its N. Thus the relative density character in \eqref{pl:o:eq:7.5} is trivial.

\paragraph{Closure of the patch}\label{pl:o:sec:7.3}

Pull the Levi descendant back to $\GL_2\times\GL_6$. By \S\ref{pl:o:sec:2} and \S\ref{pl:o:sec:3} this is injective, has nilpotent wavefront, is transpose-odd, and is equivariant under the subgroup \eqref{pl:o:eq:7.4} on both sides.

On the first-factor Gauss open, extract the two N-characters. Write the remaining diagonal matrix as $\diag(a_1,a_2)$. The reductive parameters $a,b\in\GL_2$ in the two copies of \eqref{pl:o:eq:7.4} act by
\[
(a_1,a_2)\longmapsto
(\det a\det b)^{-1}(a_1,a_2)
\]
and on $\GL_{6}$ by the corresponding two-sided GR action. The ratio $a_1/a_2$ is passive. Frobenius descent on the nonzero coordinate $a_2$, normalized to $1$, has the stabilizer
\begin{equation}
\det a\det b=1.
\label{pl:o:eq:7.6}
\end{equation}
Its phase is the balanced-GR phase, its density ratio is one, and transpose fixes the normalized torus coordinate. By \eqref{pl:o:eq:6.1}, the slice wavefront is zero on the passive ratio and nilpotent on $\GL_{6}$. Lemma~\ref{m9:balanced-gr}, including its passive-parameter version, therefore kills the entire first-factor Gauss open.

To remove the remaining Whittaker boundary, apply
Lemma~\ref{pl:whittaker-boundary} with $n=2$ to the first $\GL_{2}$
factor, retaining the wavefront of the original pulled-back distribution. Thus the pullback is zero, the Levi descendant is zero, and
\begin{equation}
{T|_{O_3}=0.}
\label{pl:o:eq:7.7}
\end{equation}
The conclusion uses NW GR vanishing with central balancing.

\subsubsection{The A3+A2+A1 patch: quotient and closure}\label{pl:o:sec:8}

\paragraph{The rational group}\label{pl:o:sec:8.1}

Use factors $\GL_{4,X}$, $\GL_{2,Y}$, $\GL_{3,A}$, with simple-root chains (7,6,5), (2), and (1,3), respectively. Then
\begin{equation}
{
M_4\simeq
(\GL_4\times \GL_2\times \GL_3)/
\{(t_XI_4,t_YI_2,t_AI_3):t_Xt_Yt_A=1\}.
}
\label{pl:o:eq:8.1}
\end{equation}
The kernel is a connected split two-dimensional torus.

For diagonal entries $X_i,Y_i,A_i$, the torus map is
\begin{equation}
\begin{aligned}
\alpha_1&=A_1/A_2,&\alpha_2&=Y_1/Y_2,&\alpha_3&=A_2/A_3,\\
\alpha_4&=X_4Y_2A_3,&\alpha_5&=X_3/X_4,&
\alpha_6&=X_2/X_3,&\alpha_7&=X_1/X_2.
\end{aligned}
\label{pl:o:eq:8.2}
\end{equation}
Together with the stated pinned root maps, this determines the homomorphism. It is F-surjective: set $X_4=Y_2=1$, solve $A_3$ from $\alpha_{4}$, and solve the remaining ratios successively. Its kernel is exactly \eqref{pl:o:eq:8.1}. The integral cocharacter matrix has seven nonzero Smith invariants all equal to one.

Each of the three scalar centers maps to $z_4$. The pullback of $\chi_4$ is
\begin{equation}
(\det X)^3(\det Y)^6(\det A)^4.
\label{pl:o:eq:8.3}
\end{equation}
This confirms $d_4=12$. On derived groups, the kernel is $(-I_4,-I_2,I_3)$, giving precisely the $j=4$ row of \eqref{pl:o:eq:2.1}.

For completeness, the cocharacter calculation underlying \eqref{pl:o:eq:8.2} is as follows. For a factor $\GL_{n}$ with chain $\gamma_1,\ldots,\gamma_{n-1}$, the first diagonal cocharacter is
\[
\frac{z_4}{n}+\sum_{k=1}^{n-1}\frac{n-k}{n}\gamma_k^\vee;
\]
subtract successive $\gamma_k^\vee$ to obtain the others. Their $\alpha_{4}$ pairings are integral because the last coefficient $1/n$ cancels $z_4/n$ against the adjacent end root. All other pairings are the ordinary $\GL_{n}$ ones. This gives \eqref{pl:o:eq:8.2}.

\paragraph{Subgroup preimage and phase}\label{pl:o:sec:8.2}

The inverse image on rational points of H4 is
\begin{equation}
{
\widetilde H_4=
\{(cD_2(a)u,\ a,\ (c\det a)^{-1}n):
 c\in F^\times,\ a\in \GL_2,\ u\in U_{2,2},\ n\in N_3\}.
}
\label{pl:o:eq:8.4}
\end{equation}
The character is
\begin{equation}
\widetilde\xi_4=
\psi\bigl(\operatorname{Tr}u_{12}+n_{12}+n_{23}\bigr).
\label{pl:o:eq:8.5}
\end{equation}
Here the $\GL_{3}$ simple entries are the pinned $\alpha_{1}$ and $\alpha_{3}$ coordinates, so their coefficients are the physical $v_1$ and $e$ coefficients. The $\GL_{4}$ trace is the physical x-trace.

The parameter c ranges over \textbf{all} of $F^\times$. Indeed \eqref{pl:o:eq:8.4} contains every central-kernel element in \eqref{pl:o:eq:8.1}: choose $a=t_YI_2$ and $c=t_X/t_Y$. Its third scalar is then $(t_Xt_Y)^{-1}=t_A$.

For $a=\diag(t,1)$ and $c=1$, formula \eqref{pl:o:eq:8.2}
gives the simple-root values
\[(1,t,1,t^{-1},t,t^{-1},t),\]
hence the cocharacter $v$. On the simply connected diagonal $\SL_{2}$ the map is $a\mapsto(D_2(a),a,I_3)$. This checks both the $H_0$ embedding and its $\PGL_{2}$ rational form. The connected central kernel and $H_0$ have the expected dimensions, so \eqref{pl:o:eq:8.4} is the full inverse image.

Transpose in all three factors lifts $\iota$, and the phase on the opposite side is consequently the exact one in \eqref{pl:o:eq:1.8}. Both central kernel directions act trivially on their density lines. The reductive part acts by determinant-one conjugation on the $\Mat_{2}$ block and centrally on the $\GL_{3}$ Whittaker radical; the relative density character is trivial.

\paragraph{Closure}\label{pl:o:sec:8.3}

Pullback through \eqref{pl:o:eq:8.1} gives an odd NW distribution on
$\GL_4\times\GL_2\times\GL_3$. Extract the two maximal-unipotent
characters on the $\GL_{3}$ Gauss open. The remaining diagonal
entries $(a_1,a_2,a_3)$ are all multiplied by
\[
\bigl(c_1\det b_1\,c_2\det b_2\bigr)^{-1}
=\bigl(\nu_1\nu_2\bigr)^{-1},\qquad \nu(c,b)=c\det b.
\]
The two ratios $a_1/a_3$ and $a_2/a_3$ are passive. Normalizing
$a_3=1$ by Frobenius descent leaves exactly the balanced extended-Shalika subgroup $\nu_1\nu_2=1$ on $\GL_4\times\GL_2$. The character is the trace character in \eqref{pl:o:eq:5.4}, the density ratio is one, and transpose acts trivially on the two torus ratios. Formula \eqref{pl:o:eq:6.1} supplies zero wavefront in those passive directions and nilpotent wavefront on $\GL_4\times\GL_2$.

Applying the balanced extended-Shalika theorem
\eqref{pl:o:eq:5.6} clears this Gauss open. We then use
Lemma~\ref{pl:whittaker-boundary}, with $n=3$, to clear the remaining
Whittaker boundary. Pullback injectivity and \eqref{pl:o:eq:1.4} give
\begin{equation}
{T|_{O_4}=0.}
\label{pl:o:eq:8.6}
\end{equation}
The inputs are the Shalika distribution theorem and the NW descent
lemmas proved above.

\subsection{The \texorpdfstring{$D_6$}{D6} patch}
\label{m10:d6}
\begin{corollary}
\label{m10:o1}
The original $T$ vanishes on $O_1$.
\end{corollary}
\begin{proof}
Apply free descent to the semispin central form
$M_{10,1}$. After $m=xb_-$ its datum is exactly the $\kappa=1$ case of
Theorem~\ref{m9:native-semispin}, with the $F$-rational acting subgroup,
plus-plus action, NW Chevalley operation and trivial density character.
It therefore clears the entire $48$-cell patch, including the $32$ cells
assigned after $O_3,O_4$, independently of the type-A opens. Injectivity
of free descent proves the claim.
\end{proof}

\subsection{The relative \texorpdfstring{$D_5+A_1$}{D5+A1} theorem}
\label{m10:b5-section}
The derived group of $M_{10,6}$ is
\[
 D_{10,6}=
 (\Spin_{10}\times\SL_2)/\langle(-1,-I_2)\rangle,
\]
where the subscript records the omitted $E_7$ node. The Levi group is
\[
 M_{10,6}=(D_{10,6}\times\mathbb G_m)/
 \{(z_6(\zeta)^{-1},\zeta):\zeta^4=1\},\qquad z_6=\omega_6^\vee.
\]
The natural product map has fourth-power rational sectors; the proof works
directly on $M_{10,6}(\F)$. Its acting subgroup is $\PGL_2\ltimes U_6$, with $\dim U_6=18$ and
$\xi_{10,6}=\psi(v_1+e+\tr z)$. Its $D_5$ simple-root order is
$(\alpha_1,\alpha_3,\alpha_4,\alpha_2,\alpha_5)$, with external root
$\alpha_7$. Write the phase coordinates of a Levi factor as
$(s_l,t_l,B_l)\in(\F^\times)^2\times\GL_2(\F)$.
The nilpotent cone has nilpotent components
on the $\mathfrak{so}_{10}$ and $\mathfrak{sl}_2$ ideals and zero central
component.

For $i\ne j$, define the invariant open
\[
 V_i^{(j)}=\{x\in M_{10,j}(\F):xb_-\omega_{10}\in O_i\}.
\]
\begin{theorem}[B5 with its predecessor condition]
\label{m10:b5}
Let $S$ be a distribution on $M_{10,6}(\F)$, satisfying NW,
$\iota_{10}$-odd and equivariant for the NW plus-plus action with
character $\xi_{10,6}\boxtimes\xi_{10,6}$. Then
\[
 S|_{V_1^{(6)}\cup V_3^{(6)}\cup V_4^{(6)}}=0
 \quad\Longrightarrow\quad S=0.
\]
\end{theorem}

The complete compatible remainder consists of cells
$49599,51648,52209,52215,52335$, of original exceptional lengths
$42,47,50,50,52$. The descending support order is
\[
 52335\ \longrightarrow\ \{52209,52215\}
 \ \longrightarrow\ 51648\ \longrightarrow\ 49599.
\]
The first four steps use ambient conormals to the full smooth
unipotent-compatible loci, with all rational sectors retained. The final
cell has an intrinsic character-one transporter. Their coordinates,
conormal calculations and predecessor opens are proved below.

\subsubsection{The setting and matrix notation}\label{b5:o:sec:1}

Write
\[
I=\{2,5,7\},\qquad K_6=\Delta\setminus\{6\},\qquad
M=M_6,\qquad U=U_6,\qquad H=H_6=H_0\ltimes U.
\]
The parabolic of $M$ with Levi $L$ is denoted by
\[
\Pi=LU,\qquad \Pi^-=\iota(\Pi).
\]
Here $L$ is the original rank-seven Levi of derived type $A_1^3$.

We retain the group
\[
[M,M]=(\operatorname{Spin}_{10}\times \SL_2)/\langle(-1,-I_2)\rangle,
\]
with the specified central gluing, and the $F$-rational subgroup $H_0 = \PGL_{2}$. The distribution convention is exactly
\[
(h,k)x=hx\iota(k),\qquad \eta(h,k)=\xi_6(h)\xi_6(k),
\qquad \xi_6=\psi(v_1+e+\operatorname{Tr}z).
\]
The translated variable is $x=mb^{-1}$. Set
\[
V=V_1^{(6)}\cup V_3^{(6)}\cup V_4^{(6)}.
\]
The assumption throughout is $S|V = 0$.

\paragraph{The pinned D5 Lie algebra}\label{b5:o:sec:1.1}

Use the inherited correspondence
\[
(\beta_1,\beta_2,\beta_3,\beta_4,\beta_5)
=(\alpha_1,\alpha_3,\alpha_4,\alpha_2,\alpha_5).
\]
In the standard Euclidean description of the D5 root system,
\[
\beta_1=\epsilon_1-\epsilon_2,\quad
\beta_2=\epsilon_2-\epsilon_3,\quad
\beta_3=\epsilon_3-\epsilon_4,\quad
\beta_4=\epsilon_4-\epsilon_5,\quad
\beta_5=\epsilon_4+\epsilon_5.
\]
Use the ordered vector basis
\[
(1,2,3,4,5,-5,-4,-3,-2,-1)
\]
with anti-diagonal symmetric Gram matrix. Define
\[
E_{i-j}=E_{i,j}-E_{-j,-i},\qquad
E_{i+j}=E_{i,-j}-E_{j,-i}\quad(i<j),
\]
\[
H_i=E_{i,i}-E_{-i,-i}.
\]
The simple positive generators are $E_{1-2}, E_{2-3}, E_{3-4}, E_{4-5}, E_{4+5}$, with their transposes as negative generators.

These matrices realize the pinned D5 Lie algebra and its rational
adjoint action, with all group elements $l,r$ taken in $L(F)$.
On the D5 ideal we use the invariant bilinear form
\[
\langle A,B\rangle=\tfrac12\operatorname{Tr}(AB),
\]
extended by any nondegenerate invariant forms on the external $A1$ and the central ideal.

The roots of $U$ are the eighteen roots
\[
\epsilon_i\pm\epsilon_j,\qquad 1\leq i\leq3,\quad i<j\leq5.
\]
The D5 part of the Levi Lie algebra has basis
\[
H_1,\ldots,H_5,E_{4-5},E_{4+5},E_{4-5}^t,E_{4+5}^t.
\]
The remaining Levi directions are the external $\mathfrak{sl}_{2}$ and the one-dimensional center.

\paragraph{Character coordinates}\label{b5:o:sec:1.2}

For $l$ in $L(F)$, write the phase $\lambda_6\circ \operatorname{Ad}(l)$ on $u$ as
\[
p_l=(s_l,t_l,B_l),\qquad s_l,t_l\in F^\times,\quad B_l\in \GL_2(F).
\]
These are obtained from the specified coordinates $(s,t,A,B)$ by forgetting $A$. Explicitly,
\begin{equation}
\begin{aligned}
p_l(Y)={}&s_lY_{1,2}+t_lY_{2,3}
 +(B_l)_{21}Y_{3,4}+(B_l)_{22}Y_{3,5}\\
&-(B_l)_{11}Y_{3,-5}+(B_l)_{12}Y_{3,-4}.
\end{aligned}\label{b5:o:eq:1.1}
\end{equation}
The corresponding four-vector is
\[
c(B)=(B_{21},B_{22},-B_{11},B_{12}),\qquad q(c(B))=-\det B.
\]
In particular, the identity trace functional is represented by $e_5-e_{-5}$, a nonisotropic vector. In this pinned identification, the coefficient of $B_{11}$ has negative sign.

To obtain a rational section of
$L\longrightarrow (\mathbb G_m)^2\times \GL_{2}$, use the section
of $H_0\backslash L$, set the $A$ coordinate equal to the identity,
and project. Infinitesimally, the nonisotropic four-vector has its three-dimensional orthogonal orbit and an independent scaling direction; the two nonzero scalar coordinates supply the other two directions. This also verifies submersivity at every point, including over residue characteristic two.

\subsubsection{The five relative words}\label{b5:o:sec:2}\label{app:b5}

For the five $E_7$ cells listed in \S\ref{m10:b5-section}, evaluate
\[
v=ww_0,\qquad w_0=-s_2s_5s_7.
\]
The resulting external Weyl component is the identity. On the five D5 weight coordinates the results are:

\begingroup\small\setlength{\tabcolsep}{3pt}
\begin{longtable}{@{}>{\raggedright\arraybackslash}p{0.235\linewidth}>{\raggedright\arraybackslash}p{0.235\linewidth}>{\raggedright\arraybackslash}p{0.235\linewidth}>{\raggedright\arraybackslash}p{0.235\linewidth}@{}}
\toprule
E7 cell & E7 length & Images of $(\epsilon_1 ,...,\epsilon_5 )$ under $v$ & A reduced D5 word for $v$ \\
\midrule\endhead
49599 & 42 & $(-1,-2,-3,4,-5)$ & \texttt{352341\allowbreak{}235234\allowbreak{}123121} \\
51648 & 47 & $(-4,-2,-3,1,-5)$ & \texttt{534235\allowbreak{}234123\allowbreak{}2} \\
52209 & 50 & $(-4,-5,3,1,2)$ & \texttt{523412\allowbreak{}3121} \\
52215 & 50 & $(-4,5,3,1,-2)$ & \texttt{423512\allowbreak{}3121} \\
52335 & 52 & $(-4,3,2,1,-5)$ & $54123121$ \\
\bottomrule\end{longtable}\endgroup

A signed integer $-4$ denotes $-\epsilon_4$. The local word uses the beta numbering above. The corresponding pinned lift is
\[
n_v=\prod n_{\beta_i},\qquad
n_{\beta_i}=\exp(E_{\beta_i})\exp(-E_{\beta_i}^t)\exp(E_{\beta_i}).
\]
The five local word lengths are $18,13,10,10,8$.

Let
\[
X_c=\{x\in M(F):xbw_0\text{ belongs to the original cell }c\}.
\]
The specified free-patch coordinates identify this set with
\[
X_c=\Pi(F)n_v\Pi^-(F).
\]
Indeed, after multiplying by $w_0$ on the right, the outer $N_6$ factors are the two ordinary Gauss factors, and their uniqueness leaves exactly the indicated parabolic double cell in $M$.

The representative $n_vbw_0$ and the specified $E_7$ lift $n_w$ differ by an element of $T(F)$:
\[
d_c=n_vbw_0n_w^{-1}\in T(F).
\]
For $x=l n_v \iota (r)$,
\[
xbw_0=l d_c n_w\,\operatorname{Ad}((bw_0)^{-1})(\iota(r)).
\]
This gives the exact change between the local Levi-pair chart and the specified E7 chart. The right phase in these coordinates is $p_r$.

All five coarse cells are transpose-fixed. For the local Weyl representatives, $v^{-1}=a v a$ with $a=1$ for 49599 and $a=s_{\beta_4}  s_{\beta_5} $ for each of the other four cells.

\subsubsection{U-stabilizer and compatible manifolds}\label{b5:o:sec:3}

At
\[
x=l n_v\iota(r),\qquad l,r\in L(F),
\]
put $p=p_l=(s_L,t_L,B_L)$ and $q=p_r=(s_R,t_R,B_R)$.

Let
\[
\mathfrak k_v=\{Y\in\mathfrak u:\operatorname{Ad}(n_v)Y^t\in\mathfrak u\}.
\]
The $U\times U$ stabilizer is the graph
\begin{equation}
Y\longmapsto
\left(l\exp(-\operatorname{Ad}(n_v)Y^t)l^{-1},
      r\exp(Y)r^{-1}\right),\qquad Y\in\mathfrak k_v.
\label{b5:o:eq:3.1}
\end{equation}
This is a parametrization by the whole nilpotent Lie algebra, using the characteristic-zero exponential. Its character is exactly
\begin{equation}
\psi\big(q(Y)-p(\operatorname{Ad}(n_v)Y^t)\big).
\label{b5:o:eq:3.2}
\end{equation}
Thus compatibility is the equality of these two phases on all of $\mathfrak k_v$.

\paragraph{Root sets of the complete unipotent stabilizers}\label{b5:o:sec:3.1}

Use $i-j$ for $\epsilon_i-\epsilon_j$ and $i+j$ for $\epsilon_i+\epsilon_j$. These are the roots of $Y$ in \eqref{b5:o:eq:3.1}:

\begingroup\small\setlength{\tabcolsep}{3pt}
\begin{longtable}{@{}>{\raggedright\arraybackslash}p{0.313\linewidth}>{\raggedright\arraybackslash}p{0.313\linewidth}>{\raggedright\arraybackslash}p{0.313\linewidth}@{}}
\toprule
Cell & Roots & Dimension \\
\midrule\endhead
49599 & All eighteen roots of $U$ & 18 \\
51648 & $1+2, 1+3, 1-4, 2-3, 2+3, 2-4, 2-5, 2+5, 3-4, 3-5, 3+5$ & 11 \\
52209 & $1-3, 1-4, 1-5, 2-3, 2-4, 2-5, 3-4, 3-5$ & 8 \\
52215 & $1-3, 1-4, 1+5, 2-3, 2-4, 2+5, 3-4, 3+5$ & 8 \\
52335 & $1-2, 1-3, 1-4, 2-3, 2-4, 3-4$ & 6 \\
\bottomrule\end{longtable}\endgroup

\paragraph{Compatible loci in the descendant}\label{b5:o:sec:3.2}

In the \textbf{local plus-plus chart just specified}, all nonzero equations in \eqref{b5:o:eq:3.2} are equivalent to the following table. Throughout impose
\[
s_Lt_Ls_Rt_R\det B_L\det B_R\ne0.
\]
The $A_L,A_R$ coordinates of the $H_0\backslash L$ chart remain invertible and unrestricted by this unipotent test.

\begingroup\small\setlength{\tabcolsep}{3pt}
\begin{longtable}{@{}>{\raggedright\arraybackslash}p{0.470\linewidth}>{\raggedright\arraybackslash}p{0.470\linewidth}@{}}
\toprule
Cell & Complete unipotent-compatibility equations \\
\midrule\endhead
49599 & $s_R=-s_L; t_R=-t_L; B_R=B_L^t$ \\
51648 & $(B_L)_{12}=0; (B_R)_{21}=0; t_R=-t_L; (B_R)_{11}=-(B_L)_{11}; (B_R)_{22}=-(B_L)_{22}$ \\
52209 & $(B_L)_{12}=0; (B_R)_{21}=0; t_R=-(B_L)_{11}; (B_R)_{22}=t_L$ \\
52215 & $(B_L)_{12}=0; (B_R)_{21}=0; t_R=(B_L)_{22}; (B_R)_{11}=-t_L$ \\
52335 & $s_R=(B_L)_{12}; t_R=-t_L; (B_R)_{21}=-s_L$ \\
\bottomrule\end{longtable}\endgroup

These equations describe the preimage of the unipotent-compatible submanifold $C_c$ under
\begin{equation}
U\times L\times L\times U\longrightarrow X_c,
\qquad (u,l,r,z)\longmapsto u l n_v\iota(r)\iota(z).
\label{b5:o:eq:3.3}
\end{equation}
Because the map is a surjective submersion and the equations are
independent linear equations in the twelve phase coordinates,
subject to the stated nonvanishing conditions, their preimage is
smooth. Submersion coordinates then show that $C_c$ is smooth.
Compatibility is closed on the cover, so the open quotient map
\eqref{b5:o:eq:3.3} also shows that it is closed in $X_c$.

At every point outside $C_c$, the unipotent stabilizer
\eqref{b5:o:eq:3.1} has nontrivial character. The ordinary JLX
stabilizer criterion therefore eliminates this complement intrinsically.

The dimensions, computed from the root tangent spaces, are:

\begingroup\small\setlength{\tabcolsep}{3pt}
\begin{longtable}{@{}>{\raggedright\arraybackslash}p{0.188\linewidth}>{\raggedright\arraybackslash}p{0.188\linewidth}>{\raggedright\arraybackslash}p{0.188\linewidth}>{\raggedright\arraybackslash}p{0.188\linewidth}>{\raggedright\arraybackslash}p{0.188\linewidth}@{}}
\toprule
Cell & $\dim X_c$ in $M_6$ & Number of independent equations & $\dim C_c$ & $\dim\CN(C_c)$ \\
\midrule\endhead
49599 & 31 & 6 & 25 & 24 \\
51648 & 38 & 5 & 33 & 16 \\
52209 & 41 & 4 & 37 & 12 \\
52215 & 41 & 4 & 37 & 12 \\
52335 & 43 & 3 & 40 & 9 \\
\bottomrule\end{longtable}\endgroup

\subsubsection{The $H\times H$ stabilizer and its character}\label{b5:o:sec:4}

We describe the full $H\times H$ stabilizers, in addition to the
$U\times U$ kernels. The conormal proof applies to the full $C_c$.

\paragraph{A group-level presentation, including rational central conditions}\label{b5:o:sec:4.1}

Set
\[
\mathcal P_v=\Pi\cap n_v\Pi^-n_v^{-1},\qquad
k^\sharp=\iota(n_v^{-1}k^{-1}n_v).
\]
For $k$ in $\mathcal P_v(F)$, write
\[
k=u_1d_1,\qquad k^\sharp=u_2d_2,
\qquad u_i\in U(F),\ d_i\in L(F).
\]
Define
\begin{equation}
\mathcal P_v(l,r)=\{k\in\mathcal P_v(F):
 l d_1 l^{-1}\in H_0(F),\ r d_2r^{-1}\in H_0(F)\}.
\label{b5:o:eq:4.1}
\end{equation}
Then the \textbf{entire} stabilizer of $l n_v \iota (r)$ is
\begin{equation}
{
\operatorname{Stab}_{H\times H}(x)=
\{(l k l^{-1},\ r k^\sharp r^{-1}):k\in\mathcal P_v(l,r)\}.
}\label{b5:o:eq:4.2}
\end{equation}
Its character is
\begin{equation}
{
\eta_x(k)=\psi\big(p_l(\log u_1)+p_r(\log u_2)\big).
}\label{b5:o:eq:4.3}
\end{equation}
Equations \eqref{b5:o:eq:4.1} parametrize the full stabilizer,
including its connected components and finite central conditions.

For a concrete root-group description, the quotient
\[
\Gamma_v=\operatorname{im}\big(k\mapsto(d_1,d_2)\big)\subset L\times L
\]
has dimension thirteen. Its kernel is precisely the unipotent group in \eqref{b5:o:eq:3.1}.

\begin{itemize}
\item For 49599, $n_v$ normalizes $L$, and $\Gamma_v$ is the graph of $d\longrightarrow \iota (n_v^{-1}d^{-1}n_v)$ on the $L$.
\item For the other four cells, it is generated by the torus pairs
  \[
  (t,v^{-1}(t)^{-1}),\quad t\in T,
  \]
  the two left root groups $(x_{\beta_4}(a),1)$ and $(x_{\beta_5}(a),1)$, the two right root groups $(1,x_{-\beta_4}(a))$ and $(1,x_{-\beta_5}(a))$, and the external root-group graph $(h,\iota (h^{-1}))$.

\end{itemize}
Thus \eqref{b5:o:eq:4.1} involves only this specified small root-group image in $L\times L$. The torus pairs encode the finite central gluing.

\paragraph{The three unipotent-alignment cases}\label{b5:o:sec:4.2}

For 51648, 52209 and 52215, set
\begin{equation}
N_L=B_{L,11}E_{4-5}+B_{L,22}E_{4+5},\qquad
N_R=B_{R,22}E_{4-5}^t+B_{R,11}E_{4+5}^t.
\label{b5:o:eq:4.4}
\end{equation}
Both are nonzero. On the compatible locus $(B_L)_{12}=(B_R)_{21}=0$; both diagonal entries of each $B$ are nonzero. The surviving left internal subgroup is $\exp(\rho  N_L)$ and the surviving right internal subgroup is $\exp(\sigma  N_R)$.

The group-level calculation uses the matrix
\[
x_{\beta_4}(a)x_{\beta_5}(b)=
\begin{pmatrix}1&a&b&-ab\\0&1&0&-b\\0&0&1&-a\\0&0&0&1\end{pmatrix}
\]
on $(4,5,-5,-4)$. At
$c(B_L)=(B_{L,21},B_{L,22},-B_{L,11},0)$ its fixed-vector
equation is $aB_{L,22}-bB_{L,11}=0$, yielding the left group
in \eqref{b5:o:eq:4.4}; transposition gives the right group.
The torus graph in \eqref{b5:o:eq:4.1} forces the remaining
diagonal parameter to one, since both middle coordinates are
nonzero. Thus the parametrization includes every finite diagonal component.
For the expanded substitutions, see Section M10-B5H of
\path{anc/model10_calculations.pdf}.

Let $\rho_{\mathrm{ext}}:M\longrightarrow \PGL_{2}$ be the external adjoint-factor projection, and put $y=\rho_{\mathrm{ext}}(x)$. Via the fixed identification of the internal three-dimensional stabilizer with $H_0$, let $X_L,X_R$ be the external $\mathfrak{sl}_{2}$ images of the internal generators $\operatorname{Ad}(l)N_L, \operatorname{Ad}(r)N_R$. The identification is linear. Put
\begin{equation}
\mathfrak a_x=\{(\rho,\sigma)\in F^2:
 \rho X_L+\sigma\operatorname{Ad}(y)d\iota_{\rm ext}(X_R)=0\}.
\label{b5:o:eq:4.5}
\end{equation}
This vector space has dimension zero or one, and the stabilizer is
an extension of $\exp(\mathfrak k_v)$ by the additive group
$\mathfrak a_x$. Being unipotent, these groups are unimodular.

The character on the stabilizer is trivial on $\exp(\mathfrak k_v)$ and is as follows on the quotient:

\begingroup\small\setlength{\tabcolsep}{3pt}
\begin{longtable}{@{}>{\raggedright\arraybackslash}p{0.470\linewidth}>{\raggedright\arraybackslash}p{0.470\linewidth}@{}}
\toprule
Cell & Phase on $(\rho ,\sigma )$ in $\mathfrak a_x$ \\
\midrule\endhead
51648 & $0$ \\
52209 & $s_R (B_L)_{11} \rho  - s_L (B_R)_{22} \sigma $ \\
52215 & $s_R (B_L)_{22} \rho  - s_L (B_R)_{11} \sigma $ \\
\bottomrule\end{longtable}\endgroup

For example, on 52209 a left $E_{\beta_4} $ generator becomes the negative of the right $E_{\beta_1} $ after Weyl transport, yielding $+s_R$; a right $E_{-\beta_4} $ produces the compensating left phase $-s_L$. On 52215 the same statement uses $\beta_5 $ instead. All other cross-root contributions are character-zero. Exponentiation gives the displayed additive characters on the quotient, while compatibility kills the kernel contribution.

Consequently the locus compatible with the $H\times H$ stabilizer character on 52209 and 52215 is exactly the subset where the displayed linear functional annihilates $\mathfrak a_x$. At points with $\mathfrak a_x=0$ it is automatically compatible. At aligned points, a nonzero restriction gives additional intrinsic vanishing. The conormal proof applies on the entire smooth $C_c$, which contains this subset.

\paragraph{Cell 52335}\label{b5:o:sec:4.3}

Here $(B_L)_{12}$ and $(B_R)_{21}$ are nonzero. The internal stabilizers of the coefficient vectors are split one-dimensional tori rather than the unipotent groups in \eqref{b5:o:eq:4.4}.

The normalizers are
\[
\begin{aligned}
n_L&=x_{\beta_4}(-B_{L,11}/B_{L,12})x_{\beta_5}(B_{L,22}/B_{L,12}),\\
n_R&=x_{-\beta_4}(-B_{R,22}/B_{R,21})x_{-\beta_5}(B_{R,11}/B_{R,21}).
\end{aligned}
\]
Their fixed-vector equations give the two tori
\[
n_L^{-1}\diag(1,a,a^{-1},1)n_L,\qquad
n_R^{-1}\diag(1,a,a^{-1},1)n_R.
\]
The common parameter and every finite central restriction are
matched by the group torus graph in \eqref{b5:o:eq:4.1}. Section M10-B5H gives the full vector
substitutions; all those points belong to the same smooth
$C_{52335}$ used by the conormal proof.

The character \eqref{b5:o:eq:4.3} is trivial. The two cross-root groups on each side map to $1\pm 5$, on which the Bessel phase is zero. The kernel has dimension six. Its root sum in the left chart is $(3,1,-1,3,0)$, so the surviving $\epsilon_5 $ torus has determinant one on it. Its modular character is therefore trivial.

\subsubsection{Intrinsic vanishing on cell 49599}\label{b5:o:sec:5}

Let $n=n_v$ for this cell. It takes $U^-$ to $U$, normalizes $L$, and swaps the two internal leaf root groups while fixing their diagonal subgroup. Thus
\[
X_{49599}=\Pi n.
\]
Choose $s_0$ in the split adjoint torus by
\begin{equation}
(\alpha_1(s_0),\ldots,\alpha_7(s_0))=(-1,1,-1,1,1,1,1),
\qquad g_0=s_0n.
\label{b5:o:eq:5.1}
\end{equation}
The phase of $s_0$ is $(-1,-1,I_2)$, so it satisfies the compatibility requirement against the identity right phase.

\paragraph{The finite central sign}\label{b5:o:sec:5.1}

The sum of the eighteen $U$ coroots is
\[
8\beta_1^\vee+14\beta_2^\vee+18\beta_3^\vee
 +9\beta_4^\vee+9\beta_5^\vee.
\]
Multiplication of the rank-one Tits relations along the displayed word gives $n^2$ equal to this cocharacter at $-1$. Its seven simple-root values in the E7 torus are
\[
(1,1,1,1,1,-1,1).
\]
The same values are obtained for $v(s_0)s_0^{-1}$. Since the seven simple roots are a basis of the character lattice of the adjoint E7 torus, this proves the torus identity
\[
v(s_0)s_0^{-1}=n^2.
\]
Consequently
\begin{equation}
\iota(g_0)=n^{-1}s_0=g_0.\label{b5:o:eq:5.2}
\end{equation}

\paragraph{The compatible locus and its orbits}\label{b5:o:sec:5.2}

Let
\[
R=\operatorname{Stab}_{L}(\lambda_6).
\]
In the specified Levi quotient coordinates this is exactly
\[
s=t=1,\qquad B=I_2,
\]
with the $A$ coordinate unrestricted. The compatible locus is
\begin{equation}
C_{49599}=U R g_0.\label{b5:o:eq:5.3}
\end{equation}
Indeed, after the unique $U L n$ decomposition, compatibility says that the Levi phase equals the phase of $s_0$, which is equivalent to membership in $R s_0$.

The specified rational section gives an exact representative
\begin{equation}
c(C)=j\left(\operatorname{diag}(1,1,(\det C)^{-1}),
            \operatorname{diag}(C,I_2,I_2),1\right),
\qquad C\in \GL_2(F).
\label{b5:o:eq:5.4}
\end{equation}
Every element of $R(F)$ is uniquely $h c(C)$, with $h\in H_0(F)$. Conjugation by $H_0=\PGL_{2}$ acts on $C$ by ordinary rational conjugation, and
\[
\iota(c(C))=c(C^t).
\]
The connected group $R$ has derived Lie algebra equal to the
internal diagonal $\mathfrak{sl}_{2}$ plus the external
$\mathfrak{sl}_{2}$, and connected center $Z(M)$. Since the lift $n$
fixes both diagonal simple generators and the external factor, and
centralizes $Z(M)$, the element $g_0$ centralizes $R$.

It follows that the $H\times H$ orbits in \eqref{b5:o:eq:5.3} are parameterized by the \textbf{rational conjugacy classes of invertible $2\times2$ matrices}. A complete description is trace and determinant, together with the distinction between a scalar matrix and a nonscalar matrix in the repeated-eigenvalue case.

At $c(C)g_0$, the stabilizer is
\begin{equation}
U\rtimes\{\bar a\in \PGL_2(F):aCa^{-1}=C\},
\label{b5:o:eq:5.5}
\end{equation}
with the right factor determined by $\iota (h^{-1})$. The condition fixes $C$, and the stabilizer character is one.

\paragraph{Character-one transporter}\label{b5:o:sec:5.3}

Choose $a\in \GL_{2}(F)$ with $aCa^{-1}=C^t$. Such an $a$ always exists: for nonscalar $C$, put it in its rational companion form
\[
K=\begin{pmatrix}0&-\det C\\1&\operatorname{Tr}C\end{pmatrix};
\]
the invertible symmetric matrix $\diag(1,-\det C)$ intertwines $K$ with $K^t$. For scalar matrices the identity is a conjugator.

Let $h(a)$ be the corresponding element of $H_0(F)=\PGL_{2}(F)$. The pair
\[
\bigl(h(a),\ \iota(h(a)^{-1})\bigr)
\]
acts by conjugation and sends $c(C)g_0$ to its transpose. Its character is one, and this property is preserved when the pair
is transported back through the $H\times H$ normalization: the
anti-involution exchanges the two factors and preserves their
product character.

\cite[Lemma 2.7]{JLX} therefore kills the odd distribution line on each compatible orbit. The algebraic constructibility criterion and the ordinary open/closed exact sequence combine this with the incompatible complement. Thus
\begin{equation}
{\mathcal D'_{\widetilde\xi}(X_{49599})=0.}\label{b5:o:eq:5.6}
\end{equation}
This is intrinsic odd vanishing.

\subsubsection{The common semisimple-conormal calculation}\label{b5:o:sec:6}

\paragraph{Tangent and normal spaces}\label{b5:o:sec:6.1}

At $x=l n_v \iota (r)$ use right translation and then $\operatorname{Ad}(l^{-1})$ to identify the tangent space with $\mathfrak m$. For a root $Y$ in $\mathfrak k_v$, put
\[
F_Y(p,q)=p(\operatorname{Ad}(n_v)Y^t)-q(Y).
\]
For $A,B\in\mathfrak l$, its differential along $l \exp(A), r \exp(B)$ is
\begin{equation}
dF_Y(A,B)=p([A,\operatorname{Ad}(n_v)Y^t])-q([B,Y]).
\label{b5:o:eq:6.1}
\end{equation}
Let $D_v$ denote the collection of these differentials for the independent equations of Section~\ref{b5:o:sec:3}. Then the compatible tangent space in this trivialization is
\begin{equation}
{
\mathcal T_v(p,q)=
\mathfrak u+\operatorname{Ad}(n_v)\mathfrak u^-
+\{A+\operatorname{Ad}(n_v)B^t:(A,B)\in\ker D_v\}.
}\label{b5:o:eq:6.2}
\end{equation}
The external $\mathfrak{sl}_{2}$ and the center are included in the Levi variables in this formula. The normal space is $\mathfrak m/\mathcal T_v(p,q)$ and the conormal is its annihilator. Under a stabilizer pair $(h1,h2)$, the normal representation is the one induced by $\operatorname{Ad}(h1)$ in right trivialization.

We need only one explicit element of that annihilator.

\paragraph{Two matrix shapes, one family-level argument}\label{b5:o:sec:6.2}

Set
\[
H_* = H_2-H_3=\beta_2^\vee,
\qquad E_*=E_{1+2}=E_{1,-2}-E_{2,-1}.
\]
Then
\[
[H_*,E_*]=E_*,\qquad E_*^2=0.
\]
The following table gives conormals in the trivialization of \eqref{b5:o:eq:6.2}:

\begingroup\small\setlength{\tabcolsep}{3pt}
\begin{longtable}{@{}>{\raggedright\arraybackslash}p{0.313\linewidth}>{\raggedright\arraybackslash}p{0.313\linewidth}>{\raggedright\arraybackslash}p{0.313\linewidth}@{}}
\toprule
Cell & Selected compatibility equation $F$ & Conormal $C$ \\
\midrule\endhead
49599 & $-t_L-t_R=0$ & $-t_L H_\ast$ \\
51648 & $-t_L-t_R=0$ & $-t_L H_\ast$ \\
52209 & $t_L-(B_R)_{22}=0$ & $t_L H_\ast - (B_R)_{12} E_\ast$ \\
52215 & $t_L+(B_R)_{11}=0$ & $t_L H_\ast - (B_R)_{12} E_\ast$ \\
52335 & $-t_L-t_R=0$ & $-t_L H_\ast$ \\
\bottomrule\end{longtable}\endgroup

For each row, the exact identities are
\begin{equation}
\langle C,\mathfrak u\rangle=0,
\qquad \langle C,\operatorname{Ad}(n_v)\mathfrak u^-\rangle=0,
\label{b5:o:eq:6.3}
\end{equation}
and, \textbf{on the compatible locus},
\begin{equation}
{
\langle C,A+\operatorname{Ad}(n_v)B^t\rangle=dF(A,B)
\quad\text{for every }A,B\in\mathfrak l.
}\label{b5:o:eq:6.4}
\end{equation}
These identities hold for all coefficients satisfying compatibility.

For the diagonal rows, $v$ takes $\beta_2 $ to $-\beta_2 $; the Levi leaf root spaces commute with the simple $\beta_2 $ character coordinate. Equation \eqref{b5:o:eq:6.4} is the differential of $-t_L-t_R=0$.

For 52209, the stabilizer root is $Y=E_{3-5}$ and
\[
\operatorname{Ad}(n_v)Y^t=E_{2-3}.
\]
The diagonal part of the differential of $t_L-(B_R)_{22}$ is represented by $t_L H_\ast$. The right leaf-root derivative contributing $(B_R)_{12}$ is represented by $-(B_R)_{12} E_\ast$; the other possible derivative is zero because $(B_R)_{21}=0$. For 52215 use $Y=E_{3+5}$; its right phase is $-(B_R)_{11}$, giving $t_L+(B_R)_{11}$ and exactly the same correction. Direct checking on the nine displayed D5 Levi basis vectors proves \eqref{b5:o:eq:6.4}. External and central vectors contribute zero to both sides.

For \eqref{b5:o:eq:6.3}, the diagonal part pairs trivially with all root vectors. The only root vector which can pair with $E_\ast$ is $E_{-1-2}$. It occurs in neither $u$ nor $\operatorname{Ad}(n_v)u^-$ for the two rows: its inverse Weyl image is a Levi root.

Thus \eqref{b5:o:eq:6.2}--\eqref{b5:o:eq:6.4} prove that $C$ is
a conormal to $C_c$ in the ambient group.

\paragraph{Semisimplicity at all parameters}\label{b5:o:sec:6.3}

Since $t_L \ne  0$ throughout the compatible locus, the diagonal
conormal is nonzero semisimple.

For the two rows,
\begin{equation}
{
 t_LH_*-B_{R,12}E_*
 =\operatorname{Ad}\bigl(x_{\epsilon_1+\epsilon_2}(B_{R,12}/t_L)\bigr)(t_LH_*).
}\label{b5:o:eq:6.5}
\end{equation}
Hence it too is nonzero semisimple \textbf{for every value of $(B_R)_{12}$}, including zero. In all five cases,
\begin{equation}
\tfrac12\operatorname{Tr}(C^2)=2t_L^2\ne0.
\label{b5:o:eq:6.6}
\end{equation}
At the point $u l n_v \iota (r) \iota (z)$, the corresponding right-trivialized conormal is
\[
\operatorname{Ad}(ul)C.
\]
It is semisimple over $F$ and lies in the D5 ideal.

The construction is defined over $M(F)$: $H_\ast$, $E_\ast$, the root subgroup in \eqref{b5:o:eq:6.5}, and $\operatorname{Ad}(ul)$ all exist in the inherited pinning of the group. It therefore applies in every fourth-power rational sector of $M(F)$.

\paragraph{The elementary boundary lemma}\label{b5:o:sec:6.4}

\textbf{Five-cell Bessel boundary lemma.} Let $W$ be an open submanifold of $M(F)$, and suppose an ambient distribution $R$ on $W$ has nilpotent wavefront and is supported in $C_c$, for one of these five cells. Then $R=0$.

\textbf{Proof.} At every point of $C_c$, \S\ref{b5:o:sec:6.1}--\ref{b5:o:sec:6.3} supply a nonzero semisimple conormal. By \cite[Theorem 2.19]{JLX}, the wavefront of an ambient distribution supported on a smooth submanifold is invariant under conormal translation; at a support point its extended wavefront contains the conormal space. The exhibited cotangent vector would therefore belong to the wavefront of $R$, contradicting NW nilpotence. This argument is applied at every possible support point. Open locality, \cite[Proposition 2.17]{JLX}, preserves the wavefront hypothesis on $W$.

The lemma also applies if prior $H\times H$ stabilizer arguments have confined the support to a smaller subset of $C_c$: support is still contained in the same smooth ambient submanifold $C_c$.

\subsubsection{The support filtration}\label{b5:o:sec:7}

Let $\ell(x)$ be the length of the original E7 minimal double-coset representative containing $x b w_0$. Define
\begin{equation}
\mathcal U_{>d}=V\cup\{x\in M(F):\ell(x)>d\}.
\label{b5:o:eq:7.1}
\end{equation}
For a particular cell $c$ of length $d$, define
\begin{equation}
\mathcal W_c=\mathcal U_{>d}\cup X_c.
\label{b5:o:eq:7.2}
\end{equation}
Both sets are invariant opens, as follows from the ordinary Bruhat
closure relation, which only decreases length. The complement of the union of all cells of length greater than $d$
and a chosen cell of length $d$ is a downward-closed finite union. Pulling back to $M$ and adjoining the open $V$ preserves openness.

Under the hypothesis $S|V=0$, every compatible cell outside $V$ is one of the five listed cells. All other cells in the intervening length gaps vanish by intrinsic character incompatibility.

The exact sequence
\[
0\longrightarrow\mathcal D'(Z)\longrightarrow\mathcal D'(W)
 \longrightarrow\mathcal D'(W\setminus Z)
\]
identifies distributions supported on the closed remainder. The NW condition is retained on the restriction of $S$ to the ambient open.

The predecessor is always $V$ together with all cells of greater
original $E_7$ length. The five compatible removals are
\[
52335\ (52)\ \longrightarrow\ \{52209,52215\}\ (50)
\ \longrightarrow\ 51648\ (47)\ \longrightarrow\ 49599\ (42).
\]
For the first four steps, apply Section~\ref{b5:o:sec:6} on the
separate opens $\mathcal W_c$; for the last, use the intrinsic
theorem \eqref{b5:o:eq:5.6}. Between these lengths and below $42$, every
remaining locus is character-incompatible. This gives every
successive predecessor vanishing. The full support ledger is
Table M10-B5P of
\path{anc/model10_calculations.pdf}.

At each spectral step, the original $S$ is restricted only to the open $W_c$. Once its predecessor open is zero, the supported-distribution exact sequence and ordinary character obstruction place that ambient restriction on $C_c$. Section~\ref{b5:o:sec:6} then applies to the \textbf{ambient restriction}, retaining its nilpotent wavefront condition.

The two length-50 cells are isolated in their separate opens \eqref{b5:o:eq:7.2}. Original E7 lengths are removed in the order
\[
{52\ \longrightarrow\ 50,50\ \longrightarrow\ 47\ \longrightarrow\ 42.}
\]
This proves B5 by the stated open/closed support induction.

\subsubsection{Orbit parameters and a transporter test}\label{b5:o:sec:9}

Let
\[
\mathcal E_v=\{(H_0l,H_0r)\in(H_0\backslash L)^2:
             \text{the compatibility equations of \S\ref{b5:o:sec:3} hold}\}.
\]
Use the coordinates $(s,t,A,B)$ of \S\ref{b5:o:sec:1.2}, so this is an explicitly given open subset of a linear-equation locus in twenty coordinates. The small group $\Gamma_v$ of Section~\ref{b5:o:sec:4} acts by
\[
(H_0l,H_0r)\cdot(d_1,d_2)=(H_0ld_1,H_0rd_2).
\]
There is an exact orbit-parameter presentation
\begin{equation}
{\quad (H\times H)\backslash C_c\ =\ \mathcal E_v(F)/\Gamma_v(F).\quad}\label{b5:o:eq:9.1}
\end{equation}
The kernel of $\mathcal P_v\longrightarrow \Gamma_v$ is split unipotent, and its root coordinates give rational lifts of every point in \eqref{b5:o:eq:9.1}.

To write the action in coefficient coordinates, let
\[
d=j(\operatorname{diag}(a,b,c),\operatorname{diag}(g_1,g_2,g_3),\lambda),
\]
so that right multiplication changes these coordinates by
\begin{equation}
(s,t,A,B)\mapsto
\left(\frac ab s,\ b\lambda\det(g_3)t,
       g_1^t A g_2^{-t},\ g_2^t B g_3^{-t}\right).
\label{b5:o:eq:9.2}
\end{equation}
Together with the generators in \S\ref{b5:o:sec:4}, this gives a
complete rational orbit-equivalence test. On $49599$ it reduces to the
$\GL_2$ conjugacy classification of \S\ref{b5:o:sec:5}.

For a character-one transpose test, choose $a_v,b_v\in L(F)$ satisfying
\[
\iota(n_v)=a_v n_v\iota(b_v).
\]
They exist from the small leaf-Weyl comparison in Section~\ref{b5:o:sec:2}; the remaining torus multiplier is defined by this equality in the group and must be retained. A transpose transporter exists for the bare point $l n_v \iota (r)$ exactly when there is $k\in\mathcal P_v(F)$, with projections $(d_1,d_2)$, such that
\begin{equation}
r a_v=h_{01}l d_1,\qquad l b_v=h_{02}r d_2
\label{b5:o:eq:9.3}
\end{equation}
for some $h_{01},h_{02}\in H_0(F)$. Writing $k=u_1d_1$, $k^\sharp=u_2d_2$, the resulting transporter is
\[
\left(h_{01}l u_1^{-1}l^{-1},\ h_{02}r u_2^{-1}r^{-1}\right).
\]
It has character one exactly when
\begin{equation}
\psi\big(p_l(\log u_1)+p_r(\log u_2)\big)=1.
\label{b5:o:eq:9.4}
\end{equation}
The torus graph in \S\ref{b5:o:sec:4} retains every rational central
condition in \eqref{b5:o:eq:9.3}--\eqref{b5:o:eq:9.4}.

\begin{corollary}[The proper-initial open]
\label{m10:proper-open}
$T|_{\Omega_{\mathrm{proper}}}=0$.
\end{corollary}
\begin{proof}
The type-A theorem and Corollary~\ref{m10:o1} clear
$O_1\cup O_3\cup O_4$.
In the free coordinates of $O_6$, the pullback of $O_i$ is exactly
$N_6\times V_i^{(6)}\times N_6$. Thus its descendant satisfies the three-open premise of
Theorem~\ref{m10:b5}, which clears $S_6$. Free descent then clears
$O_6$, completing the four-patch cover of the stated union.
\end{proof}

\paragraph{The identity cell}\label{gk:o:sec:E.1}

Compatibility in \(P\), after removing the outer \(U\)-factors, is
\[
\lambda_\xi\circ\operatorname{Ad}(\ell)=\lambda_\xi,
\]
hence \(\ell\in H_0\). Therefore the compatible part of \(P\) is exactly \(H\).

For \(h\in H\),
\[
(\tau(h)h^{-1},1)\cdot h=\tau(h),
\]
and
\[
\eta(\tau(h)h^{-1},1)=1.
\]
This is an explicit character-one transporter, so \cite[Lemma 2.7]{JLX} applies.

\Needspace{14\baselineskip}
\subsection{The residual exceptional filtration}
\label{m10:residual}
\paragraph{Residual families and their support order}\label{app:residual}
The $28$ residual cells are exhausted by the following disjoint families.
Their full words, lengths, transpose partners and intersection dimensions
are in Table M10-S28 of \path{anc/model10_calculations.pdf} and
\path{anc/model10/residual_e7/signatures_28.json}.
\begin{center}\small
\begin{tabular}{@{}>{\raggedright\arraybackslash}p{.30\linewidth} >{\raggedright\arraybackslash}p{.64\linewidth}@{}}
\toprule
Argument & Cell IDs\\\midrule
Mixed stabilizer & $21092,37733,39171,42519,44094,46017,49001,$
$49596,49812,50122,50415,50767,50996,51205$\\
Ordinary transporter & $2526,14578,32937,49054,49466,50724,51188,51387$\\
Closed-branch sieve and transporter & $51960$\\
Opposite-weight plane and transporter & $52208$\\
Ambient conormal, predecessor $G_{>47}=0$ & $51774,51793$\\
Ambient conormal, predecessor $G_{>38}=0$ & $44670,45681$\\
\bottomrule
\end{tabular}
\end{center}
The first family has intrinsic equivariant vanishing;
the next three have intrinsic odd vanishing. The last two are ambient
support removals. The descending Bruhat induction begins after the
proper-open union has vanished, reaches length $47$ before length $38$,
and retains the same ambient distribution at both conormal steps.
The $51960$ argument includes its closed two-equation locus, and the
$52208$ argument includes both axes and their origin.

The intrinsic arguments below include the closed exceptional branch
on cell $51960$ and the whole opposite-weight plane on cell $52208$.
The next proposition states the predecessor conditions for the two
spectral pairs. 

\begin{proposition}[Conditional residual spectral removal]
\label{m10:residual-conditional}
For the original global NW odd equivariant distribution $T$, put
$\mathcal O_{>d}=\bigcup_{\ell(v)>d}G_{10,v}$. The two spectral packets obey
\[
 T|_{\mathcal O_{>47}}=0
 \Longrightarrow
 T|_{\mathcal O_{>47}\cup G_{10,51774}\cup G_{10,51793}}=0,
\]
\[
 T|_{\mathcal O_{>38}}=0
 \Longrightarrow
 T|_{\mathcal O_{>38}\cup G_{10,44670}\cup G_{10,45681}}=0.
\]
Together with the $24$ intrinsic results, identity transporter and empty
compatibility loci, the initial step
$T|_{\Omega_{\mathrm{proper}}}=0$ suffices to run the descending Bruhat
support filtration to $T=0$.
\end{proposition}

These four cells admit intrinsic odd distributions. Their smooth
compatible supports are excluded by the ambient NW condition: the
Cartan conormal $\gamma_3$, equal to $\alpha_3$ on the Cartan and zero
on root spaces, has nonzero semisimple Killing dual on all six branches.
Corollary~\ref{m10:proper-open} starts the filtration, which reaches
length $47$ before $38$ and supplies both predecessor conditions.

\subsubsection{Structural signatures and the coefficient redundancy}\label{re:o:sec:3}

\paragraph{Signatures}\label{re:o:sec:3.1}\label{resigs:o:sec:3.1}
For a residual word put
$c_U(w)=\dim(\mathfrak u\cap\Ad(n_w^{-1})\mathfrak u)$.
Table M10-S28 of \path{anc/model10_calculations.pdf} records the
lengths, coefficient dimensions, $c_U$ and signed root actions for
all $28$ words in \S\ref{app:residual}. Except at $14578$, the common
Levi-root system is empty: all three positive Levi roots map into
$U$ under both $w$ and $w^{-1}$. At $14578$ the common roots are
$\pm\alpha_7$, their pinned $\SL_2$ commutes with $n_w$, the
coefficient dimension is $12$, and $c_U=35$. These root facts determine
which Gaussian operations are independent.

\paragraph{The section and the positive-Levi operations}\label{re:o:sec:3.2}

Use the specified section \(\ell(p)\), for \(p=(s,t,A,B)\):
\begin{equation}
g_3=I,\quad g_2=B^t,\quad g_1=B^tA^t,
\quad a=st,\quad b=t,\quad
c=(st^2\det A(\det B)^2)^{-1},\quad\lambda=1.
\label{re:o:eq:3.1}
\end{equation}
After removing the outer U factors, points have the form
\begin{equation}
x(p,q)=\ell(p)n_w\ell(q)^{-1}.
\label{re:o:eq:3.2}
\end{equation}
The right coefficient in this convention is \(p(\ell_2^{-1})\), as in \S\ref{m10:datum}.

Right multiplication of a Levi representative by the three positive root groups at nodes \(7,5,2\) acts on its coefficient matrices as
\begin{equation}
(A,B)\longmapsto
\bigl(L(r_1)A L(-r_2),\ L(r_2)B L(-r_3)\bigr),
\quad L(r)=\begin{pmatrix}1&0\\r&1\end{pmatrix}.
\label{re:o:eq:3.3}
\end{equation}
It leaves s and t unchanged. This follows directly from
\(A^t=g_2^{-1}g_1\), \(B^t=g_3^{-1}g_2\).

For the 27 cells without common Levi roots, the operations \eqref{re:o:eq:3.3} are available independently on p and q. For example, inserting a positive Levi element a on the left of \(n_w\) can be moved to the other side as \(n_w^{-1}an_w\in U\); after conjugation by the remaining Levi factor it is still in U. By the section theorem, Levi representatives with the same p differ
by an H0 element. Hence these operations relate points in the same
H-double orbit.

If \(x=k y\) and \(\tau(y)=k_0y\), with \(\eta(k_0)=1\), then
\begin{equation}
\tau(x)=\theta(k)k_0k^{-1}x,
\qquad \eta(\theta(k)k_0k^{-1})=1,
\label{re:o:eq:3.4}
\end{equation}
because the specified tau automorphism preserves eta. Thus a character-one transporter on a core extends to its whole H-saturation.

The torus redundancy is
\begin{equation}
(p,q)\longmapsto
\bigl(p\circ\operatorname{Ad}(w(z)),\ q\circ\operatorname{Ad}(z)\bigr),\qquad z\in T.
\label{re:o:eq:3.5}
\end{equation}
Indeed \(\ell_p w(z)n_w(\ell_qz)^{-1}=\ell_p n_w\ell_q^{-1}\). Returning to the section introduces the $H_0$ corrections used below.

\paragraph{Root-sign convention}\label{re:o:sec:3.3}

Let \(e_i=e_{\alpha_i}\). The specified cocycle is
\[
\epsilon(a,b)=(-1)^{\sum_i a_i b_i+\sum_{i>j,\ i\sim j}a_i b_j}.
\]
For \(\beta\ne\pm\alpha_i\), evaluation of \(\operatorname{Ad}(n_i)e_{\beta}\) uses
\begin{equation}
\begin{cases}
\epsilon(\alpha_i,\beta)e_{s_i\beta},&\langle\beta,\alpha_i^\vee\rangle=-1,\\
\epsilon(-\alpha_i,\beta)e_{s_i\beta},&\langle\beta,\alpha_i^\vee\rangle=1,\\
e_{\beta},&\langle\beta,\alpha_i^\vee\rangle=0.
\end{cases}
\label{re:o:eq:3.6}
\end{equation}
At \(\beta=\pm\alpha_i\), it sends \(e_{\beta}\) to \(e_{-\beta}\). These formulas follow by the two-dimensional root string and \(E_i=e_{\alpha_i},F_i=-e_{-\alpha_i}\). Words are evaluated right to left.

The physical-to-root-basis signs for \((s,t,A_{11},A_{12},A_{21},A_{22},B_{11},B_{12},B_{21},B_{22})\) are
\begin{equation}
(1,1,-1,1,1,-1,-1,1,1,-1).
\label{re:o:eq:3.7}
\end{equation}
In particular all four diagonal trace coordinates have sign $-1$. The character identities below follow from \eqref{re:o:eq:3.6}, \eqref{re:o:eq:3.7}, and the displayed section; the coefficient data record every contributing term.

\subsubsection{Family M: fourteen intrinsic even-zero cells}\label{re:o:sec:4}

\paragraph{The mixed diagonal-unipotent construction}

Let
\[
E_H=E_2+E_5+E_7\in\mathfrak h_0.
\]
If $A,B$ are lower triangular and invertible, then
\begin{equation}
Y(p):=\operatorname{Ad}(\ell(p)^{-1})E_H
=E_2+\frac{B_{22}}{B_{11}}E_5
+\frac{A_{22}B_{22}}{A_{11}B_{11}}E_7.
\label{re:o:eq:4.1}
\end{equation}
This follows by conjugating \(\begin{psmallmatrix}0&1\\0&0\end{psmallmatrix}\) by \(g_3=I,g_2=B^t,g_1=B^tA^t\). An upper triangular matrix with diagonal (a,d) multiplies that root vector under inverse conjugation by d/a. All denominators in \eqref{re:o:eq:4.1} are nonzero on this locus.

For a \textbf{left direction}, assume \(\operatorname{Ad}(n_w^{-1})Y(p)\in\mathfrak u\). Then the following exact one-parameter subgroup fixes \(x=x(p,q)\):
\begin{equation}
h_L(z)=\exp(zE_H),\qquad
h_R(z)=\ell(q)\exp\!\left(z\operatorname{Ad}(n_w^{-1})Y(p)\right)\ell(q)^{-1}.
\label{re:o:eq:4.2}
\end{equation}
Here \(h_L\in H_0\) and \(h_R\in U\). The identity \(h_L(z)x=xh_R(z)\) is obtained by moving \(\exp(zY(p))\) across \(n_w\). Its character is
\begin{equation}
\eta(h_L(z),h_R(z))
=\psi\bigl(zc_w\bigr),\qquad
c_w=-q\bigl(\operatorname{Ad}(n_w^{-1})Y(p)\bigr).
\label{re:o:eq:4.3}
\end{equation}

For a \textbf{right direction}, use
\begin{equation}
h_R(z)=\exp(zE_H),\qquad
h_L(z)=\ell(p)\exp\!\left(z\operatorname{Ad}(n_w)Y(q)\right)\ell(p)^{-1},
\label{re:o:eq:4.4}
\end{equation}
with
\begin{equation}
c_w=p\bigl(\operatorname{Ad}(n_w)Y(q)\bigr).
\label{re:o:eq:4.5}
\end{equation}
These are polynomial unipotent exponentials, valid over every characteristic-zero field under consideration.

If \(c_w\ne0\), the character is nontrivial as z varies over F. The modular character of any algebraic stabilizer is trivial on this unipotent one-parameter subgroup. Consequently \cite[Theorem 2.5]{JLX} applies even though the acting group is not unipotent. Outside the U-compatible locus, the specified U-incompatibility result already applies. Constructibility therefore gives
\begin{equation}
\mathcal D'_\eta(G_w)=0.
\label{re:o:eq:4.6}
\end{equation}
The obstruction is supplied by the displayed subgroup of the full
$H$-stabilizer. \hfill\(\square\)

\paragraph{Finite character and transporter formulas}\label{supp:e7-formulas}
Table M10-M14 of \path{anc/model10_calculations.pdf} gives the exact
$c_w$ in \eqref{re:o:eq:4.3} or \eqref{re:o:eq:4.5} for every mixed
row of \S\ref{app:residual}. The matching equations make the required
matrices lower triangular and invertible; every displayed factor is
therefore nonzero. For example, $c_{49001}=-(B_L)_{11}\ne0$.
These are identities on the entire compatible loci. Together with
the explicit unipotent stabilizers above they prove the fourteen
intrinsic vanishings.

Table M10-T12 of the same supplement gives the Laurent seven-tuples
solving \eqref{re:o:eq:6.3} on the Gaussian cores
\[
 \begin{gathered}
 2526_{DD},14578_{DD},32937_{DD},32937_{AD},49054_{AD},49466_{DD},\\
 50724_{DA},51188_{DA},51387_{AD},51960_{AD},51960_{DD},52208_{DA}.
 \end{gathered}
\]
All denominators are required to be nonzero. The $51960_{DD}$ tuple
uses exactly \eqref{re:o:eq:7.2}; the remaining $52208$ branch is
handled by the whole-plane descent below. The exact symbolic substitutions
are also in \path{anc/model10/residual_e7/tau_torus_certificates.json}.
For instance, the $2526_{DD}$ tuple is
\[
 \left(\frac{B_{R,22}}{s_R},\frac{A_{R,22}}{A_{L,22}},-1,
 -\frac{A_{L,22}s_R}{A_{R,22}B_{R,22}},
 -\frac{A_{R,22}B_{R,22}}{A_{L,22}s_R},
 \frac{A_{L,22}^2s_R}{A_{R,22}^2B_{R,22}},
 \frac{A_{R,22}^2B_{R,22}}{A_{L,22}^2s_R}\right).
\]
Since $G$ is split adjoint, its simple roots form a character basis
of $T$, and each tuple defines an $F$-rational torus element.
The transporter identity, its character and the distributional
conclusion are proved in \S\ref{re:o:sec:6}.

\subsubsection{Gaussian core normal forms}\label{re:o:sec:5}

Use
\[
\mathsf D(a,d)=\begin{pmatrix}a&0\\0&d\end{pmatrix},\quad
\mathsf A(b,c)=\begin{pmatrix}0&b\\c&0\end{pmatrix},\quad
\mathsf V(b,c,e)=\begin{pmatrix}0&b\\c&e\end{pmatrix}.
\]
All displayed diagonal or anti-diagonal entries are nonzero. Only e in \(\mathsf V\) is allowed to vanish. A pattern AD means A is anti-diagonal and B is diagonal.

Here are explicit choices for \eqref{re:o:eq:3.3}. Write
\(A=\begin{psmallmatrix}a&b\\c&d\end{psmallmatrix}\),
\(B=\begin{psmallmatrix}e&f\\g&h\end{psmallmatrix}\).

\begingroup\small\setlength{\tabcolsep}{3pt}
\begin{longtable}{@{}>{\raggedright\arraybackslash}p{0.188\linewidth}>{\raggedright\arraybackslash}p{0.188\linewidth}>{\raggedright\arraybackslash}p{0.188\linewidth}>{\raggedright\arraybackslash}p{0.188\linewidth}>{\raggedright\arraybackslash}p{0.188\linewidth}@{}}
\toprule
Branch & r1 & r2 & r3 & Result \\
\midrule\endhead
b=f=0 & (r2d-c)/a & -g/e & 0 & DD \\
b$\ne$0, f=0 & -d/b & a/b & (g+r2e)/h & AD \\
b=0, f$\ne$0 & (r2d-c)/a & -h/f & e/f & DA \\
b$\ne$0, f$\ne$0 & -d/b & a/b & e/f & AV \\
\bottomrule\end{longtable}\endgroup

In the last row,
\begin{equation}
A'=\begin{pmatrix}0&b\\-\det A/b&0\end{pmatrix},\quad
B'=\begin{pmatrix}0&f\\-\det B/f&h+af/b\end{pmatrix}.
\label{re:o:eq:5.1}
\end{equation}
The extra affine coordinate is \((AB)_{12}/A_{12}\). These rational formulas give the normal forms and their nonzero denominators.

The required core branches are:

\begingroup\small\setlength{\tabcolsep}{3pt}
\begin{longtable}{@{}>{\raggedright\arraybackslash}p{0.470\linewidth}>{\raggedright\arraybackslash}p{0.470\linewidth}@{}}
\toprule
Cell & Left / right core \\
\midrule\endhead
2526 & DD / DD \\
14578 & DD / DD after the shared SL2 operation below \\
32937 & DD / DD if the common A12 is zero; AD / AD otherwise \\
49054 & AD / AD \\
49466 & DD / DD \\
50724 & DA / DA \\
51188 & DA / DA \\
51387 & AD / AD \\
51960 & DD / DD if the common A12 is zero; AD / AD otherwise \\
52208 & DA / DA if the common A12 is zero; AV / AV otherwise \\
44670 / 45681 & AD / DD and DD / AD \\
51774 / 51793 & AD / DA and DA / AD on the open branch; DD / DA and DA / DD on the closed branch \\
\bottomrule\end{longtable}\endgroup

\paragraph{The shared SL2 operation for 14578}

By the matching equations of \S\ref{gk:o:sec:D.2}, the second
column of \(A_R\) is minus the second column of \(A_L\), and
this column is nonzero. The signed root calculation shows that the
positive and negative node-7 root groups commute with \(n_w\):
\[
\operatorname{Ad}(n_w)E_7=E_7,
\qquad \operatorname{Ad}(n_w)F_7=F_7.
\]
The common root SL2 can therefore multiply both A matrices on the left by the same determinant-one matrix. If their common second column is \((b,d)^t\) with \(d\ne0\), use
\(\begin{psmallmatrix}1&-b/d\\0&1\end{psmallmatrix}\).
If \(d=0\), use \(J=\begin{psmallmatrix}0&1\\-1&0\end{psmallmatrix}\).
Both A matrices become lower triangular.

Keep the node-7 operation shared, set \(r_1=0\), and choose independently
\(r_{2,L}=(A_L)_{21}/(A_L)_{22}\) and
\(r_{2,R}=(A_R)_{21}/(A_R)_{22}\), and then choose each
\(r_3=(B_{21}+r_2B_{11})/B_{22}\).
These operations use nodes 5 and 2, which are available independently,
and put both pairs in DD form. Together, the two rational charts
cover the nonzero second column.

\subsubsection{Family T: character-one torus transporters}\label{re:o:sec:6}

\paragraph{The universal construction}\label{re:o:sec:6.1}

For a tau-fixed minimal representative, \(w=w^{-1}\). Set
\[
d_w=n_w^2\in T,\qquad c=b_-r.
\]
For every root \(\beta\),
\begin{equation}
\beta(d_w)=c(w,\beta)c(w,w\beta)\in\{1,-1\}.
\label{re:o:eq:6.1}
\end{equation}
The seven simple-root values determine $d_w$ in the adjoint torus.

For \(y=\ell_p n_w\ell_q^{-1}\), put
\[
a_\tau=c\ell_q d_w^{-1},\qquad b_\tau=c\ell_p.
\]
Then
\(\tau(y)=a_\tau n_w b_\tau^{-1}\).
Because left multiplication by c multiplies the whole coefficient functional by -1, the coefficient pair of these two Levi factors is
\begin{equation}
p'_{\beta}=-\beta(d_w^{-1})q_{\beta},\qquad
q'_{\beta}=-p_{\beta}.
\label{re:o:eq:6.2}
\end{equation}
The formula applies to all ten coordinates, with the physical signs cancelling on each coordinate ratio.

Choose \(z\in T(F)\) satisfying
\begin{equation}
z^{w^{-1}\beta}=p'_{\beta}/p_{\beta},
\qquad z^\beta=q'_{\beta}/q_{\beta}
\label{re:o:eq:6.3}
\end{equation}
for every nonzero core coordinate. Zero coordinates remain zero. The explicit $z$ values are Table M10-T12 of
\path{anc/model10_calculations.pdf}; see \S\ref{supp:e7-formulas}. Define
\begin{equation}
{
 h_L=a_\tau\bigl(\ell_p w(z)\bigr)^{-1},
 \qquad h_R=b_\tau\bigl(\ell_qz\bigr)^{-1}.
}
\label{re:o:eq:6.4}
\end{equation}
The coefficient equalities \eqref{re:o:eq:6.3}, together with the specified identity \(H_0\backslash L\simeq\mathcal O\), prove \(h_L,h_R\in H_0(F)\). Direct multiplication gives
\begin{equation}
h_Lyh_R^{-1}
=a_\tau w(z)^{-1}n_wz b_\tau^{-1}
=a_\tau n_w b_\tau^{-1}=\tau(y).
\label{re:o:eq:6.5}
\end{equation}
Since both elements are in H0, their character is exactly one.

The resulting rational character-one transporter extends, by the
Gaussian operations and \eqref{re:o:eq:3.4}, to every point of the
compatible branch. Applying \cite[Lemma 2.7]{JLX} orbitwise and then
constructibility proves vanishing on its algebraic saturation.

\paragraph{Rationality and the torus certificates}\label{re:o:sec:6.2}

The simple roots are a basis of \(X^*(T)\), since G is split adjoint. Thus
\[
(\alpha_1(z),\ldots,\alpha_7(z))\in(F^\times)^7
\]
can be prescribed arbitrarily. Every entry in Table M10-T12 is a Laurent monomial in coordinates required to be nonzero. 

For a transparent algebraic verification, the certificate forms the integer matrix W whose rows are \(w^{-1}\beta\) on the left and \(\beta\) on the right. It substitutes each displayed seven-tuple in \eqref{re:o:eq:6.3} and checks all equations identically. An integer Smith decomposition is used only to find these monomial formulas. All nonzero Smith invariants of the relevant matrices are 1; the resulting Laurent formulas give the F-rational solutions.

\paragraph{The resulting vanishings}\label{re:o:sec:6.3}

This proves intrinsic odd vanishing on
\begin{equation}
{2526,14578,32937,49054,49466,50724,51188,51387.}
\label{re:o:eq:6.6}
\end{equation}
For 32937 both branches of the common A12 are included. It also proves the open AD branch of 51960 and the DA branch of 52208. The two remaining branches are treated next.

\subsubsection{Cell 51960: the exceptional locus, including singular parameters}\label{re:o:sec:7}

On the closed branch \((A_L)_{12}=(A_R)_{12}=0\), Gaussian reduction gives DD/DD. Retain the following matching equations
\[
(B_L)_{11}=-t_R,\qquad (B_R)_{11}=t_L.
\]
The left and right diagonal-unipotent stabilizers have characters \(\psi(zc_L)\), \(\psi(zc_R)\), where
\begin{equation}
\begin{aligned}
c_L&=-\frac{(A_L)_{11}s_Rt_R+(A_L)_{22}(A_R)_{22}(B_L)_{22}}
{(A_L)_{11}t_R},\\
c_R&=\frac{(A_L)_{22}(A_R)_{22}(B_R)_{22}+(A_R)_{11}t_Ls_L}
{(A_R)_{11}t_L}.
\end{aligned}
\label{re:o:eq:7.1}
\end{equation}
All denominators are nonzero. Off their common zero locus, \cite[Theorem 2.5]{JLX} gives even vanishing. The remaining exact locus is
\begin{equation}
{
\begin{aligned}
s_R&=-\frac{(A_L)_{22}(A_R)_{22}(B_L)_{22}}{(A_L)_{11}t_R},\\
s_L&=-\frac{(A_L)_{22}(A_R)_{22}(B_R)_{22}}{(A_R)_{11}t_L}.
\end{aligned}}
\label{re:o:eq:7.2}
\end{equation}
On this torus locus, the $51960_{DD}$ tuple of Table M10-T12
satisfies all transporter equations after the two displayed
substitutions. Its character-one transporter therefore gives
vanishing on the remaining locus.

The branch $A_{12}\ne0$ follows from \S\ref{re:o:sec:6}. On
$A_{12}=0$, use the open complement of \eqref{re:o:eq:7.2}, followed
by its closed-support injection. Hence
\begin{equation}
{\mathcal D'_{\widetilde\eta}(G_{51960})=0.}
\label{re:o:eq:7.3}
\end{equation}

\subsubsection{Cell 52208: an opposite-weight-plane descendant}\label{re:o:sec:8}

The common $A_{12}=0$ branch is DA/DA and has the transporter
of Table M10-T12. Suppose now
\(x=(A_L)_{12}=(A_R)_{12}\ne0\).
The AV/AV core has the form
\begin{equation}
\begin{gathered}
A_L=\begin{pmatrix}0&x\\a&0\end{pmatrix},\quad
A_R=\begin{pmatrix}0&x\\b&0\end{pmatrix},\\
B_L=\begin{pmatrix}0&-s\\c&u\end{pmatrix},\quad
B_R=\begin{pmatrix}0&r\\d&v\end{pmatrix},\\
s_L=-r,\quad s_R=s,\quad t_L=t_R=t.
\end{gathered}
\label{re:o:eq:8.1}
\end{equation}
The eight variables \(x,a,b,c,d,r,s,t\) are nonzero; u and v are arbitrary elements of F. These are the precise coordinates, including the axes and the origin.

\paragraph{The quotient and the effective action}\label{re:o:sec:8.1}

The two base invariants are
\begin{equation}
I_1=abrs t,\qquad I_2=xcrds t^2.
\label{re:o:eq:8.2}
\end{equation}
The following torus, in simple-root coordinates, makes \(x,a,r,s,c,t\) all equal to 1:
\begin{equation}
\left(
\frac1s,\frac{a}{crt},\frac1t,
 xcrst^2,\frac1{axrst},\ arst,\ 1
\right).
\label{re:o:eq:8.3}
\end{equation}
The remaining nonzero entries become \(b=I_1,d=I_2\). Formula \eqref{re:o:eq:8.3} is an explicit rational section of the base quotient, and identifies the quotient on F-points.

The remaining torus has simple-root values
\begin{equation}
z(\lambda)=(1,\lambda,1,1,\lambda^{-1},1,\lambda).
\label{re:o:eq:8.4}
\end{equation}
It has weights \((-1,+1)\) on \((u,v)\) and fixes the two
nonzero base parameters.

After compensating tau on the normalized base by the torus from \eqref{re:o:eq:6.3}, its action is
\begin{equation}
(u,v)\longmapsto(-v/I_1,-I_1u).
\label{re:o:eq:8.5}
\end{equation}
In coordinates \(U=I_1u,V=-v\), this is simply
\begin{equation}
\lambda:(U,V)\mapsto(\lambda^{-1}U,\lambda V),
\qquad\sigma(U,V)=(V,U).
\label{re:o:eq:8.6}
\end{equation}
All compensating elements at this stage are Levi H0 elements from \eqref{re:o:eq:6.4}, so the character is one.

\paragraph{Frobenius and the density calculation}\label{re:o:sec:8.2}

Gaussian normalization eliminates six independent positive-Levi
coordinates, while \eqref{re:o:eq:8.3} eliminates the six base torus
coordinates. The resulting action map from
\[
K_H\times\bigl((F^\times)^2\times F^2\bigr)
\]
onto this compatible branch is a surjective submersion. In the Gaussian coordinates its fibers have dimension eight: the seven-dimensional common-U stabilizer, and the residual torus \eqref{re:o:eq:8.4}. At the plane origin the latter is in the point stabilizer; away from it it moves the plane coordinate. This keeps the total fiber dimension eight. The uniqueness of the Gaussian cores, the explicit base section, and the saturated torus weights show that these are the complete rational fibers.

More concretely, for \(Y\in\mathfrak u\cap\operatorname{Ad}(n_w^{-1})\mathfrak u\), the common stabilizer frame is
\begin{equation}
Y\longmapsto
\bigl(\operatorname{Ad}(\ell_p n_w)Y,\operatorname{Ad}(\ell_q)Y\bigr).
\label{re:o:eq:8.7}
\end{equation}
The specified compatibility equations make eta trivial on this whole unipotent fiber. Under \eqref{re:o:eq:8.4} the frame transforms by its root weights. Its seven roots, in simple-root coordinates, are
\begin{equation}
\begin{array}{lll}
(0,0,0,0,1,1,1),&(0,0,1,0,0,0,0),&(0,1,0,1,1,0,0),\\
(0,1,1,1,1,0,0),&(1,0,0,0,0,0,0),&(1,0,1,0,0,0,0),\\
(1,1,1,1,1,0,0).
\end{array}
\label{re:o:eq:8.8}
\end{equation}
Every one pairs to zero with \((0,1,0,0,-1,0,1)\). Thus the determinant character on the common-U density is one. The torus itself is unimodular, and its determinant on the plane is \(\lambda^{-1}\lambda=1\).

To check unimodularity of K\_H, note that U is unipotent and the
determinant of the H0-action on its Lie algebra is an algebraic
character of PGL2, hence trivial. These frames describe the fibers as the common unipotent group together with the residual torus, acting with the indicated section cocycle. The relative density character for twisted Frobenius descent is therefore trivial. With product Haar densities it gives an injective descent to ordinary distributions on
\[
(F^\times)^2\times F^2
\]
that are invariant under the first action in \eqref{re:o:eq:8.6} and odd under the second. The base is fixed. The same calculation shows compatibility with tau: the plane Jacobian has absolute value one, the fiber density has trivial torus weight, and the compensators have eta=1. Positive densities may be chosen tau-invariant.

Lemma~\ref{m8:s:c:plane}, with the two nonzero base parameters
passive, now kills the descendant on the whole plane, including its axes
and origin. Injectivity of Frobenius descent proves intrinsic odd
vanishing on this branch and hence on all of cell $52208$.

\subsubsection{The final two exchanged packets: stabilizers and residual loci}\label{re:o:sec:9}

Having established the 24 intrinsic conclusions, we turn to the
remaining packets:
\begin{equation}
\mathcal S_1=(44670,45681),\qquad
\mathcal S_2=(51774,51793).
\label{re:o:eq:9.1}
\end{equation}

\paragraph{Complete stabilizer calculation on the cores}\label{re:o:sec:9.1}

To determine the Lie stabilizer modulo the common U-intersection,
we use the following small linear system.

For \(Z_L,Z_R\in\mathfrak{sl}_2\), expand
\[
\operatorname{Ad}(\ell_p^{-1})Z_L,
\quad \operatorname{Ad}(\ell_q^{-1})Z_R
\]
in the six Levi root lines and the Cartan. Equate the Cartan parts after w, discard any root component whose partner lies outside \(\mathfrak p\), and equate the common Levi root components. The components mapping between a Levi root and U determine the corresponding U terms uniquely. The common U$\cap$wU variables remain free. This accounts for \textbf{every} component of
\(\mathfrak h\cap\operatorname{Ad}(x)\mathfrak h\).

The residual character is the left U phase minus the right U phase. The outcome is:

\begingroup\small\setlength{\tabcolsep}{3pt}
\begin{longtable}{@{}>{\raggedright\arraybackslash}p{0.16\linewidth}>{\raggedright\arraybackslash}p{0.16\linewidth}>{\raggedright\arraybackslash}p{0.16\linewidth}>{\raggedright\arraybackslash}p{0.46\linewidth}@{}}
\toprule
Cell/core & Common U dimension & Extra stabilizer & Extra phase \\
\midrule\endhead
44670 AD/DD & 19 & right diagonal Ga & 0 \\
45681 DD/AD & 19 & left diagonal Ga & 0 \\
51774 AD/DA & 10 & none & 0 \\
51774 DD/DA & 10 & left diagonal Ga & \(\displaystyle\frac{-t_R\bigl((A_L)_{11}s_R-(A_L)_{22}(A_R)_{22}\bigr)}{(A_L)_{11}(B_L)_{11}}\) \\
51793 DA/AD & 10 & none & 0 \\
51793 DA/DD & 10 & right diagonal Ga & \(\displaystyle\frac{(B_R)_{22}\bigl((A_L)_{22}(A_R)_{22}+(A_R)_{11}s_L\bigr)}{(A_R)_{11}(B_R)_{11}}\) \\
\bottomrule\end{longtable}\endgroup

Thus the closed branches of the second packet are further confined to
\begin{equation}
\begin{aligned}
51774:&\quad (A_L)_{12}=0,
\quad (A_L)_{11}s_R=(A_L)_{22}(A_R)_{22},\\
51793:&\quad (A_R)_{12}=0,
\quad (A_L)_{22}(A_R)_{22}+(A_R)_{11}s_L=0.
\end{aligned}
\label{re:o:eq:9.2}
\end{equation}
The complements are intrinsically even-zero. The surviving stabilizers are unipotent: the coefficient-torus weight
matrices have rank seven and saturated image, and the remaining
directions in the table are $\mathbb G_a$. Their character is trivial
on \eqref{re:o:eq:9.2}, and their modular character is one.

\paragraph{Frobenius descendants and the stabilizer character}\label{re:o:sec:9.2}

After Gaussian reduction, the surviving parameter spaces are tori.
Quotienting by the effective rank-seven torus, whose character lattice
is saturated, gives a geometric quotient of dimension

\begin{itemize}
\item one for the first packet;
\item two on the open AD/DA and DA/AD branches of the second packet;
\item one on each boundary \eqref{re:o:eq:9.2}.

\end{itemize}
Each fiber is a single H-double orbit, with equivariant density
obtained from \(\eta(k)d\dot k\) on \(K_H/(K_H)_x\). This
density exists because the stabilizer is unipotent, unimodular and
eta-trivial. A rational section of the split-torus quotient trivializes this density line. Thus the Frobenius descendants are ordinary distributions on \(F^\times\), or on \((F^\times)^2\), with the two coarse cells exchanged by tau.

Intrinsic odd nonvanishing on the packets can in fact be seen on closed orbits. For 44670 one may take
\begin{equation}
p=(1,1,\begin{pmatrix}0&1\\1&0\end{pmatrix},I),\quad
q=(1,1,\operatorname{diag}(-1,1),\operatorname{diag}(1,-1)).
\label{re:o:eq:9.3}
\end{equation}
For 51774 use the closed boundary point
\begin{equation}
p=(1,1,I,I),\quad
q=(1,1,I,\begin{pmatrix}0&-1\\1&0\end{pmatrix}).
\label{re:o:eq:9.4}
\end{equation}
These satisfy the matching equations of \S\ref{gk:o:sec:D.2} and the stated extra equations on each locus. On the first locus and on the closed locus \eqref{re:o:eq:9.2}, the Gaussian quotient is a torus and its subsequent torus quotient has closed fibers. Hence these H-double orbits are closed in the indicated compatible locus and in their coarse cell. Their nonzero orbital eta-distributions push forward by closed embeddings. Subtracting the tau image gives a nonzero odd distribution on the two-cell packet; the cells are clopen in their same-length shell.

The orbital distributions have nonnilpotent wavefront, so their
intrinsic nonvanishing is consistent with the ambient exclusion below.

\subsubsection{Conormal exclusion for the four spectral families}\label{re:o:sec:10}

\paragraph{Monomial representatives and the conormal}\label{re:o:sec:10.1}

For a monomial core write
\begin{equation}
\ell_p=t_p m_p,\qquad \ell_q=t_qm_q,
\quad x=t\,n_v,
\quad v=m_pwm_q^{-1},
\label{re:o:eq:10.1}
\end{equation}
where \(n_v=m_pn_wm_q^{-1}\). For a diagonal/anti-diagonal core, the WL factors are:

\begingroup\small\setlength{\tabcolsep}{3pt}
\begin{longtable}{@{}>{\raggedright\arraybackslash}p{0.188\linewidth}>{\raggedright\arraybackslash}p{0.188\linewidth}>{\raggedright\arraybackslash}p{0.188\linewidth}>{\raggedright\arraybackslash}p{0.188\linewidth}>{\raggedright\arraybackslash}p{0.188\linewidth}@{}}
\toprule
Cell/core & m\_p & m\_q & \(v^{-1}\alpha_3\) & \(\alpha_3(t)\) \\
\midrule\endhead
44670 AD/DD & s7 & 1 & $\alpha_{3}$ & 1 \\
45681 DD/AD & 1 & s7 & $\alpha_{3}$ & 1 \\
51774 AD/DA & s7 & s7s5 & $\alpha_{2}$+$\alpha_{4}$ & 1 \\
51774 DD/DA & 1 & s7s5 & $\alpha_{2}$+$\alpha_{4}$ & 1 \\
51793 DA/AD & s7s5 & s7 & $\alpha_{2}$+$\alpha_{4}$ & -1 \\
51793 DA/DD & s7s5 & 1 & $\alpha_{2}$+$\alpha_{4}$ & -1 \\
\bottomrule\end{longtable}\endgroup

The root maps and torus constants are exact; all seven torus coordinates appear in \texttt{\detokenize{spectral_conormal_certificates.json}}. For example, the simple-root coordinates of a diagonal Levi factor are
\begin{equation}
(a/b,\ d_5/d_6,\ b\lambda d_5d_6,\ d_4/d_5,
 d_3/d_4,\ d_2/d_3,\ d_1/d_2),
\label{re:o:eq:10.2}
\end{equation}
where \(d_1,\dots,d_6\) are its physical GL6 diagonal entries. Factoring each physical anti-diagonal block as a diagonal matrix times \(\begin{psmallmatrix}0&1\\-1&0\end{psmallmatrix}\) and using \eqref{re:o:eq:3.1} yields the displayed constants. The same constants hold on \eqref{re:o:eq:9.2}.

Let \(\gamma_3\) be the functional on \(\mathfrak g\) equal to \(\alpha_3\) on the Cartan and zero on every root space. Its Killing-dual is a nonzero multiple of the semisimple coroot \(h_3\).

The H0 torus has simple-root values
\begin{equation}
(1,z,1,z^{-1},z,z^{-1},z),
\label{re:o:eq:10.3}
\end{equation}
so \(\alpha_3\) vanishes on its Lie algebra. It also vanishes on its v-translate: the fourth column of the table is either \(\alpha_3\) or \(\alpha_2+\alpha_4\), and both vanish on \eqref{re:o:eq:10.3}.

The tangent to the H$\times$H saturation of the monomial core, after right translation to the Lie algebra, is the sum of
\begin{equation}
\mathfrak h+\operatorname{Ad}(x)\mathfrak h
\quad\text{and the tangent to the core torus.}
\label{re:o:eq:10.4}
\end{equation}
The root-space parts of the first summand are annihilated by \(\gamma_3\), its two Cartan directions are annihilated by \eqref{re:o:eq:10.3} and the displayed root maps, and the core tangent is annihilated because \(\alpha_3(t)=\pm1\) is constant. Hence
\begin{equation}
{\gamma_3\in \CN^G_{Y,x},\qquad\gamma_3\notin\mathcal N.}
\label{re:o:eq:10.5}
\end{equation}
At an arbitrary point of the saturation its transported conormal is still nonnilpotent: translations and adjoint actions preserve nilpotence.

\paragraph{Smoothness and submersion, with the dimensions checked}\label{re:o:sec:10.2}

The compatibility equations of \S\ref{gk:o:sec:D.2} are linear equations inside an open subset of coefficient space. The Gaussian changes are rational isomorphisms on their specified nonzero-pivot charts. The exceptional conditions \eqref{re:o:eq:9.2} solve a nonzero scalar variable uniquely and hence define smooth torus loci.

The coefficient cover of a P-double cell is a smooth homogeneous quotient; the saturated compatibility conditions descend because their inverse image is exactly the given coefficient condition. Thus its codimension is the coefficient codimension. The dimensions are:

\begingroup\small\setlength{\tabcolsep}{3pt}
\begin{longtable}{@{}>{\raggedright\arraybackslash}p{0.130\linewidth}>{\raggedright\arraybackslash}p{0.130\linewidth}>{\raggedright\arraybackslash}p{0.130\linewidth}>{\raggedright\arraybackslash}p{0.130\linewidth}>{\raggedright\arraybackslash}p{0.130\linewidth}>{\raggedright\arraybackslash}p{0.130\linewidth}>{\raggedright\arraybackslash}p{0.130\linewidth}@{}}
\toprule
Locus Y & Dimension of coarse cell & Coefficient codimension & dim Y & dim H$\times$H orbit & Core torus dimension & Cartan overlap \\
\midrule\endhead
S1 compatible locus & 114 & 7 & 107 & 106 & 3 & 2 \\
S2 open branch & 123 & 5 & 118 & 116 & 4 & 2 \\
S2 boundary \eqref{re:o:eq:9.2} & 123 & 7 & 116 & 115 & 3 & 2 \\
\bottomrule\end{longtable}\endgroup

For example the coarse-cell dimension is
\(2\dim P-(13+c_U)=133-c_U\) here: the parabolic intersection has its seven-dimensional Cartan, six mixed Levi-root directions, and the common U. The coefficient cover has a constant-dimensional smooth fiber, so subtracting the coefficient codimension is legitimate.

Using the full stabilizer calculation in \S\ref{re:o:sec:9} and
the independence of the two Cartan directions in
\eqref{re:o:eq:10.4}, the last three columns give exactly dim Y
at every core point. Thus the action map onto Y is a surjective
submersion. Thus \eqref{re:o:eq:10.4} is the full tangent space of the smooth support used in the wavefront argument.

\paragraph{Proposition~\ref{m10:residual-conditional} --- the conditional spectral-removal theorem}

Put
\begin{equation}
\mathcal O_{>d}=\bigcup_{\ell(v)>d}G_v.
\label{re:o:eq:10.6}
\end{equation}
This is an open set: its complement is a downward-closed finite Bruhat union, using the specified closure formula. If a packet \(\mathcal S\) has length d, then
\begin{equation}
M=\mathcal O_{>d}\ \sqcup\!\bigsqcup_{w\in\mathcal S}G_w
\label{re:o:eq:10.7}
\end{equation}
is also an open set. It is tau-stable for either packet in \eqref{re:o:eq:9.1}.

For the original global NW odd equivariant distribution T,
\begin{equation}
{
T|_{\mathcal O_{>38}}=0
\Longrightarrow
T|_{\mathcal O_{>38}\cup G_{44670}\cup G_{45681}}=0,
}
\label{re:o:eq:10.8}
\end{equation}
and
\begin{equation}
{
T|_{\mathcal O_{>47}}=0
\Longrightarrow
T|_{\mathcal O_{>47}\cup G_{51774}\cup G_{51793}}=0.
}
\label{re:o:eq:10.9}
\end{equation}
The same proof works with smaller local predecessors that isolate the indicated packet.

\textbf{Proof.} The premise and the closed-support injection first put the restriction of the original T to M on the packet. The specified U-incompatibility result confines it to the compatible locus. For S2 the extra stabilizers in \S\ref{re:o:sec:9} further remove the complement of \eqref{re:o:eq:9.2} on its closed branch. The ambient restriction $T|_M$ retains its NW condition.

For S1 the resulting support lies in the smooth Y of the first row
of the dimension table, with nonnilpotent conormal
\eqref{re:o:eq:10.5} at every point. Since the original T retains
nilpotent wavefront on M, Proposition~\ref{cf:conormal} applies with
the ambient witness \eqref{re:o:eq:10.5} at every support point and
gives zero.

For S2, remove the closed zero-A12 branches from M. On the resulting
open, support lies in the smooth open-branch Y, so the same conormal
argument gives zero. With support now confined to the zero-A12
branches, apply the extra character sieve to obtain
\eqref{re:o:eq:9.2}. These are closed smooth submanifolds of the
current M and are eliminated by the same conormal argument. Thus T
vanishes on M. \hfill\(\square\)

The exclusion uses the semisimple conormal $\gamma_3$, open locality
and the nilpotent wavefront of the same ambient $T$.

\subsubsection{Filtration and the relative residual theorem}\label{re:o:sec:11}

To insert the intrinsic results into the Bruhat filtration, use each
locally closed cell or its open/closed refinement. Within a fixed
length shell, the cells are clopen because their boundaries have
strictly smaller length. Since tau preserves length, we treat
exchanged cells as packets.

For an open embedding $j:V\hookrightarrow X$ and its closed
complement $i:Z\hookrightarrow X$, use the sequence
\begin{equation}
0\longrightarrow\mathcal D'_\chi(Z)
\xrightarrow{i_*}\mathcal D'_\chi(X)
\longrightarrow\mathcal D'_\chi(V).
\label{re:o:eq:11.1}
\end{equation}
When the image on $V$ is zero, exactness identifies the distribution
with the pushforward $i_*S$ of a distribution on $Z$.

For the relative residual theorem of \S\ref{m10:residual}, assume its stated premise and work with the increasing opens
\[
\Omega_{\rm proper}\cup\mathcal O_{>d},
\]
processing d in decreasing order and refining each same-length shell into individual fixed cells and exchanged packets. Outside \(\Omega_{\rm proper}\), every cell is, by the specified input, either U-incompatible, the already-solved identity cell, or one of these 28. The 24 intrinsic results give zero kernels at their steps. At the two spectral packets the previously processed larger-length open is already zero, so the proof of Proposition~\ref{m10:residual-conditional} applies to the \textbf{same original global T}. The identity step is the specified one. This proves T=0.

The proper-Levi theorems supply the initial premise. The individual
branch statements then enter in descending Bruhat order, including
every intervening incompatible cell.

\subsection{The distribution theorem and the product bound}
\begin{theorem}[Model-10 distribution theorem]
\label{m10:distribution}
Every $\tau_{10}$-odd, $(H_{10}\times H_{10},\eta_{10})$-equivariant
distribution $T\in\NW(G_{10})$ is zero. Here the nilpotent cone is
$\cN(\mathfrak e_7^*)$.
\end{theorem}
\begin{proof}
The identity compatible locus is $H_{10}$ and has the character-one
transporter $(\tau_{10}(h)h^{-1},1)$. Proposition~\ref{m10:finite-localization},
Corollary~\ref{m10:proper-open} and Proposition~\ref{m10:residual-conditional}
then give zero by the stated open/closed support filtration.
\end{proof}

\begin{corollary}[The GK product bound]
\label{m10:gp2}
For every irreducible admissible smooth complex representation $\pi$ of
$G_{10}(\F)$,
\[
 m_{10}(\pi)m_{10}^{\vee}(\pi)\le1,\qquad
 m_{10}(\pi)=\dim\Hom_{H_{10}}(\pi,\xi_{10}),\quad
 m_{10}^{\vee}(\pi)=\dim\Hom_{H_{10}}(\pi^\vee,\xi_{10}^{-1}).
\]
\end{corollary}
The common coefficient construction \eqref{cf:center-span} and
Theorem~\ref{cf:nilpotent} supply the NW hypothesis for the GK
application. The individual bound follows from the representation-theoretic
argument in the next section.

\part{From the \texorpdfstring{$E_7$}{E7} product bound to multiplicity one}\label{part:5}
\section{From the exceptional product bound to individual multiplicity one}
\label{ep:section}

\paragraph{Proof strategy}
We obtain the individual bound by identifying the two period
dimensions in the duality cases and proving period vanishing in
the remaining case. Hermitian duality handles unitarizable
representations. For the other Langlands quotients, rational duality
handles type-$A/D$ inducing Levis and generic tempered $E_6$ data
(\S\ref{ep:langlands}). For nongeneric tempered data on the maximal
$E_6$ Levi, we prove that the entire positive standard module has zero
Model-10 period (Theorem~\ref{ep:e6-vanishing}). A Mackey filtration
leaves an ordinary Whittaker quotient and one open contribution,
which reduces to a normalized $D_5$ descendant;
tempered exponents and the compact-boundary operator imply zero period.

Throughout this part, normalized induction is realized by
$f(pg)=\delta_P(p)^{1/2}\rho(p)f(g)$. Unless explicitly normalized,
$J_{U,\chi}(\rho)$ is the quotient by the vectors
$\rho(u)v-\chi(u)v$; a convention using $(\rho\otimes\chi)_U$
would therefore use the inverse character. A star on a vector space
means its algebraic dual. The residual $\PGL_2$ modules are treated
in the smooth category.

\subsection{Langlands classification and rational duality}
\label{ep:langlands}
For tempered inducing data with positive parameter, we use
\cite[Theorem~3.5(i), p.~396]{Konno}. We give the quotient and
duality arguments below.
Let $Q_*=M_*N_*$ be the standard maximal parabolic omitting $\alpha_7$.
In split adjoint $E_7$, the Levi $M_*$ has a simply connected
derived group of type $E_6$ and a connected split central torus. More generally,
for every standard Levi $M_J$,
\[
 X^*(Z(M_J))=Q(E_7)/\mathbb ZJ.
\]
Since the simple roots form a basis of $Q(E_7)$, this quotient is torsion-free.
Among the $127$ proper standard subsets, $\{1,2,3,4,5,6\}$ is the unique
$E_6$ subdiagram; every other simple component has type $A$ or $D$.

\begin{proposition}[Rational duality and the Langlands quotient]
\label{ep:dual-langlands}
Let $P=MN$ be a proper standard parabolic and
$\pi=J_P^{G_{10}}(\sigma e^\nu)$ its Langlands quotient, with $\sigma$
irreducible tempered and $\nu$ positive. If $n_G$ represents the longest element $w_{\mathrm{long}}=-1$ of the Weyl group of $G_{10}$ and
$\vartheta=\Ad(n_G)|_M$, then
\[
 \pi^\vee\simeq
 J_P^{G_{10}}\bigl((\sigma^\vee)^\vartheta e^\nu\bigr).
\]
If all simple factors of $M$ have type $A$ or $D$, then
$\sigma^\vartheta\simeq\sigma^\vee$. This conclusion also holds for
$M=M_*$ when $\sigma$ is ordinary nondegenerate Whittaker generic.
In these cases $\pi\simeq\pi^\vee$.
\end{proposition}

The rational-conjugacy proof and the quotient formula follow below.
\begin{proposition}[Connected-center duality in types $A$ and $D$]
\label{ep:connected-center-duality}
Let $M$ be a split connected reductive $F$-group with split connected
center and with every adjoint simple factor of type $A$ or $D$.
Let $C_M$ be an $F$-automorphism inverting a fixed split maximal torus.
Then every irreducible admissible smooth complex representation $\sigma$
of $M(F)$ satisfies $\sigma^{C_M}\simeq\sigma^\vee$.
\end{proposition}
The proof first lifts rational inverse-conjugacy through a connected
central quotient, and then supplies the conjugators in types $A$ and $D$.

\subsubsection{The rational lifting argument}\label{ep:o:sec:3}

Consider a central quotient of split connected reductive groups
\[
1\longrightarrow Z\longrightarrow M\xrightarrow{q}\overline M\longrightarrow1,
\]
with \(Z\) a split torus and compatible automorphisms \(C_M,C_{\overline M}\), with \(C_M\) inverting a fixed split maximal torus.

Choose \(g\in M(F)\) in the common open set where both \(g\) and \(q(g)\) are strongly regular and the asserted rational inverse-conjugacy downstairs holds. Write
\[
S=Z_M(g),\qquad \overline S=Z_{\overline M}(q(g)).
\]
These are connected tori and
\begin{equation}
1\longrightarrow Z\longrightarrow S\longrightarrow\overline S\longrightarrow1
\label{ep:o:eq:3.1}
\end{equation}
is exact. This is the inverse image of the maximal centralizer torus under the central quotient; the roots of the quotient and of \(M\) agree.

Torus inversion gives a geometric conjugator \(\gamma\) satisfying
\[
\gamma g^{-1}\gamma^{-1}=C_M(g).
\]
The cocycle \(a_\tau=\gamma^{-1}\tau(\gamma)\) has values in \(S\). A rational conjugator \(\overline h\) downstairs satisfies
\[
\overline h^{-1}q(\gamma)\in\overline S(\overline F).
\]
Writing this element as \(u\), one gets
\[
q(a_\tau)=u^{-1}\tau(u),
\]
so the image class in \(H^1(F,\overline S)\) is trivial.

The cohomology groups of tori are abelian. Exactness of
\eqref{ep:o:eq:3.1} and \(H^1(F,Z)=0\) give an injection
\[
H^1(F,S)\hookrightarrow H^1(F,\overline S).
\]
Thus \([a]=0\), and modifying \(\gamma\) on the right by the
corresponding element of \(S(\overline F)\) makes it F-rational.

On this dense open set,
\[
\Theta_{\sigma^{C_M}}(g)=\Theta_\sigma(C_M(g))
=\Theta_\sigma(g^{-1})=\Theta_{\sigma^\vee}(g).
\]
A proper algebraic closed subset has no p-adic open interior here. Local constancy extends the equality throughout the regular semisimple set, and local integrability identifies the character distributions. Character determination then gives \(\sigma^{C_M}\simeq\sigma^\vee\).

For the analogous lifting argument in the GSpin setting see \cite[Lemmas 3.2--3.3 and (3.5)]{MN}. Character determination is recalled in \cite[Corollary 2.3.2]{Casselman}.

\subsubsection{Rational conjugacy in types A and D}\label{ep:o:sec:4}

\paragraph{Type A}\label{ep:o:sec:4.1}

A matrix over \(F\) is similar over \(F\) to its transpose. Hence \(g^{-t}\) is rationally conjugate to \(g^{-1}\) in \(\GL_n(F)\). The quotient to \(\PGL_n\) is surjective on F-points because its kernel is \(\mathbb G_m\). Projecting the conjugacy and applying the character argument proves the needed type-A statement.

\paragraph{Type D and the determinant check}\label{ep:o:sec:4.2}

Let
\[
J=\begin{pmatrix}0&I_n\\I_n&0\end{pmatrix},\qquad
 g^tJg=\mu J,\qquad \det g=\mu^n.
\]
Work on the dense open subset of \(\GSO(J)\) where \(g\) has separable characteristic polynomial and no eigenvalue is fixed by \(\alpha\mapsto\mu/\alpha\). Then \(A=F[g]\) is an étale algebra of dimension \(2n\), and the representation space is a rank-one \(A\)-module. The adjoint involution is
\[
\overline g=\mu g^{-1}.
\]
For a suitable \(c_A\in A^\times\) fixed by this involution, the symmetric form is
\[
B(v,w)=\operatorname{Tr}_{A/F}(c_Av\overline w).
\]
Indeed, the étale trace pairing represents \(B(1,-)\), and symmetry forces \(\overline{c_A}=c_A\).

The F-linear involution \(\kappa(v)=\overline v\) preserves \(B\). Its plus and minus eigenspaces both have dimension \(n\), so
\[
\det\kappa=(-1)^n,\qquad
\kappa g\kappa^{-1}=\mu g^{-1}.
\]
The matrix \(J\) is another isometry, also with determinant \((-1)^n\). Therefore
\begin{equation}
h=J\kappa\in \SO(J)(F),
\quad hgh^{-1}=J\mu g^{-1}J=g^t,
\quad g^{-t}=hg^{-1}h^{-1}.
\label{ep:o:eq:4.1}
\end{equation}
The two determinant signs ensure that the conjugator lies in the
connected group \(\SO(J)\).

The quotient \(\GSO(J)\to \operatorname{PSO}(J)\) has the split scalar torus as kernel, so is surjective on F-points. Conjugacy descends. The construction uses division by 2 in the characteristic-zero field and applies in arbitrary residue characteristic.

\paragraph{The connected-center Levi}\label{ep:o:sec:4.3}

On a split adjoint group, two torus-inverting F-automorphisms differ by an automorphism fixing the torus pointwise. Its simple-root scaling constants determine an element of the adjoint torus over \(F\), since simple roots are a character basis. Thus every such torus-inverting automorphism has the rational duality established above on an A/D adjoint product.

Apply Section~\ref{ep:o:sec:3} to
\[
M\longrightarrow M_{\rm ad}
\]
using the split connected center of \(M\). This proves
Proposition~\ref{ep:connected-center-duality} for all the stated
central forms. In the D6 case, rational lifting uses the scalar-torus
quotient of the actual Levi onto the adjoint D6 group.

\subsubsection{Root data of the proper standard Levis}\label{ep:o:sec:5}

For adjoint E7,
\[
X^*(T)=Q(E_7),\qquad
X^*(Z(M_I))=Q(E_7)/\mathbb Z I.
\]
Since the simple roots are a basis of \(Q(E_7)\), this quotient is free for every subset \(I\). The center is therefore split and connected.

The integral root calculation gives
\[
\det A_{E_7}=2,\qquad
X_*(T)=Q^\vee(E_7)+\mathbb Z\frac{H_2+H_5+H_7}{2}.
\]
Intersecting with the rational span of the coroots in \(I\) yields
\begin{equation}
[X_*(T\cap M_I^{\rm der}):Q_I^\vee]
=\begin{cases}
2,&\{2,5,7\}\subset I,\\
1,&\text{otherwise}.
\end{cases}
\label{ep:o:eq:5.1}
\end{equation}
The lattice index records the central form of the embedded derived group.

Among the \(2^7-1=127\) proper standard subsets, those with a D or E
factor are:

\begingroup\small\setlength{\tabcolsep}{3pt}
\begin{longtable}{@{}>{\raggedright\arraybackslash}p{0.188\linewidth}>{\raggedright\arraybackslash}p{0.188\linewidth}>{\raggedright\arraybackslash}p{0.188\linewidth}>{\raggedright\arraybackslash}p{0.188\linewidth}>{\raggedright\arraybackslash}p{0.188\linewidth}@{}}
\toprule
\(I\) & Root type & Positive roots & Center dimension & Derived lattice index \\
\midrule\endhead
2,3,4,5 & D4 & 12 & 3 & 1 \\
1,2,3,4,5 & D5 & 20 & 2 & 1 \\
2,3,4,5,6 & D5 & 20 & 2 & 1 \\
2,3,4,5,7 & D4+A1 & 13 & 2 & 2 \\
1,2,3,4,5,6 & E6 & 36 & 1 & 1 \\
1,2,3,4,5,7 & D5+A1 & 21 & 1 & 2 \\
2,3,4,5,6,7 & D6 & 30 & 1 & 2 \\
\bottomrule\end{longtable}\endgroup

All other entries have only A factors or are tori. The unique
exceptional standard Levi \(M_*\) has derived index one, hence
simply connected derived group E6. The D6 derived group has semispin
central form: on nodes 2 through 7, the vector coweight has
coroot-lattice coset \( (H_2+H_3)/2 \), whereas the extra coset is
\( (H_2+H_5+H_7)/2 \), and their difference is not integral.

\subsubsection{Langlands quotients and generic E6 data}\label{ep:o:sec:6}

Let \(P=MN\), \(\sigma\) irreducible tempered, and \(\nu\) strictly P-positive. The canonical opposite-module morphism is
\[
J'_{\bar P|P}(\sigma e^\nu):i_P^G(\sigma e^\nu)
\longrightarrow i_{\bar P}^G(\sigma e^\nu).
\]
Both sides have the same inducing representation and parameter.
Proposition~\ref{app:konno-image-independent} constructs this nonzero
morphism from the normalized Jacquet projection and Frobenius
reciprocity, using a compactly supported open-cell integral, and proves
that its image is the unique irreducible quotient \(\pi\) described in
\cite[Corollary~3.2(ii), pp.~393--394]{Konno}. Hence \(\pi\) embeds into the opposite standard module.

Dualizing this embedding gives a surjection
\[
i_{\bar P}^G(\sigma^\vee e^{-\nu})\twoheadrightarrow\pi^\vee.
\]
The parameter \(-\nu\) is positive for \(\bar P\). Choose a rational longest-Weyl representative \(n_G\), and set \(\vartheta=\operatorname{Ad}(n_G)|_M\). Since the longest Weyl element of E7 is \(-1\), conjugation takes \(\bar P\) to \(P\) and \(-\nu\) to \(+\nu\). Therefore
\begin{equation}
{
\bigl(J_P^G(\sigma e^\nu)\bigr)^\vee
\simeq J_P^G((\sigma^\vee)^\vartheta e^\nu).
}
\label{ep:o:eq:6.1}
\end{equation}
Thus the quotient formula follows by dualizing the embedding into the
opposite standard module.

If \(\sigma^\vartheta\simeq\sigma^\vee\), then \((\sigma^\vee)^\vartheta\simeq\sigma\), since \(\vartheta^2\) is inner. Thus \(\pi\) is self-contragredient. Section~\ref{ep:o:sec:4} establishes this for every A/D Levi. Together with the character transport in \S\ref{ep:hermitian} and the specified product bound, it gives GP1 for all such Langlands quotients, including reducibility parameters.

For generic tempered \(\sigma\) on \(M_*(F)\), both Whittaker hypotheses of \cite[Proposition~7.5]{Yang} hold. Indeed, tempered \(\sigma\) is unitarizable, so its invariant positive Hermitian form identifies \(\overline\sigma\) with \(\sigma^\vee\). If \(\ell\) is a nonzero \(\psi\)-Whittaker functional on \(\sigma\), then \(\overline v\mapsto\overline{\ell(v)}\) is a complex-linear \(\psi^{-1}\)-Whittaker functional on \(\overline\sigma\), since \(\psi\) is unitary. The smooth-dual identification and this character calculation are detailed in \S\ref{repair:shalika-criterion}. Together with the bi-Whittaker distribution invariance proved there, \cite[Proposition~7.5 and Corollary~7.7]{Yang} gives \(\sigma^\iota\simeq\sigma^\vee\) for its F-pinning. Explicitly, for a rational longest-Weyl representative $n_M$ in $M_*$, set $C^{\prime}=\operatorname{Ad}(n_M)\circ\iota_{M_*,\mathcal P}$. On the split maximal torus, $C^{\prime}(t)=w_{0,M}(w_{0,M}(t^{-1}))=t^{-1}$. Its comparison with \(\vartheta\) is conjugation by an F-point of the adjoint torus. That point lifts to \(M_*(F)\), because the center is a split torus. Hence \(\sigma^\vartheta\simeq\sigma^\vee\) here too.

The same torus-lifting argument shows that the rational nondegenerate
Whittaker data form one orbit. Genericity is required here of the
tempered inducing representation \(\sigma\); the Langlands quotient
\(\pi\) may be nongeneric.

\subsubsection{Character transport and the tempered case}\label{ep:hermitian}
The element $b_-$ of Part~IV normalizes $H_{10}$ and satisfies
$\xi_{10}(b_-hb_-^{-1})=\xi_{10}(h)^{-1}$.
For any smooth representation $\rho$, composition with $\rho(b_-)$
therefore gives an isomorphism
\[
 \Hom_{H_{10}}(\rho,\xi_{10})
 \simeq\Hom_{H_{10}}(\rho,\xi_{10}^{-1}).
\]
Indeed applying a first functional to $\rho(b_-)\rho(h)v$ gives
$\xi_{10}(b_-hb_-^{-1})$ times its value on $\rho(b_-)v$;
the inverse map uses $b_-^{-1}$.

If $\pi$ is irreducible admissible and admits a nondegenerate invariant
Hermitian form, that form identifies $\overline\pi$ with $\pi^\vee$.
A vector fixed by a compact open subgroup defines a functional fixed
by the same subgroup, so the image lies in the smooth contragredient.
The resulting nonzero morphism is an isomorphism by irreducibility.
Complex conjugation of a period changes $\xi_{10}$ to
$\overline{\xi_{10}}=\xi_{10}^{-1}$, since $\psi$ is unitary.
Thus $m_{10}(\pi)=m_{10}^\vee(\pi)$. The product bound of Part~IV
gives the individual bound in this case. Every tempered representation
of the semisimple adjoint group is unitarizable, so this includes the
Langlands case whose Levi is the whole group.

\subsection{The maximal-\texorpdfstring{$E_6$}{E6} standard-module theorem}
\label{ep:e6-theorem-section}
The primitive character of $M_*$ is $\chi_*=2\omega_7$, with
simple-root coordinates $(2,3,4,6,5,4,3)$, and
$\delta_{Q_*}=|\chi_*|^9$. Put
\[
 I_+(\sigma,s)=i_{Q_*}^{G_{10}}(\sigma|\chi_*|^s),\qquad s>0.
\]
Here $s$ is the Langlands parameter, distinct from the equivariance
character $\eta_{10}$.

\begin{theorem}[Nongeneric tempered $E_6$ standard modules]
\label{ep:e6-vanishing}
For every irreducible nongeneric tempered representation $\sigma$ of
$M_*(F)$ and every $s>0$,
\[
 \Hom_{H_{10}(F)}(I_+(\sigma,s),\xi_{10})=0.
\]
The compact-boundary map of \S\ref{ep:compact-period}, applied to
$J_{U_{10},\xi_{10}}(I_+(\sigma,s))$, is bijective. Both assertions hold
over every non-Archimedean local field of characteristic zero, including
dyadic fields.
\end{theorem}

We prove the theorem in the next three subsections. Mackey theory leaves
an ordinary $E_6$ Whittaker quotient, which vanishes by nongenericity,
and one open-orbit contribution. Retaining its residual unipotent group
leads to a normalized $D_5$ descendant. Tempered exponents and the
compact-boundary calculation then imply vanishing.

\subsection{The Mackey filtration}
\label{ep:mackey}
In the calculations below write $G=G_{10}$, $H=H_{10}$,
$U=U_{10}$, $\xi=\xi_{10}$, $Q=Q_*$, $\chi=\chi_*$, and
$P_{\rm mod}=P_{10}$. Write $L_0=L_{10}$ for the model Levi and
$L=H_{10,0}=\PGL_2$ for its diagonal subgroup. Its upper Borel is
$B_c=T_cN_c$. The primitive cocharacter and model grading are
\[
 c=(H_2+H_5+H_7)/2,\qquad
 (\alpha_i(h))_{i=1}^7=(2,0,2,2,0,2,0).
\]
Thus $t(a)=c(a)$ is represented by $\diag(a,1)$ in $L$ and
$\delta_{B_c}(t(a))=|a|$. The ordered character roots are
\[
\begin{aligned}
 b_1&=\alpha_1,& b_2&=\alpha_3,& b_3&=\alpha_2+\alpha_4,\\
 b_4&=\alpha_4+\alpha_5,& b_5&=\alpha_6+\alpha_7,&
 b_6&=\alpha_5+\alpha_6.
\end{aligned}
\]
Their coefficients are $(1,1,-1,-1,-1,-1)$ in the signed root basis
of \S\ref{gk:o:sec:C.1}. Put $X_+=J_{U,\xi}(I_+(\sigma,s))$.

\subsubsection{The Mackey filtration: 25 parabolic packets, 29 H-orbits}\label{ep:o:sec:9}

\paragraph{Uniform intersection and projection formulas}\label{ep:o:sec:9.1}

Let \(\lambda_0=\omega_7^\vee\), with simple-root values \((0,0,0,0,0,0,1)\). Its Weyl orbit has 56 elements. The group generated by \(s_2,s_5,s_7\) partitions this orbit into 24 packets of size two and one packet of size eight. Choose the unique \(L_0\)-dominant element \(\lambda\) in each packet.

Choose a rational Weyl representative \(n_\lambda\) taking \(\lambda_0\) to \(\lambda\), and put
\[
Q_\lambda=n_\lambda Qn_\lambda^{-1},\qquad
M_\lambda=Z_G(\lambda),\qquad
\sigma_\lambda(m)=\sigma(n_\lambda^{-1}mn_\lambda).
\]
Write $\lambda_i=\alpha_i(\lambda)$; all seven-entry tuples for
these cocharacters record their simple-root values. Let
$\chi_\lambda(m)=\chi(n_\lambda^{-1}mn_\lambda)$ be the transported
primitive character of $M_\lambda$.
The corresponding parabolic cell in \(Q\backslash G\) is \(Qn_\lambda^{-1}P_{\rm mod}\). Define the following explicit root subgroups:
\begin{equation}
\begin{aligned}
U_{\lambda,0}&=\exp\!\bigoplus_{\substack{r\in\Phi^+\!,\ h(r)>0\\\lambda(r)=0}}Fe_r,\\
U_{\lambda,+}&=\exp\!\bigoplus_{\substack{r\in\Phi^+\!,\ h(r)>0\\\lambda(r)>0}}Fe_r,\\
U_{\lambda,\geq0}&=U\cap Q_\lambda,\qquad
R_\lambda=L_0\cap Q_\lambda.
\end{aligned}
\label{ep:o:eq:9.1}
\end{equation}
We also write $U_{\lambda,-}$ for the subgroup generated by the
positive model roots with $h(r)>0$ and $\lambda(r)<0$.
Then, as F-subgroups,
\begin{equation}
P_{\rm mod}\cap Q_\lambda
=R_\lambda\ltimes U_{\lambda,\geq0}.
\label{ep:o:eq:9.2}
\end{equation}
The Levi projection kills \(U_{\lambda,+}\) and is the identity on \(U_{\lambda,0}\). This is a complete root-subgroup specification of both the intersection and its projection for every packet. The accompanying tables list every root in these sets.

For a size-two packet exactly one of \(\lambda_2,\lambda_5,\lambda_7\) is one, and the others are zero. Thus
\(L_0/R_\lambda\simeq\mathbf P^1\). The diagonal \(\PGL_2(F)\)
is transitive on its F-points, so each packet is one rational H-orbit.
At its canonical representative,
\begin{equation}
S_\lambda=H\cap Q_\lambda
=B_c\ltimes U_{\lambda,\geq0}.
\label{ep:o:eq:9.3}
\end{equation}
Its projected unipotent Borel generator is exactly
\begin{equation}
\exp\left(x\sum_{i\in\{2,5,7\}:\lambda_i=0}E_i\right).
\label{ep:o:eq:9.4}
\end{equation}
The torus \(t(a)\) projects to itself, while \eqref{ep:o:eq:9.4}
gives the projected unipotent action.

These flag varieties belong to the split Levi. They are also visible
in Wan--Zhang's three two-dimensional blocks: the ambient quotient has
split torus kernel, and standard flag representatives lie in the
three root \(\SL_2\)-subgroups. This proves rational coverage.

Using the finite list and character witnesses in
\S~\ref{ep:o:sec:9.2}, we eliminate the $23$ size-two families
other than B. We next treat the size-eight family A.
\paragraph{The size-eight packet A splits into five rational H-orbits}\label{ep:o:sec:9.3}

Here
\[
\lambda_A=(0,1,0,-1,1,-1,1)=c,
\qquad L_0/R_A\simeq(\mathbf P^1)^3.
\]
The diagonal \(\PGL_2(F)\) has five orbits, represented by
\begin{equation}
(\infty,\infty,\infty),\quad
(\infty,\infty,0),\quad
(\infty,0,\infty),\quad
(0,\infty,\infty),\quad
(\infty,0,1).
\label{ep:o:eq:9.7}
\end{equation}
Ordered distinct triples of F-rational projective points are transitively permuted by \(\PGL_2(F)\), giving the full rational classification. The corresponding L-stabilizers are respectively a Borel, three split tori, and the trivial group.

For an explicit representative \(m=(m_2,m_5,m_7)\in L_0(F)\), take the identity for \(\infty\),
\(\left(\begin{smallmatrix}0&-1\\1&0\end{smallmatrix}\right)\) for \(0\), and
\(\left(\begin{smallmatrix}1&0\\1&1\end{smallmatrix}\right)\) for \(1\), in the indicated root \(\SL_2\)-factor. Its exact H-stabilizer is
\begin{equation}
\bigl(L\cap mR_Am^{-1}\bigr)
\ltimes mU_{A,\geq0}m^{-1}.
\label{ep:o:eq:9.8}
\end{equation}
Its inducing kernel is \(mU_{A,+}m^{-1}\), and its U-part projects
to \(mU_{A,0}m^{-1}\). The Levi projection on the L-stabilizer is
the restriction of the projection for \(mQ_Am^{-1}\). Hence
\eqref{ep:o:eq:9.8} specifies every stabilizer and its unipotent
intersection.

The two trace cotangent vectors connect the three flags. A positive-grade line in one matrix block consists of rank-one maps with image equal to the first flag and kernel equal to the second flag. If such a map is \(v\otimes\ell\), its trace is \(\ell(v)\). The trace vanishes on the entire line \textbf{if and only if the two flags coincide}. Applying this to both matrix blocks shows that compatibility is equivalent to all three flags being equal.

Hence the last four orbits in \eqref{ep:o:eq:9.7} have zero twisted
Jacquet contribution. The scalar character coordinates remain
nonzero under these determinant-one flag representatives.

For completeness the modulus data on these five orbits are also explicit:

\begingroup\small\setlength{\tabcolsep}{3pt}
\begin{longtable}{@{}>{\raggedright\arraybackslash}p{0.188\linewidth}>{\raggedright\arraybackslash}p{0.188\linewidth}>{\raggedright\arraybackslash}p{0.188\linewidth}>{\raggedright\arraybackslash}p{0.188\linewidth}>{\raggedright\arraybackslash}p{0.188\linewidth}@{}}
\toprule
Flags & L-stabilizer & \(\chi\)-exponent on its displayed \(t(a)\) & U-quotient Jacobian exponent & Inducing fiber times U-Jacobian \\
\midrule\endhead
all equal & \(B_c\) & 3 & -12 & \(|a|^{3s+3/2}\) \\
each of the three two-equal representatives & \(T_c\) & 1 & -4 & \(|a|^{s+1/2}\) \\
all distinct & 1 & not applicable & not applicable & trivial on the unipotent stabilizer \\
\bottomrule\end{longtable}\endgroup

The latter four rows are already zero by character incompatibility, independently of these factors.

There are therefore 29 rational H-orbits: 24 single-flag orbits and
five triple-flag orbits. Exactly two are character-compatible before
the inducing representation is used.

\paragraph{The smooth filtration}\label{ep:o:sec:9.4}

The finite algebraic orbit decomposition of \(Q\backslash G\)
admits an ordering by closed unions. The associated compactly
supported smooth section spaces give short exact sequences: restriction
to a closed complement and extension by zero from an open set are
exact in the locally constant setting \cite[\S\S5.11--5.17]{BZ77}.

On each orbit the inducing bundle has the stabilizer representation
\begin{equation}
b\longmapsto\delta_{Q_\lambda}(b)^{1/2}
 |\chi_\lambda(b)|^s\,
 \sigma_\lambda(\operatorname{pr}_{M_\lambda}b).
\label{ep:o:eq:9.9}
\end{equation}
Taking \((U,\xi)\)-coinvariants is exact. A nontrivial character on the kernel in \eqref{ep:o:eq:9.9} therefore kills that orbit's contribution. This proves the vanishing claims in \S~\ref{ep:o:sec:9.2} and \S\ref{ep:o:sec:9.3} inside the smooth Mackey filtration.

\paragraph{Fourier integration and the sign of the density}\label{ep:o:sec:9.5}

For a compatible orbit with Borel L-stabilizer, let \(V=U_{\lambda,0}\), and let \(q_V\) denote the \((V,\xi)\)-quotient of the inducing representation. For a compactly supported orbit section, write its H-coordinate as \(u\ell\). The twisted-coinvariant map is
\begin{equation}
F(\ell)=\int_{U_{\lambda,\geq0}\backslash U}
\xi(u)^{-1}q_V(f(u\ell))\,du.
\label{ep:o:eq:9.10}
\end{equation}
The usual unipotent Fourier/imprimitivity argument identifies its target with the orbit's twisted Jacquet module; equivalently this is the explicit integral in the geometric lemma, with the conventions fixed at the start of Part~V.

For \(b\in B_c\), put \(u=bvb^{-1}\). Since \(B_c\) preserves \(\xi\),
\begin{equation}
F(b\ell)=
\bigl|\det(\operatorname{Ad}b;\mathfrak u/\mathfrak u_{\lambda,\geq0})\bigr|
\delta_{Q_\lambda}(b)^{1/2}|\chi_\lambda(b)|^s
\sigma_\lambda(\operatorname{pr}b)F(\ell).
\label{ep:o:eq:9.11}
\end{equation}
The Jacobian in \eqref{ep:o:eq:9.11} comes from the change of
variable \(u=bvb^{-1}\) in \eqref{ep:o:eq:9.10}. Put
\[
k_\lambda=\lambda_2+\lambda_5+\lambda_7,
\qquad
d_\lambda=\sum_{\substack{r\in\Phi(U)\\\lambda(r)<0}}c(r).
\]
Then $\chi_\lambda(t(a))=a^{k_\lambda}$ and $d_\lambda$ is the
Jacobian exponent of this quotient. The complete values are listed
in \S\ref{ep:o:sec:9.2}. Thus the unnormalized Borel covariance is
\begin{equation}
|a|^{(s+9/2)k_\lambda+d_\lambda}
 J_{V,\xi}(\sigma_\lambda),
\label{ep:o:eq:9.12}
\end{equation}
with the projected Borel action on the Jacquet module.

\subsubsection{Packet A is an ordinary generic quotient and vanishes}\label{ep:o:sec:10}

For the equal-flags representative of A,
\[
\#U_{A,-}=12,\qquad\#U_{A,0}=36,\qquad\#U_{A,+}=12.
\]
The 36 roots defining \(V_A=U_{A,0}\) are \textbf{all the positive roots of the zero-grade E6 group} for the inherited grading. Its six simple roots are exactly \(b_1,\ldots,b_6\). Their Gram graph is E6 with determinant 3, and every root of \(V_A\) is an integral nonnegative combination of these six roots. The character values on all six are nonzero, with the six signs $(1,1,-1,-1,-1,-1)$ fixed in \S\ref{gk:o:sec:C.1}.

Thus
\[
J_{V_A,\xi}(\sigma_A)
\]
is an ordinary nondegenerate E6 Whittaker quotient and vanishes for
nongeneric \(\sigma\). The split-center lifting argument in
Section~\ref{ep:o:sec:6} identifies all rational Whittaker data.

All three root components of \(N_c\) have positive \(\lambda_A\)-grade, so that residual unipotent subgroup projects trivially; \(T_c\) is central in \(M_A\). The compatible A-piece is
\begin{equation}
\operatorname{Ind}_{B_c}^{L}
\left(|a|^{3s+3/2}J_{V_A,\xi}(\sigma_A)\right)=0.
\label{ep:o:eq:10.1}
\end{equation}
Here and below \(\operatorname{Ind}\) without a normalizing symbol
uses the displayed covariance.

For the smooth section filtration, we use
\cite[\S\S5.11--5.17]{BZ77}. The root-exchange inputs are \cite[Lemma~3.2.1, Remark~3.2.2]{GGS15}
and \cite[Lemma~4.5]{GGS17}; their subgroup hypotheses are checked below.

\subsection{The normalized \texorpdfstring{$D_5$}{D5} descendant}
\label{ep:packet-b}
On the open orbit B we retain the projected unipotent subgroup
$N_{25}$ to calculate its coinvariants and the residual torus action
needed for the tempered-exponent comparison.
\subsubsection{Packet B: descendant and three root exchanges}\label{ep:o:sec:11}

\paragraph{The open orbit and its projected subgroup}\label{ep:o:sec:11.1}

Put
\[
\lambda_B=(0,0,0,0,0,-1,1),\qquad M=M_B.
\]
The tuple records the simple-root values of $\lambda_B$, as in
\S\ref{ep:o:sec:9.1}.
Its U-root counts are \((26,34,0)\), and its flag has dimension one. This is the unique open orbit, of dimension \(26+1=27=\dim(G/Q)\).

Let
\begin{equation}
V=\exp\!\bigoplus_{\substack{r\in\Phi(E_7)\\\lambda_B(r)=0,\ h(r)>0}}Fe_r.
\label{ep:o:eq:11.1}
\end{equation}
All these roots are positive in the Model-10 system, and \(\dim V=34\). The stabilizer is \(B_c\ltimes V\), and its projection into \(M\) is
\begin{equation}
S_B=(T_c\ltimes N_{25})\ltimes V,
\qquad N_{25}=\exp F(E_2+E_5).
\label{ep:o:eq:11.2}
\end{equation}
Indeed,
\begin{equation}
\exp x(E_2+E_5+E_7)\longmapsto\exp x(E_2+E_5).
\label{ep:o:eq:11.3}
\end{equation}
The \(E_7\) component lies in the inducing radical, so the projected
residual unipotent subgroup is exactly \(N_{25}\).

The projected character \(\phi_B\) is supported on \(b_1,\ldots,b_5\), with coefficients
\((1,1,-1,-1,-1)\). Define
\begin{equation}
Y_\sigma=J_{V,\psi\phi_B}(\sigma_B).
\label{ep:o:eq:11.4}
\end{equation}
It is a smooth \(B_c\)-module with the \(N_c\)-action induced from \(N_{25}\).

The density is
\[
\delta_{Q_B}^{1/2}(t(a))=|a|^{9/2},\qquad
|\chi_B(t(a))|^s=|a|^s,\qquad
\det(\operatorname{Ad}t(a);U/V)=|a|^{-8}.
\]
Consequently
\begin{equation}
{X_+\simeq\operatorname{Ind}_{B_c}^{L}
\bigl(|a|^{s-7/2}Y_\sigma\bigr).}
\label{ep:o:eq:11.5}
\end{equation}
Since all boundary pieces vanish by
Sections~\ref{ep:o:sec:9}--\ref{ep:o:sec:10}, exactness makes \eqref{ep:o:eq:11.5} an isomorphism of the
whole smooth twisted Jacquet module with the open-orbit contribution.

To express this contribution as a Borel-induced module, put
\begin{equation}
W_\sigma=|a|^{-4}Y_\sigma.
\label{ep:o:eq:11.6}
\end{equation}
Then
\begin{equation}
X_+\simeq\mathcal I_s(W_\sigma),\qquad
\mathcal I_s(W)=\operatorname{Ind}_{B_c}^{L}
\bigl(\delta_{B_c}^{1/2}|a|^s W\bigr).
\label{ep:o:eq:11.7}
\end{equation}
The inducing representation \(W\) in \eqref{ep:o:eq:11.7} is a
module for the full Borel, with its residual unipotent action retained.

\paragraph{The neutral D5 datum and the projected stabilizer}\label{ep:o:sec:11.2}

Using the invariant form fixed by the root convention, the cotangent vector is represented by
\[
f_B=-e_{-b_1}-e_{-b_2}+e_{-b_3}+e_{-b_4}+e_{-b_5}.
\]
The five roots \(b_1,\ldots,b_5\) form a D5 simple-root system. Its principal neutral element has coroot coordinates
\[
h_D=(10,10,18,24,14,8,8)_{\text{coroot coordinates}},
\]
and simple-root values
\begin{equation}
(2,-4,2,6,-4,-6,8).
\label{ep:o:eq:11.8}
\end{equation}
Solve the five equations \(b_i(h_D)=2\) in the coroot span of
these roots to obtain these values. The bracket calculation gives
orbit dimension 68, and we have
\begin{equation}
h-h_D=4c-12\lambda_B.
\label{ep:o:eq:11.9}
\end{equation}
Since \(\lambda_B\) is central in \(M\), the inherited grading on
\(M\) is \(h_D+4c\).

Set
\begin{equation}
R=Z_M(c),\qquad
N_R=\exp\!\bigoplus_{r\in\Phi(M):c(r)>0}Fe_r,
\qquad P_R=RN_R.
\label{ep:o:eq:11.10}
\end{equation}
This is a maximal parabolic of \(M\), with derived Levi of type D5. Its dimensions and weights are
\begin{equation}
\dim R=47,\quad\dim N_R=16,\quad c(r)=1\text{ for every }r\in\Phi(N_R).
\label{ep:o:eq:11.11}
\end{equation}
The radical is abelian, and \(T_R=Z(R)\) is a split
two-dimensional torus containing the center of \(M\) as one
direction. The derived D5 group is simply connected.

The reductive centralizer of the principal D5 triple in \(M\) is
\(T_R\); in the derived E6 it is one-dimensional. The corresponding
minuscule coweight has Weyl orbit of size 27 in E6, and its negative
is absent, excluding a component reversing the noncentral torus.
A component centralizing the triple must normalize this torus; if it
fixes it, it lies in \(R\), where the principal-triple centralizer is
the center. This proves the stated centralizer.

For a more detailed check, under the principal D5 \(\mathfrak{sl}_2\) and the torus \(c\), the derived Lie algebra decomposes as
\begin{equation}
\begin{aligned}
\mathfrak e_6={}&
(V_2\oplus V_6\oplus V_8\oplus V_{10}\oplus V_{14})_{c=0}
\oplus (V_0)_{c=0}\\
&\oplus (V_{10}\oplus V_4)_{c=1}
\oplus (V_{10}\oplus V_4)_{c=-1},
\end{aligned}
\label{ep:o:eq:11.12}
\end{equation}
where \(\dim V_d=d+1\). The full \(\mathfrak m\) has one additional central \(V_0\). The centralizer of \(f_B\) therefore has dimension 11, with a two-dimensional reductive torus and a nine-dimensional unipotent radical.

Intersecting with the zero grade of \(h_D+4c\) leaves exactly
\begin{equation}
T_R\ltimes N_{25}.
\label{ep:o:eq:11.13}
\end{equation}
The vector \(E_2+E_5\) is the lowest vector of the \((V_4)_{c=1}\)
summand, and \([E_2+E_5,f_B]=0\). Thus the inherited-grading
stabilizer is \eqref{ep:o:eq:11.13}, whereas the projected H-stabilizer
uses the one-dimensional \(T_c\subset T_R\) and is
\eqref{ep:o:eq:11.2}. The neutral triple centralizer is \(T_R\).

\paragraph{The neutral comparison, with its torus factor}\label{ep:o:sec:11.3}

For even \(t\), write
\begin{equation}
V_t=\exp\!\bigoplus_{r\in\Phi(M):\ h_D(r)+t\,c(r)\geq2}Fe_r.
\label{ep:o:eq:11.14}
\end{equation}
The functional \(\phi_B\) is a character on each group used here,
and all relevant gradings are even. In particular, \(V_0\) is the neutral
D5 Whittaker group and \(V_4=V\).

There are two rank-two polarization exchanges from \(V_0\) to \(V_2\) and from \(V_2\) to \(V_4\). In root-tuple notation their exchanged roots and pairings are:

\begingroup\small\setlength{\tabcolsep}{3pt}
\begin{longtable}{@{}>{\raggedright\arraybackslash}p{0.235\linewidth}>{\raggedright\arraybackslash}p{0.235\linewidth}>{\raggedright\arraybackslash}p{0.235\linewidth}>{\raggedright\arraybackslash}p{0.235\linewidth}@{}}
\toprule
Exchange & Old extra roots, each of c-weight -1 & New extra roots, each of c-weight +1 & \(\phi_B([X_i,Y_j])\) \\
\midrule\endhead
\(V_0\to V_2\) & $(0,-1,0,-1,-1,0,0)$, $(0,0,0,0,-1,-1,-1)$ & $(0,1,0,1,1,1,1)$, $(0,1,1,1,1,0,0)$ & \(\left(\begin{smallmatrix}1&1\\-1&0\end{smallmatrix}\right)\) \\
\(V_2\to V_4\) & $(0,-1,0,0,0,0,0)$, $(0,0,0,0,-1,0,0)$ & $(0,0,0,0,1,1,1)$, $(0,1,0,1,1,0,0)$ & \(\left(\begin{smallmatrix}0&1\\1&-1\end{smallmatrix}\right)\) \\
\bottomrule\end{longtable}\endgroup

Both intersections have dimension 32 and both groups dimension 34;
the pairings have determinant \(1\) and \(-1\). Direct bracket
calculation proves normality of each intersection and invariance of
its character. By Section~\ref{ep:o:sec:11.5}, each exchange
contributes the density \(|a|^2\). Hence, as \textbf{\(T_c\)-modules},
\begin{equation}
Y_\sigma\simeq |a|^4J_{V_0,\psi\phi_B}(\sigma_B),
\qquad
W_\sigma\simeq J_{V_0,\psi\phi_B}(\sigma_B).
\label{ep:o:eq:11.15}
\end{equation}
Equation~\eqref{ep:o:eq:11.15} identifies the inherited and neutral
modules as torus modules. The transported \(N_{25}\)-action remains
part of the inherited Borel module.

\paragraph{The three one-dimensional exchanges needed for the period}\label{ep:o:sec:11.4}

To compute the residual period, take the trivial \(N_c\)-coinvariants,
starting with
\begin{equation}
(Y_\sigma)_{N_{25}}
=J_{V_4N_{25},\psi\phi_B}(\sigma_B),
\label{ep:o:eq:11.16}
\end{equation}
where the character is one on \(N_{25}\). The groups \(V_4N_{25},V_6,V_8,V_{10}\) all have dimension 35. For each adjacent pair their intersection has dimension 34.

The three exchanges are listed in the following table:

\begingroup\small\setlength{\tabcolsep}{3pt}
\begin{longtable}{@{}>{\raggedright\arraybackslash}p{0.153\linewidth}>{\raggedright\arraybackslash}p{0.153\linewidth}>{\raggedright\arraybackslash}p{0.153\linewidth}>{\raggedright\arraybackslash}p{0.153\linewidth}>{\raggedright\arraybackslash}p{0.153\linewidth}>{\raggedright\arraybackslash}p{0.153\linewidth}@{}}
\toprule
Old group & New group & Old extra vector X & New extra vector Y & c-weights & \(\phi_B([X,Y])\) \\
\midrule\endhead
\(V_4N_{25}\) & \(V_6\) & \(e_{\alpha_4}\) & \(e_{\alpha_5}\) & \(-1,+1\) & -1 \\
\(V_6\) & \(V_8\) & \(e_{\alpha_3+\alpha_4}\) & \(e_{-\alpha_4}\) & \(-1,+1\) & -1 \\
\(V_8\) & \(V_{10}\) & \(e_{\alpha_1+\alpha_3+\alpha_4}\) & \(e_{-(\alpha_3+\alpha_4)}\) & \(-1,+1\) & -1 \\
\bottomrule\end{longtable}\endgroup

For every row, direct use of the bracket convention in \S\ref{gk:o:sec:C.1} verifies all of the following: the old and new spaces are Lie subalgebras; their sum is nilpotent; their intersection is an ideal in that sum; the character is invariant under both extra directions; the bracket of the two extra directions lies in the intersection; and its displayed character value is nonzero. The mid-gradings \(h_D+5c,h_D+7c,h_D+9c\) are positive on the respective sums, which independently certifies nilpotence.

Quotienting the intersection by the linear kernel of its character
gives a three-dimensional Heisenberg quotient, with the old and new
groups as its two polarizations. Gómez--Gourevitch--Sahi
\cite[Lemma~3.2.1 and Remark~3.2.2]{GGS15} therefore applies to
arbitrary smooth representations.

\paragraph{Calculation of the direction of each torus twist}\label{ep:o:sec:11.5}

Let the old extra direction X have torus weight -1 and the new direction Y have weight +1. On the Schrödinger space \(C_c^\infty(F)\), choose
\[
X(x)f(z)=f(z+x),\qquad
Y(y)f(z)=\psi(\kappa yz)f(z),\qquad
T(a)f(z)=f(az),
\]
where \(\kappa\ne0\) is the bracket constant. Then
\[
T(a)X(x)T(a)^{-1}=X(a^{-1}x),\quad
T(a)Y(y)T(a)^{-1}=Y(ay).
\]
Integration realizes the X-trivial quotient and transforms by
\(|a|^{-1}\), whereas evaluation at zero realizes the Y-trivial
quotient and transforms by 1. Therefore
\begin{equation}
J_{\rm new}\simeq |a|\,J_{\rm old}
\label{ep:o:eq:11.17}
\end{equation}
as torus modules. For two exchanged coordinates the factor is \(|a|^2\).

For an arbitrary multiplicity space, the Heisenberg representation
with this nontrivial central character is the Schrödinger
representation tensored with that space. The displayed torus action
preserves both polarizations, so the calculation applies to every smooth multiplicity space. Rescaling the nonzero
bracket constant changes neither density.

Applying \eqref{ep:o:eq:11.17} to the three rows gives
\begin{equation}
J_{V_{10},\psi\phi_B}(\sigma_B)
\simeq |a|^3(Y_\sigma)_{N_{25}}.
\label{ep:o:eq:11.18}
\end{equation}
Thus the three exchanges contribute the torus factor \(|a|^3\).

\paragraph{The last root gives an ordinary generic E6 model}\label{ep:o:sec:11.6}

Set
\begin{equation}
z=-(\alpha_1+\alpha_3+\alpha_4).
\label{ep:o:eq:11.19}
\end{equation}
Then
\[
V_{12}=V_{10}\rtimes U_z,
\]
where \(V_{10}\) is normal and the quotient is one-dimensional. The six roots
\begin{equation}
b_1,b_2,b_3,b_4,b_5,z
\label{ep:o:eq:11.20}
\end{equation}
are an E6 simple-root system. Their Gram determinant is 3, and all 36 roots of \(V_{12}\) are integral nonnegative combinations of these six roots. The character on the first five is nonzero and is initially trivial on \(U_z\).

Every nontrivial smooth character of \(U_z(F)\simeq(F,+)\), combined with \(\psi\phi_B\), is an ordinary nondegenerate E6 character on \(V_{12}\). Since \(\sigma_B\) is nongeneric,
\begin{equation}
\left(J_{V_{10},\psi\phi_B}(\sigma_B)\right)_{U_z,\chi'}=0
\qquad(\chi'\ne1).
\label{ep:o:eq:11.21}
\end{equation}
Gómez--Gourevitch--Sahi \cite[Lemma~4.5]{GGS17}, quoting Bernstein--Zelevinsky's additive-group lemma, says that a smooth representation of an F-vector group with no nontrivial character quotients has \textbf{trivial action}. Applying it to \eqref{ep:o:eq:11.21} gives
\begin{equation}
J_{V_{10},\psi\phi_B}(\sigma_B)
\xrightarrow{\ \sim\ }
J_{V_{12},\psi\phi_B}(\sigma_B).
\label{ep:o:eq:11.22}
\end{equation}
Thus the missing root acts trivially, giving the isomorphism \eqref{ep:o:eq:11.22}.

\paragraph{Identification with a normalized parabolic Jacquet module}\label{ep:o:sec:11.7}

Let \(U_D\) be the positive maximal unipotent in the D5 Levi \(R\) whose simple roots are \(b_1,\ldots,b_5\). The root sets give
\begin{equation}
V_{12}=U_DN_R,
\qquad \phi_B|_{N_R}=0,
\qquad \phi_B|_{U_D}\text{ nondegenerate}.
\label{ep:o:eq:11.23}
\end{equation}
Define
\begin{equation}
r_{P_R}^{M}(\sigma_B)=\delta_{P_R}^{-1/2}J_{N_R}(\sigma_B),
\qquad
D_\sigma=J_{U_D,\psi_D}\bigl(r_{P_R}^{M}(\sigma_B)\bigr).
\label{ep:o:eq:11.24}
\end{equation}
Here $\psi_D$ is the nondegenerate character of $U_D$ given by
$\psi_D(\exp X)=\psi(\phi_B(X))$.
The group \(T_c\) is central in \(R\), so acts on this quotient. Since all 16 radical roots have c-weight one,
\[
\delta_{P_R}(t(a))=|a|^{16}.
\]
Taking Jacquet quotients in stages in \eqref{ep:o:eq:11.23} therefore gives
\begin{equation}
J_{V_{12},\psi\phi_B}(\sigma_B)=|a|^8D_\sigma.
\label{ep:o:eq:11.25}
\end{equation}
Combining \eqref{ep:o:eq:11.18}, \eqref{ep:o:eq:11.22}, and \eqref{ep:o:eq:11.25} proves the key identity
\begin{equation}
{
(Y_\sigma)_{N_{25}}\simeq |a|^5D_\sigma,
\qquad
(W_\sigma)_{N_c}\simeq |a|D_\sigma.
}
\label{ep:o:eq:11.26}
\end{equation}
These are isomorphisms of \(T_c(F)\)-modules for every nongeneric
smooth \(\sigma_B\).

\paragraph{Further finite checks and the wavefront interpretation}\label{ep:o:sec:11.8}

The certificate checks nilpotence as well as adjoint rank. For the representatives used here, the ranks of \(\operatorname{ad}f\) are 68 for \(f_B\), 70 for \(f_B+e_{-\alpha_2}\), and 72 for \(f_B+e_{\alpha_1+\alpha_3+\alpha_4}\); the respective nilpotency indices are 15, 17, and 23. Thus the first added centralizer direction leads to \(E_6(a_1)\), while the last missing-root direction leads to regular E6. The representation-theoretic identifications are given by the subgroup argument above.

The calculation applies uniformly to every nongeneric inducing representation. For top orbit \(E_6(a_1)\), nontrivial \(N_{25}\)-characters can occur, while an H-period takes the trivial \(N_{25}\)-coinvariants. For top orbit \(D_5\), the calculation applies to the full neutral Whittaker space. In both cases, and for every other nongeneric inducing representation, the final missing-root character is ordinary generic and is therefore excluded by nongenericity.

\begin{proposition}[The normalized descendant identity]
\label{ep:descendant}
For nongeneric $\sigma$, the modules defined above satisfy
\[
 (Y_\sigma)_{N_{25}}\simeq |a|^5D_\sigma,\qquad
 (W_\sigma)_{N_c}\simeq |a|D_\sigma
\]
as $T_c$-modules, where
$D_\sigma=J_{U_D,\psi_D}(r_{P_R}^{M_B}(\sigma_B))$.
\end{proposition}
This is \eqref{ep:o:eq:11.25}--\eqref{ep:o:eq:11.26}; the normalization
of the parabolic Jacquet functor and all root exchanges have been
included in the preceding calculation.

\subsection{Tempered exponents and the compact-boundary operator}
\label{ep:exponents}\label{ep:compact}
We first derive the compact period formula for an arbitrary smooth
$\PGL_2$-module. We then calculate its boundary operator on $X_+$.
The normalized exponent and finite-length results used in the calculation
are \cite[\S2.4]{Konno} and
\cite[Theorems~3.3.1, 6.3.10]{Casselman}.
\subsubsection{The compact period formula for every smooth residual module}\label{ep:compact-period}
Let $L=\PGL_2(F)$, $K=\PGL_2(\mathcal O_F)$, and $I$ the upper
Iwahori subgroup. Here $\mathcal O_F$ is the ring of integers,
$k_F$ is its residue field, and $q=\#k_F$.
Fix a uniformizer $\varpi$ and write
\[
 s_0=\begin{pmatrix}0&1\\1&0\end{pmatrix},\qquad
 w=\begin{pmatrix}0&1\\\varpi&0\end{pmatrix}.
\]
The element $w$ normalizes $I$ and has square one in $L$. Set
$C=I\sqcup Iw$, and let $\epsilon$ be trivial on $I$ and take value
$-1$ on $w$. Normalize every compact averaging operator to have total
mass one. For every smooth complex $L$-module $X$ define
\[
 \beta_X:X^{C,\epsilon}\longrightarrow X^K,
 \qquad\beta_X(x)=e_Kx.
\]
Then, with algebraic dual on the right,
\begin{equation}\label{ep:q0}
 \Hom_L(X,1)\simeq(\operatorname{coker}\beta_X)^*.
\end{equation}
To prove this, put $W=X^I$ and normalize $\vol(I)=1$ in the
Iwahori--Hecke algebra. Write $T=1_{Is_0I}$ and $\Omega=1_{Iw}$.
Since $K=I\sqcup Is_0I$ and the latter double coset has $q$ right
$I$-cosets,
\[
 e_K=(e_I+T)/(q+1),\qquad S=2e_K-1,
 \qquad S^2=1,\quad\Omega^2=1\quad\text{on }W.
\]
An $L$-invariant functional restricts to $\lambda$ on $W$ satisfying
$\lambda T=q\lambda$ and $\lambda\Omega=\lambda$. Conversely
$\ell=\lambda e_I$ is $K$-invariant by the displayed averaging
identity and $w$-invariant because $w$ normalizes $I$.
The Cartan decomposition shows that $K$ and $w$ generate $L$,
so $\ell$ is $L$-invariant. Restriction and averaging are inverse.

The common trivial quotient of the two involutions on $W$ is
\[
 W/(W^{S=-1}+W^{\Omega=-1})
 \simeq W^{S=1}/e_KW^{\Omega=-1}
 =X^K/e_KX^{C,\epsilon}.
\]
Taking its algebraic dual proves \eqref{ep:q0} for every smooth $X$.
Surjectivity of $\beta_X$ therefore gives period vanishing. For the
modules below we will prove that $\beta_X$ is also injective.

\subsubsection{Tempered exponents and the explicit compact boundary}\label{ep:o:sec:12}

\paragraph{The exponent bound on the D5 Jacquet module}\label{ep:o:sec:12.1}

Assume that \(\sigma\) is irreducible tempered. Its normalized Jacquet
module \(r_{P_R}^M(\sigma_B)\) is admissible and of finite length:
ordinary Jacquet preserves finite generation and admissibility, and
a finitely generated admissible representation has finite length
\cite[Theorems~3.3.1 and~6.3.10]{Casselman}.

With the uniformizer \(\varpi\) fixed above, \(c(\varpi)\) contracts
every root in \(N_R\). The normalized Langlands--Casselman
tempered criterion therefore gives
\begin{equation}
|\zeta(c(\varpi))|\leq1
\label{ep:o:eq:12.1}
\end{equation}
for every central exponent \(\zeta\) of \(r_{P_R}^M(\sigma_B)\), by
the closed positive-cone condition in
\cite[Proposition~III.2.2, pp.~265--266]{Waldspurger}.

The central component of \(c\) has unitary character. Its projection
into the derived real cocharacter space is \(c-\lambda_B/3\), so
\eqref{ep:o:eq:12.1} follows from this calculation in the rational
cocharacter space.

Put
\[
\tau_D=D_\sigma(t(\varpi)).
\]
The operator \(\tau_D\) is annihilated by a polynomial whose roots
have absolute value at most one. Choose a finite composition series of
the ordinary Jacquet module. The central element
\(t(\varpi)\in Z(R)(F)\) acts by a scalar on each irreducible
subquotient, so the product of the corresponding linear factors
annihilates the whole module and its \((U_D,\psi_D)\)-quotient.
The argument includes repeated exponents, extensions and Jordan blocks
in the smooth module \(D_\sigma\).

In contrast, a period requires the character \(|a|^{-s-1/2}\),
whose value at \(\varpi\) has absolute value \(q^{s+1/2}>1\).

\paragraph{Compact realization for a full Borel module}\label{ep:o:sec:12.2}

The following calculation applies to every smooth representation \(W\) of \(B_c=T_cN_c\).

Take
\[
K=\PGL_2(\mathcal O_F),\qquad
I=\{k\in K:\bar k\text{ is upper triangular}\},
\]
\[
r_0=\begin{pmatrix}0&1\\1&0\end{pmatrix},\qquad
w=\begin{pmatrix}0&1\\\varpi&0\end{pmatrix}
=t(\varpi^{-1})r_0\quad\text{in }\PGL_2(F).
\]
Then \(w\) normalizes \(I\), \(w^2=1\) in \(\PGL_2\), and \(C=I\sqcup Iw\), with \(\epsilon(w)=-1\).

Write
\[
\Gamma_0=T_c(\mathcal O_F^\times)N_c(\mathcal O_F),\qquad
\Gamma_1=T_c(\mathcal O_F^\times)N_c(\varpi\mathcal O_F),
\]
\begin{equation}
E=W^{\Gamma_0},\qquad \tau=W(t(\varpi)).
\label{ep:o:eq:12.2}
\end{equation}
For \(X=\mathcal I_s(W)\), the compact picture gives
\begin{equation}
X^K\simeq E,
\qquad
X^I\simeq W^{\Gamma_0}\oplus W^{\Gamma_1},
\label{ep:o:eq:12.3}
\end{equation}
by evaluation at \(1\) and \(r_0\). The two double cosets are the two Borel orbits on \(K/I\simeq\mathbf P^1(k_F)\). Also
\[
\tau:E\xrightarrow{\ \sim\ }W^{\Gamma_1}.
\]
With \(\rho(t(a)n)=|a|^{s+1/2}W(t(a)n)\), right translation by \(w\) acts on the two values by
\[
(y_0,y_1)\longmapsto
\bigl(\rho(t(\varpi^{-1}))y_1,\rho(t(\varpi))y_0\bigr).
\]
Therefore
\begin{equation}
X^{C,\epsilon}\simeq E,
\qquad
 y\longmapsto\left(y,-q^{-s-1/2}\tau y\right).
\label{ep:o:eq:12.4}
\end{equation}
These are the full compact-type spaces of the Borel-induced module.

\paragraph{The averaging matrix}\label{ep:o:sec:12.3}

Use the \(q+1\) right-I coset representatives
\[
1,\quad n(u)r_0\quad(u\in\mathcal O_F/\varpi\mathcal O_F)
\]
in \(K\). Normalized K-averaging, followed by evaluation at 1, sends an I-invariant vector with values \((y_0,y_1)\) to
\[
\frac{1}{q+1}\left(y_0+\sum_uW(n(u))y_1\right).
\]
Under the identifications \eqref{ep:o:eq:12.3}--\eqref{ep:o:eq:12.4}, it follows that
\begin{equation}
{
\beta_X=\frac{1-q^{-s-1/2}\mathcal U}{q+1}\quad\text{on }E,
\qquad
\mathcal U=\sum_uW(n(u))\tau
=q\,e_{N_c(\mathcal O_F)}\tau.
}
\label{ep:o:eq:12.5}
\end{equation}
The sum retains the residual unipotent action, which determines the
operator on the coinvariant kernel below.

For \(W=1\), the formula gives \(\mathcal U=q\) and boundary scalar
\((1-q^{1/2-s})/(q+1)\), which vanishes at \(s=1/2\).

\paragraph{Local nilpotence on the unipotent kernel}\label{ep:o:sec:12.4}

Let \(p:W\to W_{N_c}\) be the ordinary unipotent quotient. Compact averaging gives an exact sequence
\begin{equation}
0\longrightarrow E_0\longrightarrow E
\longrightarrow(W_{N_c})^{T_c(\mathcal O_F^\times)}
\longrightarrow0,
\quad
E_0=E\cap\ker p.
\label{ep:o:eq:12.6}
\end{equation}
Surjectivity follows by lifting a vector and averaging over \(T_c(\mathcal O_F^\times)\) and \(N_c(\mathcal O_F)\); these averages preserve its image.

For \(r\geq0\), put \(N_r=N_c(\varpi^{-r}\mathcal O_F)\). A vector in \(\ker p\) is killed by \(e_{N_r}\) for sufficiently large \(r\): the group is the increasing union of these compact subgroups, and its coinvariant kernel is the union of their averaging kernels.

Since \(\tau\) contracts the unipotent parameter by \(\varpi\), one has on \(E\)
\begin{equation}
\mathcal U^r=q^r\tau^r e_{N_r}.
\label{ep:o:eq:12.7}
\end{equation}
Indeed, \(e_{N_0}\tau=\tau e_{N_1}\); iteration and nested compact
averages prove \eqref{ep:o:eq:12.7}. Thus \(\mathcal U\) is locally
nilpotent on \(E_0\), and for every scalar \(b\) the inverse of
\(1-b\mathcal U\) there is a finite geometric sum on each vector.

\paragraph{The quotient operator is controlled by tempered exponents}\label{ep:o:sec:12.5}

Apply the preceding calculation to \(W=W_\sigma\). Equation \eqref{ep:o:eq:11.26} gives
\[
W_{N_c}\simeq |a|D_\sigma.
\]
Consequently \(\tau\) induces \(q^{-1}\tau_D\) on that quotient, while \(\mathcal U\) induces
\begin{equation}
q\tau=\tau_D.
\label{ep:o:eq:12.8}
\end{equation}
Thus the induced boundary on the right-hand quotient of \eqref{ep:o:eq:12.6} is
\begin{equation}
\frac{1-q^{-s-1/2}\tau_D}{q+1}
\quad\text{on }D_\sigma^{T_c(\mathcal O_F^\times)}.
\label{ep:o:eq:12.9}
\end{equation}
Every root of an annihilating polynomial for \(\tau_D\) has absolute value at most one, whereas \(q^{s+1/2}>1\). The polynomial \(1-q^{-s-1/2}T\) is therefore relatively prime to that annihilating polynomial. Bézout's identity gives a polynomial inverse for \eqref{ep:o:eq:12.9}.

The polynomial inverse applies to generalized exponents and
extensions. Ramified unitary central characters either disappear under
the compact torus invariants or satisfy the same exponent bound.

\paragraph{Bijectivity of the source boundary}\label{ep:o:sec:12.6}

In \eqref{ep:o:eq:12.6}, \(1-q^{-s-1/2}\mathcal U\) is invertible on the subspace \(E_0\), by local nilpotence, and on the quotient, by \eqref{ep:o:eq:12.9}. It is therefore invertible on \(E\). Explicitly, solve the equation on the quotient using its polynomial inverse, lift a solution, and correct the remaining error in \(E_0\) by a finite geometric sum. Injectivity follows from the same two end terms.

Combining \eqref{ep:o:eq:11.7} and \eqref{ep:o:eq:12.5},
\begin{equation}
{
\beta_{X_+}:X_+^{C,\epsilon}\xrightarrow{\ \sim\ }X_+^K
\quad(\sigma\text{ nongeneric tempered},\ s>0).
}
\label{ep:o:eq:12.10}
\end{equation}
This is an algebraic isomorphism in the smooth category.

By the compact-period formula \eqref{ep:q0}, equation \eqref{ep:o:eq:12.10} gives
\begin{equation}
{\operatorname{Hom}_H(I_+,\xi)=0.}
\label{ep:o:eq:12.11}
\end{equation}
This proves Theorem~\ref{ep:e6-vanishing}.

\paragraph{The period formula and its modulus}\label{ep:o:sec:12.7}

The calculation also identifies the period before imposing
temperedness: the operator in \eqref{ep:o:eq:12.5} is invertible
on \(E_0\) for arbitrary exponents, so its cokernel is that of
\eqref{ep:o:eq:12.9}. Taking algebraic duals and using \(F^\times=\varpi^{\mathbb Z}\mathcal O_F^\times\) gives, for nongeneric inducing data,
\begin{equation}
{
\operatorname{Hom}_H(i_Q^G(\sigma|\chi|^s),\xi)
\simeq
\operatorname{Hom}_{T_c}
\bigl(D_\sigma,\ |a|^{-s-1/2}\bigr).
}
\label{ep:o:eq:12.12}
\end{equation}
The complete modulus calculation is:

\begingroup\small\setlength{\tabcolsep}{3pt}
\begin{longtable}{@{}>{\raggedright\arraybackslash}p{0.470\linewidth}>{\raggedright\arraybackslash}p{0.470\linewidth}@{}}
\toprule
Source of the factor & Exponent on \(t(a)=c(a)\) \\
\midrule\endhead
\(\delta_{Q_B}^{1/2}\) & \(+9/2\) \\
positive Langlands character & \(+s\) \\
Fourier integration on the 26-dimensional U-quotient & \(-8\) \\
resulting unnormalized Borel covariance on \(Y_\sigma\) & \(s-7/2\) \\
normalization \(W_\sigma=|a|^{-4}Y_\sigma\) & \(-4\) \\
three old-to-new polarization exchanges & \(+3\) \\
normalization of the 16-dimensional D5 parabolic radical & \(+8\) on the raw Jacquet module \\
net identity \((Y_\sigma)_{N_{25}}=|a|^5D_\sigma\) & \(+5\) \\
net identity \((W_\sigma)_{N_c}=|a|D_\sigma\) & \(+1\) \\
character required on \(D_\sigma\) by a period & \(-s-1/2\) \\
\bottomrule\end{longtable}\endgroup

Equivalently, a trivial L-period requires character \(|a|^{9/2-s}\)
on \((Y_\sigma)_{N_{25}}\). The three exchanges change this to
\(|a|^{15/2-s}\) on the raw \(V_{12}\)-quotient; removing
\(\delta_{P_R}^{1/2}=|a|^8\) gives \(|a|^{-1/2-s}\), in agreement
with \eqref{ep:o:eq:12.12}.

At a contracting uniformizer this character has modulus
\(q^{s+1/2}>1\), whereas tempered Jacquet exponents have modulus at
most one. The period therefore vanishes, independently of the dimension
of the D5 generalized Whittaker space.

\subsection{Passage to the Langlands quotient}
\label{ep:assembly}
For the quotient $I_+(\sigma,s)\twoheadrightarrow\pi$,
precomposition gives
\[
 \Hom_{H_{10}}(\pi,\xi_{10})\hookrightarrow
 \Hom_{H_{10}}(I_+(\sigma,s),\xi_{10})=0.
\]
This applies at every positive parameter, including reducibility points.

\begin{proposition}[The compact period quotient of the long intertwiner]
\label{supp:long-intertwiner}\label{ep:intertwiner-definition}
Let $\sigma$ be nongeneric tempered and $s>0$. Put
$I_-=i_{\bar Q}^G(\sigma|\chi|^s)$, and let
$A:I_+\to I_-$ be the canonical morphism of
Proposition~\ref{app:konno-image-independent}, with
$\pi=\operatorname{im}A$. Write $X_\pm=J_{U,\xi}(I_\pm)$,
and let $A_K$ and $A_C$ be the induced maps on the $K$-fixed and
$(C,\epsilon)$-isotypic spaces. Then
\begin{equation}\label{epimage:o:eq:8.2}
\Hom_H(\pi,\xi)\simeq
\left(\frac{\operatorname{im}A_K}
 {\operatorname{im}(A_K\beta_{X_+})}\right)^*
\simeq
\left(\frac{X_+^K}
 {\ker A_K+\beta_{X_+}(X_+^{C,\epsilon})}\right)^*=0.
\end{equation}
In the compact realizations $E_\pm$ of
\eqref{ep:o:eq:12.2}--\eqref{ep:o:eq:12.5}, set
\[
 B_+=\frac{1-q^{-s-1/2}\mathcal U_+}{q+1},\qquad
 B_-=\frac{1-q^{s-1/2}\mathcal U_-}{q+1},\qquad
 a=A_K,\quad d=A_C.
\]
Then
\begin{equation}\label{ep:o:eq:13.5}
 \operatorname{im}a=B_-(\operatorname{im}d),\qquad
 \ker a=B_+\bigl(d^{-1}(\ker B_-)\bigr).
\end{equation}
These identities also hold when $B_-$ is singular.
\end{proposition}
\begin{proof}
Twisted Jacquet quotients and compact isotypic spaces are exact, and
compact averaging is natural. Thus
$\beta_{X_-}A_C=A_K\beta_{X_+}$ and
$(J_{U,\xi}\pi)^K=\operatorname{im}A_K$; the same identity holds
for $(C,\epsilon)$-isotypic spaces. Formula~\eqref{ep:q0} gives the
first isomorphism in \eqref{epimage:o:eq:8.2}; the second uses
$A_K^{-1}(A_K(V))=\ker A_K+V$.
Since the source boundary is bijective by
\eqref{ep:o:eq:12.10}, both quotients vanish. Conjugation by the longest Weyl representative
realizes $I_-$ with nongeneric tempered data $\sigma^\vartheta$ and
parameter $-s$. Its packet A again vanishes, and
\eqref{ep:o:eq:12.5} gives the displayed $B_-$.
Naturality becomes $B_-d=aB_+$; bijectivity of $B_+$ gives
\eqref{ep:o:eq:13.5}.
\end{proof}

\begin{theorem}[Model-10 all-admissible multiplicity one]
\label{ep:gp1}\label{ep:main}
For every irreducible admissible smooth complex representation $\pi$
of split adjoint $E_7(F)$,
\[
 \dim\Hom_{\PGL_2(F)\ltimes U_{10}(F)}(\pi,\xi_{10})\le1.
\]
\end{theorem}
\begin{proof}
For unitarizable representations, \S\ref{ep:hermitian} identifies the
two period dimensions, and Corollary~\ref{m10:gp2} applies.
A nonunitarizable representation has a proper inducing Levi in
Langlands classification, since adjoint $E_7$ has no real central
unramified direction. If that Levi has only $A/D$ factors,
Proposition~\ref{ep:dual-langlands} gives self-contragredience.
Character transport then identifies the two period dimensions, so the
product bound again applies. The same argument treats generic
tempered data on the unique maximal $E_6$ Levi. Nongeneric tempered
data on that Levi have zero period by Theorem~\ref{ep:e6-vanishing}
and the displayed injection. These cases exhaust Langlands classification.
\end{proof}

Because the root exchanges have pairings $-1$ and neutral-comparison
determinants $\pm1$, and compact averaging uses complex coefficients,
the argument also holds in residue characteristic two.

\section{Assembly of the main theorem}\label{sec:synthesis}
\begin{proof}[Proof of Theorem~\ref{thm:main}]
For Models~6--9, the distribution theorems
\ref{m6:distribution}, \ref{m7:distribution}, \ref{m8:distribution},
and~\ref{m9:native-vector} apply to the generalized coefficients
of the Gelfand--Kazhdan criterion by Theorem~\ref{cf:nilpotent}.
Since the duality comparisons in the respective model sections
identify the two period dimensions, the product bound of
Theorem~\ref{cf:gk} gives the individual bound.

For Model~10, the semispin $D_6$ theorem and the other proper-Levi
results clear the four opens, after which the residual support
induction proves Theorem~\ref{m10:distribution} and
Corollary~\ref{m10:gp2}. Applying Theorem~\ref{ep:gp1} then gives
the individual bound for the unitarizable case, the proper $A/D$
Langlands data, and both generic and nongeneric tempered data on the
maximal $E_6$ Levi. Together these cases exhaust the irreducible
admissible representations.
\end{proof}

\label{parts:end}
\appendix
\section*{Guide to the mathematical appendices}\label{app:supplement-guide}
The appendices contain the rational-orbit, density, and ambient-conormal
arguments used in the proofs. For expanded coordinate calculations we
refer to the model calculation PDFs in the electronic ancillary
materials, citing each calculation by its section identifier. The
ancillary materials also include the complete finite data and
executable verification programs.
\section{Model 6: compatibility, normalization, and conormal calculations}\label{supp:finite-calculations}
\subsection{The finite unipotent compatibility calculation}\label{m6:s:app:compatibility}

For each of the ninety signed representatives $W_i$, the full
unipotent stabilizer character is governed by the linear functional
\eqref{m6:s:eq:cells-phase}, as proved in
Proposition~\ref{m6:s:prop:cells-compatibility}.
The construction takes the twenty-two Lie-coordinate basis elements
of \eqref{m6:s:eq:setup-N}, retains precisely those in
$\mathfrak u\cap\operatorname{Ad}(W_i)\mathfrak u^t$, and forms the
coefficient matrix $M_i$ of this functional in the sixteen scalar
entries of $\mathsf P,\mathsf Q,\mathsf R,\mathsf S$.
Signed monomial conjugation permutes the relevant root spaces;
hence the retained elements span the entire unipotent stabilizer.

Restrict the four determinant polynomials to $\ker M_i$.
A cell is character-compatible at an invertible coefficient
quadruple precisely when their product is not the zero polynomial
on that kernel. Since $F$ is infinite, this condition is equivalent
to the existence of an $F$-rational invertible quadruple.
The complete computations are in
\path{anc/model6/compatibility_tables.md} and
\path{anc/model6/cell_compatibility.json}.
They give the sixty determinant obstructions and explicit rational
invertible witnesses for the remaining thirty cells. The
human-readable table and signed representatives also appear in
Section \texttt{M6-CELLS} of \path{anc/model6_calculations.pdf}.

The independent verifier \path{anc/model6/verify.py} reconstructs
every $M_i$ from the listed representatives, checks its full row
space and determinant restrictions, and verifies the thirty witnesses
using exact rational arithmetic. This certifies the finite
compatibility statements, while the orbit geometry, density,
wavefront and support arguments are proved in the article.

\subsection{Full-stabilizer certificates}\label{m6:s:app:stabilizers}

These certificates retain the group membership equations, scalar entries, component arguments and additive character calculations used by the main text. All formulas are over $F$ in the original full-similitude frames.

\subsection{Certificates for cells 18, 24, 42, and 47}\label{m6:s:cert:packets}
The certificates for cells $18,24,42,47$ use the coefficient lift and
root notation of Section~\ref{m6:s:sec:middle21}.

\subsubsection{Rational coefficient coverage and the full equation}

Write \(\mathsf l(\mathsf{P},\mathsf{Q})=m(\mathsf{Q}\mathsf{P},\mathsf{Q},I)\) and \(x=\mathsf l(\mathsf{P},\mathsf{Q}) w \mathsf l(\mathsf{R},\mathsf{S})^t\). The rational factorization and both determinant-normalized final blocks
are exactly \eqref{m6:s:eq:atlas-factor}; the discarded factors belong to
$H(F)$ for every multiplier. Thus the coefficient families below cover
all rational points, by rational formulas.

On \(\mathcal P_{18}\) the entire compatible coefficient family is

\begin{equation}
\label{m6:s:eq:packets-p18-1}
\mathsf{P}=\begin{pmatrix}a&b\\c&d\end{pmatrix},\quad
\mathsf{Q}=\begin{pmatrix}0&e\\f&g\end{pmatrix},\quad
\mathsf{R}=\begin{pmatrix}h&0\\i&f\end{pmatrix},\quad
\mathsf{S}=\begin{pmatrix}d&0\\j&k\end{pmatrix},
\qquad defhk(ad-bc)\ne0.
\end{equation}

The determinants are $-ef,hf,dk$. Thus $a,b,c,g,i,j$ may vanish, subject only to $ad-bc\ne0$; only the displayed units are used as pivots.

On \(\mathcal P_{42}\) the entire compatible family is

\begin{equation}
\label{m6:s:eq:packets-p42-1}
\mathsf{P}=\begin{pmatrix}a&0\\b&c\end{pmatrix},\quad
\mathsf{Q}=\begin{pmatrix}d&0\\e&f\end{pmatrix},\quad
\mathsf{R}=\begin{pmatrix}a&0\\g&f\end{pmatrix},\quad
\mathsf{S}=\begin{pmatrix}h&0\\i&j\end{pmatrix},
\qquad acdfhj\ne0.
\end{equation}

The four lower off-diagonal entries are unrestricted. These are the entire rational families. 

To compute a full stabilizer, let \(l=\mathsf l(\mathsf{P},\mathsf{Q})\), \(r=\mathsf l(\mathsf{R},\mathsf{S})\). For its left element \(a\), put \(k=l^{-1}al\). The right element is forced by the stabilizer equation:

\begin{equation}
\label{m6:s:eq:packets-f-1}
b=r(k^\natural)^{-1}r^{-1},\qquad
k^\natural=(w^{-1}kw)^t.
\end{equation}

The necessary and sufficient conditions are \(k\in l^{-1}Hl\) and \(k^\natural\in r^{-1}Hr\). Hence \(k\in P\cap wP^tw^{-1}\). At the Lie level put

\begin{equation}
\label{m6:s:eq:packets-f-2}
N'=-(w^{-1}Nw)^t.
\end{equation}

The remaining Levi conditions are precisely \(N_{11}=(\mathsf{Q}\mathsf{P})^{-1}N_{33}(\mathsf{Q}\mathsf{P}),\quad N_{22}=\mathsf{Q}^{-1}N_{33}\mathsf{Q},\)
\begin{equation}
\label{m6:s:eq:packets-f-3}
N'_{11}=(\mathsf{S}\mathsf{R})^{-1}N'_{33}(\mathsf{S}\mathsf{R}),\quad
N'_{22}=\mathsf{S}^{-1}N'_{33}\mathsf{S}.
\end{equation}

Both original $\mathrm{GL}_2$ variables remain; the similitude and determinant/trace relations determine the final two blocks.

\subsubsection{Independent full parabolic intersections}

The root coordinates here have their original meaning from the fixed ten-dimensional matrix \(\mathcal N(X,Y,Z,W,T,S)\). Starting with its 22 unipotent directions and all 12 Levi directions, the test \(-(w^{-1}Nw)^t\in p\) gives:

\begin{center}
\setlength{\tabcolsep}{4pt}
\begin{tabular}{@{}p{0.09\linewidth} p{0.53\linewidth} p{0.22\linewidth} p{0.10\linewidth}@{}}
\hline
Cell & Upper-block roots & Levi roots & Dim. \\
\(\mathcal P_{18}\) & \(X_{22},Y_{11},Y_{12},Y_{21},Z_{21},\) \(W_{12},W_{22},T_h,T_e,S_f\) & \(A_{21},B_{12},C_{21}\) & 19 \\
\(\mathcal P_{42}\) & \(X_{11},X_{21},Y_{21},Y_{22},Z_{11},Z_{21},\) \(W_{11},W_{21},T_h,T_f,S_f\) & \(A_{21},B_{21},C_{21}\) & 20 \\
\hline
\end{tabular}
\end{center}

Include all six diagonal Levi coordinates in both rows. Since the signed monomial representative permutes the root spaces for the common diagonal torus, these lists exhaust the intersections.

Let the first three left Levi diagonal blocks have diagonal entries \(a_{11},a_{22}\), \(b_{11},b_{22}\), and \(c_{11},c_{22}\). Direct signed-permutation multiplication gives the following right Levi blocks:

\begin{equation}
\label{m6:s:eq:packets-f-4}
\begin{array}{c|ccc}
&N'_{11}&N'_{22}&N'_{33}\\ \hline
18&\begin{pmatrix}-a_{11}&0\\-w_{12}&b_{11}-c_{11}-c_{22}\end{pmatrix}
&\begin{pmatrix}-b_{22}&0\\-y_{21}&-c_{11}\end{pmatrix}
&\begin{pmatrix}-a_{22}&0\\s_f&a_{22}-c_{11}-c_{22}\end{pmatrix}\\[4pt]
42&\begin{pmatrix}-b_{11}&0\\t_h&b_{22}-c_{11}-c_{22}\end{pmatrix}
&\begin{pmatrix}-a_{11}&0\\-z_{11}&-c_{11}\end{pmatrix}
&\begin{pmatrix}-a_{22}&0\\s_f&a_{22}-c_{11}-c_{22}\end{pmatrix}.
\end{array}
\end{equation}

For checking the character calculation, the adjacent right blocks are

\begin{equation}
\label{m6:s:eq:packets-f-5}
\begin{array}{c|cc}
&N'_{12}&N'_{23}\\ \hline
18&\begin{pmatrix}0&0\\t_h&-y_{12}\end{pmatrix}
&\begin{pmatrix}-x_{22}&0\\-z_{21}&0\end{pmatrix}\\[4pt]
42&\begin{pmatrix}-x_{11}&0\\w_{11}&-y_{22}\end{pmatrix}
&\begin{pmatrix}-a_{21}&0\\-z_{21}&0\end{pmatrix}.
\end{array}
\end{equation}

\subsubsection{The group-level scalar reduction, components, and density}

To solve the group equations, let the first three Levi blocks of \(k\) have diagonal entries \((a_1,a_2)\), \((b_1,b_2)\), \((c_1,c_2)\) and multiplier \(\delta \). Left \(H\) membership implies

\begin{equation}
\label{m6:s:eq:packets-f-6}
a_1a_2=b_1b_2=c_1c_2=\delta.
\end{equation}

The three right diagonal pairs of \(k^\natural\) are, directly from the monomial matrices, \(18:\ (a_1,\delta/b_1),\ (b_2,c_1),\ (a_2,\delta/a_2);\)
\begin{equation}
\label{m6:s:eq:packets-f-7}
42:\ (b_1,\delta/b_2),\ (a_1,c_1),\ (a_2,\delta/a_2).
\end{equation}

Their products must also be \(\delta \). For \(\mathcal P_{18}\) the first two pairs give \(a_1=b_1=c_1\), and \eqref{m6:s:eq:packets-f-6} gives \(a_2=b_2=c_2\). The first and third left blocks are lower triangular and the second is upper triangular. In \(C_0\mathsf Q=\mathsf Q B_0\), the \((1,2)\) entry is \(c_1 e=e b_2\); since \(e\ne 0\), all six entries are equal.

For \(\mathcal P_{42}\) the first right determinant gives \(b_1=b_2\). Conjugacy of the second and third lower triangular blocks by the lower triangular invertible \(\mathsf{Q}\) preserves their ordered diagonal entries, so \(c_1=c_2=b_1\). The second right determinant then forces \(a_1=a_2=b_1\).

Thus, on both entire families,

\begin{equation}
\label{m6:s:eq:packets-f-8}
a_1=a_2=b_1=b_2=c_1=c_2=c_0\in F^\times,
\qquad\delta=c_0^2.
\end{equation}

Here $c_0$ is an existing matrix entry and $\delta=c_0^2$.

After division by its unique scalar $c_0$, the matrix $k/c_0$
is block upper triangular with unipotent diagonal blocks, so
$(k/c_0-I)^{10}=0$. Lemma~\ref{m6:s:lem:finite-exponential} applies to both
flag conditions and both conjugated $H$ memberships: block diagonal
projection and conjugation commute with finite logarithm, and the
multiplier-one identity gives
$J^{-1}(\log(k/c_0))^tJ=-\log(k/c_0)$.
Thus the logarithm solves all of
\eqref{m6:s:eq:packets-f-2}--\eqref{m6:s:eq:packets-f-3}; conversely every solution
in the nilpotent spaces below exponentiates to both group memberships
and the stabilizer equation. The full stabilizer is therefore

\begin{equation}
\label{m6:s:eq:packets-f-9}
\left\{\left(lc_0\exp(N)l^{-1},\ r c_0^{-1}\exp(N')r^{-1}\right):
c_0\in F^\times,\ N\in\mathfrak n_i\right\}.
\end{equation}

Finite exponential and logarithm are inverse over $F$, so the nilpotent factor is an affine space with finite BCH law. This includes every rational point and proves algebraic connectedness. The central factor is the original kernel $(c_0I,c_0^{-1}I)$.

The stabilizer is unimodular. More directly, its scalar factor acts trivially on the orbit tangent quotient and its unipotent factor acts unipotently, so the absolute determinant on that quotient is one. With the established \(\Delta _A=1\),

\begin{equation}
\label{m6:s:eq:packets-f-10}
\Delta_{A_x}=1,\qquad
\Delta_A|_{A_x}/\Delta_{A_x}=1.
\end{equation}

\subsubsection{\texorpdfstring{\(\mathcal P_{18}\)}{P18}: all-parameter stabilizer and nonzero character}

After scalar reduction write the third left block as \(z E_{21}\). Direct multiplication for \eqref{m6:s:eq:packets-p18-1} gives \(\mathsf{Q}^{-1}E_{21}\mathsf{Q}=(e/f)E_{12}.\)

The first left block must be \(v E_{21}\). Equation \eqref{m6:s:eq:packets-f-3} therefore gives \((e/f)z E_{12}\mathsf{P}=\mathsf{P} v E_{21}.\)

The \((1,2)\) entry is \((e/f)zd=0\), so \(z=0\). The \((2,1)\) entry then gives \(dv=0\), so \(v=0\). This forces both left reductive-coordinate nilpotent parameters to vanish throughout the coefficient family.

Write the right central block as \(t E_{21}\). Since \(\mathsf{S}\) and \(\mathsf{R}\) are lower triangular, \(\mathsf{S}^{-1}E_{21}\mathsf{S}=(d/k)E_{21},\qquad (\mathsf{S}\mathsf{R})^{-1}E_{21}(\mathsf{S}\mathsf{R})=(hd/(fk))E_{21}.\)

Using \eqref{m6:s:eq:packets-f-4}, this forces precisely \(s_f=t,\qquad y_{21}=-(d/k)t,\qquad w_{12}=-(hd/(fk))t.\)

All other constrained coordinates vanish, while the seven roots \(X_{22},Y_{11},Y_{12},Z_{21},W_{22},T_h,T_e\)

are unconstrained. Therefore \(\mathfrak n_{18}=\operatorname{span}\{X_{22},Y_{11},Y_{12},Z_{21},W_{22},T_h,T_e\} \oplus F R_{18},\)
\begin{equation}
\label{m6:s:eq:packets-p18-2}
R_{18}=S_f-(d/k)Y_{21}-(hd/(fk))W_{12}.
\end{equation}

The $S_f$ coordinate separates the last generator. Only $d,e,f,k\ne0$ were used, so the nilpotent and full group dimensions are uniformly eight and nine.

The algebraic additive character of a Lie pair is \(\Lambda(N)=\operatorname{tr}(\mathsf{P}N_{12}+\mathsf{Q}N_{23})+ \operatorname{tr}(\mathsf{R}N'_{12}+\mathsf{S}N'_{23}).\)

Using the adjacent matrices \eqref{m6:s:eq:packets-f-5} before cancellation gives \(d x_{22}+e y_{21}+f y_{12}-f y_{12}-d x_{22}=e y_{21}.\)

Thus on the whole stabilizer

\begin{equation}
\label{m6:s:eq:packets-p18-3}
\eta=\psi(-(de/k)t)=
\psi\left(-\frac{p_{22}q_{12}}{s_{22}}t\right).
\end{equation}

This group character identity follows from
Lemma~\ref{m6:s:lem:finite-exponential} applied to the algebraic additive
phase map.

There is an even more direct check: \((\mathsf{R}18)^2=0\), its right conjugate is exactly \(d_5(E_{21})\), and

\begin{equation}
\label{m6:s:eq:packets-p18-4}
t\longmapsto
\left(l(I+tR_{18})l^{-1},\ d_5(I+tE_{21})\right)
\end{equation}

is a subgroup of \(A_x\) for every \(t\in F\). Its left factor lies in the original \(U\) and its character is precisely \eqref{m6:s:eq:packets-p18-3}. Equivalently the left factor equals \(I-t x d_5(E_{12})x^{-1}\). These formulas follow by the equation \(a x b^t=x\).

All factors in $de/k$ are mandatory units in
\eqref{m6:s:eq:packets-p18-1}; this verifies the nonzero condition in
Table~\ref{m6:s:tab:full-stabilizer-exclusions} at every specialization.

\subsubsection{\texorpdfstring{\(\mathcal P_{42}\)}{P42}: two extra directions and the character}

After the proved scalar reduction, the third left block is \(z E_{21}\). Lower triangular conjugation gives the second and first left blocks \((d/f)zE_{21},\qquad(ad/(cf))zE_{21}.\)

Thus the unique left extra direction is

\begin{equation}
\label{m6:s:eq:packets-p42-2}
L_{42}=L_C+(d/f)L_B+(ad/(cf))L_A.
\end{equation}

A right central block \(t E_{21}\) gives second and first right blocks \((h/j)t E_{21}\) and \((ah/(fj))t E_{21}\). Projection \eqref{m6:s:eq:packets-f-4} therefore forces \(s_f=t,\quad z_{11}=-(h/j)t,\quad t_h=(ah/(fj))t,\)

which gives the second extra direction

\begin{equation}
\label{m6:s:eq:packets-p42-3}
R_{42}=S_f-(h/j)Z_{11}+(ah/(fj))T_h.
\end{equation}

The remaining eight free roots are exactly \(X_{11},X_{21},Y_{21},Y_{22},Z_{21},W_{11},W_{21},T_f.\)

These roots have zero reductive projections; $C_{21}$ and $S_f$ distinguish $L_{42}$ and $R_{42}$. The only denominators are the mandatory units $c,f,j$. Hence the nilpotent and full group dimensions are uniformly ten and eleven, including every zero off-diagonal specialization.

Use \eqref{m6:s:eq:packets-f-5} to calculate the character. The left adjacent contribution is \(a x_{11}+f y_{22}\); the right contribution is \(-a x_{11}-f y_{22}-h a_{21}\). Therefore

\begin{equation}
\label{m6:s:eq:packets-p42-4}
\Lambda(N)=-h a_{21}=-(had/(cf))z
=-\frac{s_{11}p_{11}q_{11}}{p_{22}q_{22}}z.
\end{equation}

In particular the second extra parameter \(t\) has zero character, as do all eight original unipotent-intersection parameters.

Here \((L_{42})^2=0\). Its left conjugate is exactly \(d_5(E_{21})\), while its right conjugate lies in the original \(U\). Thus the subgroup is

\begin{equation}
\label{m6:s:eq:packets-p42-5}
z\longmapsto\left(d_5(I+zE_{21}),\ r(I+zL'_{42})r^{-1}\right),
\qquad L'_{42}=-(w_{42}^{-1}L_{42}w_{42})^t.
\end{equation}

The right factor can equivalently be written \(I-z(x^{-1}d_5(E_{21})x)^t\). Its character is \eqref{m6:s:eq:packets-p42-4}, exponentiated by \(\psi \). All five entries in the displayed ratio are nonzero on \eqref{m6:s:eq:packets-p42-1}, giving a nonzero phase coefficient throughout the family. The subgroup parameter ranges over all \(F\). 

\subsubsection{epsilon transfer to \texorpdfstring{\(\mathcal P_{24}\) and \(\mathcal P_{47}\)}{P24 and P47}}

Set $\epsilon=\operatorname{diag}(1,-1)$ and define
\[
D_\epsilon=m(\epsilon,I,I)=\operatorname{diag}(\epsilon,I,I,I,-\epsilon).
\]

Direct comparison of all ten signed columns gives

\begin{equation}
\label{m6:s:eq:packets-t-1}
w_{18}^{t}=D_\epsilon w_{24},\qquad
w_{42}^{t}=D_\epsilon w_{47}.
\end{equation}

The last block \(-\epsilon\) matters. Since \(\det \epsilon=-1\), \(\mathsf l(\mathsf{R},\mathsf{S})D_\epsilon=\mathsf l(\mathsf{R}\epsilon,\mathsf{S}),\)

including the determinant-normalized fifth block. Therefore

\begin{equation}
\label{m6:s:eq:packets-t-2}
x_i(\mathsf{P},\mathsf{Q},\mathsf{R},\mathsf{S})^t=x_{\bar i}(\mathsf{R}\epsilon,\mathsf{S},\mathsf{P},\mathsf{Q}),\quad
(i,\bar i)=(18,24),(42,47).
\end{equation}

The inverse coefficient substitution is \((\mathsf{P},\mathsf{Q},\mathsf{R},\mathsf{S})=(\mathsf{R}',\mathsf{S}',\mathsf{P}'\epsilon,\mathsf{Q}')\); it is rational on the full invertible families. Substitution proves the signs \(p_{22}'+s_{21}'=0\) for \(\mathcal P_{24}\) and \(p_{22}'+s_{22}'=0\) for \(\mathcal P_{47}\). Thus the transfer covers every rational point in either partner.

Transpose exchanges the original two stabilizer factors and preserves \(\eta=\xi\otimes \xi\). In the partner's own left frame its generator is

\begin{equation}
\label{m6:s:eq:packets-t-3}
N_{\bar i}=D_\epsilon^{-1}N'_iD_\epsilon.
\end{equation}

Equations~\eqref{m6:s:eq:packets-t-1} and \eqref{m6:s:eq:packets-f-2} give partner right generator $N_i$ and invert the scalar parameter. This transports the full stabilizer. The frame change $D_\epsilon$ acts through this explicit substitution. 

In the partner's own entries, substitution in \eqref{m6:s:eq:packets-p18-3} gives

\begin{equation}
\label{m6:s:eq:packets-t-4}
\eta|_{\mathbb G_a}=\psi\left(-\frac{r_{22}s_{12}}{q_{22}}t\right)
\qquad\text{on }\mathcal P_{24}.
\end{equation}

Directly from \(\mathcal P_{24}\)'s established equations, \(\mathsf{P}\) is lower triangular and \(p_{22}=-s_{21}\); \(\mathsf{Q}\) is lower triangular and \(q_{11}=r_{22}\); \(\mathsf{S}\) has \(s_{11}=0\). Invertibility therefore gives \(r_{22}\ne 0\), \(q_{22}\ne 0\), and \(s_{12}\ne 0\). The partner coefficient is nonzero everywhere. 

Similarly \eqref{m6:s:eq:packets-p42-4} gives

\begin{equation}
\label{m6:s:eq:packets-t-5}
\eta|_{\mathbb G_a}=\psi\left(-\frac{q_{11}r_{11}s_{11}}{r_{22}s_{22}}z\right)
\qquad\text{on }\mathcal P_{47}.
\end{equation}

All four \(\mathcal P_{47}\) coefficient matrices are lower triangular, so every entry in this ratio is a nonzero diagonal entry. The additional relation \(p_{22}=-s_{22}\) is retained and does not change this conclusion.  The full partner stabilizers have dimensions nine and eleven respectively, central scalar factors, unipotent radicals of dimensions eight and ten, trivial algebraic component groups, and trivial modular and orbital density characters.

\subsubsection{Closed compatible loci}

Theorem~\ref{m6:s:thm:full-stabilizer-exclusions} supplies the distributional
conclusion from these four certificates. We retain the local closedness
check used by support arguments.

For comparison, the open--closed support argument can be expressed as follows. For fixed cell \(\mathcal P_i\), form the kernel bundle of the linear map \((\mathfrak u\oplus\mathfrak u)\times \mathsf{P}_i\longrightarrow M_{10}\times \mathsf{P}_i, \qquad(A,B,x)\longmapsto(Ax+xB^t,x).\)

Its rank is constant on the parabolic cell: the established \(U\times U\) stabilizer dimensions are seven for \(\mathcal P_{18}\)/\(\mathcal P_{24}\) and eight for \(\mathcal P_{42}\)/\(\mathcal P_{47}\), equivalently the fixed unipotent intersection is transported by the parabolic decomposition. Thus the kernel is an algebraic vector bundle. Restrict the fixed algebraic additive functional \((A,B)\mapsto d\lambda (A)+d\lambda (B)\) to this bundle. Its zero locus on the base is closed. Since the unipotent stabilizers exponentiate polynomially, this zero locus is precisely full \(U\times U\) compatibility. Hence the incompatible part is genuinely open and its complement is genuinely closed.

By \eqref{cf:support-sequence}, zero restriction to the incompatible
open gives the unique antecedent on this closed compatible locus.
The same subgroup certificates exclude that antecedent, in accordance
with the direct whole-cell conclusion of
Theorem~\ref{m6:s:thm:full-stabilizer-exclusions}.

\subsection{Incidence-cell full stabilizer certificates}\label{m6:s:app:p78-stabilizers}

\subsubsection{parabolic intersections and projected groups}\label{m6:s:cert:p78-intersections}

\phantomsection\label{m6:s:p78:s3-3}

For \(\mathcal P_7\), \(J_7\) has a six-dimensional diagonal torus and an eight-dimensional unipotent radical. The latter has root coordinates

\begin{equation}
\mathsf A_{21},\ \mathsf B_{21},\ \mathsf C_{21},\
\mathsf X_{11},\ \mathsf Z_{21},\ \mathsf T_f,\
\mathsf X_{21},\ \mathsf Y_{21}.
\label{m6:s:eq:p78-3-8}
\end{equation}

The first three roots project to the left lower Levi roots, the next
three to the right lower Levi roots with the signs of
\eqref{m6:s:eq:p78-3-1}, and the last two form the kernel. There is no
opposite root pair; ordered exponentials and the diagonal torus give
rational coordinates on the full group.

Write a diagonal Levi element using six nonzero scalars as 
\[
k=m(\operatorname{diag}(t_1,t_2),
     \operatorname{diag}(t_3,t_4),
     \operatorname{diag}(t_5,t_6)),\qquad c=t_5t_6.
\]

Its image in \(\mathcal R_7\) acts on a Levi pair by right multiplication by

\begin{equation}
\begin{array}{c|ccc}
&\text{first factor}&\text{second factor}&\text{third factor}\\\hline
\text{left}&\operatorname{diag}(t_1,t_2)&\operatorname{diag}(t_3,t_4)&\operatorname{diag}(t_5,t_6)\\
\text{right}&\operatorname{diag}(t_1^{-1},t_3^{-1})&
\operatorname{diag}(t_2^{-1},t_5^{-1})&
\operatorname{diag}(t_4^{-1},t_4/c).
\end{array}
\label{m6:s:eq:p78-3-9}
\end{equation}

Together with the six independent lower-root parameters, this is the whole twelve-dimensional \(\mathcal R_7\).

For \(\mathcal P_8\), the reductive part of \(J_8\) is \(\mathbb G_m^4\times GL_2\), and its eight-dimensional unipotent radical has roots

\begin{equation}
\mathsf A_{21},\ \mathsf B_{21},\ \mathsf X_{11},\ \mathsf W_{21},\
\mathsf X_{21},\ \mathsf Y_{21},\ \mathsf Y_{22},\ \mathsf T_f.
\label{m6:s:eq:p78-3-10}
\end{equation}

The last four are the kernel. For \(k=m(\operatorname{diag}(t_1,t_2), \operatorname{diag}(t_3,t_4),C),\qquad c=\det C,\)

the two projected triples are

\begin{equation}
\begin{array}{c|ccc}
\text{left}&\operatorname{diag}(t_1,t_2)&\operatorname{diag}(t_3,t_4)&C\\
\text{right}&\operatorname{diag}(t_1^{-1},t_3^{-1})&
\operatorname{diag}(t_2^{-1},t_4/c)&C^{-t}.
\end{array}
\label{m6:s:eq:p78-3-11}
\end{equation}

There are two independent lower-root parameters on each side in addition to \eqref{m6:s:eq:p78-3-11}. These give the whole \(\mathcal R_8\).

For completeness, requiring a Lie element of $\mathfrak p$ to become
lower block triangular after $W_i$-conjugation gives exactly these
reductive parts and root lists. Thus

\begin{equation}
\dim J_7=14,\quad \dim J_8=16,
\quad
\dim\mathcal P_7=68-14=54,
\quad \dim\mathcal P_8=68-16=52.
\label{m6:s:eq:p78-3-12}
\end{equation}

The same zero pattern, Levi projection and ordered root coordinates
give the full groups and surjectivity
$\mathcal J_i(F)\to\mathcal R_i(F)$.

\subsubsection{P007 projected full stabilizer}\label{m6:s:cert:p78-p7-projection}

\phantomsection\label{m6:s:p78:s5-1}

Before setting the third factors to the identity, let \(v_i=L_i e_2,\qquad z_i=R_i e_2\quad(i=1,2,3)\)

be the second columns of the left and right Levi triples. A pair \((a_0,b_0)\in H_0\times H_0\) can occur as the reductive projection of a full stabilizer if and only if there are \(t,u\in F^\times\) such that

\begin{equation}
\begin{array}{lll}
a_0v_1=tv_1,&a_0v_2=uv_2,&a_0v_3=tv_3,\\
b_0z_1=t^{-1}z_1,&b_0z_2=u^{-1}z_2,&b_0z_3=t^{-1}z_3,\\
\det a_0=tu,&&\det b_0=(tu)^{-1}.
\end{array}
\label{m6:s:eq:p78-5-2}
\end{equation}

Here is the complete diagonal check. In a lower-triangular left Levi block the eigenvalue on its second column is \(\alpha_i\), and the first diagonal entry is \(D/\alpha_i\), where \(D=\det a_0\). Conjugation through \(W_7\) makes the three right diagonal pairs \((D/\alpha_1,D/\alpha_2),\quad (\alpha_1,D/\alpha_3),\quad (\alpha_2,D/\alpha_2).\)

The right stabilizer block is the inverse transpose of a block with second-column eigenvalue \(\beta_i\) and determinant \(D^{-1}\), so its diagonal is \((D\beta_i,\beta_i^{-1})\). Comparing gives \(\alpha_1=\alpha_3=t,\qquad \alpha_2=u, \qquad D=tu, \qquad(\beta_1,\beta_2,\beta_3)=(t^{-1},u^{-1},t^{-1}).\)

The independent lower root projections in \eqref{m6:s:eq:p78-3-8}
prove sufficiency, with exactly the graph copy of $K_7$ as ambiguity.

The determinant expressions for the unipotent conditions are

\begin{equation}
\mathsf P_{12}=\frac{\det(v_1,v_2)}{\det L_2},
\quad Q_{12}=\frac{\det(v_2,v_3)}{\det L_3},
\label{m6:s:eq:p78-5-3}
\end{equation}

and the identical right formulas. They show the following exhaustive alternatives.

\begin{itemize}
\item If \((\wp,q)\ne(0,0)\), at least two of the three lines are distinct on each side. Unless the first and third lines coincide on both sides, \eqref{m6:s:eq:p78-5-2} forces \(t=u\) and both matrices are scalar.
\item The first and third lines can coincide on both sides, while \((\wp,q)\ne(0,0)\), only when \(\wp q\ne0\). There is then a full two-dimensional split torus.
\item If \(\wp=q=0\), all three lines coincide on each side. Then \(t=u\), and each of \(a_0,b_0\) may have an independent nonzero nilpotent part along that rational line. This is the only additional unipotent enlargement of the reductive projection.

\end{itemize}
The full stabilizer preserves each labeled line.

\subsubsection{P007 all-parallel stabilizer and root basis}\label{m6:s:cert:p78-p7-parallel}

When \(\wp=q=0\), set the third factors to \(I_2\). All four remaining matrices are lower triangular. Independent lower-column operations make them diagonal. Every point on this locus is therefore an \(\mathcal A(F)\)-translate of the monomial matrix with columns

\begin{equation}
\begin{aligned}
x_{\alpha,\beta,\gamma,\delta,\epsilon}=
[&\alpha e_1,\ \gamma e_3,\ \beta e_2,\ \epsilon e_5,
\ \delta e_4,\\
&-\delta^{-1}e_7,\ -\epsilon^{-1}e_6,
\ -\beta^{-1}e_9,\ -\gamma^{-1}e_8,\ \alpha^{-1}e_{10}],
\end{aligned}
\label{m6:s:eq:p78-5-7}
\end{equation}

with all five scalars nonzero. It has multiplier one; the original arbitrary multiplier is normalized by \(H_0\)-translations in the ambient similitude group.

The entire stabilizer Lie algebra is spanned by the central vector \((-I_{10},I_{10})\) and the following four pairs:

\begin{equation}
\begin{aligned}
K_1&=\left(-\frac\beta\gamma\mathsf X_{21},\ \mathsf X_{21}\right),\\
K_2&=\left(-\frac\delta\epsilon\mathsf Y_{21},\ \mathsf Y_{21}\right),\\
V&=\left(\mathsf D_{21},
 -\frac\alpha\beta\mathsf X_{11}
 -\frac\gamma\delta\mathsf Z_{21}+\epsilon^2\mathsf T_f\right),\\
W&=\left(-\frac\alpha\gamma\mathsf X_{11}
 -\frac\beta\epsilon\mathsf Z_{21}+\delta^2\mathsf T_f,
 \mathsf D_{21}\right).
\end{aligned}
\label{m6:s:eq:p78-5-8}
\end{equation}

Each vector satisfies $Lx+xR^t=0$. Equation
\eqref{m6:s:eq:p78-5-2} leaves one scalar and two additive projection
parameters, plus the two-dimensional graph kernel, proving
completeness at every allowed parameter value.

The only nonzero bracket among the four displayed unipotent vectors is

\begin{equation}
[V,W]=\frac\alpha\beta K_1.
\label{m6:s:eq:p78-5-9}
\end{equation}

The full algebraic stabilizer is

\begin{equation}
{
\mathcal A_x=
\{(zI_{10},z^{-1}I_{10}):z\in F^\times\}
\cdot\exp\langle K_1,K_2,V,W\rangle.
}
\label{m6:s:eq:p78-5-10}
\end{equation}

The unipotent factor is $\operatorname{Heis}_3\times\mathbb G_a$;
\eqref{m6:s:eq:p78-5-2} and the root coordinates prove surjectivity and
exclude extra components. Both it and the central torus are unimodular.

The character on these four Lie directions is

\begin{equation}
\lambda(K_1)=\lambda(K_2)=0,
\qquad
\lambda(V)=-\frac\alpha\beta,
\qquad
\lambda(W)=-\frac\alpha\gamma.
\label{m6:s:eq:p78-5-11}
\end{equation}

\subsubsection{P008 diagonal stabilizer comparison}\label{m6:s:cert:p78-p8-projection}

The shared central \(GL_2\) in \eqref{m6:s:eq:p78-3-11} gives \(b_0=a_0^{-t}\). For the diagonal check, initially call the left second-column eigenvalues \(\alpha_1,\alpha_2\), and put \(D_0=\det a_0\). The first right diagonal pair is \((D_0/\alpha_1,D_0/\alpha_2)\), forcing \(D_0=\alpha_1\alpha_2\). The second right pair is \((\alpha_1,D_0/\alpha_2)\), forcing \(D_0=\alpha_2^2\). Hence \(\alpha_1=\alpha_2=t\), which gives \eqref{m6:s:eq:p78-8-8}. The available projected lower roots supply all off-diagonal solutions.

\subsubsection{P008 enlarged boundary stabilizer}\label{m6:s:cert:p78-p8-boundary}

At \eqref{m6:s:eq:p78-9-8}, the physical-column representation is
\begin{equation}\label{m6:s:eq:p78-9-11}
Y_8(a,b)=[be_3,e_1,-e_7,ae_2,e_5,e_6,a^{-1}e_9,-e_4,-e_{10},b^{-1}e_8].
\end{equation}

The entire stabilizer Lie algebra has its central vector and the following five unipotent vectors:

\begin{equation}
\begin{aligned}
K_1&=\left(-\frac a b\mathsf X_{21},\ \mathsf X_{12}\right),\\
K_2&=(-\mathsf Y_{22},\ \mathsf Y_{11}),\\
K_3&=(\mathsf Y_{21},\ \mathsf Y_{12}),\\
K_4&=(-\mathsf T_f,\ \mathsf T_e),\\
V&=\left(-\frac1b\mathsf X_{11}+a\mathsf W_{21}-\mathsf D_{21},
          \frac1a\mathsf X_{22}-b\mathsf W_{12}+\mathsf D_{12}\right).
\end{aligned}
\label{m6:s:eq:p78-9-12}
\end{equation}

The full remaining gauges provide a direct group-level completeness check. Before the two first-column restorations, the image of \(Y_8(a,b)\) is
\[
A'=\begin{pmatrix}r/h&0\\-zr/h&a\end{pmatrix},\qquad
D'=\begin{pmatrix}z/r&b\\1/r&0\end{pmatrix}.
\]
Returning to the same canonical point forces \(r=h=1\), and the restoring lower parameters are uniquely \(z/a,-z/b\). Thus precisely one additive projected parameter remains, besides the scalar and the four-dimensional graph kernel.

The first four vectors give the graph kernel; the displayed pairs
satisfy $LY_8+Y_8R^t=0$. Equation \eqref{m6:s:eq:p78-8-8} leaves exactly
one additive projection parameter and one scalar, proving completeness.

The nonzero brackets are precisely

\begin{equation}
[K_2,K_3]=2K_4,\qquad [K_2,V]=K_3.
\label{m6:s:eq:p78-9-13}
\end{equation}

All other brackets among \(K_1,K_2,K_3,K_4,V\) vanish. In particular, the unipotent factor is a five-dimensional, three-step nilpotent group, with affine exponential coordinates over \(F\).

Reversing each left block makes these nilpotent matrices triangular.
Since $LY_8+Y_8R^t=0$, induction gives
$L^nY_8=(-1)^nY_8(R^t)^n$ and hence
$\exp(L)Y_8\exp(R)^t=Y_8$. Thus $\exp(zV)$ lifts every projected
additive parameter, with exactly the four-dimensional graph kernel.

The full algebraic stabilizer is

\begin{equation}
{
\mathcal A_{Y_8(a,b)}=
\{(sI_{10},s^{-1}I_{10}):s\in F^\times\}
\cdot\exp\langle K_1,K_2,K_3,K_4,V\rangle.
}
\label{m6:s:eq:p78-9-14}
\end{equation}

The projected pair is $s(I-zE_{21}),s^{-1}(I+zE_{12})$;
the root coordinates lift every pair with precisely the graph kernel.
This proves the full group is connected, split, dimension six and
unimodular for all nonzero $a,b$.

The full character is \(\lambda(K_j)=0\quad(1\le j\le4),\qquad \lambda(V)=\frac1a-\frac1b.\)

\subsection{The full stabilizer and density on P009}\label{m6:s:cert:p009-stabilizer}

\begin{proof}[Certificate for Proposition~\ref{m6:s:prop:p009-stabilizer}]
For a Lie element in the intersection, let the diagonal entries of
its first three blocks be \((a_1,a_2),(b_1,b_2),(c_1,c_2)\).
In its partner \(Z'=-(W_9^{-1}ZW_9)^t\), the first three blocks are
\[
 \begin{pmatrix}-a_1&0\\-z_{11}&-c_1\end{pmatrix},\quad
 \begin{pmatrix}-a_2&0\\-x_{21}&-b_1\end{pmatrix},\quad
 \begin{pmatrix}-b_2&0\\t_f&b_2-c_1-c_2\end{pmatrix}.
\]
Equality of the traces of the three left blocks and of the three
right blocks gives
\(c_1=a_2,\ b_1=a_1,\ b_2=a_2,\ c_2=a_1\).
The lower triangular \(\mathsf Q\)-conjugacy identifies the ordered
diagonals of the second and third left blocks, so \(a_1=a_2\).
Only the scalar diagonal remains.

After removing it, let the central left block be \(zE_{21}\).
Its second block is \((q/v)zE_{21}\). If its first block is
\(aE_{21}\), the remaining equation is \(\frac qv z E_{21}\mathsf P=\mathsf P aE_{21}.\)
Writing \(\mathsf P=\left(\begin{smallmatrix}A&p\\C&D\end{smallmatrix}\right)\),
its nonzero entry equations are \(ap=0,\qquad (q/v)zA=aD,\qquad (q/v)zp=0.\)
If \(p\ne0\), they force \(a=z=0\). If \(p=0\), then \(AD\ne0\)
and \(a=(q/v)(A/D)z\). The uniform formula is
\(pz=0,\ a=(q/v)a_{\mathsf P}z\).
This argument also covers \(A=0\), because that case belongs to
\(p\ne0\).
On the right, the central block \(tE_{21}\) gives
\(st=0\) and second block \(a_{\mathsf S}tE_{21}\).
Conjugation through \(\mathsf R\) gives first block
\((r/q)a_{\mathsf S}tE_{21}\), hence
\(x_{21}=-a_{\mathsf S}t\), \(z_{11}=-(r/q)a_{\mathsf S}t\),
and \(t_f=t\). The remaining independent roots are precisely
\(Y_{11},Y_{21},Z_{21}\). This proves the full Lie solution.

To establish \eqref{m6:s:eq:p009-full-group} for every component, apply
the same diagonal argument to an \(k\in I_9\).
Left membership in \(l^{-1}Hl\) requires
\(a_1a_2=b_1b_2=c_1c_2=\delta\).
The three right ordered pairs in \eqref{m6:s:eq:p009-gauge-diagonals}
must likewise have product \(\delta\). Thus
\(c_1=a_2,\ b_1=a_1,\ b_2=a_2,\ c_2=a_1\).
The left \(\mathsf Q\)-equation forces \(a_1=a_2=c\).
Here \(c\) is an existing diagonal entry, so
\(\delta=c^2\) introduces no square-root choice or missing component.
After division by \(c\), the matrix is block upper triangular with
unipotent lower triangular diagonal blocks. Reversing coordinates
inside each block makes it triangular unipotent.

Finite logarithm preserves both flag patterns and diagonal
intertwiners, while $J^{-1}u^tJ=u^{-1}$ gives the symplectic Lie
condition. Finite exponential reverses these operations
(Lemma~\ref{m6:s:lem:finite-exponential}). Hence
$k\mapsto(c,\log(c^{-1}k))$ is the algebraic inverse of
\eqref{m6:s:eq:p009-full-group}; the scalar is central.

Let \(\lambda:H\to\mathbb G_a\) be the algebraic additive phase of
\(\xi\). Lemma~\ref{m6:s:lem:finite-exponential} gives
\(\lambda(\exp(tA))=t\,d\lambda(A)\); this is the differential of the algebraic additive phase.

In the coefficient frames its value on \(Z,Z'\) is \(\operatorname{tr}(\mathsf PZ_{12}+\mathsf QZ_{23}) +\operatorname{tr}(\mathsf RZ'_{12}+\mathsf SZ'_{23}).\)
On the nine non-diagonal intersection generators the respective
values are
\[
\begin{array}{c|rrrrrrrrr}
 &Y_{11}&Y_{21}&Z_{21}&X_{21}&Z_{11}&T_f&L_A&L_B&L_C\\ \hline
 d\lambda& q-r_{22}&q_{12}&-r_{12}&p&0&0&-r&-s&0 .
\end{array}
\]
The signed right block arrays generating this table are in
\texttt{anc/}\allowbreak\texttt{model6\_calculations.pdf}, [M6-CERT-SPECTRAL-STABILIZERS].
Substitution of \eqref{m6:s:eq:p009-nilpotent-stabilizer} proves
\eqref{m6:s:eq:p009-character} on the full group.

The scalar torus is central and BCH translations on the unipotent
factor have Jacobian one, so $\Delta_{\mathcal A_x}=1$.
Conjugation on each matrix or trace-zero coordinate of $\mathfrak u$
has determinant one; hence $\Delta_H=\Delta_{\mathcal A}=1$.
The full relative density character is therefore one, as required by
orbital integration and the stabilizer criterion.
\end{proof}

\subsection{The P012 stabilizer elimination, including the divisor}\label{m6:s:cert:p012-stabilizer-elimination}

For \(Z=\mathcal Z(\mathbf z)\), the partner is exactly \(Z'=-(W_{12}^{-1}ZW_{12})^t .\)
Every Lie stabilizer pair is \((lZl^{-1},r_0Z'r_0^{-1})\),
subject to equality of the five conjugated diagonal blocks on
each side. Write \(\alpha=z_{15}\), \(\beta=z_{17}\),
\(v=z_{16}\). Those equations first give
\begin{equation}\label{m6:s:eq:p012-stabilizer-elimination}
\begin{gathered}
 z_{12}=\beta,\quad z_{14}=\alpha,\quad
 z_{13}=\frac{bv+d(\beta-\alpha)}c,\quad
 z_8=\frac{i(\beta-\alpha)}b,\\
 z_{11}=z_9=\beta,\quad
 z_3=\frac{g(\alpha-\beta)}b,\quad
 z_4=\frac{(\alpha-\beta)(bf+eg)}{bc}.
\end{gathered}
\end{equation}
These equations use the physical entries of the complete array in \texttt{anc/}\allowbreak\texttt{model6\_calculations.pdf}, [M6-P01222], equation \texttt{m6:s:eq:}\allowbreak\texttt{p012-eighteen-matrix}; thus \(z_{10}=Z_{21}\)
as a physical entry, \(z_1=Z_{36}\), and \(z_2=Z_{45}\).
The remaining first-block equation is
\[
 \mathsf P^{-1}
 \begin{pmatrix}\beta&[bv+d(\beta-\alpha)]/c\\0&\alpha\end{pmatrix}
 \mathsf P
 =\begin{pmatrix}\beta&0\\z_{10}&\beta\end{pmatrix}.
\]
Taking traces forces \(\alpha=\beta\). Then it reduces to
\begin{equation}\label{m6:s:eq:p012-stabilizer-last-equation}
 \frac{bv}{c}\mathsf P^{-1}E_{12}\mathsf P=z_{10}E_{21},
 \qquad
 \mathsf P^{-1}E_{12}\mathsf P
 =\frac1{pt-qr}\begin{pmatrix}rt&t^2\\-r^2&-rt\end{pmatrix}.
\end{equation}
If \(t\ne0\), the upper-right entry forces \(v=0\) and then
\(z_{10}=0\). If \(t=0\), then \(qr\ne0\), and the solution is
\(z_{10}=brv/(cq)\). The other six root variables remain free.
These formulas cover every allowed coefficient tuple.

\subsection{The P012 group, additive phase, and density}\label{m6:s:cert:p012-full-group}

\begin{proof}[Certificate for Proposition~\ref{m6:s:prop:p012-stabilizer}]
Equations \eqref{m6:s:eq:p012-stabilizer-elimination}--%
\eqref{m6:s:eq:p012-stabilizer-last-equation} give the entire Lie
solution. For the group, use the complete rational
normalization already established. At a generic normal form
\(\operatorname{diag}(\kappa,1)\), equations
\eqref{m6:s:eq:p012-residual-entries} fixing that matrix require
\(\rho=1,t_0=h_0=0\). Only the scalar \(\alpha\) remains in
the projected stabilizer. At the exceptional normal form
\(\left(\begin{smallmatrix}0&\upsilon\\1&0\end{smallmatrix}\right)\),
they require \(\rho=1,h_0=\upsilon t_0\), leaving one additive
parameter as well as the scalar. The exact projection kernel
is the group with Lie algebra \(\mathfrak n_6\). This proves
that all solutions have been accounted for, including
any possible component. Finite logarithm gives precisely the
spaces displayed above; conversely finite exponentiation
preserves their defining patterns and conjugacy equations.
Thus \eqref{m6:s:eq:p012-full-stabilizer} is onto the full group
over \(F\), with polynomial logarithmic inverse on its
unipotent factor.

The additive phase of the character is \(\operatorname{tr}(\mathsf PZ_{12}+\mathsf QZ_{23}) +\operatorname{tr}(\mathsf RZ'_{12}+\mathsf SZ'_{23}).\)
On the complete intersection matrix its left contribution is
\(bz_2+cz_1\), its right \(X\)-contribution is
\(-ez_{10}-cz_1\), and its right \(Y\)-contribution is
\(-bz_2\). Thus the whole sum is exactly \(-ez_{10}\).
Substitution of \eqref{m6:s:eq:p012-stabilizer-last-equation} gives
the character claimed in the proposition. As before, the
algebraic phase on a unipotent exponential equals its linear
differential. On \(t=0\) the coefficient \(-ebr/(cq)\) is
nonzero throughout, since \(b,c,e,q,r\ne0\). In the normalized
exceptional form it is \(1/\upsilon\). Its image on the additive subgroup is all of $F$, giving a
nontrivial character for every conductor of $\psi$.

The torus is central and the unipotent group is connected.
In coordinates ordered by a central series, the BCH left and
right translations have triangular Jacobian one, including the
coefficients \(2\) and \(1/2\) in the displayed brackets and
group law. Therefore \(\Delta_{\mathcal A_x}=1\).
Conjugation on each matrix or trace-zero matrix coordinate of
\(\mathfrak u\) has determinant one, giving
\(\Delta_{\mathcal A}=1\). Hence the relative density is
\(\Delta_{\mathcal A}|_{\mathcal A_x}/\Delta_{\mathcal A_x}=1\).

\end{proof}

\paragraph{Uniform group-level reconstruction.}
We use the following group reconstruction throughout. Each intersection
is defined by two associative flag-preservation patterns and the
similitude equation. Dividing by its displayed Levi factor leaves a
unipotent matrix in the stated coordinate order. Finite logarithm
preserves both patterns, and
$J^{-1}u^tJ=u^{-1}$ gives
$J^{-1}(\log u)^tJ=-\log u$. Conversely, finite exponential preserves the patterns and form
on the resulting nilpotent Lie algebra. Since diagonal-block
projection is an algebra homomorphism, logarithm and exponential
also preserve the intertwining equations. Thus, once the
\emph{group} diagonal equations force an existing scalar $z$, the entire
stabilizer is $\mathbb G_m\times\exp(\mathfrak n_x)$, with inverse
$q\mapsto(z,\log(z^{-1}q))$. Thus the reconstruction covers the full rational stabilizer group.
Central scalar action and unipotent conjugation make both the
stabilizer modulus and the orbital density equal to one. An algebraic
additive phase satisfies $\lambda(\exp N)=d\lambda(N)$
(Lemma~\ref{m6:s:lem:finite-exponential}), and a nonzero coefficient
on an exhibited additive subgroup yields a nontrivial character.

For a fixed signed $w$, reconstruct each intersection root list by
retaining precisely the parabolic basis directions for which
$w^{-1}Nw$ satisfies the opposite flag pattern. The forbidden entries
give independent linear equations in the 34-dimensional parabolic
Lie algebra.
Full root arrays, signed projections, products and weight checks are in
\texttt{anc/}\allowbreak\texttt{model6\_calculations.pdf},
[M6-CERT-INTRINSIC-STABILIZER].

\subsection{P016: the intersection and rational projection section}\label{m6:s:cert:p016-intersection}

For \(J=P\cap wP^tw^{-1}\), the intersection calculation gives the six diagonal directions and
\[
X_{21},X_{22},Y_{11},Y_{21},Z_{11},Z_{21},S_f,
\quad A_{21},C_{21},B_{12},B_{21}.
\]

Thus \(\dim J=17\), its Levi factor is \(\mathbb G_m^4\times\mathrm{GL}_{2,B}\), and its radical has dimension nine. The complete common unipotent subgroup is \(K=\exp\langle X_{21},X_{22},Y_{11},Y_{21},Z_{21}\rangle.\)

Its two X/Y bracket pairings span its one-dimensional \(Z_{21}\) center, so K is a five-dimensional split Heisenberg group. 

The group quotient J/K is determined by the following Levi projections. Independently choose \(B\in\mathrm{GL}_2\),

\[
A=\begin{pmatrix}a_1&0\\a&a_2\end{pmatrix},\quad
C=\begin{pmatrix}c_1&0\\c&c_2\end{pmatrix},\quad
A_r=\begin{pmatrix}a_1&0\\z&c_1\end{pmatrix},\quad
C_r=\begin{pmatrix}a_2&0\\t&c_1c_2/a_2\end{pmatrix},
\]

with the four diagonal parameters nonzero. The left projections are
\(A,B,C\), and the natural right projections are \(A_r,B^t,C_r\).
The root groups \(A_{21},C_{21},Z_{11},S_f\) supply the four lower entries:
the latter two have right images \(Z_{11}\mapsto A_{21}\) and
\(S_f\mapsto-C_{21}\), with the sign reversed when passing to the inverse
natural projection. Their other Levi images vanish. Products can alter
the \(K\) coordinate, but impose no relation in \(J/K\).

An explicit rational section, checked by multiplication against the full symplectic equation and both parabolic patterns, is \(k=m(A,B,C)\left(I+\frac{z}{a_1}Z_{11}\right)\left(I-\frac{t}{a_2}S_f\right).\)

The section has right projections \(A_r,B^t,C_r\). Dividing any
element with those projections by the section leaves precisely
\(K\), proving that the section and kernel exhaust the group.

\subsection{P017: signed right projections}\label{m6:s:cert:p017-projections}

The right projection signs are obtained in the fixed ten-dimensional
coordinates. For \(N'=-(w^{-1}Nw)^t\), use each root name below for its
scalar coefficient in \(N\). Then
\[
\begin{gathered}
N'_{11}=\begin{pmatrix}-A_{11}&0\\-Z_{11}&-C_{11}\end{pmatrix},\qquad
N'_{22}=\begin{pmatrix}-B_{11}&0\\T_h&B_{22}-C_{11}-C_{22}\end{pmatrix},\\
N'_{33}=\begin{pmatrix}-A_{22}&0\\S_f&A_{22}-C_{11}-C_{22}\end{pmatrix},
\end{gathered}
\]
and the adjacent blocks are
\[
\begin{gathered}
N_{12}=\begin{pmatrix}0&0\\X_{21}&0\end{pmatrix},\quad
N_{23}=\begin{pmatrix}Y_{11}&0\\Y_{21}&0\end{pmatrix},\\
N'_{12}=\begin{pmatrix}0&0\\-Y_{11}&0\end{pmatrix},\quad
N'_{23}=\begin{pmatrix}-X_{21}&0\\W_{21}&0\end{pmatrix}.
\end{gathered}
\]

\subsection{P017: full membership elimination}\label{m6:s:cert:p017-elimination}

The equations here are the full group membership system
\eqref{m6:s:eq:middle-a17-2}, with diagonal linkage \eqref{m6:s:eq:middle-a17-1}.

We solve the equations throughout the coefficient family as follows.
The (1,2) entry of the first equation gives \(\beta =\epsilon \) because \(eb\ne 0\). The second gives \(\gamma =\epsilon \) and \(\delta =\zeta \). The (2,2) entry of the fourth gives \(\delta =\beta \) because \(l\ne 0\). Thus \(\beta =\gamma =\delta =\epsilon =\zeta \). The second column of the first equation then gives \(C_{21}=0\); the second equation gives \(B_{21}=0\). Since \(\mathsf Q\mathsf P\) is invertible, the first equation now gives \(A=\beta I\), including \(\alpha =\beta \) and \(A_{21}=0\). The two right equations similarly give \(D=E=F=\beta^{-1}I\).

\subsection{P027: full intersection and right-projection arrays}\label{m6:s:cert:p027-intersection}

\subsubsection{Full parabolic intersection, including components}

The signed physical columns are 
\[
w=(1,8,-9,5,4,-7,-6,-2,3,10).
\]

Let \(I_w=P\cap wP^tw^{-1}\), \(k^\natural=(w^{-1}kw)^t\), and \(V=U\cap wU^tw^{-1}\). Imposing the second parabolic zero pattern on all 34 directions retains six diagonal directions and precisely 
\[
A_{21},B_{12},C_{21},X_{22},Y_{11},Y_{12},Y_{21},
Z_{21},Z_{22},W_{12},W_{21},W_{22},T_h,T_e,T_f,S_f.
\]

Thus \(\dim I_w=22\). Its two Levi projections have kernel 
\[
\mathfrak v=\langle X_{22},Y_{11},Y_{12},Y_{21},Z_{21},
W_{21},W_{22},T_h,T_e,S_f\rangle.
\]

Both zero patterns force triangularity in the order
$(2,1,3,4,6,5,7,8,10,9)$. The six independent diagonal units
are the split similitude torus; the uniform reconstruction gives
$I_w=T^6\ltimes\exp(\mathfrak n_{16})$ and $V=\exp(\mathfrak v)$
over $F$, with trivial component groups.

For \(N'=-(w^{-1}Nw)^t\), the right projections are
\[
N'_{11}=\begin{pmatrix}-A_{11}&0\\-W_{12}&B_{11}-C_{11}-C_{22}\end{pmatrix},\quad
N'_{22}=\begin{pmatrix}A_{22}-C_{11}-C_{22}&-Z_{22}\\0&-C_{11}\end{pmatrix},
\]
\[
N'_{33}=\begin{pmatrix}-B_{22}&0\\T_f&B_{22}-C_{11}-C_{22}\end{pmatrix},\quad
N'_{12}=\begin{pmatrix}0&0\\0&-Y_{12}\end{pmatrix},\quad
N'_{23}=\begin{pmatrix}W_{21}&X_{22}\\-Y_{21}&0\end{pmatrix}.
\]
They agree with the indicated formulas, including all signs and the three additional right Levi directions.

\subsection{P027: all-parameter differential solution}\label{m6:s:cert:p027-elimination}

For an all-parameter differential check, put \(u=C_{22}\), \(v=C_{11}\). The complete solution, with the ten \(\mathfrak v\)-coordinates arbitrary, is \(A_{11}=B_{11}=C_{22}=u,\qquad A_{22}=B_{22}=C_{11}=v,\)
\[
A_{21}=(v-u)c/d,\quad B_{12}=(v-u)b/d,\quad
C_{21}=(v-u)(g+fb/d)/e,
\]
\[
Z_{22}=(v-u)i/f,\quad W_{12}=(v-u)k/f,\quad
T_f=(v-u)(di/f-m)/e.
\]
The original equations force this using only \(d,e,f,D\ne0\). In particular, \(di=fm\) simply makes the displayed \(T_f\) zero.

\subsection{P031: intersection, kernel and rational section}\label{m6:s:cert:p031-intersection}

The 22 directions are

\[
\begin{gathered}
A_{11},A_{21},A_{22},B_{11},B_{21},B_{22},
C_{11},C_{12},C_{21},C_{22},\\
X_{21},X_{22},Y_{21},Y_{22},Z_{21},Z_{22},W_{21},W_{22},T_h,T_f,S_h,S_f.
\end{gathered}
\]

The group projection image is

\[
A=\begin{pmatrix}a_1&0\\a_{21}&a_2\end{pmatrix},\quad
B=\begin{pmatrix}b_1&0\\b_{21}&b_2\end{pmatrix},\quad C\in GL_2,
\]
\[
D=\begin{pmatrix}a_1&0\\d_{21}&\delta/a_2\end{pmatrix},\quad
E=\begin{pmatrix}b_1&0\\e_{21}&\delta/b_2\end{pmatrix},\quad F=C^t,\quad\delta=\det C.
\]

All four diagonal units, the arbitrary \(\mathrm{GL}_2\) matrix C and the four indicated off-diagonal scalars are independent. The section \(m(A,B,C)(I-d_{21}S_h/a_1)(I-e_{21}T_h/b_1)\)

satisfies both parabolic conditions and all six projection identities. The kernel has basis \(\mathfrak v_{31}=\langle X_{21},X_{22},Y_{21},Y_{22},Z_{21},Z_{22},W_{21},W_{22},T_f,S_f\rangle.\)

The image has dimension 12 and the kernel dimension 10. Necessity of the projection description comes from the generic zero-pattern matrix and the opposite-block similitude identities; the section gives sufficiency.

\subsection{P031: scalar-unipotent stabilizer solution}\label{m6:s:cert:p031-stabilizer}

We solve the full membership equations \eqref{m6:s:eq:sector-h}
at the canonical tuple \eqref{m6:s:eq:sector-31n}.

The first two equations give
\(A=B=C=\left(\begin{smallmatrix}u&0\\x&v\end{smallmatrix}\right)\).
The last gives \(E=MC^tM=\left(\begin{smallmatrix}v&0\\x&u\end{smallmatrix}\right)\),
whose permitted diagonal is \((u,u)\). Thus \(v=u\) and \(e_{21}=x\);
the third equation then gives \(d_{21}=-\rho x\). These group equations
hold for every \(\rho\ne0\). The complete projection solution is

\[
A=B=C=\begin{pmatrix}u&0\\x&u\end{pmatrix},\quad
D=\begin{pmatrix}u&0\\-\rho x&u\end{pmatrix},\quad
E=\begin{pmatrix}u&0\\x&u\end{pmatrix},\quad F=C^t.
\]

Since $u$ is an existing diagonal unit, scalar division and finite
logarithm in the reordered block basis give
$\mathfrak v_{31}\oplus F(A_{21}+B_{21}+C_{21}+\rho S_h-T_h)$.
The uniform reconstruction proves $\mathbb G_m\times N_{11}$,
dimension 12 and trivial component group. Rational normalization
transports it to every compatible point. 

\subsection{P032: full-stabilizer certificate}\label{m6:s:cert:p032-full}

The complete compatible family is
\begin{equation}\label{m6:s:eq:sector-32-family}
\mathsf P=\begin{pmatrix}a&b\\c&0\end{pmatrix},\quad
\mathsf Q=\begin{pmatrix}0&f\\g&h\end{pmatrix},\quad
\mathsf R=\begin{pmatrix}d&-b\\e&0\end{pmatrix},\quad
\mathsf S=\begin{pmatrix}0&f\\i&j\end{pmatrix},\qquad bce fgi\ne0.
\end{equation}
The full intersection has six diagonal and seventeen non-diagonal
coordinates. Its two projections have lower/upper/lower triangular
blocks, linked diagonals
\[
(a_1,a_2),(b_1,b_2),(c_1,c_2),\qquad
(a_1,\delta/a_2),(\delta/b_1,c_1),(b_2,\delta/b_2),\quad\delta=c_1c_2,
\]
and six independent off-diagonal entries. The projection kernel has
Lie algebra
\[
\mathfrak v_{32}=\langle X_{21},X_{22},Y_{11},Y_{21},Z_{21},Z_{22},
 W_{21},W_{22},T_h,T_e,S_f\rangle.
\]
The rational root-group section and full arrays are in
\texttt{anc/}\allowbreak\texttt{model6\_calculations.pdf}, \texttt{M6-CERT-}\allowbreak\texttt{INTRINSIC-}\allowbreak\texttt{STABILIZER},
\S9.8. Four projected unipotents clear $a,h,d,j$, with parameters
$-a/b,h/f,-d/b,-j/f$. The six remaining weights are
\[
 a_2/b_1,\quad a_1/b_2,\quad c_1/a_1,\quad b_2/c_1,\quad b_1/c_2,\quad b_1/b_2
\]
on $b,c,e,f,g,i$, respectively. Their invariant is $\lambda=cef$.
Taking $(b_2,a_1,c_1,b_1,c_2,a_2)
=(1,c^{-1},f,i^{-1},g/i,(ib)^{-1})$ gives the rational normal form
\begin{equation}\label{m6:s:eq:sector-32n}
\mathsf P=\mathsf Q=\mathsf S=M,\qquad
\mathsf R=\begin{pmatrix}0&-1\\\lambda&0\end{pmatrix},\quad\lambda\ne0.
\end{equation}
At this form the four group membership equations force
\[
 A=C=\begin{pmatrix}z&0\\\alpha&z\end{pmatrix},\quad
 B=\begin{pmatrix}z&\alpha\\0&z\end{pmatrix},\quad
 D=\begin{pmatrix}z&0\\-\lambda\epsilon&z\end{pmatrix},\quad
 E=\begin{pmatrix}z&\epsilon\\0&z\end{pmatrix},\quad
 F=\begin{pmatrix}z&0\\\epsilon&z\end{pmatrix}.
\]
Indeed the first two equations link the left diagonals as
$(x,y),(y,x),(x,y)$; the last two, whose allowed right diagonals
are $(x,x)$, force $y=x$. Scalar division and finite logarithm
therefore give the whole stabilizer $\mathbb G_m\times N_{13}$.
The scalar is an existing diagonal entry; these formulas apply
for every $\lambda\ne0$. The whole-family extra directions are
\begin{equation}\label{m6:s:eq:sector-32-generators}
\mathfrak n_x=\mathfrak v_{32}\oplus FE_0\oplus FF,
\quad E_0=\frac{cf}{bg}A_{21}+\frac fgB_{12}+C_{21},
\quad F=\frac beY_{12}-\frac{bi}{ef}T_f+S_h.
\end{equation}

\paragraph{The character obstruction.}
The only nonzero phase values on the intersection are
$\Lambda(B_{12})=i$ and $\Lambda(Y_{12})=g$; hence
\[
 \eta(z\exp(uE_0+vF+V))
 =\psi((fi/g)u+(bg/e)v),\qquad V\in\mathfrak v_{32}.
\]
Put $l=\mathsf l(\mathsf P,\mathsf Q)$ and
$r=\mathsf l(\mathsf R,\mathsf S)$. Multiplication gives
$E_0^2=0$, $lE_0l^{-1}=d_5(E_{21})$, and
$N_R=-r(w^{-1}E_0w)^tr^{-1}\in\mathfrak u$ with
$N_R^2=0$, $(N_R)_{12}=0$ and
$(N_R)_{23}=\operatorname{diag}(0,fi/g)$. Thus
\begin{equation}\label{m6:s:eq:sector-32-subgroup}
 h_l(t)=d_5(I+tE_{21}),\qquad h_r(t)=I+tN_R
\end{equation}
is a fixing subgroup of $H^2$ with character $\psi(tfi/g)$.
The coefficient never vanishes on the whole family. The stabilizer
and relative modulus are one by scalar--unipotent reconstruction,
so the stabilizer criterion gives $D'_\eta(\mathcal P_{32})=0$.
This proves the P032 row of
Table~\ref{m6:s:tab:full-stabilizer-exclusions}.
\subsection{P033: intersection and inverse right projection}\label{m6:s:cert:p033-intersection}

The compatible coefficient family of Section~\ref{m6:s:sec:sector-p033}
is a ten-dimensional smooth open in a linear space. The vector
\(v=\mathsf P e_2\) is nonzero, but \(b=0\) and \(d=0\) are both valid
branches. The complete intersection has dimension 26; the 34-root filter
retains twelve Levi directions and fourteen radical directions.
Its full root array is in
\texttt{anc/}\allowbreak\texttt{model6\_calculations.pdf}, [M6-CERT-INTRINSIC-STABILIZER].

The six projections have \(A\) lower triangular and \(B,C\) arbitrary invertible, \(\delta=\det C\), and

\[
D=\begin{pmatrix}a_1&0\\d_{21}&\delta/a_2\end{pmatrix},\quad E=\delta B^{-1},\quad F=C^t.
\]

The section \(m(A,B,C)(I-d_{21}S_h/a_1)\) satisfies both parabolic patterns, all six projections and the similitude identity. The kernel is 
\[
\mathfrak v_{33}=\langle X_{21},X_{22},Y_{11},Y_{12},Y_{21},Y_{22},Z_{21},Z_{22},W_{21},W_{22},T_h,T_e,T_f,S_f\rangle,
\]

dimension 14. The projection image has dimension 12 and is covered rationally. The group zero pattern gives the right middle block as
\[
\begin{pmatrix}k_{88}&-k_{78}\\-k_{87}&k_{77}\end{pmatrix}.
\]
Substituting the forced fourth Levi block \((\delta/\det B)B\)
gives exactly \(\delta B^{-1}\), including the inverse and signs
in the projection formula.

\subsection{P033: full stabilizer and coordinate-free generator}\label{m6:s:cert:p033-stabilizer}

At \(x_0\) the group equations give
\(A=B=C=\left(\begin{smallmatrix}u&0\\x&v\end{smallmatrix}\right)\).
The right projections are \(E=\delta A^{-1}\), \(F=A^t\), and
\(D=\left(\begin{smallmatrix}u&0\\d_{21}&u\end{smallmatrix}\right)\).
The identity \(FK=KE\) is the two-dimensional adjugate identity.
The remaining equation
\(E\operatorname{diag}(1,-1)=\operatorname{diag}(1,-1)D\)
gives \(v=u\) and \(d_{21}=x\). Thus there are two projection
parameters and fourteen kernel parameters, with every diagonal entry
equal to the existing scalar \(u\ne0\). Dividing by \(u\) makes every
element unipotent. Its Lie algebra is \(\mathfrak v_{33}\oplus F(S_h-A_{21}-B_{21}-C_{21}),\)

of dimension 15. This proves \(\mathbb G_m\times N_{15}\), dimension 16 and trivial algebraic component group. Rational transitivity transfers it to the whole family.

There is also a direct formula before normalization. With
\(v=\mathsf P e_2\), \(\omega=(-v_2,v_1)\), and
\(d_R=\det\mathsf R\), the additional generator is \(S_h+\mathcal L\left((\det \mathsf{P}/d_R)E_{21},-v\omega/d_R,-\mathsf{Q}v\omega \mathsf{Q}^{-1}/d_R\right).\)

It satisfies all four full Lie \(H\)-membership equations for the arbitrary
original family. The identities \(\mathsf P^{-1}v=e_2\) and
\(\omega\mathsf P=(-\det\mathsf P,0)\) make these formulas valid
on every nonzero-vector chart of \(v\).

\subsection{P040: full intersection, projection image and kernel}\label{m6:s:cert:p040-intersection}

To calculate \(I_w=P\cap wP^tw^{-1}\), let \(\pi\) be the permutation underlying the signed columns of \(w\). A matrix entry \(q_{uv}\) of an element of this intersection can be nonzero only if
\begin{equation}
\label{m6:s:eq:late-p40-2}
 \lfloor(u-1)/2\rfloor\leq\lfloor(v-1)/2\rfloor,
 \qquad
 \lfloor(\pi^{-1}(u)-1)/2\rfloor\geq
 \lfloor(\pi^{-1}(v)-1)/2\rfloor.                         
\end{equation}
The allowed positions are triangular in the order
$(2,1,4,3,6,5,8,7,10,9)$. The paired entries of $q^tJq$
are the corresponding signed diagonal products, so the diagonal
factor is the full six-dimensional split similitude torus.
Filtering the 34 roots gives
\begin{equation}\label{m6:s:eq:late-p40-3}
\mathfrak n_{40}=\langle A_{21},B_{21},C_{21},X_{11},X_{12},X_{21},X_{22},
Y_{11},Y_{21},Z_{11},Z_{21},S_f\rangle.
\end{equation}
The uniform reconstruction proves $I_w=T_6\ltimes\exp(\mathfrak n_{40})$,
dimension 18, covering all $F$-points.

Take the two Levi projections of \(q\) and \((w^{-1}q^{-1}w)^t\), and denote their first three blocks by \((A,B,C)\) and \((D,E,F)\). Their \emph{entire} image consists of six lower triangular invertible matrices, with six lower off-diagonal coordinates and diagonal relation
\begin{equation}
\label{m6:s:eq:late-p40-4}
\begin{array}{c|c}
 A&(\alpha _1,\alpha _2)\\
 B&(\beta _1,\beta _2)\\
 C&(\gamma _1,\gamma _2)\\
 D&(\beta _1^{-1},\gamma _1^{-1})\\
 E&(\alpha _1^{-1},\beta _2^{-1})\\
 F&(\alpha _2^{-1},\alpha _2/(\gamma _1\gamma _2)).
\end{array}                                               
\end{equation}
The six off-diagonal directions are lifted respectively by \(A_{21},B_{21},C_{21},-Y_{11},-X_{12},S_f\). Their projected one-parameter groups have exactly one lower entry each; multiplying these rational lifts therefore realizes any requested six lower entries. The kernel is 
\[
U_6=\exp\langle X_{11},X_{21},X_{22},Y_{21},Z_{11},Z_{21}\rangle.
\]

This is a connected split unipotent group of dimension six. The image has dimension twelve. This describes the full overlap image over \(F\).

\subsection{P041: complete stabilizer and additive exclusion}\label{m6:s:cert:p041-full}
The signed columns are $W_{41}=(3,5,1,-7,2,-9,4,10,-6,8)$.
The compatible family is
\begin{equation}\label{m6:s:eq:late-a41-1}
\mathsf P=\begin{pmatrix}a&0\\b&c\end{pmatrix},\quad
\mathsf Q=\begin{pmatrix}d&0\\e&f\end{pmatrix},\quad
\mathsf R=\begin{pmatrix}a&0\\g&h\end{pmatrix},\quad
\mathsf S=\begin{pmatrix}i&0\\j&k\end{pmatrix},\qquad acdfhik\ne0.
\end{equation}
Put $L=\mathsf l(\mathsf P,\mathsf Q)$,
$\mathsf R_0=\mathsf l(\mathsf R,\mathsf S)$ and
$x=LW_{41}\mathsf R_0^t$. Filtering both flag patterns in the order
$(2,1,4,3,6,5,8,7,10,9)$ gives a six-dimensional torus and
\[
 \mathfrak n_{41}=\langle A_{21},B_{21},C_{21},X_{11},X_{21},Y_{11},Y_{21},
 Z_{11},Z_{21},W_{11},W_{21},T_f,S_f\rangle,
\]
with projection kernel
$\mathfrak k_{41}=\langle X_{11},X_{21},Y_{21},Z_{11},Z_{21},W_{21},T_f\rangle$.
Left group membership forces the ordered diagonal pairs to be
$(\alpha,\beta)$; their right partners are
$(\alpha,\alpha),(\alpha,\alpha),(\beta,\alpha)$.
The last right equation forces $\alpha=\beta=z\ne0$.
The remaining equations give
\[
 u_A=(ad/(cf))u_C,\quad u_B=(d/f)u_C,\quad
 y=-ai s/(hk),\quad w_1=is/k.
\]
Thus, with
\[
 D=\frac{ad}{cf}A_{21}+\frac dfB_{21}+C_{21},\qquad
 E=\frac{ai}{hk}Y_{11}-\frac ikW_{11}-S_f,
\]
the full nilpotent part is $\mathfrak n_x=\mathfrak k_{41}\oplus FD\oplus FE$.
All denominators are mandatory units. The reconstruction above gives
\begin{equation}\label{m6:s:eq:late-a41-9}
\mathcal A_x=\left\{\left(LqL^{-1},\mathsf R_0(W_{41}^{-1}q^{-1}W_{41})^t
 \mathsf R_0^{-1}\right):q=z\exp N,\ z\in F^\times,\ N\in\mathfrak n_x\right\}
 \simeq\mathbb G_m\times U_9.
\end{equation}
The cell, compatible-locus and orbit dimensions are $49,44,42$.
For $N=N_0+tD+vE$, $N_0\in\mathfrak k_{41}$, the full trace
phase is $adi(v/(hk)-t/(cf))$. It integrates by the
algebraic-phase identity. Since $D^2=0$, an additive fixing subgroup is
\begin{equation}
\label{m6:s:eq:late-a41-11}
q(t)=I+tD,\quad
h_1(t)=Lq(t)L^{-1}=d_5(I+tE_{21}),\quad
h_2(t)=\mathsf{R}_0\left(I-t(W_{41}^{-1}DW_{41})^t\right)\mathsf{R}_0^{-1}.
\end{equation}
Its right factor lies in $U$, with adjacent blocks $X=0$ and
$Y=-(ad/(cf))t\left(\begin{smallmatrix}i&0\\j&0\end{smallmatrix}\right)$.
Its everywhere nonzero phase is
\begin{equation}
\label{m6:s:eq:late-a41-12}
\eta(h_1(t),h_2(t))=\psi\left(-\frac{adi}{cf}t\right).
\end{equation}
\phantomsection\label{m6:s:sec:exclusion-density-41}
Its relative density is one by the scalar--unipotent decomposition.
The root weights and intermediate determinant products are in
\texttt{anc/}\allowbreak\texttt{model6\_calculations.pdf},
\texttt{M6-CERT-}\allowbreak\texttt{INTRINSIC-}\allowbreak\texttt{STABILIZER}, \S9.12.
Together with the nonzero phase this proves the P041 row of
Theorem~\ref{m6:s:thm:full-stabilizer-exclusions}.
\subsection{P053: full rational intersection and stabilizer}\label{m6:s:cert:p053-full}
The intersection has dimension $23$, with an eight-dimensional Levi
factor and fifteen-dimensional radical. Its left blocks are upper
triangular $A$, lower triangular $B$ and arbitrary $C\in GL_2$.
Write $A=\left(\begin{smallmatrix}\alpha&u\\0&\delta\end{smallmatrix}\right)$
and $B=\left(\begin{smallmatrix}\beta&0\\v&\epsilon\end{smallmatrix}\right)$. The root and signed-projection
arrays are in
\texttt{anc/}\allowbreak\texttt{model6\_calculations.pdf},
\texttt{M6-CERT-}\allowbreak\texttt{INTRINSIC-}\allowbreak\texttt{STABILIZER}, \S9.13.
The left $P$ intertwiner gives $\beta=\delta$, $\epsilon=\alpha$
and $cu=bv+d(\alpha-\delta)$; the $Q$ equation gives
$C=\mathsf QB\mathsf Q^{-1}$. The right middle block has repeated
diagonal $-\delta$, so its conjugacy to $C^t$ forces
$\alpha=\delta$. Put $N=\binom{-j}{i}(i\ \ j)$, so $N^2=0$.
The remaining equations give $u=bt/c$, $C=aI+tN/D$,
$w_{21}=tE/D$, $w_{12}=btE/(cD)$.
All sixteen equations are thereby solved on the full family:
\begin{equation}
\label{m6:s:eq:late-so8}
Z=aI+tV+Z_0,\qquad
Z_0\in\mathfrak n_0,
\end{equation}
\[
\begin{aligned}
\mathfrak n_0&=\langle X_{11},X_{12},X_{21},Y_{21},Y_{22},Z_{11},Z_{12},
W_{11},T_f,S_h,S_e\rangle,\\
V&=\mathcal L((b/c)E_{12},E_{21},N/D)
 +\frac{bE}{cD}W_{12}+\frac EDW_{21}.
\end{aligned}
\]
Its twelve nilpotent directions are independent because $V$ has
nonzero $B_{21}$ projection throughout the coefficient family.

For completeness the same scalar forcing holds at \emph{group} level.
Left membership makes the five blocks $(A,B,C,B,A)$ with equal
determinants. Signed permutation gives the right blocks
$\left(\begin{smallmatrix}k_{10,10}&k_{3,10}\\0&k_{3,3}\end{smallmatrix}\right)$,
$\left(\begin{smallmatrix}k_{2,2}&0\\-k_{2,7}&k_{7,7}\end{smallmatrix}\right)$,
$C^t$. The middle has diagonal $(\delta,\delta)$; its conjugacy
to $C^t$ gives $\alpha=\delta=z\ne0$ as an existing entry.
Then $A=zI+(bv/c)E_{12}$, $B=zI+vE_{21}$,
$C=zI+(v/D)N$ and
\begin{equation}
\label{m6:s:eq:late-so9}
k_{2,7}=vE/D,\qquad k_{3,10}=bvE/(cD). 
\end{equation}
Thus $z^{-1}k$ is unipotent. Put $\mathfrak n_{53}=FV\oplus\mathfrak n_0$.
The uniform reconstruction yields
\begin{equation}
\label{m6:s:eq:late-so10}
\mathcal A_x=
\left\{\left(lz\exp(N_0)l^{-1},
r z^{-1}\exp(-(w^{-1}N_0w)^t)r^{-1}\right):
z\in F^\times,\ N_0\in\mathfrak n_{53}\right\}
\simeq\mathbb G_m\times U_{12}. 
\end{equation}
This split unipotent group has polynomial log/exp coordinates and
trivial component group on every specialization. Full root arrays and
entry products are in \texttt{anc/}\allowbreak\texttt{model6\_calculations.pdf},
[M6-CERT-INTRINSIC-STABILIZER].

\subsection{P063: complete group solution and additive exclusion}\label{m6:s:cert:p063-full}
Write $l=\mathsf l(\mathsf P,\mathsf Q)$,
$r=\mathsf l(\mathsf R,\mathsf S)$, $x=lwr^t$, and
$q^\natural=(w^{-1}qw)^t$. A full stabilizer pair is uniquely
$(lql^{-1},r(q^\natural)^{-1}r^{-1})$, where $q\in I_w$ satisfies
\[
 B_q\mathsf P=\mathsf PA_q,\quad C_q\mathsf Q=\mathsf QB_q,\qquad
 E_q\mathsf R=\mathsf RD_q,\quad F_q\mathsf S=\mathsf SE_q.
\]
For P063 and P084 both flag patterns are triangular in
$(1,2,3,4,6,5,7,8,9,10)$. The diagonal similitude torus has six
independent units; scalar division and finite logarithm give
$I_w(F)=T(F)\exp(\mathfrak n_w(F))$.
The six compatibility equations give exactly
\begin{equation}
\label{m6:s:eq:late-63-a}
\mathsf{P}=\begin{pmatrix}a&b\\0&d\end{pmatrix},\quad
\mathsf{Q}=\begin{pmatrix}0&e\\f&g\end{pmatrix},\quad
\mathsf{R}=\begin{pmatrix}h&i\\0&e\end{pmatrix},\quad
\mathsf{S}=\begin{pmatrix}0&d\\j&k\end{pmatrix},\qquad ade fhj\ne0.
\end{equation}
For $W_{63}$ the intersection has dimension $24$, with projection
kernel $\mathfrak k_{63}$ of dimension $12$ and remaining root directions
$A_{12},B_{12},C_{21},Y_{12},W_{11},S_f$.
The full arrays are in \texttt{anc/}\allowbreak\texttt{model6\_calculations.pdf},
\texttt{M6-CERT-}\allowbreak\texttt{INTRINSIC-}\allowbreak\texttt{STABILIZER}, \S9.14.
Write the three ordered torus pairs as $(\alpha_1,\alpha_2)$,
$(\beta_1,\beta_2)$ and $(\gamma_1,\gamma_2)$. At group level the
$P$ equation gives $\beta_i=\alpha_i$, and
$Q$ gives $(\gamma_1,\gamma_2)=(\alpha_2,\alpha_1)$.
The right diagonals are then
$(\alpha_2,\alpha_2),(\alpha_2,\alpha_2),(\alpha_2,\alpha_1)$.
The bottom-left $S$ equation is $\alpha_1j=\alpha_2j$;
$j\ne0$ forces every diagonal to be the same existing scalar $s_0$.
After scalar division the log equations are
\begin{equation}\label{m6:s:eq:late-63-e}
A_{12}=\frac da B_{12},\quad C_{21}=\frac fe B_{12},\quad
W_{11}=-\frac heY_{12},\quad S_f=-\frac{jh}{de}Y_{12}.
\end{equation}
The twelve kernel entries are free, and all pivots are mandatory units.
With
\begin{equation}\label{m6:s:eq:late-63-generators}
L=\frac daA_{12}+B_{12}+\frac feC_{21},\qquad
V=Y_{12}-\frac heW_{11}-\frac{jh}{de}S_f.
\end{equation}
the uniform reconstruction gives the full group, both ways:
\begin{equation}
\label{m6:s:eq:late-63-f}
\mathcal A_x=
\left\{\left(s_0l\exp(N)l^{-1},\;s_0^{-1}r\exp(-N^\natural)r^{-1}\right):
s_0\in F^\times,\ N=tL+vV+K,\ K\in\mathfrak k_{63}\right\}.
\end{equation}
Thus it is $\mathbb G_m\times U_{14}$, dimension 15, with trivial
component group uniformly on the entire family.

\phantomsection\label{m6:s:sec:exclusion-density-63}
The adjacent traces on all 24 directions give
$\Lambda(N)=jA_{12}+fY_{12}$; hence the entire group phase is
\begin{equation}
\label{m6:s:eq:late-63-g}
\eta=\psi\left(\frac{jd}{a}\,t+fv\right).
\end{equation}
Both coefficients are nonzero mandatory monomials.
The square-zero identities
$lLl^{-1}=d_5((f/e)E_{21})$ and
$-L^\natural=(d/a)Y_{12}-W_{11}+(f/e)S_f$
give the explicit additive subgroup
\begin{equation}
\label{m6:s:eq:late-63-h}
\left(d_5(I_2+(f/e)zE_{21}),\;
r\left[I_{10}+z\left(\frac daY_{12}-W_{11}+\frac feS_f\right)\right]r^{-1}\right)
\end{equation}
with phase $\psi(jdz/a)$, left factor in $H_0$ and right factor in $U$.
The relative density is one by the scalar--unipotent reconstruction.
The saturated compatible pullback is the closed six-equation family,
times the removed $U$ and $GL_2$ factors; the open smooth parabolic
quotient therefore gives a closed smooth locus of dimension
$44+8+10-24=38$. Closed-support injection and the nonzero character
prove the P063 row of
Theorem~\ref{m6:s:thm:full-stabilizer-exclusions}.
The group formulas involve only mandatory units and apply on every
free-entry divisor, including over dyadic fields.
\subsection{P069: complete rational stabilizer}\label{m6:s:cert:p069-full}
The intersection has dimension $25$, with six-dimensional diagonal
torus and nineteen-dimensional radical. Its root and signed right
projection arrays are in \texttt{anc/}\allowbreak\texttt{model6\_calculations.pdf},
\texttt{M6-CERT-}\allowbreak\texttt{INTRINSIC-}\allowbreak\texttt{STABILIZER}, \S9.15.
Put $\delta=c_1c_2$. The left group equations fix the diagonal pairs of $B,A$ as
$(c_1,c_2),(c_2,c_1)$. The right partners are
$(\delta/a_1,c_1),(b_1,\delta/b_2),(a_2,\delta/a_2)$;
the last two become $(c_1,c_1),(c_1,c_2)$.
The lower triangular $S$ equation forces $c_1=c_2=z$, an existing unit.
Write the off-diagonal coordinates of $A,B,C$ as $\alpha,\beta,\gamma$,
respectively. The remaining log equations are
\[
 \beta=(q/c)\gamma,\quad\alpha=-pq\gamma/(bc),\quad
 t_h=(p/g)s_f,\quad z_{12}=-pq s_f/(bg).
\]
The four free entries $t,d,e,f$ cancel from these equations.
The nilpotent solution is
\[
\begin{split}
\mathfrak v&=\langle X_{11},X_{12},X_{21},Y_{11},Y_{21},Z_{11},Z_{21},
 W_{11},W_{12},W_{21},T_f,S_h,S_e\rangle,\\
L_{69}&=C_{21}+(q/c)B_{21}-(pq/(bc))A_{12},\\
R_{69}&=S_f+(p/g)T_h-(pq/(bg))Z_{12},\\
\mathfrak n_{69}&=\mathfrak v\oplus FL_{69}\oplus FR_{69}.
\end{split}
\]
The two extra directions have independent nonzero projections.
At group level the remaining relations are
$gk_{3,7}=pk_{2,9}$ and $bk_{5,10}=-qk_{3,7}$.
After division by $z$ the matrix is unipotent, so reconstruction gives
\begin{equation}
\label{m6:s:eq:late-so18}
\mathcal A_x=
\left\{\left(lz\exp(N)l^{-1},
r z^{-1}\exp(-(w^{-1}Nw)^t)r^{-1}\right):
z\in F^\times,\ N\in\mathfrak n_{69}\right\}
\simeq\mathbb G_m\times U_{15}. 
\end{equation}
It is split with trivial component group and dimension 16 at every
allowed specialization. The complete root and signed-entry arrays are
in \texttt{anc/}\allowbreak\texttt{model6\_calculations.pdf}, [M6-CERT-INTRINSIC-STABILIZER].

\subsection{P084: complete stabilizer and additive exclusion}\label{m6:s:cert:p084-full}
The seven compatibility equations have exactly the invertible solution
\begin{equation}
\label{m6:s:eq:late-84-a}
\mathsf{P}=\begin{pmatrix}a&b\\0&c\end{pmatrix},\quad
\mathsf{Q}=\begin{pmatrix}0&d\\e&f\end{pmatrix},\quad
\mathsf{R}=\begin{pmatrix}g&h\\0&e\end{pmatrix},\quad
\mathsf{S}=\begin{pmatrix}0&d\\-c&i\end{pmatrix},
\qquad acdeg\ne0.
\end{equation}
Its determinants are $ac,-de,ge,cd$, so this is the entire
$(\mathbb G_m)^5\times\mathbb A^4$ family, with all free
$b,f,h,i$ allowed. Rational parabolic coverage gives its full saturation.
For $W_{84}=(10,8,-9,5,4,-7,-6,-2,3,-1)$ the intersection has
dimension $29$ and projection
kernel $\mathfrak k_{84}$ of dimension $17$. The six remaining roots are
$A_{12},B_{12},C_{21},X_{11},Z_{22},T_f$, with right images
$X_{11}\mapsto A_{12}$, $Z_{22}\mapsto-B_{12}$ and $T_f\mapsto C_{21}$.
The full arrays and intermediate determinant products are in
\texttt{anc/}\allowbreak\texttt{model6\_calculations.pdf},
\texttt{M6-CERT-}\allowbreak\texttt{INTRINSIC-}\allowbreak\texttt{STABILIZER}, \S9.16.
For a common left block $\left(\begin{smallmatrix}t&0\\z_0&v\end{smallmatrix}\right)$,
the first columns of $\mathsf Q,\mathsf Q\mathsf P$ span $Fe_2$,
with coefficients $e,ae\ne0$. Left membership gives diagonal pairs
$(v,t),(v,t),(t,v)$. The first columns of
$\mathsf S,\mathsf S\mathsf R$ similarly span $Fe_2$, with
coefficients $-c,-cg\ne0$. Equality of the first two right pairs
therefore gives $(t^{-1},t^{-1})=(v^{-1},t^{-1})$, hence $t=v=s_0$.
After scalar division the projected log equations are
\begin{equation}
\label{m6:s:eq:late-84-g}
\begin{array}{c|ccc}
&1&2&3\\\hline
N&\dfrac{cd}{ae}zE_{12}&\dfrac de zE_{12}&zE_{21}\\[3pt]
N'&-\dfrac{de}{cg}z'E_{12}&-\dfrac dc z'E_{12}&z'E_{21}.
\end{array}
\end{equation}
\begin{equation}
\label{m6:s:eq:late-84-h}
X_{11}=-\frac{de}{cg}T_f,\qquad
Z_{22}=\frac dc T_f,\qquad
A_{12}=\frac{cd}{ae}C_{21},\qquad
B_{12}=\frac de C_{21}.
\end{equation}
All 17 kernel coordinates are free. The four pivots are mandatory
units, so $b,f,h,i$ are arbitrary. With
\begin{equation}\label{m6:s:eq:late-84-generators}
L=\frac{cd}{ae}A_{12}+\frac deB_{12}+C_{21},\qquad
\mathsf{R}_0=-\frac{de}{cg}X_{11}+\frac dcZ_{22}+T_f.
\end{equation}
the uniform reconstruction gives the full group, both ways:
\begin{equation}
\label{m6:s:eq:late-84-i}
\mathcal A_x=
\left\{\left(s_0l\exp Nl^{-1},\;s_0^{-1}r\exp(-N^\natural)r^{-1}\right):
s_0\in F^\times,\ N=N_0+zL+z'\mathsf{R}_0,\ N_0\in\mathfrak k_{84}\right\}.
\end{equation}
It is $\mathbb G_m\times U_{19}$ of dimension 20 at every allowed
specialization.

\phantomsection\label{m6:s:sec:exclusion-density-84}
The complete trace phase is $\Lambda(N)=gA_{12}+aX_{11}$, so
\begin{equation}
\label{m6:s:eq:late-84-j}
\eta=\psi\left(\frac{cdg}{ae}z-\frac{ade}{cg}z'\right).
\end{equation}
Both coefficients are everywhere nonzero. Matrix multiplication gives
$lLl^{-1}=d_5(E_{21})$, $L'=-L^\natural=(cd/(ae))X_{11}+(d/e)Z_{22}+T_f$
and $L^2=(L')^2=0$; hence
\begin{equation}
\label{m6:s:eq:late-84-k}
\left(d_5(I_2+tE_{21}),\;r(I_{10}+tL')r^{-1}\right),\qquad t\in F,
\end{equation}
is an additive stabilizer subgroup. Its right factor is in $U$,
with adjacent blocks
$X=\left(\begin{smallmatrix}0&0\\-cgi/(ae)&cdg/(ae)\end{smallmatrix}\right)$,
$Y=0$, and phase $\psi(tcdg/(ae))$.
The second square-zero direction has phase $\psi(-tade/(cg))$.
The relative density is one by scalar--unipotent reconstruction.
As for P063, the saturated compatible pullback descends to a closed
smooth locus, of dimension $44+8+9-29=32$. Closed-support injection
and the nonzero phase prove the P084 row of
Theorem~\ref{m6:s:thm:full-stabilizer-exclusions}.
Every denominator is a mandatory unit, so all free-entry divisors
and dyadic fields are included.
\subsection{P090: unipotent graph and signed coefficient dictionary}\label{m6:s:cert:p090-compatibility}

\subsubsection{Solving the entire unipotent stabilizer}\label{m6:s:p090:s4-1}

Define \(\sigma(q)=w_0q^tw_0\). On \(U\), the two adjacent blocks transform as 
\[
(X,Y)\longmapsto\bigl(-(X^\sharp)^t,-(Y^\sharp)^t\bigr),
\qquad X^\sharp=j^{-1}X^tj.
\]

Since \(\operatorname{tr}X^\sharp=\operatorname{tr}X\), it follows that

\begin{equation}
\lambda(\sigma(q))=-\lambda(q).
\label{m6:s:eq:p090-4-1}
\end{equation}

At \(mw_0=m(A,B,C)w_0\), solve \(v\,mw_0\,z^t=mw_0\) for \(z\). The complete \(U\times U\) stabilizer is

\begin{equation}
K^U_{mw_0}
=\left\{\left(v,\sigma(m^{-1}v^{-1}m)\right):v\in U\right\}.
\label{m6:s:eq:p090-4-2}
\end{equation}

Every \(v\in U\) gives a solution because \(m\) normalizes \(U\) and \(\sigma(U)=U\); conversely the matrix equation forces this unique second coordinate. It also proves that the full unipotent stabilizer is connected, smooth, and of dimension 22 throughout \(\mathcal P_{90}\).

By \eqref{m6:s:eq:p090-4-1}, the \(F\)-valued additive phase on \eqref{m6:s:eq:p090-4-2} is

\begin{equation}
\lambda(v)+\lambda(m^{-1}vm)
=\operatorname{tr}\bigl((I_2+BA^{-1})X_v\bigr)
 +\operatorname{tr}\bigl((I_2+CB^{-1})Y_v\bigr).
\label{m6:s:eq:p090-4-3}
\end{equation}

The adjacent blocks \(X_v,Y_v\) range independently over all of \(M_2(F)\). Nondegeneracy of the trace pairing therefore gives the necessary and sufficient conditions

\begin{equation}
B=-A,\qquad C=-B=A.
\label{m6:s:eq:p090-4-4}
\end{equation}

 If either coefficient in \eqref{m6:s:eq:p090-4-3} is nonzero, choose an elementary matrix on which the corresponding trace pairing is nonzero and scale it by arbitrary \(t\in F\). The phase then ranges over all of \(F\). Since \(\psi\) is nontrivial, the stabilizer character is nontrivial. Thus the stabilizer character is trivial precisely when its algebraic phase coefficients vanish.

At a general point \(u_0mw_0\), the stabilizer \eqref{m6:s:eq:p090-4-2} is conjugated by \((u_0,1)\). The character is unchanged under this conjugation. Equations \eqref{m6:s:eq:p090-4-4} are consequently valid for all points in the unique parametrization \eqref{m6:s:eq:p090-3-4}.

\subsubsection{The four coefficient matrices in signed coordinates}\label{m6:s:p090:s4-2}

For the signed representative \(W_{90}\), use the coefficients \(\mathsf P=B_l^{-1}A_l,\quad Q=C_l^{-1}B_l,\quad R=B_r^{-1}A_r,\quad S=C_r^{-1}B_r.\)

The upright letter \(\mathsf P\) denotes a coefficient matrix, distinct from the parabolic \(P\). The compatibility equations are

\begin{equation}
R=-\mathsf P,\qquad Q=jS.
\label{m6:s:eq:p090-4-5}
\end{equation}

Indeed, writing \(c_r=\det C_r\), direct multiplication yields \(lW_{90}r^t=m(A_0,B_0,C_0)w_0,\)

\begin{equation}
A_0=c_rA_lA_r^{-1}j,\qquad
B_0=c_rB_lB_r^{-1}j,\qquad
C_0=C_lC_r^t.
\label{m6:s:eq:p090-4-6}
\end{equation}

The identity \(jD^t=(\det D)D^{-1}j\) for \(D\in\mathrm{GL}_2\) proves these expressions, including the similitude scalar \(c_r\). The equation \(B_0=-A_0\) is equivalent to \(R=-\mathsf P\). Under that equality, \(C_0=A_0\) is equivalent to \(C_l^{-1}B_l=jC_r^{-1}B_r,\)

which is \(Q=jS\). Thus \eqref{m6:s:eq:p090-4-5} and \eqref{m6:s:eq:p090-4-4} describe exactly the same compatible locus in their respective coordinates.

There are eight independent scalar equations in \eqref{m6:s:eq:p090-4-5}. Its entire invertible coefficient family is

\begin{equation}
(\mathsf P,Q,R,S)=(\mathsf P,jS,-\mathsf P,S),\qquad \mathsf P,S\in\mathrm{GL}_2(F).
\label{m6:s:eq:p090-4-7}
\end{equation}

There is no determinant divisor beyond the original invertibility conditions: \(\det R=\det\mathsf P,\quad\det Q=\det S,\quad \det\mathsf P\det Q\det R\det S=(\det\mathsf P)^2(\det S)^2\ne0.\)

Every parameter pair in \eqref{m6:s:eq:p090-4-7} is realized rationally, for example by \(C_l=C_r=I_2\), \(B_l=Q\), \(A_l=Q\mathsf P\), \(B_r=S\), \(A_r=SR\).

In the alternative unique coordinates \(g=u\,m(A,B,C)W_{90}\), \eqref{m6:s:eq:p090-3-1} gives \(m(A,B,C)W_{90}=m(Aj,Bj,C)w_0.\)

Accordingly the compatible equations in those coordinates are

\begin{equation}
B=-A,\qquad C=Aj,
\quad m(A,-A,Aj)W_{90}=d_5(Aj)c.
\label{m6:s:eq:p090-4-8}
\end{equation}

Thus the coefficient equations in block-reversal coordinates are \(B=-A,C=A\), whereas those in the unique signed coordinates \(u\,m(A,B,C)W_{90}\) are \(B=-A,C=Aj\). The relation \eqref{m6:s:eq:p090-3-1} is required when passing between them.

\subsection{Rational normalization and density certificates}\label{m6:s:app:normalizations}

The rational pivot domains and absolute Jacobians below are part of each certificate. Passive parameters and nontrivial intermediate density characters remain present throughout the indicated changes of variables.

\subsection{Rational normalization of both P009 families}\label{m6:s:cert:p009-normalization}

\begin{proof}[Certificate for Proposition~\ref{m6:s:prop:p009-rational-orbits}]
First perform \eqref{m6:s:eq:p009-global-normalization}.
On the generic branch let \(p=\widehat p_{12}\),
\(s=\widehat s_{12}\).
With diagonal parameters one, take
\(\alpha=-\widehat p_{11}/p\),
\(\beta=\widehat p_{22}/p\),
\(\varepsilon=\widehat s_{11}/s\), and
\(\varphi=-\widehat s_{22}/s\). The two resulting matrices are
\[
 \begin{pmatrix}0&p\\-\det\widehat{\mathsf P}/p&0\end{pmatrix},
 \qquad
 \begin{pmatrix}0&s\\-\det\widehat{\mathsf S}/s&0\end{pmatrix}.
\]
Diagonal action multiplies their four displayed nonzero entries
by \(A,A/B,B,A^{-1}\), respectively. Choosing \(A=1/p\),
\(B=1/s\) gives the asserted normal form and
\eqref{m6:s:eq:p009-generic-invariants}. The weight formulas prove that
these invariants are constant under the entire gauge. They are
therefore both necessary and sufficient.

On the exceptional branch the normalized equation is
\(\widehat p_{11}+\widehat s_{12}\widehat p_{22}=0\).
All three displayed factors are nonzero. With all diagonal
parameters one, take \(\alpha=-\widehat p_{21}/\widehat p_{22},\qquad \beta=0,\qquad \varepsilon=\widehat s_{11}/\widehat s_{12},\qquad \varphi=-\widehat s_{22}/\widehat s_{12}.\)
These remove \(\widehat p_{21}\) and the two diagonal entries
of \(\widehat{\mathsf S}\). Then take \(B=1/\widehat s_{12},\qquad A=1/(\widehat s_{12}\widehat p_{22}).\)
The resulting \(\mathsf P\) is \(\operatorname{diag}(-1,1)\),
and the lower-left entry of \(\mathsf S\) is exactly
\eqref{m6:s:eq:p009-exceptional-invariant}.
The same weight formulas prove necessity. Every stated invariant
value occurs.

These reductions and their lifts to \(I_9(F)\) and the two
\(H(F)\)-factors are rational over \(F\). They apply to
nonsplit, repeated-root and nonsemisimple coefficients, for all
\(\rho,\sigma,\kappa\in F^\times\).

At \(x_{\rho,\sigma}\) the nilpotent stabilizer is spanned by the
commuting \(Y_{11},Y_{21},Z_{21}\). At \(y_\kappa\) add
\(L_*=-L_A+L_B+L_C\). The only nonzero bracket among these four
generators is \([L_*,Y_{11}]=Y_{21}\). In logarithmic coordinates
\((u_1,u_2,u_3,z)\), multiplication is \((u,z)*(u',z')= \left(u_1+u_1',\ u_2+u_2'+\tfrac12(zu_1'-u_1z'),\ u_3+u_3',\ z+z'\right).\)
This proves the asserted Heisenberg factor and, again, its
translation Jacobian one for every \(\kappa\ne0\).
\end{proof}

\subsection{P006 Gaussian isomorphisms and their densities}\label{m6:s:cert:p006-gaussian}

\subsubsection{Semisimple chart}\label{m6:s:p006:s4-1}

Every point of \(Y_p\) has a unique expression

\begin{equation}
u(z e_1^t)\,x_p(A,b,c)\,l(e_1v),
\qquad z\in F^2,\quad v\in(F^2)^*.
\label{m6:s:eq:p006-4-4}
\end{equation}

In ordinary block coordinates this expression is

\begin{equation}
\begin{pmatrix}
A+pzv&[pz,b]\\
\begin{bmatrix}pv\\c\end{bmatrix}&pE_{11}
\end{pmatrix}.
\label{m6:s:eq:p006-4-5}
\end{equation}

Thus \(z=B_{:1}/p\), \(v=C_{1:}/p\), and the normalized upper block is
\(A_{\rm raw}-B_{:1}C_{1:}/p\). This proves rational surjectivity and uniqueness on the entire fiber.

The four Gaussian directions carry

\begin{equation}
\psi(z_1+v_1)\,dz\,dv.
\label{m6:s:eq:p006-4-6}
\end{equation}

Their block-coordinate Jacobian is \(|p|^4\), a positive factor in the fixed invariant parameter \(p\). Under a diagonal centralizer, the weights on these four coordinates are \(1,h^{-1},1,h\), where \(h=s/t\); their product is one. Transpose exchanges the column and row coordinates with absolute determinant one and preserves their phase.

\subsubsection{Nilpotent chart}\label{m6:s:p006:s4-2}

Every point of \(Y_N\) has a unique expression

\begin{equation}
u(z e_1^t)\,x_N(A,b,c)\,l(e_2v).
\label{m6:s:eq:p006-4-7}
\end{equation}

Its raw blocks are

\begin{equation}
A_{\rm raw}=A+zv,\qquad B=[b,z],\qquad
C=\begin{bmatrix}v\\c\end{bmatrix},\qquad d=N.
\label{m6:s:eq:p006-4-8}
\end{equation}

The inverse is \(z=B_{:2}\), \(v=C_{1:}\), and \(A=A_{\rm raw}-zv\). Its coordinate Jacobian is one. The character-weighted Gaussian measure is

\begin{equation}
\psi(z_1+v_2)\,dz\,dv.
\label{m6:s:eq:p006-4-9}
\end{equation}

Corrected transpose acts on these Gaussian coordinates by
\((z,v)\mapsto(wv^t,z^tw)\), with absolute determinant one, preserving \eqref{m6:s:eq:p006-4-9}.

\subsection{P006 nonzero-leg pivot normalization}\label{m6:s:cert:p006-nonzero-normalization}

Write \(b=(r,q)^t\), \(c=(q,s)\), and take

\begin{equation}
h=\frac{s-r}{2q},\qquad m=\frac{r+s}{2}.
\label{m6:s:eq:p006-6-10}
\end{equation}

The operation \eqref{m6:s:eq:p006-6-4} takes the two vectors to \(b'=(m,q)^t,\qquad c'=(q,m).\)

Moreover its phase is one, because

\begin{equation}
h(r-s)+2qh^2=0.
\label{m6:s:eq:p006-6-11}
\end{equation}

Let \(A'\) be its new upper block. The residual translations

\begin{equation}
u=\begin{pmatrix}-A'_{11}/q\\-A'_{21}/q\end{pmatrix},
\qquad
v=\left(0,-A'_{22}/q+mA'_{21}/q^2\right)
\label{m6:s:eq:p006-6-12}
\end{equation}

take \(A'\) to \(\kappa E_{12}\), where

\begin{equation}
\kappa=A'_{12}-\frac{m}{q}(A'_{11}+A'_{22})
          +\frac{m^2}{q^2}A'_{21}.
\label{m6:s:eq:p006-6-13}
\end{equation}

\subsection{Incidence-cell rational normalization and density certificates}\label{m6:s:app:p78-normalization}

\subsubsection{Incidence root and Levi notation}\label{m6:s:cert:p78-notation}

\phantomsection\label{m6:s:p78:s1-2}

All block matrices have five blocks of size two unless otherwise stated. Set

\[
j=\begin{pmatrix}0&-1\\1&0\end{pmatrix},
\qquad J=w_5\otimes j,
\qquad G=GSp(J)(F).
\]

For a two-by-two matrix, \(A^\sharp=\operatorname{tr}(A)I_2-A\). The fixed Levi is

\begin{equation}
m(A,B,C)=\operatorname{diag}
\left(A,B,C,\frac{\det C}{\det B}B,
                 \frac{\det C}{\det A}A\right).
\label{m6:s:eq:p78-1-2}
\end{equation}

Its three factors are unrestricted \(GL_2\)'s. Its multiplier is \(\det C\).

We use the following already fixed root-coordinate notation:

\begin{equation}
\mathcal N(X,Y,Z,W,T,S)=
\begin{pmatrix}
0&X&Z&W&S\\
0&0&Y&T&-W^\sharp\\
0&0&0&-Y^\sharp&-Z^\sharp\\
0&0&0&0&-X^\sharp\\
0&0&0&0&0
\end{pmatrix},
\label{m6:s:eq:p78-1-3}
\end{equation}

where \(X,Y,Z,W\in M_2\) and \(T,S\in\mathfrak{sl}_2\). Thus \(U=\exp\mathfrak u\), \(\dim U=22\), and \(H=U\rtimes d_5(GL_2),\qquad \xi(u d_5(a))=\psi(\operatorname{tr}X_u+\operatorname{tr}Y_u).\)

The symbol \(\mathsf X_{ij}\) means \eqref{m6:s:eq:p78-1-3} with \(X=E_{ij}\) and all other coordinates zero; similarly for \(\mathsf Y,\mathsf Z,\mathsf W\). Put \(\mathsf T_f=\mathcal N(T=E_{21})\), \(\mathsf T_e=\mathcal N(T=E_{12})\). The symbols \(\mathsf A_{21},\mathsf B_{21},\mathsf C_{21}\) mean the Levi Lie elements obtained by differentiating the first, second, or third factor of \eqref{m6:s:eq:p78-1-2} in direction \(E_{21}\). Finally, \(\mathsf D_{ij}=\operatorname{diag}(E_{ij},E_{ij},E_{ij},E_{ij},E_{ij})\)

is a Lie element of \(H_0\). These distinctions matter in the stabilizer tables below.

\subsubsection{Affine fiber coordinates and rational uniqueness}\label{m6:s:cert:p78-affine}

The kernel $K_8$ in \eqref{m6:s:eq:p78-3-4}, with brackets
\eqref{m6:s:eq:p78-3-5}, is a one-dimensional additive group times a
three-dimensional Heisenberg group.

Let \(\mathfrak u_i^\circ\) be the span of all the root coordinates in \eqref{m6:s:eq:p78-1-3} except those listed for \(K_i\) in \eqref{m6:s:eq:p78-3-4}, and put \(U_i^\circ=\exp\mathfrak u_i^\circ\). Direct brackets give

\begin{equation}
[\mathfrak u,\mathfrak u_i^\circ]\subset\mathfrak u_i^\circ,
\qquad
U=U_i^\circ\rtimes K_i.
\label{m6:s:eq:p78-3-6}
\end{equation}

For $i=7$ the omitted roots have degree one. For $i=8$ the only
omitted higher root $T_{21}$ is produced by the two omitted second-row
$Y$ roots; a retained first-row $Y$ root contributes only retained
$T_{11}-T_{22},T_{12}$, and an $X$ bracket contributes only retained
higher blocks. This proves the ideal assertion.

The dimensions are \(\dim U_7^\circ=20\), \(\dim U_8^\circ=18\). Exponential coordinates and the semidirect product in \eqref{m6:s:eq:p78-3-6} give polynomial rational inverse maps. Once a rational local section \((l,r)\) of \eqref{m6:s:eq:p78-3-3} is chosen, every point of the fiber has the unique form

\begin{equation}
{g=l u^\circ W_i v^t r^t,
\qquad u^\circ\in U_i^\circ,\quad v\in U.}
\label{m6:s:eq:p78-3-7}
\end{equation}

To check uniqueness, split the first unipotent as \(u^\circ k\) and move \(k\) through \(W_i\) into the second unipotent. The graph kernel in \eqref{m6:s:eq:p78-3-1} is exactly the ambiguity so removed. This proves that the fibers of \(\mathcal P_i\to\mathscr Y_i\) are rational affine spaces of dimensions 42 and 40, respectively.

\subsubsection{Root products for the intermediate modular lines}\label{m6:s:cert:p78-root-densities}

Multiplying the eight weights of
\eqref{m6:s:eq:p78-3-8} gives the intermediate modular character
used in the main text. For $\mathcal P_8$, the eight dimensions are
five one-dimensional roots, a two-dimensional $C$-summand of determinant
$t_4^2/c$, and one further root of weight $t_4^2/c$.
The exact weights and both products are in
\texttt{anc/}\allowbreak\texttt{model6\_calculations.pdf}, [M6-CERT-INCIDENCE-NORMALIZATION].
These are intersection characters; their restriction to the full
stabilizer determines the orbital densities.

\subsubsection{Rational transporters for normalization}\label{m6:s:cert:p78-transporters}

\phantomsection\label{m6:s:p78:s4-2}

If \(k\in J_i\), write \(k=k_Uk_L,\qquad k'=W_i^{-1}kW_i=k'_Uk'_L,\)

where \(k_U\in U\), \(k'_U\in U^t\), and the other factors are Levi factors. Replacing

\begin{equation}
(l,r)\longmapsto(lk_L,r(k'_L)^{-t})
\label{m6:s:eq:p78-4-3}
\end{equation}

is realized on the group point by

\begin{equation}
{
 h_1=l k_U^{-1}l^{-1},\qquad
 h_2=r(k'_U)^t r^{-1},
\qquad x_i' =h_1x_i h_2^t.
}
\label{m6:s:eq:p78-4-4}
\end{equation}

Subsequent normalization of the two third Levi factors uses their inverses in the two \(H_0\) factors and has character one. Formula \eqref{m6:s:eq:p78-4-4} is a rational formula for the original \(H\times H\) transporter at every elementary normalization below.

\subsubsection{Normalized diagonal and shared-central actions}\label{m6:s:cert:p78-gauges}

\phantomsection\label{m6:s:p78:s4-4}

For \(\mathcal P_7\), \eqref{m6:s:eq:p78-3-9}, followed by third-factor normalization, gives

\begin{equation}
\begin{aligned}
A&\mapsto\operatorname{diag}(t_5^{-1},t_6^{-1})A\operatorname{diag}(t_1,t_2),\\
B&\mapsto\operatorname{diag}(t_5^{-1},t_6^{-1})B\operatorname{diag}(t_3,t_4),\\
D&\mapsto\operatorname{diag}(t_4,t_5t_6/t_4)D\operatorname{diag}(t_1^{-1},t_3^{-1}),\\
E&\mapsto\operatorname{diag}(t_4,t_5t_6/t_4)E\operatorname{diag}(t_2^{-1},t_5^{-1}).
\end{aligned}
\label{m6:s:eq:p78-4-6}
\end{equation}

For \(\mathcal P_8\), \eqref{m6:s:eq:p78-3-11} gives

\begin{equation}
\begin{aligned}
A&\mapsto C^{-1}A\operatorname{diag}(t_1,t_2),&
B&\mapsto C^{-1}B\operatorname{diag}(t_3,t_4),\\
D&\mapsto C^tD\operatorname{diag}(t_1^{-1},t_3^{-1}),&
E&\mapsto C^tE\operatorname{diag}(t_2^{-1},t_4/\det C).
\end{aligned}
\label{m6:s:eq:p78-4-7}
\end{equation}

These formulas, the independent lower roots, and \eqref{m6:s:eq:p78-4-4} provide both rational existence and a complete test for equivalence of normal forms. None uses an algebraic-closure orbit classification.

\subsubsection{P007 normalization on the nonzero q branch}\label{m6:s:cert:p78-p7-q}

Let \(q=B_{12}=E_{12}\ne0\). A lower operation on the third left factor, followed by normalizing it back to \(I_2\), replaces \(B\) by \(n_-(-u)B\). Take \(u=B_{22}/q\), and then use a lower-column operation on \(B\) with parameter \(-B_{11}/q\). This gives

\[
B=\begin{pmatrix}0&q\\c&0\end{pmatrix},\qquad c=-\frac{\det B_{\rm old}}q\ne0.
\]

Do the same independently on the right to obtain

\[
E=\begin{pmatrix}0&q\\d&0\end{pmatrix},\qquad d=-\frac{\det E_{\rm old}}q\ne0.
\]

These four operations have zero phase. In \eqref{m6:s:eq:p78-4-6} choose

\begin{equation}
(t_1,t_2,t_3,t_4,t_5,t_6)
=(1,qd,c^{-1},q^{-1},1,1).
\label{m6:s:eq:p78-6-1}
\end{equation}

\subsubsection{Completeness of the P007 nonzero-q residual gauge}\label{m6:s:cert:p78-p7-gauge-test}

To verify completeness, a general left third-factor lower operation and left second-factor lower operation act on \(f\) by

\[
n_-(-u)f n_-(v)=\begin{pmatrix}v&1\\1-uv&-u\end{pmatrix}.
\]

Being anti-diagonal forces \(u=v=0\). The same holds on the right. The remaining diagonal constraints are
\(t_2=t_3=t_6\), \(t_4=t_5\); their two effective ratios are exactly \(h=t_6/t_5\), \(r=t_1/t_6\). These ratios give the entire effective action on \eqref{m6:s:eq:p78-6-3}.

On this remaining torus the density and relative modular character are
both one, as recorded in \eqref{m6:s:eq:p78-6-4}.

\subsubsection{P007 invariants before Gaussian normalization}\label{m6:s:cert:p78-p7-invariants}

For an explicit description of the invariant functions before Gaussian normalization, with \(C_l=C_r=I_2\), put \(q=B_{12}\). Then

\begin{equation}
\begin{aligned}
p&=\frac{\det B\det E}{q}(B^{-1}A)_{12},\\
k&=\frac{\det B\det E}{q}\det A\det D,\\
be&=\frac{\det B\det E}{q^2}A_{12}D_{12}.
\end{aligned}
\label{m6:s:eq:p78-6-9}
\end{equation}

The axis conditions are respectively $A_{12}=0,D_{12}=0$.
These formulas are invariant under every residual gauge and agree
on the rational chart overlaps by the complete normal-form comparison.

\subsubsection{P007 affine phase calculation}\label{m6:s:cert:p78-p7-phase}

Before eliminating the two lower first-column entries, write

\begin{equation}
A=\begin{pmatrix}(\alpha+bc)/p&b\\c&p\end{pmatrix},
\qquad
D=\begin{pmatrix}(\delta+ef)/p&e\\f&p\end{pmatrix}.
\label{m6:s:eq:p78-6-11}
\end{equation}

Here \(\alpha=\det A\ne0\), \(\delta=\det D\ne0\). The first-factor lower-column operations change \(c\) by \(pv\) and \(f\) by \(pw\). Formula \eqref{m6:s:eq:p78-4-5} gives the total phase change

\begin{equation}
vf+wc+pvw
=\frac{(c+pv)(f+pw)-cf}{p}.
\label{m6:s:eq:p78-6-12}
\end{equation}

The phase potential is $\psi(cf/p)$ and clearing $c,f$ uses its
inverse. It is symmetric in $c,f$; their weights $r,r^{-1}$ preserve
$dc\,df$. The raw substitutions $A_{21}=pv,D_{21}=pw$ have
Jacobian $|p|^2$, a positive transpose-invariant factor.

\subsubsection{P007 rational normalization on the zero-q branch}\label{m6:s:cert:p78-p7-qzero}

First make \(B,E\) diagonal by their lower-column operations. Write their diagonal entries as \((b_1,b_2)\), \((e_1,e_2)\). All are nonzero. The choice

\begin{equation}
(t_1,t_2,t_3,t_4,t_5,t_6)
=(1,e_1/b_2,1/b_1,1/b_2,1,1)
\label{m6:s:eq:p78-6-19}
\end{equation}

gives \(B=I_2,\qquad E=\operatorname{diag}(1,z), \qquad z=b_2e_2\ne0, \qquad A_{12}=D_{12}=p\ne0.\)

A third-left-factor lower operation with parameter \(A_{22}/p\), compensated by a second-left-factor lower operation, sets \(A_{22}=0\) while restoring \(B=I_2\). On the right use the parameter \(D_{22}/p\) and compensate in the second factor with that parameter divided by \(z\). These operations have zero phase. Then the first-factor lower operations clear \(A_{11},D_{11}\). The phase before this last clearing is \(\psi(A_{11}D_{11}/p)\), with the exact difference formula \eqref{m6:s:eq:p78-6-12}.

The remaining diagonal action scales the three nonzero entries of

\[
A=\begin{pmatrix}0&p\\a&0\end{pmatrix},\qquad
D=\begin{pmatrix}0&p\\d&0\end{pmatrix}
\]

by \(p\mapsto hp,\qquad a\mapsto ra,\qquad d\mapsto (hr)^{-1}d.\)

Setting $h=p^{-1},r=a^{-1}$ gives \eqref{m6:s:eq:p78-6-20}.

\subsubsection{P008 global compatible normalization}\label{m6:s:cert:p78-p8-global}

In \eqref{m6:s:eq:p78-4-7}, take \(C=B\) and all four scalars initially one. This makes \(B=I_2\). Since \(B^t(-jBe_2)=\det B\,e_1,\)

the new \(E\) has second column \(e_1\). It is therefore

\[
E=\begin{pmatrix}e&1\\d&0\end{pmatrix},\qquad d\ne0.
\]

A right second-factor lower-column operation clears \(e\). Taking \(t_2=d\), with the other scalars one and \(C=I_2\), gives \eqref{m6:s:eq:p78-8-3}.

\subsubsection{Completeness of the P008 residual gauge}\label{m6:s:cert:p78-p8-gauge-test}

For \(h,r\in F^\times\), \(z\in F\), realize \eqref{m6:s:eq:p78-8-4} using \eqref{m6:s:eq:p78-4-7} with \(t_1=r,t_2=t_4=1,t_3=h,C=C_{h,z}\). Restore \(B\) by a left second-factor lower operation of parameter \(z\) and \(E\) by a right one of parameter \(-z\); both have phase zero.

Conversely, preserving \(B=I_2\) forces the shared matrix \(C\) to be lower triangular after its independent second-factor lower operation. Preserving \(E=f\) then forces \(t_2=t_4\) and determines the right second-factor lower parameter. Removing the common scalar gives exactly \eqref{m6:s:eq:p78-8-4}. Thus the remaining list is complete.

\subsubsection{P008 nonzero-branch normalization}\label{m6:s:cert:p78-p8-nonzero}

After \eqref{m6:s:eq:p78-8-3}, let \(a=A_{22},b=D_{12}\), \(\alpha=\det A\). In \eqref{m6:s:eq:p78-8-4} choose

\begin{equation}
h=p,\qquad r=p/\alpha,\qquad z=(a-b)/2.
\label{m6:s:eq:p78-9-1}
\end{equation}

The result has \(p=1\), equal two remaining second-column entries
\(m=(a+b)/2\), and \(\det A=1\). The first-factor lower operations now clear \(A_{11}\) and \(D_{21}\). Invertibility forces the other entries, giving \eqref{m6:s:eq:p78-9-2}.

\subsubsection{P008 nonzero-branch invariant comparison}\label{m6:s:cert:p78-p8-invariants}

For example, after the first operation \(C=B\), the determinant of the resulting \(E\) remains \(\det E\), so its lower-left entry is \(-\det E\). Taking \(t_2=-\det E\) gives

\begin{equation}
p=-\det E\,(B^{-1}A)_{12},\quad
 a=-\det E\,(B^{-1}A)_{22},\quad
 b=(B^tD)_{12}.
\label{m6:s:eq:p78-9-4}
\end{equation}

Equations \eqref{m6:s:eq:p78-8-5}--\eqref{m6:s:eq:p78-8-6} give
\eqref{m6:s:eq:p78-9-3}. Between canonical forms $p=1$ forces $h=1$;
equal second-column entries force $z=0$ and fix $m$; $\det A=1$
forces $r=1$; the first-column zeros force both lower parameters
to vanish. Thus $(k,m)$ are complete.

\subsubsection{P008 nonzero-branch transpose restoration}\label{m6:s:cert:p78-p8-transpose}

Starting from the transpose Levi pair gives
\((A,B;D,E)\mapsto(D,E\operatorname{diag}(1,-1);A,B)\).
For \eqref{m6:s:eq:p78-9-2}, the first shared-central correction uses \(C=j\), followed by \(t_2=-1\). The resulting two first-factor matrices are

\[
\begin{pmatrix}0&-1\\-k&m\end{pmatrix},\qquad
\begin{pmatrix}-1&m\\0&-1\end{pmatrix}.
\]

The restoration \(h=-1,r=1/k,z=0\) gives \eqref{m6:s:eq:p78-9-2}.
The total original reductive
transporter is \((-f,-f)\), equivalent to \((f,f)\) on the group
point. Its character is one, and every density factor in this
restoration has absolute value one.

\subsubsection{P008 boundary normalization}\label{m6:s:cert:p78-p8-zero}

In \eqref{m6:s:eq:p78-8-3}, invertibility now gives

\[
A=\begin{pmatrix}c&0\\d&a\end{pmatrix},\qquad
D=\begin{pmatrix}e&b\\f_0&0\end{pmatrix},
\qquad a b c f_0\ne0.
\]

Use the first lower-column parameters \(-d/a\) and \(-e/b\), then use \eqref{m6:s:eq:p78-8-4} with \(h=cf_0,\qquad r=f_0,\qquad z=0.\)
This gives exactly the complete normal form \eqref{m6:s:eq:p78-9-8}.

\subsubsection{P008 boundary transpose and density comparison}\label{m6:s:cert:p78-p8-density-swap}

For a precise comparison of stabilizer coordinates, the automorphism induced by \(\theta\) maps the basis in \eqref{m6:s:eq:p78-9-12} at \((a,b)\) to the basis at \((b,a)\) by

\begin{equation}
K_1\mapsto-\frac a b K_1,\quad
K_2\mapsto-K_2,\quad K_3\mapsto K_3,\quad
K_4\mapsto-K_4,\quad V\mapsto-V,
\qquad s\mapsto s^{-1}.
\label{m6:s:eq:p78-9-18}
\end{equation}


\subsection{P016: intermediate determinant weights}\label{m6:s:cert:p016-weights}

For diagonal pairs $(\alpha,\beta),(\gamma,\delta),(\epsilon,\zeta)$,
the determinant on the nine radical roots is $(\beta/\epsilon)^6$,
and on $K$ it is $(\beta/\epsilon)^3$. These follow by multiplying
the weights of the exact roots in
Appendix~\ref{m6:s:cert:p016-intersection}; the individual factors are in
\texttt{anc/}\allowbreak\texttt{model6\_calculations.pdf}, [M6-CERT-INTRINSIC-NORMALIZATION].
Thus $q=a_2/c_1=\beta/\epsilon$ gives exactly the modulus and
orbit-Jacobian formulas of Section~\ref{m6:s:sec:middle-p016}.

\subsection{P016: rational quotient fibers and sections}\label{m6:s:cert:p016-fibers}

Suppose two tuples have the same \((N,\kappa)\). Their nonzero factors
satisfy \(v'=a_2v\) and \(\phi'=r\phi\) for unique nonzero scalars.
Equality of \(N\) gives \(d'=a_2rd\). Put \(c_1=r^{-1}\).
Then \(A=\mathsf P^{-1}\mathsf P'\) is lower triangular with second
diagonal entry \(a_2\), and \(C=\mathsf Q(\mathsf Q')^{-1}\) is lower
triangular with first diagonal entry \(c_1\). Set
\(A_r=(\mathsf R')^{-1}\mathsf R\). Since
\(\mathsf R e_2=\phi^t\) and \(\mathsf R'e_2=r\phi^t\), its second
column is \(c_1e_2\). The identity \(d'=(a_2/c_1)d\) forces its first
diagonal entry to equal that of \(A\). Finally
\(C_r=\mathsf S'\mathsf S^{-1}\) has first row \((a_2,0)\), and equality
of \(\kappa\) gives \(\det C_r=c_1c_2\). These are exactly the linked
projection conditions. Each gauge parameter is uniquely determined,
proving rational transitivity and freeness on every fiber.

For surjectivity, factor a nonzero rank-one matrix as \(N=v\phi\)
on any nonzero-entry chart. Complete \(v\) to the second column of
\(\mathsf P\) with determinant one, and \(\phi^t\) to the second column
of \(\mathsf R\) with determinant one. Complete \(\phi\) to the first
row of \(\mathsf Q\) with determinant \(\kappa\), and \(v^t\) to the
first row of \(\mathsf S\) with determinant one. For example, a column
\((v_1,v_2)^t\) is completed by \((v_2^{-1},0)^t\) if \(v_2\ne0\),
or by \((0,-v_1^{-1})^t\) if \(v_1\ne0\); the corresponding row
formulas follow by transpose and a prescribed determinant. Thus only a
chosen nonzero vector coordinate is divided by. These rational local
sections have \(d=1\), cover nilpotents as well, and identify each chart
preimage with \(K_8\) times the chart. The quotient is a surjective
analytic submersion.

\subsection{P017: radical, kernel and orbit determinants}\label{m6:s:cert:p017-weights}

For the exact twelve radical roots, the radical and kernel determinants are
\[
\frac{\beta^6\delta^6}{\epsilon^9\zeta^3},\qquad
\frac{\beta^3\delta^4}{\epsilon^5\zeta^2}.
\]
The complete weight list is recorded in [M6-CERT-INTRINSIC-NORMALIZATION] of the Model~6 ancillary calculations.
The raw 38-dimensional $U^2/K$ orbit Jacobian is therefore
$j=|\beta^6\delta^6/(\epsilon^9\zeta^3)|$ and $\Delta_J=j^{-1}$.

\subsection{P017: six rational pivot operations}\label{m6:s:cert:p017-normalization}

The complete coefficient gauge is \((\mathsf{P},\mathsf{Q},\mathsf{R},\mathsf{S})\mapsto(B^{-1}\mathsf{P}A,C^{-1}\mathsf{Q}B,E^{-1}\mathsf{R}D,F^{-1}\mathsf{S}E).\)

The six projected lower operations use parameters
$-a/b,d/b,(f+gd/b)/e,-h/e,j/e,(k+lj/e)/b$ in the
$A,B,C,D,E,F$ slots. They lift rationally through the six retained roots
and give
\[
\mathsf P_u=\begin{pmatrix}0&b\\-\det\mathsf P/b&0\end{pmatrix},\quad
\mathsf Q_u=\operatorname{diag}(e,g),\quad
\mathsf R_u=\begin{pmatrix}0&e\\-\det\mathsf R/e&0\end{pmatrix},\quad
\mathsf S_u=\operatorname{diag}(b,l).
\]
The torus choice
$(\beta,\gamma,\epsilon,\delta,\zeta,\alpha)
=(1,b,eb,l^{-1},g/l,-b/(\det\mathsf P\,l))$
in \eqref{m6:s:eq:middle-a17-1} gives
\[
\mathsf P_\kappa=\begin{pmatrix}0&1\\1&0\end{pmatrix},\quad
\mathsf Q_\kappa=\mathsf S_\kappa=I,\quad
\mathsf R_\kappa=\begin{pmatrix}0&1\\\kappa&0\end{pmatrix},\quad
\kappa=\det\mathsf P\det\mathsf R\,gl.
\]
The full gauge preserves $\kappa$ and every divisor is forced nonzero.
The expanded projected matrices are in
\texttt{anc/}\allowbreak\texttt{model6\_calculations.pdf}, [M6-CERT-INTRINSIC-NORMALIZATION]. These are explicit rational transporters for every \(F\)-point. Thus \(\kappa\) is a complete rational orbit invariant.

\subsection{P027: full rational gauge and invariant}\label{m6:s:cert:p027-normalization}

For \(k\in I_w\), the replacement \((l,r)\mapsto(lk^{-1},rk^\natural)\) preserves \(lwr^t\). Absorbing the resulting unipotent parts and common third Levi blocks into the original \(H\times H\) gives \((\mathsf{P},\mathsf{Q},\mathsf{R},\mathsf{S})\mapsto (B_0\mathsf{P}A_0^{-1},C_0\mathsf{Q}B_0^{-1},(B'_0)^{-1}\mathsf{R}A'_0,(C'_0)^{-1}\mathsf{S}B'_0).\)
Thus this is a rational orbit equivalence.

The permitted triangular left blocks
\[
A_0=\begin{pmatrix}D/d&0\\c&d\end{pmatrix},\quad
B_0=\begin{pmatrix}1&-b/d\\0&1\end{pmatrix},\quad
C_0=\begin{pmatrix}1/e&0\\-(g+fb/d)/(ef)&1/f\end{pmatrix}
\]
give an \(m(A_0,B_0,C_0)\in I_w(F)\), of multiplier \(1/(ef)\). Direct ten-by-ten multiplication yields
\[
\mathsf{P}_0=I,\quad \mathsf{Q}_0=s=\begin{pmatrix}0&1\\1&0\end{pmatrix},\quad
\mathsf{R}_0=\begin{pmatrix}efhD&di\\ekD/d&1\end{pmatrix},\quad
\mathsf{S}_0=\begin{pmatrix}0&1\\-1&fm\end{pmatrix}.
\]
Write \(i_0=di,\ k_0=ekD/d,\ m_0=fm\). The root-subgroup product \((I-k_0N_{W_{12}})(I+i_0N_{Z_{22}}) (I+(i_0-m_0)N_{T_f})\)
has identity left Levi and right Levi blocks
\[
\begin{pmatrix}1&0\\-k_0&1\end{pmatrix},\quad
\begin{pmatrix}1&i_0\\0&1\end{pmatrix},\quad
\begin{pmatrix}1&0\\m_0-i_0&1\end{pmatrix}.
\]
These identities include the full matrices' commutator terms. The resulting tuple is \(\mathsf{P}=I,\quad \mathsf{Q}=s,\quad \mathsf{R}=\operatorname{diag}(t,1),\quad \mathsf{S}=-j, \qquad t=eDD_R\in F^\times.\)
Only forced nonzero \(d,e,f,D\) occur as denominators. This covers \(a=h=0\), \(di=fm\), vanishing free entries and all square classes of \(t\).

Invariance must hold under the entire gauge. On its six-dimensional torus write the first three diagonal pairs as \((\alpha,\beta),(\gamma,\delta),(\epsilon,\zeta)\). The right pairs are \((\alpha,\epsilon\zeta/\gamma),\quad (\epsilon\zeta/\beta,\epsilon),\quad (\delta,\epsilon\zeta/\delta).\)
The factors \(e,\det \mathsf{P},\det \mathsf{R}\) acquire respectively \(\epsilon/\delta,\qquad \gamma\delta/(\alpha\beta),\qquad \alpha\beta/(\gamma\epsilon),\)
whose product is one. Each of the 16 radical root groups has determinant-one Levi projections and leaves the linked coefficient \(e\) unchanged. The full group decomposition above proves invariance under all gauges. Consequently \(t:C_{27}\to F^\times\) is a complete rational invariant and each fiber is a single \(\mathcal A(F)\)-orbit.

\subsection{P027: kernel-weight certificate}\label{m6:s:cert:p027-weights}

The ten kernel weights consist of six $1$'s, two $\lambda/\mu$'s
and two $\mu/\lambda$'s, so their product is one. Directions with
nontrivial weights are $(Y_{11},T_e)$ and $(W_{21},S_f)$,
respectively; the full table is in
\texttt{anc/}\allowbreak\texttt{model6\_calculations.pdf}, [M6-CERT-INTRINSIC-NORMALIZATION].

\subsection{P031: rational normalization}\label{m6:s:cert:p031-normalization}

For sufficiency take first \(A=C=D=I\), \(B=\mathsf P\), and

\[
E=\begin{pmatrix}a&0\\-e/(dc)&1/c\end{pmatrix};
\]

This gives \(\mathsf P=I\) and \(\mathsf R=\operatorname{diag}(\rho,-1)\). Next choose \(A=B=I\), take \(C\) equal to the resulting \(\mathsf Q\), and put \(D=E=\operatorname{diag}(1,\det\mathsf Q)\). The new matrices have \(\mathsf Q=I\) and

\[
\mathsf{S}=\begin{pmatrix}\nu&1\\\kappa&0\end{pmatrix},\quad\kappa\ne0.
\]

The upper-right entry is one because \(\mathsf S e_2=-j\mathsf Qe_2\)
and \(\mathsf Q^tj\mathsf Q=(\det\mathsf Q)j\). Now choose
\[
A=B=C=\operatorname{diag}(1,\kappa^{-1}),\qquad
E=\begin{pmatrix}1&0\\\nu&1\end{pmatrix},\qquad
D=\begin{pmatrix}1&0\\-\rho\nu&1\end{pmatrix}.
\]
This yields the normal form

\begin{equation}
\label{m6:s:eq:sector-31n}
\mathsf{P}=\mathsf{Q}=I,\qquad \mathsf{R}=\operatorname{diag}(\rho,-1),\qquad \mathsf{S}=M=\begin{pmatrix}0&1\\1&0\end{pmatrix}.
\end{equation}

All denominators are forced nonzero on the whole family. In particular the construction covers \(g=0\) and \(h=0\), arbitrary original Jordan forms and nonsplit coefficient matrices. Equality of \(\rho\) is necessary and sufficient over \(F\).

\subsection{P040: rational pivots and diagonal weights}\label{m6:s:cert:p040-normalization}

Let \(L(z)=\left(\begin{smallmatrix}1&0\\z&1\end{smallmatrix}\right)\). Use the six available lower unipotents
\[
\begin{split}
 B&=L(-d/e),\quad C=L(g/e),\quad A=L(-(b+ad/e)/c),\\
 E&=L(-i/c),\quad F=L(k/c),\quad D=L(-(h+ai/c)/e).
\end{split}
\]
Direct multiplication gives
\[
 \mathsf{P}=\operatorname{diag}(a,c),\quad \mathsf{R}=\operatorname{diag}(a,e),
 \quad \mathsf{Q}=\begin{pmatrix}0&e\\f_0&0\end{pmatrix},\quad
 \mathsf{S}=\begin{pmatrix}0&c\\j_0&0\end{pmatrix},
\]
where \(f_0=-\det \mathsf{Q}/e\ne0\) and \(j_0=-\det \mathsf{S}/c\ne0\). The only denominators so far are the forced nonzero entries \(c,e\).

Under \eqref{m6:s:eq:late-p40-4}, the five remaining entries \((a,c,e,f_0,j_0)\) have weights 
\[
\alpha _1/\beta _1,\quad\alpha _2/\beta _2,\quad
 \beta _2/\gamma _1,\quad\beta _1/\gamma _2,\quad
 \gamma _1\gamma _2/(\alpha _1\alpha _2).
\]

Their product is invariant. The rational choice 
\[
\alpha _1=a^{-1},\ \alpha _2=(ce)^{-1},\ \beta _1=1,
 \ \beta _2=e^{-1},\ \gamma _1=1,\ \gamma _2=f_0
\]

gives the canonical tuple \eqref{m6:s:eq:late-p40-6}.

\subsection{P034: finite rational charts before full compatibility}\label{m6:s:cert:p034-pair-charts}

\subsubsection{Finite rational charts for the remaining pair}\label{m6:s:p034:s2-3}

The following cyclic-vector charts normalize the pair over $F$. For a nonscalar matrix \(B\), at least one of

\begin{equation}
v=e_1,\quad e_2,\quad e_1+e_2
\label{m6:s:eq:p034-2-8}
\end{equation}

is cyclic, meaning \(\det(v,Bv)\ne0\). Indeed, if both coordinate vectors are noncyclic, \(B\) is diagonal; if it is nonscalar, their sum is cyclic.

On the corresponding chart set \(P_B=(v,Bv),\quad t=\operatorname{tr}B,\quad \delta=\det B.\)

Cayley--Hamilton gives

\begin{equation}
P_B^{-1}BP_B=J(t,\delta):=
\begin{pmatrix}0&-\delta\\1&t\end{pmatrix}.
\label{m6:s:eq:p034-2-9}
\end{equation}

The pure reductive element \(a=P_B^{-1}\) normalizes the pair to

\begin{equation}
\bigl(J(t,\delta),\ P_B^{-1}CP_B\bigr).
\label{m6:s:eq:p034-2-10}
\end{equation}

These are rational, finite-chart formulas on the whole nonscalar-\(B\) locus. Two normalized pairs with the same first matrix are equivalent precisely under its full centralizer. If \(B=bI\) and \(C\) is nonscalar, use the same three charts for \(C\). If both are scalar, no conjugation normalization remains.

The finitely many rational cyclic charts differ by the appropriate centralizer action.

\subsection{Tangent, submersion, and conormal certificates}\label{m6:s:app:conormals}

The bases below annihilate the full indicated tangent spaces, including all varying support parameters and the ambient similitude direction. The support dimensions and nonnilpotent witnesses are stated at their uses in the main text.

\subsection{The six- and eight-dimensional P009 conormals}\label{m6:s:cert:p009-conormals}

Use \eqref{m6:s:eq:p009-trace-trivialization}. The ambient Lie algebra is
\(J^{-1}\operatorname{Sym}_{10}\oplus FI_{10}\), of dimension 56.
The summands are trace-orthogonal and nondegenerate, with
\(\operatorname{tr}(I_{10}^2)=10\ne0\).

The physical columns of the two normal forms are
\begin{equation}\label{m6:s:eq:p009-physical-normal-forms}
\begin{aligned}
x_{\rho,\sigma}={}&[e_5,\rho\sigma e_2,e_3,\rho e_1,e_4,
 -e_7,-\rho^{-1}e_{10},e_8,(\rho\sigma)^{-1}e_9,e_6],\\
y_\kappa={}&[e_5,-\kappa e_1,e_3,\kappa e_2,e_4,
 -e_7,\kappa^{-1}e_9,e_8,\kappa^{-1}e_{10},e_6].
\end{aligned}
\end{equation}
Put
\[
\begin{aligned}
 Q_0&=\operatorname{diag}(1,1,0,0,0,0,0,0,-1,-1),\\
 Q_1&=\operatorname{diag}(0,1,0,0,0,0,0,0,-1,0),\\
 D_9&=\operatorname{diag}(0,0,1,0,-1,1,0,-1,0,0),\\
 D_{9,0}&=\operatorname{diag}(1,-1,-1,1,0,0,-1,1,1,-1).
\end{aligned}
\]
Differentiating the columns gives \((\partial_\rho x)x^{-1}=\rho^{-1}Q_0,\qquad (\partial_\sigma x)x^{-1}=\sigma^{-1}Q_1,\qquad (\partial_\kappa y)y^{-1}=\kappa^{-1}Q_0.\)

\begin{lemma}\label{m6:s:cert:p009-bases}
The entire ambient conormals are
\begin{equation}\label{m6:s:eq:p009-generic-conormal}
 CN^G_{\mathcal C_9^\circ,x_{\rho,\sigma}}=
 \operatorname{span}_F
 \{Y_{11},Y_{21},Z_{11},Z_{21}-\sigma^{-1}X_{11},
       L_C-\tfrac12L_B,D_9\},
\end{equation}
and
\begin{equation}\label{m6:s:eq:p009-exceptional-conormal}
\begin{split}
 CN^G_{\mathcal E_9,y_\kappa}=\operatorname{span}_F
 \{&Y_{11},Y_{21},Z_{11}+X_{21},Z_{21},
 L_B-L_A,D_{9,0},\\
 &L_C-\tfrac12L_A,D_9\}.
\end{split}
\end{equation}
At every point of either manifold this full conormal contains a
nonnilpotent element.
\end{lemma}

\begin{proof}
The full tangent dimensions are 50 and 48, as proved in
Proposition~\ref{m6:s:prop:p009-conormals}. The differentials of the invariant
parameters \(\rho,\sigma\), respectively \(\kappa\), are independent
on the corresponding normal-form family and vanish on orbit tangents.
We compute the complete tangent annihilators from
\eqref{m6:s:eq:p009-root-matrix} and \eqref{m6:s:eq:p009-physical-normal-forms}.

Let \(\mathcal B_{\mathfrak h}\) consist of the sixteen
\(X_{ij},Y_{ij},Z_{ij},W_{ij}\), the six
\(T_h,T_e,T_f,S_h,S_e,S_f\), and the four \(d_5(E_{ij})\).
For a cotangent matrix \(Z\), the equations are exactly
\begin{equation}\label{m6:s:eq:p009-conormal-equations}
 \operatorname{tr}(ZA)=0,\quad
 \operatorname{tr}(x^{-1}ZxA^t)=0\quad(A\in\mathcal B_{\mathfrak h}),
 \qquad
 \operatorname{tr}(ZQ_0)=\operatorname{tr}(ZQ_1)=0.
\end{equation}
For \(y\) omit the last \(Q_1\) equation.
Here all entries of the equations are explicit: for a physical
matrix unit \(\operatorname{tr}(ZE_{ab})=Z_{ba}\), while if
the \(j\)-th column in \eqref{m6:s:eq:p009-physical-normal-forms} is
\(c_je_{\pi(j)}\), then \((x^{-1}Zx)_{ij}=c_i^{-1}c_jZ_{\pi(i),\pi(j)}.\)
The entry formula verifies each of the six, respectively eight,
displayed generators against every equation in
\eqref{m6:s:eq:p009-conormal-equations}; the full zero-pairing tables are
in \texttt{anc/}\allowbreak\texttt{model6\_calculations.pdf}, [M6-CERT-SPECTRAL-CONORMALS].
The generic cancellations are $Z_{21}-\sigma^{-1}X_{11}$ and
$L_C-\tfrac12L_B$; the exceptional ones are $Z_{11}+X_{21}$,
$L_B-L_A$ and $L_C-\tfrac12L_A$. Independent root, lower Levi and
diagonal entries give dimensions 6 and 8, and the previously proved
full tangent codimensions prove completeness. The scalar tangent
requires $\operatorname{tr}Z=0$, retaining the full $GSp_{10}$ ambient.

For the nonnilpotent witness, $D_9$ annihilates the left action since
its odd and even diagonal sums vanish and upper roots pair trivially.
The physical columns give
$x^{-1}D_9x=y^{-1}D_9y=\operatorname{diag}(-1,0,1,0,0,0,0,-1,0,1)$,
whose odd and even sums also vanish. Its zero entries on the supports
of $Q_0,Q_1$ prove the parameter equations.
Finally \(\det(tI-D_9)=t^6(t-1)^2(t+1)^2,\qquad \operatorname{tr}(D_9^2)=4\ne0.\)
At \(h_1xh_2^t\), right-trivialized tangents and their trace duals
transform by \(\operatorname{Ad}_{h_1}\). Thus
\(h_1D_9h_1^{-1}\) is a nonnilpotent conormal at every point,
including every similitude and every rational parameter.
\end{proof}

\subsection{The full P012 tangent and ten-dimensional conormal}\label{m6:s:cert:p012-conormals}

We retain the full right-trivialized trace pairing
\eqref{m6:s:eq:p009-trace-trivialization}. At \(x=lW_{12}r_0^t\),
conjugate these tangent and cotangent coordinates by \(l^{-1}\).
The entire tangent image of the saturation map is
\begin{equation}\label{m6:s:eq:p012-full-tangent-span}
\begin{split}
\mathcal T={}&\mathfrak u+\operatorname{Ad}_{W_{12}}\mathfrak u^t\\
&+\operatorname{span}\{l^{-1}d_5(a)l:a\in M_2(F)\}\\
&+\operatorname{span}\{W_{12}r_0^td_5(a)^tr_0^{-t}W_{12}^{-1}:
                         a\in M_2(F)\}\\
&+\operatorname{span}\{l^{-1}\partial_vl+
 W_{12}(\partial_vr_0)^tr_0^{-t}W_{12}^{-1}:v\in\mathcal V\},
\end{split}
\end{equation}
where
\(\mathcal V=(p,q,r,t,b,c,d,e,f,g,h,i)\).
It contains all forty-four unipotent, eight reductive-action,
and twelve coefficient directions. In particular \(b\) and
\(c\) are differentiated simultaneously in every block where
they occur. This is the ambient tangent space.

All derivative matrices are specified by the following
logarithmic derivative:
\begin{equation}\label{m6:s:eq:p012-logarithmic-derivative}
\begin{split}
 m(A,B,I)^{-1}\delta m(A,B,I)=\operatorname{diag}(
 &A^{-1}\delta A,B^{-1}\delta B,0,\\
 &B^{-1}\delta B-\operatorname{tr}(B^{-1}\delta B)I,\\
 &A^{-1}\delta A-\operatorname{tr}(A^{-1}\delta A)I).
\end{split}
\end{equation}
For the left factor take \(A=\mathsf Q\mathsf P,B=\mathsf Q\);
for the right factor take \(A=\mathsf S\mathsf R,B=\mathsf S\)
and transpose after taking the same left logarithmic derivative,
since
\((\delta r_0)^tr_0^{-t}=(r_0^{-1}\delta r_0)^t\).
The variations are exactly
\[
\begin{gathered}
\delta(\mathsf Q\mathsf P)=(\delta\mathsf Q)\mathsf P+
                           \mathsf Q\delta\mathsf P,\qquad
\delta(\mathsf S\mathsf R)=(\delta\mathsf S)\mathsf R+
                           \mathsf S\delta\mathsf R,\\
\delta\mathsf Q=\begin{pmatrix}0&\delta b\\\delta c&\delta d\end{pmatrix},
\quad
\delta\mathsf R=\begin{pmatrix}\delta e&0\\\delta f&\delta c\end{pmatrix},
\quad
\delta\mathsf S=\begin{pmatrix}\delta g&\delta b\\\delta h&\delta i\end{pmatrix}.
\end{gathered}
\]

The intersection
\(\mathfrak u\cap\operatorname{Ad}_{W_{12}}\mathfrak u^t\)
is the six-dimensional \(\mathfrak n_6\).
Hence the first line of
\eqref{m6:s:eq:p012-full-tangent-span} has dimension \(22+22-6=38\).
Its entire annihilator in the full fifty-six-dimensional
\(\mathfrak{gsp}_{10}\) consists precisely of the matrices
\(\mathcal Z(\mathbf z)\) in the complete physical array in
\texttt{anc/}\allowbreak\texttt{model6\_calculations.pdf}, [M6-P01222],
original equation \texttt{m6:s:eq:}\allowbreak\texttt{p012-eighteen-matrix}. Indeed the equations \(Z^tJ+JZ=\lambda J,\quad \operatorname{tr}(ZN)=0,\quad \operatorname{tr}(ZW_{12}N^tW_{12}^{-1})=0 \quad(N\text{ a root basis vector of }\mathfrak u)\)
give exactly the displayed entries. All eighteen free entries
are independent. Their number is \(56-38\), which proves
completeness of this initial annihilator.

We now solve every remaining equation in
\eqref{m6:s:eq:p012-full-tangent-span}.
Put \(\lambda=z_{15}+z_{17}\). Vary \(\mathsf P\)
arbitrarily; \(a=\mathsf P^{-1}\delta\mathsf P\) ranges over
all of \(M_2(F)\). The resulting pairing is
\[
 \operatorname{tr}
 \left[\begin{pmatrix}2z_9-\lambda&0\\2z_{10}&2z_{11}-\lambda\end{pmatrix}a\right].
\]
It vanishes identically exactly when
\(z_{10}=0,\ z_9=z_{11}=\lambda/2\).
The scalar left \(H_0\) direction has pairing
\(\operatorname{tr}Z=5\lambda\); it forces \(\lambda=0\).
Thus the central similitude coordinate has been retained and
eliminated by its tangent equation:
\begin{equation}\label{m6:s:eq:p012-conormal-first-four}
 z_9=z_{10}=z_{11}=0,\qquad z_{17}=-z_{15}.
\end{equation}
The remaining left \(d_5(E_{11})\) and \(d_5(E_{12})\)
pairings are \(-z_{12}+z_{14}+z_{15},\qquad z_{16}+\frac{2c}{b}z_{13} -\frac{2d}{b}(z_{12}-z_{14}).\)
They give
\begin{equation}\label{m6:s:eq:p012-conormal-next-two}
 z_{12}=z_{14}+z_{15},\qquad
 z_{16}=\frac{2d}{b}z_{15}-\frac{2c}{b}z_{13}.
\end{equation}
Define the two remaining expressions \(\Phi=z_3+\frac ce z_4+ \left(\frac fe+\frac gb\right)z_{14}+\frac fe z_{15}, \qquad \Psi=z_8+\frac{2i}{b}z_{14}.\)
After \eqref{m6:s:eq:p012-conormal-first-four} and
\eqref{m6:s:eq:p012-conormal-next-two}, all four left reductive
pairings vanish. The four right pairings, in order
\(E_{11},E_{12},E_{21},E_{22}\), are exactly
\begin{equation}\label{m6:s:eq:p012-right-pairing-table}
 \frac{2bi}{\Delta}\Phi,\qquad
 \frac{2i^2}{\Delta}\Phi-\Psi,\qquad
 -\frac{2b^2}{\Delta}\Phi,\qquad
 -\frac{2bi}{\Delta}\Phi.
\end{equation}
Among the twelve coefficient pairings the only possibly
nonzero ones are
\begin{equation}\label{m6:s:eq:p012-parameter-pairing-table}
 \langle Z,\partial_b\rangle=\frac{2i}{\Delta}\Phi,\qquad
 \langle Z,\partial_i\rangle=-\frac{2b}{\Delta}\Phi.
\end{equation}
Since \(b\Delta\ne0\), the \(i\)-derivative forces \(\Phi=0\)
even when \(i=0\). The right \(E_{12}\) direction then gives
\(\Psi=0\). Conversely the displayed complete pairing lists
prove sufficiency. The final two equations are
\begin{equation}\label{m6:s:eq:p012-conormal-last-two}
 z_3=-\frac ce z_4-\left(\frac fe+\frac gb\right)z_{14}
                    -\frac fe z_{15},\qquad
 z_8=-\frac{2i}{b}z_{14}.
\end{equation}
The eight pivot variables are distinct and all pivot
coefficients are one. The only variable denominators in the
final formulas are \(b,e\), already nonzero.
Thus the full differential has rank forty-six at every point,
including \(t=0\).

\begin{lemma}\label{m6:s:cert:p012-basis}
Put \(U_0=E_{25}+E_{69}\), and define the following ten physical
matrices in the \(l^{-1}\)-conjugated right trivialization:
\begin{equation}\label{m6:s:eq:p012-conormal-basis}
\begin{aligned}
C_1&=E_{35}-E_{68},\\
C_2&=E_{36}+E_{58},\\
C_3&=E_{45}+E_{67},\\
C_4&=E_{18}+E_{3,10}-(c/e)U_0,\\
C_5&=E_{28}-E_{39},\\
C_6&=E_{37}-E_{48},\\
C_7&=E_{38},\\
C_8&=E_{34}+E_{78}-(2c/b)E_{65},\\
C_9&=\operatorname{diag}(0,0,1,1,0,0,-1,-1,0,0)
       \\
   &\qquad -(f/e+g/b)U_0-(2i/b)E_{47},\\
C_{10}&=\operatorname{diag}(0,0,1,0,1,-1,0,-1,0,0)
       \\
   &\qquad -(f/e)U_0+(2d/b)E_{65}.
\end{aligned}
\end{equation}
At \(x=lW_{12}r_0^t\), the full ambient conormal of
\(\mathcal C_{12}\) is
\(\operatorname{Ad}_l\operatorname{span}_F(C_1,\ldots,C_{10})\).
At \(h_1xh_2^t\) it is the corresponding
\(\operatorname{Ad}_{h_1l}\)-translate.
It contains a nonnilpotent element at every point.
\end{lemma}

\begin{proof}
Solving
\eqref{m6:s:eq:p012-conormal-first-four},
\eqref{m6:s:eq:p012-conormal-next-two}, and
\eqref{m6:s:eq:p012-conormal-last-two} in the physical array of
[M6-P01222] in the Model~6 ancillary calculations,
original equation \texttt{m6:s:eq:}\allowbreak\texttt{p012-eighteen-matrix} gives exactly the span in
\eqref{m6:s:eq:p012-conormal-basis}. Its ten coefficients are read,
in order, from entries \[(3,5),(3,6),(4,5),(1,8),(2,8),(3,7),(3,8),(3,4),(4,4),(5,5).\]
The coefficient matrix on these entries is the identity.
Thus all ten vectors are independent uniformly, including
when the optional coefficients vanish. The preceding
elimination proved necessity and sufficiency against every
action and coefficient direction, so this is the entire
ambient conormal. It has dimension \(56-46=10\).

The entire tangent is the trace annihilator of this ten-dimensional
space. In the same conjugated trivialization its exact equations are
\begin{equation}\label{m6:s:eq:p012-tangent-equations}
\operatorname{tr}(BC_j)=0\quad(1\le j\le10),\qquad B\in\mathfrak{gsp}_{10}.
\end{equation}
The ten entry equations follow by
$\operatorname{tr}(BE_{ab})=B_{ba}$. The expanded calculation is recorded as
\texttt{M6-CERT-SPECTRAL-CONORMALS} in
\path{anc/model6_calculations.pdf}.
This is the full support tangent, because the preceding elimination
included every action and coefficient direction.

Take \(D_{12}=C_9\). Its diagonal part, \(U_0\), and \(E_{47}\)
all satisfy \(A^tJ+JA=0\), so \(D_{12}\in\mathfrak{sp}_{10}\).
Its off-diagonal entries are strictly above the physical
diagonal. Hence
\begin{equation}\label{m6:s:eq:p012-witness-polynomial}
 \det(sI-D_{12})=s^6(s-1)^2(s+1)^2,\qquad
 \operatorname{tr}(D_{12}^2)=4\ne0.
\end{equation}
This holds on the entire coefficient envelope, including
\(t=0\) and \(f=g=i=0\). More generally a linear combination
of \(C_1,\ldots,C_{10}\) whose last two coefficients are
\(\omega,\nu\) has characteristic polynomial \(s^4(s^2-\omega^2)(s^2-\nu^2)(s^2-(\omega+\nu)^2).\)
Indeed in physical order \(1,2,3,4,6,5,7,8,9,10\), the matrix
is upper triangular with diagonal
\(0,0,\omega+\nu,\omega,-\nu,\nu,-\omega,-\omega-\nu,0,0\).
Thus its nilpotent part is exactly the eight-dimensional
subspace \(\omega=\nu=0\).

Undoing the conjugated trivialization gives the asserted
\(\operatorname{Ad}_l\) formula. The action
\(x\mapsto h_1xh_2^t\) sends right-trivialized tangent and
cotangent matrices by \(\operatorname{Ad}_{h_1}\), proving
the statement at every saturated point. Rational coverage
was established above, so this includes all rational points
and all multipliers.
\end{proof}

The scalar equation and the divisor tangent are retained in
Proposition~\ref{m6:s:prop:p012-conormal}; the full partner conormal is
transported by \eqref{m6:s:eq:p009-cotangent-transpose}.

\subsection{The full five-dimensional residual P006 conormal}\label{m6:s:cert:p006-conormal}

\begin{lemma}\label{m6:s:cert:p006-basis}
The full conormal of \(\mathscr S\) in \(\mathscr W\) is the five-dimensional space \eqref{m6:s:eq:p006-10-5}. At every point it contains a nonnilpotent cotangent vector.
\end{lemma}
\begin{proof}

Identify a cotangent vector at \(x\in GL_4\) with a matrix \(Z\) by

\begin{equation}
\lambda_Z(\dot x)=\operatorname{tr}(Z\dot x x^{-1}).
\label{m6:s:eq:p006-10-1}
\end{equation}

The convention agrees with Lemma~\ref{m6:s:lem:p006-conormal}.
\phantomsection\label{m6:s:p006:s10-1}

For arbitrary \(\mathsf P,Q,A_0\in M_2\), the full orbit tangent directions are

\begin{equation}
\begin{pmatrix}0&\mathsf P\\0&0\end{pmatrix}x,
\qquad
x\begin{pmatrix}0&0\\Q&0\end{pmatrix},
\qquad
[d(A_0),x],
\quad d(A_0)=\operatorname{diag}(A_0,A_0).
\label{m6:s:eq:p006-10-2}
\end{equation}

Together with \(\partial_\beta x\) and \(\partial_\gamma x\), these span \eqref{m6:s:eq:p006-9-6}. The conormal is therefore obtained by imposing

\begin{equation}
\begin{aligned}
\operatorname{tr}\left(Z\begin{pmatrix}0&\mathsf P\\0&0\end{pmatrix}\right)&=0,\\
\operatorname{tr}\left(x^{-1}Zx\begin{pmatrix}0&0\\Q&0\end{pmatrix}\right)&=0,\\
\operatorname{tr}\bigl((Z-x^{-1}Zx)d(A_0)\bigr)&=0,\\
\operatorname{tr}(Z(\partial_\beta x)x^{-1})&=0,\\
\operatorname{tr}(Z(\partial_\gamma x)x^{-1})&=0.
\end{aligned}
\label{m6:s:eq:p006-10-3}
\end{equation}

At the canonical point the inverse has only four nonzero entries:
\begin{equation}\label{m6:s:eq:p006-10-4}
x^{-1}=-\frac1{\beta\gamma}E_{12}+\frac1\gamma E_{24}
             +\frac1\beta E_{31}+E_{43}.
\end{equation}

Solving \eqref{m6:s:eq:p006-10-3} gives the entire conormal:

\begin{equation}
{
CN^{\mathscr W}_{\mathscr S,x_{\beta,\gamma}}=
\left\{
\begin{pmatrix}
\lambda&-\rho+(\mu-\nu)/\beta&\mu&\sigma\\
0&-\lambda&0&\nu\\
0&0&-\lambda&\rho\\
0&0&0&\lambda
\end{pmatrix}:
\lambda,\mu,\nu,\rho,\sigma\in F
\right\}.
}
\label{m6:s:eq:p006-10-5}
\end{equation}

For a direct check of completeness, the resulting linear conditions on \(Z=(z_{ij})\) are

\begin{equation}
\begin{gathered}
z_{21}=z_{23}=z_{31}=z_{32}=z_{41}=z_{42}=z_{43}=0,\\
z_{11}=z_{44},\qquad z_{22}=z_{33}=-z_{44},\\
z_{12}+z_{34}=(z_{13}-z_{24})/\beta.
\end{gathered}
\label{m6:s:eq:p006-10-6}
\end{equation}

These eleven independent conditions are equivalent to
\eqref{m6:s:eq:p006-10-3}: substitution of \eqref{m6:s:eq:p006-10-5} gives
zero in every equation, and \(16-11=5\) proves completeness.
\phantomsection\label{m6:s:p006:s10-2}

For every matrix in \eqref{m6:s:eq:p006-10-5},

\begin{equation}
\operatorname{tr}(Z^2)=4\lambda^2.
\label{m6:s:eq:p006-10-7}
\end{equation}

Setting \(\lambda=1\), \(\mu=\nu=\rho=\sigma=0\) gives the
main-text witness \(Z_*\); the nilpotent part is precisely
\(\lambda=0\), when the matrix is strictly upper triangular.

At an arbitrary point \(y=LxR\) obtained by a fixed \(K_S\)-action, right-translated cotangent matrices transform by

\begin{equation}
Z\longmapsto LZL^{-1},\qquad L=u(X)d(a).
\label{m6:s:eq:p006-10-8}
\end{equation}

Thus every point of \(\mathscr S\) has a nonnilpotent conormal witness, with the same nonzero trace-square.

For reference, corrected transpose transforms a cotangent matrix by

\begin{equation}
Z\longmapsto d(w)(x^{-1}Zx)^t d(w).
\label{m6:s:eq:p006-10-9}
\end{equation}

For the diagonal witness this is again \(Z_*\) at \(x_{\gamma,\beta}\). Thus the witness is compatible with the transpose exchange of the residual parameters.

\end{proof}

\subsection{The P006 ambient radical and full Levi trace pairing}\label{m6:s:cert:p006-ambient-pairing}

At \(g_0=\iota(x,I_2)\), identify cotangent matrices using

\begin{equation}
\Lambda_Z(\dot g)=\operatorname{tr}(Z\dot g g_0^{-1}).
\label{m6:s:eq:p006-11-4}
\end{equation}

The tangent directions from the two radical factors in \eqref{m6:s:eq:p006-11-1} are \(\mathfrak n'\) and \(\mathfrak n'^t\), after right translation. Their annihilator under the trace form is the Levi Lie algebra. A general Levi Lie element is

\begin{equation}
Z=\operatorname{diag}\left(
 B_0,C_0,\upsilon I_4-J_4^{-1}B_0^tJ_4\right),
\qquad\upsilon=\operatorname{tr}C_0.
\label{m6:s:eq:p006-11-5}
\end{equation}

Pairing this element with the Levi tangent parametrized by \((A_0,D_0)\in\mathfrak{gl}_4\oplus\mathfrak{gl}_2\) gives

\begin{equation}
2\operatorname{tr}(B_0A_0)-\upsilon\operatorname{tr}A_0
+\operatorname{tr}(C_0D_0)
+(4\upsilon-\operatorname{tr}B_0)\operatorname{tr}D_0.
\label{m6:s:eq:p006-11-6}
\end{equation}

Hence its cotangent vector matrices on the two product factors are

\begin{equation}
E=2B_0-\upsilon I_4,
\qquad F_2=C_0+(4\upsilon-\operatorname{tr}B_0)I_2.
\label{m6:s:eq:p006-11-7}
\end{equation}

\subsection{Incidence-cell tangent certificates}\label{m6:s:app:p78-tangents}

\subsubsection{Incidence tangents, submersion and full conormals}\label{m6:s:cert:p78-tangents}

\phantomsection\label{m6:s:p78:s11}

\paragraph{The ambient tangent space of each entire parabolic cell}\label{m6:s:p78:s11-1}

We calculate the full parabolic-cell conormal separately from
the model-orbit annihilator and compatible-locus dimensions.

Identify a cotangent vector at \(g\in G\) with \(Z\in\mathfrak{gsp}_{10}\) using

\begin{equation}
\lambda_Z(\dot g)=\operatorname{tr}(Z\dot g g^{-1}).
\label{m6:s:eq:p78-11-1}
\end{equation}

The trace form is nondegenerate on \(\mathfrak{gsp}_{10}\) in characteristic zero. The scalar-center tangent is included. The Lie algebra \(\mathfrak p\) has dimension \(12+22=34\). For the smooth parabolic cell \(C_i=P W_i P^t\),

\begin{equation}
T_g C_i=\mathfrak p g+g\mathfrak p^t,
\qquad
(T_g C_i)g^{-1}=\mathfrak p+\operatorname{Ad}(g)\mathfrak p^t.
\label{m6:s:eq:p78-11-2}
\end{equation}

The stabilizer of the \(P\times P\)-action is the graph of \(J_i\). Thus \eqref{m6:s:eq:p78-3-12} gives the ranks of the full tangent map:

\begin{equation}
\begin{array}{c|c|c|c}
&\text{source Lie dimension}&\text{kernel dimension}&\dim T_g C_i\\\hline
\mathcal P_7&68&14&54\\
\mathcal P_8&68&16&52.
\end{array}
\label{m6:s:eq:p78-11-3}
\end{equation}

Strict upper-block triangularity gives $\mathfrak u\perp\mathfrak p$;
the equality $56-34=22=\dim\mathfrak u$ then proves

\begin{equation}
\mathfrak p^\perp=\mathfrak u,
\qquad (\mathfrak p^t)^\perp=\mathfrak u^t.
\label{m6:s:eq:p78-11-4}
\end{equation}

Taking the annihilator of \eqref{m6:s:eq:p78-11-2} therefore gives the complete ambient conormal

\begin{equation}
{CN^G_{C_i,g}=\mathfrak u\cap\operatorname{Ad}(g)\mathfrak u^t.}
\label{m6:s:eq:p78-11-5}
\end{equation}

All parabolic tangent directions, including all Levi and similitude directions, have been included in this calculation.

\paragraph{Explicit matrices for both full conormals}\label{m6:s:p78:s11-2}

At the two representatives, \eqref{m6:s:eq:p78-3-4} and \eqref{m6:s:eq:p78-11-5} give

\begin{equation}
CN^G_{C_7,W_7}=F\mathsf X_{21}\oplus F\mathsf Y_{21},
\label{m6:s:eq:p78-11-6}
\end{equation}

and

\begin{equation}
CN^G_{C_8,W_8}
=F\mathsf X_{21}\oplus F\mathsf Y_{21}
 \oplus F\mathsf Y_{22}\oplus F\mathsf T_f.
\label{m6:s:eq:p78-11-7}
\end{equation}

In the root coordinates \eqref{m6:s:eq:p78-1-3}, these whole
matrix spaces are equivalently
\begin{equation}\label{m6:s:eq:p78-11-8}
\mathcal N(aE_{21},bE_{21},0,0,0,0),\qquad a,b\in F,
\end{equation}
and, putting $Y=\left(\begin{smallmatrix}0&0\\b&c\end{smallmatrix}\right)$,
\begin{equation}\label{m6:s:eq:p78-11-9}
\mathcal N(aE_{21},Y,0,0,dE_{21},0),\qquad a,b,c,d\in F.
\end{equation}
The expanded block arrays are in
\texttt{anc/}\allowbreak\texttt{model6\_calculations.pdf}, [M6-CERT-INCIDENCE-TANGENTS].

The displayed matrices are independent, symplectic, and annihilate
$\mathfrak p$; their $W_i^{-1}$-conjugates lie in $\mathfrak u^t$,
annihilating $\operatorname{Ad}(W_i)\mathfrak p^t$.
The ambient codimensions $2,4$ prove completeness at every point,
including the central direction.

In fact every matrix in \eqref{m6:s:eq:p78-11-8} or \eqref{m6:s:eq:p78-11-9} has square zero. For \eqref{m6:s:eq:p78-11-9}, the only potentially nonzero products of successive blocks are killed by

\begin{equation}
E_{21}Y=0,\quad E_{21}^2=0,\quad
Y(-Y^\sharp)=-\det(Y)I_2=0,\quad
(-Y^\sharp)E_{21}=0.
\label{m6:s:eq:p78-11-10}
\end{equation}

The assertion for \eqref{m6:s:eq:p78-11-8} is its simpler specialization. This proves, in particular, that these are entirely nilpotent conormal spaces.

At an arbitrary rational point \(g=u_0l W_i r^t v_0^t,\)

both unipotent radicals are preserved by their respective parabolics. Formula \eqref{m6:s:eq:p78-11-5} becomes

\begin{equation}
{CN^G_{C_i,g}=\operatorname{Ad}(u_0l)\operatorname{Lie}K_i.}
\label{m6:s:eq:p78-11-11}
\end{equation}

Conjugation preserves the square-zero property. Thus the entire ambient conormal is nilpotent at every point of either entire parabolic cell.

\paragraph{The tangent check for the common incidence submersion}\label{m6:s:p78:s11-3}

The additional middle symplectic directions in Section~\ref{m6:s:p78:s2} are genuinely tangent directions of the single group \(G_*\). At a representative, one may vary the middle factor as \(\exp(tZ)W_i\), with \(Z\in\mathfrak g_*\), because \(w_o\) commutes with \(G_*\). Right translation identifies the new tangent direction with \(Z\).

For \(\mathcal P_7\), the trace-pairing matrix of the conormal basis \((\mathsf X_{21},\mathsf Y_{21})\) with the added direction \(\mathsf Y_{21}^t\) is \((0,2)^t\). For \(\mathcal P_8\), pairing the basis in \eqref{m6:s:eq:p78-11-7} with the three added directions \((\mathsf Y_{21}^t,\mathsf Y_{22}^t,\mathsf T_f^t)\) gives

\begin{equation}
\begin{pmatrix}
0&0&0\\
2&0&0\\
0&2&0\\
0&0&1
\end{pmatrix}.
\label{m6:s:eq:p78-11-12}
\end{equation}

The annihilator is exactly $F\mathsf X_{21}$, so the rank is $55$,
including every middle tangent; the factors $2$ are nonzero also
dyadically. At either representative $g_{55}=E_{22}$ and
$\mathsf X_{21}^tW_i$ changes its $(1,1)$ entry by $-1$, giving
derivative $-1$ for $\det g_{55}$. Parabolic translation and the
established first incidence locus prove smoothness everywhere.

\paragraph{Model-orbit tangents and compatible dimensions, kept separate}\label{m6:s:p78:s11-4}

The full model-orbit tangent is

\begin{equation}
T_g(\mathcal A g)=\{Ag+gB^t:A,B\in\mathfrak h\}.
\label{m6:s:eq:p78-11-13}
\end{equation}

Its dimension is \(52-\dim\mathcal A_g\). The full stabilizers already solved give the following exact table. ``Locus dimension'' includes the unipotent fibers and all listed invariant parameters.

\begin{center}\small\setlength{\tabcolsep}{3pt}
\begin{tabular}{>{\raggedright\arraybackslash}p{0.44\linewidth}>{\raggedright\arraybackslash}p{0.16\linewidth}>{\raggedright\arraybackslash}p{0.16\linewidth}>{\raggedright\arraybackslash}p{0.16\linewidth}}
\hline
Cell and locus & Full stabilizer dimension & Model-orbit dimension & Locus dimension \\ \hline
\(\mathcal P_7\), \(q\ne0,p\ne0\), away from the plane origin & 3 & 49 & 52 for the whole open branch \\
\(\mathcal P_7\), plane origin \(b=e=0\), with \(p,k\) variable & 4 & 48 & 50 \\
\(\mathcal P_7\), \(q\ne0,p=0\) & 3 & 49 & 51 \\
\(\mathcal P_7\), \(q=0,\wp\ne0\) & 3 & 49 & 51 \\
\(\mathcal P_7\), \(q=\wp=0\) & 5 & 47 & 50 \\
\(\mathcal P_8\), \(p\ne0\) & 5 & 47 & 49 \\
\(\mathcal P_8\), \(p=0\), before the full-character equation & 6 & 46 & 48 \\
\(\mathcal P_8\), \(p=0,a=b\) & 6 & 46 & 47 \\
\hline
\end{tabular}
\end{center}

For example, the compatible base dimensions are \(24-12-2=10\) and \(24-12-3=9\), from the complete Levi quotient and the two or three independent compatibility equations. Restoring the affine fibers gives 52 and 49. On \(\mathcal P_8\) the independent invariant coordinates \((k,m)\) supply the two directions transverse to the 47-dimensional orbit. On its boundary \((a,b)\) supply two directions transverse to the 46-dimensional orbit; after imposing \(a=b\), exactly one remains. The rational normalizations prove these dimension assertions also at \(m=0\), repeated characteristic roots, and \(a=b\).

The table records orbit tangent dimensions and the dimensions of their saturated compatible families.

\paragraph{Intrinsic vanishing and the coarse conormal}\label{m6:s:p78:s11-5}

The entire-cell conormals are square-zero. Both cell vanishings follow
from intrinsic character and transpose arguments;
Proposition~\ref{m6:s:prop:p006} supplies their predecessor open.


\subsection{P031--P033: incidence cone, submersion and compatibility jump}\label{m6:s:cert:three-cell-cone}

\subsubsection{Polynomial matrix and exhaustive rank equations}

Put

\[
Y=\begin{pmatrix}a&b\\c&d\end{pmatrix},\quad T=\begin{pmatrix}h&e\\f&-h\end{pmatrix},\qquad
g(Y,T)=\exp\bigl(\mathcal N(0,Y,0,0,T,0)^t\bigr)W_{33}.
\]

Multiplication in the fixed ten-dimensional form gives the following entire middle six-by-six matrix, on physical indices 3 through 8:

\[
\begin{pmatrix}
0&0&0&0&0&-1\\
0&0&0&0&1&0\\
0&0&1&0&c&-a\\
0&0&0&1&d&-b\\
0&-1&-d&c&f&(ad-bc)/2-h\\
1&0&b&-a&-(ad-bc)/2-h&-e
\end{pmatrix}.
\]

The outer entries are those of \(W_{33}\). Define

\[
C=Y^tj=\begin{pmatrix}c&-a\\d&-b\end{pmatrix},\qquad
L=-\operatorname{adj}(Y)^t=\begin{pmatrix}-d&c\\b&-a\end{pmatrix},
\]
\[
D=\begin{pmatrix}f&\det(Y)/2-h\\-\det(Y)/2-h&-e\end{pmatrix}.
\]

Multiplication gives \(LC=-\det(Y)j\). The incidence cone is

\begin{equation}
\label{m6:s:eq:degeneration-vii-1}
\mathcal Z=\left\{(Y,D):\det Y=0,\ D=D^t,
\ \operatorname{rank}\binom C D\leq1,\quad \operatorname{rank}(L\ D)\leq1\right\}.
\end{equation}

This is an explicit polynomial description: the equations are \(ad-bc=0\), \(D_{12}-D_{21}=0\), and all two-by-two minors of the displayed four-by-two and two-by-four matrices. For \(D=\left(\begin{smallmatrix}d_{11}&d_{12}\\d_{12}&d_{22}\end{smallmatrix}\right)\), the corresponding \(T\) is uniquely \(\left(\begin{smallmatrix}-d_{12}&-d_{22}\\d_{11}&d_{12}\end{smallmatrix}\right)\).

There are precisely three rank strata. A rank-two \(D\) contradicts either incidence condition. If \(D=0\), \(Y\) has rank zero or one. The southeast four-by-four middle block has rank

\[
\operatorname{rank}\begin{pmatrix}I_2&C\\L&D\end{pmatrix}
=2+\operatorname{rank}(D-LC)=2+\operatorname{rank}D
\]

on the cone; the two neighboring rectangular ranks are the ranks in \eqref{m6:s:eq:degeneration-vii-1}. The outer pivot contributes one to the relevant full corner ranks. Thus the three strata and their full matrices are

\[
\begin{array}{c|c}
\operatorname{rank}D=1&
R_{31}=\begin{pmatrix}7&5&3&1\\5&4&2&1\\3&2&2&1\\1&1&1&1\end{pmatrix}\\[3pt]
D=0,\ \operatorname{rank}Y=1&
R_{32}=\begin{pmatrix}7&5&3&1\\5&3&2&1\\3&2&1&1\\1&1&1&1\end{pmatrix}\\[3pt]
Y=D=0&
R_{33}=\begin{pmatrix}7&5&3&1\\5&3&1&1\\3&1&1&1\\1&1&1&1\end{pmatrix}.
\end{array}
\]

Conversely, within this seven-coordinate slice the corner-rank requirements of any of these three cells force the two rectangular ranks to be at most one. In particular \(\operatorname{rank}C\leq1\), so \(\det Y=0\), and the displayed formula for \(D\) then makes it symmetric. The rank conditions therefore imply exactly \eqref{m6:s:eq:degeneration-vii-1}. This proves exhaustion inside the slice.

\subsubsection{Rational coverage, dimensions and local geometry}

If \(Y\neq0\), its rank is one and it admits a factorization \(Y=u v^t\) with \(u,v\in F^2\setminus\{0\}\), obtained by dividing by any nonzero entry. The row space of \(C\) is \(F(ju)^t\), and the column space of \(L\) is \(Fju\). The incidence equations imply

\begin{equation}
\label{m6:s:eq:degeneration-vii-2}
D=s(ju)(ju)^t,\qquad s\in F.
\end{equation}

Conversely \eqref{m6:s:eq:degeneration-vii-2} satisfies all the equations. The redundancy is precisely \((u,v,s)\sim(zu,z^{-1}v,z^{-2}s),\quad z\in F^\times.\)

If \(Y=0\), a nonzero allowed \(D\) is symmetric of rank one. Such a matrix has a nonzero diagonal entry: otherwise its nonzero off-diagonal entry would give nonzero determinant. If \(D_{11}\neq0\), for example,

\[
D=D_{11}\begin{pmatrix}1\\D_{12}/D_{11}\end{pmatrix}
\begin{pmatrix}1&D_{12}/D_{11}\end{pmatrix}.
\]

There is an analogous second chart. The scalar \(D_{11}\) ranges over \(F^\times\), and the two charts cover every rational point.

On \(\operatorname{rank}D=1\), the symmetric rank-one \(D\) varies in a smooth two-dimensional space. Its image line determines the allowed image line of \(Y\), and the two column coefficients of \(Y\) vary freely. This stratum is smooth of dimension four. When \(D=0\) and \(Y\) has rank one, it is the usual smooth three-dimensional rank-one matrix variety. The last stratum is a point. The \(Y\neq0\) parametrization has dimension \(2+2+1-1=4\); its closure contains every \(Y=0,D\neq0\) point by scaling a permitted nonzero \(Y\) to zero. The cone is irreducible of dimension four. In the seven-coordinate ambient space with \(D\) already symmetric, all remaining equations are quadratic, giving tangent dimension seven at the origin; the origin is singular.

Each stratum contains rational points of its indicated cell. The established rational \(P\)-double-cell classification, applied to these three ranks, implies that its \(P(F)^2\)-saturation is exactly that cell. This determines the global rational saturation of each stratum.

There is also a local analytic chart with this incidence cone. Let \(\mathsf{Q}_{\mathrm{out}}\) preserve the first isotropic two-space. Its Lie algebra is \(\mathfrak q_{\mathrm{out}}=\mathfrak p\oplus\mathfrak v_7, \qquad \mathfrak v_7=\langle Y_{11}^t,Y_{12}^t,Y_{21}^t,Y_{22}^t,T_h^t,T_e^t,T_f^t\rangle.\)

The tangent spaces at \(W_{33}\) give \(\dim T(P W_{33}P^t)=42,\qquad \dim\bigl(T(P W_{33}P^t)+\mathfrak v_7W_{33}\bigr)=49,\)
\[
\dim T(\mathsf{Q}_{\mathrm{out}}W_{33}\mathsf{Q}_{\mathrm{out}}^t)=49.
\]

Direct conjugation also gives \(\operatorname{Ad}(W_{33})(\mathfrak v_7^t)=\mathfrak v_7\). Selecting a 42-dimensional complement to the parabolic stabilizer in the source, the action map with these seven additional coordinates has invertible differential onto the outer orbit. The non-Archimedean inverse function theorem gives a local analytic chart. Since the parabolic factors preserve corner ranks, its three-cell part is a sufficiently small neighborhood in the smooth 42-dimensional factor times a sufficiently small neighborhood of zero in \eqref{m6:s:eq:degeneration-vii-1}. The degeneration curves below meet the other strata arbitrarily near the origin; parabolic transitivity transports such charts to their rational points. This gives the required local analytic model.

\subsubsection{degenerations and the compatibility jump}

Substitution gives the following polynomial curves and rank behavior:

\begin{center}
\begin{tabular}{@{}p{0.307\linewidth} p{0.307\linewidth} p{0.307\linewidth}@{}}
\hline
Curve & Cell for \(t\neq0\) & Cell at \(t=0\) \\
\(\exp(tY_{11}^t)W_{33}\) & \(\mathcal P_{32}\) & \(\mathcal P_{33}\) \\
\(\exp(tT_f^t)W_{33}\) & \(\mathcal P_{31}\) & \(\mathcal P_{33}\) \\
\(\exp((Y_{11}+tT_e)^t)W_{33}\) & \(\mathcal P_{31}\) & \(\mathcal P_{32}\) \\
\hline
\end{tabular}
\end{center}

They establish the stated algebraic and \(F\)-analytic incidences.

For the compatible \(\mathcal P_{33}\) limit, put \(D_0=\operatorname{diag}(1,-1)\), \(K_0=-j\), \(x_0=W_{33}\mathsf l(D_0,K_0)^t\). The curve \(x(t)=\exp(tT_e^t)x_0\) is polynomial at zero. For \(t\neq0\), matrix multiplication verifies

\begin{equation}
\label{m6:s:eq:degeneration-vii-3}
\exp(tT_e^t)W_{33}
=(I+t^{-1}E_{3,8})m(I,\operatorname{diag}(-t^{-1},1),I)
W_{31}(I-t^{-1}E_{3,8}^t).
\end{equation}

Its coefficient matrices are \(\mathsf{P}=\operatorname{diag}(-t,1),\quad \mathsf{Q}=\operatorname{diag}(-t^{-1},1),\quad \mathsf{R}=D_0,\quad \mathsf{S}=K_0.\)

All \(\mathcal P_{31}\) compatibility equations hold and \(\rho=-t\). The full \(\mathcal P_{31}\) stabilizer calculation establishes full compatibility. Hence \(x_0\) is a compatible limit, despite the denominators in this auxiliary cell factorization.

For the \(\mathcal P_{32}\) limit, fix any \(\lambda\in F^\times\), put \(M=\left(\begin{smallmatrix}0&1\\1&0\end{smallmatrix}\right)\), and use

\[
w(t)=\exp(tY_{21}^t)W_{32},\quad \mathsf{P}=M,\quad \mathsf{R}=\begin{pmatrix}0&-1\\\lambda&0\end{pmatrix},
\]
\begin{equation}
\label{m6:s:eq:degeneration-vii-4}
\mathsf{Q}(t)=\begin{pmatrix}-t&1\\1&0\end{pmatrix},\quad
\mathsf{S}(t)=\begin{pmatrix}t&1\\1&0\end{pmatrix},\quad
z_\lambda(t)=\mathsf l(\mathsf{P},\mathsf{Q}(t))w(t)\mathsf l(\mathsf{R},\mathsf{S}(t))^t.
\end{equation}

All determinant conditions hold for every \(t\), and the curve is regular. Independently solving the original 44-variable unipotent stabilizer equation at the moving representative gives dimension ten for \(t\neq0\). On a full kernel basis its possibly nonzero character values are \(p_{12}+r_{12},\quad r_{22},\quad -tq_{21}+s_{11},\quad p_{22},\quad tp_{22},\quad -q_{11}-ts_{21}.\)

These are zero for \eqref{m6:s:eq:degeneration-vii-4}. The other four values are zero identically. This verifies compatibility along the moving matrix.

The following factorization identifies the invariant on this curve. Put
\[
\begin{aligned}
q_0&=m(I,M,I),\\
l_0&=m(I,\operatorname{diag}(1,-t^{-1}),\operatorname{diag}(t,t^{-1})),\\
u&=\exp(t^{-1}Y_{21}).
\end{aligned}
\]

one has \(\exp(tY_{21}^t)W_{32}=u l_0q_0W_{31}q_0u^t\). Taking Levi projections in \eqref{m6:s:eq:degeneration-vii-4} gives

\[
\mathsf{P}'=\operatorname{diag}(-t,1),\quad \mathsf{Q}'=\begin{pmatrix}-t^{-2}&-1\\0&t\end{pmatrix},\quad
\mathsf{R}'=\operatorname{diag}(\lambda,-1),\quad \mathsf{S}'=\begin{pmatrix}1&t\\0&1\end{pmatrix}.
\]

They obey the five \(\mathcal P_{31}\) equations and give \(\rho=-\lambda t\). For \(t\neq0\), full compatibility follows from the entire \(\mathcal P_{31}\) family. At zero, \eqref{m6:s:eq:degeneration-vii-4} is the \(\mathcal P_{32}\) normal form with this arbitrary \(\lambda\), whose full character is nontrivial. Translating the curve by \(\mathcal A(F)\) elements covers every unipotently compatible \(\mathcal P_{32}\) point as a limit. The stabilizer gains new character-bearing directions at that boundary; full compatibility really is not closed across this jump.

The boundary \(\mathcal P_{32}\) retains the parameter \(\lambda\in F^\times\), whereas \(\mathcal P_{33}\) has one orbit.

\subsection{P040, P041 and P053: degeneration curves and minors}\label{m6:s:cert:outer-degenerations}

Write physical matrix units with ten-dimensional indices in this section. For the signed \(W_{41},W_{53}\), define \(g_1(t)=(I+tE_{74})W_{41},\qquad g_2(t)=(I-t(E_{61}+E_{10,5}))W_{53}.\)
Both root matrices square to zero and satisfy \(Z^tJ+JZ=0\); hence both curves lie in \(\operatorname{Sp}(J)\) for every \(t\). Their full rank matrices, computed, are \(\mathsf{R}(g_1(t))=R_{40},\quad \mathsf{R}(g_1(0))=R_{41},\qquad \mathsf{R}(g_2(t))=R_{41},\quad \mathsf{R}(g_2(0))=R_{53}\quad(t\ne0).\)
The corresponding torus matrices are \(D_1(t)=\operatorname{diag}(1,1,1,1,1,t,t,t,t,t)\) and \(D_2(t)=\operatorname{diag}(1,1,1,1,t^{-1},t,1,1,1,1)\). They satisfy the Levi identities: the first has multiplier \(t\) and scales \(E_{74}\) by \(t\); the second has multiplier one and scales both \(E_{61}\) and \(E_{10,5}\) by \(t\). Rational parabolic translation then puts every lower-cell point in the closure of the higher cell, for the analytic and Zariski topologies.

For the stronger compatible degeneration, set \(x_{53}=(9,3,-7,1,5,6,-10,-4,-8,2),\quad Z=E_{52}+E_{96}-E_{74},\quad b(t)=(I+tZ)x_{53}.\)
Again \(Z^2=0\) and \(Z^tJ+JZ=0\). Multiplication gives exactly \(b(t)=x_{53}+t(E_{5,10}+E_{7,8}+E_{9,6}),\qquad b(0)=x_{53}.\)
Its ranks are \(R_{40}\) for every nonzero \(t\), and \(R_{53}\) at zero.

Now return to the nilpotent-root names. Put
\[
 N_L=t^{-1}(T_f-Z_{21}),\quad N_R=t^{-1}(T_e-Z_{12}),\quad
 a_t=\begin{pmatrix}0&-t\\t^{-1}&0\end{pmatrix},
\]
\[
                 h_1=d_5(a_t)(I+N_L),\qquad h_2=I+N_R.
\]
Independent exact ten-by-ten multiplication verifies
\begin{equation}
\label{m6:s:eq:late-p40-10}
 N_L^2=N_R^2=0,\quad h_i^tJh_i=J,\qquad
                    h_1b(t)h_2^t=x_{1/t}.                
\end{equation}
The nilpotent terms have zero adjacent trace, so both characters are exactly one. This is a rational transporter on all \(t\ne0\), with denominators only powers of \(t\). It proves \(b(t)\in C_{40}\) and \(\kappa(b(t))=t^{-1}\). The closure contains \(x_{53}\), and therefore its full \(H(F)\times H(F)\)-orbit. The rational and geometric transitivity proved for \(\mathcal P_{53}\) identify this orbit with all of \(C_{53}\). Therefore this gives \(C_{53}\subset\overline{C_{40}}^{\,G(F)}\)
and the algebraic closure inclusion.

For explicit rank certificates, the following minors use physical one-based rows and columns. The displayed ranks are the generic polynomial ranks, so larger minors vanish identically. Each listed full-rank minor is a nonzero monomial in \(t\), proving the same ranks at every \(t\ne0\). The first two rows of the rank matrix determine the other two by the symplectic complementary-rank identities; all sixteen were also directly recalculated.

\begin{longtable}{@{}p{0.07\linewidth} p{0.08\linewidth} p{0.07\linewidth} p{0.23\linewidth} p{0.24\linewidth} p{0.09\linewidth}@{}}
\hline
Curve & Corner & Rank & Rows & Columns & Minor \\
\(b\) & 11 & 6 & 4,5,6,7,8,10 & 3,5,6,7,8,9 & -1 \\
\(b\) & 12 & 5 & 4,5,6,8,10 & 5,6,7,8,9 & -1 \\
\(b\) & 13 & 4 & 4,5,8,10 & 7,8,9,10 & -t \\
\(b\) & 14 & 2 & 5,8 & 9,10 & t \\
\(b\) & 21 & 5 & 5,6,7,8,10 & 3,5,6,7,9 & 1 \\
\(b\) & 22 & 5 & 5,6,7,8,10 & 5,6,7,8,9 & t \\
\(b\) & 23 & 4 & 5,7,8,10 & 7,8,9,10 & -\(t^2\) \\
\(b\) & 24 & 2 & 5,8 & 9,10 & t \\
\(g_1\) & 11 & 6 & 4,6,7,8,9,10 & 4,6,7,8,9,10 & 1 \\
\(g_1\) & 12 & 5 & 4,6,8,9,10 & 6,7,8,9,10 & -1 \\
\(g_1\) & 13 & 4 & 4,6,8,10 & 7,8,9,10 & -1 \\
\(g_1\) & 14 & 2 & 6,8 & 9,10 & -1 \\
\(g_1\) & 21 & 5 & 6,7,8,9,10 & 4,6,8,9,10 & 1 \\
\(g_1\) & 22 & 5 & 6,7,8,9,10 & 6,7,8,9,10 & t \\
\(g_1\) & 23 & 4 & 6,7,8,10 & 7,8,9,10 & t \\
\(g_1\) & 24 & 2 & 6,8 & 9,10 & -1 \\
\(g_2\) & 11 & 6 & 4,5,6,7,8,9 & 4,5,6,7,8,9 & -1 \\
\(g_2\) & 12 & 5 & 4,5,6,8,9 & 5,6,7,8,9 & -1 \\
\(g_2\) & 13 & 4 & 4,6,8,9 & 7,8,9,10 & -t \\
\(g_2\) & 14 & 2 & 6,8 & 9,10 & t \\
\(g_2\) & 21 & 5 & 5,6,7,8,9 & 4,5,6,8,9 & 1 \\
\(g_2\) & 22 & 4 & 5,6,8,9 & 5,6,8,9 & -1 \\
\(g_2\) & 23 & 3 & 6,8,9 & 8,9,10 & -t \\
\(g_2\) & 24 & 2 & 6,8 & 9,10 & t \\
\hline
\end{longtable}

In the matrix presentation, \(M=tE_{11}\to0\). Its conjugate
\[
        \begin{pmatrix}t&1\\0&0\end{pmatrix}\longrightarrow E_{12}
\]
has a nonzero nilpotent limit; the conjugator \(\left(\begin{smallmatrix}1&-1/t\\0&1\end{smallmatrix}\right)\), \eqref{m6:s:eq:late-p40-9}, and \eqref{m6:s:eq:late-p40-10} are singular at \(t=0\). The normalized section \(\sigma\) has parameter domain \(\operatorname{tr}M\ne0\); the displayed curves establish the compatible closure relation in the original ambient group.

\subsection{Rank-zero invariant fibers and the Jordan boundary}
\label{m6:s:p034:s8}\label{m6:s:app:p034-fibers}
\begin{proposition}\label{m6:s:prop:p034-fibers}
Every conjugation-invariant ordinary distribution on $GL_2(F)$ is
transpose-even. A distinct-root characteristic-polynomial fiber has
one-dimensional invariant distribution space. On a double-root fiber
with root $\lambda\ne0$, this space is
$\mathbb C\mu_\lambda\oplus\mathbb C\delta_{\lambda I}$,
where $\mu_\lambda$ is the locally finite Jordan orbital extension.
\end{proposition}
\begin{proof}
A distinct-root fiber is a single closed rational conjugacy orbit,
with split or nonsplit unimodular torus centralizer. Its invariant
line is generated by a positive orbital measure, which transpose fixes
by the transporter identities
\eqref{m6:s:eq:p034-7-1}--\eqref{m6:s:eq:p034-7-3}.

For a double root, translate by $\lambda I$ and write
\[
 N=\begin{pmatrix}z&r\\s&-z\end{pmatrix},\qquad q(N)=z^2+rs.
\]
The charts $r\ne0$ and $s\ne0$ cover the nonzero nilpotent orbit.
Its invariant residue density is
\[
 \frac{|dz\,dr|}{|r|},\qquad \frac{|dz\,ds|}{|s|}.
\]
They agree on the overlap since $s=-z^2/r$ has absolute Jacobian
$|s/r|$ with $z$ fixed. Conjugation preserves $q$ and has determinant
one on $\mathfrak{sl}_2$, so preserves this density. Under scalar
dilation it satisfies $\mu_{\rm nil}(cE)=|c|\mu_{\rm nil}(E)$.
For $\Sigma=\{q=0:\max(|z|,|r|,|s|)=1\}$ and a uniformizer $\varpi$,
\[
 \mu_{\rm nil}(0<\|N\|\le1)
 =\mu_{\rm nil}(\Sigma)\sum_{k\ge0}|\varpi|^k<\infty.
\]
Thus assigning zero mass to the origin gives a locally finite
invariant measure on the whole cone. If additive Haar measure gives
$\mathcal O_F$ volume one and the residue cardinality is $q_F$, the
disjoint unit charts $|r|=1,|z|\le1$ and $|s|=1,|z|<1$ have masses
$1-q_F^{-1}$ and $q_F^{-1}(1-q_F^{-1})$. Hence
$\mu_{\rm nil}(\Sigma)=1-q_F^{-2}$ and the displayed ball has mass
$1+q_F^{-1}$. The same charts and mass formula apply
in residue characteristic two. Transpose exchanges the charts
and preserves their density; the extended measure is even.

Translate this measure by $\lambda I$ to obtain $\mu_\lambda$.
Any invariant distribution on the double-root fiber restricts on its
open Jordan orbit to $c\mu_\lambda$. Subtracting this extension
leaves a distribution supported on the scalar point. The ordinary
closed-support sequence gives
\[
 D'(\chi^{-1}(2\lambda,\lambda^2))^{GL_2}
 =\mathbb C\mu_\lambda\oplus\mathbb C\delta_{\lambda I}.
\]
The two lines are independent by restriction to the Jordan orbit,
and both are even. Formula~\eqref{m6:s:eq:p034-5-13} carries this
extension to $X_0$ by product with the smooth character-weighted
$M_2$ measure; its unitary phase preserves local finiteness.

Finally $D'(GL_2)^{\operatorname{Ad}GL_2}$ is a weak-* closed module
under functions pulled back by $\chi:GL_2\to F\times F^\times$.
Bernstein localization \cite[Theorem 2.8]{JLX} topologically generates
it by the subspaces supported on closed fibers. These are the even
lines and two-dimensional spaces just described. Continuity of
transpose implies global evenness. This is topological generation
by fiber-supported distributions. Combining with
\eqref{m6:s:eq:p034-5-16} proves Proposition~\ref{m6:s:prop:p034}.
\end{proof}

\section{Model 6: rational orbits and boundary calculations}\label{supp:m6-orbital-families}
\subsection{Independent Shalika classifications}\label{supp:m6-sharpness}
\subsubsection{Shalika coordinates and the rank dictionary}
\label{supp:m6-shalika-coordinates}
We continue the Gauss decomposition of Proposition~\ref{m6:rank-dictionary}.
The following external-factor descent is used by the independent
rank-one and rank-zero calculations.

Under the combined chart below, the cells
$\mathcal P_1,\mathcal P_6,\mathcal P_{34}$ correspond respectively to
\[
 X_r=\left\{\begin{pmatrix}A&B\\C&d\end{pmatrix}\in\GL_4:
             \rank d=r\right\},\qquad r=2,1,0.
\]
The correspondence retains the stabilizer, character, transpose, and
trivial residual density. It commutes with open restrictions and
closed-support injections. These assertions follow from the
averaging and rank calculation below.

\paragraph{\texorpdfstring{External \(GL_2\) Frobenius descent}{External GL2 Frobenius descent}}\label{m6:s:p006:s2-2}

The external coordinate transforms by \(C\mapsto a_1Ca_2^t\). Its base is a single rational orbit, with base point \(I_2\), section
\(z\mapsto((d(z),z),1)\), and stabilizer condition \(a_2=a_1^{-t}\).

The global coordinate map

\begin{equation}
GL_2\times GL_4\longrightarrow GL_4\times GL_2,
\qquad (z,x)\longmapsto(d(z)x,z)
\label{m6:s:eq:p006-2-7}
\end{equation}

is a rational isomorphism, with inverse \((A,C)\mapsto(C,d(C)^{-1}A)\). The Frobenius averaging is

\begin{equation}
\mathcal I f(x)=\int_{GL_2}f(d(z)x,z)\,dz.
\label{m6:s:eq:p006-2-8}
\end{equation}

Combining \eqref{m6:s:eq:p006-2-2} and \eqref{m6:s:eq:p006-2-7} gives the common chart
\begin{equation}
\begin{aligned}
\Phi:N'\times GL_2\times GL_4\times N'&\longrightarrow\Omega',\\
(n_1,z,x,n_2)&\longmapsto n_1\iota(d(z)x,z)n_2^t.
\end{aligned}
\label{m6:s:eq:extraction-424-chart}
\end{equation}
Its validity on the whole big cell, including both boundary ranks, follows from block Gaussian inversion and the displayed inverse of \eqref{m6:s:eq:p006-2-7}.

The full residual group is

\[
K_S=(M_2\oplus M_2)\rtimes GL_2.
\]

To avoid repeated transposes in the calculations, write the lower right-unipotent parameter as \(V=X_2^t\), and put

\[
u(X)=\begin{pmatrix}I&X\\0&I\end{pmatrix},\qquad
l(V)=\begin{pmatrix}I&0\\V&I\end{pmatrix}.
\]

Then

\begin{equation}
{(X,V,a)\cdot x=u(X)d(a)x d(a)^{-1}l(V),\qquad
\eta_S(X,V,a)=\psi(\operatorname{tr}X+\operatorname{tr}V).}
\label{m6:s:eq:p006-2-9}
\end{equation}

Both matrix spaces in the semidirect product transform by \(Z\mapsto aZa^{-1}\). In these coordinates its multiplication is

\begin{equation}
(X,V,a)(X',V',a')=
(X+aX'a^{-1},\ V+aV'a^{-1},\ aa').
\label{m6:s:eq:p006-2-10}
\end{equation}

The transpose involution and its automorphism of \(K_S\) are

\begin{equation}
\tau(x)=x^t,\qquad
\sigma_\tau(X,V,a)=(V^t,X^t,a^{-t}).
\label{m6:s:eq:p006-2-11}
\end{equation}

The character is preserved. Extend it by \(\widetilde\eta_S(\tau)=-1\).

The determinants of conjugation on both copies of \(M_2\) are one. Hence \(H_S\), \(H_S\times H_S\), and \(K_S\) are unimodular. On the external base the stabilizer acts on tangent matrices by conjugation, also of determinant one. Transpose has absolute determinant one. Thus Frobenius descent has trivial relative density, including at the odd element.

The averaging and transpose are compatible on equivariant distributions: the conjugation introduced by the global coordinates \eqref{m6:s:eq:p006-2-7} is absorbed by equivariance.

\paragraph{Why the Shalika object is representation-valued}
Let $N_S=\{n(X):X\in M_2(F)\}\subset\GL_4$, with
$\psi_S(n(X))=\psi(\operatorname{tr}X)$. For a representation
$\sigma$ of $\GL_4$, its twisted Jacquet space
$J_{N_S,\psi_S}(\sigma)$ is the quotient by the span of
$\sigma(n)v-\psi_S(n)v$. The diagonal $d(\GL_2)$ normalizes
$(N_S,\psi_S)$, so this quotient is itself a $\GL_2$-representation.
For an independent representation $\rho$ of the outer $\GL_2$,
the graph subgroup $H_S$ gives, by this quotient definition,
\[
 \Hom_{H_S}(\sigma\boxtimes\rho,\xi_S)
 \simeq\Hom_{\GL_2}
       (J_{N_S,\psi_S}(\sigma)\otimes\rho,\mathbf1).
\]
The internal diagonal group acts on the Jacquet module, while the
external factor acts on $\rho$; the period pairs them diagonally.
Thus extraction retains a representation-valued Shalika object.
The ordinary-distribution proof here uses the action, density
and rank strata in this proposition.

\paragraph{The rank dictionary and its support maps}\label{m6:s:p006:s2-3}

Write

\begin{equation}
x=\begin{pmatrix}A&B\\C&d\end{pmatrix}\in GL_4.
\label{m6:s:eq:p006-2-12}
\end{equation}

In \eqref{m6:s:eq:p006-2-2}, the last \(4\times4\) block is unaffected by either radical factor. Its last \(2\times2\) block has the same rank as the upper-left \(2\times2\) block of \(A_{\mathrm{Levi}}^{-1}\). For any invertible \(4\times4\) block matrix, this rank equals the rank of its lower-right block.

Here is a direct nullity proof of the latter assertion. If \(M\in GL_4\), a vector \(v\) lies in the kernel of \((M^{-1})_{11}\) precisely when
\(M^{-1}(v,0)=(0,w)\). This is equivalent to
\(M_{22}w=0\) and \(v=M_{12}w\). Invertibility makes \(w\mapsto v\) an isomorphism of these two kernels. Both blocks are two-dimensional, so their ranks agree. Conjugation by \(j\), transpose, and multiplication by \(\det C\) do not change rank.

Finally, in \eqref{m6:s:eq:p006-2-7}, the lower-right block of \(A_{\mathrm{Levi}}=d(z)x\) is \(zd\). Therefore throughout the big cell

\begin{equation}
\operatorname{rank}g_{55}=\operatorname{rank}d.
\label{m6:s:eq:p006-2-13}
\end{equation}

Together with the established three-cell description of \(\Omega'\), this proves the exact identifications

\begin{equation}
\mathcal P_6\longleftrightarrow
X_1:=\left\{\begin{pmatrix}A&B\\C&d\end{pmatrix}\in GL_4:
\operatorname{rank}d=1\right\},
\label{m6:s:eq:p006-2-14}
\end{equation}

and

\begin{equation}
V_6:=\Omega'\cap\{g_{55}\ne0\}
\longleftrightarrow
\mathscr W:=\{x\in GL_4:d\ne0\}.
\label{m6:s:eq:p006-2-15}
\end{equation}

Both \(V_6\) and \(\mathscr W\) are opens in their respective reductive groups. In these opens, \(\mathcal P_6\) and \(X_1\) are closed, with respective open complements \(\mathcal P_1\) and \(\{\det d\ne0\}\).

The distribution correspondences restrict to

\begin{equation}
D'_{\widetilde\eta}(\mathcal P_6)
\simeq D'_{\widetilde\eta_S}(X_1),\qquad
D'_{\widetilde\eta}(V_6)
\simeq D'_{\widetilde\eta_S}(\mathscr W).
\label{m6:s:eq:p006-2-16}
\end{equation}

They commute with the indicated open restrictions and closed-support injections.

The same identity gives the remaining rank-zero equality
\begin{equation}
\Phi^{-1}(\mathcal P_{34})=N'\times GL_2\times X_0\times N',
\qquad X_0=\left\{\begin{pmatrix}A&B\\C&0\end{pmatrix}\in GL_4\right\}.
\label{m6:s:eq:p034-1-9}
\end{equation}
Thus the rank-zero blocks range over the entire invertible locus. The corresponding complements are
\begin{equation}
GL_4=\mathscr W\sqcup X_0,
\label{m6:s:eq:p034-1-10}
\end{equation}
\begin{equation}
V_6=\mathcal P_1\sqcup\mathcal P_6,
\qquad\Omega'\setminus V_6=\mathcal P_{34}.
\label{m6:s:eq:p034-1-11}
\end{equation}
Here \(X_0\) is closed in \(GL_4\), whereas \(X_1\) is closed only in \(\mathscr W\). Averaging therefore also gives
\begin{equation}
D'_{\widetilde\eta}(\mathcal P_{34})
\simeq D'_{\widetilde\eta_S}(X_0),
\label{m6:s:eq:p034-1-12}
\end{equation}
with the same assertion for the even character. Rank zero has a one-point lower-block base, hence contributes no new base-orbit density. These are intrinsic correspondences or maps of closed-support antecedents.

\paragraph{Full original stabilizers are retained}\label{m6:s:p006:s2-4}

At \(\iota(x,I_2)\), an element of the original \(H\times H\)-stabilizer must have trivial \(N'\times N'\) coordinates by uniqueness of \eqref{m6:s:eq:p006-2-2}. Its remaining coordinates are exactly the \(K_S\)-stabilizer of \(x\).

Explicitly, \((X,V,a)\) embeds as the pair of Levi-intersection elements corresponding to

\begin{equation}
(n(X)d(a),a),\qquad
(n(V^t)d(a^{-t}),a^{-t}).
\label{m6:s:eq:p006-2-17}
\end{equation}

Thus the calculations below determine the full stabilizers. Stabilizers at every other point of each of the three cells are conjugate to these through the displayed radical/external coordinate changes. This exhausts the original rational group stabilizer on both boundary ranks.

\subsubsection{Sharpness on the rank-one stratum}\label{m6:s:sec:p006}

The rank-one cell has intrinsic odd orbital distributions, whose
ambient conormals exclude their occurrence in a global NW distribution.
The calculation also supplies the spectral coordinates used on the
remaining boundary. Expanded algebra and Jacobians are in
\path{anc/model6_calculations.pdf}, section
\texttt{M6-P006}. Coordinates are those of
Proposition~\ref{m6:rank-dictionary}; the Gaussian and conormal
certificates are \S\ref{m6:s:cert:p006-gaussian},
\S\ref{m6:s:cert:p006-nonzero-normalization}, and
\S\ref{m6:s:cert:p006-conormal}--\ref{m6:s:cert:p006-ambient-pairing}.

\begin{proposition}\label{m6:s:prop:p006}
Let \(T\in\NW(G)\) be an \(\eta\)-equivariant, transpose-odd ordinary distribution on \(G\). If \(T|_\Omega=0\), then \(T|_{\mathcal O_6}=0\), where \(\mathcal O_6=\Omega\sqcup\mathcal P_6\). Nevertheless \(D'_{\widetilde\eta}(\mathcal P_6)\ne0\).
\end{proposition}

\paragraph{The entire rational rank-one base.}\label{m6:s:p006:s3}
The block projection $p_d:X_1\to\{\operatorname{rank}d=1\}$ is a
surjective submersion for conjugation: $\left(\begin{smallmatrix}0&I\\I&d\end{smallmatrix}\right)$
is invertible for every $d$. For $p=\operatorname{tr}d\ne0$ the base is
one rational orbit represented by $pE_{11}$, with full diagonal-torus
centralizer. For $p=0$ it is one rational nilpotent orbit represented by
$N=E_{12}$, with centralizer
\begin{equation}
 C_N=\{s(I+hN):s\in F^\times,h\in F\}.
 \label{m6:s:eq:p006-3-2}
\end{equation}
Local rational sections use an image vector and a kernel vector in the
first case, and
\begin{equation}
 a=(de_j,e_j),\qquad aNa^{-1}=d,\qquad de_j\ne0,
 \label{m6:s:eq:p006-3-3}
\end{equation}
in the second. They cover
all rational points.
Localize first by $p$ and descend to the fixed base. The even fiber
groups are $(M_2\oplus M_2)\rtimes T_2$ and
$(M_2\oplus M_2)\rtimes C_N$, respectively. They and their base
centralizers are unimodular: semisimple tangent weights are reciprocal,
and the nilpotent base action is a unipotent shear. All Frobenius
relative densities are one. Transpose fixes $pE_{11}$; over $N$ use
\begin{equation}
 w=\begin{pmatrix}0&1\\1&0\end{pmatrix},\qquad
 \theta(x)=d(w)x^td(w).
 \label{m6:s:eq:p006-3-6}
\end{equation}
It is a character-one translate of transpose, has square one and
absolute base Jacobian one, so the odd sign remains $-1$.

\paragraph{Gaussian extraction and the complete additive sieve.}\label{m6:s:p006:s4}
Use the canonical matrices
\begin{align}
 x_p(A,b,c)&=\begin{pmatrix}
 A_{11}&A_{12}&0&b_1\\A_{21}&A_{22}&0&b_2\\
 0&0&p&0\\c_1&c_2&0&0\end{pmatrix},
 \label{m6:s:eq:p006-4-1}\\
 x_N(A,b,c)&=\begin{pmatrix}
 A_{11}&A_{12}&b_1&0\\A_{21}&A_{22}&b_2&0\\
 0&0&0&1\\c_1&c_2&0&0\end{pmatrix}.
 \label{m6:s:eq:p006-4-2}
\end{align}
Invertibility is exactly that of
$\mathcal B(A,b,c)=\left(\begin{smallmatrix}A&b\\c&0\end{smallmatrix}\right)$,
so $b,c\ne0$. The global Gaussian charts
\eqref{m6:s:eq:p006-4-4}, \eqref{m6:s:eq:p006-4-7} have factors
$\psi(z_1+v_1)dz\,dv$, Jacobian $|p|^4$ in the semisimple case, and
$\psi(z_1+v_2)dz\,dv$, Jacobian one in the nilpotent case. Their
centralizer and transpose absolute determinants are one; corrected
transpose sends $(z,v)$ to $(wv^t,z^tw)$.
The remaining action is
\begin{equation}
 A\mapsto A+uc+bv,\quad
 \eta=\begin{cases}\psi(u_2+v_2),&d=pE_{11},\\
                         \psi(u_2+v_1),&d=N.
          \end{cases}
 \label{m6:s:eq:p006-4-10}
\end{equation}
The equation $uc+bv=0$ has exactly the solutions
\begin{equation}
 u=kb,\quad v=-kc\quad(k\in F),
 \label{m6:s:eq:p006-4-12}
\end{equation}
because its image line forces $u\in Fb$. Thus the image is
three-dimensional and full unipotent compatibility is
\begin{equation}
 b_2=c_2\quad(d=pE_{11}),\qquad b_2=c_1\quad(d=N).
 \label{m6:s:eq:p006-4-13}
\end{equation}
With $\ell_b=(-b_2,b_1)$, $k_c=(-c_2,c_1)^t$, the missing scalar is
\begin{equation}
 \det\mathcal B=-\ell_bAk_c\ne0.
 \label{m6:s:eq:p006-4-14}
\end{equation}
These statements solve the entire affine quotient, including degenerate
$A$.

\paragraph{Semisimple branch: the plane and its origin.}\label{m6:s:p006:s5}
\begin{proposition}\label{m6:s:prop:p006-semisimple}
Every intrinsic odd equivariant distribution on \(\{x\in X_1:\operatorname{tr}d\ne0\}\) is zero.
\end{proposition}
\begin{proof}
Fix $p\ne0$. A diagonal centralizer has effective parameter $h=s/t$
and acts by $(A,b,c)\mapsto(R_hAR_h^{-1},R_hb,cR_h^{-1})$,
$R_h=\operatorname{diag}(h,1)$. Its entire stabilizer is
\begin{equation}
 (h-1)b_1=(h-1)c_1=0,\quad
 uc+bv=A-R_hAR_h^{-1},\quad\eta=\psi(u_2+v_2).
 \label{m6:s:eq:p006-5-2}
\end{equation}
Gaussian preservation forces all extracted coordinates to zero.
If $(b_1,c_1)\ne(0,0)$ this is $\mathbb G_m\times\mathbb G_a$;
otherwise normalization makes $A$ a nonzero multiple of $E_{11}$ and
it is $T_2\times\mathbb G_a$. In both cases the stabilizer is
unimodular, and its character is trivial exactly on
$b_2=c_2=q$. Remove the incompatible open by its displayed subgroup.

\subparagraph{The branch $q\ne0$.}\label{m6:s:p006:s5-2}\label{m6:s:p006:s5-3}
Write $b=(r,q)^t$, $c=(s,q)$. The unique chart
\begin{equation}
 A=\alpha E_{11}+uc+bv,\quad u=(u_1,u_2)^t,\quad v=(v_1,0)
 \label{m6:s:eq:p006-5-3}
\end{equation}
has inverse $u_1=A_{12}/q$, $u_2=A_{22}/q$,
$v_1=A_{21}/q-sA_{22}/q^2$ and
$\alpha=A_{11}-sA_{12}/q-rA_{21}/q+rsA_{22}/q^2$.
It gives $\det x_p=-pq^2\alpha$, so $\alpha q\ne0$, and a
surjective product chart even at $r=s=0$. The three additive factors
are $\psi(u_2)du_1du_2dv_1$, with Jacobian $|q|^3$; together with
$|p|^4$ this depends only on invariant parameters. The quotient is
\begin{equation}
 h:(\alpha,q,r,s)\mapsto(\alpha,q,hr,h^{-1}s),\qquad
 \tau:(\alpha,q,r,s)\mapsto(\alpha,q,s,r).
 \label{m6:s:eq:p006-5-8}
\end{equation}
Here every $h\in F^\times$ occurs. On the additive coordinates $h$
acts as $(hu_1,u_2,h^{-1}v_1)$ and transpose as
$(v_1+(s/q)u_2,u_2,u_1-(r/q)u_2)$; their absolute determinants and
phase ratios are one. Hence this quotient is untwisted.
For $rs=c_0\ne0$, the fiber is one rational orbit parametrized by
$r\mapsto(r,c_0/r)$, and transpose preserves multiplicative Haar.
For $c_0=0$ retain both punctured axes and their common origin. Their
odd difference cannot extend invariantly across the origin:
Lemma~\ref{m8:s:c:plane} applies to the whole zero-product fiber.
Localization in $rs$, and then in $\alpha,q$, proves odd vanishing on
the entire $q\ne0$ branch.

\subparagraph{The branch $q=0$ and its density.}\label{m6:s:p006:s5-4}
Now $b=(r,0)^t$, $c=(s,0)$, $rs\ne0$, and the unique chart is
\begin{equation}
 A=\begin{pmatrix}su_1&rv_2\\su_2&\alpha\end{pmatrix},\qquad
 \det x_p=-p\alpha rs\ne0.
 \label{m6:s:eq:p006-5-14}
\end{equation}
The torus acts by $(hu_1,u_2,v_2)$ with determinant $|h|$.
The normalized fiber measure is therefore
\begin{equation}
 |s|\psi(u_2+v_2)du_1du_2dv_2.
 \label{m6:s:eq:p006-5-16}
\end{equation}
Transpose exchanges $r,s$ and sends the additive coordinates to
$((s/r)u_1,v_2,u_2)$, with determinant $|s/r|$; pushforward changes
$|s|$ to $|r|$, exactly the target normalization. The raw Jacobian
$|s|^2|r|$ differs only by invariant $|rs|$.
Thus the residual action is again untwisted interchange, and each
$rs$ fiber is a single rational torus orbit with character-one
transporter $\operatorname{diag}(s/r,1)$ and unimodular stabilizer.
This gives odd vanishing on $q=0$.
The order is incompatible open, compatible $q\ne0$ open, then its
$q=0$ closed complement. Frobenius descent after localization in $p$
proves the intrinsic proposition.
\end{proof}

\paragraph{Nilpotent branch with nonzero legs.}\label{m6:s:p006:s6}
\begin{proposition}\label{m6:s:prop:p006-nonzero-leg}
On the nilpotent rank-one branch, the full unipotent compatibility equation is \(b_2=c_1\). The open locus where this common value is nonzero has intrinsic odd vanishing.
\end{proposition}
\begin{proof}
Corrected transpose sends $(A,b,c)$ to $(wA^tw,wc^t,b^tw)$.
For $a_h=I+hN$ the Gaussian correction is
$X_h=h(a_hb)e_1^t$, $V_h=-he_2(ca_h^{-1})$, giving exactly
\begin{equation}
 u(X_h)d(a_h)x_N(A,b,c)d(a_h)^{-1}l(V_h)
 =x_N(a_h(A+h^2bc)a_h^{-1},a_hb,ca_h^{-1}).
 \label{m6:s:eq:p006-6-4}
\end{equation}
Its phase is $\psi(h(b_1-c_2)+h^2(b_2+c_1))$; the quadratic term
is retained. The full stabilizer master equations are
\begin{equation}
 \begin{gathered}
 a=s(I+hN),\quad hb_2=hc_1=0,\quad
 X=hbe_1^t+ue_2^t,\quad V=e_1v-he_2c,\\
 uc+bv=A-a_hAa_h^{-1}-h^2bc,\qquad
 \eta=\psi(u_2+v_1+h(b_1-c_2)).
 \end{gathered}
 \label{m6:s:eq:p006-6-6}
\end{equation}
When either leg is nonzero these force $h=0$, leaving the whole
stabilizer $\mathbb G_m\times\mathbb G_a$, with character
$\psi(k(b_2-c_1))$. On compatibility write
$b=(r,q)^t$, $c=(q,s)$, $q\ne0$. Taking
$h=(s-r)/(2q)$ in the corrected centralizer operation gives phase one
and legs $(m,q)^t,(q,m)$, $m=(r+s)/2$.
The explicit affine inverse in
\S\ref{m6:s:cert:p006-nonzero-normalization} sends $A$ to
$\kappa E_{12}$, with $\kappa=-\det x/q^2\ne0$.
Thus every rational orbit has the unique invariant triple
$(q,m,\kappa)$ and representative
\begin{equation}
 y_{q,m,\kappa}=x_N(\kappa E_{12},(m,q)^t,(q,m)),\quad
 \theta(y_{q,m,\kappa})=y_{q,m,\kappa}.
 \label{m6:s:eq:p006-6-14}
\end{equation}
Coverage follows from the full centralizer operation and the entire
three-dimensional affine quotient for arbitrary $A$. The stabilizer is unimodular with trivial character. Additive
Jacobian is $|q|^3$, and $dr\,ds=|2q|dm\,dh$; these positive
invariant factors create no character, including over dyadic $F$.
For any normalization $g$ to this representative, ordinary transpose
is reached by $t_w\sigma_\theta(g)^{-1}g$, whose character is one
because $\sigma_\theta$ preserves $\eta_S$. Positive quotient Haar
then gives transpose sign $+1$ on every orbit. Constructiveness gives
odd vanishing on the whole branch. Division by $2q$ is legal; no
square root or division by $m$ occurs.
\end{proof}

\paragraph{Nilpotent zero-leg branch: all rational orbits.}\label{m6:s:p006:s7}
\begin{proposition}\label{m6:s:prop:p006-zero-leg}
On the nilpotent branch with \(b_2=c_1=0\), the ordered parameters \((\alpha,\beta,\gamma)\in(F^\times)^3\) classify rational orbits. The full stabilizer is \(\mathbb G_m\times\mathbb G_a^2\), its character is trivial precisely when \((\beta-\gamma)(\alpha+\beta\gamma)=0\), and the locus \(\beta=\gamma\) has intrinsic odd vanishing.
\end{proposition}
\begin{proof}
Invertibility forces $b=\beta e_1$, $c=\gamma e_2^t$,
$\beta\gamma\ne0$. The translations
$u_1=-A_{12}/\gamma$, $u_2=-A_{22}/\gamma$,
$v_1=-A_{11}/\beta$, $v_2=0$ give
\begin{equation}
 x_{\alpha,\beta,\gamma}=
 \begin{pmatrix}0&0&\beta&0\\\alpha&0&0&0\\0&0&0&1\\0&\gamma&0&0\end{pmatrix},
 \quad\det x=-\alpha\beta\gamma\ne0.
 \label{m6:s:eq:p006-7-3}
\end{equation}
The corrected centralizer leaves $\beta,\gamma,A_{21}=\alpha$
unchanged; together with the displayed normalization this proves
complete rational coverage, including all continuous parameter values.
Solving the group equations, with $a=s(I+hN)$, gives the full stabilizer
\begin{equation}
 X=\begin{pmatrix}h\beta&-(\beta/\gamma)z+
              ((\alpha-\beta\gamma)/\gamma)h^2\\0&(\alpha/\gamma)h\end{pmatrix},
 \quad V=\begin{pmatrix}-(\alpha/\beta)h&z\\0&-h\gamma\end{pmatrix}.
 \label{m6:s:eq:p006-7-6}
\end{equation}
Indeed preservation of $B,C$ fixes the first column of $X$ and second
row of $V$; the remaining three independent equations determine the
other entries with free $z$. The noncentral law is
$(h,z)(h',z')=(h+h',z+z'+(\alpha/\beta-\gamma)hh')$;
$\zeta=z-\tfrac12(\alpha/\beta-\gamma)h^2$ makes it additive.
The full group is therefore connected and unimodular of dimension three.
Its entire phase is
\begin{equation}
 \operatorname{tr}X+\operatorname{tr}V
 =h(\beta-\gamma)(1+\alpha/(\beta\gamma)),\qquad
 (\beta-\gamma)(\alpha+\beta\gamma)=0
 \label{m6:s:eq:p006-7-11}
\end{equation}
for compatibility. If this fails, the genuine additive subgroup
$\zeta=0$ detects nontrivial character, so removes the entire open.

Transpose interchanges the ordered orbit invariants:
\begin{equation}
 \theta(x_{\alpha,\beta,\gamma})=x_{\alpha,\gamma,\beta}.
 \label{m6:s:eq:p006-7-13}
\end{equation}
On stabilizers it sends $(s,h,\zeta)$ to
$(s^{-1},-h,-(\beta/\gamma)\zeta)$; the comparison Jacobian is
$|\beta/\gamma|$, equal to one on a fixed orbit.
For completeness the additive quotient chart is
$A=\left(\begin{smallmatrix}\beta v_1&\gamma u_1\\\alpha&\gamma u_2\end{smallmatrix}\right)$
with normalized measure
\begin{equation}
 |\gamma|\psi(u_2+v_1)du_1du_2dv_1.
 \label{m6:s:eq:p006-7-17}
\end{equation}
Its raw Jacobian is $|\beta\gamma^2|$; corrected transpose sends the
coordinates to $((\gamma/\beta)u_1,v_1,u_2)$, pushing the coefficient
to $|\beta|$, exactly the target normalization. The residual interchange
has no density twist. On $\beta=\gamma$ the representative is fixed and
the ordinary-transpose transporter $t_w$ has character one, so positive
orbital Haar gives sign $+1$ and odd vanishing. On
$\alpha=-\beta\gamma$, $\beta\ne\gamma$, the two rational orbits are
exchanged and no such transporter exists.
\end{proof}

\paragraph{Intrinsic support and a whole-cell odd example.}\label{m6:s:p006:s8}
\begin{proposition}\label{m6:s:prop:p006-support}
Every intrinsic odd distribution on \(X_1\) is supported on the closed locus \(\mathscr S\) defined below. The space \(D'_{\widetilde\eta_S}(X_1)\) is nonzero.
\end{proposition}
\begin{proof}
\subparagraph{Intrinsic closed loci.}\label{m6:s:p006:s8-1}
In $\mathscr W=\{d\ne0\}$ define
\begin{equation}
 \begin{aligned}
 \mathscr N&=\{\det d=\operatorname{tr}d=0\},&
 q_B&=\operatorname{tr}(dB),\quad q_C=\operatorname{tr}(dC),\\
 \mathscr C&=\{x\in\mathscr N:q_B=q_C\},&
 \mathscr C_0&=\{x\in\mathscr N:q_B=q_C=0\}.
 \end{aligned}
 \label{m6:s:eq:p006-8-1}
\end{equation}
These are invariant closed loci; $d^2=0$ proves invariance of $q_B,q_C$
under both additive actions. On $\mathscr C_0$ there are unique regular
scalars $\beta,\gamma$ with $Bd=\beta d$, $dC=\gamma d$:
writing $d=u\phi$, $\phi(u)=0$, the trace-zero conditions force
$Bu\in Fu$ and $\phi C\in F\phi$. Entry ratios on $d_{ij}\ne0$
give regular functions agreeing on overlaps. Normalization proves
$\beta\gamma\ne0$, and
$\alpha=-\det x/(\beta\gamma)$ is regular, nonzero and invariant.
Transpose exchanges $\beta,\gamma$ and fixes $\alpha$.
Hence
\begin{equation}
 \mathscr E=\{\beta=\gamma\},\quad
 \mathscr S=\{\alpha=-\beta\gamma\}
           =\{\det x=(\beta\gamma)^2\}\subset\mathscr C_0
 \label{m6:s:eq:p006-8-4}
\end{equation}
are closed in $\mathscr W$; the scalar equation is imposed only where
the regular scalars exist.
The support order is semisimple open, then
$\mathscr N\setminus\mathscr C$, then
$\mathscr C\setminus\mathscr C_0$, then
$\mathscr C_0\setminus(\mathscr E\cup\mathscr S)$, and finally
$\mathscr E\setminus\mathscr S$ open in $\mathscr E\cup\mathscr S$.
The preceding vanishing results, applied to each ordinary closed-support
antecedent, leave
\begin{equation}
 \operatorname{supp}S\subset\mathscr S
       \quad(S\in D'_{\widetilde\eta_S}(X_1)).
 \label{m6:s:eq:p006-8-5}
\end{equation}
The residual support $\mathscr S$ is treated by ambient conormal exclusion below.

\subparagraph{Closed rational orbits.}\label{m6:s:p006:s8-3}
The regular map $x\mapsto(\alpha,\beta,\gamma)$ on $\mathscr C_0$
has one rational orbit per fiber, by the base pivot sections, Gaussian
chart and zero-leg translations. Every such orbit is closed in $X_1$
and $\mathscr W$. For $x_*=x_{-2,1,2}$, $\det x_*=4$, and its full
stabilizer phase is zero. Thus its equivariant positive orbital Haar
pushes forward to all of $X_1$: compact test support meets the closed
orbit in a compact set. Its transpose orbit has distinct ordered
invariants $(-2,2,1)$, so
\begin{equation}
 \mathcal T_{x_*}-\tau_*\mathcal T_{x_*}\ne0
             \quad\text{in }D'_{\widetilde\eta_S}(X_1).
 \label{m6:s:eq:p006-8-8}
\end{equation}
Common extraction carries this to the whole $\mathcal P_6$ cell,
with full radical/external product coordinates and full stabilizer.
The following ambient calculation gives a nonnilpotent conormal
direction on its closed-support extension.
The construction works over dyadic $F$, since $2\ne0$ in the field.
\end{proof}

\paragraph{The residual support is a smooth closed submanifold}\label{m6:s:p006:s9}

\begin{lemma}\label{m6:s:lem:p006-smooth}
The locus \(\mathscr S\) is a smooth closed submanifold of \(\mathscr W\), of dimension eleven, including along \(\beta=\gamma\).
\end{lemma}
\begin{proof}

Fix \(d=N\). On \(\mathscr S\), the Gaussian coordinates are

\begin{equation}
B=\begin{pmatrix}\beta&z_1\\0&z_2\end{pmatrix},\qquad
C=\begin{pmatrix}v_1&v_2\\0&\gamma\end{pmatrix},\qquad
A_{\rm raw}=\widehat A+zv,
\label{m6:s:eq:p006-9-1}
\end{equation}

where

\begin{equation}
\widehat A=\begin{pmatrix}a_{11}&a_{12}\\
-\beta\gamma&a_{22}\end{pmatrix},
\qquad \beta\gamma\ne0.
\label{m6:s:eq:p006-9-2}
\end{equation}

This is a polynomial graph, with rational inverse, in the nine independent coordinates

\begin{equation}
(z_1,z_2,v_1,v_2,a_{11},a_{12},a_{22},\beta,\gamma).
\label{m6:s:eq:p006-9-3}
\end{equation}

For varying nonzero nilpotent \(d\), choose any of the rational local sections \(a(d)\) from \eqref{m6:s:eq:p006-3-3}. Conjugation by \(d(a(d))\) gives local product charts with the two-dimensional nilpotent conjugacy orbit. All transitions are regular on the relevant pivot charts. Thus

\begin{equation}
{\dim\mathscr S=9+2=11.}
\label{m6:s:eq:p006-9-4}
\end{equation}

This proves smoothness on all of \(\mathscr S\), including \(\beta=\gamma\). Closedness in \(\mathscr W\) was proved by the intrinsic equations in Section~\ref{m6:s:p006:s8-1}.

At \(x_{\beta,\gamma}=x_{-\beta\gamma,\beta,\gamma}\), the full orbit dimension is

\begin{equation}
\dim K_S-\dim K_x=12-3=9.
\label{m6:s:eq:p006-9-5}
\end{equation}

The two remaining parameter directions are genuinely transverse: the invariant functions \(\beta,\gamma\) have independent derivatives on the displayed canonical family and vanish on orbit tangents. Therefore

\begin{equation}
{T_x\mathscr S=T_x(K_Sx)+
 F\,\partial_\beta x_{\beta,\gamma}+
 F\,\partial_\gamma x_{\beta,\gamma},\qquad\dim T_x\mathscr S=11.}
\label{m6:s:eq:p006-9-6}
\end{equation}

The relation \(\alpha=-\beta\gamma\) leaves two transverse parameters \(\beta,\gamma\).

\end{proof}

\paragraph{\texorpdfstring{The full conormal in \(GL_4\)}{The full conormal in GL4}}\label{m6:s:p006:s10}

\begin{lemma}\label{m6:s:lem:p006-conormal}
The full conormal of \(\mathscr S\) in \(\mathscr W\) is the
five-dimensional space \eqref{m6:s:eq:p006-10-5} in
\S\ref{m6:s:cert:p006-conormal}. At every point it contains a
nonnilpotent cotangent vector.
\end{lemma}
\begin{proof}
Use \(\lambda_Z(\dot x)=\operatorname{tr}(Z\dot x x^{-1})\).
The full tangent is \eqref{m6:s:eq:p006-9-6}, including both
\(\beta,\gamma\) parameter directions. The certificate
\eqref{m6:s:eq:p006-10-2}--\eqref{m6:s:eq:p006-10-6} tests all orbit and
parameter directions, gives eleven independent equations, and
solves their entire five-dimensional annihilator. Its uniform
witness is
\[
Z_*=\operatorname{diag}(1,-1,-1,1),\qquad
\operatorname{tr}(Z_*^2)=4\ne0.
\]
At \(y=LxR\), right-trivialized cotangents transform by
\(Z\mapsto LZL^{-1}\). Thus the witness exists everywhere,
including \(\beta=\gamma\), with the same trace-square.
The corrected-transpose cotangent formula is
\eqref{m6:s:eq:p006-10-9}; it sends this witness to itself at
\(x_{\gamma,\beta}\).
\end{proof}

\paragraph{\texorpdfstring{Full lift to \(\operatorname{GSp}_{10}\) and the spectral exclusion}{Full lift to GSp10 and the spectral exclusion}}\label{m6:s:p006:s11}

This section makes the ambient spectral step explicit. It is possible to use a descended nilpotent-wavefront statement on \(GL_4\), but the independent proof here applies the conormal obstruction directly to the ambient NW distribution.

\subparagraph{The ambient support manifold}\label{m6:s:p006:s11-1}

Combine the established big-cell coordinates with \eqref{m6:s:eq:p006-2-7}:

\begin{equation}
\Phi:N'\times GL_2\times GL_4\times N'\longrightarrow\Omega',
\quad
(n_1,z,x,n_2)\longmapsto
n_1\iota(d(z)x,z)n_2^t.
\label{m6:s:eq:p006-11-1}
\end{equation}

This is an analytic isomorphism defined over \(F\). On \(V_6\) its \(GL_4\) coordinate lies in \(\mathscr W\). Define

\begin{equation}
\widehat{\mathscr S}=
\Phi(N'\times GL_2\times\mathscr S\times N')\subset V_6.
\label{m6:s:eq:p006-11-2}
\end{equation}

Since \(\mathscr S\) is smooth and closed in \(\mathscr W\), this is a smooth closed submanifold of the open ambient manifold \(V_6\subset G\), with

\begin{equation}
\dim\widehat{\mathscr S}=18+4+11+18=51,
\qquad\dim G=56.
\label{m6:s:eq:p006-11-3}
\end{equation}

The radical and external-base Haar factors in the distribution correspondence are smooth and nowhere zero. Thus support on \(\mathscr S\) after the intrinsic reductions means support of the same \(T|_{V_6}\) on \(\widehat{\mathscr S}\).

\subparagraph{An explicit lift of every conormal direction}\label{m6:s:p006:s11-2}

At \(g_0=\iota(x,I_2)\), use
\(\Lambda_Z(\dot g)=\operatorname{tr}(Z\dot g g_0^{-1})\).
\S\ref{m6:s:cert:p006-ambient-pairing} computes the annihilator
of both radical factors in the full similitude Lie algebra:
\eqref{m6:s:eq:p006-11-5} is the complete Levi element and
\eqref{m6:s:eq:p006-11-6}--\eqref{m6:s:eq:p006-11-7} give its trace pairing
and the two product-factor cotangent vectors \(E,F_2\).

For a \(4\times4\) block matrix \(E\), put

\begin{equation}
\operatorname{pr}(E)=E_{11}+E_{22}\in M_2,
\label{m6:s:eq:p006-11-8}
\end{equation}

where the subscripts in \eqref{m6:s:eq:p006-11-8} denote diagonal \(2\times2\) blocks. The external \(z\)-direction in \eqref{m6:s:eq:p006-11-1}, at \(z=I\), has tangent coordinates \((d(A_0),A_0)\). Annihilating that direction is exactly

\begin{equation}
F_2+\operatorname{pr}(E)=0.
\label{m6:s:eq:p006-11-9}
\end{equation}

Taking the trace in \eqref{m6:s:eq:p006-11-9} and using \eqref{m6:s:eq:p006-11-7} yields \(5\upsilon=0\). Since the field has characteristic zero, \(\upsilon=0\). It follows that the unique lift of a cotangent vector \(E\) on the \(x\)-coordinate, annihilating the radical and external-base directions, is

\begin{equation}
{
\operatorname{Lift}_6(E)=\operatorname{diag}\left(
 \frac12E,
 -\operatorname{pr}(E)+\frac12\operatorname{tr}(E)I_2,
 -\frac12J_4^{-1}E^tJ_4
 \right).
}
\label{m6:s:eq:p006-11-10}
\end{equation}

The central \(2\times2\) block has trace zero, and direct substitution gives
\(\operatorname{Lift}_6(E)^tJ_{10}+J_{10}\operatorname{Lift}_6(E)=0\). Thus the lift lies in the exact \(\mathfrak{sp}_{10}\) inside the chosen similitude Lie algebra.

In particular, the entire ambient conormal is

\begin{equation}
{
CN^{V_6}_{\widehat{\mathscr S},\iota(x_{\beta,\gamma},I)}
=\operatorname{Lift}_6\left(CN^{\mathscr W}_{\mathscr S,x_{\beta,\gamma}}\right).
}
\label{m6:s:eq:p006-11-11}
\end{equation}

The five-dimensional space on the right is given explicitly by \eqref{m6:s:eq:p006-10-5}. The radical/base annihilator conditions in \eqref{m6:s:eq:p006-11-11} determine the unique lift of each of the five residual cotangent vectors, matching the ambient codimension five.

For \eqref{m6:s:eq:p006-10-5}, \(\operatorname{tr}E=0\) and

\begin{equation}
\operatorname{pr}(E)=\frac{\mu-\nu}{\beta}E_{12}.
\label{m6:s:eq:p006-11-12}
\end{equation}

Thus all blocks of every lifted conormal are specified explicitly by \eqref{m6:s:eq:p006-10-5}, \eqref{m6:s:eq:p006-11-10}, and \eqref{m6:s:eq:p006-11-12}.

\subparagraph{The uniform ambient witness}\label{m6:s:p006:s11-3}

For \(E=Z_*\), one has \(\operatorname{pr}(E)=0\) and
\(J_4^{-1}E^tJ_4=E\). Therefore

\begin{equation}
2\operatorname{Lift}_6(Z_*)=
\operatorname{diag}(1,-1,-1,1,0,0,-1,1,1,-1)=\widehat Z_*.
\label{m6:s:eq:p006-11-13}
\end{equation}

It satisfies

\begin{equation}
\widehat Z_*^tJ_{10}+J_{10}\widehat Z_*=0,
\qquad\operatorname{tr}(\widehat Z_*^2)=8\ne0.
\label{m6:s:eq:p006-11-14}
\end{equation}

Thus it is not in the nilpotent cone under the invariant trace-form identification with the dual Lie algebra.

Every point of \(\widehat{\mathscr S}\) is in an original \(H\times H\)-orbit of some \(\iota(x_{\beta,\gamma},I)\), by \eqref{m6:s:eq:p006-11-1} and the rational normalizations above. Under left--right multiplication the right-translated cotangent matrix is conjugated by the left multiplier. Consequently \eqref{m6:s:eq:p006-11-13} supplies a nonnilpotent conormal at every point of \(\widehat{\mathscr S}\).

\subparagraph{Wavefront transfer}\label{m6:s:p006:s11-4}

The coordinate formula \eqref{m6:s:eq:p006-11-1} writes the original equivariant distribution as smooth, nowhere-zero, character-weighted Haar factors in the two radical variables, smooth Haar in the external \(z\)-variable, and the descended distribution in the \(x\)-variable.

A cotangent vector in the wavefront set of the last factor gives a cotangent vector in the original distribution with zero cotangent components in the three smooth factors. This elementary product assertion follows in local product charts by taking product test functions: the Fourier transform of each compactly localized smooth factor is nonzero at zero, and its Fourier transform introduces no singular direction. Analytic coordinate changes transport wavefront sets by the differential \cite[Lemma 2.15]{JLX}.

The resulting right-translated lift is exactly \eqref{m6:s:eq:p006-11-10}. If it is nilpotent, its first \(4\times4\) block \(E/2\) is nilpotent, hence \(E\) is nilpotent. Thus the descended \(GL_4\) distribution does inherit nilpotent wavefront fibers on its open domain.

The next step applies the ambient witness
\eqref{m6:s:eq:p006-11-13} directly to the original distribution.

\subparagraph{Exclusion for the same original global distribution}\label{m6:s:p006:s11-5}

By the established Model 7 open result, \(T|_\Omega=0\). In particular the descended distribution vanishes on \(\{\det d\ne0\}\), corresponding to \(\mathcal P_1\) in \(V_6\). The intrinsic support reductions in Section~\ref{m6:s:p006:s8} therefore imply

\begin{equation}
\operatorname{supp}(T|_{V_6})\subset\widehat{\mathscr S}.
\label{m6:s:eq:p006-11-15}
\end{equation}

Proposition~\ref{cf:conormal} now applies to the original \(T\), the open \(V_6\), and the smooth closed manifold \(\widehat{\mathscr S}\). The support inclusion is \eqref{m6:s:eq:p006-11-15}; the full conormal contains the transported witness \eqref{m6:s:eq:p006-11-13} at every point, with trace-square \(8\ne0\). The assumed NW condition of the original ambient distribution supplies the nilpotent wavefront fibers.
Hence

\begin{equation}
{T|_{V_6}=0.}
\label{m6:s:eq:p006-11-16}
\end{equation}

The closed-support extensions of the intrinsic odd orbital distributions in \S\ref{m6:s:p006:s8-3} have the nonnilpotent conormal direction just computed.

\paragraph{\texorpdfstring{Proof of Proposition~\ref{m6:s:prop:p006}}{Proof of the rank-one boundary proposition}}\label{m6:s:p006:s12}
\begin{proof}[Completion of the proof]
The opens satisfy
\begin{equation}
 \mathcal O_6=\Omega\cup V_6,\qquad
 \Omega\cap V_6=\mathcal P_1,\qquad
 V_6\subset\Omega'.
 \label{m6:s:eq:p006-12-1}
\end{equation}
The NW conclusion of Proposition~\ref{m6:rank-dictionary} gives
$T|_{V_6}=0$. Together with the stated predecessor $T|_\Omega=0$,
open gluing gives
\begin{equation}
 T|_{\mathcal O_6}=0.
 \label{m6:s:eq:p006-12-6}
\end{equation}
The independent support and conormal proof above gives the same
conclusion by the independent rank-stratum argument. The intrinsic
odd examples in \S\ref{m6:s:p006:s8-3} prove the remaining nonvanishing
assertion.
\end{proof}

\subsubsection{The rank-zero outer cell}\label{m6:s:sec:outer}

The rank-zero cell admits a complete distributional description,
including its nonsplit and Jordan fibers. Expanded coordinate formulas
are in \path{anc/model6_calculations.pdf}, \texttt{M6-RANKZERO};
\S\ref{m6:s:app:p034-fibers} describes the invariant fiber measures.

\begin{proposition}\label{m6:s:prop:p034}
Let
\begin{equation}
X_0=\left\{\begin{pmatrix}A&B\\C&0\end{pmatrix}\in GL_4(F)\right\}.
\label{m6:s:eq:p034-0-1}
\end{equation}
Then \(D'_{\widetilde\eta_S}(X_0)=0\), and under the parabolic extraction one has
\begin{equation}
D'_{\widetilde\eta}(\mathcal P_{34})=0.
\label{m6:s:eq:p034-0-2}
\end{equation}
Every even-character equivariant distribution on \(X_0\) is supported on
\begin{equation}
\mathscr Y_0=\left\{\begin{pmatrix}A&B\\B&0\end{pmatrix}:A\in M_2(F),\ B\in GL_2(F)\right\},
\label{m6:s:eq:p034-0-3}
\end{equation}
and there is an explicit transpose-compatible isomorphism
\begin{equation}
D'_{\eta_S}(X_0)\simeq D'(GL_2(F))^{\operatorname{Ad}GL_2(F)}.
\label{m6:s:eq:p034-0-4}
\end{equation}
\end{proposition}

\paragraph{Rational chart and the full stabilizer.}\label{m6:s:p034:s1}\label{m6:s:p034:s2}\label{m6:s:p034:s3}
Use the common extraction of \S\ref{supp:m6-shalika-coordinates};
its radical, external $GL_2$, full-similitude density and open/closed
dictionary remain unchanged. A kernel vector for $B$, or a left kernel
vector for $C$, would give a kernel for $x$, so both are invertible.
Conversely
\begin{equation}
 x^{-1}=\begin{pmatrix}0&C^{-1}\\B^{-1}&-B^{-1}AC^{-1}\end{pmatrix},
 \qquad\det x=\det B\det C.
 \label{m6:s:eq:p034-2-1}
\end{equation}
Thus $X_0\simeq M_2\times GL_2\times GL_2$ is smooth, closed in $GL_4$, of
dimension $12$. The action is
\begin{equation}
 (A,B,C)\longmapsto(aAa^{-1}+XaCa^{-1}+aBa^{-1}V,
                         aBa^{-1},aCa^{-1}).
 \label{m6:s:eq:p034-2-4}
\end{equation}
The global normalization $(-AC^{-1},0,I)$ sends $x$ to $y_{B,C}$ with
upper block zero. Equivalently $A=ZC$, with
$dA=|\det C|^2dZ$. It preserves the simultaneous-conjugacy problem for
$(B,C)$: independent left/right normalization of these blocks is not
allowed. Finite rational companion charts for nonscalar $B$ use one of
$e_1,e_2,e_1+e_2$ as a cyclic vector; if $B$ is scalar apply them to $C$.
The formulas are given in \S\ref{m6:s:cert:p034-pair-charts}.

\begin{proposition}\label{m6:s:prop:p034-stabilizer}
The full stabilizer at every point of \(X_0\) is given by \eqref{m6:s:eq:p034-3-2}. It is connected and unimodular, with dimensions eight, six, or five according to the common centralizer. Its full character is trivial precisely on \(C=B\).
\end{proposition}
\begin{proof}
The entire solved stabilizer is
\begin{equation}
 a\in Z(B,C),\quad X\in M_2,\quad
 V=B^{-1}(A-aAa^{-1}-XC),\qquad
 Z(B,C)=Z_{GL_2}(B)\cap Z_{GL_2}(C).
 \label{m6:s:eq:p034-3-2}
\end{equation}
The action gives $aB=Ba$, $aC=Ca$, and
$XC+BV=A-aAa^{-1}$, proving the formula at the group level.
At $y_{B,C}$ the stabilizer is $M_2\rtimes Z(B,C)$, with conjugation on
$M_2$; the front-end extraction embeds it into the original
$H\times H$ and uniqueness of Gauss coordinates excludes extra radicals.
A nonscalar $2\times2$ matrix has matrix centralizer $F[B]$, since it is
cyclic. Hence the common centralizer is $GL_2$ if both matrices are
scalar, the units of a quadratic algebra if they commute and at least
one is nonscalar, and the scalar $\mathbb G_m$ otherwise: a common
nonscalar centralizer element would put both matrices in its commutative
quadratic algebra. The quadratic algebras are $F\times F$, a quadratic
field, and $F[\epsilon]/\epsilon^2$. Their unit groups are respectively
$\mathbb G_m^2$, $\operatorname{Res}_{E/F}\mathbb G_m$, and
$\mathbb G_m\times\mathbb G_a$. These are connected and unimodular,
and $\det(\operatorname{Ad}(a);M_2)=1$. Stabilizer dimensions are
$8,6,5$ according to common-centralizer dimension $4,2,1$;
$K$ is also unimodular, so the relative density is one.
Taking trace in the solved $V$ equation, the $A$ terms cancel because
$a$ commutes with $B$, giving
\begin{equation}
 \eta_S|_{K_x}=\psi\!\left(\operatorname{tr}((I-CB^{-1})X)\right).
 \label{m6:s:eq:p034-3-10}
\end{equation}
If $Q=I-CB^{-1}\ne0$, choose $Q_{ij}\ne0$; the genuine subgroup
$(tE_{ji},-tB^{-1}E_{ji}C,I)$ has nontrivial character $\psi(tQ_{ij})$.
Conversely $C=B$ makes the character one on every stabilizer component.
\end{proof}

\paragraph{Closed support and exact descent.}\label{m6:s:p034:s4}\label{m6:s:p034:s5}
On the stable open $B\ne C$, the preceding subgroup and trivial
relative modulus give zero by the stabilizer-character criterion
\cite[Theorem 2.5 and Remark 2.6]{JLX}. The ordinary closed-support
sequence therefore identifies $D'_{\eta_S}(X_0)$ with
$D'_{\eta_S}(\mathscr Y_0)$ and preserves transpose parity. This uses the left-exact support sequence.

\begin{proposition}\label{m6:s:prop:p034-correspondence}
The averaging map \eqref{m6:s:eq:p034-5-13} gives the isomorphism \eqref{m6:s:eq:p034-0-4}; its transpose action is ordinary transpose on \(GL_2\), with no residual density character.
\end{proposition}
\begin{proof}
The global rational product chart and its action are
\begin{align}
 \Theta(Z,B)&=\begin{pmatrix}ZB&B\\B&0\end{pmatrix},
 \label{m6:s:eq:p034-5-1}\\
 (Z,B)&\longmapsto(aZa^{-1}+X+B'VB'^{-1},B'),\quad B'=aBa^{-1}.
 \label{m6:s:eq:p034-5-2}
\end{align}
The quotient is $B$, with section $Z=0$. The additive fiber stabilizer
is $(X,-B^{-1}XB)$, with character one, and the full point stabilizer is
$M_2\rtimes Z_{GL_2}(B)$. Translation changes $\operatorname{tr}Z$ by
$\operatorname{tr}X+\operatorname{tr}V$, with absolute Jacobian one.
The fixed normalization of fiber measure is
$\psi(\operatorname{tr}Z)dZ$, or
$|\det B|^{-2}\psi(\operatorname{tr}(AB^{-1}))dA$ in the original chart.
For $f\in D(X_0)$ the resulting formula is
\begin{equation}
 S_\Lambda(f)=\Lambda\!\left[B\longmapsto\int_{M_2}
 f\!\left(\begin{pmatrix}ZB&B\\B&0\end{pmatrix}\right)
 \psi(\operatorname{tr}Z)dZ\right].
 \label{m6:s:eq:p034-5-13}
\end{equation}
The integrand is compactly supported in the global product chart.
Averaging is onto: choose $\rho\in D(M_2)$ with
$\int\rho(Z)\psi(\operatorname{tr}Z)dZ=1$ and use
$\rho(Z)\varphi(B)$. On $X_0$, an explicit extension is
$\rho(AB^{-1})\varphi(B)\kappa(CB^{-1})$, with
$\kappa\in D(GL_2)$, $\kappa(I)=1$; compactness follows from the
invertible product change $(A,B,C)\leftrightarrow(AB^{-1},B,CB^{-1})$.
Uniqueness of translation-equivariant Haar on $M_2$, and
$D(M_2\times GL_2)=D(M_2)\otimes D(GL_2)$, show that every equivariant
distribution has this unique form. Equivariance under $a$ is exactly
conjugation invariance of $\Lambda$; the second additive factor imposes
no further condition, by the trace and Jacobian calculation above.
Finally transpose sends $(Z,B)$ to $(B^tZ^tB^{-t},B^t)$, preserving
trace and absolute fiber Jacobian. Thus, writing the inner averaging as
$\mathcal J_0$,
\begin{equation}
 \mathcal J_0(f\circ\tau)(B)=\mathcal J_0f(B^t).
 \label{m6:s:eq:p034-5-16}
\end{equation}
\end{proof}

\paragraph{All rational orbits and their transpose signs.}\label{m6:s:p034:s6}\label{m6:s:p034:s7}
For fixed $(t,\delta)=(\operatorname{tr}B,\det B)$, $\delta\ne0$,
if $t^2-4\delta\ne0$ the entire fiber is one rational orbit by a cyclic
companion basis. It is closed and smooth of dimension two, with full
centralizer a split torus or the multiplicative group of the quadratic
field. This includes ramified dyadic extensions. For a repeated root
$\lambda=t/2\ne0$ the fiber is
\begin{equation}
 \chi^{-1}(2\lambda,\lambda^2)=\{\lambda I+N:N^2=0\},
 \label{m6:s:eq:p034-6-2}
\end{equation}
with exactly the scalar orbit and the nonzero Jordan orbit.
For $N\ne0$, choose $Ne_j\ne0$; the basis $(Ne_j,e_j)$ conjugates $N$
to $E_{12}$. These two rational pivot charts prove transitivity over $F$. The Jordan centralizer is $\{s(I+hN)\}$ and its boundary
is the scalar orbit, as $\lambda I+\varpi^kN\to\lambda I$.
The full compatible $K$-orbits are $M_2$ times these residual orbits,
and their homogeneous topology follows from the product chart and
rational companion/pivot sections.

With $J(t,\delta)$ the companion matrix and
$S_t=\left(\begin{smallmatrix}0&1\\1&t\end{smallmatrix}\right)$,
\begin{equation}
 S_tJ=J^tS_t,\qquad S_t^t=S_t,\qquad\det S_t=-1.
 \label{m6:s:eq:p034-7-1}
\end{equation}
If $B=P_BJP_B^{-1}$, put $a_B=P_B^{-t}S_tP_B^{-1}$; for scalar $B$
put $a_B=I$. Then
\begin{equation}
 a_B^t=a_B,\qquad a_BBa_B^{-1}=B^t.
 \label{m6:s:eq:p034-7-3}
\end{equation}
The transporter $k_B=(0,0,a_B)$ has character one. At
$x=\Theta(Z,B)$ the transporter is
$(-a_BZa_B^{-1},Z^t,a_B)$, still of character one.
Both $K$ and its full stabilizer are unimodular, so the orbital line is
one-dimensional and generated by positive quotient Haar weighted by
$\eta_S$. Lemma~\ref{m8:s:c:transporter} gives transpose scalar $+1$.
Equivalently, $k_B^{-1}\tau$ is an order-two odd stabilizer component:
its absolute Haar modulus is one while its prescribed character is $-1$.
All rational orbits thus have zero odd orbital space. Constructiveness
\cite[Theorem 2.3 and Remark 2.4]{JLX} yields
\begin{equation}
 D'_{\widetilde\eta_S}(\mathscr Y_0)=0.
 \label{m6:s:eq:p034-7-15}
\end{equation}
Closed support and common extraction prove
Proposition~\ref{m6:s:prop:p034}.

\paragraph{The open consequence.}\label{m6:s:p034:s10}
The exact rank dictionary gives
$\mathcal U=\Omega\cup\Omega'$,
$\Omega\cap\Omega'=\mathcal P_1$,
$V_6=\mathcal P_1\sqcup\mathcal P_6$, and
$\mathcal U\setminus\mathcal O_6=\mathcal P_{34}$; $V_6$ uses a
nonzero corner matrix. Model-7 extraction and
Proposition~\ref{m6:rank-dictionary} kill the same original NW odd
ambient distribution on $\mathcal U$. The intrinsic rank-zero proof gives the independent local even
classification. The rank-one step in
Proposition~\ref{m6:s:prop:p006} uses ambient nilpotent wavefront.

\subsection{Boundary calculations}
\subsection{Rational quotient calculations for the three middle families}\label{m6:supp:middle-orbits}

For each intersection, refine the two flags simultaneously. Gaussian
elimination using the root groups below removes the maps that
strictly increase the common flag order, while the symplectic
equation pairs opposite entries according to \(\mathcal N\).
The only opposite-root pair for \(W_{16}\) is \(B_{12},B_{21}\); it gives the
entire \(\mathrm{GL}_2\) factor, including its Weyl elements. For \(W_{17}\)
and \(W_{27}\) there is no opposite-root pair. All other diagonal entries
are arbitrary elements of \(F^\times\), and the kernel of the common graded
action is the split unipotent group given by the listed roots. This proves
the group-level exhaustion used below.

\paragraph{The cell \texorpdfstring{$\mathcal P_{16}$}{P16}}\label{m6:s:sec:middle-p016}

The signed representative is
\(W_{16}=(1,5,3,4,2,-9,7,8,-6,10)\), with \(W_{16}^t=W_{16}=W_{16}^{-1}\).

\subparagraph{Entire family and intersection}

The compatibility equations give exactly

\[
 \mathsf{P}=\begin{pmatrix}p_{11}&p_{12}\\p_{21}&p_{22}\end{pmatrix},\quad
 \mathsf{Q}=\begin{pmatrix}q_{11}&q_{12}\\q_{21}&q_{22}\end{pmatrix},\quad
 \mathsf{R}=\begin{pmatrix}r_{11}&q_{11}\\r_{21}&q_{12}\end{pmatrix},\quad
 \mathsf{S}=\begin{pmatrix}p_{12}&p_{22}\\s_{21}&s_{22}\end{pmatrix},
\]

with all four determinants nonzero. This is a smooth, 12-dimensional open subset of affine 12-space. The column \((p_{12},p_{22})^t\) and row \((q_{11},q_{12})\) are nonzero; each individual entry may vanish. The identity \(m(A,B,C)=d_5(C)\mathsf l(B^{-1}A,C^{-1}B)\)

and the established rational parabolic-cell presentation cover all original rational points, including arbitrary similitude multipliers.

For \(J=P\cap wP^tw^{-1}\), the complete intersection has dimension \(17\), Levi factor
\(\mathbb G_m^4\times\mathrm{GL}_{2,B}\), and nine-dimensional radical.
Its common unipotent kernel is the five-dimensional split Heisenberg group
\(K=\exp\langle X_{21},X_{22},Y_{11},Y_{21},Z_{21}\rangle\).
The exact projection matrices \(A,C,A_r,C_r\), their linked diagonal
entries and the rational group section are recorded in
Appendix~\ref{m6:s:cert:p016-intersection}. They exhaust \(J/K\) over \(F\),
including every component and the full \(B\in\mathrm{GL}_2(F)\) factor.

Consequently its complete coefficient action is

\begin{equation}
\label{m6:s:eq:middle-a16-1}
(\mathsf{P},\mathsf{Q},\mathsf{R},\mathsf{S})\mapsto
(B^{-1}\mathsf{P}A,\ C^{-1}\mathsf{Q}B,\ B^t\mathsf{R}A_r^{-1},\ C_r\mathsf{S}B^{-t}). 
\end{equation}

Substitution shows that the action preserves all four compatibility
equations. Setting \(B=I\) gives an eight-dimensional subgroup
\(K_8\), which commutes with the \(B\) factor in \(J/K\).

\subparagraph{Character and all components of the stabilizer}

The adjacent blocks in conjugation by the displayed signed representative give the exact additive functional

\begin{equation}
\label{m6:s:eq:middle-a16-2}
\lambda=(p_{12}-s_{11})x_{21}+(p_{22}-s_{12})x_{22}
 +(q_{11}-r_{12})y_{11}+(q_{12}-r_{22})y_{21}. 
\end{equation}

All remaining roots in \(J\) have zero character phase. After
the coefficient-frame change, the non-kernel root generators and
Levi generators have zero additive character on both sides; for
\(K\), cancellation is given by \eqref{m6:s:eq:middle-a16-2}. Since
these generators preserve compatibility, their products have
trivial moving-frame gauge cocycle.

Put \(v=\mathsf{P}e_2\), \(\phi=e_1^t\mathsf{Q}\), \(d=\det \mathsf{P}\det \mathsf{R}\), and \(N=v\phi/d\). The full group membership equations reduce to \(Bv=\beta v,\qquad \phi B=\beta\phi,\)

which are equivalent to \(BN=NB\), including when \(\phi v=0\). Indeed equality of the two rank-one products first preserves the nonzero image line, and then equates the row eigenvalue to the same scalar. The four lower projection entries are uniquely forced by \(B\). Their required diagonal pairs are \((\alpha,\beta)\) and \((\beta,\alpha)\). Every such tuple has a rational lift in \(J\) by the preceding group projection computation. Therefore, on every point of this entire family,

\begin{equation}
\label{m6:s:eq:middle-a16-3}
1\longrightarrow K\longrightarrow\mathcal A_x
\longrightarrow Z_{\mathrm{GL}_2}(N)\longrightarrow1, \qquad
\dim\mathcal A_x=7. 
\end{equation}

To describe the full group on both branches, take the uniform matrix \(N_\nu=\begin{pmatrix}\nu&1\\0&0\end{pmatrix}\). Its centralizer is explicitly

\begin{equation}
\label{m6:s:eq:middle-a16-4}
\left\{\begin{pmatrix}u&v\\0&u-\nu v\end{pmatrix}:
u(u-\nu v)\ne0\right\}. 
\end{equation}

For \(\nu \ne 0\) this is a split two-dimensional torus; for \(\nu =0\) it is \(\mathbb G_m\times\mathbb G_a\). Thus \eqref{m6:s:eq:middle-a16-4} gives dimension two on both
branches. The explicit projection lifts make the quotient map in
\eqref{m6:s:eq:middle-a16-3} surjective on \(F\)-points. For \(\nu \ne 0\) the torus lifts already in a rational diagonal coefficient chart, so the group is a split torus by \(K\). For \(\nu =0,\) divide by the scalar entry \(u\): the remaining ten-by-ten matrix is unipotent. Thus the full stabilizer is the central scalar \(\mathbb G_m\) times a six-dimensional split unipotent group. Its finite logarithm includes every rational point.

The algebraic additive phase vanishes on the stabilizer, so the
entire unipotently compatible locus is full-character-compatible.

\subparagraph{Density calculation}

Write the diagonal entries of \(A,B,C\) as
\((\alpha,\beta),(\gamma,\delta),(\epsilon,\zeta)\).
The radical and kernel weight products are respectively
\((\beta/\epsilon)^6\) and \((\beta/\epsilon)^3\), as computed in
Appendix~\ref{m6:s:cert:p016-weights}. Put \(q=a_2/c_1=\beta/\epsilon\). Then \(\Delta_J=|q|^{-6},\quad \Delta_{J/K}=\Delta_{K_8}=|q|^{-3},\quad j_{U^2/K}=|q|^6.\)

On a full stabilizer, fixing \(d\) forces \(q=1\), so the
action on \(K\) has determinant one. Both centralizers in
\eqref{m6:s:eq:middle-a16-4} are unimodular, and hence so is the
full stabilizer. Since \(\mathcal A=H^2\) is unimodular as well,
the relative orbital density character is one on all of
\(\mathcal A_x\).

For the finer coefficient reduction, the test-function convention gives \(D(f\circ k)=|q|^6D(f).\)

To determine the density character, integrate over the right \(J\)-fiber with left Haar measure. Right translation has factor \(\Delta_J^{-1}\). Formula \eqref{m6:s:eq:middle-a16-1} gives \(d^\prime=qd\). Multiplication of a distribution by \(\lvert d\rvert^6\) then gives the invariant distribution \(D_0=|d|^6D\), because \(|d|^6\circ k^{-1}=|q|^{-6}|d|^6\). This establishes the correction with the stated test-function convention.

\subparagraph{residual quotient, including all fibers}

Set \(\kappa=\det \mathsf{Q}\det \mathsf{S}\). Direct substitution in the full gauge gives

\begin{equation}
\label{m6:s:eq:middle-a16-5}
d'=qd,\qquad \kappa'=\kappa,\qquad N'=B^{-1}NB. 
\end{equation}

To identify the map
\((\mathsf{P},\mathsf{Q},\mathsf{R},\mathsf{S})\mapsto(N,\kappa)\)
as the quotient by the free \(K_8\) action, we determine its fibers.
Appendix~\ref{m6:s:cert:p016-fibers} proves that each fiber is a free transitive
\(K_8(F)\)-orbit and gives rational sections on the nonzero-entry charts
of a rank-one matrix. Only a chosen nonzero vector coordinate is inverted.
The sections have \(d=1\), include the nonzero nilpotent locus, and make
the quotient a surjective analytic submersion.

After the density correction, the free right \(K_8\) action is extracted with right Haar measure on \(K_8\). Changes between these \(d=1\) sections are left translations by elements with \(q=1\). Their \(K_8\) modular factor is therefore one. Transport by \(B\) keeps \(d=1\), and transpose also keeps \(d=1\); restoring a local section again requires only \(q=1\) gauges. The resulting residual problem is the untwisted \(\mathrm{GL}_2\)-conjugation action on \(\{N\in M_2(F):\operatorname{rank}N=1\}\times F^\times_\kappa.\)

\subparagraph{Rational orbits, nilpotent test, transpose and topology}

A rank-one \(N\) has characteristic polynomial \(T(T-\nu)\), with
\(\nu=\operatorname{tr}N\). If \(\nu\ne0\), its rational image and kernel
give a split diagonal normal form. If \(\nu=0\), choose \(u\) with
\(Nu\ne0\); the basis \((Nu,u)\) gives a nonzero nilpotent normal form.
There is consequently one rational \(\mathrm{GL}_2\)-orbit for each
\(\nu\in F\), with its split semisimple and nonzero nilpotent
branches as described above.

The canonical coefficient section is

\[
\mathsf{P}_0=\begin{pmatrix}0&1\\-1&0\end{pmatrix},\quad
\mathsf{Q}_{\nu,\kappa}=\begin{pmatrix}\nu&1\\-\kappa&0\end{pmatrix},\quad
\mathsf{R}_\nu=\begin{pmatrix}1&\nu\\0&1\end{pmatrix},\quad \mathsf{S}_0=I_2.
\]

It gives exactly \(N_\nu\) and \(\kappa\), including \(\nu=0\), in the ten-dimensional representative \(\mathsf l(\mathsf{P}_0,\mathsf{Q}_{\nu,\kappa})w\mathsf l(\mathsf{R}_\nu,I)^t\). The entire rational orbit parameter space in the original compatible locus is therefore \(F\times F^\times\), and every fiber has dimension \(52-7=45\). The compatible locus has dimension \(52+12-17=47\). Surjective local quotient charts from the parabolic presentation descend its smoothness and closedness in the cell. More explicitly, its preimage is the saturated closed four-equation subset of the open coefficient space, and the surjective submersion is an open quotient map.

The invariant map \((\nu,\kappa)\) is continuous on the compatible
locus, and each fiber is the single rational orbit just classified.
The fibers are closed there, hence closed in the whole cell.

Transpose swaps \((\mathsf P,\mathsf Q,\mathsf R,\mathsf S)\) with \((\mathsf R,\mathsf S,\mathsf P,\mathsf Q)\). Hence \(d\) and \(\kappa\) are fixed and \(N\) goes to \(N^t\). The uniform transporter

\[
b_\nu=\begin{pmatrix}\nu&1\\1&0\end{pmatrix},\qquad
\det b_\nu=-1,\qquad b_\nu N_\nu=N_\nu^t b_\nu
\]

is valid at \(\nu =0\) as well as all \(\nu \ne 0\). The gauge phase is one by \eqref{m6:s:eq:middle-a16-2}, and the previously checked line through extraction has no sign. The involution on the original \(J\)-fiber preserves Haar measure because its positive Haar modulus has square one; the same is true of its induced positive orbital density after the character phase has been removed.

The continuous invariant map to \(F\times F^\times\) is fixed by
conjugation and transpose. On every fiber, including \(\nu=0\),
the full unimodular character-trivial stabilizer gives one orbital line,
and the character-one transporter makes it even by
Lemma~\ref{m8:s:c:transporter}. Proposition~\ref{cf:localization} and the
open--closed sequence, after the incompatible open is removed, give
\[
{D'_{\mathrm{odd}}(\mathcal P_{16})=0.}
\]
Every compatible orbit is closed in the cell, so its nonzero weighted
Haar line has a closed pushforward. Hence
\(D'_\eta(\mathcal P_{16})\ne0\).

\paragraph{The cell \texorpdfstring{$\mathcal P_{17}$}{P17}}\label{m6:s:sec:middle-p017}

Here
\(W_{17}=(1,5,3,-7,2,-9,4,8,-6,10)\).

\subparagraph{Entire family and full intersection}

The compatibility equations give

\[
\mathsf{P}=\begin{pmatrix}a&b\\c&d\end{pmatrix},\quad
\mathsf{Q}=\begin{pmatrix}e&0\\f&g\end{pmatrix},\quad
\mathsf{R}=\begin{pmatrix}h&e\\i&j\end{pmatrix},\quad
\mathsf{S}=\begin{pmatrix}b&0\\k&l\end{pmatrix},
\quad b e g l\det \mathsf{P}\det \mathsf{R}\ne0.
\]

All eight of a,c,d,f,h,i,j,k are unrestricted. The intersection calculation gives the six diagonal roots, \(A_{21}\),\(B_{21}\),\(C_{21}\), and \(X_{21},Y_{11},Y_{21},Z_{11},Z_{21},W_{21},T_h,T_f,S_f.\)

Thus the full parabolic intersection has dimension 18, a split six-dimensional torus Levi, and radical dimension 12. Its complete common unipotent subgroup is

\[
K=\exp\langle X_{21},Y_{11},Y_{21},Z_{21},W_{21},T_f\rangle,
\qquad \dim K=6.
\]

The full quotient J/K has six lower root coordinates. The roots \(A_{21}\),\(B_{21}\),\(C_{21}\) give the three left lower entries; \(Z_{11}\),\(T_h\),\(S_f\) give the three right lower entries and no other Levi entries. The torus and these root groups therefore realize all projection tuples over \(F\). Their products can change only the already retained kernel coordinate.

For left diagonal pairs (\(\alpha ,\beta ),(\gamma ,\delta ),(\epsilon ,\zeta ),\) the inverse-natural right diagonal pairs are

\begin{equation}
\label{m6:s:eq:middle-a17-1}
(\alpha^{-1},\epsilon^{-1}),\quad
(\gamma^{-1},\delta/(\epsilon\zeta)),\quad
(\beta^{-1},\beta/(\epsilon\zeta)). 
\end{equation}

This matches the directly recomputed ten-by-ten conjugation.

\subparagraph{group stabilizer and character}

The signed right-projection and adjacent-root arrays are
Appendix~\ref{m6:s:cert:p017-projections}; they use
\(N'=-(w^{-1}Nw)^t\), with all signs in the fixed ten-dimensional form.

Let \(A,B,C,D,E,F\) be the six lower-triangular projection matrices, with the linkage \eqref{m6:s:eq:middle-a17-1}. Membership in the original \(H^2\) is exactly

\begin{equation}
\label{m6:s:eq:middle-a17-2}
\mathsf{Q}\mathsf{P}A=C\mathsf{Q}\mathsf{P},\qquad \mathsf{Q}B=C\mathsf{Q},\qquad \mathsf{S}\mathsf{R}D=F\mathsf{S}\mathsf{R},\qquad \mathsf{S}E=F\mathsf{S}. 
\end{equation}

The complete elimination of \eqref{m6:s:eq:middle-a17-2} is
Appendix~\ref{m6:s:cert:p017-elimination}. It uses only the mandatory
\(ebl\ne0\) and invertibility, and forces
\(A=B=C=\beta I, D=E=F=\beta^{-1}I\).

Consequently

\[
\mathcal A_x\simeq\mathbb G_m\times K,
\qquad \dim\mathcal A_x=7,
\]

with the scalar pair \((zI,z^{-1}I)\) and the full six-dimensional split unipotent kernel conjugated by the two coefficient lifts. All rational solutions, including all group components, are included. This dimension does not depend on any free coefficient entry or on the eventual parameter \(\kappa \).

The adjacent-root blocks give \(\lambda=bX_{21}+eY_{11}-eY_{11}-bX_{21}=0.\)

The potential \(W_{21}\) contribution on the right multiplies \(s_{12}=0\). The other kernel directions have no character contribution. The scalar factor has character one. Thus the full stabilizer has exact additive character zero. The entire unipotently compatible family remains fully compatible.

\subparagraph{Density}

Appendix~\ref{m6:s:cert:p017-weights} gives the complete radical-weight array.
Its product is \(\beta^6\delta^6/(\epsilon^9\zeta^3)\), and the kernel
product is \(\beta^3\delta^4/(\epsilon^5\zeta^2)\). Thus the raw
38-dimensional \(U^2/K\) Jacobian and intersection modulus are
\[
j=\left|\frac{\beta^6\delta^6}{\epsilon^9\zeta^3}\right|,
\qquad\Delta_J=j^{-1}.
\]

 On the full stabilizer all six diagonal parameters coincide, so this factor is one. More intrinsically, the stabilizer is the direct product of a central split torus and a unipotent group; it is unimodular. Therefore the full relative orbital density line is trivial.

\subparagraph{Rational normalization}

The complete gauge and its six lower-unipotent operations, followed by
the linked diagonal operation, are given in
Appendix~\ref{m6:s:cert:p017-normalization}. They produce
\[
\mathsf P_\kappa=\begin{pmatrix}0&1\\1&0\end{pmatrix},\quad
\mathsf Q_\kappa=\mathsf S_\kappa=I,\quad
\mathsf R_\kappa=\begin{pmatrix}0&1\\\kappa&0\end{pmatrix},\qquad
\kappa=\det\mathsf P\det\mathsf R\,gl\in F^\times.
\]
The divisions are only by \(b,e,g,l,\det\mathsf P,\det\mathsf R\), all
forced nonzero. The invariant is complete under the entire rational
gauge, so these are \(F\)-transporters.

Every \(\kappa\in F^\times\) occurs, including nonsquare values for which \(R_\kappa\) has nonsplit characteristic polynomial \(T^2-\kappa\).

The compatible locus has dimension 46. The map \(H^2\times F^\times\to C_{17}\), using the displayed canonical section, is a surjective analytic submersion: the orbit derivative has rank 45 from the full group stabilizer, and \(\kappa \) of the section is the identity parameter. Its fibers are the complete stabilizers above. Smoothness and closedness also descend from the saturated four-equation coefficient locus under the original open submersion. Every \(\kappa \) fiber is a single closed rational orbit in the whole cell.

\subparagraph{homogeneous conclusion and its rational meaning}

Let \(V=U^2\) be the normal unipotent subgroup of \(H^2\),
with Levi quotient \(K_0\)=\(\mathrm{GL}_2\times \mathrm{GL}_2\). On each homogeneous orbit \(H^2/ \mathcal A_x\), the projection of the \emph{full} stabilizer is \(Z_- =\{(zI,z^{-1}I):z\in F^\times\},\)

and its intersection with V is the complete six-dimensional unipotent K. Quotienting the \(V\)-orbit directions therefore gives exactly \(K_0/Z_-\). The 38-dimensional V-fiber and its character/density have already been accounted for. \(Z_-\) acts trivially on its own stabilizer and on the residual tangent line; the resulting \(K_0\)-homogeneous density and character are trivial.

The residual space is consequently \(\bigl((\mathrm{GL}_2\times\mathrm{GL}_2)/Z_-\bigr)\times F^\times.\)

The map \([(a,b)]\mapsto(ab^t,[b])\) identifies it with \(\mathrm{GL}_2\times \mathrm{PGL}_2\times F^\times\). This is rationally surjective: \(\mathrm{PGL}_2\) is the nonzero-determinant open subset of projective 3-space, and an \(F\)-rational projective point has an \(F\)-valued matrix representative. Equivalently, the elementary scalar ambiguity is exactly \(Z_-\).

In these coordinates the action is \((A,B):(c,[b],\kappa)\mapsto(AcB^t,[Bb],\kappa).\)

Frobenius descent over the identity in the first residual factor has stabilizer \(\{(B^{-t},B)\}\); its tangent action has determinant one. This reduces the distribution problem to the regular left \(\mathrm{GL}_2\) action on \(\mathrm{PGL}_2\), with the passive parameter \(\kappa \) and trivial density.

\subparagraph{Transpose, Haar sign, localization and extension}

Writing \(d_0=\operatorname{diag}(1,-1)\), multiplication gives the original ten-by-ten polynomial identity for arbitrary \(\kappa \): \(x_\kappa^t=d_5(d_0)x_\kappa d_5(d_0)^t.\)

Both transporter factors have character one and the associated extended stabilizer element has order two. The coefficient transpose sends \((\mathsf{P},\mathsf{Q},\mathsf{R},\mathsf{S})\) to \((d_0\mathsf{R},\mathsf{S}d_0,\mathsf{P},\mathsf{Q})\) and preserves \(\kappa \). On the double quotient it is \([(a,b)]\mapsto[(bd_0,ad_0)]\); on the final \(\mathrm{PGL}_2\) conclusion it is \([b]\mapsto[b^{-t}d_0].\)

Inverse transpose is an involutive automorphism, hence its positive Haar modulus is one. Right translation by \([d_0]\) preserves Haar because \(\mathrm{PGL}_2\) is unimodular. Thus the Haar line has sign +1. The residual base \(F^\times\) is fixed, and every fiber has zero odd-distribution space. Proposition~\ref{cf:localization} is applicable and yields zero odd distributions on the compatible locus.

The incompatible open has unipotent-character vanishing; the
closed-support injection therefore gives
\[
{D'_{\mathrm{odd}}(\mathcal P_{17})=0.}
\]
Each compatible \(\kappa\)-orbit is closed in the cell, with trivial full
character and relative density, so its weighted Haar line has a nonzero
closed pushforward. Thus \(D'_\eta(\mathcal P_{17})\ne0\).

\paragraph{The cell \texorpdfstring{$\mathcal P_{27}$}{P27}}\label{m6:s:sec:middle-p027}
\subparagraph{Entire compatible family}

Put \(\mathsf l(\mathsf{P},\mathsf{Q})=m(\mathsf{Q}\mathsf{P},\mathsf{Q},I)\), \(\mathcal A=H\times H\), and \(x=\mathsf l(\mathsf{P},\mathsf{Q})w\mathsf l(\mathsf{R},\mathsf{S})^t\). The five compatibility equations give exactly
\[
\mathsf{P}=\begin{pmatrix}a&b\\c&d\end{pmatrix},\quad
\mathsf{Q}=\begin{pmatrix}0&e\\f&g\end{pmatrix},\quad
\mathsf{R}=\begin{pmatrix}h&i\\k&f\end{pmatrix},\quad
\mathsf{S}=\begin{pmatrix}0&e\\-d&m\end{pmatrix},
\]
with exactly \(def\ne0,\qquad D=ad-bc\ne0,\qquad D_R=fh-ik\ne0.\)
This is an open subset of affine 11-space. The remaining entries \(a,h,i,k\) and the value \(di-fm\) are unrestricted.

Rational coverage follows from the rational cell factorization and 
\[
m(A,B,C)=d_5(C)\mathsf l(B^{-1}A,C^{-1}B).
\]

This factorization uses \(F\)-matrices and preserves arbitrary multipliers. The following gauge calculation resolves the redundancy of the coefficient tuples.

\subparagraph{Full parabolic intersection}
For \(w=(1,8,-9,5,4,-7,-6,-2,3,10)\), put
\(I_w=P\cap wP^tw^{-1}\), \(k^\natural=(w^{-1}kw)^t\), and
\(V=U\cap wU^tw^{-1}\). Appendix~\ref{m6:s:cert:p027-intersection} gives
the complete 22-dimensional intersection, its two projection arrays and
the ten-dimensional kernel \(V=\exp(\mathfrak v)\).
The full group is \(T^6\ltimes\exp(\mathfrak n_{16})\); the fixed
triangular order proves rational coverage and excludes extra components.

\subparagraph{rational gauge and invariant}

Appendix~\ref{m6:s:cert:p027-normalization} gives the overlap action,
its rational triangular normalization and its three right-root corrections.
They send every compatible tuple to
\[
\mathsf P=I,\quad \mathsf Q=s=\begin{pmatrix}0&1\\1&0\end{pmatrix},
\quad \mathsf R=\operatorname{diag}(t,1),\quad \mathsf S=-j,
\qquad t=eDD_R\in F^\times.
\]
Only the forced units \(d,e,f,D\) are divided by, and the initial
triangular element has multiplier \(1/(ef)\). Thus the branches
\(a=h=0\), \(di=fm\), arbitrary free zero entries and all square classes
of \(t\) are covered. The certificate verifies invariance under the full
six-dimensional torus and all sixteen radical roots; equality of \(t\)
is necessary and sufficient for rational orbit equivalence.

\subparagraph{full stabilizer and character-zero locus}

For every tuple the complete stabilizer is represented uniquely by 
\[
(h_1,h_2)=(lkl^{-1},r(k^\natural)^{-1}r^{-1}),\qquad k\in I_w,
\]

with both original \(H\)-membership conditions. At the canonical tuple, the left blocks \(A,B,C\) have lower/upper/lower triangular shapes. Membership gives \(A=B=sCs\); hence all are diagonal. The right blocks have the same three shapes. Since \(\mathsf{R}=\operatorname{diag}(t,1)\), right membership gives \(A'=\mathsf{R}^{-1}B'\mathsf{R}\), again forcing all three diagonal. Their remaining equalities are exactly
\[
A=B=\operatorname{diag}(\lambda,\mu),\qquad
C=\operatorname{diag}(\mu,\lambda).
\]
Dividing by the unique such torus element leaves both Levi projections equal to the identity, hence leaves exactly \(V\). Conversely these elements satisfy the original ten-by-ten stabilizer equation. Therefore
\[
\mathcal A_{x_t}\simeq T^2\ltimes V,\quad
T^2=\{m(\operatorname{diag}(\lambda,\mu),
\operatorname{diag}(\lambda,\mu),\operatorname{diag}(\mu,\lambda))\}.
\]
Its original reductive projection is \((a,a^{-1})\), \(a=\operatorname{diag}(\mu,\lambda)\). This proves the full group has dimension 12, is split and algebraically connected, and has trivial component group. Rational normalization transfers it to every compatible point.

The all-parameter solution of the membership equations is recorded in
Appendix~\ref{m6:s:cert:p027-elimination}. Its only divisors are the mandatory
\(d,e,f,D\ne0\); in particular \(di=fm\) introduces no jump.

The character functional is identically zero even on the entire intersection:
\[
\operatorname{tr}(\mathsf{P}N_{12}+\mathsf{Q}N_{23}+\mathsf{R}N'_{12}+\mathsf{S}N'_{23})=0.
\]
Group triviality has a separate justification: the canonical torus is in \(H_0\times H_0\), where the character is one. For each of the ten generators of \(V\), finite-exponential matrices satisfy both \(H\)-membership conditions and the stabilizer equation, and their adjacent-block additive character coordinates sum to zero. Finite exponential and the adjacent-block additive homomorphism give trivial character on all of \(V(F)\). Thus the entire previously \(U\times U\)-compatible locus is full-character-compatible.

\subparagraph{modular and density characters}

The ten kernel weights in Appendix~\ref{m6:s:cert:p027-weights} consist of
six weights \(1\), two \(\lambda/\mu\), and two \(\mu/\lambda\).

The product is one, so the full semidirect stabilizer is unimodular.
Since \(\mathcal A\) is also unimodular, its relative orbital line is exactly \(\Delta_{\mathcal A}|_{\mathcal A_{x_t}}\Delta_{\mathcal A_{x_t}}^{-1}=1.\)
The full frame torus has the following coefficient-density determinants:
\[
\det\operatorname{Ad}(k)|_{\mathfrak v}
=\beta^6\gamma^6\delta^2\epsilon^{-8}\zeta^{-6},
\qquad
\det\operatorname{Ad}(k)|_{\mathfrak n_{16}}
=\beta^8\gamma^8\delta^3\epsilon^{-11}\zeta^{-8}.
\]
Direct computation on the unipotent orbit quotient gives
\[
\frac{\det\operatorname{Ad}(k)|_{\mathfrak u}
\det\operatorname{Ad}((k^\natural)^{-1})|_{\mathfrak u}}
{\det\operatorname{Ad}(k)|_{\mathfrak v}}
=\beta^8\gamma^8\delta^3\epsilon^{-11}\zeta^{-8}.
\]
Its absolute value is the forward Jacobian; the density pullback and group modular function use its inverse. On the stabilizing torus,
\((\beta,\gamma,\delta,\epsilon,\zeta)=(\mu,\lambda,\mu,\mu,\lambda)\),
so this equals one. The nontrivial raw line has therefore been retained until the legitimate homogeneous reduction.

\subparagraph{idempotent conclusion: action, fibers and transpose}

Set \(l_0=\mathsf l(I,s)\), \(r_t=\mathsf l(\operatorname{diag}(t,1),-j)\), \(x_t=l_0wr_t^t\). Multiplication gives
\[
\begin{array}{c|rrrrrrrrrr}
(i,j)&(1,7)&(7,1)&(2,2)&(3,5)&(5,3)&(4,10)&(10,4)&(6,8)&(8,6)&(9,9)\\\hline
(x_t)_{ij}&-1&-1&-t&1&1&-1&-1&1&1&-t^{-1}.
\end{array}
\]
In particular \(x_t=x_t^t\), also for nonsquare \(t\).

The surjective action map \(\mathcal A\times F^\times\to C_{27}\) has the complete stabilizers above. Quotienting by the normal subgroup \(U\times U\) yields for each \(t\) \(K_0/T_\Delta,\qquad K_0=GL_2\times GL_2,\qquad T_\Delta=\{(d,d^{-1}):d\text{ invertible diagonal}\}.\)
The unipotent fibers are the 34-dimensional spaces \((U\times U)/V_t\). Split unipotent root coordinates give rational local sections and character-weighted Haar factors. The stabilizing torus has trivial determinant on both the full unipotent pair and \(V_t\), so this extraction line is trivial. Fixed coordinates on \(V\), transported through \(l_0,r_t\), vary analytically with \(t\); thus they also trivialize the line across the parameter family.

The map \([(a,b)]\longmapsto(z=ab^t,E=aE_{11}a^{-1})\)
is an exact \(F\)-rational isomorphism \(K_0/T_\Delta\simeq GL_2\times\mathcal E,\qquad \mathcal E=\{E:E^2=E,\operatorname{rank}E=1\}.\)
Given \(E\), a nonzero column of \(E\) and a nonzero column of \(I-E\) lie in its distinct eigenspaces and form an \(F\)-basis \(a\). Put \(b=z^ta^{-t}\). All choices differ exactly by \(T_\Delta\). The finitely many column choices give rational analytic sections. This establishes rational surjectivity, identifies the fibers, and
gives smooth local quotient charts.

The precise action and transpose are \((a_0,b_0):(z,E)\mapsto(a_0zb_0^t,a_0Ea_0^{-1}),\qquad \tau(z,E)=(z^t,z^tE^tz^{-t}).\)
Frobenius on \(z\in GL_2\) has the global rational section \((z,I)\), stabilizer \((a,a^{-t})\), and gives conjugation on \(\mathcal E\). Both groups in this Frobenius step are unimodular, the base tangent conjugation determinant is one, and transpose has absolute determinant one. Thus the resulting character and density are trivial and the induced involution is \(E\mapsto E^t\).

The residual problem has dimension \(34+4+2+1=41\). Every rank-one idempotent has characteristic polynomial \(T(T-1)\). Its eigenvalues are distinct and in \(F\), and \(\mathcal E(F)\) is one rational conjugacy orbit with full split diagonal-torus stabilizer. Nilpotent, repeated-root, zero-rank and nonsplit residual branches are excluded by the exact isomorphism. The passive invariant \(t\in F^\times\) retains every square class.

\subparagraph{Closedness, extension and localization}

The \(U\times U\)-stabilizer Lie algebras form a constant-rank kernel bundle on the parabolic cell, with rank ten. Restricting the fixed additive character functional gives a regular section of its dual. Its zero set \(C_{27}\) is therefore closed in the whole intrinsic cell.

The parabolic parameter map is a surjective submersion with smooth 22-dimensional fibers. Its compatible preimage is the 52-dimensional acting group times the smooth 11-variable coefficient family. Local quotient charts give smooth \(C_{27}\) of dimension \(52+11-22=41\). The fully gauge-invariant \(t=eDD_R\) descends continuously and analytically. Its derivative in \(h\) on the coefficient chart is \(efD\ne0\), so it is a submersion. Every fiber is the single rational orbit already classified, of dimension \(52-12=40\). A fiber is closed in \(C_{27}\), hence in the cell. The action map has forty orbit directions and the additional parameter direction, making it a surjective analytic submersion; the homogeneous quotient topology is therefore the orbit topology.

The incompatible open has even-character vanishing by the unipotent
criterion; the support sequence identifies the remaining distributions
with their unique antecedents on \(C_{27}\). The invariant
\(t:C_{27}\to F^\times\) is fixed by the extended action. Each fiber is
the rational homogeneous orbit with trivial full character and density;
transpose fixes \(x_t\), so Lemma~\ref{m8:s:c:transporter} makes its Haar line
even. Proposition~\ref{cf:localization} and closed-support injection give
\(D'_{\mathrm{odd}}(\mathcal P_{27})=0\).

For completeness, at \(h_1x_th_2^t\) the transporter \((h_2h_1^{-1},h_1h_2^{-1})\)
has character one by cancellation.

Each compatible orbit carries a nonzero even equivariant Haar
distribution. Because the orbit is closed in the intrinsic cell,
its closed pushforward gives a nonzero even distribution on the cell.

These calculations verify the middle-family rows of Proposition~\ref{m6:orbit-verification}.

\subsection{Full orbit calculations in the rank-one sector}\label{m6:supp:sector-orbits}

\paragraph{The cell \texorpdfstring{$\mathcal P_{31}$}{P31}}\label{m6:s:sec:sector-p031}

\subparagraph{Entire coefficient family}

The five compatibility equations give exactly

\[
\mathsf{P}=\begin{pmatrix}a&0\\b&c\end{pmatrix},\quad
\mathsf{R}=\begin{pmatrix}d&0\\e&-c\end{pmatrix},\quad
\mathsf{Q}=\begin{pmatrix}f&-h\\i&g\end{pmatrix},\quad
\mathsf{S}=\begin{pmatrix}j&g\\k&h\end{pmatrix},
\]

with \(acd(fg+hi)(jh-gk)\ne 0\). This is an eleven-dimensional smooth determinant-open subset of an affine linear space. Either \(g\) or \(h\) may vanish; the normalization below applies to both cases. Every choice is realized by \(\mathsf l(\mathsf{P},\mathsf{Q}) W_{31} \mathsf l(\mathsf{R},\mathsf{S})^t\), and the identity \(m(A,B,C)=d_5(C)\mathsf l(B^{-1}A,C^{-1}B)\) retains every original multiplier.

\subparagraph{Full intersection and gauge}

Appendix~\ref{m6:s:cert:p031-intersection} gives the 22-direction intersection
and its rational projection section. Its image has dimension twelve
and kernel \(\exp(\mathfrak v_{31})\) dimension ten. The projections
have \(A,B,D,E\) lower triangular, \(C\in\mathrm{GL}_2\), and
\(F=C^t\); all linked diagonals and all four lower entries are retained
there. Necessity from the two zero patterns and sufficiency of the
section determine the full group.

\subparagraph{rational invariant and all stabilizers}

In \eqref{m6:s:eq:sector11-gauge}, the diagonal entries transform as \(p_{11}\mapsto (a_1/b_1)p_{11}\), \(r_{11}\mapsto (b_1/a_1)r_{11}\); hence \(\rho=p_{11}r_{11}\) is invariant under \emph{every} gauge.

The rational operations in Appendix~\ref{m6:s:cert:p031-normalization}
produce the canonical tuple
\[
\mathsf P=\mathsf Q=I,\qquad
\mathsf R=\operatorname{diag}(\rho,-1),\qquad
\mathsf S=M=\begin{pmatrix}0&1\\1&0\end{pmatrix}.
\]
Only forced units are divided by. Thus \(g=0\), \(h=0\), Jordan or
nonsplit original coefficient matrices are all included. Equality of \(\rho\) is sufficient
over \(F\) as well as necessary.

The full group stabilizer at \eqref{m6:s:eq:sector-31n} is obtained by imposing

\begin{equation}
\label{m6:s:eq:sector-h}
B\mathsf{P}=\mathsf{P}A,\quad C\mathsf{Q}=\mathsf{Q}B,\quad E\mathsf{R}=\mathsf{R}D,\quad F\mathsf{S}=\mathsf{S}E.
\end{equation}

Appendix~\ref{m6:s:cert:p031-stabilizer} solves these equations for every
\(\rho\ne0\).
The full stabilizer is \(\mathbb G_m\times N_{11}\), dimension twelve,
with trivial algebraic component group. Its unipotent Lie algebra is \(\mathfrak v_{31}\oplus F(A_{21}+B_{21}+C_{21}+\rho S_h-T_h).\)
Rational normalization transfers this full group to every compatible point.

\subparagraph{Full character and density}

For Z in the complete intersection set \(Z'=-(w^{-1}Zw)^t\). Evaluating

\[
\Lambda(Z)=\operatorname{tr}(\mathsf{P}Z_{12}+\mathsf{Q}Z_{23})+
\operatorname{tr}(\mathsf{R}Z'_{12}+\mathsf{S}Z'_{23})
\]

on the 22 basis vectors gives exactly the following potentially nonzero values:
\[
\begin{aligned}
\Lambda(X_{21})&=p_{12},&\Lambda(X_{22})&=p_{22}+r_{22},\\
\Lambda(Y_{21})&=q_{12}+s_{22},&\Lambda(Y_{22})&=q_{22}-s_{12},\\
\Lambda(W_{22})&=-r_{12}.
\end{aligned}
\]
Every other basis value is zero, and all five displayed values vanish on
the coefficient family. Thus the algebraic additive character vanishes
on the full stabilizer Lie algebra. It vanishes directly on the scalar
subgroup. On a finite exponential its value is the linear coefficient
\(\Lambda(N)\), by the polynomial identity for an algebraic homomorphism
to \(\mathbb G_a\). Therefore the algebraic additive homomorphism itself
is zero on the full group, and the character is one on every rational point.

Since the scalar factor is central, the stabilizer and \(H^2\) have modulus one. The full relative orbital density character is exactly one.

\subparagraph{Closedness, smoothness, transpose and localization statement}

The parabolic orbit map \(q:P(F)^2\longrightarrow \mathcal P_{31}(F)\) is surjective by the rational classification. Its stabilizer has the dimension 22 and is smooth; characteristic-zero orbit maps provide local analytic sections. Thus q is an open quotient. Its entire inverse image of the compatible locus is \(H^2\) times the closed five-equation coefficient family. Consequently \(C_{31}\) is closed. Smoothness descends along this smooth quotient, giving dimension \(52+11-22=41\).

The regular invariant \(\rho\) descends through the coefficient presentation; equivalently it is regular on the local sections. A fiber is closed in \(C_{31}\), which is closed in the cell. The rational classification identifies the fiber with one rational orbit of dimension \(52-12=40\). Thus every orbit \(O_\rho\) is closed in the whole cell. The differential of \(\rho=ad\) is nonzero, since both \(a\) and \(d\) are nonzero. In particular every invariant fiber is smooth.

The matrix identity \(W_{31}^t=\mathsf l(I,D_\tau)W_{31}\), with \(D_\tau=\operatorname{diag}(1,-1)\), holds. It induces \((\mathsf{P},\mathsf{Q},\mathsf{R},\mathsf{S})\mapsto(D_\tau \mathsf{R} D_\tau,\mathsf{S} D_\tau,\mathsf{P},\mathsf{Q})\) and preserves \(\rho\). A rational transporter follows from the \emph{full} normalization and its group section. If \(\tau(x)=a\cdot x\), swapping the two factors gives \(\theta(a)a\in \operatorname{Stab}(x)\). The algebraic additive character \(\lambda \) satisfies \(\lambda \circ\theta=\lambda \) and vanishes on the entire group stabilizer by the character calculation above. Therefore \(2\lambda (a)=0\) in \(F\), hence \(\lambda (a)=0\) because \(\operatorname{char}F=0\). This holds in every characteristic-zero field considered here, including dyadic \(F\).

The character-one transporter and positive density give an even
orbital line by Lemma~\ref{m8:s:c:transporter}. Every fiber of
\(\rho:C_{31}\to F^\times\) is the closed rational orbit above and the
extended group fixes the base. Proposition~\ref{cf:localization}, followed
by the support sequence after removing the incompatible open, proves
intrinsic odd vanishing on \(\mathcal P_{31}\).

\paragraph{The cell \texorpdfstring{$\mathcal P_{33}$}{P33}}\label{m6:s:sec:sector-p033}

\subparagraph{Entire family, intersection and gauge}

Put \(K=-j\). The six compatibility equations give

\[
\mathsf{P}=\begin{pmatrix}a&b\\c&d\end{pmatrix},\quad
\mathsf{Q}\in\mathrm{GL}_2,\quad \mathsf{R}=\begin{pmatrix}i&-b\\j&-d\end{pmatrix},\quad \mathsf{S}=K\mathsf{Q},
\qquad\det \mathsf{P}\det \mathsf{Q}\det \mathsf{R}\ne0.
\]

The coefficient family is smooth of dimension ten. The vector
\(v=\mathsf P e_2\) is nonzero; both \(b=0\) and \(d=0\) are permitted.
The full intersection has dimension 26. Its rational projection image
has dimension twelve and its unipotent kernel
\(\exp(\mathfrak v_{33})\) dimension fourteen. The exact projections are
\(A\) lower triangular, \(B,C\in\mathrm{GL}_2\), and
\[
D=\begin{pmatrix}a_1&0\\d_{21}&\delta/a_2\end{pmatrix},\quad
E=\delta B^{-1},\quad F=C^t,\qquad\delta=\det C.
\]
Appendix~\ref{m6:s:cert:p033-intersection} gives the root basis, rational
section, and direct verification of the inverse and signs in \(E\).

\subparagraph{A normal form for every rational point}

Since \(\mathsf R e_2=-\mathsf P e_2\),

\[
\mathsf{P}^{-1}\mathsf{R}=\begin{pmatrix}u&0\\v&-1\end{pmatrix},\quad u\ne0.
\]

Take \(A=I\), \(B=\mathsf P\), \(C=\mathsf Q\mathsf P\),
\(\delta=\det(\mathsf Q\mathsf P)\), and \(D_{21}=-\delta v\).
Substitution in the full gauge gives

\[
(\mathsf{P},\mathsf{Q},\mathsf{R},\mathsf{S})\mapsto(I,I,\operatorname{diag}(\kappa,-1),K),\quad\kappa=\delta u\ne0.
\]

Next take \(A=B=C=\operatorname{diag}(1,\kappa^{-1})\) and \(D=I\).
This is a legal projection point because \(\delta=\kappa^{-1}\) and
\(\delta/a_2=1\). It sends \(\kappa\) to \(1\). Thus every compatible
\(F\)-point belongs to the orbit of

\[
x_0=W_{33}\mathsf l(\operatorname{diag}(1,-1),K)^t.
\]

These identities apply to all the original variables. The cases \(\mathsf P=I\) and \(\mathsf P=M\)
illustrate the two extreme branches \(b=0\) and \(d=0\), respectively;
the same formulas cover both and all intermediate points. The construction supplies a rational transporter for every input, with arbitrary nonzero
\(\kappa\) normalized by a diagonal matrix already over \(F\).

\subparagraph{Full stabilizer and density}

Solving the four group membership equations at \(x_0\) gives the full
stabilizer \(\mathbb G_m\times N_{15}\), dimension sixteen, with trivial
component group throughout the family. Its unipotent Lie algebra is \(\mathfrak v_{33}\oplus F(S_h-A_{21}-B_{21}-C_{21}).\)
Appendix~\ref{m6:s:cert:p033-stabilizer} gives both the full group solution
and the all-parameter generator; its denominators are only determinants,
so either entry of \(\mathsf P e_2\) may vanish.

On the 26 intersection directions, the potentially nonzero character
values are precisely
\[
\begin{aligned}
\Lambda(X_{21})&=p_{12}+r_{12},&\Lambda(X_{22})&=p_{22}+r_{22},\\
\Lambda(Y_{11})&=q_{11}+s_{21},&\Lambda(Y_{12})&=q_{21}-s_{11},\\
\Lambda(Y_{21})&=q_{12}+s_{22},&\Lambda(Y_{22})&=q_{22}-s_{12}.
\end{aligned}
\]
All six expressions vanish on the entire coefficient family, and every
other basis value is zero. The scalar and finite-unipotent argument
therefore gives algebraic character zero and \(\eta=1\) on the full
rational stabilizer. Both stabilizer and acting group are unimodular,
so the relative orbital density character is exactly one.

\subparagraph{Closed smooth orbit, extraction, and transpose}

The entire inverse image of the compatible set under the open smooth parabolic quotient is \(H^2\) times the closed ten-variable compatible family. It is saturated by the complete gauge. Hence the compatible set is closed in \(\mathcal P_{33}\). Rational transitivity makes it the single smooth orbit, of dimension \(52-16=36=52+10-26\). The cell has dimension 42.

The intersection with \(U^2\) is exactly \(\exp(\mathfrak v_{33})\),
dimension 14. For the full stabilizer projection to
\(K_0=\mathrm{GL}_2^2\), take \(\Gamma=S_h-A_{21}-B_{21}-C_{21}\).
Its left \(H_0\) Lie projection is \(-E_{21}\). After conjugation by
the right coefficient lift \(r\), its right \(H_0\) Lie projection is
\(E_{12}\). Renaming the parameter gives precisely \(J_0=\{(c(I+tE_{21}),c^{-1}(I-tE_{12})):c\in F^\times,t\in F\}.\)

This sign and the \(E_{12}\) direction were checked by multiplication as ten-by-ten conjugation identities. Every element lifts by the finite exponential. The homogeneous map \(C_{33}\longrightarrow K_0/J_0\) has the \(U^2\)-orbits as its fibers, each of dimension 30; the base has dimension 6. The action maps are analytic submersions. Weighted Haar extraction is legitimate because the character is trivial on the fiber stabilizer; scalar conjugation is trivial and unipotent conjugation has determinant one on the relevant unipotent Lie spaces. Thus the residual density line at \(J_0\) is trivial.

The full intermediate stabilizer is \(J_0=\mathbb G_m\times\mathbb G_a\). Direct transpose settles the distribution problem on the extracted homogeneous quotient.

The identity \(x_0^t=d_5(M)x_0d_5(M)^t\) follows by multiplication. Both translators lie in \(H_0\), their characters are one, and the corresponding extended stabilizer element is an involution with odd character -1. Its positive density scalar is one. Therefore the homogeneous compatible orbit has zero odd distribution space. Since it is closed in the cell and the incompatible complement is open, the open restriction followed by closed-support injection yields \(D'_{\mathrm{odd}}(\mathcal P_{33})=0\).

\subsection{Full orbit calculations on the outer interval}\label{m6:supp:outer-orbits}

\paragraph{Group equations, densities, and rational homogeneous spaces}

For a signed representative \(w\), set \(I_w=P\cap wP^t w^{-1}\). At \(x=lwr^t\), its full stabilizer is in bijection with those \(q\in I_w\) for which

\begin{equation}
\label{m6:s:eq:late-c1}
\left(lql^{-1},\ r(w^{-1}q^{-1}w)^tr^{-1}\right)\in H\times H.
\end{equation}

These are group equations. On either side, the equalities of the first three Levi blocks force equality of all five, by the similitude relation. Membership can equivalently be tested using \(q\), rather than \(q^{-1}\), on the right, since \(H\) is a group. This distinction fixes the sign in the character after differentiation.

For the solvable intersections, the two coordinate-flag conditions
force a common triangular matrix pattern. Their diagonal entries give
an arbitrary split torus and the diagonal-normalized remainder is
unipotent. Lemma~\ref{m6:s:lem:finite-exponential} identifies this entire
remainder using the imposed flag, symplectic, and intertwining
equations. The $\mathcal P_{53}$ intersection has a $\mathrm{GL}_2$ Levi
block and is treated separately below, before the scalar reduction of
its full stabilizer.

Conjugation by \(d_5(a)\) on each \(M_2\) or \(\mathfrak{sl}_2\) coordinate has determinant one. Hence \(H\) and \(\mathcal A\) are unimodular. We compare characters with the modular ratio

\begin{equation}
\label{m6:s:eq:late-c2}
\Delta_{\mathcal A}|_{\mathcal A_x}\Delta_{\mathcal A_x}^{-1}.
\end{equation}

We use the conventions of \S\ref{m6:s:sec:setup} for the parabolic
determinant \(\delta_P\) and the group modular function, keeping
preliminary extraction determinants separate from
\eqref{m6:s:eq:late-c2}. For the surviving compatible orbits below, all full-stabilizer
ratios are one. Any absolute constant from a coordinate chart fixes
the normalization of Haar measure.

\subparagraph{Support and orbital distributions}

The support exact sequence of \cite[(2.2)--(2.3)]{JLX} gives, for a closed invariant \(Z\subset X\), only the needed left-exact equivariant sequence

\begin{equation}
\label{m6:s:eq:late-c3}
0\longrightarrow D'_\chi(Z)\xrightarrow{i_*}D'_\chi(X)
\longrightarrow D'_\chi(X\setminus Z).
\end{equation}

This is \eqref{cf:support-sequence}. The unique closed-support
antecedent and preservation of equivariance and oddness are those of
Lemma~\ref{cf:support}.

We use Proposition~\ref{cf:stabilizer} and its algebraic rational-point
orbit assembly, also for the finite extension by transpose. The initial
incompatible complement is excluded by the $U\times U$ sieve; the
additional P041 exclusion is a row of
Theorem~\ref{m6:s:thm:full-stabilizer-exclusions}.

For a transitive orbit \(O=\mathcal A(F)/K(F)\), if both \(\eta|_K=1\) and \eqref{m6:s:eq:late-c2} are trivial, positive quotient Haar gives

\begin{equation}
\label{m6:s:eq:late-c4}
\mathcal O(f)=\int_{\mathcal A(F)/K(F)} f(a x_0)\eta(a)\,d\mu(a).
\end{equation}

The density criterion \eqref{cf:orbital} makes this line
one-dimensional. The orbit submersion gives its homogeneous-space topology
and locally finite quotient Haar measure; on a closed orbit the integral
pushes forward, and a compact open patch of constant weight proves
nonvanishing.

A character-one transpose transporter makes the orbital line even by
Lemma~\ref{m8:s:c:transporter}. For P040 this applies in every fiber of its
regular complete invariant, so Proposition~\ref{cf:localization} gives
vanishing on the whole compatible family.

\subparagraph{Compatible loci and rational transitivity}

The \(U\times U\) graph kernels have constant dimension on each fixed parabolic cell: \(U\times U\) is normal in \(P\times P\), so the parabolic action transports them. They therefore form an algebraic vector bundle. The additive character differential is a regular section of its dual, and compatibility is its zero locus. It is closed in the cell. Vanishing of this differential is equivalent to triviality of the character on the entire unipotent group, by the finite characteristic-zero logarithm.

The smooth parabolic orbit map \(P^2\to \mathcal P_i\) has fiber \(I_{w_i}\). Its inverse image of the compatible locus is, in the Levi/unipotent coordinates, \(U^2\times \mathrm{GL}_2^2\times Z_i\), where \(Z_i\) is the complete open coefficient family. Its dimension is \(52+\dim Z_i\). Smoothness and geometric irreducibility descend from this smooth surjective base change, giving the compatible-locus dimensions recorded below. This argument uses the entire saturated family.

On the atlas fiber product the two compatible pullbacks have the same
geometric support, since the all-parameter equations are intrinsic over
every characteristic-zero extension. They are smooth and reduced, so their
radical ideals agree. Faithfully flat descent gives the closed reduced
compatible subscheme, and smoothness descends.

For \(\mathcal P_{53}\) and \(\mathcal P_{69}\) every geometric orbit has the full dimension of that geometrically irreducible compatible variety. Each orbit is locally closed; its closure must therefore be the whole irreducible variety, and the orbit is open. Distinct orbits would be disjoint nonempty opens, impossible in an irreducible variety. This proves geometric transitivity only after uniformity at every geometric point has been established.

For the separate rational-transitivity step, a rational base point
has stabilizer $K=\mathbb G_m\times U_n$. Multiplicative Hilbert 90 and
a split additive central series give $H^1(F,K)=1$. Additive Hilbert 90
follows by averaging over a finite Galois splitting extension in
characteristic zero; induction through the series proves the unipotent
assertion. Lemma~\ref{m6:s:lem:rational-homogeneous} then gives one rational
orbit.

\paragraph{The cell \texorpdfstring{$\mathcal P_{40}$}{P40}}\label{m6:s:sec:p040}

\subparagraph{atlas and full intersection}

In the fixed five-block realization write \(w=(3,5,1,4,2,-9,7,10,-6,8),\quad \mathsf l(\mathsf P,\mathsf Q)=m(\mathsf{Q}\mathsf{P},\mathsf{Q},I).\)
The signed-column convention gives \(w^t=w\) and \(w^tJw=J\). Every rational point of the parabolic cell is an \(H(F)\times H(F)\)-translate of \(\mathsf l(\mathsf P,\mathsf Q)w\mathsf l(\mathsf R,\mathsf S)^t\), with all four matrices invertible. This retains both exterior diagonal \(\mathrm{GL}_2\)'s, including all multipliers.

The five compatibility equations are \(p_{11}=r_{11},\ p_{12}=r_{12}=0,\ p_{22}=s_{12},\ q_{12}=r_{22}.\)
Their complete solution, with precisely the invertibility inequalities, is
\begin{equation}
\label{m6:s:eq:late-p40-1}
 \mathsf{P}=\begin{pmatrix}a&0\\b&c\end{pmatrix},\quad
 \mathsf{Q}=\begin{pmatrix}d&e\\f&g\end{pmatrix},\quad
 \mathsf{R}=\begin{pmatrix}a&0\\h&e\end{pmatrix},\quad
 \mathsf{S}=\begin{pmatrix}i&c\\j&k\end{pmatrix},\qquad
 ace\det \mathsf{Q}\det \mathsf{S}\ne0.                                      
\end{equation}
This is an open subset of an eleven-dimensional linear space. None of
\(b,d,g,h,i,j,k\), \(ad+be\), or \(ai+ch\) is assumed nonzero.

The two flag conditions and the similitude equation give the full group
\(I_w=T_6\ltimes\exp(\mathfrak n_{40})\), of dimension eighteen,
with twelve-dimensional radical. Appendix~\ref{m6:s:cert:p040-intersection}
contains the matrix zero pattern, root basis and rational projection
lifts. The entire projection image consists of six lower triangular
blocks \((A,B,C,D,E,F)\) with the diagonal linkage
\eqref{m6:s:eq:late-p40-4} and six independent lower entries. Its kernel is \(U_6=\exp\langle X_{11},X_{21},X_{22},Y_{21},Z_{11},Z_{21}\rangle.\)
The common triangular order and finite logarithm establish all rational
points and exclude extra components.

\subparagraph{Gauge, invariant and divisors}

The exact ambiguity in the four coefficients is
\begin{equation}
\label{m6:s:eq:late-p40-5}
 (\mathsf{P},\mathsf{Q},\mathsf{R},\mathsf{S})\mapsto(B^{-1}\mathsf{P}A,C^{-1}\mathsf{Q}B,E^{-1}\mathsf{R}D,F^{-1}\mathsf{S}E).     
\end{equation}
It follows by multiplying the two Levi representatives by the two Levi projections above and refactoring out their common third blocks into the exterior \(H_0\) factors. Thus it is the ambiguity of rational cell representatives.

The six lower-unipotent operations and the diagonal weights are
Appendix~\ref{m6:s:cert:p040-normalization}. Their first divisions are only
by the mandatory \(c,e\ne0\); thereafter
\(f_0=-\det\mathsf Q/e\) and \(j_0=-\det\mathsf S/c\) are units.
The final rational diagonal choice gives

\begin{equation}
\label{m6:s:eq:late-p40-6}
 \mathsf{P}=\mathsf{R}=I,\qquad \mathsf{Q}=\begin{pmatrix}0&1\\1&0\end{pmatrix},\qquad
 \mathsf{S}=\begin{pmatrix}0&1\\\kappa&0\end{pmatrix},
 \qquad\kappa=a\det \mathsf{Q}\det \mathsf{S}\in F^\times.                 
\end{equation}
All \(\kappa\ne0\) occur. Under the lower-unipotent portion of \eqref{m6:s:eq:late-p40-5}, \(a\), \(\det \mathsf{Q}\), and \(\det \mathsf{S}\) are individually unchanged; under its diagonal portion their weights cancel. Since \eqref{m6:s:eq:late-p40-5} exhausts the overlap, \(\kappa\) is a complete invariant of \(H(F)\times H(F)\)-orbits.

The divisors \(ad+be=0\) and \(ai+ch=0\), separately and together, are fully included. The canonical parameters \(a=c=e=f=1\), \(j=\kappa\), \(b=d=g=h=i=k=0\) lie on both. The normalizations are regular there. Both divisors are therefore included in the same normalization.

\subparagraph{Full stabilizer, phase and relative density}

Put \(E_0=E_{11}\), \(F_0=E_{22}\). Exact multiplication of the normalized Levi matrices and \(w\) gives
\begin{equation}
\label{m6:s:eq:late-p40-7}
 x_\kappa=\begin{pmatrix}
 0&\kappa F_0&E_0&0&0\\
 \kappa F_0&E_0&0&0&0\\
 E_0&0&0&0&F_0\\
 0&0&0&F_0&\kappa^{-1}E_0\\
 0&0&F_0&\kappa^{-1}E_0&0
 \end{pmatrix}.                                          
\end{equation}
It satisfies \(x_\kappa^t=x_\kappa\), \(x_\kappa^tJx_\kappa=J\) and the \(\mathcal P_{40}\) ranks for every \(\kappa\ne0\).

To determine the full stabilizer, use the complete overlap image. Fixing \(\mathsf{P}=\mathsf{R}=I\) forces \(A=B\), \(D=E\). Fixing the exchange matrix \(\mathsf{Q}\) forces \(C=\mathsf{Q}B\mathsf{Q}^{-1}\). As both \(B,C\) are lower triangular, \(B\) must be diagonal. Fixing \(\mathsf{S}\) likewise forces \(E\) diagonal. In \eqref{m6:s:eq:late-p40-4}, the remaining equations leave exactly two arbitrary nonzero diagonal entries. Every lower off-diagonal image coordinate is zero. The preimage is the full connected kernel \(U_6\). These equations describe the full group, including every component
and parameter value.

Let \(l_0=m(\mathsf{Q},\mathsf{Q},I)\), \(V=l_0U_6l_0^{-1}\), and let \(T_2\) be the diagonal torus in \(\mathrm{GL}_2\). The full stabilizer is
\begin{equation}
\label{m6:s:eq:late-p40-8}
 K_\kappa=\left\{
 \left(h,(x_\kappa^{-1}h^{-1}x_\kappa)^t\right):
 h\in V\rtimes d_5(T_2)\right\}.                         
\end{equation}
On the torus the pair is exactly
\[
 \left(d_5(\operatorname{diag}(\alpha,\beta)),
 d_5(\operatorname{diag}(\alpha^{-1},\beta^{-1}))\right).
\]
The unipotent Lie algebra in intersection coordinates consists of
\[
 \mathcal N(X,Y,Z,0,0,0),\quad
 X=\begin{pmatrix}x_1&0\\x_2&x_3\end{pmatrix},\quad
 Y=\begin{pmatrix}0&0\\y&0\end{pmatrix},\quad
 Z=\begin{pmatrix}z_1&0\\z_2&0\end{pmatrix}.
\]
The only nonzero bracket of its indicated coordinate basis is \([X_{22},Y_{21}]=Z_{21}\). After conjugation by \(l_0\), the left character phase is \(x_1+x_3+y\); the right phase is its negative. These phases follow from the two adjacent blocks in the stabilizer formula. On finite exponentials the adjacent blocks equal their Lie adjacent blocks; the trace-sum character therefore has these exact phases on the whole group. The torus contributes character one. Hence \(\eta|_{K_\kappa}=1.\)
The entire locus \eqref{m6:s:eq:late-p40-1} is full-character-compatible.

In the order \((x_1,x_2,x_3,y,z_1,z_2)\), the torus adjoint weights on \(V\) are exactly 
\[
1,\ \alpha/\beta,\ 1,\ 1,\ \beta/\alpha,\ 1.
\]

Their product is one. The unipotent part and torus are unimodular, so \(\Delta_{K_\kappa}=1\). As \(\Delta_{H\times H}=1\), the \emph{relative} orbital Frobenius density character is exactly one. It works unchanged when \(|2|\ne1\).

Thus every compatible stabilizer has dimension eight and is a split two-torus semidirect a connected six-dimensional unipotent group. The rational normalization and \eqref{m6:s:eq:late-p40-8} give this full group at every compatible point.

\subparagraph{Closedness, smoothness, orbital distributions and oddness}

The dimension of the parabolic cell is \(2\dim P-\dim I_w=68-18=50\). The \(U\times U\) stabilizer dimension is the constant six, so its Lie kernels form an algebraic vector bundle on the cell. The vanishing of the trace-sum functional on that bundle is a closed condition; call its zero set \(C_{40}\). By the preceding full-character calculation this is also the fully compatible locus.

The smooth parabolic atlas pulls \(C_{40}\) back to \eqref{m6:s:eq:late-p40-1}, with the freely varying two exterior \(U\)'s and two exterior \(\mathrm{GL}_2\)'s. This pullback is smooth of dimension \(44+8+11=63.\)
The atlas has smooth eighteen-dimensional fibers. Smooth descent proves that \(C_{40}\) is smooth of dimension 45. 

The regular function \(a\det \mathsf{Q}\det \mathsf{S}\) is invariant under the full overlap relation, so it descends through the smooth atlas to a regular map \(\kappa:C_{40}\longrightarrow\mathbb G_m.\)
This invariance holds on the smooth reduced atlas fiber product, since the same algebraic overlap formulas hold over every characteristic-zero extension. Its fibers are exactly the rational and geometric orbits, by the explicit field-rational normalization. Each fiber is consequently closed in \(C_{40}\), and \(C_{40}\) is closed in \(\mathcal P_{40}\). Every \(\kappa\)-orbit is therefore closed \emph{in \(\mathcal P_{40}\)}, of dimension \(52-8=44\).

At an incompatible point the additive-character criterion applies. At a compatible point, transport from the symmetric point \(x_\kappa\) supplies an odd stabilizer element conjugate to transpose itself. Its extended character is \(-1\), of the character of the transporter: the two transporter factors cancel under conjugation. Its modular character is one, because it has order two. The stabilizer distribution criterion for the extended action therefore yields
\[
                  D'_{\widetilde\eta}(\mathcal P_{40})=0.
\]

Even vanishing is false. On the closed orbit \(O_\kappa\), the integral
\[
 T_\kappa(f)=\int_{(H\times H)/K_\kappa}
                 \eta(a)f(a\cdot x_\kappa)\,d\bar a
\]
is well-defined because \(\eta|K_\kappa=1\) and the relative modular character is one. The quotient is the orbit with its analytic topology: its orbit map is a smooth surjective submersion, hence locally open, and the induced homogeneous-space bijection is a homeomorphism. Quotient Haar is a Radon measure. Because the orbit is closed in the cell, every compactly supported locally constant cell test function restricts to a compactly supported orbit test function. Its closed pushforward is a locally finite nonzero distribution. A test supported where the character is constant gives a nonzero value.

Transpose fixes the base point and preserves the character. On the one-dimensional orbital density line it acts by a positive scalar whose square is one, hence by \(+1\). Thus these are transpose-even distributions. This proves \(D'_\eta(\mathcal P_{40})\ne0\).

\subparagraph{rank-one presentation and the nilpotent test}

For
\[
 Y=\{M\in M_2(F):\operatorname{rank}M=1,\ t=\operatorname{tr}M\ne0\},
\]
Cayley--Hamilton gives \(E=M/t\) idempotent of rank one, and \(F_0=I-E\). Define
\begin{equation}
\label{m6:s:eq:late-p40-9}
 \sigma(M)=\begin{pmatrix}
 0&t^{-1}F_0&E&0&0\\
 t^{-1}F_0&E&0&0&0\\
 E&0&0&0&F_0\\
 0&0&0&F_0&M\\
 0&0&F_0&M&0
 \end{pmatrix}.                                           
\end{equation}
The two eigenspaces of \(M\) are the rational lines for eigenvalues \(t\) and zero. Hence every rational \(M\in Y\) is \(\mathrm{GL}_2(F)\)-conjugate to \(\operatorname{diag}(t,0)\). At that diagonal matrix, \eqref{m6:s:eq:late-p40-9} is exactly \(x_{1/t}\). The formula satisfies \(\sigma(aMa^{-1})=d_5(a)\sigma(M)d_5(a)^{-1},\quad \sigma(M)^t=\sigma(M^t),\quad \sigma(M)_{45}=M.\)
Conjugation by the similitude \(d_5(a)\) preserves the multiplier-one symplectic group. These identities prove the form equation and cell membership for every \(M\). The precise embedded conjugation group is \(a\mapsto(d_5(a),d_5(a^{-t}))\subset H\times H;\)
its character is one. The invariant is exactly \(t=\kappa^{-1}\). At the diagonal normal representative its centralizer is the split diagonal torus; at a general point it is its rational conjugate, and both it and \(\mathrm{GL}_2\) are unimodular, giving relative density one. Transpose fixes the diagonal base point and is positive and involutive on quotient Haar, so its sign is \(+1\).

The variety \(Y\) is smooth of dimension three: it is the nonzero-trace open of the smooth rank-one locus in the determinant hypersurface. The map \((H\times H)\times Y\longrightarrow C_{40},\quad (a,M)\mapsto a\cdot\sigma(M)\)
is surjective on \(F\)-points by \eqref{m6:s:eq:late-p40-6}. Its differential contains all 44 orbit directions; the trace direction contributes a further direction, since \(d\kappa=-t^{-2}\,dt\ne0.\)
It is a submersion onto the smooth 45-dimensional target at every point. This proves saturation and submersion everywhere.

The full stabilizer retains the six-dimensional unipotent group,
whose determinant-one torus action gives trivial orbital density. 

 A trace-zero rank-one point does not occur in this specified residual domain: \(t=\kappa^{-1}\), and \(ace\det \mathsf{Q}\det \mathsf{S}\ne0\) forces \(\kappa\ne0\). However, a literal raw block of an arbitrary unnormalized compatible \(\mathcal P_{40}\) point can be nilpotent. The curve in the next section has \(b(t)_{45}=b(t)_{54}=-E_{21}\qquad(t\ne0).\)
Thus \(b(1)\) is an explicit compatible \(\mathcal P_{40}\) point with a trace-zero rank-one raw (4,5) block. The residual statement refers to the normalized matrix: after the stated normalization, the residual matrix is \(tE_{11}\), of nonzero trace.

\paragraph{The cell \texorpdfstring{$\mathcal P_{53}$}{P53}}\label{m6:s:sec:p053}

\subparagraph{Full-stabilizer equations}

For \(E_w=P\cap wP^t w^{-1},\qquad Z'=-(w^{-1}Zw)^t,\)
the stabilizer Lie equations at such a coefficient point are
\begin{equation}
\label{m6:s:eq:late-so1}
\mathsf{P}Z_{11}=Z_{22}\mathsf{P},\quad \mathsf{Q}Z_{22}=Z_{33}\mathsf{Q},\quad
\mathsf{R}Z'_{11}=Z'_{22}\mathsf{R},\quad \mathsf{S}Z'_{22}=Z'_{33}\mathsf{S}. 
\end{equation}
Here subscripts denote two-by-two blocks. They come from equality of diagonal blocks in both conjugates of H. The same equations with positive \(k^\natural=(w^{-1}kw)^t\) describe group membership, since membership of the second factor may be tested before inversion. The pair associated to \(k\) is
\begin{equation}
\label{m6:s:eq:late-so2}
\left(lkl^{-1},r(k^\natural)^{-1}r^{-1}\right),
\quad l=\mathsf l(\mathsf{P},\mathsf{Q}),\ r=\mathsf l(\mathsf{R},\mathsf{S}).
\end{equation}

The calculated algebraic phase is
\begin{equation}
\label{m6:s:eq:late-so3}
\Lambda(Z)=\operatorname{tr}(\mathsf{P}Z_{12}+\mathsf{Q}Z_{23})
 +\operatorname{tr}(\mathsf{R}Z'_{12}+\mathsf{S}Z'_{23}). 
\end{equation}
This is the differential of the algebraic additive phase underlying \(\eta=\xi\boxtimes\xi\).

\subparagraph{\(\mathcal P_{53}\): compatible family}

The signed representative is \(W_{53}=(10,3,2,-7,5,6,4,-9,-8,-1).\)
Recomputing \eqref{m6:s:eq:late-so3} on the full parabolic intersection gives only
\[
\begin{array}{c|cccccc}
Z&X_{11}&X_{12}&X_{21}&Y_{21}&Y_{22}&S_h\\ \hline
\Lambda(Z)&p_{11}&p_{21}+r_{21}&p_{12}-r_{12}
&q_{12}+s_{22}&q_{22}-s_{12}&r_{11}.
\end{array}
\]
All other allowed root values are zero. The entire solution in the four invertible coefficient matrices is
\begin{equation}
\label{m6:s:eq:late-so4}
\mathsf{P}=\begin{pmatrix}0&b\\c&d\end{pmatrix},\quad
\mathsf{Q}=\begin{pmatrix}f&-j\\h&i\end{pmatrix},\quad
\mathsf{R}=\begin{pmatrix}0&b\\-c&e\end{pmatrix},\quad
\mathsf{S}=\begin{pmatrix}k&i\\l&j\end{pmatrix}, 
\end{equation}
\[
bcDE\ne0,\qquad D=fi+hj,\qquad E=kj-il.
\]
Indeed the four determinants are \(-bc,D,bc,E\). These are the only restrictions. In particular \(d,e,f,h,k,l\) have no individual nonvanishing conditions; \(i=0\) and \(j=0\) are both permitted. The simultaneous condition \(i=j=0\) is already excluded by \(D\ne0\), rather than by an extra pivot choice.

As a nonempty principal open subset of a ten-dimensional affine
linear space, the family is smooth and geometrically irreducible.
Together with rational parabolic factorization, it covers all
compatible rational points. The calculation below accounts for
the overlap group.

\subparagraph{Full stabilizer}
Appendix~\ref{m6:s:cert:p053-full} gives the complete 23-dimensional
intersection \(L_{53}\ltimes\exp(\mathfrak e_{53,u})\), with
\(L_{53}\simeq\mathbb G_m^4\times\mathrm{GL}_2\) and radical
dimension fifteen. Its full membership elimination yields
\[
\mathcal A_x=\{(lz\exp(N_0)l^{-1},\,
rz^{-1}\exp(-(w^{-1}N_0w)^t)r^{-1}):
z\in F^\times, N_0\in\mathfrak n_{53}\}
\simeq\mathbb G_m\times U_{12}.
\]
The basis of \(\mathfrak n_{53}\), all-parameter solution
\eqref{m6:s:eq:late-so8}, and group equations
\eqref{m6:s:eq:late-so9}--\eqref{m6:s:eq:late-so10} are retained there.
The required units are \(b,c,D,E\), so these formulas apply when \(i=0\) or \(j=0\).
This is the full group on every branch, with central scalar, split
unipotent factor, trivial component group and dimension thirteen.

\subparagraph{character, density, closedness and rational transitivity}

The six phase values in the preceding table vanish on the entire family \eqref{m6:s:eq:late-so4}, and every other intersection basis value is zero. Thus the algebraic phase vanishes even before imposing the full \(H\) conditions. On the entire unipotent factor of \eqref{m6:s:eq:late-so10}, the algebraic additive homomorphism consequently vanishes by finite exponentiation; the scalar pair has character one directly. This proves \(\eta|_{\mathcal A_x}=1\) as a group statement.

The scalar factor acts trivially under conjugation, and every unipotent element acts unipotently on both \(\mathfrak a\) and \(\mathfrak a_x\). Hence

\begin{equation}
\label{m6:s:eq:late-so53d}
\Delta_{\mathcal A}|_{\mathcal A_x}\Delta_{\mathcal A_x}^{-1}=1,
\qquad
\left|\det(\operatorname{Ad}(k)|_{\mathfrak a/\mathfrak a_x})\right|=1.
\end{equation}

This computes the full orbit density.

Let \(C_{53}\) be the compatible subset. The complete smooth parabolic atlas has fiber dimension 23, and its compatible inverse image is \(U^2\times \mathrm{GL}_2^2\times Z_{53}\), of dimension \(44+8+10=62\). It is a closed subset inside the invertible source coordinates and is smooth and geometrically irreducible. Compatibility is saturated under all changes of parabolic representative. The descent argument in the common topology section therefore proves that \(C_{53}\) is closed in \(\mathcal P_{53}\) and smooth and geometrically irreducible of dimension \(62-23=39\). The cell itself has dimension \(68-23=45\).

Every geometric stabilizer has dimension 13, by the uniform equations, so every geometric orbit in \(C_{53}\) has dimension \(52-13=39\). By local closedness, each such orbit is open in this irreducible variety. They cannot be disjoint nonempty opens. Thus the whole compatible variety is one geometric orbit.

The rational base point below has full stabilizer \(\mathbb G_m\times U_{12}\). A split central additive filtration of this unipotent group and Hilbert 90 give \(H^1(F,U_{12})=1\) and \(H^1(F,\mathbb G_m)=1\). The explicit transporter-cocycle argument in the common section then gives exactly \(C_{53}(F)=\mathcal A(F)x_{53,0}\). This establishes rational transitivity before using any degeneration.

\subparagraph{Explicit transpose and the unique orbital line}

Put

\[
F_0=\begin{pmatrix}0&1\\1&0\end{pmatrix},\qquad
J_0=\begin{pmatrix}0&1\\-1&0\end{pmatrix},\qquad
(\mathsf{P},\mathsf{Q},\mathsf{R},\mathsf{S})=(F_0,I,J_0,F_0).
\]

The exact product \(x_{53,0}=\mathsf l(F_0,I)W_{53}\mathsf l(J_0,F_0)^t\) has multiplier one and signed columns \((9,3,-7,1,5,6,-10,-4,-8,2).\)

Multiplication gives

\begin{equation}
\label{m6:s:eq:late-so53t}
{d_5(F_0)x_{53,0}d_5(F_0)=x_{53,0}^t.}
\end{equation}

Both factors have \(\xi=1\). For \(x=a x_{53,0}\), the transporter \(\theta(a)(d_5(F_0),d_5(F_0))a^{-1}\) has character one. The positive invariant quotient density changes under transpose by a positive scalar whose square is one; its exact sign is therefore \(+1\).

The rationally surjective orbit submersion has rank 39 and local analytic
sections. Thus \eqref{m6:s:eq:late-c4} is a locally finite nonzero generator
of the orbital line, and closedness of \(C_{53}\) permits pushforward.
The incompatible open has even vanishing; \eqref{m6:s:eq:late-c3} identifies
all cell distributions with this line.
Consequently

\[
{\dim D'_\eta(\mathcal P_{53})=1,\qquad D'_{\widetilde\eta}(\mathcal P_{53})=0.}
\]

For the two boundary branches, \(i=0,j=h=k=1\) gives \(D=E=1\), leaving \(f,l,d,e\) free. Alternatively \(j=0,i=f=1,l=-1\) gives \(D=E=1\), leaving \(h,k,d,e\) free. All permitted free-zero and dyadic values are covered by the same equations.

\subsection{The late scalar-unipotent orbit}\label{m6:supp:late-orbits}
\paragraph{The cell $\mathcal P_{69}$.}\label{m6:s:sec:p069}
The complete compatible coefficient family is
\[
\mathsf P=\begin{pmatrix}0&p\\-b&t\end{pmatrix},\quad
\mathsf Q=\begin{pmatrix}q&0\\d&c\end{pmatrix},\quad
\mathsf R=\begin{pmatrix}0&q\\b&e\end{pmatrix},\quad
\mathsf S=\begin{pmatrix}p&0\\f&g\end{pmatrix},\qquad pbqcg\ne0.
\]
The four free entries $t,d,e,f$ are unrestricted.
Section~\ref{m6:s:cert:p069-full} gives the full intersection of dimension
$25$ and the stabilizer $\mathbb G_m\times U_{15}$, uniformly on this
family. Its adjacent-trace phase is zero on every intersection root:
the potentially nonzero values, in the order
$X_{11},X_{12},X_{21},Y_{11},Y_{21},W_{12},W_{21}$, are
\[
 p_{11},\ p_{21}+r_{21},\ p_{12}-s_{11},\ q_{11}-r_{12},\ q_{12},\ -r_{11},\ s_{12},
\]
and vanish on the displayed family. Scalar--unipotent reconstruction
therefore gives trivial full character and relative orbital density.

The compatible pullback is closed and saturated in the parabolic
coefficient source. The smooth open quotient gives a closed smooth,
geometrically irreducible locus $C_{69}$ of dimension
$44+8+9-25=36$. Every geometric orbit has dimension $52-16=36$;
its local closedness makes it open, so irreducibility implies a single
geometric orbit. For the rational point with coefficients
$(-j,I,M,I)$, the split scalar--unipotent stabilizer has trivial
$H^1$ by multiplicative and additive Hilbert 90. The transporter
cocycle then gives one rational orbit as well. The orbit submersion
supplies its homogeneous-space topology.

With $D_0=\operatorname{diag}(1,-1)$, matrix multiplication gives
\[
 x_{69,0}^t=d_5(D_0)x_{69,0}d_5(D_0)^t.
\]
The two factors have character one and multiplier $-1$.
Consequently the positive orbital line is transpose-even.
Its locally finite measure pushes forward from the closed compatible
orbit; the incompatible complement vanishes by the character criterion.
Thus
\[
 \dim D'_\eta(\mathcal P_{69})=1,\qquad
 D'_{\widetilde\eta}(\mathcal P_{69})=0.
\]

\subsection{The complete support table}
\subsubsection{Global support assembly}\label{m6:s:sec:assembly}

We now prove Theorem~\ref{m6:distribution}. Fix one ordinary distribution
$T\in D'_{\widetilde\eta}(G)\cap\NW(G)$. Every conditional application
in this section concerns this same distribution on $G$.

\paragraph{The opens and the filtration}
Retain $\mathcal O_j,\mathcal S_j,V_i,B_i$ from
\eqref{m6:s:eq:cells-opens}--\eqref{m6:s:eq:cells-boundary}. Put
$V_{i,k}=V_i\cup V_k$ and $B_{i,k}=B_i\cup B_k$ for the paired targets.
Besides $\Omega=\mathcal O_5$, the following determinant open and union of opens occur:
\begin{equation}\label{m6:s:eq:assembly-M}
\begin{split}
\Omega'&=\mathcal P[1,6,34],\\
\mathcal M&=\mathcal O_8\cup\Omega'\cup V_{17}\\
&=\mathcal P[1,\ldots,9,11,14,15,16,17,34].
\end{split}
\end{equation}
The determinant description of $\Omega'$ is given in \S\ref{m6:s:sec:extraction-424}. These are opens in $G$.

For the outer interval set
\begin{equation}\label{m6:s:eq:assembly-interval}
\mathcal B=V_{53}\cap\bigcup_{a,b}\{R_{ab}>(R_{40})_{ab}\},
\qquad \mathcal L=V_{53}\setminus\mathcal B.
\end{equation}
Both $\mathcal B$ and $V_{53}$ are invariant transpose-stable opens.
The interval and boundary are
\begin{equation}\label{m6:s:eq:assembly-interval-members}
\begin{split}
\mathcal L&=\mathcal P[40,41,42,47,48,51,52,53],\\
\mathcal B&=\mathcal P[1,\ldots,18,21,22,23,24,25,29,30,31,34,35,36,39,46].
\end{split}
\end{equation}
In particular $\mathcal B$ differs from the individual rank boundary
$B_{53}$. The intrinsic odd theorem for the entire interval is proved
with its own relative support filtration in Section~\ref{m6:s:sec:outerlate}.

\begin{proposition}\label{m6:s:prop:assembly-filtration}
The sets $\mathcal W_0=\varnothing,\mathcal W_1,\ldots,\mathcal W_{23}=G$
in Table~\ref{m6:s:tab:assembly-filtration} form an increasing sequence of
invariant transpose-stable opens. The third column is exactly
the closed difference $\mathcal W_j\setminus\mathcal W_{j-1}$.
Its twenty-three index sets partition $\{1,\ldots,90\}$.
\end{proposition}
\begin{proof}
Each expression in the second column is a union of opens already
defined by rank inequalities or determinant nonvanishing. Thus it is
open in the original group. The index sets follow by inserting the
rank table in \eqref{m6:s:eq:cells-rank-open}--\eqref{m6:s:eq:cells-boundary};
the complete boundary memberships are recorded in
\S~\ref{m6:s:app:support}. Taking their successive differences gives
the third column. In particular, the cumulative cell counts are
\[
5,6,8,10,14,16,19,21,23,27,28,31,34,38,46,54,55,60,61,82,83,89,90.
\]
Every difference is nonempty, the differences are disjoint, and
their union has all ninety indices. Transpose stability can be seen
both from the open definitions and from the pairs in
\S~\ref{m6:s:app:ledger}. The prefixes used without symmetrization
are $\mathcal O_5,\mathcal O_6,\mathcal O_8,\mathcal O_{26},
\mathcal O_{27},\mathcal O_{30},\mathcal O_{34},\mathcal O_{89}$;
each is closed under that pairing. The other early opens use
$\mathcal S_j$. The sets $\Omega',\mathcal M,\mathcal B,V_{53}$
and the late $B_i,V_i$ in the table are also stable. The difference
of the preceding open from the next one is therefore closed in the
next open, invariant, and transpose-stable.
\end{proof}

Let $S_T=\operatorname{supp}(T)$. In the table, each antecedent has
already vanished, and each last entry is the resulting support inclusion
for this same global $T$. An ``even'' entry uses ordinary equivariance; its relative cell
filtration applies to the individual cells with their ordinary actions. The closed-support passage in every row is
\eqref{cf:support-sequence}.

\begingroup\small\setlength{\tabcolsep}{3pt}\renewcommand{\arraystretch}{1.2}
\begin{longtable}{@{}r >{\raggedright\arraybackslash}p{2.80cm} >{\raggedright\arraybackslash}p{3.45cm} >{\raggedright\arraybackslash}p{2.75cm} >{\raggedright\arraybackslash}p{2.75cm} >{\raggedright\arraybackslash}p{2.18cm}@{}}
\caption{The twenty-three-step support induction.}
\label{m6:s:tab:assembly-filtration}\\
\toprule Step&Current open $\mathcal W_j$&New closed difference&Already-cleared antecedent&Local theorem&Resulting support\\\midrule
\endfirsthead
\toprule Step&Current open $\mathcal W_j$&New closed difference&Already-cleared antecedent&Local theorem&Resulting support\\\midrule\endhead
\midrule\multicolumn{6}{r}{Continued on the next page}\\\endfoot
\bottomrule\endlastfoot
1&$\Omega=\mathcal O_5$&1,2,3,4,5&$\varnothing$&Prop.~\ref{m6:s:prop:model7-open}&$S_T\subset G\setminus\mathcal W_{1}$\\
2&$\mathcal O_6$&6&$\Omega$&Prop.~\ref{m6:rank-dictionary}&$S_T\subset G\setminus\mathcal W_{2}$\\
3&$\mathcal O_8$&7,8&$\mathcal O_6$&Prop.~\ref{m6:s:prop:p007-p008}&$S_T\subset G\setminus\mathcal W_{3}$\\
4&$\mathcal S_9$&9,14&$B_{9,14}\subset\mathcal O_8$&Prop.~\ref{m6:s:prop:p009-spectral}&$S_T\subset G\setminus\mathcal W_{4}$\\
5&$\mathcal S_{11}$&10,11,15,21&$\mathcal S_9$&Prop.~\ref{m6:s:prop:incompatible}; even&$S_T\subset G\setminus\mathcal W_{5}$\\
6&$\mathcal S_{12}$&12,22&$B_{12,22}\subset\mathcal S_{11}$&Prop.~\ref{m6:s:prop:p012-spectral}&$S_T\subset G\setminus\mathcal W_{6}$\\
7&$\mathcal W_6\allowbreak\cup\mathcal M$&16,17,34&$\mathcal W_6$; see \eqref{m6:s:eq:assembly-Bmiddle}&Props.~\ref{m6:s:prop:middle21}, \ref{m6:rank-dictionary}&$S_T\subset G\setminus\mathcal W_{7}$\\
8&$\mathcal S_{17}\allowbreak\cup\Omega'$&13,23&$\mathcal W_7$&Prop.~\ref{m6:s:prop:incompatible}; even&$S_T\subset G\setminus\mathcal W_{8}$\\
9&$\mathcal S_{18}\allowbreak\cup\Omega'$&18,24&$\mathcal W_8$&Thm.~\ref{m6:s:thm:full-stabilizer-exclusions}; even&$S_T\subset G\setminus\mathcal W_{9}$\\
10&$\mathcal O_{26}\allowbreak\cup\Omega'$&19,20,25,26&$\mathcal W_9$&Prop.~\ref{m6:s:prop:incompatible}; even&$S_T\subset G\setminus\mathcal W_{10}$\\
11&$\mathcal O_{27}\allowbreak\cup\Omega'$&27&$B_{27}=\mathcal O_{26}$&Prop.~\ref{m6:s:prop:middle21}&$S_T\subset G\setminus\mathcal W_{11}$\\
12&$\mathcal O_{30}\allowbreak\cup\Omega'$&28,29,30&$\mathcal W_{11}$&Prop.~\ref{m6:s:prop:incompatible}; even&$S_T\subset G\setminus\mathcal W_{12}$\\
13&$\mathcal O_{34}$&31,32,33&$\mathcal O_{30}\cup\Omega'$&Prop.~\ref{m6:s:prop:sector11}&$S_T\subset G\setminus\mathcal W_{13}$\\
14&$\mathcal O_{34}\allowbreak\cup\mathcal B$&35,36,39,46&$\mathcal O_{34}$&Prop.~\ref{m6:s:prop:incompatible}; even&$S_T\subset G\setminus\mathcal W_{14}$\\
15&$\mathcal O_{34}\allowbreak\cup V_{53}$&40,41,42,47,48,51,52,53&$\mathcal B\subset\mathcal W_{14}$&Sec.~\ref{m6:s:sec:outerinterval}&$S_T\subset G\setminus\mathcal W_{15}$\\
16&$\mathcal W_{15}\allowbreak\cup B_{63}$&37,43,44,49,59,60,61,62&$\mathcal W_{15}$&Prop.~\ref{m6:s:prop:incompatible}; even&$S_T\subset G\setminus\mathcal W_{16}$\\
17&$\mathcal W_{15}\allowbreak\cup V_{63}$&63&$B_{63}\subset\mathcal W_{16}$&Thm.~\ref{m6:s:thm:full-stabilizer-exclusions}; even&$S_T\subset G\setminus\mathcal W_{17}$\\
18&$\mathcal W_{17}\allowbreak\cup B_{69}$&54,56,66,67,68&$\mathcal W_{17}$&Prop.~\ref{m6:s:prop:incompatible}; even&$S_T\subset G\setminus\mathcal W_{18}$\\
19&$\mathcal W_{17}\allowbreak\cup V_{69}$&69&$B_{69}\subset\mathcal W_{18}$&Sec.~\ref{m6:s:sec:lateboundary}&$S_T\subset G\setminus\mathcal W_{19}$\\
20&$\mathcal W_{19}\allowbreak\cup B_{84}$&38, 45, 50, 55, 57, 58, 64, 65, 70, 71, 72, 73, 74, 75, 76, 78, 79, 80, 81, 82, 83&$\mathcal W_{19}$&Prop.~\ref{m6:s:prop:incompatible}; even&$S_T\subset G\setminus\mathcal W_{20}$\\
21&$\mathcal W_{19}\allowbreak\cup V_{84}$&84&$B_{84}\subset\mathcal W_{20}$&Thm.~\ref{m6:s:thm:full-stabilizer-exclusions}; even&$S_T\subset G\setminus\mathcal W_{21}$\\
22&$\mathcal O_{89}$&77,85,86,87,88,89&$\mathcal W_{21}$&Prop.~\ref{m6:s:prop:incompatible}; even&$S_T\subset G\setminus\mathcal W_{22}$\\
23&$G$&90&$\mathcal O_{89}$&Prop.~\ref{m6:s:prop:p090}&$S_T\subset G\setminus\mathcal W_{23}$\\
\end{longtable}\endgroup

\paragraph{Verification of the nonautomatic predecessors.}
The table retains all twenty-three predecessors and results; expanded
rank comparisons and row-by-row applications are in
\texttt{anc/model6\_calculations.pdf}, section \texttt{M6-SUPPORT}.
We check here the steps whose support hypothesis is substantive.
The two initial extractions clear $\Omega\cup\Omega'$, hence
$\mathcal O_6=\Omega\cup V_6$; the intrinsic two-cell theorem then
clears $\mathcal O_8$. The spectral boundaries are exactly
\begin{align}
 B_9=B_{14}&=\mathcal P[1,2,6,7]\subset\mathcal O_8,
 \label{m6:s:eq:assembly-B9}\\
 B_{12}&=\mathcal P[1,2,3,4,6,7,8,9,10,11],\quad
 B_{22}=\mathcal P[1,2,3,4,6,7,8,14,15,21].
 \label{m6:s:eq:assembly-B12}
\end{align}
The first paired conormal theorem clears $\mathcal S_9$. Ordinary
equivariant vanishing on $10,11,15,21$, filtered by lexicographic
opens, clears $\mathcal S_{11}$; both second-pair boundaries now lie
there. Applying its conormal theorem to the same original $T$ clears
$\mathcal S_{12}$.
For step seven the complete middle boundaries are
\begin{equation}\label{m6:s:eq:assembly-Bmiddle}
 B_{16}=\mathcal P[1,2,6,7,9,14]=V_9\cup V_{14},\qquad
 B_{17}=\mathcal P[1,2,3,4,6,7,8,9,11,14,15,16].
\end{equation}
First clear $V_{16}$ from its intrinsic theorem and its already-zero
boundary. The only further points of $B_{17}$ outside
$\mathcal O_8\cup V_{16}$ are cells $11,15$, removed at step five;
hence clear $V_{17}$. Gluing with the already-zero $\Omega'$ gives
$T|_{\mathcal M}=0$ for the original global distribution.
The $27$ predecessor is $B_{27}=\mathcal O_{26}$, and step thirteen
uses the proved whole union $\mathcal O_{33}\setminus\mathcal O_{30}$
for $31,32,33$.

For step fourteen $\mathcal B\setminus\mathcal O_{34}$ consists
exactly of the incompatible cells $35,36,39,46$.
Thus $T|_{\mathcal B}=0$ and the intrinsic odd theorem for the entire
interval $\mathcal L$ gives $T|_{V_{53}}=0$ by closed-support injection.
The late boundary differences are
\begin{equation}\label{m6:s:eq:assembly-late-differences}
\begin{split}
 B_{63}\setminus\mathcal W_{15}&=\mathcal P[37,43,44,49,59,60,61,62],\\
 B_{69}\setminus\mathcal W_{17}&=\mathcal P[54,56,66,67,68],\\
 B_{84}\setminus\mathcal W_{19}&=\mathcal P[38,45,50,55,57,58,64,65,70,71,72,73,74,75,76,\\
 &\hspace{4em}78,79,80,81,82,83].
\end{split}
\end{equation}
All these cells are incompatible. After clearing each entire boundary,
use even vanishing on $63,84$ and odd vanishing on $69$ inside $V_i$.
The resulting $\mathcal W_{21}$ excludes $77$;
$\mathcal O_{89}\setminus\mathcal W_{21}$ is exactly
$\mathcal P[77,85,86,87,88,89]$, all incompatible.
Finally the unique odd antecedent on the closed $\mathcal P_{90}$
is zero by Proposition~\ref{m6:s:prop:p090}. This proves
Theorem~\ref{m6:distribution}.
At every intrinsic row, lexicographic relative opens and the ordinary
left-exact support sequence give the antecedent only after its
predecessor is zero. At the spectral rows, the theorem applies to the
original global $T$ with the displayed boundary.

\subsection{Incidence and spectral boundary calculations}\label{supp:m6-boundary-calculations}
To supply the boundary arguments of \S\ref{m6:packets}, the next
three calculations determine the rational gauges, complete group
stabilizers, and relative densities.
\subsubsection{Rational incidence calculations for $\mathcal P_7$ and $\mathcal P_8$}\label{supp:m6-incidence-calculations}

The following calculations supply the rational reductions used in \S\ref{m6:s:sec:p007-p008}. Expanded arrays are also available in \path{anc/model6_calculations.pdf}.

\paragraph{Conventions}\label{m6:s:p78:s1}

We use the fixed five-by-five two-block realization of $G=GSp(J)(F)$,
$J=w_5\otimes j$, $j=\left(\begin{smallmatrix}0&-1\\1&0\end{smallmatrix}\right)$,
and $A^\sharp=\operatorname{tr}(A)I_2-A$. The unrestricted Levi
$m(A,B,C)$ has multiplier $\det C$; its explicit matrix and the root
coordinates $\mathcal N(X,Y,Z,W,T,S)$ are recorded in
\eqref{m6:s:eq:p78-1-2}--\eqref{m6:s:eq:p78-1-3} of
\S\ref{m6:s:cert:p78-notation}. In particular, $\dim U=22$ and
$\xi(ud_5(a))=\psi(\operatorname{tr}X_u+\operatorname{tr}Y_u)$.
The symbols $\mathsf X_{ij},\mathsf Y_{ij},\ldots$ designate precisely
those root matrices; $\mathsf D_{ij}=d_5(E_{ij})$ denotes a Levi Lie
direction, as specified in that calculation.

\subparagraph{Actions, phases, and measures}\label{m6:s:p78:s1-3}

Let

\begin{equation}
\mathcal A=H\times H,
\qquad (h_1,h_2)g=h_1gh_2^t,
\qquad \eta(h_1,h_2)=\xi(h_1)\xi(h_2).
\label{m6:s:eq:p78-1-4}
\end{equation}

Transpose swaps the two factors and the odd extension has
$\widetilde\eta(\tau)=-1$. Our convention is
\begin{equation}\label{m6:s:eq:p78-1-5}
D'_\eta(X)=\{S:S(f\circ k)=\eta(k)^{-1}S(f)\},
\end{equation}
so weighted integration uses $\eta$. Both $H$ and $\mathcal A$
are unimodular, while
\begin{equation}\label{m6:s:eq:p78-1-6}
\delta_P(m(A,B,C))=|\det A|^9|\det B|^5|\det C|^{-14},
\qquad\delta_P|_{H_0}=1.
\end{equation}
With $\int f(xk)\,dx=\Delta_K(k)^{-1}\int f(x)\,dx$,
$\Delta_P=\delta_P^{-1}$.

\subparagraph{Distributional tools}\label{m6:s:p78:s1-4}
We use Theorems~\ref{cf:stabilizer}, \ref{cf:localization}, and
\ref{cf:frobenius}, the left-exact support sequence
\eqref{cf:support-sequence}, and Lemma~\ref{m8:s:c:transporter}.
Frobenius reciprocity retains the relative modular line, also for
nonreductive acting groups and families of unipotent orbital
distributions \cite[Theorem~2.5.7]{AGdescent}. A character-one transpose
transporter acts by $+1$ on an existing orbital line because its measure
scalar is positive and the action is involutive. The entire-plane
Lemma~\ref{m8:s:c:plane}, including its origin, supplies the separate
two-axis argument below.

\paragraph{Common geometric front end and envelopes}\label{m6:s:p78:s2}

\subparagraph{The target rank matrices}\label{m6:s:p78:s2-1}

The rank invariants of Section~\ref{m6:s:sec:cells} give

\begin{equation}
R_7=\begin{pmatrix}
7&6&4&2\\6&5&4&2\\4&4&3&2\\2&2&2&1
\end{pmatrix},
\qquad
R_8=\begin{pmatrix}
7&6&4&2\\6&5&3&2\\4&3&3&2\\2&2&2&1
\end{pmatrix}.
\label{m6:s:eq:p78-2-1}
\end{equation}

Their signed physical-column representatives are

\begin{equation}
\begin{aligned}
W_7&=[1,3,2,5,4,-7,-6,-9,-8,10],\\
W_8&=[1,3,2,-7,5,6,4,-9,-8,10].
\end{aligned}
\label{m6:s:eq:p78-2-2}
\end{equation}

With $\varepsilon=m(I_2,\operatorname{diag}(1,-1),I_2)$, the
signed-column convention gives
\begin{equation}\label{m6:s:eq:p78-2-3}
W_7^t=W_7,\qquad W_8^t=\varepsilon W_8=W_8\varepsilon.
\end{equation}

\subparagraph{The incidence boundary}\label{m6:s:p78:s2-2}

Both cells have $\operatorname{rank}g_{55}=1$ and $R_{33}=3$.
They lie outside both big cells for $Q=(2,6,2)$ and $Q'=(4,2,4)$;
the $\mathcal P_6$ $GL_4$ Shalika coordinates therefore do not
apply to either whole cell.

\subparagraph{\texorpdfstring{A framed \(Sp_4\) incidence slice}{A framed Sp4 incidence slice}}\label{m6:s:p78:s2-3}

Put $w_o=[1,3,2,4,5,6,7,-9,-8,10]$.
The subgroup $G_*\simeq Sp_4$ acts on
$\langle e_4,e_5,e_6,e_7\rangle$, trivially on its orthogonal
complement, and commutes with $w_o$.

Its intersection \(P_*=G_*\cap P\) is the parabolic stabilizing the line \(Fe_4\). Relative to coordinates \((e_4,e_5,e_6,e_7)\), its three opposite-parabolic incidence loci are

\begin{equation}
\begin{array}{ll}
\Sigma_0:&x_{44}\ne0,\\
\Sigma_7:&x_{44}=0,\quad (x_{24},x_{34})\ne(0,0),\\
\Sigma_8:&x_{24}=x_{34}=x_{44}=0.
\end{array}
\label{m6:s:eq:p78-2-4}
\end{equation}

On $\Sigma_8$, invertibility forces $x_{14}\ne0$.
These are the three rational relative positions of an isotropic line
and its opposite. Each rational locus is covered by extending its indicated nonzero pivot
to a symplectic basis with pairing one; the construction is rational over $F$.

The respective representatives in \(G_*\) are the identity and the signed four-column matrices \([2,1,-4,-3],\qquad [-4,2,3,1].\)

Multiplying the last two by \(w_o\) gives exactly \(W_7,W_8\). Consequently the single map

\begin{equation}
\Phi:U\times L\times G_*\times L\times U\longrightarrow G,
\qquad
(u,l,x,r,v)\longmapsto u l w_o x r^t v^t
\label{m6:s:eq:p78-2-5}
\end{equation}

has image

\begin{equation}
\mathcal C=\mathcal P_6\sqcup\mathcal P_7\sqcup\mathcal P_8.
\label{m6:s:eq:p78-2-6}
\end{equation}

The three loci give the three cells. The differential
\eqref{m6:s:eq:p78-11-12} has rank $55$: the added $G_*$ directions
remove every cell conormal except $F\mathsf X_{21}$, with nonzero
pairings $2,2,1$. Absorbing the two $P_*$ factors transports this
test to every source point, including the first locus. Thus $\Phi$
is a surjective analytic submersion onto $\mathcal C$.

\subparagraph{The ambient open and the limitation of the common slice}\label{m6:s:p78:s2-4}

For a symmetric rank matrix \(R_i\), define the open set

\begin{equation}
V_i=\bigcap_{a,b=1}^4\{g:R_{ab}(g)\ge(R_i)_{ab}\},
\qquad i=7,8.
\label{m6:s:eq:p78-2-7}
\end{equation}

Set
\begin{equation}\label{m6:s:eq:p78-2-8}
\mathcal O_7=\Omega\cup V_7,\qquad\mathcal O_8=\Omega\cup V_8.
\end{equation}
The rank classification and $R_6\ge R_7\ge R_8$ give exactly the
first seven and eight cells: dominance implies no later lexicographic
position, and every non-$\Omega$ predecessor dominates the target.
Consequently $\mathcal C=\mathcal O_8\setminus\Omega$ is closed in
$\mathcal O_8$. It is the smooth dimension-$55$ part of the outer
$(r,s)=(2,1)$ incidence locus; locally $\det g_{55}=0$ has nonzero
derivative along $\mathsf X_{21}^tW_i$.
The submersion has varying $U^2$ stabilizers, which
have dimensions two and four on the two cells. The following
quotients retain both different isotropy groups.

\paragraph{rational cell factorization, residual groups, character, and density}\label{m6:s:p78:s3}

\subparagraph{Homogeneous factorization and its rational meaning}\label{m6:s:p78:s3-1}

For \(i=7,8\), let

\begin{equation}
J_i=P\cap W_iP^t W_i^{-1},
\qquad
\mathcal J_i=
\{(k,(W_i^{-1}k^{-1}W_i)^t):k\in J_i\}.
\label{m6:s:eq:p78-3-1}
\end{equation}

Then the rational factorization is

\begin{equation}
\mathcal P_i\simeq(P\times P)/\mathcal J_i,
\qquad (p_1,p_2)\longmapsto p_1W_ip_2^t.
\label{m6:s:eq:p78-3-2}
\end{equation}

The rational parabolic classification represents every $F$-point;
two pairs represent the same point exactly when their ratio belongs
to $\mathcal J_i(F)$.

Put \(\mathcal U=U\times U\), \(G_0=H_0\times H_0\simeq GL_2\times GL_2\), and let \(\mathcal R_i\) be the image of \(\mathcal J_i\) in \(L\times L\). The resulting residual space and group are

\begin{equation}
\mathscr Y_i=(L\times L)/\mathcal R_i,
\qquad G_0\text{ acting by left multiplication}.
\label{m6:s:eq:p78-3-3}
\end{equation}

Both spaces have dimension 12. The residual ordinary character on \(G_0\) is trivial; the unipotent orbital-distribution line carries the calculated density.

\subparagraph{The two unipotent kernels and affine fiber coordinates}\label{m6:s:p78:s3-2}

The single-group kernels are

\begin{equation}
\begin{aligned}
K_7&=U\cap W_7U^tW_7^{-1}
 =\exp\langle\mathsf X_{21},\mathsf Y_{21}\rangle,\\
K_8&=U\cap W_8U^tW_8^{-1}
 =\exp\langle\mathsf X_{21},\mathsf Y_{21},\mathsf Y_{22},\mathsf T_f\rangle.
\end{aligned}
\label{m6:s:eq:p78-3-4}
\end{equation}

The kernel bracket is
\begin{equation}\label{m6:s:eq:p78-3-5}
[\mathsf Y_{21},\mathsf Y_{22}]=-2\mathsf T_f.
\end{equation}
Deleting the kernel roots gives ideals $\mathfrak u_i^\circ$ and
polynomial splittings $U=U_i^\circ\rtimes K_i$. The affine fibers
have dimensions $42,40$; their polynomial inverses and uniqueness,
including the higher $T_{21}$ root, are
\eqref{m6:s:eq:p78-3-7} in \S\ref{m6:s:cert:p78-affine}.

\subparagraph{intersection groups}\label{m6:s:p78:intersections-main}
The intersection calculation \S\ref{m6:s:cert:p78-intersections} gives
\[
\begin{array}{c|c|c|c}
i&\text{reductive part of }J_i&\dim R_u(J_i)&\dim\mathcal P_i\\\hline
7&\mathbb G_m^6&8&54\\
8&\mathbb G_m^4\times GL_2&8&52
\end{array}
\]
Thus $\dim J_7=14$, $\dim J_8=16$ and $\dim\mathcal R_i=12$.
The projected diagonal actions are
\eqref{m6:s:eq:p78-3-9}, \eqref{m6:s:eq:p78-3-11}, with six
independent lower roots for $i=7$ and two on each side for $i=8$.
The root-factor lifts give $\mathcal J_i(F)\twoheadrightarrow
\mathcal R_i(F)$. For $i=7$, $c=t_5t_6$; for $i=8$ the shared
central matrix is $C\in GL_2$, $c=\det C$.

\subparagraph{Character-weighted extraction and its measure normalization}\label{m6:s:p78:s3-4}

Write \(l=m(A_l,B_l,C_l)\), \(r=m(A_r,B_r,C_r)\), and retain

\begin{equation}
\mathsf P=B_l^{-1}A_l,\quad Q=C_l^{-1}B_l,
\qquad R=B_r^{-1}A_r,\quad S=C_r^{-1}B_r.
\label{m6:s:eq:p78-3-13}
\end{equation}

The upright letter \(\mathsf P\) denotes a two-by-two coefficient matrix; it is distinct from the parabolic subgroup \(P\).

On a local chart of \(\mathscr Y_i\), fix Haar measures on \(U_i^\circ,U,K_i\) compatible with \eqref{m6:s:eq:p78-3-6}. The fiber integration can be normalized as

\begin{equation}
\begin{split}
\mathcal J_i f(l,r)
={}&\delta_P(l)\delta_P(r)
\int_{U_i^\circ\times U}
 f(lu^\circ W_i v^t r^t)\\
&\quad\cdot
\psi\!\left(\operatorname{tr}(\mathsf P X_{u^\circ}+QY_{u^\circ})
 +\operatorname{tr}(R X_v+S Y_v)\right)
\,du^\circ\,dv.
\end{split}
\label{m6:s:eq:p78-3-14}
\end{equation}

The graph-stabilizer Haar measure is pushed forward from $W_i$.
Conjugation by $(l,r)$ on $U\times U$ has determinant
$\delta_P(l)\delta_P(r)$, giving the displayed positive prefactor
by Weil quotient integration.

Splitting $u=u^\circ k$ makes the integral independent of the
kernel exactly when its character is trivial, namely
\eqref{m6:s:eq:p78-3-15}
This is the closed zero section of the kernel-character bundle;
off it a one-parameter stabilizer kills distributions.
On it weighted Haar uniqueness and the rational affine charts
give the exact $G_0$ correspondence valued in the orbital line.
Weighted-integral-one test functions prove local surjectivity;
the correspondence commutes with open restriction and
closed-support injection.

\subparagraph{All intermediate density characters}\label{m6:s:p78:s3-5}

Let \(d_i(k)\) be the absolute determinant of conjugation by \(k\in J_i\) on \(\operatorname{Lie}K_i\). On the reductive coordinates above,

\begin{equation}
{
 d_7=\left|\frac{t_2t_4}{t_3t_5}\right|,
\qquad
 d_8=\left|\frac{t_2t_4^4}{t_3c^2}\right|.
}
\label{m6:s:eq:p78-3-16}
\end{equation}

The $i=7$ kernel weights are $t_2/t_3,t_4/t_5$. For $i=8$ the
$X_{21}$ weight is $t_2/t_3$, the second-row $Y$ determinant is
$t_4^2/c$, and the $T_f$ weight is $t_4^2/c$; every unipotent has
determinant one. These products give \eqref{m6:s:eq:p78-3-16}.

The homogeneous-space check is
\begin{equation}\label{m6:s:eq:p78-3-17}
\begin{array}{c|c|c|c}
i&\delta_{J_i}&\delta_P(k_L)\delta_P(k_R)&
\Delta_{P^2}|_{\mathcal J_i}\Delta_{\mathcal J_i}^{-1}\\\hline
7&|t_2^3t_4^4/(t_3^3t_5^4)|&|t_2^4t_4^5/(t_3^4t_5^5)|&d_7^{-1}\\
8&|t_2^3t_4^6/(c^3t_3^3)|&|t_2^4t_4^{10}/(c^5t_3^4)|&d_8^{-1}
\end{array}
\end{equation}
Here $\delta_{J_i}$ is the modulus on its unipotent radical;
the eight root weights are in \S\ref{m6:s:cert:p78-root-densities}.

In the anchor convention \eqref{m6:s:eq:p78-3-14}, a projected intersection
gauge multiplies the fiber measure by $d_i$ and the inverse transporter
character. Its dual relative modular line is $d_i^{-1}$. These are the cell-quotient densities.

\paragraph{Normalization calculus and character bookkeeping}\label{m6:s:p78:s4}

\subparagraph{Removing the two third Levi factors}\label{m6:s:p78:s4-1}

Independent character-one \(H_0\)-translations let us use

\begin{equation}
l=m(A,B,I_2),\qquad r=m(D,E,I_2),
\qquad A=Q\mathsf P,\ B=Q,\ D=SR,\ E=S.
\label{m6:s:eq:p78-4-1}
\end{equation}

We write

\begin{equation}
x_i(A,B,D,E)=m(A,B,I_2)W_i m(D,E,I_2)^t.
\label{m6:s:eq:p78-4-2}
\end{equation}

The representatives are taken modulo the projected intersection just computed; their invariants below are constant under that full gauge group.

Put

\[
f=\begin{pmatrix}0&1\\1&0\end{pmatrix},
\qquad n_-(v)=\begin{pmatrix}1&0\\v&1\end{pmatrix}.
\]

Every projected gauge has the rational $H\times H$ transporter
\eqref{m6:s:eq:p78-4-4}, proved in \S\ref{m6:s:cert:p78-transporters}.
Normalizing the two third Levi factors uses character-one $H_0$
translations.

\subparagraph{Which elementary operations carry a phase}\label{m6:s:p78:s4-3}

There are independent lower-column operations on the first left and first right matrices: \(A\mapsto A n_-(v),\qquad D\mapsto D n_-(w).\)

Their respective transporter characters have exponents

\begin{equation}
{v(E^{-1}D)_{11},\qquad w(B^{-1}A)_{11}.}
\label{m6:s:eq:p78-4-5}
\end{equation}
These positive signs come from the $\mathsf A_{21}$ and
$\mathsf X_{11}$ lifts. Second-factor lower operations, the
third-factor higher-root lifts for $\mathcal P_7$, the shared
$GL_2$ restorations for $\mathcal P_8$, and all diagonal gauges
have phase zero. The complete gauge actions
\eqref{m6:s:eq:p78-4-6}--\eqref{m6:s:eq:p78-4-7}, with the rational
transporters \eqref{m6:s:eq:p78-4-4}, give both existence and the
full equivalence test; the residual branch actions follow below.

\paragraph{\texorpdfstring{\(\mathcal P_7\): stabilizer on the unipotent-compatible family}{P7: stabilizer on the unipotent-compatible family}}\label{m6:s:p78:s5}

\begin{proposition}\label{m6:s:prop:p007-stabilizers}
On the unipotent-compatible locus of \(\mathcal P_7\), the stabilizers have dimensions three, four, or five as described below. They are connected and unimodular. Compatibility of the character on the stabilizer is precisely \eqref{m6:s:eq:p78-5-13}.
\end{proposition}
\par\noindent\textit{Proof.}

In \eqref{m6:s:eq:p78-4-1}, write

\begin{equation}
\wp=(B^{-1}A)_{12}=(E^{-1}D)_{12},
\qquad q=B_{12}=E_{12}.
\label{m6:s:eq:p78-5-1}
\end{equation}

The zero/nonzero conditions on \(\wp,q\) are gauge invariant, although the individual nonzero values before normalization need not be.

\subparagraph{Stabilizer branches}
The eigenline equations \eqref{m6:s:eq:p78-5-2} and determinant
conditions \eqref{m6:s:eq:p78-5-3} give three exhaustive cases.
For $(\wp,q)\ne(0,0)$ the projection is scalar, except when
the first and third lines coincide on both sides; this requires
$\wp q\ne0$ and gives a split two-torus. For $\wp=q=0$ all
three lines coincide and two nilpotent projection parameters
appear. The lines are individually labeled over $F$, excluding
permutation components and nonsplit stabilizer tori.

\subparagraph{Stabilizers away from the all-parallel locus}\label{m6:s:p78:s5-2}

The graph kernel at the normal forms is
\begin{equation}\label{m6:s:eq:p78-5-4}
\{(l e^N l^{-1},r(W_7^{-1}e^{-N}W_7)^t r^{-1}):
N=s\mathsf X_{21}+t\mathsf Y_{21},\ s,t\in F\}.
\end{equation}
Adjoining the central pair $(zI_{10},z^{-1}I_{10})$ gives
\begin{equation}\label{m6:s:eq:p78-5-5}
\mathcal A_x\simeq\mathbb G_m\times\mathbb G_a^2,\quad\dim=3.
\end{equation}
At the plane origin the extra pure reductive pair is
$(d_5(\operatorname{diag}(u,t)),d_5(\operatorname{diag}(u^{-1},t^{-1})))$.
Both kernel weights in \eqref{m6:s:eq:p78-3-16} become one, so
\begin{equation}\label{m6:s:eq:p78-5-6}
\mathcal A_x\simeq\mathbb G_m^2\times\mathbb G_a^2,\quad\dim=4.
\end{equation}
These complete split connected stabilizers are unimodular;
\eqref{m6:s:eq:p78-3-15} kills the kernel character and the
reductive character is trivial. Their orbital density is one.

\subparagraph{The all-parallel locus: enlarged group and nontrivial character}\label{m6:s:p78:s5-3}

On $\wp=q=0$, lower-column operations give the monomial forms
\eqref{m6:s:eq:p78-5-7} with five nonzero parameters
$\alpha,\beta,\gamma,\delta,\epsilon$, after legitimate $H_0$
normalization of the arbitrary multiplier. The whole stabilizer
is $\mathbb G_m\times\operatorname{Heis}_3\times\mathbb G_a$.
The generators, bracket and group law
\eqref{m6:s:eq:p78-5-8}--\eqref{m6:s:eq:p78-5-10} give
$\lambda(K_1)=\lambda(K_2)=0$,
$\lambda(V)=-\alpha/\beta$, $\lambda(W)=-\alpha/\gamma$.

Here \(\lambda(L,R)\) is the sum of the two adjacent-trace Lie characters. Since these are algebraic additive characters on the unipotent factor,

\begin{equation}
\eta\bigl(\exp(sK_1+tK_2+vV+wW)\bigr)
=\psi\!\left(-\frac\alpha\beta v-\frac\alpha\gamma w\right).
\label{m6:s:eq:p78-5-12}
\end{equation}

Since $\alpha/\beta\ne0$, the subgroup $\exp(vV)$ has nontrivial
character at every all-parallel point.

\subparagraph{\texorpdfstring{\(\mathcal P_7\) stabilizer-character compatibility}{P7 stabilizer-character compatibility}}\label{m6:s:p78:s5-4}

The result is \eqref{m6:s:eq:p78-5-13}.

All three stabilizer branches, of dimensions $3,4,5$, are
unimodular and have orbital density character one. The incompatible
all-parallel locus is closed in the unipotent-compatible family and
is removed as a closed-support layer after the other branches.

\hfill\(\square\)\par

\paragraph{\texorpdfstring{\(\mathcal P_7\): exhaustive rational normal forms and their densities}{P7: exhaustive rational normal forms and their densities}}\label{m6:s:p78:s6}

\subparagraph{\texorpdfstring{The branch \(q\ne0\): a global Gaussian normalization}{The branch q is nonzero: a global Gaussian normalization}}\label{m6:s:p78:s6-1}

For $q=B_{12}=E_{12}\ne0$, clear the two diagonal entries of each $B,E$ by left and right
lower operations, leaving the nonzero lower-left pivots
$-\det B/q,-\det E/q$, then rescale diagonally. These zero-phase operations
\eqref{m6:s:eq:p78-6-1} give $B=E=f$ and
\begin{equation}\label{m6:s:eq:p78-6-2}
A_{22}=D_{22}=p.
\end{equation}

The complete remaining gauge transformations are independent lower-column operations on \(A,D\) and

\begin{equation}
{
A\mapsto\operatorname{diag}(h,1)A\operatorname{diag}(r,1),
\qquad
D\mapsto\operatorname{diag}(h^{-1},1)D\operatorname{diag}(r^{-1},1),
\quad h,r\in F^\times.
}
\label{m6:s:eq:p78-6-3}
\end{equation}

\S\ref{m6:s:cert:p78-p7-gauge-test} proves that this is the complete
entire remaining action.

\begin{equation}
d_7=1,
\qquad \Delta_{P\times P}|_{\mathcal J_7}\Delta_{\mathcal J_7}^{-1}=1.
\label{m6:s:eq:p78-6-4}
\end{equation}

This is the density calculation relevant to the plane conclusion.

\subparagraph{\texorpdfstring{The branch \(q\ne0,\ p\ne0\): the whole opposite-weight plane}{The branch q is nonzero, p is nonzero: the whole opposite-weight plane}}\label{m6:s:p78:s6-2}

Use lower-column operations to set \(A_{21}=D_{21}=0\). Then use \(r\) in \eqref{m6:s:eq:p78-6-3} to normalize \(\det A=1\). The normal forms are

\begin{equation}
{
X_7(p,k,b,e)=x_7\left(
\begin{pmatrix}p^{-1}&b\\0&p\end{pmatrix},\ f,
\begin{pmatrix}k/p&e\\0&p\end{pmatrix},\ f\right),
\quad p,k\in F^\times,\ b,e\in F.
}
\label{m6:s:eq:p78-6-5}
\end{equation}

Every point on this branch reaches the displayed family rationally,
and every allowed parameter gives invertible Levi matrices. Its only
remaining effective equivalence is

\eqref{m6:s:eq:p78-6-6}

It is realized on the original group with character one by

\begin{equation}
\left(d_5(\operatorname{diag}(1,h^{-1})),
      d_5(\operatorname{diag}(1,h))\right).
\label{m6:s:eq:p78-6-7}
\end{equation}

Here $h$ ranges over all $F^\times$ and transpose is
\begin{equation}\label{m6:s:eq:p78-6-8}
X_7(p,k,b,e)^t=X_7(p,k,e,b).
\end{equation}
The reciprocal first-column renormalization is an intersection
gauge of phase and density one. The pre-normalization invariants
are \eqref{m6:s:eq:p78-6-9}; in particular the two axes are
$A_{12}=0$ and $D_{12}=0$ separately. The complete orbit list is
\begin{equation}\label{m6:s:eq:p78-6-10}
\begin{array}{c|c|c}
(b,e)&\text{invariant beyond }(p,k)&\text{transpose}\\\hline
be\ne0&be\in F^\times&\text{fixed orbit}\\
(0,0)&\text{origin}&\text{fixed orbit}\\
(b,0),\ b\ne0&\text{first axis}&\text{second axis}\\
(0,e),\ e\ne0&\text{second axis}&\text{first axis}
\end{array}
\end{equation}
For $be\ne0$, the transporter \eqref{m6:s:eq:p78-6-7} with
$h=e/b$ has character one; the origin is symmetric. The two
punctured axes remain distinct rational orbits.

\subparagraph{The affine phase and the constant-rank distributional slice}\label{m6:s:p78:s6-3}

First-column clearing has potential $\psi(cf/p)$ with difference
$vf+wc+pvw$ (\eqref{m6:s:eq:p78-6-11}--\eqref{m6:s:eq:p78-6-12}).
It is symmetric in $c,f$; their opposite $r$-weights preserve
$dc\,df$, and the raw Jacobian $|p|^2$ is positive and
transpose invariant.

For the remaining quotient, let

\begin{equation}
B_0=\left\{
\bigl(\operatorname{diag}(s,t),\operatorname{diag}(s^{-1},t^{-1})\bigr):s,t\in F^\times
\right\}\subset G_0.
\label{m6:s:eq:p78-6-13}
\end{equation}

The kernel-free residual compatible space on this branch has the associated-bundle description

\begin{equation}
{
\mathscr Y_{7,q\ne0,p\ne0}
\simeq G_0\times^{B_0}
\bigl((F^\times)^2_{p,k}\times F^2_{b,e}\bigr).
}
\label{m6:s:eq:p78-6-14}
\end{equation}

Here $h=s/t$ gives \eqref{m6:s:eq:p78-6-6} and the central diagonal
acts trivially. The normalization proves surjectivity; the complete
gauge comparison, $A_{21}=D_{21}=0$ and $\det A=1$ leave exactly
$B_0$, proving injectivity. Rational column-pivot sections of $G_0/B_0$
and the normalization inverses give analytic charts also at $b=e=0$.

Both $G_0,B_0$ are unimodular. The orbital line restricted to
$B_0$ is trivial by \eqref{m6:s:eq:p78-6-4} and the character-one
transporter \eqref{m6:s:eq:p78-6-7}. Transpose swaps the two base
factors, preserves their positive quotient Haar measure, and
swaps $b,e$. Frobenius reciprocity and affine extraction therefore
give precisely
\eqref{m6:s:eq:p78-6-15}
The bundle dimension is $8-2+4=10$; restoring $42$ affine
directions gives $52$. Off the origin the orbit contributes
$49$ dimensions with one overlapping parameter; at the origin
it contributes $48$ with none. Hence action-plus-slice rank is
$52$ across the stabilizer jump.

\subparagraph{\texorpdfstring{The branch \(q\ne0,\ p=0\)}{The branch q is nonzero, p=0}}\label{m6:s:p78:s6-4}

After \(B=E=f\), invertibility forces \(A_{12}D_{12}\ne0\). The remaining lower-column operations set \(A_{11}=D_{11}=0\). Using both parameters in \eqref{m6:s:eq:p78-6-3} gives

\begin{equation}
{
Y_7(c,k)=x_7\left(f,f,
\begin{pmatrix}0&c\\k&0\end{pmatrix},f\right),
\qquad c,k\in F^\times.
}
\label{m6:s:eq:p78-6-16}
\end{equation}

Before the torus step, write
$A=\left(\begin{smallmatrix}0&b\\a&0\end{smallmatrix}\right)$,
$D=\left(\begin{smallmatrix}0&e\\d&0\end{smallmatrix}\right)$.
Taking $h=b^{-1},r=a^{-1}$ gives the complete invariants $c=be,k=ad$,
with every nonzero pair realized. The first-column potential is
\begin{equation}\label{m6:s:eq:p78-6-17}
\psi\left((d/b)A_{11}+(a/e)D_{11}\right).
\end{equation}
Its coefficients are \eqref{m6:s:eq:p78-4-5}; its raw Jacobian
is $|c|$ and the residual density is one.

The explicit character-one transpose transporter at \eqref{m6:s:eq:p78-6-16} is

\begin{equation}
{
\left(d_5(\operatorname{diag}(1,c^{-1})),
      d_5(\operatorname{diag}(1,c))\right).
}
\label{m6:s:eq:p78-6-18}
\end{equation}

It fixes the complete invariant pair \((c,k)\). The stabilizer is \eqref{m6:s:eq:p78-5-5}, so the positive-measure argument proves intrinsic odd vanishing on the whole branch.

\subparagraph{\texorpdfstring{The branch \(q=0,\ \wp\ne0\)}{The branch q=0, wp is nonzero}}\label{m6:s:p78:s6-5}

On $q=0$, $B,E$ are lower triangular. The nonzero diagonal
pivots give $B=I_2,E=\operatorname{diag}(1,z)$ with $z\ne0$ and
$A_{12}=D_{12}=p\ne0$; the explicit diagonal choices are
\eqref{m6:s:eq:p78-6-19}. Third-factor lower operations, compensated
in the second factors, clear $A_{22},D_{22}$ with zero phase.
First-factor clearing has potential $\psi(A_{11}D_{11}/p)$.
The residual action $(p,a,d)\mapsto(hp,ra,(hr)^{-1}d)$ and
$h=p^{-1},r=a^{-1}$ give

\begin{equation}
{
Z_7(z,k)=x_7\left(f,I_2,
\begin{pmatrix}0&1\\k&0\end{pmatrix},
\operatorname{diag}(1,z)\right),
\qquad z,k\in F^\times,
}
\label{m6:s:eq:p78-6-20}
\end{equation}

where \(k=pad\). The complete invariant formulas before normalization are

\begin{equation}
z=B_{22}E_{22},\qquad
k=\frac{\det A\det D}{(B^{-1}A)_{12}}.
\label{m6:s:eq:p78-6-21}
\end{equation}

Here $B,E$ are lower triangular, and both formulas are invariant
under every lower operation and \eqref{m6:s:eq:p78-4-6}.

The stabilizer is again \eqref{m6:s:eq:p78-5-5}. The character-one transporter to transpose is

\begin{equation}
{
\left(d_5(\operatorname{diag}(1,z)),
      d_5(\operatorname{diag}(1,z^{-1}))\right).
}
\label{m6:s:eq:p78-6-22}
\end{equation}

The intermediate fiber density is $d_7=|h|^2$ on the diagonal
action preserving $B=I_2,E=\operatorname{diag}(1,z)$.
Transpose normalization contributes $|z|^{-2}$ then $|z|^2$;
unipotent restorations have Jacobian one. On the point stabilizer
$h=1$, so both its modular character and the total transpose-density
ratio are one. Lemma~\ref{m8:s:c:transporter} gives odd vanishing.

\label{m6:s:p78:s6-6}
All normalizations use rational pivot charts and only the stated
nonzero divisors. They cover quotient topologies, arbitrary
nonsemisimple coefficients (including $p^2=1,k=p^2$), and all
$h\in F^\times$. The individually labeled rational eigenlines
exclude an additional nonsplit stabilizer torus.

\paragraph{\texorpdfstring{\(\mathcal P_7\): the two-axis conclusion and intrinsic assembly}{P7: the two-axis conclusion and intrinsic assembly}}\label{m6:s:p78:s7}

\subparagraph{A reusable opposite-weight-plane theorem}\label{m6:s:p78:s7-1}

Lemma~\ref{m8:s:c:plane} applies to \eqref{m6:s:eq:p78-6-15}, with
passive $(p,k)$. The two punctured axes are distinct compatible
orbits exchanged by transpose, so their union admits an odd orbital
difference. However $X_7(p,k,b,0)$ remains in the same cell at $b=0$
because both Levi determinants stay nonzero, and likewise on the
other axis. Both closures therefore meet the included origin
$X_7(p,k,0,0)$. The plane density is untwisted, so the lemma's annulus
argument prohibits invariant extension of that odd difference.
\label{m6:s:p78:s7-2}

\subparagraph{\texorpdfstring{Intrinsic vanishing on $\mathcal P_7$}{Intrinsic vanishing on P7}}\label{m6:s:p78:s7-3}

The intrinsic conclusion and its support argument are given in
Proposition~\ref{m6:s:prop:p007}.

\paragraph{\texorpdfstring{\(\mathcal P_8\): one global compatible normalization and the stabilizer}{P8: one global compatible normalization and the stabilizer}}\label{m6:s:p78:s8}

\label{m6:s:p78:s8-1}\label{m6:s:p78:s8-2}
The compatible column equations are
\begin{equation}\label{m6:s:eq:p78-8-1}
(B^{-1}A)_{12}=(E^{-1}D)_{12},\qquad
B_{12}=-E_{22},\quad B_{22}=E_{12}.
\end{equation}
Equivalently, for the nonzero $v=Be_2$,
\begin{equation}\label{m6:s:eq:p78-8-2}
Ee_2=-jv.
\end{equation}
Use the entire invertible $C=B$ in
\eqref{m6:s:eq:p78-4-7}.
Then clear the right lower-column entry and scale its mandatory
nonzero lower-left pivot. These zero-phase operations give
\begin{equation}\label{m6:s:eq:p78-8-3}
B=I_2,\qquad E=f,\qquad A_{12}=D_{22}=p.
\end{equation}
This covers every compatible $F$-point, even when either
component of $Be_2$ vanishes (\S\ref{m6:s:cert:p78-p8-global}).

\subparagraph{The remaining group action and its nontrivial density}\label{m6:s:p78:s8-3}

The transformations preserving \(B=I_2,E=f\), apart from the first-factor lower-column operations, are

\eqref{m6:s:eq:p78-8-4}

Here $h,r\in F^\times$ and $z\in F$.
Its realization and complete converse are recorded in
\S\ref{m6:s:cert:p78-p8-gauge-test}: both restoring second-factor
operations have zero phase and exhaust the orbit equivalence.

For $a=A_{22},b=D_{12}$ the complete remaining action gives
\begin{equation}\label{m6:s:eq:p78-8-5}
p\mapsto p/h,\quad a\mapsto a-zp/h,\quad b\mapsto b+zp/h,
\end{equation}
\begin{equation}\label{m6:s:eq:p78-8-6}
\det A\mapsto(r/h)\det A,\qquad\det D\mapsto r^{-1}\det D.
\end{equation}
The retained fiber and dual density characters are
\begin{equation}\label{m6:s:eq:p78-8-7}
d_8=|h|^{-3},\qquad d_8^{-1}=|h|^3.
\end{equation}
Indeed $t_2=t_4,c=t_3t_4$ in
\eqref{m6:s:eq:p78-3-16}; the other determinants are one.
Point stabilizers have $h=1$ and transpose restoration $h=-1$.

\subparagraph{Reductive projection of the stabilizer}\label{m6:s:p78:s8-4}

For a normalized general Levi pair let \(v_1=Ae_2,v_2=Be_2\), \(z_1=De_2,z_2=Ee_2\). The complete projected stabilizer equations are

\begin{equation}
{
 b_0=a_0^{-t},\quad
 a_0v_1=tv_1,\quad a_0v_2=tv_2,\quad
 b_0z_1=t^{-1}z_1,\quad b_0z_2=t^{-1}z_2,
\quad\det a_0=t^2.
}
\label{m6:s:eq:p78-8-8}
\end{equation}

The shared $GL_2$ forces $b_0=a_0^{-t}$.
For $(B^{-1}A)_{12}\ne0$, $v_1,v_2$ are independent,
so the projection is $a_0=tI_2,b_0=t^{-1}I_2$ and
\begin{equation}\label{m6:s:eq:p78-8-9}
\mathcal A_x\simeq\mathbb G_m\times
(\mathbb G_a\times\operatorname{Heis}_3),\quad\dim=5.
\end{equation}
This is the full unimodular character-trivial stabilizer for
all $(k,m)$, independent of spectra. On the zero branch both
left vectors lie on one line and the right vectors on its
annihilator; exactly one extra unipotent parameter appears.
The lower-root sufficiency test is \S\ref{m6:s:cert:p78-p8-projection};
its group character is computed below.

\paragraph{\texorpdfstring{\(\mathcal P_8\): all rational normal forms, transpose, and intrinsic vanishing}{P8: all rational normal forms, transpose, and intrinsic vanishing}}\label{m6:s:p78:s9}

The resulting statement is Proposition~\ref{m6:s:prop:p008}.
\par\noindent\textit{Proof of Proposition~\ref{m6:s:prop:p008}.}

\subparagraph{\texorpdfstring{The branch \(p\ne0\)}{The branch p is nonzero}}\label{m6:s:p78:s9-1}

For $p\ne0$, let $a=A_{22}$, $b=D_{12}$ and $\alpha=\det A$.
The zero-phase choice $h=p$, $r=p/\alpha$, $z=(a-b)/2$,
followed by first-column clearing, gives $p=1$,
$m=(a+b)/2$, $\det A=1$ and the normal family below.
The precise operations are \eqref{m6:s:eq:p78-9-1} in
\S\ref{m6:s:cert:p78-p8-nonzero}.

\begin{equation}
{
X_8(k,m)=x_8\left(
\begin{pmatrix}0&1\\-1&m\end{pmatrix},I_2,
\begin{pmatrix}k&m\\0&1\end{pmatrix},f\right),
\qquad k\in F^\times,\ m\in F.
}
\label{m6:s:eq:p78-9-2}
\end{equation}

Its two complete invariants, already before \eqref{m6:s:eq:p78-8-3}, are

\begin{equation}
{
k=\frac{\det A\det D}{(B^{-1}A)_{12}},
\qquad
m=\frac{(B^tD)_{12}-\det E\,(B^{-1}A)_{22}}2.
}
\label{m6:s:eq:p78-9-3}
\end{equation}

The invariant derivation is \eqref{m6:s:eq:p78-9-4}.
Between canonical forms, $p=1$ forces $h=1$, equal second-column
entries force $z=0$ and fix $m$, and $\det A=1$ forces $r=1$;
both first-column lower parameters vanish. Thus $(k,m)$ are
complete rational invariants.
The clearing potential is
\begin{equation}\label{m6:s:eq:p78-9-5}
\psi(A_{11}D_{21})\quad(p=1),
\end{equation}
and for general $p$ it is $\psi(A_{11}D_{21}/p)$, with difference
$vD_{21}+wA_{11}+pvw$ as in \eqref{m6:s:eq:p78-4-5}.
The shared-central normalization has phase zero.
For $a=m+z,b=m-z$, $da\,db=|2|\,dm\,dz$; this positive scalar
is retained over dyadic fields, along with \eqref{m6:s:eq:p78-8-7}.

Define the fixed character-one translate of transpose

\begin{equation}
\theta(g)=d_5(f)g^t d_5(f).
\label{m6:s:eq:p78-9-6}
\end{equation}

It has square one and normalizes the original action and character. Direct multiplication gives

\begin{equation}
{\theta(X_8(k,m))=X_8(k,m).}
\label{m6:s:eq:p78-9-7}
\end{equation}

The multiplication test \S\ref{m6:s:cert:p78-p8-transpose}
retains the $\varepsilon$ frame. Its reductive transporter
$(-f,-f)$ equals $(f,f)$ on the group point; its phase is zero
and restoration densities have absolute value one.
The whole stabilizer \eqref{m6:s:eq:p78-8-9} is unchanged for all
$(k,m)$, so its positive orbital measure is $\theta$-even.
Localization proves odd vanishing on the entire $p\ne0$ branch.

\subparagraph{\texorpdfstring{The branch \(p=0\): rational normal forms before imposing the extra condition}{The branch p=0: rational normal forms before imposing the extra condition}}\label{m6:s:p78:s9-2}

For $p=0$, write
$A=\left(\begin{smallmatrix}c&0\\d&a\end{smallmatrix}\right)$,
$D=\left(\begin{smallmatrix}e&b\\f_0&0\end{smallmatrix}\right)$;
invertibility forces $abcf_0\ne0$.
The clearing parameters and the final $h=cf_0,r=f_0,z=0$
normalization are in \S\ref{m6:s:cert:p78-p8-zero}.

The complete normal forms are

\begin{equation}
{
Y_8(a,b)=x_8\left(
\operatorname{diag}(1,a),I_2,
\begin{pmatrix}0&b\\1&0\end{pmatrix},f\right),
\qquad a,b\in F^\times.
}
\label{m6:s:eq:p78-9-8}
\end{equation}

The two ordered invariants are

\begin{equation}
{a=-\det E\,(B^{-1}A)_{22},\qquad b=(B^tD)_{12}}
\label{m6:s:eq:p78-9-9}
\end{equation}

On $p=0$, \eqref{m6:s:eq:p78-8-4} and the first-column operations
fix each invariant. The normalization reaches every $Y_8(a,b)$,
so the ordered pair is a complete invariant.

The phase before the first-column clearing is

\begin{equation}
\psi\left(\frac{d f_0}{a}+\frac{e c}{b}\right).
\label{m6:s:eq:p78-9-10}
\end{equation}

The phase differences are $vf_0,wc$ as in \eqref{m6:s:eq:p78-4-5}.
The raw translation Jacobian is $|ab|$ and the diagonal restoration
density is $|cf_0|^{-3}$. The following stabilizer and transpose
calculations retain these positive normalization factors.

\subparagraph{The enlarged stabilizer, its brackets, and its character}\label{m6:s:p78:s9-3}

The complete boundary group
\eqref{m6:s:eq:p78-9-11}--\eqref{m6:s:eq:p78-9-14}
is $\mathbb G_m$ times a split connected five-dimensional
three-step unipotent group, hence unimodular. The residual
stabilizer has $h=r=1$ and one parameter $z$, lifted by
$\exp(zV)$; its fiber is the four-dimensional graph kernel.
Its phase is $\lambda(K_j)=0$, $\lambda(V)=1/a-1/b$.

Equivalently,

\eqref{m6:s:eq:p78-9-15}

There is also a direct group-normalization check of the coefficient. In \eqref{m6:s:eq:p78-8-4} take \(h=r=1\). Before restoring the two first columns,

\[
A=\begin{pmatrix}1&0\\-z&a\end{pmatrix},\qquad
D=\begin{pmatrix}z&b\\1&0\end{pmatrix}.
\]

The restoring lower parameters are \(z/a\), \(-z/b\). Formula \eqref{m6:s:eq:p78-4-5} gives \(z/a-z/b\), exactly \eqref{m6:s:eq:p78-9-15}. The shared-central and second-factor restorations have zero phase. This computation checks the potentially hidden stabilizer character beyond the unipotent sieve.

Because \(a,b\ne0\), the additional compatibility condition is

\begin{equation}
{a=b.}
\label{m6:s:eq:p78-9-16}
\end{equation}

The stabilizer is unimodular throughout $ab\ne0$. At $a\ne b$ its
additive character is nontrivial, giving the stabilizer-character
exclusion.

\subparagraph{Transpose on all invariants and the comparison of transported densities}\label{m6:s:p78:s9-4}

The same fixed involution \eqref{m6:s:eq:p78-9-6} acts by

\begin{equation}
{\theta(Y_8(a,b))=Y_8(b,a).}
\label{m6:s:eq:p78-9-17}
\end{equation}

Thus $a=b$ is exactly the transpose-stable boundary and every
exchanged pair is already incompatible with the stabilizer character.
The basis comparison \eqref{m6:s:eq:p78-9-18} has raw determinant
$|a/b|$, compensated by the transported kernel vectors $-bK_1$
and $-aK_1$. On $a=b$ both comparisons already have determinant one.

Hence the raw and transported densities agree on $a=b$.
Transpose normalization has zero phase and shared-central matrix
$j$, followed by scalar $-1$ and $h=-1$; all density factors have
absolute value one. On $a=b$ the ordinary transpose transporter is
\begin{equation}\label{m6:s:eq:p78-9-19}
(d_5(f),d_5(f)).
\end{equation}
It has character one and acts by $+1$ on the positive orbital line.

\subparagraph{Stabilizer-character compatibility and rational coverage}\label{m6:s:p78:s9-5}

On the whole \(\mathcal P_8\) unipotent-compatible family, the stabilizer has dimension five when \(\mathsf P_{12}\ne0\), and dimension six when \(\mathsf P_{12}=0\). All are unimodular. The locus compatible with the stabilizer character is

\begin{equation}
{
\begin{cases}
\mathsf P_{12}\ne0,&\text{no additional equation};\\
\mathsf P_{12}=0,&(Q^tSR)_{12}=-\det S\,\mathsf P_{22},
\end{cases}
}
\label{m6:s:eq:p78-9-20}
\end{equation}

in the coefficient notation of \eqref{m6:s:eq:p78-3-13}, always together with the three equations \eqref{m6:s:eq:p78-3-15}. The second line is \eqref{m6:s:eq:p78-9-16} using \eqref{m6:s:eq:p78-9-9}. The extra incompatible locus is open inside the closed \(\mathsf P_{12}=0\) layer.

Every $(k,m)\in F^\times\times F$ and every ordered
$(a,b)\in(F^\times)^2$ occurs, with the complete displayed
invariants. Split, nonsplit and repeated spectra, including
$m=\pm2$, are covered by the same framed-column calculation.
The density calculation retains the $|2|$ Jacobian.

\subparagraph{\texorpdfstring{Intrinsic assembly for \(\mathcal P_8\)}{Intrinsic assembly for P8}}\label{m6:s:p78:s9-6}

First remove the unipotent-incompatible open by Proposition~\ref{cf:stabilizer}.
On the resulting closed compatible antecedent, the open
$\mathsf P_{12}\ne0$ branch vanishes by \eqref{m6:s:eq:p78-9-7},
\eqref{m6:s:eq:p78-8-9} and localization. On the remaining closed
$\mathsf P_{12}=0$ layer, remove the open $a\ne b$ part by
\eqref{m6:s:eq:p78-9-15}, then the closed $a=b$ part by
\eqref{m6:s:eq:p78-9-19} and Lemma~\ref{m8:s:c:transporter}.
The support sequence \eqref{cf:support-sequence} gives
\eqref{m6:s:eq:p78-9-21}
\hfill\(\square\)\par

\paragraph{Transport away from normal forms and comparison of conclusions}\label{m6:s:p78:s10}

\subparagraph{Character-one transporters at arbitrary points}\label{m6:s:p78:s10-1}
The factorization, affine phases and diagonal operations above give a
rational normalizer $n_xx=y$ on each stated branch. Applying
Lemma~\ref{m6:s:lem:normal-transport}, with $\sigma$ interchanging the two
factors of $\mathcal A$, gives
\begin{equation}\label{m6:s:eq:p78-10-1}
k_x=\sigma(n_x)^{-1}k_y n_x,\qquad k_xx=\tau(x),\qquad\eta(k_x)=1
\end{equation}
whenever $k_yy=\tau(y)$ and $\eta(k_y)=1$.
The normalizer may have nontrivial character. On the two punctured
axes of $\mathcal P_7$, the invariants distinguish the exchanged
orbits; only the complete plane lemma applies there.

\subparagraph{The two residual arguments}\label{m6:s:p78:s10-2}
The residual problems, with their preceding density and phase checks, are
\begin{equation}\label{m6:s:eq:p78-10-2}
\begin{array}{c|c}
\mathcal P_7&(F^\times)^2\times F^2,\quad h(b,e)=(hb,h^{-1}e)\\
\mathcal P_8&F^\times\times F,\quad\text{then ordered }(a,b)\in(F^\times)^2, a=b.
\end{array}
\end{equation}

\subsubsection{Rational calculations for $\mathcal P_9$ and $\mathcal P_{14}$}\label{supp:m6-p009-calculations}

The following calculations supply the rational reductions used in \S\ref{m6:s:sec:p009-p014}. Expanded arrays are also available in \path{anc/model6_calculations.pdf}.

\paragraph{The coefficient family and the intersection}

In this section the upright letters
\(\mathsf P,\mathsf Q,\mathsf R,\mathsf S\) denote \(2\)-by-\(2\)
coefficient matrices, distinct from the parabolic \(P\). Put \(\ell_0(\mathsf P,\mathsf Q)=m(\mathsf Q\mathsf P,\mathsf Q,I_2).\)
Thus a representative with coefficients \((\mathsf P,\mathsf Q,\mathsf
R,\mathsf S)\) is
\(x=lW_9r_0^t\), where \(l=\ell_0(\mathsf P,\mathsf Q)\) and
\(r_0=\ell_0(\mathsf R,\mathsf S)\).
In signed-column notation the two representatives are \(W_9=[1,5,2,3,4,-7,-8,-9,-6,10],\qquad W_{14}=[1,3,4,5,2,-9,-6,-7,-8,10].\)
Here, as in the cell table, a signed integer \(a\) denotes the column
\(\operatorname{sgn}(a)e_{|a|}\). These matrices satisfy \(W_{14}=W_9^t\).

The unipotent compatibility equations are
\(q_{11}=r_{22}\), \(q_{12}=r_{12}=0\). Their entire invertible solution is \eqref{m6:s:eq:p009-family},
with \(\mathsf P,\mathsf S\) otherwise arbitrary. It is a smooth
thirteen-dimensional coefficient space. In particular, neither
\(p_{11}\) nor \(s_{11}\) is required to be nonzero.

The rational parabolic decomposition gives every point of
\(\mathcal P_9(F)\) in the form \(h_1lW_9r_0^th_2^t\).
For clarity, the normalization of the central Levi blocks uses the identity
\begin{equation}\label{m6:s:eq:p009-levi-normalization}
 m(A,B,C)=d_5(C)\ell_0(B^{-1}A,C^{-1}B).
\end{equation}
The fourth and fifth blocks are $(\det C/\det B)B$ and
$(\det C/\det A)A$. Both normalizers $d_5(C)$ lie in $H$,
with arbitrary $C\in GL_2(F)$, so all similitudes are retained.

We give the root notation used in the stabilizer and conormal
calculations. The matrix
\begin{equation}\label{m6:s:eq:p009-root-matrix}
\mathcal N(X,Y,Z,W,T_0,S_0)=
\begin{pmatrix}
0&X&Z&W&S_0\\
0&0&Y&T_0&-W^\sharp\\
0&0&0&-Y^\sharp&-Z^\sharp\\
0&0&0&0&-X^\sharp\\
0&0&0&0&0
\end{pmatrix}
\end{equation}
has \(X,Y,Z,W\in M_2(F)\) and \(T_0,S_0\in\mathfrak{sl}_2(F)\).
The symbol \(X_{ij}\), respectively \(Y_{ij},Z_{ij},W_{ij}\),
means that the corresponding coordinate is \(E_{ij}\) and all others
are zero. In a \(T_0\) or \(S_0\) coordinate, the suffixes \(h,e,f\)
mean \(\operatorname{diag}(1,-1),E_{12},E_{21}\).
Define also
\[
 \mathcal L(A,B,C)=\operatorname{diag}
 \bigl(A,B,C,B+(\operatorname{tr}C-\operatorname{tr}B)I_2,
 A+(\operatorname{tr}C-\operatorname{tr}A)I_2\bigr),
\]
and \(L_A=\mathcal L(E_{21},0,0)\),
\(L_B=\mathcal L(0,E_{21},0)\), \(L_C=\mathcal L(0,0,E_{21})\).
All these symbols denote \(10\)-by-\(10\) matrices.

Let \(I_9=P\cap W_9P^tW_9^{-1}\). To determine its Lie algebra,
start with the twelve entries of \(\mathcal L(A,B,C)\), the sixteen
entries of \(X,Y,Z,W\), and the six trace-zero entries of \(T_0,S_0\).
Conjugation by the displayed signed permutation and the second
parabolic zero pattern retain precisely the six diagonal Levi
coordinates and
\begin{equation}\label{m6:s:eq:p009-intersection-roots}
 L_A,L_B,L_C,X_{21},Y_{11},Y_{21},Z_{11},Z_{21},T_f.
\end{equation}
Indeed signed-permutation conjugation sends each $E_{ab}$ to
$\epsilon_a\epsilon_bE_{\pi^{-1}(a),\pi^{-1}(b)}$; requiring each
linked root entry and diagonal block to be lower block triangular
retains exactly these roots and removes nineteen coordinates.
Hence $\dim I_9=15$.

For use at the group level, write the first three left Levi
projections of \(k\in I_9\) as \(a,b,c\), and those of
\(k^\natural=(W_9^{-1}kW_9)^t\) as \(d,e,f\).
All six are lower triangular, and their ordered diagonal entries are
\begin{equation}\label{m6:s:eq:p009-gauge-diagonals}
\begin{aligned}
 \operatorname{diag}(a)&=(a_1,a_2),&
 \operatorname{diag}(b)&=(b_1,b_2),&
 \operatorname{diag}(c)&=(c_1,c_2),\\
 \operatorname{diag}(d)&=(a_1,c_1),&
 \operatorname{diag}(e)&=(a_2,b_1),&
 \operatorname{diag}(f)&=(b_2,\delta/b_2),\quad \delta=c_1c_2.
\end{aligned}
\end{equation}
All six diagonal and six lower parameters are independent.
A rational section is $U_LT_9U_R$: $T_9$ has diagonal
$(a_1,a_2,b_1,b_2,c_1,c_2,\delta/b_2,\delta/b_1,\delta/a_2,\delta/a_1)$;
$U_L$ is the ordered product of $I+sL_A,I+sL_B,I+sL_C$
with $s=a_{21}/a_1,b_{21}/b_1,c_{21}/c_1$, respectively;
$U_R$ uses $Z_{11},X_{21},T_f$ with
$s=d_{21}/a_1,e_{21}/a_2,-f_{21}/b_2$.
These square-zero roots have trivial opposite Levi projections.
Their common kernel is $FY_{11}\oplus FY_{21}\oplus FZ_{21}$.
Thus this section realizes the whole projection image over $F$;
the expanded products are in \texttt{M6-P00914}.

\paragraph{Stabilizers, their characters, and densities}

\begin{proposition}\label{m6:s:prop:p009-stabilizer}
For \eqref{m6:s:eq:p009-family}, put \(p=p_{12}\), \(s=s_{12}\), and
\[
 a_{\mathsf P}=\frac{p_{11}^2}{\det\mathsf P},\qquad
 a_{\mathsf S}=\frac{s_{11}^2}{\det\mathsf S},\qquad
 L_9=L_C+\frac qvL_B+\frac qv a_{\mathsf P}L_A,\qquad
 \mathscr R_9=T_f-a_{\mathsf S}X_{21}-\frac rq a_{\mathsf S}Z_{11}.
\]
Let
\begin{equation}\label{m6:s:eq:p009-nilpotent-stabilizer}
 \mathfrak n_x=
 \{u_1Y_{11}+u_2Y_{21}+u_3Z_{21}+zL_9+t\mathscr R_9:
       pz=0,\ st=0\}.
\end{equation}
The stabilizer of \(x=lW_9r_0^t\) is
\begin{equation}\label{m6:s:eq:p009-full-group}
\mathcal A_x=
\left\{\left(cl\exp(N)l^{-1},
c^{-1}r_0\exp(N')r_0^{-1}\right):
c\in F^\times,\ N\in\mathfrak n_x,\quad
N'=-(W_9^{-1}NW_9)^t\right\}.
\end{equation}
It is a central \(\mathbb G_m\) times a connected unipotent group,
has trivial algebraic component group, and has dimension \(4+\mathbf1_{p=0}+\mathbf1_{s=0}.\)
Its character is \(\psi(\Lambda_9(N))\), where
\begin{equation}\label{m6:s:eq:p009-character}
 \Lambda_9(N)=-p\,a_{\mathsf S}t
              -\frac qv(ra_{\mathsf P}+s)z.
\end{equation}
Both its modular function and its relative orbital density character
are one.
\end{proposition}

\begin{proof}
The six diagonal-product equations on both sides and the lower
$\mathsf Q$ conjugacy force one existing scalar $c$, with multiplier
$c^2$, with $c$ an existing matrix entry.
After removing it, the left $E_{21}$ equation is
$(q/v)zE_{21}\mathsf P=\mathsf P aE_{21}$, which gives
$pz=0$ and $a=(q/v)a_{\mathsf P}z$. The right equation gives
$st=0$ and the coefficients in $\mathscr R_9$; the common kernel
is $FY_{11}\oplus FY_{21}\oplus FZ_{21}$.
The remainder is triangular unipotent. Finite logarithm and
exponential preserve its patterns, conjugacy and symplectic equations,
so the formula describes the whole group over $F$.
The additive phase on the intersection roots is
$(0,0,0,p,0,0,-r,-s,0)$ in the order
$(Y_{11},Y_{21},Z_{21},X_{21},Z_{11},T_f,L_A,L_B,L_C)$;
substitution gives $\Lambda_9$.
The scalar is central and BCH translations in a central-series
basis have determinant one. Thus the stabilizer and $\mathcal A$
are unimodular, giving relative density one.
The expanded elimination and group inverse are
\S\ref{m6:s:cert:p009-stabilizer}.
\end{proof}

Write \(F_0=rp_{11}+sp_{22}\). The complete character table is
\eqref{m6:s:eq:p009-character-table}
The nonzero coefficients in the second and fourth rows cannot
vanish, by invertibility of the triangular matrices. A nonzero
coefficient times \(F\) is all of \(F\), regardless of the conductor
of \(\psi\). It follows that the locus where the character on the stabilizer is trivial is
\eqref{m6:s:eq:p009-full-compatible}
In the second set \(s\ne0\) follows from invertibility.
This union is not locally closed. It contains the dense open
\(\{ps\ne0\}\) of the irreducible coefficient space, but is not open
at a point of \(p=F_0=0\): varying \(p_{11}\) slightly while retaining
\(p=0\) changes \(F_0\) to a nonzero value without changing
invertibility. A dense locally closed subset would be open.
Pullback through the coefficient map proves the same assertion in
the cell.

\paragraph{Rational orbit classification}

Changing the parabolic factorization by \(k\in I_9\) replaces
\((l,r_0)\) by \((lk,r_0(k^\natural)^{-1})\). Normalizing the central
Levi blocks gives the complete coefficient action
\begin{equation}\label{m6:s:eq:p009-gauge-action}
 (\mathsf P,\mathsf Q,\mathsf R,\mathsf S)\longmapsto
 (b^{-1}\mathsf Pa,c^{-1}\mathsf Qb,e\mathsf Rd^{-1},f\mathsf Se^{-1}).
\end{equation}
Conversely, an equality of \(\mathcal A(F)\)-orbits, after the two
parabolic factors are moved to opposite sides of \(W_9\), produces
an element of \(I_9(F)\). Thus this action gives necessity as well
as sufficiency for orbit equivalence.

Its relevant weights on the entire coefficient space are
\[
\begin{gathered}
 p'=(a_2/b_1)p,\quad s'=(b_2/b_1)s,\quad
 q'=(b_1/c_1)q,\quad v'=(b_2/c_2)v,\quad r'=(a_2/a_1)r,\\
 \det\mathsf P'=\frac{a_1a_2}{b_1b_2}\det\mathsf P,\qquad
 \det\mathsf S'=\frac{c_1c_2}{a_2b_1}\det\mathsf S.
\end{gathered}
\]
On \(p=0\), moreover,
\(p_{11}'=(a_1/b_1)p_{11}\),
\(p_{22}'=(a_2/b_2)p_{22}\), and
\(F_0'=(a_2/b_1)F_0\).
These identities prove the invariance of \(p=0\), \(s=0\), and
\(p=F_0=0\), including the assertions about descended closed loci
used below.

The gauge choice
\[
 a=\operatorname{diag}(r,1),\quad b=e=I,\quad
 c=\mathsf Q,\quad d=\mathsf R,\quad f=\operatorname{diag}(1,qv)
\]
is allowed by \eqref{m6:s:eq:p009-gauge-diagonals} and has the rational
lift already displayed. It gives
\begin{equation}\label{m6:s:eq:p009-global-normalization}
 \mathsf Q'=\mathsf R'=I,\qquad
 \widehat{\mathsf P}=\mathsf P\operatorname{diag}(r,1),\qquad
 \widehat{\mathsf S}=\operatorname{diag}(1,qv)\mathsf S .
\end{equation}
Thus \(u,w\) introduce no orbit invariant. The group preserving
\(\mathsf Q=\mathsf R=I\) consists of
\[
\begin{aligned}
 a&=\begin{pmatrix}a_0&0\\\alpha&a_0\end{pmatrix},&
 b=c&=\begin{pmatrix}b_1&0\\\beta&b_2\end{pmatrix},\\
 d=e&=\begin{pmatrix}a_0&0\\\varepsilon&b_1\end{pmatrix},&
 f&=\begin{pmatrix}b_2&0\\\varphi&b_1\end{pmatrix}.
\end{aligned}
\]
All four lower entries are free and the three diagonal parameters
are nonzero. The effective parameters \(A=a_0/b_1\) and \(B=b_2/b_1\)
range over all of \((F^\times)^2\).

The resulting statement is Proposition~\ref{m6:s:prop:p009-rational-orbits}.

\begin{proof}[Proof of Proposition~\ref{m6:s:prop:p009-rational-orbits}]
First apply \eqref{m6:s:eq:p009-global-normalization}.
On $ps\ne0$, the four lower parameters clear the two diagonal
entries of each remaining matrix, giving
$\left(\begin{smallmatrix}0&p\\-\det\widehat{\mathsf P}/p&0\end{smallmatrix}\right)$
and
$\left(\begin{smallmatrix}0&s\\-\det\widehat{\mathsf S}/s&0\end{smallmatrix}\right)$.
Their entries have diagonal weights $A,A/B,B,A^{-1}$;
$A=1/p,B=1/s$ gives $x_{\rho,\sigma}$.
On $p=F_0=0$, the normalized equation is
$\widehat p_{11}+\widehat s_{12}\widehat p_{22}=0$, with
all three factors nonzero. Lower operations clear
$\widehat p_{21}$ and the diagonal entries of $\widehat{\mathsf S}$;
$B=1/\widehat s_{12}$ and
$A=1/(\widehat s_{12}\widehat p_{22})$ give $y_\kappa$.
The preceding weight equations prove necessity of each invariant;
the reductions prove sufficiency and realization of every value.
The explicit section of $I_9(F)$ lifts all operations to $H(F)^2$.
These rational reductions cover nonsplit, repeated-root and
nonsemisimple coefficients.
At $x_{\rho,\sigma}$ the three kernel roots commute.
At $y_\kappa$ add $L_*=-L_A+L_B+L_C$, with sole bracket
$[L_*,Y_{11}]=Y_{21}$, proving the stated Heisenberg factor.
The pivot parameters and BCH law are recorded in
\S\ref{m6:s:cert:p009-normalization}.
\end{proof}

\paragraph{Closed orbits and intrinsic odd distributions}

Let $\mathcal C_9$ be the unipotent-compatible envelope.
The smooth parabolic quotient $P^2\to\mathcal P_9$, with
fifteen-dimensional fibers, has saturated inverse image
$H^2\times\{\eqref{m6:s:eq:p009-family}\}$.
Compatibility is gauge invariant. Local analytic sections make
this an open quotient on $F$-points; closedness and smoothness
therefore descend, with
\begin{equation}\label{m6:s:eq:p009-support-dimension}
\dim\mathcal C_9=68-3-15=50.
\end{equation}

For fixed \(\rho,\sigma\ne0\), the inverse image in the coefficient space
of the orbit of \(x_{\rho,\sigma}\) is
\begin{equation}\label{m6:s:eq:p009-closed-generic}
 qvp\det\mathsf S+\rho s=0,\qquad
 rqv\det\mathsf P\det\mathsf S-\rho\sigma p=0.
\end{equation}
The second equation forces \(p\ne0\), and the first forces \(s\ne0\).
Thus these saturated equations add no double-zero boundary.
They are equivalent to \eqref{m6:s:eq:p009-generic-invariants}.
The orbit of \(y_\kappa\) similarly has inverse image
\begin{equation}\label{m6:s:eq:p009-closed-exceptional}
 p=0,\qquad rp_{11}+sp_{22}=0,\qquad
 \kappa+qvp_{22}\det\mathsf S=0.
\end{equation}
Invertibility forces \(s\ne0\) here as well.
Quotient topology now proves that each individual rational orbit
is closed in \(\mathcal C_9\) and hence in the whole cell.
The smooth orbit map \(\mathcal A(F)\to\mathcal A(F)x\) has local
analytic sections, so the subspace topology of the orbit is its
homogeneous quotient topology.

For either family the character and relative density on the
stabilizer are trivial. With a positive invariant Radon measure,
\begin{equation}\label{m6:s:eq:p009-orbital-distribution}
\mathscr T_x(f)=\int_{\mathcal A(F)/\mathcal A_x(F)}
f(ax)\eta(a)\,d\mu_x(a\mathcal A_x)
\end{equation}
is finite on compact test supports: their intersection with the
closed orbit is compact in its quotient topology, and $|\eta|=1$.
Changing variables gives
$\mathscr T_x(f\circ a_0)=\eta(a_0)^{-1}\mathscr T_x(f)$.

Transpose swaps coefficients to
$(\mathsf R,\mathsf S,\mathsf P,\mathsf Q)$ and swaps the full
stabilizer factors. Thus it transports the complete compatible
family, its character, invariant parameters and relative density
to $\mathcal P_{14}$; the expanded partner formulas are in
\texttt{M6-P00914}. The rank matrices are
\begin{equation}\label{m6:s:eq:p009-ranks}
R_9=\begin{pmatrix}7&6&4&2\\5&5&4&2\\4&4&3&1\\2&2&2&1\end{pmatrix},
\qquad R_{14}=R_9^t.
\end{equation}
They are incomparable, so each cell is open and closed in their
pair and no $\mathcal A$-transporter reaches its partner.
Transport the positive measure and phase by transpose.
The transported logarithmic Jacobian is one; scalar inversion
preserves multiplicative Haar measure, and $\tau^2=1$.

The intrinsic conclusion and its support argument are given in
Proposition~\ref{m6:s:prop:p009-intrinsic}.

\subsubsection{Rational calculations for $\mathcal P_{12}$ and $\mathcal P_{22}$}\label{supp:m6-p012-calculations}

The following calculations supply the rational reductions used in \S\ref{m6:s:sec:p012-p022}. Expanded arrays are also available in \path{anc/model6_calculations.pdf}.

\paragraph{The entire coefficient family and rational gauge}

Retain the upright coefficient notation and the section
\(\ell_0(\mathsf P,\mathsf Q)=m(\mathsf Q\mathsf P,\mathsf Q,I_2)\)
of Section~\ref{m6:s:sec:p009-p014}. The signed representative is
\begin{equation}\label{m6:s:eq:p012-representative}
 W_{12}=[e_1,e_8,e_2,e_5,e_4,-e_7,-e_6,-e_9,e_3,e_{10}].
\end{equation}
Every rational point of its cell is \(h_1\ell_0(\mathsf P,\mathsf Q)W_{12} \ell_0(\mathsf R,\mathsf S)^th_2^t,\qquad h_1,h_2\in H(F).\)
Identity \eqref{m6:s:eq:p009-levi-normalization} proves this from the
rational parabolic factorization and accounts for every multiplier.
The four unipotent compatibility equations are \(q_{11}=0,\qquad q_{12}=s_{12},\qquad q_{21}=r_{22},\qquad r_{12}=0.\)
Their entire invertible solution is \eqref{m6:s:eq:p012-family}.
Write \(\Delta=gi-bh\), and denote this coefficient space by
\(\mathcal Q_{12}\). It has twelve independent coordinates:
\(\det\mathsf Q=-bc\), \(\det\mathsf R=ec\), and the remaining
determinants are \(pt-qr,\Delta\).
Thus \(\mathcal Q_{12}\) is smooth, closed relative to
\(\operatorname{GL}_2^4\), and open in its affine twelve-dimensional
linear space. The optional entries \(d,f,g,h,i\) may vanish; so
may \(t\).

Put \(I_{12}=P\cap W_{12}P^tW_{12}^{-1}\).
Its six diagonal parameters give the torus
\begin{equation}\label{m6:s:eq:p012-intersection-torus}
 \operatorname{diag}
 (a,b_0,c_0,d_0,e_0,\mu/e_0,\mu/d_0,\mu/c_0,\mu/b_0,\mu/a).
\end{equation}
The twelve off-diagonal directions of its Lie algebra, in physical
matrix units, are
\begin{equation}\label{m6:s:eq:p012-intersection-roots}
\begin{aligned}
 &E_{21}+E_{10,9},\ E_{34}+E_{78},\ E_{65},\
 E_{18}+E_{3,10},\ E_{25}+E_{69},\ -E_{47},\\
 &E_{38},\ -E_{28}+E_{39},\ -E_{37}+E_{48},\
 E_{36}+E_{58},\ E_{45}+E_{67},\ -E_{35}+E_{68}.
\end{aligned}
\end{equation}
These lists are obtained by imposing block upper triangularity,
block lower triangularity after conjugation by \(W_{12}^{-1}\),
and \(Z^tJ+JZ=\lambda J\). The eighteen free parameters consist of the six diagonal
coordinates and these twelve roots. Equivalently,
\begin{equation}\label{m6:s:eq:p012-eighteen-matrix}
\mathfrak i_{12}=\mathfrak t_6\oplus
\operatorname{span}_F\{\text{the twelve roots in
\eqref{m6:s:eq:p012-intersection-roots}}\}.
\end{equation}
The expanded $10$-by-$10$ linear array and all zero-pattern
conditions are in \texttt{M6-P01222}, equation
\texttt{m6:s:eq:p012-eighteen-matrix}.

In the coordinate
order 
\[
(2,1,3,4,6,5,7,8,10,9),
\]

the common zero pattern is
triangular. The diagonal of an invertible member is exactly
\eqref{m6:s:eq:p012-intersection-torus}; dividing it out leaves a
unipotent matrix. Finite logarithm and exponential preserve the
two patterns and the symplectic equation, as in the proof of
Proposition~\ref{m6:s:prop:p009-stabilizer}. Thus \(I_{12}\) is the
split six-dimensional torus times a connected twelve-dimensional
unipotent radical.

Write the first three Levi blocks of \(k\in I_{12}\) as \(A,B,C\),
and those of \(W_{12}^{-1}kW_{12}\in P^t\) as \(D,E,F_0\).
Their complete projection image is
\begin{equation}\label{m6:s:eq:p012-gauge-group}
\begin{aligned}
A&=\begin{pmatrix}a&0\\u&b_0\end{pmatrix},&
B&=\begin{pmatrix}c_0&v\\0&d_0\end{pmatrix},&
C&=\begin{pmatrix}e_0&0\\z&\mu/e_0\end{pmatrix},\\
D&=\begin{pmatrix}a&k_0\\0&\mu/c_0\end{pmatrix},&
E&=\begin{pmatrix}b_0&l_0\\0&e_0\end{pmatrix},&
F_0&=\begin{pmatrix}d_0&n\\0&\mu/d_0\end{pmatrix}.
\end{aligned}
\end{equation}
All six diagonal entries are nonzero and all six off-diagonal
parameters are free. The first six roots in
\eqref{m6:s:eq:p012-intersection-roots} project to the six respective
entries; the last six form the kernel. Multiplying the torus by
$I+sN$ for each square-zero first root realizes any desired entry
after division by a nonzero diagonal. Their commutators affect
only the kernel. Thus every projected tuple has an explicit
rational lift over $F$.

The complete coefficient equivalence is
\begin{equation}\label{m6:s:eq:p012-gauge-action}
 (\mathsf P,\mathsf Q,\mathsf R,\mathsf S)\longmapsto
 (B^{-1}\mathsf PA,\ C^{-1}\mathsf QB,\
 E^t\mathsf RD^{-t},\ F_0^t\mathsf SE^{-t}).
\end{equation}
It follows by replacing the two parabolic factors using \(k\)
and \((W_{12}^{-1}k^{-1}W_{12})^t\), then normalizing their
central blocks. Conversely an equality of two \(H(F)^2\)
orbits gives exactly such an element of \(I_{12}(F)\).
The equivalence relation is therefore the entire rational one.

\paragraph{Explicit rational normal forms}

For \eqref{m6:s:eq:p012-family}, choose the following tuple in
\eqref{m6:s:eq:p012-gauge-group}:
\begin{equation}\label{m6:s:eq:p012-normalizing-gauge}
\begin{aligned}
A&=\operatorname{diag}(-c\Delta e,-c\Delta),&
B&=I,& C&=\begin{pmatrix}b&0\\d&-c\end{pmatrix},\\
D&=\begin{pmatrix}-c\Delta e&-ge-bf\\0&-bc\end{pmatrix},&
E&=\begin{pmatrix}-c\Delta&g\\0&b\end{pmatrix},&
F_0&=\begin{pmatrix}1&ci\\0&-bc\end{pmatrix}.
\end{aligned}
\end{equation}
Here \(\mu=-bc\), and every required inverse exists. Substitution
in \eqref{m6:s:eq:p012-gauge-action} gives
\begin{equation}\label{m6:s:eq:p012-normalized-triple}
 \mathsf Q'=\mathsf S'=J_0=\begin{pmatrix}0&1\\-1&0\end{pmatrix},
 \qquad
 \mathsf R'=A_0=\operatorname{diag}(1,-1),\qquad
 \mathsf P'=\mathsf P\operatorname{diag}(-c\Delta e,-c\Delta).
\end{equation}
These formulas include \(d=0\), \(f=0\), \(g=0\), or \(i=0\)
whenever permitted by the determinants.

For the normalized triple, its $\mathsf Q$ and $\mathsf R$
equations give $C=J_0BJ_0^{-1}$ and $E=A_0DA_0$;
the $\mathsf S$ equation makes the off-diagonal entries of
$D,E,F_0$ zero. The complete remaining tuple is
\[
A=\begin{pmatrix}\alpha&0\\u&\alpha\end{pmatrix},\quad
B=\begin{pmatrix}\alpha&v\\0&\delta\end{pmatrix},\quad
C=\begin{pmatrix}\delta&0\\-v&\alpha\end{pmatrix},\quad
D=E=\operatorname{diag}(\alpha,\delta),\quad
F_0=\operatorname{diag}(\delta,\alpha).
\]

Put \(\rho=\delta/\alpha\), \(h_0=v/\alpha\),
\(t_0=u/\alpha\). Its effective action is
\eqref{m6:s:eq:p012-residual-action}
the additional \(\alpha\in F^\times\) is the scalar action kernel.
For \(\mathsf P=\left(\begin{smallmatrix}p&q\\r&t\end{smallmatrix}\right)\),
the four transformed entries are
\begin{equation}\label{m6:s:eq:p012-residual-entries}
 t'=t/\rho,\quad r'=(r+tt_0)/\rho,\quad
 q'=q-h_0t/\rho,\quad
 p'=p+qt_0-h_0(r+tt_0)/\rho.
\end{equation}

The resulting statement is Proposition~\ref{m6:s:prop:p012-orbits}.

\begin{proof}[Proof of Proposition~\ref{m6:s:prop:p012-orbits}]
First use \eqref{m6:s:eq:p012-normalizing-gauge}. If the resulting
lower-right entry \(t\ne0\), take
\(\rho=t,\ t_0=-r/t,\ h_0=q\) in
\eqref{m6:s:eq:p012-residual-entries}. The result is
\(\operatorname{diag}((pt-qr)/t,1)\).
Under the residual action, determinant and \(t\) are both
multiplied by \(\rho^{-1}\), so their ratio is invariant.
Returning through \eqref{m6:s:eq:p012-normalized-triple} gives the
stated formula for \(\kappa\). If \(t=0\), invertibility gives
\(q,r\ne0\); choose \(\rho=r\) and
\(h_0=p+qt_0\). The normal form is
\(\left(\begin{smallmatrix}0&q\\1&0\end{smallmatrix}\right)\),
and its upper-right entry is invariant under all residual
transformations on this branch. Returning to the original
coefficients gives \(\upsilon\).

Any equivalence between the normalized triples must have the
form \eqref{m6:s:eq:p012-residual-action}, so the invariants are
necessary as well as sufficient. The displayed gauge has a
rational lift to the intersection, hence to rational
\(H^2\)-transporters. This proves \(F\)-rational coverage.
The left-right operations apply to nonsplit, repeated-root, Jordan
and nilpotent coefficient configurations allowed by invertibility.
\end{proof}

Write the compatible normal representative as
\begin{equation}\label{m6:s:eq:p012-physical-normal}
\begin{aligned}
x_\kappa
&=\ell_0(\operatorname{diag}(\kappa,1),J_0)
  W_{12}\ell_0(A_0,J_0)^t\\
&=[-e_7,\kappa e_2,e_5,-e_1,e_3,e_8,e_{10},e_6,
                     \kappa^{-1}e_9,-e_4].
\end{aligned}
\end{equation}
The word compatible here refers to the character on the stabilizer; the next
calculation proves that precisely this branch survives it.

\paragraph{The group stabilizer and its character}

Set \(l=\ell_0(\mathsf P,\mathsf Q)\),
\(r_0=\ell_0(\mathsf R,\mathsf S)\), and \(x=lW_{12}r_0^t\).
The physical-entry elimination
\eqref{m6:s:eq:p012-stabilizer-elimination}--\eqref{m6:s:eq:p012-stabilizer-last-equation}
has six free roots for $t\ne0$ and one extra for $t=0$, where
$qr\ne0$.

Using the root matrices in \eqref{m6:s:eq:p009-root-matrix}, put
\[
 \mathfrak n_6=
 \left\{\mathcal N(0,Y,0,W,T_0,0):
 Y=\begin{pmatrix}y_0&y_1\\y_2&0\end{pmatrix},\
 W=\omega_0E_{22},\
 T_0=t_0\operatorname{diag}(1,-1)+t_1E_{12}\right\}.
\]
Its ordered basis is
\((Y_{11},Y_{12},Y_{21},W_{22},T_h,T_e)\).
The only nonzero brackets, up to antisymmetry, are
\begin{equation}\label{m6:s:eq:p012-kernel-brackets}
 [Y_{11},Y_{12}]=2T_e,\qquad [Y_{12},Y_{21}]=T_h.
\end{equation}
Thus it is two-step nilpotent with center
\(FW_{22}\oplus FT_h\oplus FT_e\).
On \(t=0\) define
\begin{equation}\label{m6:s:eq:p012-extra-generator}
 N=\operatorname{diag}
 \left(\frac{br}{cq}E_{21},\frac bcE_{12},E_{21},
       \frac bcE_{12},\frac{br}{cq}E_{21}\right),
 \qquad \mathfrak n_7=\mathfrak n_6+FN.
\end{equation}
It satisfies \(N^2=0\), \(lNl^{-1}=d_5(E_{21})\), and
\(-(W_{12}^{-1}NW_{12})^t\in\mathfrak u\).
Its additional brackets are
\begin{equation}\label{m6:s:eq:p012-extra-brackets}
 [N,Y_{12}]=-Y_{11},\qquad
 [N,Y_{21}]=(b/c)Y_{11},\qquad
 [N,T_h]=-2(b/c)T_e.
\end{equation}
All its matrices are nilpotent and are strictly triangular in
a common compatible coordinate order.

\begin{proposition}\label{m6:s:prop:p012-stabilizer}
The stabilizer at \(x=lW_{12}r_0^t\) is
\begin{equation}\label{m6:s:eq:p012-full-stabilizer}
 \left\{\left(zl\exp(Z)l^{-1},
 z^{-1}r_0\exp(-(W_{12}^{-1}ZW_{12})^t)r_0^{-1}\right):
 z\in F^\times,\ Z\in\mathfrak n_j\right\},
\end{equation}
where \(j=6\) for \(t\ne0\) and \(j=7\) for \(t=0\).
It is \(\mathbb G_m\times N_j\), with \(N_j\) a split connected
unipotent group. Its dimensions are seven and eight, respectively;
its scalar torus is central and its component group is trivial.
The restriction of $\eta$ to the stabilizer is
\[
 \eta|_{\mathcal A_x}=1\quad(t\ne0),\qquad
 \eta(z,\exp(Z_0+vN))
       =\psi\left(-\frac{ebr}{cq}v\right)\quad(t=0).
\]
Both modular functions and the relative orbital density character
are one. In particular the locus where this restriction is trivial is exactly
\(t\ne0\).
\end{proposition}

\begin{proof}
The complete residual action fixing $\operatorname{diag}(\kappa,1)$
requires $\rho=1,t_0=h_0=0$, leaving only the scalar $\alpha$.
Fixing $\left(\begin{smallmatrix}0&\upsilon\\1&0\end{smallmatrix}\right)$
requires $\rho=1,h_0=\upsilon t_0$, leaving one additive parameter.
The projection kernel is precisely $N_6$; finite logarithm and
exponential preserve its patterns and conjugacy equations.
This proves the full group formula over $F$, including all
components and its polynomial logarithmic inverse.
On the intersection matrix, the left phase is $bz_2+cz_1$;
the right contributions are $-ez_{10}-cz_1$ and $-bz_2$.
The total is $-ez_{10}$, and the extra direction on $t=0$ has
$z_{10}=brv/(cq)$. Its coefficient is nonzero, so its image is
all of $F$ independently of the conductor of $\psi$.
The torus is central and BCH translations in a central-series
basis have Jacobian one, retaining the factors $2$ and $1/2$.
The same conjugation determinant on each root matrix coordinate
of $H$ is one. Thus both modular functions and their relative
ratio are one. The expanded group law and Jacobians are
\S\ref{m6:s:cert:p012-full-group}.
\end{proof}

Thus the compatible orbits have dimension \(52-7=45\);
the exceptional, character-incompatible orbits have dimension
\(52-8=44\). The rational orbit parameter on the compatible
branch is all of \(F^\times\), rather than a finite set.

\paragraph{transpose coordinates and intrinsic nonvanishing}

The designated representative of the partner cell differs from
\(W_{12}^t\) by a Levi sign. With \(A_0=\operatorname{diag}(1,-1)\),
put
\[
 d_{12}=m(A_0,I,I)
 =\operatorname{diag}(1,-1,1,1,1,1,1,1,-1,1).
\]
Comparison of the signed physical columns gives
\begin{equation}\label{m6:s:eq:p012-transpose-sign}
 W_{22}=d_{12}W_{12}^t,\qquad d_{12}^2=I.
\end{equation}
Accordingly the coefficient transpose, when both designated
representatives are used, is
\begin{equation}\label{m6:s:eq:p012-coefficient-transpose}
 (\mathsf P,\mathsf Q,\mathsf R,\mathsf S)_{12}
 \longmapsto(\mathsf RA_0,\mathsf S,\mathsf P,\mathsf Q)_{22}.
\end{equation}
Indeed \(r_0W_{12}^tl^t=(r_0d_{12})W_{22}l^t\), and
\(r_0d_{12}=m(\mathsf S\mathsf RA_0,\mathsf S,I)\).
The inverse map is
\((\mathsf P,\mathsf Q,\mathsf R,\mathsf S)_{22}
\mapsto(\mathsf R,\mathsf S,\mathsf PA_0,\mathsf Q)_{12}\).

Substitution gives the complete partner family and the equations
$p_{12}=0$, $p_{22}+s_{21}=0$, $q_{12}=s_{12}$, $s_{11}=0$,
with the determinant conditions \eqref{m6:s:eq:p012-family}.
Its compatible branch is $r_{22}\ne0$ and has invariant
\begin{equation}\label{m6:s:eq:p012-partner-invariants}
\kappa_{22}=-s_{21}\det(\mathsf Q)p_{11}\det(\mathsf R)/r_{22}.
\end{equation}
The incompatible invariant is
$\upsilon_{22}=-s_{21}\det(\mathsf Q)r_{12}$.
Under \eqref{m6:s:eq:p012-coefficient-transpose}, these equal
$\kappa,\upsilon$ without inversion or sign change.
The normal triples become
$(I,J_0,\operatorname{diag}(\kappa,1),J_0)$ and
$(I,J_0,\left(\begin{smallmatrix}0&\upsilon\\1&0\end{smallmatrix}\right),J_0)$;
the groups, characters and densities are factor-swapped
\eqref{m6:s:eq:p012-full-stabilizer}.
The expanded substitution is in \texttt{M6-P01222}.

Let $\mathcal C_{12}$ be the saturation of $\mathcal Q_{12}$.
The smooth parabolic quotient $q_{12}:P^2\to\mathcal P_{12}$
has eighteen-dimensional fibers and local analytic sections,
so it is an open quotient on rational points. Gauge-invariant
compatibility gives the saturated inverse image
\begin{equation}\label{m6:s:eq:p012-saturated-envelope}
q_{12}^{-1}(\mathcal C_{12})=H(F)^2\times\mathcal Q_{12}(F).
\end{equation}
Its closedness therefore descends to the cell.

For fixed \(\kappa\ne0\), the inverse image in the coefficient space of
the orbit of \(x_\kappa\) is
\begin{equation}\label{m6:s:eq:p012-closed-orbit}
 \kappa t+c\Delta e\det\mathsf P=0
\end{equation}
inside \(\mathcal Q_{12}\). Its second term never vanishes,
so it forces \(t\ne0\).
The rational normal forms identify every solution with this
one orbit. Hence its inverse image is closed and saturated,
so the orbit is closed in the whole cell, algebraically and
on \(F\)-points. The rational reduction and rational lifting
exclude a smaller analytic rational orbit inside this fiber.
The smooth orbit map \(H(F)^2\to H(F)^2x_\kappa\) has local
analytic sections as well; its quotient topology is the
subspace topology. In the partner coefficients the closed
equation is \(\kappa r_{22}+s_{21}\det(\mathsf Q)p_{11}\det(\mathsf R)=0.\)

With positive invariant quotient measure, set
\begin{equation}\label{m6:s:eq:p012-orbital-distribution}
\mathscr T_\kappa(f)=\int_{\mathcal A(F)/\mathcal A_{x_\kappa}(F)}
f(ax_\kappa)\eta(a)\,d\mu_\kappa(a\mathcal A_{x_\kappa}).
\end{equation}
Character triviality makes the phase well defined, and density one
supplies the measure. Closedness and equality of quotient and
subspace topologies make the integral finite on compact supports.
It satisfies $\mathscr T_\kappa(f\circ a_0)=\eta(a_0)^{-1}\mathscr T_\kappa(f)$.

The rank matrices are
\begin{equation}\label{m6:s:eq:p012-rank-matrices}
 R_{12}=\begin{pmatrix}7&6&4&2\\5&4&3&1\\3&3&2&1\\2&2&2&1\end{pmatrix},
 \qquad R_{22}=R_{12}^t .
\end{equation}
Their \((1,2)\) and \((2,1)\) entries are respectively \(6,5\)
and \(5,6\), so each cell is open and closed in their pair.

The intrinsic conclusion and its support argument are given in
Proposition~\ref{m6:s:prop:p012-intrinsic}.

\section{Model 8: rational families and stabilizers}\label{supp:m8-families}
\subsection{Supplementary stabilizer catalogues}\label{m8:s:app:B}
\subsubsection{Ordinary stabilizer lifts on the invertible-final-block stratum}
At the Levi point in \eqref{m8:s:s2:e5-5}, the ordinary \(K\)-stabilizer is exactly
\begin{equation}
{
K_m=\{(v(s),v(-\kappa s),h):h\in Z_{GL_2}(a,C)\}.
}
\label{m8:s:s2:e5-11}
\end{equation}
Uniqueness of the factorization makes this the full stabilizer. At a general point \(x=n_Lmn_R^t\), its stabilizer is the conjugate of \eqref{m8:s:s2:e5-11} by \((n_L,n_R,1)\in K\).

\subsubsection{General ordinary stabilizer on the zero-block Siegel boundary}
For completeness, the full stabilizer before imposing \eqref{m8:s:s3:eq:8-5} consists of
\begin{equation}
{
h\in Z_{\mathrm{GL}_2}(B)\cap Z_{\mathrm{GL}_2}(C),\quad P\in M_2(F),\quad
R=B^{-1}(A-hAh^{-1}-PC).
}
\label{m8:s:s3:eq:8-6}
\end{equation}
This follows directly from \eqref{m8:s:s3:eq:3-3} for every \(A\).

\subsubsection{S4 intersection group and parameter equivalence}
\paragraph{The full equivalence relation}\label{m8:s:s45:sub3.3}

Let \(P=LU\) be the fixed nonmaximal parabolic. The relevant intersection is
\[
\mathcal J_w=P\cap wP^tw^{-1}.
\]
Every element of \(\mathcal J_w\), and only such an element, has the form
\begin{equation}\label{m8:s:s45:eq3.6}
{
j=L\left(
\begin{pmatrix}\alpha&0\\\gamma&\delta\end{pmatrix},
\begin{pmatrix}\epsilon&\phi\\0&\zeta\end{pmatrix},\omega
\right)
 u\left(0,y,\begin{pmatrix}0&p\\q&r\end{pmatrix},0\right),
}
\end{equation}
with \(\alpha,\delta,\epsilon,\zeta,\omega\ne0\) and the six other parameters arbitrary. This follows by imposing the lower block-triangular zero conditions on \(w^{-1}jw\); the form equations then determine its other entries. It is an eleven-dimensional split solvable group.

Let \(\pi_L:P\to L\), \(\pi_L^-:P^t\to L\) be the block-diagonal projections. For \eqref{m8:s:s45:eq3.6},
\[
\begin{aligned}
\pi_L(j)&=L(A_0,B_0,\omega),\\
\pi_L^-(w^{-1}jw)&=L(C_0,E_0,\omega),
\end{aligned}
\]
where
\begin{equation}\label{m8:s:s45:eq3.7}
{
C_0=\begin{pmatrix}\alpha&\alpha p\\0&\omega/\epsilon\end{pmatrix},
\qquad
E_0=\begin{pmatrix}\delta&\delta q\\0&\omega/\zeta\end{pmatrix}.
}
\end{equation}

Write a normal point as \(x(l,r_0)=lwr_0\), with \(l,r_0\in L\). Two normal points are in the same \textbf{even} \(K(F)\)-orbit if and only if there exist \(h\in \operatorname{GL}_2(F)\), \(j\in\mathcal J_w(F)\) such that
\begin{equation}\label{m8:s:s45:eq3.8}
{
\begin{aligned}
l'&=D(h)l\,\pi_L(j)^{-1},\\
r_0'&=\pi_L^-(w^{-1}jw)\,r_0D(h)^{-1}.
\end{aligned}}
\end{equation}
This is the residual equivalence relation on all the parameters of \eqref{m8:s:s45:eq3.2}.

Necessity and sufficiency, including the unipotents, are Lemma~\ref{m8:s:c:normal-family} with this $w$ and these projections.

Formula \eqref{m8:s:s45:eq3.8} gives the exact orbit redundancy of the parameters in \eqref{m8:s:s45:eq3.2}. The vanishing proof in \S\ref{m8:s:sec:s4-s5} uses the two-unipotent stabilizer on the original stratum.

\subsubsection{Full S4 diagonal stabilizer catalogue}
\label{m8:s:s45:sub5}
\begin{proposition}\label{m8:s:s45:full-stabilizer}
The full stabilizer of \eqref{m8:s:s45:eq3.2} consists exactly of the triples \eqref{m8:s:s45:eq5.2}, for \(h\) satisfying \eqref{m8:s:s45:eq5.1}. Its character is \eqref{m8:s:s45:eq5.3}.
\end{proposition}
\begin{proof}

We record all diagonal stabilizer directions and their characters.

\paragraph{Projection to the diagonal group}\label{m8:s:s45:sub5.1}

Let \(\mathcal H_x\subset \operatorname{GL}_2(F)\) consist of the matrices \(h\) for which there exist \(\alpha,\delta\ne0\) and \(\gamma,\phi,p,q\in F\) satisfying
\begin{equation}\label{m8:s:s45:eq5.1}
{
\begin{aligned}
A^{-1}hA&=\begin{pmatrix}\alpha&0\\\gamma&\delta\end{pmatrix},&
B^{-1}hB&=\begin{pmatrix}\alpha&\phi\\0&\delta\end{pmatrix},\\
ChC^{-1}&=\begin{pmatrix}\alpha&\alpha p\\0&\delta\end{pmatrix},&
EhE^{-1}&=\begin{pmatrix}\delta&\delta q\\0&\alpha\end{pmatrix}.
\end{aligned}}
\end{equation}
For a given \(h\), these parameters are unique. This is the exact projection of \(K_x\) to the diagonal \(\operatorname{GL}_2\).

Indeed, the stabilizer case of \eqref{m8:s:s45:eq3.8} requires
\[
\pi_L(j)=l^{-1}D(h)l,
\qquad
\pi_L^-(w^{-1}jw)=r_0D(h)r_0^{-1}.
\]
In \eqref{m8:s:s45:eq3.6}, this forces \(\omega=\det h=\alpha\delta=\epsilon\zeta\). Taking determinants of the two matrices in \eqref{m8:s:s45:eq3.7} gives \(\epsilon=\alpha\) and \(\zeta=\delta\). The remaining equations are exactly \eqref{m8:s:s45:eq5.1}. The lower-block equations follow from these equations and the common similitude.

\paragraph{All stabilizer matrices}\label{m8:s:s45:sub5.2}

For \(h\in\mathcal H_x\) and arbitrary \(y,r\in F\), put
\[
Z_h(r)=\begin{pmatrix}0&p\\q&r\end{pmatrix},\qquad
M_h(r)=\begin{pmatrix}\gamma/\alpha&r\\0&\phi/\delta\end{pmatrix}.
\]
Then all elements of \(K_x\), and only those elements, are
\begin{equation}\label{m8:s:s45:eq5.2}
{
\begin{aligned}
u&=u\left(0,\kappa y,
       \kappa hAZ_h(r)B^{-1}h^{-1},0\right),\\
v&=u\left(-C^tM_h(r)^tE^{-t},0,
          -\frac{\det E}{s}yF_0,0\right),\\
k&=(u,v,h).
\end{aligned}}
\end{equation}
These are the \(u(X,y,Z,z)\)-coordinates of every stabilizer
element, with both \(z\)-coordinates zero. The quadratic entry
in each of \(u\) and \(v\) is prescribed by \eqref{m8:s:p:U}.

We recover these matrices from the projection equations and
\eqref{m8:s:c:normal-transporter}. Conversely, substitution verifies
every triple in the normal-form equation. The two free kernel
coordinates $y,r$ therefore complete the stabilizer. Expanded inverse-block recovery is
in \texttt{anc/model8\_calculations.pdf}, calculation M8-C02,
under ``All stabilizer matrices''.

The complete character is
\begin{equation}\label{m8:s:s45:eq5.3}
{
\begin{aligned}
\eta_K(u,v,h)
&=\psi\left(2\kappa y-\operatorname{tr}(T_0M_h(r))\right)\\
&=\psi\left(2\kappa y
 -(T_0)_{11}\frac\gamma\alpha
 -(T_0)_{21}r
 -(T_0)_{22}\frac\phi\delta\right).
\end{aligned}}
\end{equation}
This includes every contribution from the diagonal stabilizer's possible unipotent directions. Its restriction to \(h=I_2\) recovers \eqref{m8:s:s45:eq4.5}, so triviality on \(K_x\) is impossible.

\paragraph{Additional diagonal unipotents}\label{m8:s:s45:sub5.3}

Define four rational lines
\begin{equation}\label{m8:s:s45:eq5.4}
\mathscr A_1=BF e_1,\quad \mathscr A_2=C^{-1}Fe_1,
\qquad
\mathscr D_1=AFe_2,\quad \mathscr D_2=E^{-1}Fe_1.
\end{equation}
Equations \eqref{m8:s:s45:eq5.1} say that \(h\) acts with eigenvalue \(\alpha\) on the two \(\mathscr A\)-lines and eigenvalue \(\delta\) on the two \(\mathscr D\)-lines. There are exactly three possibilities:

\begin{enumerate}
\item If \(\mathscr A_1=\mathscr A_2\), \(\mathscr D_1=\mathscr D_2\), and these two lines are distinct, \(\mathcal H_x\) is the split torus acting by separate scalars on the two lines.
\item If all four lines coincide, \(\mathcal H_x\) is the product of scalar matrices and the one-parameter unipotent group fixing that line and acting trivially on its quotient.
\item Otherwise \(\mathcal H_x\) consists only of scalar matrices.
\end{enumerate}

For justification, a nonscalar matrix with distinct eigenvalues has exactly its two eigenlines; a nonscalar matrix with repeated eigenvalue has exactly one eigenline. These alternatives applied to \eqref{m8:s:s45:eq5.4} give the list over \(F\).

An explicit nonscalar unipotent example, including its exact
character $\psi(2y/t-2k)$, is in
\texttt{anc/model8\_calculations.pdf}, calculation M8-C02,
under ``Additional diagonal unipotents''. The $y$-direction in \eqref{m8:s:s45:eq5.3} proves vanishing
at every point.

\end{proof}

\subsubsection{General S6 full stabilizer catalogue}
\paragraph{Full \(S_6\) projection and formula}\label{m8:s:v:s6-1}

For a stabilizer of \(lw_6r\), \eqref{m8:s:c:normal-equations} with \(l'=l,r'=r\) gives
\[
\pi_L(j)=l^{-1}D(h)l,
\qquad
\pi_L^-(w_6jw_6)=rD(h)r^{-1},
\qquad \omega=\det h.
\]
The determinants of the two leading blocks in \eqref{m8:s:v:e4-3} force \(\epsilon=\alpha\) and \(\zeta=\delta\). Hence the complete projection equations are
\begin{equation}
\begin{aligned}
A^{-1}hA&=\begin{pmatrix}\alpha&0\\\gamma&\delta\end{pmatrix},&
B^{-1}hB&=\begin{pmatrix}\alpha&\phi\\0&\delta\end{pmatrix},\\
ChC^{-1}&=\begin{pmatrix}\alpha&\alpha q_0\\0&\delta\end{pmatrix},&
EhE^{-1}&=\begin{pmatrix}\alpha&0\\-\delta r_0&\delta\end{pmatrix},
\end{aligned}
\label{m8:s:v:e6-1}
\end{equation}
with \(\alpha\delta\ne0\). This describes every determinant and sign locus.

Set
\[
Z_0=\begin{pmatrix}0&q_0\\r_0&z\end{pmatrix},
\qquad
Z_1=\begin{pmatrix}0&-\gamma/\alpha\\\phi/\delta&py-z\end{pmatrix}.
\]
For every \(h\) satisfying \eqref{m8:s:v:e6-1}, and arbitrary \(p,y,z\), the full stabilizer element is
\begin{equation}
\begin{aligned}
u&=u(hApE_{22}B^{-1}h^{-1},\kappa y,
       \kappa hAZ_0B^{-1}h^{-1},0),\\
v&=u(-C^tyE_{22}E^{-t},-\ell p,
       \ell C^tZ_1E^{-t},0),\\
k&=(u,v,h).
\end{aligned}
\label{m8:s:v:e6-2}
\end{equation}
Conversely every full stabilizer element has this form. In particular the right factor was obtained by multiplying
\[
\bigl(\pi_L^-(w_6jw_6)(w_6jw_6)^{-1}\bigr)^t
=u(-yE_{22},-p,Z_1,0)
\]
after the determinant constraints above were imposed.

Taking traces gives exactly the compatibility equations
\eqref{m8:s:v:e6-3}.

\subsection{Terminal factorization and stabilizer calculations}\label{m8:supp:terminal}

\paragraph{Unipotent coordinates and matrix identities}\label{m8:s:t:sub1}

Use the split datum, $J$, $A^\#$ and $Q$ of
\S\ref{m8:s:sec:preliminaries}. Every block has size two. Retain
$W=\left(\begin{smallmatrix}0&1\\1&0\end{smallmatrix}\right)$ and
$\Sigma=\operatorname{diag}(1,-1)$, and put $A^*=A/\det A$ for
invertible $A$. The complete matrix $u(X,y,Z,z)$ is
\eqref{m8:s:p:U}; its character and Levi notation are

\begin{equation}\label{m8:s:t:eq1.2}
\xi(u(X,y,Z,z))=\psi(\operatorname{tr}X+2y),\qquad
\ell(A,B,t)=\operatorname{diag}(A,B,tB^*,tA^*).
\end{equation}

\subparagraph{Upper unitriangular isometries}\label{m8:s:t:sub1.1}
For a general block upper-unitriangular matrix, applying the
isometry equations to blocks $(2,4)$, $(3,3)$, $(3,4)$, and $(4,4)$
gives, successively,
\[
 F_0=-X^\#,\quad Y=yI_2,\quad H=Z^\#-yX^\#,
 \quad T=zI_2+\tfrac12(XZ^\#-ZX^\#).
\]
Conversely these equations satisfy every block of $g^tQg=Q$.
They give exactly the ten affine coordinates of
\eqref{m8:s:p:U}.
The entrywise form calculation is preserved in
\texttt{anc/model8\_calculations.pdf}, calculation M8-C03.

\subparagraph{Group law, inversion, and Levi conjugation}\label{m8:s:t:sub1.2}
The block product gives
\begin{equation}\label{m8:s:t:eq1.3}
\begin{aligned}
&(X,y,Z,z)(X',y',Z',z')
 =(X+X',y+y',Z+Z'+y'X,z+z'+c),\\
c&=\tfrac12\bigl(\operatorname{tr}(XZ'^\#)-\operatorname{tr}(ZX'^\#)
-y'\operatorname{tr}(XX'^\#)\bigr).
\end{aligned}
\end{equation}
The $(1,3)$ block gives the third coordinate and half the trace
of the $(1,4)$ block gives $c$. Hence
\begin{equation}\label{m8:s:t:eq1.4}
u(X,y,Z,z)^{-1}=u(-X,-y,-Z+yX,-z).
\end{equation}
The first two coordinates add, so $\xi$ is a character. Levi conjugation is
\begin{equation}\label{m8:s:t:eq1.5}
\operatorname{Ad}(\ell(A,B,t))u(X,y,Z,z)
=u\left(AXB^{-1},\frac{\det B}{t}y,
\frac{\det B}{t}AZB^{-1},\frac{\det A}{t}z\right).
\end{equation}
These are identities of the complete matrices in
\eqref{m8:s:p:U}; in particular the inverse includes $yX$.

\subparagraph{Equivariance conventions}\label{m8:s:t:sub1.3}

Let \(D(h)=\operatorname{diag}(h,h,h,h)\). Retain

\[
K=(U\times U)\rtimes\mathrm{GL}_2(F),\qquad
(u,v,h)x=uD(h)xD(h)^{-1}v^t,
\]

\begin{equation}\label{m8:s:t:eq1.6}
\eta_K(u,v,h)=\xi(u)\xi(v),\qquad
\iota(u,v,h)=(v,u,h^{-t}).
\end{equation}

The action of \(h\) on the second \(U\) is by \(D(h)^{-t}\). Transpose satisfies \((kx)^t=\iota(k)x^t\), and \(\eta_K\circ\iota=\eta_K\). The convention throughout is

\begin{equation}\label{m8:s:t:eq1.7}
T(f\circ k^{-1})=\eta_K(k)T(f).
\end{equation}

The extension has \(\widetilde\eta_K(\tau)=-1\). 

\paragraph{Terminal rank geometry}\label{m8:s:t:sub2}
\label{m8:s:t:sub2.1}\label{m8:s:t:sub2.2}
The inverse-corner argument of \S\ref{m8:s:sec:ranks} gives
\begin{equation}\label{m8:s:t:eq2.1}
 r_{m,n}=m+n-8+r_{8-m,8-n}.
\end{equation}
Thus the terminal rows of Table~\ref{m8:s:r:table} have the
following equivalent block conditions. Each listed pivot is invertible.
\begin{equation}\label{m8:s:t:eq2.4}
\begin{array}{c|l|l}
 i&\text{zero blocks}&\text{pivots}\\\hline
9&x_{34},x_{44}&x_{43},x_{24}\\
13&x_{43},x_{44}&x_{34},x_{42}\\
14&x_{33},x_{34},x_{43},x_{44}&x_{24},x_{42}\\
17&x_{24},x_{34},x_{42},x_{43},x_{44}&x_{33},x_{14},x_{41}\\
18&x_{24},x_{33},x_{34},x_{42},x_{43},x_{44}&x_{14},x_{41},x_{23},x_{32}.
\end{array}
\end{equation}
For example, in row 9 the zero blocks give $r_{4,2}=0$,
$x_{43}$ gives $r_{2,4}=2$, and the complementary identity
forces $r_{6,2}=2$, hence invertibility of $x_{24}$.
The same argument, its transpose, and the invertibility of the
last block row/column give all five rows; conversely they recover
every rank in Table~\ref{m8:s:r:table}. Complete block-by-block
rank verification is in \texttt{anc/model8\_calculations.pdf},
calculation M8-C03.

\subparagraph{The locally closed families}\label{m8:s:t:sub2.4}
Put $\Omega_b=\{r_{6,2}=2\}$, an open transpose-stable set,
and intersect it with the closed equations
$x_{34}=x_{44}=0$, respectively $x_{43}=x_{44}=0$.
By the table their two loci satisfy
\begin{equation}\label{m8:s:t:eq2.5}
 \mathcal B_L=S_9\sqcup S_{14},\quad
 \mathcal B_R=S_{13}\sqcup S_{14},\quad
 \mathcal B_L\cap\mathcal B_R=S_{14}.
\end{equation}
For the intersection, the invertible $x_{24}$ already accounts
for $r_{6,4}=2$; this forces $x_{33}=0$. Consequently
$\mathcal B=\mathcal B_L\cup\mathcal B_R$ is closed in $\Omega_b$,
and $S_9,S_{13}$ are separate clopen pieces of
$\mathcal B\setminus S_{14}$. Also
\begin{equation}\label{m8:s:t:eq2.6}
 \mathcal D=\{r_{6,2}=0\}=S_{17}\sqcup S_{18}
\end{equation}
is closed, with $S_{18}=\mathcal D\cap\{x_{33}=0\}$ closed in it.
These are the locally closed spaces used in the support sequences.

\paragraph{Rational factorization on the whole stratum}\label{m8:s:t:sub3}
\begin{proposition}\label{m8:s:t:factorization}
For \(i\in\{9,13,14,17,18\}\), multiplication gives an algebraic isomorphism over \(F\), \(\mathcal V_i\times L\times U\simeq S_i\), where \(\mathcal V_i=U\cap w_iUw_i^{-1}\). Its inverse has rational formulas with denominators nonzero on the whole stratum. The Levi parameter of a normal form is constant on its \(U\times U\)-orbit.
\end{proposition}
\begin{proof}\leavevmode

\subparagraph{Representatives and Levi normalization}\label{m8:s:t:sub3.1}

The representatives, with every block specified, are

\[
w_9=\begin{pmatrix}0&I&0&0\\0&0&0&W\\W&0&0&0\\0&0&I&0\end{pmatrix},\qquad
w_{13}=w_9^t=\begin{pmatrix}0&0&W&0\\I&0&0&0\\0&0&0&I\\0&W&0&0\end{pmatrix},
\]

\begin{equation}\label{m8:s:t:eq3.1}
w_{14}=\begin{pmatrix}0&0&W&0\\0&0&0&W\\W&0&0&0\\0&W&0&0\end{pmatrix},\quad
w_{17}=\begin{pmatrix}0&0&0&W\\0&I&0&0\\0&0&I&0\\W&0&0&0\end{pmatrix},\quad
w_{18}=\begin{pmatrix}0&0&0&W\\0&0&W&0\\0&W&0&0\\W&0&0&0\end{pmatrix}.
\end{equation}

Equivalently, their eight rows of nonzero entries, specified by the row occupied by each successive column, are

\[
\begin{array}{c|c}
9&(6,5,1,2,7,8,4,3)\\
13&(3,4,8,7,2,1,5,6)\\
14&(6,5,8,7,2,1,4,3)\\
17&(8,7,3,4,5,6,2,1)\\
18&(8,7,6,5,4,3,2,1).
\end{array}
\]

All nonzero scalar entries are \(+1\). The scalar-coordinate reconstruction is in \S\ref{m8:s:t:expanded-representatives}. Direct multiplication verifies \(w_i^tQw_i=Q\), \(\det w_i=1\), and the rank matrices in Table~\ref{m8:s:r:table}.

The full normalization formulas are as follows; the third Levi parameter stays \(t\):

\begin{equation}\label{m8:s:t:eq3.2}
\begin{array}{c|c}
i& w_i\ell(A,B,t)w_i^{-1}=\ell(A',B',t)\\\hline
9&(A',B')=(B,\ tWAW/\det A)\\
13&(A',B')=(tWBW/\det B,\ A)\\
14&(A',B')=(tWBW/\det B,\ tWAW/\det A)\\
17&(A',B')=(tWAW/\det A,\ B)\\
18&(A',B')=(tWAW/\det A,\ tWBW/\det B).
\end{array}
\end{equation}

Thus each \(w_i\) normalizes the fixed Levi. The matrices \(w_i\) specify normal forms; the orbit action throughout remains the action of \(K\).

\subparagraph{Elimination over the ground field}\label{m8:s:t:sub3.2}
Let $N^+$ be the block upper-unitriangular group in $\GL_8$,
$L_{\rm GL}$ the block diagonal group, and $p_i(a)$ the column
of the nonzero block in row $a$ of $w_i$. The permitted left
positions are $a<b$ with $p_i(a)<p_i(b)$:
\[
\begin{array}{c|c|c}
i&(p_i(1),p_i(2),p_i(3),p_i(4))&
 C_i^{\rm GL}=N^+\cap w_iN^+w_i^{-1}\\\hline
9&(2,4,1,3)&12,14,34\\
13&(3,1,4,2)&13,23,24\\
14&(3,4,1,2)&12,34\\
17&(4,2,3,1)&23\\
18&(4,3,2,1)&\text{none}.
\end{array}
\]
Process rows $a=4,3,2,1$. If row $b>a$ has been reduced to its
pivot in column $p_i(b)>p_i(a)$, eliminate the current entry by
\begin{equation}\label{m8:s:t:eq3.3}
R_a\longleftarrow R_a-Y_{a,p_i(b)}Y_{b,p_i(b)}^{-1}R_b.
\end{equation}
The suffix ranks then force every entry to the right of $p_i(a)$
to vanish: a remaining entry would add a row independent of the
processed pivot rows. At $p_i(a)$ the suffix rank increases by two,
so
\begin{equation}\label{m8:s:t:eq3.4}
P_a=Y_{a,p_i(a)}\in\GL_2(F).
\end{equation}
Eliminate entries to its left by
\begin{equation}\label{m8:s:t:eq3.5}
C_j\longleftarrow C_j-C_{p_i(a)}P_a^{-1}Y_{a,j},\qquad j<p_i(a).
\end{equation}
These right lower-unitriangular operations preserve all processed
rows, whose distinct pivots lie elsewhere. Both sets of operations
preserve the terminal ranks. Thus every point, over $F$ or any
extension, has a factorization
\begin{equation}\label{m8:s:t:eq3.6}
x=c\,d\,w_i\,v^t,\qquad
c\in C_i^{\rm GL},\quad d\in L_{\rm GL},\quad v\in N^+.
\end{equation}
\label{m8:s:t:sub3.3}
It is unique: comparison of diagonal blocks first identifies $d$;
the remaining ambiguity lies in
$C_i^{\rm GL}\cap w_i(N^+)^tw_i^{-1}=1$, since their permitted
positions have opposite signs in the $p_i$ ordering.

\subparagraph{Flag minors and denominators}\label{m8:s:t:sub3.4}
Let $\Delta_a$ use row blocks $a,\ldots,4$ and column blocks
$\{p_i(a),\ldots,p_i(4)\}$, sorted increasingly, and put
$\Delta_5=1$. Restricting the two unitriangular factors of
\eqref{m8:s:t:eq3.6} to these rows and columns gives
\begin{equation}\label{m8:s:t:eq3.7}
\Delta_a=\prod_{b=a}^4\det P_b,\qquad
\det P_a=\Delta_a/\Delta_{a+1}.
\end{equation}
Because all blocks have size two, the block permutation has scalar
sign $+1$. The pivots are invertible by elimination, so these
original minors are nonzero throughout the stratum and every
division in the inverse algorithm is valid on its whole domain. The expanded column-set table is in
\texttt{anc/model8\_calculations.pdf}, M8-C03,
``Flag minors and denominators''.

\subparagraph{Orthogonality of the factors}\label{m8:s:t:sub3.5}

Define the algebraic group automorphism

\[
\Theta(g)=Qg^{-t}Q.
\]

It preserves block diagonals and block upper-unitriangular matrices: inverse transpose reverses upper to lower, and the block antidiagonal \(Q\) reverses the order back. Since \(\Theta(w_i)=w_i\), it also preserves \(C_i^{\mathrm{GL}}\). Moreover

\[
\Theta(v^t)=\Theta(v)^t,\qquad
\Theta(x)=\mu(x)^{-1}x.
\]

Apply \(\Theta\) to \eqref{m8:s:t:eq3.6}. A second factorization of \(\mu(x)^{-1}x\) is obtained by leaving \(c,v\) unchanged and replacing \(d\) by \(\mu(x)^{-1}d\). Both expressions belong to the same unique GL factorization. Thus

\begin{equation}\label{m8:s:t:eq3.8}
\Theta(c)=c,\qquad \Theta(v)=v,\qquad
\Theta(d)=\mu(x)^{-1}d.
\end{equation}

The first two equalities say \(c^tQc=v^tQv=Q\). \S\ref{m8:s:t:sub1.1} identifies these factors with the \(U\) coordinates. The diagonal equality says

\[
d_4=\mu(x)d_1^*,\qquad d_3=\mu(x)d_2^*,
\]

hence \(d=\ell(d_1,d_2,\mu(x))\). The scalar \(\mu(x)^{-1}\) affects only the diagonal factor. Both unitriangular factors are uniquely determined.

Set

\[
\mathcal V_i=U\cap w_iUw_i^{-1}.
\]

We have proved the asserted global algebraic isomorphism

\begin{equation}\label{m8:s:t:eq3.9}
{\mathcal V_i\times L\times U\ \simeq\ S_i,
\qquad(c,\ell,v)\mapsto c\ell w_iv^t.}
\end{equation}

The inverse is regular on the entire locally closed stratum by the nonvanishing denominators above.  

\subparagraph{Explicit inverse recovery}\label{m8:s:t:sub3.6}
Once the row/column algorithm recovers $c,A,B,t$, we recover
the final radical on the entire stratum by
\begin{equation}\label{m8:s:t:eq3.15}
 v^t=w_i^{-1}\ell(A,B,t)^{-1}c^{-1}x\in U^t.
\end{equation}
For the five expanded block formulas see
\texttt{anc/model8\_calculations.pdf}, calculation M8-C03,
under ``Explicit inverse formulas''. They implement
\eqref{m8:s:t:eq3.3}--\eqref{m8:s:t:eq3.8} using the nonzero
denominators established above.

\subparagraph{Levi parameters on a two-sided orbit}\label{m8:s:t:sub3.7}

A further property is needed, beyond uniqueness within the restricted left factor \(\mathcal V_i\). Suppose arbitrary \(u,v\in U\) satisfy

\[
u\ell w_iv^t=\ell'w_i.
\]

Then

\[
\ell'^{-1}u\ell=w_iv^{-t}w_i^{-1}.
\]

Comparing diagonals gives \(\ell'^{-1}\ell=1\), because
\(w_i\) is a block permutation with invertible blocks. This
proves Levi uniqueness for all five representatives, as needed
for the residual orbit relations and full stabilizers.

\end{proof}
\paragraph{Complete unipotent stabilizers}\label{m8:s:t:sub4}
\begin{proposition}\label{m8:s:t:unipotent-stabilizers}
The complete unipotent stabilizers are the groups in \eqref{m8:s:t:eq4.10}, with all equations and coordinates given in \eqref{m8:s:t:eq4.2}--\eqref{m8:s:t:eq4.9}. Their dimensions are \(5,5,6,9,10\), respectively. The stabilizer for \(i=17\) is Heisenberg; that for \(i=18\) is \(U\).
\end{proposition}
\begin{proof}
At $x=\ell w_i$, put $n=\ell^{-1}u\ell$. The full equation is
\begin{equation}\label{m8:s:t:eq4.1}
 u\ell w_iv^t=\ell w_i\iff
 n\in\mathcal N_i=U\cap w_iU^tw_i^{-1},\quad
 v=r_i(n)=(w_i^{-1}n^{-1}w_i)^t.
\end{equation}
The input is the entire ten-coordinate matrix
\eqref{m8:s:p:U}, with $H=Z^\#-yX^\#$ and
$T_0=zI+\tfrac12(XZ^\#-ZX^\#)$. Conjugate by the five
specified block permutations and set every upper off-diagonal
block to zero. The complete equations are
\begin{equation}\label{m8:s:t:eq4.2}
\begin{array}{c|c|c}
 i&\mathcal N_i&\mathcal V_i=U\cap w_iUw_i^{-1}\\\hline
9&u(0,y,Z,0)&u(X,0,0,z)\\
13&u(X,0,0,z)&u(0,y,Z,0)\\
14&u(0,y,Z,z)&u(X,0,0,0)\\
17&u(X,0,Z,z)&u(0,y,0,0)\\
18&U&1.
\end{array}
\end{equation}
For row 9 the forbidden blocks contain $X^\#$ and $T_0$,
forcing $X=z=0$; row 13 forces $y=Z=0$; rows 14,17 force
only $X=0$, respectively $y=0$; row 18 has no forbidden block.
Substitution makes every forbidden block zero, proving sufficiency
as group equations. Testing the lower blocks instead gives the
last column. The five complete conjugated matrices and the
corresponding equation ideals are in
\texttt{anc/model8\_calculations.pdf}, calculation M8-C03,
under ``The full stabilizer equations''.
\begin{equation}\label{m8:s:t:eq4.7}
 \dim\mathcal V_i=(5,5,4,1,0),\quad
 \dim\mathcal N_i=(5,5,6,9,10).
\end{equation}
\subparagraph{The nonlinear group structures}\label{m8:s:t:sub4.3}
The group law \eqref{m8:s:t:eq1.3} makes the first three
$\mathcal N_i$ additive. On $\mathcal N_{17}$ it is the
Heisenberg law with alternating pairing
\begin{equation}\label{m8:s:t:eq4.8}
 \omega((X,Z),(X',Z'))=
 \operatorname{tr}(XZ'^\#)-\operatorname{tr}(ZX'^\#).
\end{equation}
The pairing $\operatorname{tr}(AB^\#)=ah+de-bg-cf$ on
$A=(\begin{smallmatrix}a&b\\c&d\end{smallmatrix})$,
$B=(\begin{smallmatrix}e&f\\g&h\end{smallmatrix})$ is nondegenerate;
hence $\omega$ is nondegenerate and the center is the $z$-line.

\subparagraph{Right stabilizer coordinates}\label{m8:s:t:sub4.4}
Applying the complete inverse \eqref{m8:s:t:eq1.4}, including
its $yX$ correction, gives
\[
\begin{array}{c|c}
i&(X_v,y_v,Z_v,z_v)\\\hline
9&(-WZ^t,0,0,-y)\\
13&(0,-z,-X^tW,0)\\
14&(0,-z,-WZ^tW,-y)\\
17&(\Sigma ZJ,0,-\Sigma XJ,-z)\\
18&(\Sigma X\Sigma,-y,-\Sigma Z\Sigma,-z).
\end{array}
\]
Every row satisfies the whole matrix identity
\begin{equation}\label{m8:s:t:eq4.9}
 nw_ir_i(n)^t=w_i.
\end{equation}
In row 18 the inverse correction cancels during coordinate
recovery. Thus
\begin{equation}\label{m8:s:t:eq4.10}
 (U\times U)_{\ell w_i}
 =\{(\ell n\ell^{-1},r_i(n)):n\in\mathcal N_i\}.
\end{equation}
Inversion and transpose are both anti-homomorphisms, so this
is a group isomorphism, with inverse recovered from the first
component. It proves completeness on every $F$-point.
\end{proof}

\paragraph{Residual actions, full stabilizers, and modular functions}\label{m8:s:t:sub6}
\begin{proposition}\label{m8:s:t:full-stabilizers}
The residual rational orbit relation on each of the five strata is simultaneous conjugation in \eqref{m8:s:t:eq6.2}, with \(t\) fixed. The full stabilizer is \eqref{m8:s:t:eq6.3}, with action \eqref{m8:s:t:eq6.4}. Both \(K\) and all these stabilizers are unimodular.
\end{proposition}
\begin{proof}\leavevmode

\subparagraph{Rational residual equivalence}\label{m8:s:t:sub6.1}

Set

\[
(\varepsilon_A,\varepsilon_B)=
\begin{cases}
(0,1)&i=9,\\
(1,0)&i=13,\\
(1,1)&i=14,\\
(1,0)&i=17,\\
(1,1)&i=18.
\end{cases}
\]

Inspection of the nonzero blocks of \(D(h)\ell w_iD(h)^{-1}\) gives

\begin{equation}\label{m8:s:t:eq6.1}
A_h=hAW^{\varepsilon_A}h^{-1}W^{\varepsilon_A},\qquad
B_h=hBW^{\varepsilon_B}h^{-1}W^{\varepsilon_B},\qquad t_h=t.
\end{equation}

The lower two Levi blocks satisfy the same formula automatically, since \(\det A_h=\det A\), \(\det B_h=\det B\), and the appropriate pivot \(W\)'s match.

With

\[
M_1=AW^{\varepsilon_A},\qquad M_2=BW^{\varepsilon_B},
\]

this is exactly

\begin{equation}\label{m8:s:t:eq6.2}
(M_1,M_2,t)\longmapsto(hM_1h^{-1},hM_2h^{-1},t).
\end{equation}

If two normal forms are \(K\)-equivalent, apply \eqref{m8:s:t:eq6.1} to the diagonal factor of a transporter. \S\ref{m8:s:t:sub3.7} then forces equality of the resulting Levi parameters. Conversely, simultaneous conjugation is realized by the element \((1,1,h)\). This proves \eqref{m8:s:t:eq6.2} over $F$.

\subparagraph{\texorpdfstring{Full \(K\)-stabilizers}{Full K-stabilizers}}\label{m8:s:t:sub6.2}

Let

\[
H_i(A,B)=Z_{\mathrm{GL}_2}(M_1,M_2).
\]

If \((u,v,h)\) fixes \(\ell w_i\), Levi uniqueness forces \(h\in H_i\). Such an \(h\) already fixes the normal-form matrix without unipotent corrections. The remaining pair is precisely \eqref{m8:s:t:eq4.10}. Hence

\begin{equation}\label{m8:s:t:eq6.3}
{K_{\ell w_i}
=\{(\ell n\ell^{-1},r_i(n),h):n\in\mathcal N_i,\ h\in H_i(A,B)\}
\simeq\mathcal N_i\rtimes H_i(A,B).}
\end{equation}

The action in this semidirect product is

\[
\rho_h(n)=\operatorname{Ad}(w_iD(h)w_i^{-1})n.
\]

Indeed \(h\) fixes \(\ell w_i\), so \(\ell^{-1}D(h)\ell=w_iD(h)w_i^{-1}\). Put

\[
a_0=W^{\varepsilon_A}hW^{\varepsilon_A},\qquad
b_0=W^{\varepsilon_B}hW^{\varepsilon_B}.
\]

Then \(w_iD(h)w_i^{-1}=\ell(a_0,b_0,\det h)\), and \eqref{m8:s:t:eq1.5} yields

\begin{equation}\label{m8:s:t:eq6.4}
{\rho_h:u(X,y,Z,z)\mapsto
u(a_0Xb_0^{-1},y,a_0Zb_0^{-1},z).}
\end{equation}

Its restriction to each displayed \(\mathcal N_i\) is the full stabilizer action. This includes the unipotent elements in every repeated-root centralizer.

The stabilizer character is

\[
\eta_K(\ell n\ell^{-1},r_i(n),h)=\psi(\Lambda_i(n)),
\]

with trivial character on \(h\).

\subparagraph{\texorpdfstring{Centralizers over \(F\)}{Centralizers over F}}\label{m8:s:t:sub6.3}

For a nonscalar two-by-two matrix \(M\), a cyclic basis proves

\begin{equation}\label{m8:s:t:eq6.5}
Z_{M_2(F)}(M)=F[M].
\end{equation}

This is two-dimensional and commutative. Therefore the common centralizer types are exactly:

\begin{center}\small
\begin{tabular}{@{}>{\raggedright\arraybackslash}p{.47\linewidth}>{\raggedright\arraybackslash}p{.47\linewidth}@{}}
\toprule
Situation & Algebraic group \(H_i\) \\
\midrule
Both \(M_1,M_2\) scalar & \(\mathrm{GL}_2\) \\
They commute; a nonscalar member has distinct split roots & \(\mathbb G_m\times\mathbb G_m\) \\
They commute; a nonscalar member has irreducible quadratic polynomial & \(\operatorname{Res}_{E/F}\mathbb G_m\), with \(E/F\) the quadratic field \\
They commute; a nonscalar member has a repeated root & \((F[\epsilon]/\epsilon^2)^\times\simeq\mathbb G_m\times\mathbb G_a\) \\
They do not commute & scalar \(\mathbb G_m\) \\
\bottomrule
\end{tabular}\end{center}

For the last row, a nonscalar element in the common centralizer would put both \(M_1,M_2\) in its commutative two-dimensional polynomial algebra, a contradiction. In the repeated-root row, writing an element as \(a+b\epsilon\), the isomorphism is \((a+b\epsilon)\mapsto(a,b/a)\), with \(a\ne0\).

All these algebraic groups are connected. The nonsplit quadratic torus has rational-point group $E^\times$.

\subparagraph{Unimodularity}\label{m8:s:t:sub6.4}
In the coordinate order $X,y,Z,z$, left and right translations
in \eqref{m8:s:t:eq1.3} have triangular Jacobian with diagonal one.
Thus $dX\,dy\,dZ\,dz$ is bi-invariant, as are the corresponding
product measures on the additive and Heisenberg intersections.
For \eqref{m8:s:t:eq6.4},
\[
 \det(X\mapsto a_0Xb_0^{-1})
 =(\det a_0)^2(\det b_0)^{-2}=1.
\]
The same determinant calculation applies to $Z$, while $y,z$
are fixed. It also applies on each $\mathcal N_i$, where the
matrix coordinates are either present in full or zero. Since the
centralizer groups are commutative or $\GL_2$, they are unimodular. Therefore
\begin{equation}\label{m8:s:t:eq6.6}
\Delta_{K_x}=1.
\end{equation}
The two diagonal actions on $U$ conjugate $X,Z$ and fix $y,z$,
so both preserve Haar measure and
\begin{equation}\label{m8:s:t:eq6.7}
\Delta_K=1.
\end{equation}
Consequently, at every rational point,
\begin{equation}\label{m8:s:t:eq6.8}
\Delta_{U\times U}|_{(U\times U)_x}\Delta_{(U\times U)_x}^{-1}=1,
\qquad \Delta_K|_{K_x}\Delta_{K_x}^{-1}=1.
\end{equation}

\subparagraph{The modulus of the full Levi action}\label{m8:s:t:sub6.5}

Let \(d_A=\det A\), \(d_B=\det B\). From \eqref{m8:s:t:eq1.5}, the determinants on \(X,y,Z,z\), respectively, are

\[
d_A^2d_B^{-2},\qquad d_B/t,\qquad
(d_B/t)^4d_A^2d_B^{-2},\qquad d_A/t.
\]

Their product is

\begin{equation}\label{m8:s:t:eq6.9}
{\det(\operatorname{Ad}\ell|_{\mathfrak u})=d_A^5d_Bt^{-6},\qquad
\text{Haar modulus}=|d_A|_F^5|d_B|_F|t|_F^{-6}.}
\end{equation}

For the diagonal actions considered above the Levi modulus equals one, and it also equals one on the compatible \(S_{18}\) locus \(A=-B,t=\det B\).

\end{proof}

\paragraph{\texorpdfstring{The compatible \(S_{18}\) chart and its density}{The compatible S18 chart and its density}}\label{m8:s:t:sub8}
\begin{proposition}\label{m8:s:t:chart}
The map \((M,v)\mapsto x(M)v^t\) in \eqref{m8:s:t:eq8.3} identifies \(\operatorname{GL}_2\times U\) with \(\mathcal C_{18}\). Transpose acts by \((M,v)\mapsto(M^t,\theta_{M^t}(v))\), where \(\theta_M\) is the anti-automorphism \eqref{m8:s:t:eq8.7}. Each \(\theta_M\) preserves \(\xi\) and product Haar measure.
\end{proposition}
\begin{proof}\leavevmode

\subparagraph{The chart and its inverse}\label{m8:s:t:sub8.1}
Since $\mathcal V_{18}=1$, rational factorization gives
\begin{equation}\label{m8:s:t:eq8.1}
 S_{18}\simeq L\times U.
\end{equation}
Compatibility forces $A=-B$, $t=\det B$. Setting $M=BW$
gives exactly $x(M)$ of \eqref{m8:s:t:eq8.2}, with
$\mu(x(M))=-\det M$ and $\det x(M)=(\det M)^4$.
Consequently \eqref{m8:s:t:eq8.3} is a global algebraic
isomorphism, with inverse
\begin{equation}\label{m8:s:t:eq8.4}
 M=x_{23}=-x_{14}=x_{32}=-x_{41},\qquad v^t=x(M)^{-1}x.
\end{equation}
Since $M$ is invertible, the inverse is regular, and the image
is smooth of dimension $4+10=14$. On the base,
\begin{equation}\label{m8:s:t:eq8.5}
 D(h)x(M)D(h)^{-1}=x(hMh^{-1}),\qquad x(M)^t=x(M^t).
\end{equation}
Scalars act ineffectively, and distinct scalar $M$ remain distinct base points.

\subparagraph{The fiber anti-automorphism}\label{m8:s:t:sub8.2}
Block reversal defines the algebraic anti-automorphism
\begin{equation}\label{m8:s:t:eq8.6}
 \theta_M(a)=(x(M)^{-1}ax(M))^t.
\end{equation}
For $C=MJ$, its complete coordinates are
\eqref{m8:s:t:eq8.7}; the conversion uses
$M^tJM=\det(M)J$. In particular
\begin{equation}\label{m8:s:t:eq8.8}
 \xi\circ\theta_M=\xi,\qquad
 \theta_{M^t}\theta_M=\theta_M\theta_{M^t}=1.
\end{equation}
The second identity follows from $x(M^t)=x(M)^t$.

\subparagraph{The four chart actions}\label{m8:s:t:sub8.3}
From $ax(M)=x(M)\theta_M(a)^t$ and the definitions,
\begin{align}
 (1,b,1):(M,v)&\mapsto(M,bv),\label{m8:s:t:eq8.10}\\
 (a,1,1):(M,v)&\mapsto(M,v\theta_M(a)),\label{m8:s:t:eq8.11}\\
 (1,1,h):(M,v)&\mapsto(hMh^{-1},D(h)^{-t}vD(h)^t),\label{m8:s:t:eq8.12}\\
 \tau:(M,v)&\mapsto(M^t,\theta_{M^t}(v)).\label{m8:s:t:eq8.13}
\end{align}
These give $\tau^2=1$ and residual matrix transpose.

\subparagraph{The density calculation}\label{m8:s:t:sub8.4}
The group law makes $dv=dX\,dy\,dZ\,dz$ bi-invariant.
In this coordinate order, the derivative of $\theta_M$ is block
triangular with diagonal blocks
$\operatorname{Ad}(C^{-1}),1,-\operatorname{Ad}(C^{-1}),1$.
Conjugation on $M_2$ has determinant one, and $(-1)^4=1$;
the $yX$ term is off the diagonal. Thus
\begin{equation}\label{m8:s:t:eq8.14}
 (\theta_M)_*dv=dv.
\end{equation}
Conjugation in the diagonal chart action fixes $y,z$ and preserves
both $dv$ and $\xi$ on $X,Z$. Thus every fiber multiplier in
the ordinary-distribution descent is one.
The expanded block identities, coordinate conversion and Jacobian
are in \texttt{anc/model8\_calculations.pdf}, calculation M8-C03,
under ``The compatible $S_{18}$ chart and its density''.
\end{proof}

\subparagraph{The algebraic transpose extension}\label{m8:s:t:sub10.3}
The automorphism $\iota(u,v,h)=(v,u,h^{-t})$ defines the split
algebraic group $\widetilde{\mathbf K}=
\mathbf K\rtimes\underline{\mathbb Z/2\mathbb Z}$, whose rational
points are $K\rtimes\langle\tau\rangle$ and whose multiplication is
\begin{equation}\label{m8:s:t:eq10.6}
(k,\epsilon)(k',\epsilon')
 =(k\iota^\epsilon(k'),\epsilon+\epsilon').
\end{equation}
It acts algebraically by $(k,\epsilon)x=k\tau^\epsilon(x)$, with
character $\eta_K(k)(-1)^\epsilon$. For the transporter
\eqref{m8:s:t:eq10.4},
\begin{equation}\label{m8:s:t:eq10.7}
b_x=k_x^{-1}\tau\in\widetilde K_x,\qquad
\widetilde\eta_K(b_x)=-1.
\end{equation}
The even subgroups have index two. Their modular characters are
one by \eqref{m8:s:t:eq6.6}--\eqref{m8:s:t:eq6.7}; positivity and
$b_x^2\in K_x$ then give
\begin{equation}\label{m8:s:t:eq10.8}
\Delta_{\widetilde K}=\Delta_{\widetilde K_x}=1.
\end{equation}
Odd vanishing on $\mathcal C_{18}$ is already proved by the exact
descent \eqref{m8:s:t:eq9.9} and rational conjugacy localization.

\subparagraph{Six one-parameter degenerations}\label{m8:s:t:sub12.1}
Put $e_X=u(X,0,0,0)$, $q_y=u(0,y,0,0)$ and
$c_z=u(0,0,0,z)$. The following exact product formulas give the
six degenerations. Each product preserves $Q$ and has determinant
one for every $r\in F$.
\begin{center}\small
\begin{tabular}{c p{.18\linewidth} p{.49\linewidth} cc}
\toprule
Family&Product&Decisive terminal blocks&$r\ne0$&$r=0$\\
\midrule
1&$c_r^tw_{14}$&$x_{43}=rW$, $x_{34}=x_{44}=0$&$S_9$&$w_{14}$\\
2&$w_{14}c_r$&Transpose of Family 1&$S_{13}$&$w_{14}$\\
3&$e_{rI}^tw_{17}$&$x_{24}=rW$, $x_{43}=-rI$, $x_{34}=x_{44}=0$&$S_9$&$w_{17}$\\
4&$w_{17}e_{rI}$&Transpose of Family 3&$S_{13}$&$w_{17}$\\
5&$e_{rI}^tw_{18}$&$x_{24}=rW$, $x_{42}=-rW$; terminal $4\times4$ block zero&$S_{14}$&$w_{18}$\\
6&$q_r^tw_{18}$&$x_{33}=rW$; deeper terminal blocks zero&$S_{17}$&$w_{18}$\\
\bottomrule
\end{tabular}\end{center}
The necessary and sufficient block conditions in
\eqref{m8:s:t:eq2.4} prove each stated rank profile for every
$r\ne0$, since $rI$ and $rW$ are then invertible. In the order
of the table, the complete rank matrices are
\begin{equation}\label{m8:s:t:eq12.1}
 (R_9,R_{13},R_9,R_{13},R_{14},R_{17})\qquad(r\ne0).
\end{equation}
The entries are those of Table~\ref{m8:s:r:table}.

\subparagraph{Whole-stratum closure relations}\label{m8:s:t:sub12.2}

Every point of a target stratum is \(c\ell w_iv^t\) by the global chart. Multiplication by fixed elements of \(P=LU\) on the left and \(P^t\) on the right preserves each terminal rank. Applying these fixed translates to a displayed family transports its limit to that arbitrary target point and keeps its nonzero-parameter members in the source stratum.

Thus, both for p-adic closure and for Zariski closure, the nontrivial inclusions among these five strata are

\[
\overline{S_9}\supset S_{14}\cup S_{17}\cup S_{18},\qquad
\overline{S_{13}}\supset S_{14}\cup S_{17}\cup S_{18},
\]

\begin{equation}\label{m8:s:t:eq12.2}
\overline{S_{14}}\supset S_{18},\qquad
\overline{S_{17}}\supset S_{18}.
\end{equation}

Rank semicontinuity excludes every other inclusion between distinct members of this five-element collection. Explicitly, \(S_{9}\) and \(S_{13}\) are incomparable because each has a rectangular rank smaller than the other's. \(S_{14}\) cannot specialize to \(S_{17}\) because \(r_{4,4}=0\) on \(S_{14}\) but equals two on \(S_{17}\). \(S_{17}\) cannot specialize to \(S_{14}\) because \(r_{6,2}=0\) on \(S_{17}\) but equals two on \(S_{14}\). \(S_{18}\) has the minimal ranks in the collection and cannot specialize to any of the other four. The remaining prohibited directions similarly require an increase in one of Table~\ref{m8:s:r:table}'s ranks.

\subparagraph{Topology of the five-stratum union}\label{m8:s:t:sub12.4}
For $X=\operatorname{diag}(a,b)$ the closed affine embedding
\begin{equation}\label{m8:s:t:eq12.7}
f(a,b)=e_X^tw_{18}
=\begin{pmatrix}
0&0&0&W\\0&0&W&XW\\0&W&0&0\\W&-X^\#W&0&0
\end{pmatrix}
\end{equation}
has $r_{6,2}=r_{2,6}=\operatorname{rank}X$. Thus
\begin{equation}\label{m8:s:t:eq12.8}
f^{-1}(\mathcal E)=(F^\times)^2\cup\{(0,0)\},\qquad
\mathcal E=\mathcal B\cup\mathcal D.
\end{equation}
This dense subset of $F^2$ is not open in its closure at zero.
It is also not locally compact there: a compact neighborhood would
contain $(a,b_n)$ for a fixed sufficiently small $a\ne0$ and
$b_n\ne0$ tending to zero, but would omit their ambient limit
$(a,0)$. Since the plane image is closed, $\mathcal E$ is neither
locally closed nor an l-space with the inherited topology. The
support argument uses the two families $\mathcal B,\mathcal D$
separately.

\subsection{Middle-cell rational flags, intersections and boundary calculations}\label{m8:supp:middle}

Throughout this section $n_+(s)=\begin{pmatrix}1&s\\0&1\end{pmatrix}$ and $n_-(s)=\begin{pmatrix}1&0\\s&1\end{pmatrix}$; $N_+$ and $N_-$ denote the corresponding one-dimensional unipotent subgroups. We write $G_8^\circ=\Geight$ and $I_i=P\cap w_iP^tw_i^{-1}$.

\paragraph{Rational parabolic cells and representatives}\label{m8:s:m:cells}
\begin{proposition}[Rational parabolic cells]\label{m8:s:m:prop:cells}
For $i\in\{10,11,12,15,16\}$, the representatives below satisfy $w_i\in\mathrm{SO}(Q)(F)$ and $S_i(F)=P(F)w_iP(F)^t$. Every point of $S_i(F)$ therefore has the form $u_0mw_inv_0^t$, with $u_0,v_0\in U(F)$ and $m,n\in M(F)$.
\end{proposition}
\par\smallskip\noindent\textit{Proof.}\par

\subparagraph{Monomial representatives}\label{m8:s:m:representatives}

We consider the five rank strata $S_{10},S_{11},S_{12},S_{15},S_{16}$, whose dimensions are respectively $25,24,23,23,22$. Their rank matrices are those of the complete rank table.

For the monomial representatives, the row positions of columns $1,\ldots,8$ are

\begingroup
\small
\setlength{\tabcolsep}{3pt}
\renewcommand{\arraystretch}{1.18}
\begin{longtable}{p{0.16\linewidth}p{0.36\linewidth}p{0.39\linewidth}}
\toprule
i & Row positions & Column coefficients \\
\midrule
\endfirsthead
i & Row positions & Column coefficients \\
\midrule
\endhead
10 & 3,5,1,7,2,8,4,6 & 1,1,1,1,1,1,1,1 \\
11 & 8,3,2,5,4,7,6,1 & 1,1,1,1,1,-1,-1,1 \\
12 & 8,6,2,4,5,7,3,1 & 1,1,1,1,1,-1,-1,1 \\
15 & 8,3,7,4,5,2,6,1 & 1,1,1,1,1,-1,-1,1 \\
16 & 8,6,7,5,4,2,3,1 & 1,1,1,1,1,-1,-1,1 \\
\bottomrule
\end{longtable}
\endgroup

For each one, direct multiplication gives
\[
w_i^tQw_i=Q,\qquad \det w_i=1,
\]
and direct corner-rank computation gives the displayed \(R_i\). Thus every representative belongs to the selected connected special orthogonal group.

\subparagraph{Relative flags and joint multiplicities}\label{m8:s:m:flags}

Let \(E_a\) be the span of the first a standard coordinate vectors. Let \(F_b^-\) be the span of the last $b$. Then
\[
r_{a,b}(x)
 =b-\dim\bigl(xF_b^-\cap E_{8-a}\bigr),
\qquad a,b\in\{2,4,6\}.
\]
This is simply rank--nullity for projection to the last a coordinates. Including the automatic dimensions at 0 and 8, the nine ranks determine all the relevant intersection dimensions of the two self-dual flags.

Taking successive differences gives the following relative block multiplicities, displayed in the original column-block order. Reversing their columns gives the order of the increasing opposite flag:
\[
J_{10}=\begin{pmatrix}0&1&1&0\\1&0&0&1\\1&0&0&1\\0&1&1&0\end{pmatrix},\quad
J_{11}=\begin{pmatrix}0&1&0&1\\1&0&1&0\\0&1&0&1\\1&0&1&0\end{pmatrix},
\]
\[
J_{12}=\begin{pmatrix}0&1&0&1\\0&1&0&1\\1&0&1&0\\1&0&1&0\end{pmatrix},\quad
J_{15}=\begin{pmatrix}0&0&1&1\\1&1&0&0\\0&0&1&1\\1&1&0&0\end{pmatrix},
\]
\[
J_{16}=\begin{pmatrix}0&0&1&1\\0&0&1&1\\1&1&0&0\\1&1&0&0\end{pmatrix}.
\]
Every nonzero joint graded piece is one-dimensional. The grading is paired with its opposite under \((i,j)\mapsto(5-i,5-j)\). Every graded piece is paired with a distinct opposite piece.

\subparagraph{Simultaneous Witt bases over the ground field}\label{m8:s:m:witt}
The displayed multiplicities satisfy Lemma~\ref{m8:simultaneous-witt}.
Apply it to the standard flag and the $x$-image of the opposite
flag, and to the corresponding pair for $w_i$. Matching their
Witt bases gives $p\in P(F)$ with
$w_i^{-1}p^{-1}x\in P(F)^t$. The original similitude remains in
this right factor. Hence $S_i(F)=P(F)w_iP(F)^t$, and Levi--unipotent
factorization yields the asserted normal forms. The minor choices
in the lemma also give the required rational pivot sections.
\par\hfill$\square$\par\smallskip

\paragraph{Intersection groups and their Levi images}\label{m8:s:m:intersection}

\subparagraph{Root conjugations}\label{m8:s:m:roots}

Write \(X_{ij},Y,Z_{ij},Z_0\) for the root matrices obtained by setting, respectively, \(X=E_{ij}\), $y = 1$, \(Z=E_{ij}\), $z = 1$ in
\[
N(X,y,Z,z)=\begin{pmatrix}
0&X&Z&zI\\0&0&yI&Z^\#\\0&0&0&-X^\#\\0&0&0&0
\end{pmatrix}.
\]
We denote the scalar $z$-root by \(Z_0\) and the four Levi
roots by \(A_+,A_-,B_+,B_-\). 

If $w_ie_j=c_je_{\rho_i(j)}$ is the signed permutation above, then
\[
 (w_i^{-1}Rw_i)_{jk}=c_jc_kR_{\rho_i(j),\rho_i(k)}.
\]
A root lies in $\operatorname{Lie}I_i$ exactly when every nonzero
entry of this matrix satisfies
$\operatorname{block}(j)\ge\operatorname{block}(k)$.
The roots with zero image in both Levi projections give the kernels
in \S\ref{m8:s:m:log-stab}. The four mixed roots in each row are:
\[
\begin{array}{c|l|l}
i&U\text{-root to right Levi}&\text{left Levi to }U^t\\\hline
10&Y\mapsto A_+,\ Z_0\mapsto B_+&A_-\mapsto Y^t,\ B_-\mapsto Z_0^t\\
11&Z_{12}\mapsto-A_-,\ Z_{21}\mapsto B_+&A_+\mapsto-Z_{12}^t,\ B_-\mapsto Z_{21}^t\\
12&X_{11}\mapsto-A_-,\ X_{22}\mapsto B_+&A_+\mapsto-Z_{12}^t,\ B_+\mapsto Z_{21}^t\\
15&Z_{12}\mapsto-A_-,\ Z_{21}\mapsto-B_-&A_+\mapsto X_{11}^t,\ B_-\mapsto X_{22}^t\\
16&X_{11}\mapsto-A_-,\ X_{22}\mapsto-B_-&A_+\mapsto X_{11}^t,\ B_+\mapsto X_{22}^t.
\end{array}
\]
For example, row 10 sends the $y$-root to $A_+$ and the left
$A_-$-root to $Y^t$; their exponential subgroups realize the two
independent Levi directions. All four Levi roots map into the opposite unipotent radical.
The complete fourteen-root conjugation lists are in
\texttt{anc/model8\_calculations.pdf}, M8-C04, ``Root conjugations'';
they implement the displayed membership criterion.

\subparagraph{The full algebraic intersection groups}\label{m8:s:m:components}

There is a direct matrix argument. An element of $I$ preserves the standard block flag and the monomially transported opposite block flag. If \(c_i(j)\) is the column of \(w_i\) whose nonzero entry is in row j, a matrix entry in position (j,$k$) can be nonzero only if
\[
\operatorname{block}(j)\le\operatorname{block}(k),\qquad
\operatorname{block}(c_i(j))\ge\operatorname{block}(c_i(k)).
\]
The joint multiplicities in Section~\ref{m8:s:m:flags} are at most one. Therefore the equal-grade blocks have dimension at most one, and
the conditions extend to a full triangular order. Explicit self-dual orders of the standard basis are:

\begingroup
\small
\setlength{\tabcolsep}{3pt}
\renewcommand{\arraystretch}{1.18}
\begin{longtable}{p{0.16\linewidth}p{0.75\linewidth}}
\toprule
i & Basis order in which $I$ is upper triangular \\
\midrule
\endfirsthead
i & Basis order in which $I$ is upper triangular \\
\midrule
\endhead
10 & 2,4,1,6,3,8,5,7 \\
11 & 1,2,4,6,3,5,7,8 \\
12 & 1,2,3,4,5,6,7,8 \\
15 & 1,2,6,4,5,3,7,8 \\
16 & 1,2,3,4,5,6,7,8 \\
\bottomrule
\end{longtable}
\endgroup

Each order is self-dual under the $Q$-pairing of coordinates j and 9-j. Thus $I$ lies in the corresponding split Borel of $\mathrm{GSO}_8$.

Projecting $I$ onto its diagonal gives the full torus
\[
T_M=\{\operatorname{diag}(d_1,\ldots,d_8):
 d_1d_8=d_2d_7=d_3d_6=d_4d_5=\mu\},
\]
with five independent nonzero parameters. The diagonal matrices give a section, and the kernel is
unitriangular. Because polynomial logarithm and exponential
preserve both flags, they identify this kernel with the root-space
span in Section~\ref{m8:s:m:roots}. Thus the kernel is an affine
connected unipotent group over $F$ and after scalar extension.

It follows that $I$ is a split connected solvable algebraic group.
Every $F$-point is captured by its torus part and its unipotent root
coordinates. Connectedness here is understood algebraically.

\subparagraph{Levi images and dimensions}\label{m8:s:m:levi-images}

Set \(t_i=w_i^{-1}tw_i\). The image of
\[
p\longmapsto
\bigl(\pi_M(p),\pi_M(w_i^{-1}pw_i)\bigr)
\]
is exactly
\[
{\{(lt,t_ir):l\in L_i,\ t\in T_M,\ r\in R_i\}.}
\]
The groups are

\begingroup
\small
\setlength{\tabcolsep}{3pt}
\renewcommand{\arraystretch}{1.18}
\begin{longtable}{p{0.055\linewidth}p{0.245\linewidth}p{0.245\linewidth}p{0.17\linewidth}p{0.195\linewidth}}
\toprule
i & \(L_i\), in the A,$B$ factors & \(R_i\), in the A,$B$ factors & Kernel dimension & \(\dim I_i\) \\
\midrule
\endfirsthead
i & \(L_i\), in the A,$B$ factors & \(R_i\), in the A,$B$ factors & Kernel dimension & \(\dim I_i\) \\
\midrule
\endhead
10 & \(N_-\times N_-\) & \(N_+\times N_+\) & 4 & 13 \\
11 & \(N_+\times N_-\) & \(N_-\times N_+\) & 5 & 14 \\
12 & \(N_+\times N_+\) & \(N_-\times N_+\) & 6 & 15 \\
15 & \(N_+\times N_-\) & \(N_-\times N_-\) & 6 & 15 \\
16 & \(N_+\times N_+\) & \(N_-\times N_-\) & 7 & 16 \\
\bottomrule
\end{longtable}
\endgroup

The left and right mixed roots project to \((l,1)\) and
\((1,r)\), respectively, while the torus projects to
\((t,t_i)\). Taking their products realizes arbitrary $l,t,r$
over $F$; any commutator corrections lie in the kernel. Hence
the projection is surjective on $F$-points.

The kernel is exactly
\[
U\cap w_iU^tw_i^{-1}.
\]
The dimensions give a separate consistency check:
\[
\dim S_i=2\dim P-\dim I_i
 =38-(13,14,15,15,16)
 =(25,24,23,23,22).
\]

\paragraph{Residual equivalence and transition characters}\label{m8:s:m:transitions}

For two representatives \(mw_in\) and \(m'w_in'\), equality of even $K$-orbits is equivalent to
\begin{equation}
{
 m'=D(a)mlt,\qquad
 n'=r^{-1}t_i^{-1}nD(a)^{-1}
}
\label{m8:s:m:eq:transition}
\end{equation}
for some \(a\in\mathrm{GL}_2(F),l\in L_i,r\in R_i,t\in T_M\).

Necessity. An even transporter gives a relation between two $P$--$P^t$ presentations. The element relating those presentations lies in \(I_i\); applying its two Levi projections and using Section~\ref{m8:s:m:levi-images} gives \eqref{m8:s:m:eq:transition}.

Sufficiency, including the character. Choose \(p\in I_i(F)\) with projections \((lt,t_ir)\), and write
\[
p=u_p(lt),\qquad
w_i^{-1}pw_i=v_p^t(t_ir).
\]
Then
\[
(lt)w_i(t_ir)^{-1}=u_p^{-1}w_iv_p^t.
\]
Thus, before the final diagonal conjugation, a transporter has
\[
u=m u_p^{-1}m^{-1},\qquad
v=n^t v_p n^{-t}.
\]
Its character is
\begin{equation}
{\xi(m u_p^{-1}m^{-1})\,
       \xi(n^t v_p n^{-t}).}
\label{m8:s:m:eq:transition-character}
\end{equation}
Diagonal conjugation has character one and preserves this value. This gives both the parameter equivalence and the character of its transporter. 

If two lifts $p$ have the same Levi images, the resulting transporters differ by a full $U\times U$ stabilizer element. On the compatible locus its character is one, so \eqref{m8:s:m:eq:transition-character} is independent of the lift.

\subparagraph{Minor charts}\label{m8:s:m:charts}

The two charts for a right \(N_-\)-quotient, with fixed second column \(v=(v_1,v_2)^t\) and determinant \(\Delta\ne0\), are
\[
S_-^{(1)}=\begin{pmatrix}\Delta/v_2&v_1\\0&v_2\end{pmatrix},\quad v_2\ne0;
\qquad
S_-^{(2)}=\begin{pmatrix}0&v_1\\-\Delta/v_1&v_2\end{pmatrix},\quad v_1\ne0.
\]
For a right \(N_+\)-quotient, with fixed first column, they are
\[
S_+^{(1)}=\begin{pmatrix}v_1&0\\v_2&\Delta/v_1\end{pmatrix},\quad v_1\ne0;
\qquad
S_+^{(2)}=\begin{pmatrix}v_1&-\Delta/v_2\\v_2&0\end{pmatrix},\quad v_2\ne0.
\]
Equality of the fixed column and determinant forces the other columns to differ by the appropriate root operation. On the overlap,
\[
S_\pm^{(2)}=S_\pm^{(1)}
 n_\pm\!\left(-\frac{\Delta}{v_1v_2}\right).
\]
Transposition gives the left-quotient charts, interchanging the labels + and - and preserving the root parameter; the inverse in \eqref{m8:s:m:eq:transition} changes its sign.

The product of these covers gives sixteen normal-form parameter charts. Their denominators are the indicated nonzero coordinates, determinants, pairing pivots and similitudes.

\paragraph{The two-sided unipotent stabilizers}\label{m8:s:m:stab}

\subparagraph{Group law and exponential coordinates}\label{m8:s:m:group-law}

Use the complete group law and inverse in
\eqref{m8:s:t:eq1.3}--\eqref{m8:s:t:eq1.4}. Direct finite
exponentiation gives
\begin{equation}
{\exp N(X,y,Z,z)
 =u\left(X,y,Z+\tfrac y2X,z-\tfrac y6\det X\right).}
\label{m8:s:m:eq:exponential}
\end{equation}

\subparagraph{Quadratic corrections in group coordinates}\label{m8:s:m:corrections}

To read these as $U$-coordinates, keep the displayed $X$ and $y$, replace $Z$ by \(Z+yX/2\), and retain the displayed $z$. The latter is legitimate because \(y\det X=0\) and \(y'\det X'=0\) in every row.

The nonzero $Z$-corrections are
\[
\begin{array}{c|c}
\text{stratum; side}&\text{group }Z\text{-coordinate}\\\hline
S_{12},\ \text{left}&\left(\begin{smallmatrix}a&b+qy/2\\0&d\end{smallmatrix}\right)\\
S_{15},\ \text{right}&\left(\begin{smallmatrix}-q&z-ha/2\\0&r\end{smallmatrix}\right)\\
S_{16},\ \text{left}&\left(\begin{smallmatrix}a&b+qy/2\\c&d\end{smallmatrix}\right)\\
S_{16},\ \text{right}&\left(\begin{smallmatrix}-a&z-cq/2\\-y&d\end{smallmatrix}\right).
\end{array}
\]
All other coordinates are unchanged. The tables have
$N^3=(N')^3=0$ and
$(I+N+N^2/2)w_i(I+N'+(N')^2/2)^t=w_i$.

\subparagraph{Completeness over the ground field}\label{m8:s:m:full-stab}

If a pair of $U$-elements stabilizes \(mw_in\), conjugate by $m$ and \(n^t\) and take their finite polynomial logarithms. The stabilizer identity forces exactly the relation in \eqref{m8:s:m:eq:stabilizer}. Conversely, exponentiating any solution gives a stabilizer element. These operations are mutual inverses on every $F$-point.

We obtain affine connected unipotent stabilizers of dimensions
4,5,6,6,7, respectively, with polynomial exponential and logarithm
covering each entire $F$-point group.

\paragraph{\texorpdfstring{Canonical forms on $S_{10}$ and $S_{11}$}{Canonical forms on S10 and S11}}\label{m8:s:m:normalization}

Let
\[
W=\begin{pmatrix}0&1\\1&0\end{pmatrix},\qquad
H=\operatorname{diag}(1,-1),\qquad
m_0=m(W,I,1),\qquad V=m_0w_{11}m_0.
\]
Direct multiplication gives
\[
V^t=V,\qquad V^2=I_8,\qquad V\in\mathrm{SO}(Q)(F).
\]
Define
\[
F_\epsilon(E)=w_{10}m(\operatorname{diag}(1,\epsilon)E,E,
                            \epsilon\det E),
\qquad \epsilon\in F^\times,
\]
\[
G(E)=m_0w_{11}m(WE,E,\det E)=VD(E).
\]
The family \(F_\epsilon(E)\) has the form \(V_\epsilon D(E)\), where \(V_\epsilon=F_\epsilon(I)\) is symmetric. The similitudes are \(\epsilon\det E\) and \(\det E\), respectively. All these points belong to their indicated rank strata, since multiplication by the displayed Levi factors preserves the corner ranks.

\subparagraph{The two rational normalizations}\label{m8:s:m:normalization10}
\label{m8:s:m:normalization11}
On $S_{10}$ write
$R=(\begin{smallmatrix}\alpha&0\\\beta&\gamma\end{smallmatrix})$,
$S=(\begin{smallmatrix}\alpha&0\\\delta&\zeta\end{smallmatrix})$,
with $\alpha\gamma\zeta\ne0$. The following allowed transitions
reach the canonical family:
\begin{enumerate}
\item Right-multiply the left $A$-factor by
$n_-(-\beta/\gamma)$, then use
$t=m(\operatorname{diag}(\alpha^{-1},\gamma^{-1}),I,\kappa)$
and diagonal conjugation by $B^{-1}$. The left factor becomes $1$.
The first root move has character
$\psi(-2\lambda\beta/\gamma)$; the torus reparameterization
and diagonal conjugation have character one.
\item Compatibility now makes
$CD^{-1}=(\begin{smallmatrix}1&z_1\\0&\epsilon\end{smallmatrix})$,
$\epsilon\ne0$. Left-multiply $C$ by $n_+(-z_1/\epsilon)$,
with character $\psi(-2z_1/\epsilon)$. This gives
\begin{equation}\label{m8:s:m:eq:pregauge10}
 w_{10}m(\operatorname{diag}(1,\epsilon)E,E,\rho\det E),
 \qquad \rho\ne0.
\end{equation}
\item The pure diagonal move $a=\operatorname{diag}(\rho/\epsilon,1)$
sends $E$ to $E\operatorname{diag}(1,\rho/\epsilon)$ and
$\rho$ to $\epsilon$, with character one. It gives $F_\epsilon(E)$.
\end{enumerate}
The roots and characters follow from
\eqref{m8:s:m:eq:transition}--\eqref{m8:s:m:eq:transition-character}.
Every divisor used is among the nonzero
$\alpha,\gamma,\epsilon,\rho,\kappa$ or $\det B$;
the normalization characters are as computed above.

On $S_{11}$ write
$R=(\begin{smallmatrix}0&b\\c&d_0\end{smallmatrix})$,
$S=(\begin{smallmatrix}0&b\\c&e_0\end{smallmatrix})$, $bc\ne0$.
Use successively the left $A$-root $n_+(-d_0/c)$,
$t=m(\operatorname{diag}(c^{-1},b^{-1}),I,\kappa)$,
and diagonal conjugation by $B^{-1}$. The left factor becomes $m_0$.
The right ratio is now
$(\begin{smallmatrix}0&1\\1&e_1\end{smallmatrix})$;
the right $A$-root $n_-(-e_1)$ gives
\begin{equation}\label{m8:s:m:eq:pregauge11}
 m_0w_{11}m(WE,E,\rho\det E).
\end{equation}
Finally $a=\operatorname{diag}(\rho,1)$ sends
$E$ to $E\operatorname{diag}(1,\rho)$ and $\rho$ to $1$,
giving $G(E)$. Each move has character one because these mixed roots lift through
$Z$-coordinates, whose Levi conjugates have zero $X$-trace and
$y$-coordinate. Since we divide only by the forced nonzero
$b,c,\kappa,\rho$ and invertible Levi matrices, both constructions
are defined over $F$.
Their expanded steps are in \texttt{anc/model8\_calculations.pdf},
calculation M8-C04, under ``Canonical forms''.

\paragraph{Complete residual equivalence}\label{m8:s:m:residual}
\label{m8:s:m:left-equations}\label{m8:s:m:right-equations}
Apply the full projection equations
\eqref{m8:s:m:eq:transition} to the pre-gauge forms.
Keeping the left canonical $B$-factor forces the diagonal component
$a=(\begin{smallmatrix}\alpha&0\\\gamma&\delta\end{smallmatrix})$,
$\alpha\delta\ne0$. Put $t=\alpha/\delta$, $h=\gamma/\delta$.
The remaining left projections are uniquely
\[
\begin{array}{c|c|c|c}
i&l_A,l_B&t_A&t_B\\\hline
10&n_-(-h),n_-(-h)&\operatorname{diag}(\alpha^{-1},\delta^{-1})&
\operatorname{diag}(\alpha^{-1},\delta^{-1})\\
11&n_+(-h),n_-(-h)&\operatorname{diag}(\delta^{-1},\alpha^{-1})&
\operatorname{diag}(\alpha^{-1},\delta^{-1}).
\end{array}
\]
In both cases $\mu_t=(\alpha\delta)^{-1}$ and
$w_i^{-1}tw_i=\operatorname{diag}(\alpha^{-1}I_2,
\alpha^{-1}I_2,\delta^{-1}I_2,\delta^{-1}I_2)$.
For the two right-root parameters $u,v$, the ratios become
\[
\begin{array}{c|c|c}
i&\text{new ratio}&\text{forced equation}\\\hline
10&n_+(-u)\operatorname{diag}(1,\epsilon)n_+(v)&v=\epsilon u, \epsilon'=\epsilon\\
11&n_-(-u)Wn_+(v)&v=u.
\end{array}
\]
These equations hold independently of $E$, including its singular
residual loci. They force the complete pre-gauge relation
\begin{equation}\label{m8:s:m:eq:pre-residual}
 E'=n_+(s)E\begin{pmatrix}1&0\\-h&t\end{pmatrix},\qquad
 \rho'=\rho/t.
\end{equation}
The scalar relation also follows from preservation of similitude
and $\det E'=t\det E$. Using an intersection lift
$p=l\,t_M\,w_i r w_i^{-1}$ in
\eqref{m8:s:m:eq:transition-character} gives the entire characters
\begin{equation}\label{m8:s:m:pre-character}
 \chi_{10}^{\rm pre}=\psi(2s/\epsilon-2h/\rho),\qquad
 \chi_{11}^{\rm pre}=1.
\end{equation}
Their left/right exponents are exactly $2s/\epsilon$ and $-2h/\rho$. The full block substitution
is preserved in \texttt{anc/model8\_calculations.pdf}, calculation
M8-C04, under ``Characters before the similitude gauge''.

\subparagraph{After the gauge}\label{m8:s:m:gauge}
Fixing $\rho=\epsilon$, respectively $\rho=1$, forces $t=1$.
Putting $v=-h$ gives precisely
\eqref{m8:s:m:eq:residual}--\eqref{m8:s:m:eq:residual-character}.
Necessity follows from all the projection equations just solved;
the following whole-group lifts prove sufficiency with the exact
signs. These equations characterize equivalence between gauged points,
including through intermediate ungauged points.

\subparagraph{\texorpdfstring{Nonlinear lifts to $K$}{Nonlinear lifts to K}}\label{m8:s:m:lifts}

Set
\[
E'=n_+(s)En_-(v),\qquad a_v=n_-(-v).
\]
The $S_{10}$ lift is
\[
{
k_{10}(s,v)=
\bigl(u(0,s/\epsilon,0,s),\,
      u(0,v/\epsilon,0,v),\,a_v\bigr).
}
\]
It satisfies
\[
k_{10}(s,v)F_\epsilon(E)=F_\epsilon(E'),\qquad
\eta_K(k_{10}(s,v))=\psi(2(s+v)/\epsilon).
\]
The $S_{11}$ lift is
\[
{
k_{11}(s,v;E)=
\bigl(u(0,0,s a_vHa_v^{-1},0),\,
      u(0,0,vE'^tHE'^{-t},0),\,a_v\bigr).
}
\]
Here
\[
s a_vHa_v^{-1}
 =s\begin{pmatrix}1&0\\-2v&-1\end{pmatrix}.
\]
The second $Z$-coordinate uses the target $E'$ and its inverse transpose. All these inverses exist because \(\det E'=\det E\ne0\). Multiplication of the full $U$ matrices gives
\[
k_{11}(s,v;E)G(E)=G(E'),\qquad \eta_K(k_{11}(s,v;E))=1.
\]
Composites differ from the displayed composite lift by a $U\times U$ stabilizer, whose character is one on the compatible locus. Hence these lifts realize the stated residual action and character even though the $S_{11}$ lift depends on $E$.

\paragraph{Rank topology and the paired even-vanishing family}\label{m8:s:m:topology}

\subparagraph{Closure incidences}\label{m8:s:m:closures}

Write
\[
c=r_{2,2},\quad a=r_{4,4},\quad
b=r_{6,2}=r_{2,6},\quad d=r_{4,2},\quad e=r_{2,4}.
\]
In the rank notation,
\[
R=\begin{pmatrix}4+c&2+e&b\\2+d&a&d\\b&e&c\end{pmatrix}.
\]
For the five labels the parameter tuples are
\[
\begin{array}{c|c}
10&(0,2,2,1,1)\\
11&(0,2,1,1,1)\\
12&(0,2,1,0,1)\\
15&(0,2,1,1,0)\\
16&(0,0,1,0,0).
\end{array}
\]
The induced diagram is
\[
S_{10}\longrightarrow S_{11}\longrightarrow
\begin{matrix}S_{12}\\S_{15}\end{matrix}
\longrightarrow S_{16}.
\]
The arrow indicates containment of the lower stratum in the closure of the upper. It is the induced closure diagram on these five labels.

\subparagraph{One-parameter degenerations}\label{m8:s:m:degenerations}

Using the root notation of Section~\ref{m8:s:m:intersection}, the families are

\begingroup
\small
\setlength{\tabcolsep}{3pt}
\renewcommand{\arraystretch}{1.18}
\begin{longtable}{p{0.16\linewidth}p{0.36\linewidth}p{0.39\linewidth}}
\toprule
Family & Stratum for $t \ne 0$ & Limit at $t = 0$ \\
\midrule
\endfirsthead
Family & Stratum for $t \ne 0$ & Limit at $t = 0$ \\
\midrule
\endhead
\((I+tX_{12}^t)w_{11}\) & $S_{10}$ & \(w_{11}\) \\
\((I+tY^t)w_{12}\) & $S_{11}$ & \(w_{12}\) \\
\((I+tX_{22}^t)w_{15}\) & $S_{11}$ & \(w_{15}\) \\
\((I+tZ_{21}^t)w_{16}\) & $S_{12}$ & \(w_{16}\) \\
\((I+tY^t)w_{16}\) & $S_{15}$ & \(w_{16}\) \\
\bottomrule
\end{longtable}
\endgroup

The root matrices square to zero, and each family satisfies \(x(t)^tQx(t)=Q\) and \(\det x(t)=1\) identically. The corner-rank identities hold for every $t\ne0$. The complete nine-corner pivot-minor lists are in
\texttt{anc/model8\_calculations.pdf}, M8-C04,
``One-parameter degenerations''. Every listed maximal minor is
$\pm1$, $\pm t$ or $\pm t^2$.
The symbolic upper rank bounds are polynomial identities; these monomial minors achieve those ranks at every nonzero $t$. The rational double-cell result propagates the incidences from the representatives to every $F$-point of the lower stratum, by fixed left and right $P$ factors.

Upper semicontinuity excludes the reverse incidences. It also excludes either of $S_{12}$,$S_{15}$ lying in the closure of the other: the relevant zero block cannot become a nonzero rank-one block in a limit. 

\subparagraph{A locally closed union}\label{m8:s:m:familyB}

The rank table gives
\[
\mathcal B=S_{12}\sqcup S_{15}\sqcup S_{16}
 =\{c=0,b=1,\ x_{34}=0\text{ or }x_{43}=0\}.
\]
The locus $c = 0, b = 1$ is locally closed; the stated disjunction is a union of two closed block-vanishing conditions inside it. Therefore $B$ is locally closed. Transpose interchanges these two conditions and preserves $B$.

Inside $B$, $S_{16}$ is cut out by both block vanishings. It is closed. Its complement is $S_{12} \sqcup S_{15}$; these are disjoint closed subsets of that complement and hence are each open-and-closed there.

The closed-support sequence
\[
0\longrightarrow D'_{\eta_K}(S_{16})
\longrightarrow D'_{\eta_K}(\mathcal B)
\longrightarrow D'_{\eta_K}(S_{12}\sqcup S_{15})
\]
has zero ends by Proposition~\ref{m8:s:m:prop:even}. Hence
$D'_{\eta_K}(\mathcal B)=D'_{\widetilde\eta_K}(\mathcal B)=0$;
the same proof gives $D'_{\xi\boxtimes\xi}(\mathcal B)=0$.

\subparagraph{A union which is not locally closed}\label{m8:s:m:skipped}

Let
\[
\mathcal A=S_{10}\sqcup S_{12}\sqcup S_{15}\sqcup S_{16}.
\]
The displayed degenerations imply \(S_{11}\subset\overline{S_{10}}\subset\overline{\mathcal A}\), and every point of $S_{12}$ is in \(\overline{S_{11}}\). If A were locally closed, it would be open in its closure. At a point of $S_{12}$ this would give an ambient neighborhood whose intersection with \(\overline A\) is contained in A. But every such neighborhood meets $S_{11}$, which is contained in \(\overline A\setminus A\). This is a contradiction.

The argument uses $F$-topological degenerations as well as algebraic closure incidences. Consequently this union cannot replace the separate locally closed spaces in an ordinary distribution argument.

\subparagraph{Transpose on the paired strata}\label{m8:s:m:transpose-pair}

With
\[
L=m(H,I,1)=\operatorname{diag}(1,-1,1,1,1,1,-1,1),
\]
direct matrix multiplication gives
\[
{w_{12}^t=Lw_{15},\qquad w_{15}^t=w_{12}L,}
\]
\[
{w_{16}^t=m(H,-H,1)w_{16}.}
\]
Consequently the redundant parameter maps are
\[
(m,n)_{12}\longmapsto(n^tL,m^t)_{15},\qquad
(m,n)_{15}\longmapsto(n^t,Lm^t)_{12}.
\]
Their composition is the identity because \(L^t=L\) and \(L^2=I\). Even $U\times U$-equivariance already gives their vanishing.

\subparagraph{Rank envelopes and predecessor degenerations}\label{m8:s:m:predecessor-degenerations}
For a rank matrix $R_i$, let
\[
\Omega_i=\{x\in G_8^\circ:R(x)\ge R_i\text{ entrywise}\}.
\]
Writing $\mathcal O_{\rm cleared}=S_1\sqcup S_2\sqcup S_3\sqcup S_4\sqcup S_5\sqcup S_8$, the rank table gives
\begin{align*}
\Omega_{10}\setminus\mathcal O_{\rm cleared}&=S_6\sqcup S_{10},\\
\Omega_{11}\setminus\mathcal O_{\rm cleared}&=S_6\sqcup S_7\sqcup S_{10}\sqcup S_{11},\\
\Omega_{12}\setminus\mathcal O_{\rm cleared}&=S_6\sqcup S_7\sqcup S_9\sqcup S_{10}\sqcup S_{11}\sqcup S_{12},\\
\Omega_{15}\setminus\mathcal O_{\rm cleared}&=S_6\sqcup S_7\sqcup S_{10}\sqcup S_{11}\sqcup S_{13}\sqcup S_{15},\\
\Omega_{16}\setminus\mathcal O_{\rm cleared}&=S_6\sqcup S_7\sqcup S_9\sqcup S_{10}\sqcup S_{11}\sqcup S_{12}\sqcup S_{13}\sqcup S_{14}\sqcup S_{15}\sqcup S_{16}.
\end{align*}
Each $S_i$ is closed in $\Omega_i$: within the lower-rank-bound open, impose the corresponding upper rank bounds. Transpose fixes $\Omega_{10},\Omega_{11},\Omega_{16}$ and interchanges $\Omega_{12},\Omega_{15}$.

There are also the explicit predecessor degenerations
\[
(I+tX_{11}^t)w_{10}\in S_6\quad(t\ne0),\qquad
(I+tX_{21}^t)w_{11}\in S_7\quad(t\ne0).
\]
Both products preserve $Q$ and have determinant one. In the row-major order of the nine corners, nonzero maximal pivot-minor witnesses are respectively
\[
(t,1,1,1,-1,1,1,1,-t),\qquad
(t,-1,-1,-1,-1,-1,-1,-1,-t).
\]
The first family and the rational parabolic-cell description put every point of $S_{10}$ in the closure of $S_6$. The degenerations of Section~\ref{m8:s:m:degenerations} then give the same containment for all five middle strata. The second family places $S_{11}$ directly in the closure of $S_7$. Consequently the interior of
\[
\mathcal O_{\rm cleared}\sqcup S_{10}\sqcup S_{12}\sqcup S_{15}\sqcup S_{16}
\]
in $G_8^\circ$ is precisely $\mathcal O_{\rm cleared}$. The support argument therefore uses the predecessor opens, including the prior removal of $S_6$ and $S_7$, as in Section~\ref{m8:s:sec:assembly}.

\subsection{Rank-one incidence: rational coverage, fibers and densities}\label{m8:supp:s6s7}
Use the fixed $U$ and $Q$ of Model 8 and write $L(A,B,t)=m(A,B,t)$.
For the normal form $lw_ir$ put $l=L(A,B,t)$, $r=L(C,E,s)$,
$M=B^{-1}A$, $N=CE^{-1}$, $\kappa=\det B/t$, and
$\ell=\det E/s$. The representatives and compatible loci are those
of the main rank-one incidence argument; their intersection calculations
are given in \S\ref{m8:supp:s6s7-details}.

\paragraph{Rational parabolic coverage}\label{m8:s:v:s3}

\begin{lemma}\label{m8:s:v:rational-coverage}
For \(i=6,7\), one has \(S_i(F)=P(F)w_iP(F)^t\). The map from the two parabolic factors to \(S_i(F)\) admits rational local sections on pivot charts.
\end{lemma}
\begin{proof}\leavevmode

\subparagraph{Intersection dimensions and all graded pieces}\label{m8:s:v:s3-1}

Let \(F_i\) be the span of the first \(2i\) standard basis vectors, let \(F_j^-\) be the span of the last \(2j\), and put
\(G_j=xF_j^-\) \((0\leq i,j\leq4).\)
These are self-dual flags: \(F_i^\perp=F_{4-i}\) and \(G_j^\perp=G_{4-j}\). The first half is isotropic; the second half is coisotropic.

Projection to \(F/F_i\) proves, over \(F\),
\(\dim(F_i\cap G_j)=2j-r_{8-2i,2j}(x).\)
Set \(d_{ij}=\dim(F_i\cap G_j)\), with the evident zero boundary values, and let
\(m_{ij}=d_{ij}-d_{i-1,j}-d_{i,j-1}+d_{i-1,j-1}.\)
For the two prescribed rank patterns, the complete matrices of two-flag graded dimensions are
\begin{equation}
m^{(6)}=\begin{pmatrix}
0&1&0&1\\1&0&1&0\\0&1&0&1\\1&0&1&0
\end{pmatrix},
\qquad
m^{(7)}=\begin{pmatrix}
1&0&0&1\\0&1&1&0\\0&1&1&0\\1&0&0&1
\end{pmatrix}.
\label{m8:s:v:e3-2}
\end{equation}
Every nonzero graded piece is one-dimensional. The complementary piece has degree \((5-i,5-j)\). 

\subparagraph{Simultaneous Witt bases and pivot sections}\label{m8:s:v:s3-2}
Choose a basis $V=(v_1,\ldots,v_4)$ of $F_2$ adapted to
$F_1$ and $F_2\cap G_j$, with grades
\begin{equation}\label{m8:s:v:e3-3}
\begin{array}{c|cccc}
&v_1&v_2&v_3&v_4\\\hline
S_6&(1,4)&(1,2)&(2,1)&(2,3)\\
S_7&(1,4)&(1,1)&(2,3)&(2,2).
\end{array}
\end{equation}
For $v_a$ of degree $(i_a,j_a)$, Lemma~\ref{m8:simultaneous-witt}
solves for its dual column in
$E_a=F_{5-i_a}\cap G_{5-j_a}$, using
\begin{equation}\label{m8:s:v:e3-4}
F_2\cap E_a^\perp=F_{i_a-1}+(F_2\cap G_{j_a-1}).
\end{equation}
Writing the dual columns as $W$, it gives
\begin{equation}\label{m8:s:v:e3-5}
V^tQW=I_4.
\end{equation}
\label{m8:s:v:s3-3}
The flag-preserving Gram correction is
\begin{equation}\label{m8:s:v:e3-6}
W'=W-\tfrac12V(W^tQW),\qquad V^tQW'=I_4,\quad W'^tQW'=0.
\end{equation}
\label{m8:s:v:s3-4}
The reference basis
$B_{\rm ref}=(e_1,e_2,e_3,e_4,-e_8,e_7,-e_6,e_5)$ has the
same grades and Gram matrix. Thus
\begin{equation}\label{m8:s:v:e3-8}
p=B_{\rm ref}(V,W')^{-1}\in P(F),\qquad w_i^{-1}px\in P(F)^t.
\end{equation}
Preservation of $F_2$ gives determinant one, proving coverage in
the selected connected group and rational local sections on the
pairing-minor charts. The construction retains the original similitude, including for
$\mu(w_6)=-1$. Factoring both parabolics gives
$x=u_Llw_ir u_R^t$; the character of the outside unipotents is
handled by \eqref{m8:s:c:transport-back}.
\end{proof}

\paragraph{intersection groups and normal-family equivalence}\label{m8:s:v:s4}

\begin{lemma}\label{m8:s:v:intersections}
The groups \(J_i=P\cap w_iP^tw_i\) are given by \eqref{m8:s:v:e4-2} and \eqref{m8:s:v:e4-5}, with second Levi projections \eqref{m8:s:v:e4-3} and \eqref{m8:s:v:e4-6}. Their dimensions are \(12\) and \(13\). Two normal-family points are \(K\)-equivalent if and only if \eqref{m8:s:c:normal-equations} holds.
\end{lemma}
\begin{proof}\leavevmode

\subparagraph{Solving from a general parabolic element}\label{m8:s:v:s4-1}

Start with \(j=L(A_0,B_0,\omega)u(X,y,Z,z_0)\), with completely general invertible two-by-two \(A_0,B_0\), general \(X,Z\), and \(\omega\ne0\). Impose the upper-block zero equations on \(w_i j w_i\).

For \(S_6\) the complete system, after localizing at the Levi determinants and \(\omega\), is equivalent to
\begin{equation}
(A_0)_{12}=(B_0)_{21}=X_{11}=X_{12}=X_{21}=Z_{11}=z_0=0.
\label{m8:s:v:e4-1}
\end{equation}
For example, the ordinary entries \((1,6)\) and \((4,7)\) of \(w_6jw_6\) first give \((A_0)_{12}=(B_0)_{21}=0\). Invertibility then makes all four triangular diagonal entries nonzero. Entries \((1,7)\), \((1,4)\), \((2,3)\), \((1,5)\), \((1,3)\), in this order, force the remaining five equations in \eqref{m8:s:v:e4-1}. Substitution satisfies every upper-block zero equation.

Thus
\(J_6=P\cap w_6P^tw_6\)
consists exactly of
\begin{equation}
j=L\left(
\begin{pmatrix}\alpha&0\\\gamma&\delta\end{pmatrix},
\begin{pmatrix}\epsilon&\phi\\0&\zeta\end{pmatrix},\omega
\right)
 u\left(pE_{22},y,\begin{pmatrix}0&q_0\\r_0&z\end{pmatrix},0\right),
\label{m8:s:v:e4-2}
\end{equation}
where \(\alpha\delta\epsilon\zeta\omega\ne0\). Its second Levi projection is
\begin{equation}
\pi_L^-(w_6jw_6)=L\left(
\begin{pmatrix}\alpha&\alpha q_0\\0&\omega/\epsilon\end{pmatrix},
\begin{pmatrix}\omega/\delta&0\\-r_0\zeta&\zeta\end{pmatrix},\omega
\right).
\label{m8:s:v:e4-3}
\end{equation}
There are seven Levi and five unipotent parameters: \(\dim J_6=12\).

For \(S_7\) the complete system is
\((A_0)_{12}=(B_0)_{12}=X_{11}=X_{12}=Z_{11}=Z_{12}=0.\)
Ordinary entries \((1,7),(3,5)\) first give the two Levi equations. Then \((1,3),(1,5),(1,4),(1,6)\) give the four unipotent equations. All remaining upper-block entries vanish identically after substitution. In particular, \emph{both} \(y\) and the scalar \(z_0\) remain free in \(J_7\).

The complete answer is
\begin{equation}
j=L\left(
\begin{pmatrix}\alpha&0\\\gamma&\delta\end{pmatrix},
\begin{pmatrix}\epsilon&0\\\phi&\zeta\end{pmatrix},\omega
\right)
 u\left(\begin{pmatrix}0&0\\p&q\end{pmatrix},y,
          \begin{pmatrix}0&0\\r&z\end{pmatrix},z_0\right),
\label{m8:s:v:e4-5}
\end{equation}
with second projection
\begin{equation}
\pi_L^-(w_7jw_7)=L\left(
\begin{pmatrix}\alpha&\alpha z_0\\0&\omega/\delta\end{pmatrix},
\begin{pmatrix}\epsilon&\epsilon y\\0&\omega/\zeta\end{pmatrix},\omega
\right).
\label{m8:s:v:e4-6}
\end{equation}
There are seven Levi and six unipotent parameters: \(\dim J_7=13\). Subsection~\ref{m8:s:v:sA} gives the full matrices \(w_i j w_i\), including the lower blocks and all nonlinear terms.

The orbit dimensions are consequently
\(2\dim P-\dim J_6=38-12=26,\) \(2\dim P-\dim J_7=38-13=25.\)

\subparagraph{Necessity and sufficiency of the residual equivalence}\label{m8:s:v:s4-2}

Apply Lemma~\ref{m8:s:c:normal-family} with $w=w_i=w_i^{-1}$ and the two computed projections. Thus \eqref{m8:s:c:normal-equations} is the necessary and sufficient residual relation and \eqref{m8:s:c:normal-transporter} gives its unipotents, including the inverse on the right.

Writing \(\pi_L(j)=L(A_0,B_0,\omega)\) and \(\pi_L^-(w_ijw_i)=L(C_0,E_0,\omega)\), the component equations are
\begin{equation}
\begin{aligned}
A'&=hAA_0^{-1},& B'&=hBB_0^{-1},& t'&=(\det h)t/\omega,\\
C'&=C_0Ch^{-1},& E'&=E_0Eh^{-1},& s'&=\omega s/\det h.
\end{aligned}
\label{m8:s:v:e4-10}
\end{equation}
In particular \(t's'=ts\), consistent with preservation of the multiplier by \(K\).

\end{proof}

\paragraph{Canonical compatible slices and residual equivalence}\label{m8:s:v:s7}

\begin{proposition}\label{m8:s:v:slices}
The compatible loci admit the normal forms \(X_6(\rho,V)\) and \(X_7(\rho,V)\) of \eqref{m8:s:v:e7-3} and \eqref{m8:s:v:e7-8}, with \((\rho,V)\in F^\times\times\GL_2(F)\). Their complete residual equivalences are diagonal conjugation and the action \(V\mapsto n(a)V\bar n(b)\), respectively. The residual characters are \(1\) and \(\psi(2a+2b)\). Both slices satisfy \eqref{m8:s:v:e7-14}.
\end{proposition}
\begin{proof}\leavevmode

\subparagraph{$S_6$ normalization}\label{m8:s:v:s7-1}
Assume \eqref{m8:s:v:e5-12}. Choose
\[
 B_0=\begin{pmatrix}1&-M_{12}/M_{22}\\0&1\end{pmatrix},
 \quad A_0=B_0M,\quad \omega=N_{22}^{-1},\quad
 q_0=-N_{12}/N_{22},\quad r_0=-N_{21}/N_{22}.
\]
Then $A_0$ is lower triangular with diagonal
$(\det M/M_{22},M_{22})$, and $B_0$ is upper triangular.
Use \eqref{m8:s:v:e4-2} with these entries, $p=y=z=0$,
and $h=B_0B^{-1}$ in \eqref{m8:s:c:normal-equations}.
The complete output is
\begin{equation}\label{m8:s:v:e7-2}
\begin{aligned}
 A'=B'&=I,\quad t'=N_{22}/\kappa=2,\quad\rho=\det M\det N,\\
 V=E'&=\begin{pmatrix}(M_{22}N_{22})^{-1}&0\\N_{21}/N_{22}&1\end{pmatrix}
 EBB_0^{-1},\\
 C'E'^{-1}&=\operatorname{diag}(\rho,1),\quad
 \det V=\frac{\det E\det B}{M_{22}N_{22}},\quad s'=2\det V.
\end{aligned}
\end{equation}
It is $X_6(\rho,V)$, of multiplier $-4\det V=-ts$.
All divisors are nonzero because $M_{22}=2\ell$,
$N_{22}=2\kappa$, and the original matrices are invertible.
Both corrections from \eqref{m8:s:c:normal-transporter} have
$X=y=0$, so the normalization has character one.

\subparagraph{$S_6$ residual relation}\label{m8:s:v:s7-2}
Equal left canonical factors force $\pi_L(j)=D(h)$.
The opposite triangular conditions imply
$h=\operatorname{diag}(\alpha,\delta)$. The new right ratio is
\begin{equation}\label{m8:s:v:e7-4}
 \begin{pmatrix}\rho+q_0r_0&\alpha q_0/\delta\\
 \delta r_0/\alpha&1\end{pmatrix}.
\end{equation}
To remain diagonal it forces $q_0=r_0=0$ and $\rho'=\rho$;
the entire relation is therefore
\begin{equation}\label{m8:s:v:e7-5}
 V'=hVh^{-1},\qquad h\in T_2,\qquad\rho'=\rho.
\end{equation}
Conversely $(1,1,h)$ realizes every such relation with character
one; the remaining freedom is exactly the Heisenberg kernel.

\subparagraph{$S_7$ normalization}\label{m8:s:v:s7-3}
Under \eqref{m8:s:v:e5-16}, $M$ is lower and $N$ upper
triangular, with $N_{22}=-M_{22}\ne0$.
Choose \eqref{m8:s:v:e4-5} with
$A_0=M$, $B_0=I$, $\omega=\kappa^{-1}$,
$X=Z=0$, $y=0$, $z_0=N_{12}/M_{22}$, $h=B^{-1}$.
It gives left factor $1$ and right factors
\[
 V_1=\operatorname{diag}(1,\kappa^{-1})EB,\quad
 C_1=\operatorname{diag}(M_{11}N_{11},-1)V_1,\quad
 s_1=\det V_1/\ell.
\]
The upper-right cancellation is exactly
$N_{12}+(N_{12}/M_{22})N_{22}=0$.
The pure diagonal move $h_2=\operatorname{diag}(1,\ell)$,
whose second projection has leading blocks $I,I$ and multiplier
$\ell$, now yields
\begin{equation}\label{m8:s:v:e7-7}
 \rho=M_{11}N_{11},\quad
 V=\operatorname{diag}(1,\kappa^{-1})EB
       \operatorname{diag}(1,\ell^{-1}).
\end{equation}
The result is $X_7(\rho,V)$ of multiplier $\det V=ts$.
The two corrections have zero $X$-trace and $y$, including any
scalar central correction, so their characters are one.
All divisions are by forced nonzero entries or determinants.
Expanded substitutions for both normalizations are in
\texttt{anc/model8\_calculations.pdf}, calculation M8-C05,
under ``Canonical compatible slices and residual equivalence''.

\subparagraph{The \(S_7\) residual subgroup and necessity of its action}\label{m8:s:v:s7-4}

Write
\[
n(a)=\begin{pmatrix}1&a\\0&1\end{pmatrix},\qquad
\bar n(b)=\begin{pmatrix}1&0\\b&1\end{pmatrix},
\qquad u_\rho(a)=u(0,a,0,-\rho a).
\]
All nonidentity blocks of \(u_\rho(a)\) are scalar. It commutes with every \(D(h)\), and \(u_\rho(a)u_\rho(a')=u_\rho(a+a')\). Therefore
\begin{equation}
k_\rho(a,b)=\bigl(u_\rho(a),u_\rho(b),\bar n(-b)\bigr)
\label{m8:s:v:e7-9}
\end{equation}
is genuinely a homomorphism from the additive group \(F^2\) to the established \(K\). Multiplication, with the semidirect action from \eqref{m8:s:v:e1-3}, gives
\begin{equation}
k_\rho(a,b)X_7(\rho,V)=X_7(\rho,n(a)V\bar n(b)),
\qquad
\eta_K(k_\rho(a,b))=\psi(2a+2b).
\label{m8:s:v:e7-10}
\end{equation}

For necessity, equality of the first Levi factors again forces \(\pi_L(j)=D(h)\), with \(h=\begin{psmallmatrix}\alpha&0\\\gamma&\delta\end{psmallmatrix}\). Its second projection has leading blocks
\(\alpha n(z_0),\) \(\alpha n(y_0).\)
The new scalar multiplier remains \(s'=\det V\), while
\(\det E'=(\alpha/\delta)\det V.\)
The canonical equality \(s'=\det E'\) forces \(\delta=\alpha\). The new right ratio is
\[
n(z_0)\operatorname{diag}(\rho,-1)n(-y_0)
=\begin{pmatrix}\rho&-\rho y_0-z_0\\0&-1\end{pmatrix}.
\]
It follows that \(\rho'=\rho\) and \(z_0=-\rho y_0\). Finally
\(V'=n(y_0)V\bar n(-\gamma/\alpha).\)
Thus the complete residual relation is
\begin{equation}
\rho'=\rho,\qquad V'=n(a)V\bar n(b),\qquad (a,b)\in F^2.
\label{m8:s:v:e7-11}
\end{equation}
Its remaining transporter ambiguity consists of the four-dimensional kernel and an ineffective scalar \(F^\times\). 

\subparagraph{Full matrices and transpose}\label{m8:s:v:s7-5}

The product definitions \eqref{m8:s:v:e7-3} and \eqref{m8:s:v:e7-8} give
\begin{equation}
D(V^{-t})X_i(\rho,V)^tD(V^t)=X_i(\rho,V^t),\qquad i=6,7.
\label{m8:s:v:e7-14}
\end{equation}
The compensator
\(a(V)=(1,1,V^{-t})\)
has character one. The fiber and density calculations above, together with the transpose
and distribution descent in \S\ref{m8:supp:s6s7-details}, determine
its effect on distributions.

\end{proof}

\paragraph{fiber actions and surjective submersions}\label{m8:s:v:s8}

\begin{proposition}\label{m8:s:v:fiber-torsors}
The maps \(\pi_i:K\times\Sigma\to\mathcal C_i\) in \eqref{m8:s:v:e8-1} are surjective analytic submersions. Their fibers are free torsors for the unimodular groups \(H_6=H_3\times T_2\) and \(H_7=(H_4\rtimes_\beta\mathbb G_a^2)\times\mathbb G_m\), with the actions below and with relative dimensions \(5\) and \(7\).
\end{proposition}
\begin{proof}\leavevmode

Let
\begin{equation}
\Sigma=F^\times\times GL_2(F),\qquad
\pi_i:K\times\Sigma\longrightarrow\mathcal C_i,
\quad \pi_i(k,s)=kX_i(s).
\label{m8:s:v:e8-1}
\end{equation}
The coverage and normalization above prove surjectivity over \(F\). The rational pivot constructions also give local sections. Here are explicit fiber actions, which supply the remaining smoothness and descent information.

\subparagraph{Canonical kernel embeddings}\label{m8:s:v:s8-1}

Let \(t_{6,s}(p,y,z)\in K\) be the stabilizer \eqref{m8:s:v:e5-10} evaluated at
\[l=L(I,I,2),\qquad r=L(\operatorname{diag}(\rho,1)V,V,2\det V).\]
Its left component is
\(\operatorname{Ad}(L(I,I,2))n_6(p,y,z).\)
This left component is independent of \(s=(\rho,V)\) and determines the three parameters uniquely. The right component is determined by the stabilizer equation. Each \(t_{6,s}\) is an embedding of the fixed Heisenberg group \(H_3\), and its character is one.

Likewise let \(t_{7,s}(p,q,r,z)\) be \eqref{m8:s:v:e5-14} at the canonical \(S_7\) Levi parameters. Its left component is
\[
u\left(\begin{pmatrix}0&0\\p&q\end{pmatrix},0,
       \begin{pmatrix}0&0\\r&z\end{pmatrix},0\right),
\]
again independent of \(s\) and injective in the four parameters. Its character is one. Denote this fixed abstract kernel by \(H_4=\mathbb G_a^4\).

\subparagraph{Fiber covariance on \(S_6\)}\label{m8:s:v:s8-2}

Let \(g(h)=(1,1,h)\), \(h\in T_2\), and let \(h\cdot(\rho,V)=(\rho,hVh^{-1})\). Direct conjugation gives
\begin{equation}
g(h)t_{6,s}(n)g(h)^{-1}=t_{6,h\cdot s}(n).
\label{m8:s:v:e8-2}
\end{equation}
The torus fixes all three Heisenberg parameters. Thus the fiber group is the direct product
\(H_6=H_3\times T_2.\)

\subparagraph{Fiber covariance on \(S_7\)}\label{m8:s:v:s8-3}

Write \(b=(a,b_2)\in F^2\) for the residual additive parameters, to avoid confusing them with matrix entries. On \(H_4\), conjugation by \(k_\rho(a,b_2)\) is
\begin{equation}
\beta_{a,b_2}(p,q,r,z)=
\bigl(p+b_2q,\ q,\ r+b_2z-a(p+b_2q),\ z-aq\bigr).
\label{m8:s:v:e8-4}
\end{equation}
More precisely,
\begin{equation}
k_\rho(a,b_2)t_{7,s}(n)k_\rho(a,b_2)^{-1}
=t_{7,(a,b_2)\cdot s}(\beta_{a,b_2}n).
\label{m8:s:v:e8-5}
\end{equation}
The matrix of \(\beta\), in the order \((p,q,r,z)\), is
\begin{equation}
\begin{pmatrix}
1&b_2&0&0\\
0&1&0&0\\
-a&-ab_2&1&b_2\\
0&-a&0&1
\end{pmatrix},
\qquad \det\beta_{a,b_2}=1.
\label{m8:s:v:e8-6}
\end{equation}
It is an additive-group representation; this also follows from the homomorphism property of \eqref{m8:s:v:e7-9}. The fiber group is therefore explicitly
\(H_7=(H_4\rtimes_\beta\mathbb G_a^2)\times\mathbb G_m.\)
Its noncentral part is six-dimensional and unipotent, while the
final factor is the ineffective scalar subgroup
\(z(\lambda)=(1,1,\lambda I_2)\subset K.\)
The embeddings depend algebraically on \(\rho\) and realize this group law for every \(\rho\).

\subparagraph{Free fiber actions}\label{m8:s:v:s8-4}

This gives a uniform construction for both cases. Let \(B_6=T_2\), \(B_7=\mathbb G_a^2\), with residual transporters \(b_s=g(h)\) or \(b_s=k_\rho(a,b_2)\). For a fiber-group element \(j=(n,b)\), and also \(\lambda\) in the \(S_7\) case, define
\begin{equation}
\sigma_s(j)=t_{i,b\cdot s}(n)\,b_s
\quad\text{or}\quad
\sigma_s(j)=t_{7,b\cdot s}(n)\,k_\rho(b)\,z(\lambda).
\label{m8:s:v:e8-8}
\end{equation}
It sends \(X_i(s)\) to \(X_i(b\cdot s)\). Covariance \eqref{m8:s:v:e8-2} or \eqref{m8:s:v:e8-5} shows the cocycle identity
\(\sigma_{j'\cdot s}(j)\sigma_s(j')=\sigma_s(jj').\)
Consequently
\begin{equation}
j\cdot(k,s)=\bigl(k\sigma_s(j)^{-1},\ b\cdot s\bigr)
\label{m8:s:v:e8-9}
\end{equation}
is an analytic group action and leaves \(\pi_i\) unchanged.

It is \emph{free}. On \(S_6\), equality of the \(K\)-components in \eqref{m8:s:v:e8-9} forces the diagonal parameter \(h\) to be the identity and then forces the left Heisenberg matrix to be the identity. On \(S_7\), the diagonal projection \(\lambda\bar n(-b_2)\) first determines \(\lambda,b_2\); the left \(y\)-coordinate then determines \(a\); finally the four left kernel coordinates determine \(n\). Thus a fiber-group element fixing a point of \(K\times\Sigma\) is the identity.

Conversely, two points of \(K\times\Sigma\) with the same image have slice parameters related by exactly \eqref{m8:s:v:e7-5} or \eqref{m8:s:v:e7-11}. After choosing the residual transporter, their remaining difference is exactly the stabilizer described on
pp.~\pageref{m8:s:v:s6}--\pageref{m8:s:v:s6-3}. The diagonal projection and left kernel coordinates recover a unique element in \eqref{m8:s:v:e8-8}. Hence every fiber is one torsor for \(H_i\).

This proves freeness throughout, including at residual stabilizer jumps.

Combining the pivot sections with the explicit recovery of the
fiber coordinates gives local analytic trivializations
\begin{equation}
\pi_i^{-1}(W)\simeq H_i\times W.
\label{m8:s:v:e8-10}
\end{equation}
Thus $\pi_i$ is a surjective submersion. Its relative dimensions are
\(\dim H_6=3+2=5,\) \(\dim H_7=4+2+1=7,\)
and
\(24+5-5=24=\dim\mathcal C_6,\) \(24+5-7=22=\dim\mathcal C_7.\)

\end{proof}

\paragraph{Modular functions and the relative density line}\label{m8:s:v:s9}

\begin{proposition}\label{m8:s:v:density}
The total density \(dk\,d^\times\rho\,d^\times V\), divided by the fiber Haar density of \(H_i\), descends to a smooth nowhere-zero \(K\)-invariant density \(\omega_i\) on \(\mathcal C_i\). Its relative density character is trivial.
\end{proposition}
\begin{proof}\leavevmode

\subparagraph{Restriction of the parabolic modulus}\label{m8:s:v:s9-1}
The four-factor calculation \eqref{m8:s:f:general-parabolic-modulus} gives
$\delta_P(L(A,B,t))=|\det A|^5|\det B||t|^{-6}$, including the scalar central coordinate. Its restriction to $D(h)=L(h,h,\det h)$ is one, so $\Delta_U=\Delta_{U\times U}=\Delta_K=1$. The present fiber map is $K$-equivariant and its relative density
is determined by $K$ and the fiber groups.

\subparagraph{Fiber Haar measures and their normalizers}\label{m8:s:v:s9-2}

For \(H_3\), the group law \eqref{m8:s:v:e5-7} preserves \(dp\,dy\,dz\) under both left and right translation. For \(H_4\), Haar measure is \(dp\,dq\,dr\,dz\). Use multiplicative Haar measure on each torus or scalar factor.

For \(S_6\), \eqref{m8:s:v:e8-2} fixes all three kernel
parameters, giving determinant one for the residual torus action.
For \(S_7\), \eqref{m8:s:v:e8-6} gives the full normalizer
determinant one, and the scalar factor acts trivially; we use
\(da\,db_2\) as Haar measure on the additive residual group.
Thus both complete fiber groups \(H_6,H_7\) are unimodular. 

\subparagraph{Quotient density: construction and invariance}\label{m8:s:v:s9-3}

Choose Haar measure \(dk\) on \(K\), and put
\(d\sigma=d^\times\rho\,d^\times V,\) \(d^\times V=|\det V|^{-2}\,da\,db\,dc\,dd.\)
The measure \(d^\times V\) is the usual bi-invariant Haar density on \(GL_2\). It is preserved by diagonal conjugation and by left or right unipotent translation. Thus all residual actions preserve \(d\sigma\).

For fixed fiber parameters, \eqref{m8:s:v:e8-9} is a right translation on the \(K\)-coordinate by a matrix depending analytically on \(s\), together with a measure-preserving residual transformation on \(\Sigma\). Since \(K\) is unimodular, it preserves \(dk\,d\sigma\). Integration first in $k$ proves invariance; equivalently the $s$-derivatives of the translation occupy only off-diagonal blocks of a triangular differential.

Apply the tangent density identity \eqref{m8:s:c:density-line}. Here the explicit free torsors \eqref{m8:s:v:e8-9} and local trivializations \eqref{m8:s:v:e8-10} identify the vertical density with the just-computed $dh_i$. Their transition translations and normalizers have modulus one, so
\begin{equation}\label{m8:s:v:e9-6}
 \pi_i^*\omega_i=\frac{dk\,d\sigma}{dh_i}
\end{equation}
is independent of the chosen fiber coordinate and descends to a smooth nowhere-zero density on $\mathcal C_i$.

Left \(K\)-translation commutes with the fiber action and preserves both factors in \eqref{m8:s:v:e9-6}, so \(\omega_i\) is \(K\)-invariant. 

The ineffective scalar on \(S_7\) has density
\(d^\times\lambda\), character one and trivial action on every
other density factor. Its contribution to the relative density
character is therefore trivial.

\end{proof}

\subsection{Rank-one incidence: intersection and transpose calculations}
\label{m8:supp:s6s7-details}
The following calculations prove the compatible-support and transpose
statements used in \S\ref{m8:s:sec:s6-s7}. Throughout,
$L(A,B,t)=m(A,B,t)$ and the final scalar coordinate of $U$ is
written $z_0$, distinct from the entry $Z_{22}$. Thus
$\dim L=9$, $\dim U=10$, $\dim P=19$ and $\dim K=24$;
we use the faithful two-factor description \eqref{m8:s:v:e1-3}.
The rational parabolic coverage, canonical normalizations and fiber
torsors are proved in \S\ref{m8:supp:s6s7}.

\paragraph{Rank geometry and representatives}\label{m8:s:v:s2}

\subparagraph{Complementary ranks}\label{m8:s:v:s2-1}

Let \(r_{a,b}(x)\) denote the rank of the lower-right \(a\)-by-\(b\) submatrix. The block sizes used below are even. The complementary-rank proof in \S\ref{m8:s:sec:ranks} gives
\begin{equation}\label{m8:s:v:e2-1}r_{a,b}=a+b-8+r_{8-a,8-b},\qquad
 r_{6,6}=4+r_{2,2},\quad r_{6,4}=2+r_{2,4},\quad
 r_{4,6}=2+r_{4,2},\quad r_{2,6}=r_{6,2}.
\end{equation}

Put
\[
d=x_{44},\quad H=\begin{pmatrix}x_{33}&x_{34}\\x_{43}&d\end{pmatrix},
\quad S=(x_{43}\ d),\quad C_*=\binom{x_{34}}d.
\]
On both target strata,
\begin{equation}
\operatorname{rank}d=1,\quad \operatorname{rank}H=2,
\quad \operatorname{rank}S=\operatorname{rank}C_*=1.
\label{m8:s:v:e2-3}
\end{equation}

For rank-one \(d\), \(\ker d^\#=\operatorname{im}d\), and the corresponding row-space identity holds. Therefore
\(\operatorname{rank}(x_{43}\ d)=1 \iff d^\#x_{43}=0,\)
\begin{equation}
\operatorname{rank}\binom{x_{34}}d=1
\iff x_{34}d^\#=0.
\label{m8:s:v:e2-4}
\end{equation}

\subparagraph{The rank-one Schur expression}\label{m8:s:v:s2-2}

On \(d_{ij}\ne0\), set
\[
d^{(1)}=\frac{e_je_i^t}{d_{ij}}.
\]
Writing \(d=uv^t\), we have \(d_{ij}=u_iv_j\), and hence \(dd^{(1)}d=u(v_ju_i/d_{ij})v^t=uv^t=d\). The incidence conditions give
\(x_{34}=av^t,\) \(x_{43}=ub^t,\)
so
\(x_{34}d^{(1)}x_{43}=ab^t.\)
More generally this equality holds for every generalized inverse satisfying \(dd^{(1)}d=d\). Hence
\(F_0=x_{33}-x_{34}d^{(1)}x_{43}\)
is independent of the selected pivot on the incidence locus. Elementary row and column operations at the one-dimensional pivot give
\(\operatorname{rank}H=1+\operatorname{rank}F_0.\)
Thus \(\operatorname{rank}F_0=1\) in \eqref{m8:s:v:e2-3}. 

Let \(\mathcal Y\) be the locally closed rank geometry \eqref{m8:s:v:e2-3}. In \(\mathcal Y\),
\(r_{2,6}=\operatorname{rank}(x_{42}\ x_{43}\ d) =1\iff d^\#x_{42}=0.\)
The transposed statement, combined with \eqref{m8:s:v:e2-1}, gives
\begin{equation}
\begin{aligned}
S_7&=\{x\in\mathcal Y:d^\#x_{42}=0\}
=\{x\in\mathcal Y:x_{24}d^\#=0\},\\
S_6&=\{x\in\mathcal Y:d^\#x_{42}\ne0\}
=\{x\in\mathcal Y:x_{24}d^\#\ne0\}.
\end{aligned}
\label{m8:s:v:e2-6}
\end{equation}
The other entries of the rank matrix follow from \eqref{m8:s:v:e2-1}. The rank table gives exactly \(\mathcal Y=S_6\sqcup S_7\). \(S_6\) is open and \(S_7\) is closed in \(\mathcal Y\).

\subparagraph{Representatives}\label{m8:s:v:s2-3}

For the representatives \eqref{m8:s:v:e2-7},
direct multiplication gives
\(w_i^t=w_i,\) \(w_i^2=I_8,\) \(\det w_i=1,\)
\(w_6^tQw_6=-Q,\) \(w_7^tQw_7=Q.\)
Both therefore satisfy the defining determinant condition for the stipulated connected GSO group, with
\(\mu(w_6)=-1,\) \(\mu(w_7)=1.\)

The nine lower-right ranks are
\begin{equation}
R(w_6)=\begin{pmatrix}5&3&2\\3&2&1\\2&1&1\end{pmatrix},
\qquad
R(w_7)=\begin{pmatrix}5&3&1\\3&2&1\\1&1&1\end{pmatrix}.
\label{m8:s:v:e2-9}
\end{equation}
These are the rank matrices of \(S_6\) and \(S_7\).

\paragraph{Unipotent stabilizers and compatible support}\label{m8:s:v:s5}
\begin{proof}[Proof of Lemma~\ref{m8:s:v:compatible}]\leavevmode

Put
\(M=B^{-1}A,\) \(N=CE^{-1},\) \(\kappa=\det B/t,\) \(\ell=\det E/s.\)
All four matrices \(A,B,C,E\) are invertible, and \(\kappa\ell\ne0\).

\subparagraph{Lie intersections from all ten coordinates}\label{m8:s:v:s5-1}

The Lie algebra is
\begin{equation}
n(X,y,Z,z_0)=
\begin{pmatrix}
0&X&Z&z_0I\\
0&0&yI&Z^\#\\
0&0&0&-X^\#\\
0&0&0&0
\end{pmatrix}.
\label{m8:s:v:e5-2}
\end{equation}
To belong to \(\mathfrak u^t\), \(w_i n w_i\) must have zero diagonal blocks as well as zero upper blocks. Solving all those equations gives:
\begin{equation}
\begin{array}{c|l}
S_6 &X_{11}=X_{12}=X_{21}=Z_{11}=Z_{12}=Z_{21}=z_0=0,\\
S_7 &X_{11}=X_{12}=Z_{11}=Z_{12}=y=z_0=0.
\end{array}
\label{m8:s:v:e5-3}
\end{equation}
Thus
\begin{equation}
\begin{aligned}
S_6:&\quad X=pE_{22},\quad y\in F,\quad Z=zE_{22},\quad z_0=0;\\
S_7:&\quad X=\begin{pmatrix}0&0\\p&q\end{pmatrix},\quad y=0,
\quad Z=\begin{pmatrix}0&0\\r&z\end{pmatrix},\quad z_0=0.
\end{aligned}
\label{m8:s:v:e5-4}
\end{equation}
The distinction from the \(J_i\) computation is important: diagonal Levi blocks were allowed in \S\ref{m8:s:v:s4} but must vanish in \eqref{m8:s:v:e5-3}.

The finite exponential/logarithm principle in \S\ref{sec:common} identifies these complete Lie intersections with the entire group intersections.

For \(S_6\), if the Lie coordinate of \(Z_{22}\) is \(z\), then
\begin{equation}
\exp n(pE_{22},y,zE_{22},0)
=u\bigl(pE_{22},y,(z+py/2)E_{22},0\bigr).
\label{m8:s:v:e5-5}
\end{equation}
Since \(n^3=0\), exponentiation gives the displayed
\(py/2\) shift. For \(S_7\), \(n^3=0\) as well, and
\(\exp n(X,0,Z,0)=u(X,0,Z,0)\)
for the row shapes in \eqref{m8:s:v:e5-4}. The upper-right block still contains the quadratic entry $pz-qr$ prescribed by \eqref{m8:s:p:U}.

\subparagraph{\(S_6\): group law, inverse-transpose identity, and character}\label{m8:s:v:s5-2}

Set \(n_6(p,y,z)=u(pE_{22},y,zE_{22},0)\). Direct multiplication gives
\begin{equation}
(p,y,z)(p',y',z')=(p+p',y+y',z+z'+py'),
\label{m8:s:v:e5-7}
\end{equation}
so this is the three-dimensional Heisenberg group. Its inverse is
\((p,y,z)^{-1}=(-p,-y,py-z).\)
The inverse-transpose identity is
\begin{equation}
(w_6n_6(p,y,z)^{-1}w_6)^t
=u(-yE_{22},-p,(py-z)E_{22},0).
\label{m8:s:v:e5-8}
\end{equation}

The Levi conjugation formula, valid in the coordinates, is
\begin{equation}
\operatorname{Ad}(L(A,B,t))u(X,y,Z,z_0)
=u\left(AXB^{-1},\kappa y,\kappa AZB^{-1},\frac{\det A}{t}z_0\right).
\label{m8:s:v:e5-9}
\end{equation}
Conjugating \eqref{m8:s:v:e5-8} by \(l\) and \(r^t\) gives every \(U\times U\)-stabilizer element at \(lw_6r\):
\begin{equation}
\begin{aligned}
u&=u(ApE_{22}B^{-1},\kappa y,\kappa AzE_{22}B^{-1},0),\\
v&=u(-C^tyE_{22}E^{-t},-\ell p,
       \ell C^t(py-z)E_{22}E^{-t},0).
\end{aligned}
\label{m8:s:v:e5-10}
\end{equation}
Cyclicity of trace gives
\(\operatorname{tr}(ApE_{22}B^{-1})=M_{22}p,\) \(\operatorname{tr}(-C^tyE_{22}E^{-t})=-N_{22}y.\)
Adding the two \(2y\)-terms therefore gives
\eqref{m8:s:v:e5-11}.
The correction $py-z$ lies in a character-trivial $Z$-coordinate. Both minus signs come from the right inverse \eqref{m8:s:v:e5-8}.

Because \(p,y\) independently run over all of \(F\), compatibility is exactly
\eqref{m8:s:v:e5-12}.
To check independence, use the normal-family coordinates
\((M,N,B,E,\kappa,\ell)\). The Jacobian with respect to
\((\ell,\kappa)\) is diagonal with nonzero entries \(-2,-2\),
including in residue characteristic two.

\subparagraph{\(S_7\): the additive kernel and its character}\label{m8:s:v:s5-3}

For the row-shaped \(X,Z\) in \eqref{m8:s:v:e5-4}, put
\[
\widehat X=\begin{pmatrix}0&0\\-z&q\end{pmatrix},
\qquad
\widehat Z=\begin{pmatrix}0&0\\r&-p\end{pmatrix}.
\]
The identity is
\((w_7u(X,0,Z,0)^{-1}w_7)^t=u(\widehat X,0,\widehat Z,0).\)
Multiplication adds $(p,q,r,z)$, so the kernel is $\mathbb G_a^4$ despite its quadratic matrix entries.

The complete stabilizer is
\begin{equation}
\begin{aligned}
u&=u(AXB^{-1},0,\kappa AZB^{-1},0),\\
v&=u(C^t\widehat XE^{-t},0,\ell C^t\widehat ZE^{-t},0).
\end{aligned}
\label{m8:s:v:e5-14}
\end{equation}
The two trace contributions are
\(\operatorname{tr}(MX)=M_{12}p+M_{22}q,\) \(\operatorname{tr}(N^t\widehat X)=-N_{21}z+N_{22}q.\)
The character is therefore \eqref{m8:s:v:e5-15},
and compatibility is exactly \eqref{m8:s:v:e5-16}.
The coordinate \(r\) is character-trivial. The three equations are independent in the invertible-matrix normal family. On their common zero set, \(M\) is lower triangular and \(N\) upper triangular, so \(M_{22}\) and \(N_{22}\) are nonzero.

\subparagraph{Intrinsic compatible loci and distribution-theoretic support reduction}\label{m8:s:v:s5-4}

Let \(\mathcal C_i\) be the compatible locus in \(S_i\). The kernel of the infinitesimal \(U\times U\)-action is a vector bundle of constant rank 3 on \(S_6\) and 4 on \(S_7\). The algebraic additive functional
\(\Lambda(u,v)=\operatorname{tr}X_u+2y_u+\operatorname{tr}X_v+2y_v\)
restricts to zero on a stabilizer exactly when the complex character \(\psi\circ\Lambda\) is trivial on that stabilizer. Equivalently, the algebraic differential of $\Lambda$ vanishes
on the stabilizer Lie algebra.

This zero condition is algebraic, closed and preserved by \(K\)
and transpose. In the surjective normal-family charts it is defined
by the independent equations \eqref{m8:s:v:e5-12} or
\eqref{m8:s:v:e5-16}. Descending smoothness through these charts gives
\(\dim\mathcal C_6=26-2=24,\) \(\dim\mathcal C_7=25-3=22,\)
and both loci are smooth.

On each algebraic rational-point complement $S_i\setminus\mathcal C_i$, the complete unipotent stabilizer is unimodular and its character is nontrivial by \eqref{m8:s:v:e5-11} or \eqref{m8:s:v:e5-15}. The criterion \eqref{cf:orbital} therefore gives $\De(S_i\setminus\mathcal C_i)=0$. Lemma~\ref{cf:support}, applied to the closed stable pair $(S_i,\mathcal C_i)$, identifies every equivariant distribution with its unique closed-support distribution on $\mathcal C_i$.

\end{proof}

\paragraph{Full stabilizer characters}\label{m8:s:v:s6}
\begin{proof}[Proof of Lemma~\ref{m8:s:v:full-stabilizers}]\leavevmode

\subparagraph{The S6 stabilizer character}
For a stabilizer element, the left and right matrix traces are $M_{22}p$ and $-N_{22}y$; the scalar terms are $2\kappa y$ and $-2\ell p$. Conjugation by its diagonal component leaves these traces unchanged. Hence the complete character is
\begin{equation}\label{m8:s:v:e6-3}
 \eta_K(k)=\psi\bigl((M_{22}-2\ell)p+(2\kappa-N_{22})y\bigr).
\end{equation}
The diagonal stabilizer preserves the compatibility conditions
\eqref{m8:s:v:e5-12}. Its complete general-coordinate formula is in Appendix~\S\ref{m8:s:app:B}; the canonical stabilizer and all density data are given below. On $\mathcal C_6$, the lines $Ae_2$ and $Be_1$ are distinct because $M_{22}=2\ell\ne0$; the projection makes them eigenlines with eigenvalues $\delta,\alpha$. Thus the projection is contained in the split torus preserving these two lines. The exact scalar or split-torus alternatives follow from \eqref{m8:s:v:e12-1}--\eqref{m8:s:v:e12-2}.

\subparagraph{\(S_7\) stabilizer projection: determinant and trace checks}\label{m8:s:v:s6-2}

Use the same stabilizer equations with \eqref{m8:s:v:e4-6}. The first two determinant comparisons give
\(\alpha=\delta,\) \(\epsilon=\zeta.\)
The two conjugates of \(h\) have equal traces, so
\(2\alpha=2\epsilon.\)
Since $\operatorname{char}F=0$, $\epsilon=\alpha$; the trace comparison excludes the minus-sign alternative left by determinants.

The complete projection is therefore
\begin{equation}
\begin{aligned}
A^{-1}hA&=\lambda\begin{pmatrix}1&0\\g&1\end{pmatrix},&
B^{-1}hB&=\lambda\begin{pmatrix}1&0\\f&1\end{pmatrix},\\
ChC^{-1}&=\lambda\begin{pmatrix}1&z_0\\0&1\end{pmatrix},&
EhE^{-1}&=\lambda\begin{pmatrix}1&y_0\\0&1\end{pmatrix},
\end{aligned}
\label{m8:s:v:e6-4}
\end{equation}
with \(\lambda\ne0\).

The complete matrices are reconstructed by
\eqref{m8:s:c:normal-transporter}, with the projections
\eqref{m8:s:v:e6-4} and the four free kernel parameters
$p,q,r,z$ in \eqref{m8:s:v:e5-4}. This also determines the central scalar coordinates.
The expanded $u,v,h$ arrays are in
\texttt{anc/model8\_calculations.pdf}, calculation M8-C01,
under ``$S_7$ stabilizer projection: determinant and trace checks''.
Their complete character, obtained by taking the $X$ traces and
adding both scalar $2y$ terms, is
\begin{equation}\label{m8:s:v:e6-5}
 \eta_K(k)=\psi\bigl(M_{12}p+(M_{22}+N_{22})q
 +N_{21}(qy_0-z)+2\kappa y_0-2\ell f\bigr).
\end{equation}
On $\mathcal C_7$ only $2\kappa y_0-2\ell f$ remains.

\subparagraph{The extra \(S_7\) unipotent and its boundary character}\label{m8:s:v:s6-3}

Put \(V_0=EB=\begin{psmallmatrix}a&b\\c&d\end{psmallmatrix}\). Comparing the second and fourth equations in \eqref{m8:s:v:e6-4} yields
\(fV_0E_{21}V_0^{-1}=y_0E_{12}.\)
The matrix on the left before multiplication by \(f\) is
\begin{equation}
V_0E_{21}V_0^{-1}
=\frac1{ad-bc}\begin{pmatrix}bd&-b^2\\d^2&-bd\end{pmatrix}.
\label{m8:s:v:e6-8}
\end{equation}
If \(d\ne0\), the lower-left entry forces \(f=0\), hence \(y_0=0\). Thus \(h\) is scalar.

If \(d=0\), invertibility forces \(bc\ne0\), and \eqref{m8:s:v:e6-8} becomes \((b/c)E_{12}\). Therefore
\(y_0=(b/c)f.\)
On the two-unipotent-compatible locus the remaining projection equations give
\(g=(M_{11}/M_{22})f,\) \(z_0=(N_{11}/N_{22})y_0.\)
The extra direction therefore has character
\[
\psi\bigl(2(\kappa b/c-\ell)f\bigr),
\]
and its exact compatibility condition is
\begin{equation}
\kappa b=\ell c.
\label{m8:s:v:e6-10}
\end{equation}
This direction belongs to the $K$-stabilizer but not to the \(U\times U\) kernel. A non-diagonal example illustrating its character is given in Subsection~\ref{m8:s:v:sB-2}.

\end{proof}

\paragraph{Lifted transpose and vertical Jacobians}\label{m8:s:v:s10}
\begin{proof}[Proof of Lemma~\ref{m8:s:v:transpose-density}]\leavevmode

\subparagraph{The lifted formula}\label{m8:s:v:s10-1}

Since \(a(V)X_i(\rho,V)^t=X_i(\rho,V^t)\), transpose lifts to
\eqref{m8:s:v:e10-1}.
Indeed \(\pi_i\circ\Phi_i=\tau\circ\pi_i\). The lift is itself involutive: \(\theta(a(V))^{-1}a(V^t)^{-1}=1\).

Since \(\theta\) is involutive, its positive Haar scaling has
square one and it preserves Haar measure on \(K\). In the
coordinates of \eqref{m8:s:v:e1-3}, this also follows by exchanging
the unipotent factors and applying inverse-transpose on \(GL_2\).
Right multiplication by \(a(V)^{-1}\) preserves Haar measure,
as does transpose on \(d^\times V\). The dependence on \(V\)
contributes only off-diagonal differential blocks, so \(\Phi_i\)
preserves \(dk\,d\sigma\).

\subparagraph{Kernel maps in the canonical coordinates}\label{m8:s:v:s10-2}

A vertical stabilizer element transforms by
\(t_{i,s}(n)\longmapsto a(V)\theta(t_{i,s}(n))a(V)^{-1}.\)
Using the left canonical kernel coordinates from \S\ref{m8:s:v:s8}, matrix multiplication gives
\begin{equation}
F_6(p,y,z)=(-y,-p,py-z),
\label{m8:s:v:e10-3}
\end{equation}
\begin{equation}
F_7(p,q,r,z)=(z,-q,-r,p).
\label{m8:s:v:e10-4}
\end{equation}
For \(S_6\), this calculation includes the factors \(1/2\) in the left canonical \(y\)- and \(Z\)-coordinates arising from \(L(I,I,2)\). They cancel exactly when the transformed element is expressed again in the same Heisenberg coordinates.

The Jacobian of \eqref{m8:s:v:e10-3} is
\[
\begin{pmatrix}0&-1&0\\-1&0&0\\y&p&-1\end{pmatrix},
\qquad \det=1.
\]
The map \eqref{m8:s:v:e10-4} has determinant \(-1\), hence absolute Jacobian 1. Both preserve the chosen additive Haar measures over every field in scope.

\subparagraph{Residual and scalar directions}\label{m8:s:v:s10-3}

On \(S_6\) the residual torus parameter is \(h\mapsto h^{-1}\). Including the compensator gives
\(a(hVh^{-1})\,\theta(g(h))\,a(V)^{-1}=g(h^{-t})=g(h^{-1}).\)
Multiplicative Haar measure is preserved.

On \(S_7\), put \(V'=n(a)V\bar n(b)\). Then the complete transporter identity is
\[a(V')\theta(k_\rho(a,b))a(V)^{-1}=k_\rho(b,a).\]
Thus transpose exchanges the residual additive parameters and inverts the scalar parameter, preserving $da\,db$ and $d^\times\lambda$.

These maps preserve the fiber Haar densities, including the kernel normalizer actions. By the line identity \eqref{m8:s:c:density-line}, the quotient density \(\omega_i\) is therefore genuinely transpose-invariant. 

\end{proof}

\paragraph{Ordinary distribution descent}\label{m8:s:v:s11}
\begin{proof}[Proof of Proposition~\ref{m8:s:v:distribution-correspondence}]\leavevmode

\subparagraph{Verification of the fiber-descent hypotheses}
Apply Lemma~\ref{m8:s:c:fiber} to the surjective submersions \eqref{m8:s:v:e8-1}, with the free fiber torsors \eqref{m8:s:v:e8-9} and their rational local trivializations \eqref{m8:s:v:e8-10}. The chosen total density is $dk\,d\sigma$, and the quotient is the concrete $\omega_i$ of \eqref{m8:s:v:e9-6}. Thus its injective pullback has the unique form
\begin{equation}\label{m8:s:v:e11-1}
 \pi_i^*T=\eta_K(k)dk\otimes S_\Sigma.
\end{equation}
The kernel embeddings have character one, as does the ineffective scalar. On the residual transporters,
\(\eta_K(g(h))=1,\) \(\eta_K(k_\rho(a,b))=\psi(2a+2b).\)
Hence $\eta_K(\sigma_s(j))$ is independent of the slice point and the residual equation is
\(S_\Sigma(f\circ b^{-1})=\eta_K(b_s)S_\Sigma(f).\)
The plus sign follows from the exact relation $(k,s)\sim(kb_s^{-1},b\cdot s)$ in Lemma~\ref{m8:s:c:fiber}. Conversely this equation gives invariance under the entire free action, since the kernel and scalar characters are one. The overlap maps are the translations and normalizers of \eqref{m8:s:v:e8-2} and \eqref{m8:s:v:e8-5}, with determinant one as computed in \eqref{m8:s:v:e8-6}. The lemma's local Haar descent therefore applies on every chart and agrees on overlaps, proving surjectivity as well as injectivity.

\subparagraph{Transpose and the isomorphisms}\label{m8:s:v:s11-4}

Under \eqref{m8:s:v:e10-1}, the character factor in \eqref{m8:s:v:e11-1} is unchanged:
\(\eta_K(\theta(k)a(V)^{-1})=\eta_K(k).\)
The total and quotient densities are transpose-invariant by the transpose-density calculation above. Hence transpose on the descended parameter distribution is precisely
\((\rho,V)\longmapsto(\rho,V^t).\)
We have proved the exact even correspondences and, after imposing oddness, the isomorphisms
\eqref{m8:s:v:e11-4} and \eqref{m8:s:v:e11-5}.
Here the actions are exactly \eqref{m8:s:v:e7-5} and \eqref{m8:s:v:e7-10}, and \(\rho\) is passive. 

\end{proof}

\paragraph{Canonical stabilizers and opposite-weight axes}\label{m8:s:v:s12}
The following details give the stabilizers and transpose transporters
used in Proposition~\ref{m8:s:v:s6-vanishing}.

\subparagraph{All stabilizers on the compatible canonical locus}\label{m8:s:v:s12-1}

Let \(V=\begin{psmallmatrix}a&b\\c&d\end{psmallmatrix}\). The residual calculation gives
\begin{equation}
K_{X_6(\rho,V)}\simeq H_3\times T_V,
\qquad T_V=\{h\in T_2:hVh^{-1}=V\}.
\label{m8:s:v:e12-1}
\end{equation}
The Heisenberg embedding is the complete one in \eqref{m8:s:v:e5-10}, with canonical Levi parameters, and \(T_V\) is embedded by \((1,1,h)\). Equation \eqref{m8:s:v:e8-2} shows that the latter centralizes the kernel. Explicitly,
\begin{equation}
T_V=\begin{cases}
T_2,&b=c=0,\\
\{\lambda I_2:\lambda\in F^\times\},&(b,c)\ne(0,0).
\end{cases}
\label{m8:s:v:e12-2}
\end{equation}
This list includes both punctured axes. Every such stabilizer has trivial character and relative modular character.

\subparagraph{Transpose transporters}\label{m8:s:v:s12-2}

If \(bc\ne0\), then
\(h_0=\operatorname{diag}(c/b,1)\)
sends \(V\) to \(V^t\) by conjugation. If \(b=c=0\), take \(h_0=I\). In both cases
\(k_\tau=(1,1,V^th_0)\)
sends \(X_6(\rho,V)\) to its transpose and has character one.

If exactly one of \(b,c\) is zero, \eqref{m8:s:v:e7-5} proves that no even transporter to the transpose exists. This is a pair of exchanged even orbits. 

On an axis, \(\det V=ad\ne0\). Holding \(\rho,a,d\) fixed and letting the nonzero off-diagonal entry tend to zero gives the diagonal matrix, still in \(GL_2\); the product definition \eqref{m8:s:v:e7-3} shows the limit remains in \(S_6\). The diagonal matrix is not in the original even orbit by \eqref{m8:s:v:e7-5}. Thus the exchanged orbits are not closed in \(S_6\).

The axis closures contain their diagonal limit, so the distribution argument uses the whole opposite-weight plane, including its origin.

\paragraph{Canonical Whittaker stabilizers and transpose transporters}\label{m8:s:v:s13}
The following details identify the compatible boundary and transporters
used in Proposition~\ref{m8:s:v:s7-vanishing}.

\subparagraph{Residual action and additive stabilizer}\label{m8:s:v:s13-1}

The residual action and additive stabilizer are \eqref{m8:s:v:e13-1}--\eqref{m8:s:v:e13-2} in Lemma~\ref{m8:s:c:whittaker}. Its boundary condition is $b=c$.

To compare with the full stabilizer calculation above, write the original \(V_0=EB\) entries as \(b_0,c_0\). Formula \eqref{m8:s:v:e7-7} gives
\(V_{12}=b_0/\ell,\) \(V_{21}=c_0/\kappa.\)
Therefore \(V_{12}=V_{21}\) is exactly \(\kappa b_0=\ell c_0\). The residual and $K$-stabilizer compatibility equations agree.

\subparagraph{Canonical K-stabilizers and modular functions}\label{m8:s:v:s13-2}

For \(d\ne0\), the stabilizer is \(H_4\times F^\times\), with \(F^\times\) the central ineffective scalar. It has dimension 5 and is unimodular.

For $d=0$, use \eqref{m8:s:c:normal-transporter} with the projection
\eqref{m8:s:v:e6-4} and canonical parameters
\begin{equation}\label{m8:s:v:e13-3}\begin{gathered}
 A=B=I_2,\quad t=1,\quad E=V,\quad C=\diag(\rho,-1)V,\quad s=\det V,\\
 h=\lambda\bar n(k),\quad f=g=k,\quad y_0=(b/c)k,\quad z_0=-\rho y_0.
\end{gathered}
\end{equation}
Here $p,q,r,z$ remain arbitrary and $\lambda\in F^\times$. This recovers every canonical stabilizer element, including the left scalar $-\rho y_0$ and the right scalar obtained from $-(\det C/s)g=\rho k$. Its character is
\begin{equation}\label{m8:s:v:e13-5}
 \eta_K(u,v,h)=\psi\bigl(2k(b/c-1)\bigr).
\end{equation}
The additional \(k\)-direction makes the stabilizer six-dimensional. Its action on \(H_4\) is the restriction of \eqref{m8:s:v:e8-4} to residual parameters \((a,b_2)=(y_0,-k)\); its determinant is one. Thus its noncentral part is a unipotent extension of \(H_4\) by \(\mathbb G_a\), and the remaining scalar torus is central. It is unimodular. On the surviving locus \(b=c\), its character is trivial. On \(b\ne c\), the same unimodular stabilizer has a nontrivial character and kills the locus.

Thus every surviving $S_7$ stabilizer has $\Delta_{K_x}=1$ and relative modular character one.

\subparagraph{Character-one transpose transporters}\label{m8:s:v:s13-3}

For \(d\ne0\), put
\(s_0=(c-b)/d,\) \(t_0=(b-c)/d=-s_0.\)
Substitution into \eqref{m8:s:v:e13-1} gives \(n(s_0)V\bar n(t_0)=V^t\), including the upper-left entry, and
\(\chi(s_0,t_0)=\psi(0)=1.\)
The transporter in \(K\) is
\begin{equation}
k_\tau=\bigl(u_\rho(s_0),u_\rho(t_0),V^t\bar n(-t_0)\bigr).
\label{m8:s:v:e13-6}
\end{equation}
Indeed it is the product of the compensator inverse with \eqref{m8:s:v:e7-9}; the scalar-block matrices \(u_\rho\) commute with the diagonal conjugations, giving the stated conjugation formula.

On \(d=0,b=c\), the residual matrix is symmetric, and
\(k_\tau=(1,1,V^t)\)
sends \(X_7(\rho,V)\) to its transpose with character one.

Formula \eqref{m8:s:c:transport-back} carries these character-one transporters through either rational normalization, with exact cancellation of its character.

\paragraph{Stabilizers and relative modular characters}\label{m8:s:v:s14}

\begin{center}
\begin{tabular}{lll}
\toprule
Locus & Stabilizer & Character and modulus\\
\midrule
\(S_6,\ b=c=0\)&\(H_3\times T_2\)&both trivial\\
\(S_6,\ (b,c)\ne(0,0)\)&\(H_3\times F^\times\)&both trivial\\
\(S_7,\ d\ne0\)&\(H_4\times F^\times\)&both trivial\\
\(S_7,\ d=0,\ b=c\ne0\)&\((H_4\rtimes\mathbb G_a)\times F^\times\)&both trivial\\
\bottomrule
\end{tabular}
\end{center}
The torus in the first row fixes all Heisenberg coordinates. In the second and third rows the scalar factor is central. In the last row the additive factor has determinant-one action on \(H_4\), by \eqref{m8:s:v:e8-6}.

The discarded \(S_7\) locus \(d=0,b\ne c\) also has modular function one, but has the nontrivial character \eqref{m8:s:v:e13-5}. Since \(\Delta_K=1\), every relative modular character used on a surviving orbit is exactly one.

\section{Finite data for the \texorpdfstring{$D_6$}{D6} continuation}\label{supp:d6-expanded}
The rational cells, stabilizers, phase submersions and ambient support
induction are proved in \S\ref{m9:boundary-section}. Their complete
finite arrays appear in \path{anc/model9_calculations.pdf}:
M9-K48 lists the $48$ nonempty compatibility systems and their
stabilizer and manifold dimensions; M9-C35 gives all $35$ boundary
conormals; M9-V gives the tangent identities in every Levi and radical
direction. The other $444$ compatibility systems are included in the
complete $492$-row JSON/CSV in \path{anc/model9/}. The coordinate
nonvanishing arguments in Lemma~\ref{m9:conormal-family} apply to
every point of each row.

\section{Classical finite support and localization data}\label{app:classical}
\subsection{Finite cell data and reconstruction}\label{m6:s:app:ledger}
The proof of Proposition~\ref{m6:s:prop:cells-rational} gives the
complete rational classification. Enumerate the first two rows $a,b$
as compositions of two, enforce central symmetry and the nonnegative
middle row \eqref{m6:s:eq:cells-middle-row}, compute all corner ranks,
and sort them decreasingly in row-major lexicographic order. This defines the numbering P001--P090.
For each table the adapted dual-basis construction supplies a signed
symplectic monomial representative. Its physical signs matter for
compatibility, so the exact ten signed columns of the selected $W_i$
are included in the input data.

\texttt{anc/model6\_calculations.pdf}, section \texttt{M6-CELLS}, prints
the complete ninety-row table, including the first two rank rows,
all signed columns, transpose partner, full unipotent dimension,
compatibility-system rank and support step.
The same exact fields are in \path{anc/model6/cell_ledger.csv}
and \path{anc/model6/cell_compatibility.json};
\path{cell_ledger_rows.txt} preserves the row transcription for the
independent source-to-data comparison. A negative column $-k$ means
$-e_k$; every representative has multiplier one. The other two rank rows
follow from \eqref{m6:s:eq:cells-complement}.
These data classify parabolic double cells. The rational stabilizers,
character criterion and support induction are proved in the article.

\begin{table}[ht]
\centering
\begin{tabular}{@{}c r l@{}}
\toprule Outer position $(r,s)$ & Number & Cell identifiers \\\midrule
$(2,2)$ & 5 & $1$--$5$ \\
$(2,1)$ & 23 & $6$--$28$ \\
$(1,1)$ & 5 & $29$--$33$ \\
$(2,0)$ & 29 & $34$--$50$, $59$--$65$, $73$--$77$ \\
$(1,0)$ & 23 & $51$--$58$, $66$--$72$, $78$--$85$ \\
$(0,0)$ & 5 & $86$--$90$ \\
\bottomrule
\end{tabular}
\caption{Coverage of the finite cell data; the full row table is
public section \texttt{M6-CELLS}.}\label{m6:s:tab:ledger}
\end{table}
The verifier reconstructs all ninety intersections and linear systems
from $J$, the complete twenty-two-root basis and these physical
representatives. It compares every rank, transpose partner, dimension,
equation row space and witness, as specified in
\S\ref{m6:s:app:compatibility}. The separate independent verifier also
enumerates the signed Weyl placements and checks equality of the
resulting rank-table set. The displayed $W_{90}$ data
remain distinct from the block reversal $w_0$ in the terminal proof.

\subsubsection{The transpose involution}
The nontrivial pairs are listed here separately for convenience. All other
indices are fixed. Each pair is obtained by transposing the full rank matrix,
and each occurs exactly once.
\[
\begin{gathered}
(9,14),\quad (10,21),\quad (11,15),\quad (12,22),\quad (13,23)\\
(18,24),\quad (20,26),\quad (35,39),\quad (36,46),\quad (37,59)\\
(38,73),\quad (42,47),\quad (43,60),\quad (44,61),\quad (45,74)\\
(49,62),\quad (50,75),\quad (54,66),\quad (55,78),\quad (56,67)\\
(57,79),\quad (58,80),\quad (65,76),\quad (70,81),\quad (72,83)
\end{gathered}
\]
There are twenty-five pairs and forty fixed indices. The original rank
definition proves that this numerical involution is the action of transpose
on the corresponding rational cells.

\subsection{Degenerations and support boundaries}\label{m6:s:app:support}

The following tables supplement Table~\ref{m6:s:tab:assembly-filtration}.
An interval $a$--$b$ includes every integer between its endpoints.
Membership follows from the full rank matrices in
the public \texttt{M6-CELLS} rank table specified in
\S\ref{m6:s:app:ledger} and the entrywise inequalities defining $V_i$.
Every $B_i$ in this appendix is the open $V_i\setminus\mathcal P_i$,
as expressed by \eqref{m6:s:eq:cells-boundary}.

\begin{longtable}{@{}r p{12.8cm}@{}}
\caption{The early and middle rank-open boundaries.}\label{m6:s:tab:early-boundaries}\\
\toprule $i$&Cell indices in $B_i$\\\midrule\endfirsthead
\toprule $i$&Cell indices in $B_i$\\\midrule\endhead
\bottomrule\endlastfoot
6&1\\
9&1,2,6,7\\
14&1,2,6,7\\
12&1,2,3,4,6,7,8,9,10,11\\
22&1,2,3,4,6,7,8,14,15,21\\
16&1,2,6,7,9,14\\
17&1,2,3,4,6,7,8,9,11,14,15,16\\
27&1--26\\
\end{longtable}

\begin{longtable}{@{}r p{11.8cm} c@{}}
\caption{Complete late boundaries and the stage at which they vanish.}\label{m6:s:tab:late-boundaries}\\
\toprule $i$&Cell indices in $B_i$&Stage\\\midrule\endfirsthead
\toprule $i$&Cell indices in $B_i$&Stage\\\midrule\endhead
\bottomrule\endlastfoot
40&1--2,6--7,9,14,16,34--35,39&14\\
41&1--4,6--9,11,14--17,34--35,39--40&15\\
53&1--18,21--25,29--31,34--36,39--42,46--48,51--52&15\\
63&1--27,34--37,39--44,46--49,59--62&16\\
69&1--27,29--32,34--37,39--44,46--49,51--54,56,59--63,66--68&18\\
84&1--76,78--83&20\\
\end{longtable}

The interval boundary is
\[
\mathcal B=\mathcal P[1,\ldots,18,21,\ldots,25,29,30,31,34,35,36,39,46],
\]
and its complement in $V_{53}$ is
$\mathcal P[40,41,42,47,48,51,52,53]$.
Thus the only part of $\mathcal B$ outside $\mathcal O_{34}$ is
$\mathcal P[35,36,39,46]$. This is the boundary required before the
interval theorem; the individual $B_{53}$ contains cells of that interval.

\begin{longtable}{@{}>{\raggedright\arraybackslash}p{3.0cm} >{\raggedright\arraybackslash}p{5.5cm} >{\raggedright\arraybackslash}p{5.8cm}@{}}
\caption{Spectral applications and previously established antecedents.}\label{m6:s:tab:spectral-dependencies}\\
\toprule Application&Open ambient or boundary&Earlier vanishing\\\midrule\endfirsthead
\toprule Application&Open ambient or boundary&Earlier vanishing\\\midrule\endhead
\bottomrule\endlastfoot
Model $7$ initial open&$\Omega=\mathcal O_5$&Original $T$ has nilpotent wavefront; explicit extraction transfers this property.\\
$\mathcal P_6$&$B_6=\mathcal P_1$ in $V_6$&$\Omega$ is zero at stage $1$.\\
$\mathcal P_9/\mathcal P_{14}$&$B_9=B_{14}=\mathcal P[1,2,6,7]$&$\mathcal O_8$ is zero at stage $3$.\\
$\mathcal P_{12}/\mathcal P_{22}$&The two complete boundaries in Table~\ref{m6:s:tab:early-boundaries}&$\mathcal S_{11}$ is zero at stage $5$, including the reflected incompatible cells.\\
\end{longtable}

For clarity, the thirty cells with nonempty unipotent-compatible systems
divide according to the exact statements used as follows:
\[
\begin{array}{c|l}
\text{argument}&\text{indices}\\\hline
\text{initial Model }7\text{ open}&1,2,3,4,5\\
\text{intrinsic odd vanishing}&7,8,16,17,27,31,33,34,40,53,69,90\\
\text{ordinary equivariant vanishing}&18,24,32,41,42,47,63,84\\
\text{conditional nilpotent-wavefront conclusion}&6,9,14,12,22.
\end{array}
\]
The initial open has its own spectral component at the third cell.
The last row uses conditional ambient NW vanishing; the paired spectral
cells carry the intrinsic odd orbital distributions constructed above.

\subsection{Ambient tangent calculation}\label{m6:s:p78:tangent-main}
Section~\ref{m6:s:cert:p78-tangents} gives the ambient tangent and
conormal arrays, including the Levi and similitude directions.
For the entire cells,
$CN^G_{C_i,g}=\mathfrak u\cap\operatorname{Ad}(g)\mathfrak u^t$
by \eqref{m6:s:eq:p78-11-5}; the dimensions are $2,4$, and the full
matrix arrays \eqref{m6:s:eq:p78-11-8}--\eqref{m6:s:eq:p78-11-9} are square-zero
at every rational point by \eqref{m6:s:eq:p78-11-11}.
The incidence submersion has rank $55$, with the complete pairing
matrix \eqref{m6:s:eq:p78-11-12}, valid also at dyadic places.
The separate model-orbit and compatible-locus dimensions are tabulated
on p.~\pageref{m6:s:p78:s11-4}; the origin of the $\mathcal P_7$ plane
remains in the rank-$52$ action-plus-slice calculation.
The intrinsic local vanishing proved here applies globally after
the spectral removal of $\mathcal P_6$.

\subsection{The open consequence}\label{m6:s:p78:s12}

\subsubsection{The predecessor opens}\label{m6:s:p78:s12-1}

Retain the already established \(\Omega=\{g:\det g_{55}\ne0\}\) and \(\mathcal O_6\). Define \(V_7,V_8\) by \eqref{m6:s:eq:p78-2-7}, using every entry of the explicit matrices \eqref{m6:s:eq:p78-2-1}, and define \(\mathcal O_7,\mathcal O_8\) by \eqref{m6:s:eq:p78-2-8}. These are open subsets of \(G\) because rank-at-least conditions are open. They are stable under \(H\times H\), since the ranks are preserved by \(P\times P^t\), and under transpose because both matrices in \eqref{m6:s:eq:p78-2-1} are symmetric.

The identities are

\begin{equation}
\begin{aligned}
\mathcal O_7&=\mathcal O_6\sqcup\mathcal P_7,\\
\mathcal O_8&=\mathcal O_7\sqcup\mathcal P_8,\\
\mathcal C_{78}:=\mathcal O_8\setminus\mathcal O_6
 &=\mathcal P_7\sqcup\mathcal P_8.
\end{aligned}
\label{m6:s:eq:p78-12-1}
\end{equation}

The rank comparison in Section~\ref{m6:s:p78:s2-4} proves these identities:
the non-$\Omega$ predecessors dominate $R_7$ or $R_8$, and no later
lexicographic rank matrix can dominate the target.

The closed complements needed for distribution theory are precisely

\begin{equation}
\begin{array}{c|c|c}
\text{ambient space}&\text{open part}&\text{closed complement}\\\hline
\mathcal O_7&\mathcal O_6&\mathcal P_7\\
\mathcal O_8&\mathcal O_7&\mathcal P_8\\
\mathcal O_8&\mathcal O_6&\mathcal C_{78}\\
\mathcal C_{78}&\mathcal P_7&\mathcal P_8.
\end{array}
\label{m6:s:eq:p78-12-2}
\end{equation}

The last row applies because \(\mathcal P_7=\mathcal C_{78}\cap\mathcal O_7\). All four rows are transpose-stable. The common geometric slice \(\mathcal C=\mathcal O_8\setminus\Omega\) is a different closed complement; its first cell was treated in Section~\ref{m6:s:sec:p006}.

\paragraph{The relative two-cell filtration}\label{m6:s:p78:s12-2}
In $\mathcal C_{78}=\mathcal O_8\setminus\mathcal O_6$, the
$\mathcal P_7$ part is open and $\mathcal P_8$ is its closed
complement, by \eqref{m6:s:eq:p78-12-2}. Applying
\eqref{cf:support-sequence} to the intrinsic theorems
\eqref{m6:s:eq:p78-7-2} and \eqref{m6:s:eq:p78-9-21} gives
\begin{equation}\label{m6:s:eq:p78-12-6}
D'_{\widetilde\eta}(\mathcal C_{78})=0.
\end{equation}

\paragraph{The global open consequence}\label{m6:s:p78:s12-3}
Since $\mathcal C_{78}$ is closed in $\mathcal O_8$ with open
complement $\mathcal O_6$, the same support sequence gives
\begin{equation}\label{m6:s:eq:p78-12-5}
T|_{\mathcal O_6}=0\quad\Longrightarrow\quad T|_{\mathcal O_8}=0
\end{equation}
for any original global odd equivariant $T$. This completes
Proposition~\ref{m6:s:prop:p007-p008}.

\subsection{Supplementary coordinate identities}\label{m8:s:app:A}
The expanded coordinate catalogues are in
\texttt{anc/model8\_calculations.pdf}, calculation M8-C08.
They are obtained from the following input formulas.
\begin{center}
\begin{tabular}{p{.31\linewidth}p{.57\linewidth}}
\toprule
Calculation & Input and exact output used in the article\\\midrule
$S_2$ reconstruction & The global chart and two Haar descents of
\S\ref{m8:s:sec:s2}; composing them gives the twisted integral
with phase $\operatorname{tr}X_L+\operatorname{tr}X_R+2r$.\\
Siegel boundary & The affine-fiber reconstruction of
\S\ref{m8:s:s3:descent}, with density $|\det K_0|^{-2}$
and phase $\operatorname{tr}(K_0^{-1}A)$.\\
$S_4$ normal family & The product
\eqref{m8:s:s45:eq3.2}; the complete stabilizer
\eqref{m8:s:s45:eq4.2} has exponent
$2\kappa y-(T_0)_{21}r$.\\
$S_6/S_7$ matrices & The products
\eqref{m8:s:v:e7-3}, \eqref{m8:s:v:e7-8}; their transpose
compensator is \eqref{m8:s:v:e7-14}.\\
Full intersections & The entire $J_6,J_7$ parameter matrices
\eqref{m8:s:v:e4-2}, \eqref{m8:s:v:e4-5}; both second
Levi projections are \eqref{m8:s:v:e4-3}, \eqref{m8:s:v:e4-6}.\\
Nonlinear kernels & The complete group coordinates
\eqref{m8:s:v:e5-4}, group law \eqref{m8:s:v:e5-7}, and
right inverse-transpose maps \eqref{m8:s:v:e5-8},
\eqref{m8:s:v:e5-14}.\\\bottomrule
\end{tabular}
\end{center}

\paragraph{Inverse-factor reconstruction}\label{m8:s:v:sA}
For either full intersection element $j$, the right unipotent is
recovered from
\[
 \bigl(\pi_L^-(w_i j w_i)(w_i j w_i)^{-1}\bigr)^t.
\]
Together with \eqref{m8:s:c:normal-transporter}, this formula
retains every nonlinear and scalar term. Calculation M8-C08
contains the full $8\times8$ conjugated arrays, their exact inverse
factors, the expanded kernel matrices and multiplication checks.
These programs implement the whole-group identities proved above.

\paragraph{Mixed-Levi character check}\label{m8:s:v:sB-2}
The additional $S_7$ unipotent direction is essential in the
full stabilizer test. For the explicit mixed-Levi example of
calculation M8-C08, with $\kappa=2/3$, $\ell=5/7$,
$y_0=3f/5$, its entire trace-plus-scalar exponent is
\[
 2\kappa y_0-2\ell f=-22f/35.
\]
The $K$-compatibility test uses the full stabilizer in
\eqref{m8:s:v:e6-10}.

\subsection{Terminal representatives in scalar coordinates}\label{m8:s:t:expanded-representatives}
The block matrices in \eqref{m8:s:t:eq3.1} and the column-position
list in \S\ref{m8:s:t:sub3.1} specify every scalar entry: if that
list is $(\rho_i(1),\ldots,\rho_i(8))$, then
\[
 w_i=\sum_{j=1}^8 E_{\rho_i(j),j},\qquad i=9,13,14,17,18.
\]
All coefficients are $+1$, so the five $8\times8$ arrays are
recovered directly from the stated monomial data.

\subsection{Logarithmic stabilizer tables}\label{m8:s:m:log-stab}

At
\[
x=m(A,B,\mu)w_im(C,D,\nu)=mw_in,
\]
the full stabilizer consists of
\begin{equation}
{
 u=m\exp(N)m^{-1},\quad
 v=n^t\exp(N')n^{-t},\quad
 N'=-(w_i^{-1}Nw_i)^t\in\mathfrak u.
}
\label{m8:s:m:eq:stabilizer}
\end{equation}
To obtain the following tables, start with all ten coordinates in $N$ free. Requiring \(-(w_i^{-1}Nw_i)^t\) to be strictly block-upper, together with the already automatic orthogonal Lie-algebra condition, gives exactly these solutions.

$S_{10}$: four independent parameters $p,r,a,c$.
\[
\begin{array}{c|cccc}
&X&y&Z&z\\\hline
N&\begin{pmatrix}p&0\\r&0\end{pmatrix}&0&
\begin{pmatrix}a&0\\c&0\end{pmatrix}&0\\[4pt]
N'&\begin{pmatrix}-p&0\\-a&0\end{pmatrix}&0&
\begin{pmatrix}-r&0\\-c&0\end{pmatrix}&0
\end{array}
\]

$S_{11}$: five independent parameters $p,q,r,a,z$.
\[
\begin{array}{c|cccc}
&X&y&Z&z\\\hline
N&\begin{pmatrix}p&q\\r&0\end{pmatrix}&0&
\begin{pmatrix}a&0\\0&0\end{pmatrix}&z\\[4pt]
N'&\begin{pmatrix}-z&-q\\-r&0\end{pmatrix}&0&
\begin{pmatrix}-a&0\\0&0\end{pmatrix}&-p
\end{array}
\]

$S_{12}$: six independent parameters $q,y,a,b,d,z$.
\[
\begin{array}{c|cccc}
&X&y&Z&z\\\hline
N&\begin{pmatrix}0&q\\0&0\end{pmatrix}&y&
\begin{pmatrix}a&b\\0&d\end{pmatrix}&z\\[4pt]
N'&\begin{pmatrix}-z&-a\\-d&-y\end{pmatrix}&0&
\begin{pmatrix}-q&0\\0&0\end{pmatrix}&-b
\end{array}
\]

$S_{15}$: six independent parameters $p,q,r,h,a,z$.
\[
\begin{array}{c|cccc}
&X&y&Z&z\\\hline
N&\begin{pmatrix}p&q\\r&h\end{pmatrix}&0&
\begin{pmatrix}a&0\\0&0\end{pmatrix}&z\\[4pt]
N'&\begin{pmatrix}0&-a\\0&0\end{pmatrix}&h&
\begin{pmatrix}-q&z\\0&r\end{pmatrix}&-p
\end{array}
\]

$S_{16}$: seven independent parameters $q,y,a,b,c,d,z$.
\[
\begin{array}{c|cccc}
&X&y&Z&z\\\hline
N&\begin{pmatrix}0&q\\0&0\end{pmatrix}&y&
\begin{pmatrix}a&b\\c&d\end{pmatrix}&z\\[4pt]
N'&\begin{pmatrix}0&-q\\0&0\end{pmatrix}&c&
\begin{pmatrix}-a&z\\-y&d\end{pmatrix}&-b
\end{array}
\]

\subsection{Additional closure witnesses}\label{m8:s:app:C}
\subsubsection{Additional S4 and S5 closure witnesses}
\paragraph{Closure relations}\label{m8:s:s45:sub8.2}

Rank can only decrease on passing to a boundary. Thus neither \(S_4\) nor \(S_5\) lies in the closure of the other. None of \(S_1,S_2,S_3\) lies in their closures. Also \(S_8\) is incomparable with each of \(S_4,S_5\): its middle four-by-four rank is larger and its terminal two-by-two rank is smaller.

Conversely, \(S_4\) lies in the closures of \(S_1,S_2,S_3\). The following explicit degenerations exhibit the inclusions. In ordinary eight-by-eight matrix units set
\begin{equation}\label{m8:s:s45:eq8.3}
N_2=E_{62}+E_{73},\qquad N_3=E_{53}+E_{64}.
\end{equation}
These square to zero and belong to \(\mathfrak{so}(Q)\). For \(\varepsilon\ne0\), direct rank calculation gives
\begin{equation}\label{m8:s:s45:eq8.4}
(I_8+\varepsilon N_2)w\in S_2,
\qquad
(I_8+\varepsilon N_3)w\in S_3.
\end{equation}
For \(\varepsilon,\delta\ne0\),
\begin{equation}\label{m8:s:s45:eq8.5}
(I_8+\varepsilon N_2)(I_8+\delta N_3)w\in S_1.
\end{equation}
Each family tends to \(w\). The factorization of \S\ref{m8:s:s45:sub3} and stability of the rank sets under left \(P\) and right \(P^t\) transport these degenerations to every point of \(S_4\). Transpose proves the analogous statements for \(S_5\).

\subsubsection{Additional S6 and S7 closure witnesses}
\paragraph{The three degeneration paths}\label{m8:s:v:s2-4}

In ordinary matrix units put
\[
Y=E_{35}+E_{46},\qquad Z=E_{24}-E_{57}.
\]
The identities are
\[
Y^tQ+QY=0,\quad Z^tQ+QZ=0,\quad Y^2=Z^2=0.
\]
Consequently \(I+\varepsilon Y\) and \(I+\varepsilon Z\) preserve \(Q\). For nonzero \(\varepsilon\), their relevant products have ranks
\[
R(w_6(I+\varepsilon Y))=
\begin{pmatrix}5&4&2\\3&2&1\\2&2&1\end{pmatrix},
\]
\[
R((I+\varepsilon Y^t)w_6)=
\begin{pmatrix}5&3&2\\4&2&2\\2&1&1\end{pmatrix},
\]
\begin{equation}
R(w_7(I+\varepsilon Z))=
\begin{pmatrix}5&3&2\\3&2&1\\2&1&1\end{pmatrix}.
\label{m8:s:v:e2-10}
\end{equation}
These are, respectively, the \(S_4\), \(S_5\), and \(S_6\)
rank patterns. Taking the limits \(w_6,w_6,w_7\) proves the stated
boundary relations between rank strata.

\section{Rational central forms and exceptional finite data}\label{app:exceptional}
\paragraph{Rational cokernels and subgroup lifting}\label{m9:o:sec:2.6}
With $V_F=F^\times/(F^\times)^2$, the rational cokernels of the
central maps are as follows. The kernel generators use
\eqref{m9:vector-form}--\eqref{m9:common-cover}.
\begin{center}\small
\begin{tabular}{lll}
\toprule
Map & Algebraic kernel & $F$-point cokernel\\\midrule
$S\times\Gm\to G_{9,v},G_{9,s}$ & $\mu_2^2$ & $V_F^2$\\
$S\times\Gm\to\widehat G_{9,s}$ & $\langle(c_s,1)\rangle$ & $V_F$\\
$q_s:\widehat G_{9,s}\to G_{9,s}$ & $\langle(\bar c_v,-1)\rangle$ & $V_F$\\
$S\times\Gm^2\to\mathcal C_9$ & $\mu_2^2$ & $V_F^2$\\
$p_v,p_s$ & split connected $\Gm$ & $0$\\
$G_{9,v},G_{9,s}\to\operatorname{PSO}_{12}$ & split connected $\Gm$ & $0$\\
\bottomrule
\end{tabular}
\end{center}
For a finite central isogeny of split maximal tori $T'\to T$,
Hilbert~90 gives $H^1(F,T')=0$, so the connecting map onto
$H^1(F,K)$ is surjective. Applying this inside the groups proves
the finite rows. For the connected-kernel rows Hilbert~90 gives
surjectivity directly.

For $q_s$ the connecting class is $[\chi_1(m)]$ and
$s_1(a)=\alpha_1^\vee(a)$ is a rational representative. Its full
subgroup preimage is \eqref{m9:full-subgroup-preimage}; it maps onto
all of $H_{9,s}(F)$, including nonsquare classes in $\PGL_2(F)$.
The simply connected covers have smaller subgroup images:
the vector image has square determinant, and the semispin image
misses nonsquare $\PGL_2$ classes. The proof uses the natural-cover
preimage or the connected-kernel graph $R_9$, both of which retain
the full rational acting subgroup. Section M9-R of
\path{anc/model9_calculations.pdf} gives explicit independent
representatives and their cocycles for every finite row.

\paragraph{An explicit subgroup-preserving operation for literal Model 9}\label{m9:o:sec:3.4}

The literal group also admits a subgroup-preserving anti-involution. Let $w$ reverse the six two-dimensional blocks, and put
\[
B=\operatorname{diag}(I_2,-I_2,I_2,I_2,-I_2,I_2),\qquad A=Bw.
\]
Then $A\in G_v(F)$, $A^t=A$, and $A^2=1$. Define
\begin{equation}
\tau_H(g)=A g^t A^{-1}.
\label{m9:o:eq:3.7}
\end{equation}
This has square one and preserves $H_v$. Conjugating the transposed radical by $w$ changes the three phase contributions to $(-\tr X,-\tr Y,\kappa  z)$; conjugation by $B$ corrects the first two signs. Hence
\[
\xi_\kappa\circ\tau_H=\xi_\kappa.
\]
On $H_{v,0}$ it acts by $a↦a^t$. Right translation $g=xA$ gives
\[
\tau_H(xA)=x^tA,\qquad h(xA)\tau_H(k)=h xk^t A.
\]
Thus right translation by $A$ identifies its doubled distribution problem with the transpose problem. The associated automorphism $\operatorname{Ad}(A)\circ(g↦g^{-t})$ is dualizing by Proposition~\ref{m9:connected-duality}.

\subsubsection{All 48 nonempty NW compatibility loci}
Table M9-K48 of \path{anc/model9_calculations.pdf} records every
minimal word, codimension, full unipotent-stabilizer dimension,
compatible-manifold dimension, complete nontrivial matching equation,
and elimination route. On both sides impose the four matrix determinant
conditions and $p_8q_8\ne0$ of \eqref{m9:phase-open}.
The group has dimension $67$ on either form; the common cover adds
one central dimension. The table contains $13$ Siegel rows, treated together by Siegel
descent, and $35$ boundary rows, treated using the original ambient
distribution after their predecessors vanish.

The $492$-row JSON/CSV in \path{anc/model9/} contain the other $444$ rows. Their equations force $p_8=0$ or $q_8=0$, or block
both determinant monomials of a required $2$ by $2$ phase matrix.
The root-indexed rule \eqref{m9:compatibility} gives the full character
test on every row. The sufficiency of an unblocked determinant and
the rational witnesses are explained in \S\ref{m9:o:sec:5.4}; the rational coverage and ambient continuation are proved in \S\ref{m9:boundary-section}.

\subsection{Exceptional relative-position data}\label{supp:e7-enumeration}
\subsubsection{Relative Bruhat enumeration}\label{gk:o:sec:C}
The Weyl group is enumerated by its integral root-lattice action.
Projection to the unique element with no left or right $I$-descent gives
\[
 |W|=2\,903\,040,\qquad |W_I\backslash W/W_I|=52\,416.
\]
The double-coset sizes $8,32,64$ occur $1\,152,12\,096,39\,168$ times;
their weighted sum is $|W|$. Independently, the weight with simple-root
coordinates $(22,30,42,60,46,32,16)$ has stabilizer $W_I$ and an
orbit of size $362\,880$. Burnside's formula gives
\[
 (362880+3\cdot17280+4\cdot1152)/8=52416.
\]
The full records and source-only generator are in
\path{anc/model10_gk_exhaustive/model10_gk_bundle.zip.xz}. Each record
contains the word, signed root action, full compatibility equations,
and an integral witness or empty-locus certificate. The linear
nonemptiness criterion of \S\ref{gk:o:sec:D.3} therefore gives the
complete $98=69+1+28$ partition.

\paragraph{Lift and sign conventions}\label{gk:o:sec:C.1}
We use
\[
 n_i=\exp(E_i)\exp(-F_i)\exp(E_i),\qquad
 \epsilon(a,b)=(-1)^{\sum_i a_ib_i+\sum_{i>j,\ i\sim j}a_ib_j},
\]
with $[e_\alpha,e_\beta]=\epsilon(\alpha,\beta)e_{\alpha+\beta}$,
$[e_\alpha,e_{-\alpha}]=-h_\alpha$, $E_i=e_{\alpha_i}$ and
$F_i=-e_{-\alpha_i}$. The four diagonal trace root vectors are minus
the physical matrix units, so the six character coefficients are
$(1,1,-1,-1,-1,-1)$. The signed actions satisfy the square and braid
relations and the root-sum bracket identities. The count $23$ refers
to these chosen lifts; $98$ refers to the full Levi-variable loci.

\paragraph{Closure order}\label{gk:o:sec:C.2}
For minimal double-coset representatives,
\[
 \overline{Pn_wP}=\bigcup_{v\leq w}Pn_vP,
\]
where $\leq$ is ordinary Bruhat order restricted to minimal
representatives. This follows from the order-preserving
minimal/maximal projections in \cite[Theorem~1 and \S2]{HS}.
The same archive supplies the exact cover graph, with $52\,416$
vertices and $358\,608$ edges; its transitive closure has
$620\,106\,707$ strict comparable pairs. This graph records the
finite order used in the support induction, whose geometric
interpretation is given by the closure theorem.

\paragraph{The $\tau$ action}\label{gk:o:sec:C.3}
Since $b_-r\in L$, one has
\[
 \tau(Pn_wP)=Pn_{w^{-1}}P,\qquad
 \tau(n_w)=b_-r n_w^{-1}r^{-1}b_-^{-1}.
\]
Of the $98$
nonempty loci, $70$ coarse cells are fixed and $14$ pairs are exchanged.
This action preserves relative length and the full compatibility locus.

\subsection{The maximal-E6 Mackey data}
\paragraph{The full finite packet table}\label{ep:o:sec:9.2}

In this table \(U_-,U_0,U_+\) are the numbers of U-roots of grades \(-1,0,1\). The integer
\begin{equation}
d_\lambda=\sum_{r\in U:\lambda(r)<0}c(r)
\label{ep:o:eq:9.5}
\end{equation}
is the exponent of the Jacobian on \(U_{\lambda,\geq0}\backslash U\). Also
\begin{equation}
\chi_\lambda(t(a))=a^{k_\lambda},\quad
k_\lambda=\lambda_2+\lambda_5+\lambda_7,\quad
\delta_{Q_\lambda}(t(a))=|a|^{9k_\lambda}.
\label{ep:o:eq:9.6}
\end{equation}
A witness \(b_j\) means that the indicated character root has positive \(\lambda\)-grade. It belongs to the inducing unipotent kernel, on which the inducing representation is trivial but \(\xi\) is nontrivial. That entire H-orbit therefore has zero twisted Jacquet contribution.

The complete grade counts, flag roots and Jacobian exponents are in
\path{anc/model10/e6/packet_table.md}. The two compatible packets have
\[
 (U_-,U_0,U_+;d_\lambda,k_\lambda)_A=(12,36,12;-12,3),\qquad
 (U_-,U_0,U_+;d_\lambda,k_\lambda)_B=(26,34,0;-8,1).
\]
For the other $23$ size-two packets it suffices to retain one
positive-grade character root. The finite list is
\begin{center}\small
\begin{tabular}{@{}r l l@{}}
\toprule
No. & $\lambda$ & Positive-grade character root or compatible packet\\\midrule
1 & $(0,0,0,0,0,0,1)$ & $b_5$\\
2 & $(0,0,0,0,1,-1,0)$ & $b_4$\\
3 & $(0,1,1,-1,0,0,0)$ & $b_2$\\
4 & $(1,1,-1,0,0,0,0)$ & $b_1$\\
5 & $(-1,1,0,0,0,0,0)$ & $b_3$\\
6 & $(1,0,0,-1,1,0,0)$ & $b_1$\\
7 & $(-1,0,1,-1,1,0,0)$ & $b_2$\\
8 & $(0,0,-1,0,1,0,0)$ & $b_4$\\
9 & $(1,0,0,0,0,-1,1)$ & $b_1$\\
10 & $(-1,0,1,0,0,-1,1)$ & $b_2$\\
11 & $(0,1,0,-1,0,1,0)$ & $b_5$\\
12 & $(0,0,-1,1,0,-1,1)$ & $b_3$\\
13 & $(0,1,0,-1,1,-1,1)$ & $A$ (equal flags)\\
14 & $(0,1,0,0,0,-1,0)$ & $b_3$\\
15 & $(0,0,1,-1,0,0,1)$ & $b_2$\\
16 & $(1,0,-1,0,0,0,1)$ & $b_1$\\
17 & $(-1,0,0,0,0,0,1)$ & $b_5$\\
18 & $(0,0,1,-1,1,-1,0)$ & $b_2$\\
19 & $(1,0,-1,0,1,-1,0)$ & $b_1$\\
20 & $(-1,0,0,0,1,-1,0)$ & $b_4$\\
21 & $(1,1,0,-1,0,0,0)$ & $b_1$\\
22 & $(-1,1,1,-1,0,0,0)$ & $b_2$\\
23 & $(0,1,-1,0,0,0,0)$ & $b_3$\\
24 & $(0,0,0,-1,1,0,0)$ & $b_6$\\
25 & $(0,0,0,0,0,-1,1)$ & $B$\\
\bottomrule
\end{tabular}
\end{center}

For all 23 size-two packets other than B, the table supplies a character witness. Since each packet is one $H$-orbit and $L$ preserves $\xi$,
the witness at the canonical representative eliminates the whole orbit.

\section{Distributional Ginzburg--Rallis arguments}
\label{app:jlx-repairs}
Following the geometric framework of \cite{JLX}, we establish the
nilpotent-wavefront formulation used here through open-local
reductions, rational transporters, and conormal calculations. Ordinary
distribution arguments give the support reductions, while ambient
nilpotent wavefront gives the conormal exclusions. The proofs apply
over characteristic-zero local fields in every residue characteristic.
\subsection{Open-local distributional reduction following JLX \S4}
\label{repair:jlx4}

With the groups, action, character, and transpose of
Theorem~\ref{cf:gr-nw}, we have $H=S\times S$, $S=D_3(\GL_2)U_{2,2,2}$,
$\eta=\epsilon\boxtimes\epsilon\cdot\chi(\det a_1\det a_2)$, and
$\widetilde\eta(\tau)=-1$.  The field has characteristic zero; its
residue characteristic is arbitrary. The following argument concerns
distributions on an open subset itself.

Write a matrix in $2\times2$ blocks as $x=(x_{ij})_{1\leq i,j\leq3}$.
On $x_{33}=y\in\GL_2$, put
\[
 \sigma=x_{22}-x_{23}y^{-1}x_{32},\qquad N=\sigma y^{-1},\qquad
 \mathcal O_1=\{y\in\GL_2,\ \sigma\ne0\}.
\]
Block Gaussian elimination gives
\[
 \operatorname{rank}x_{\{3,4,5,6\},\{3,4,5,6\}}
     =2+\operatorname{rank}\sigma.
\]
Both rectangular lower corners have rank two.  Consequently
$\mathcal O_1$ is exactly $G^2_{\mathrm{open},1}$ of JLX, and is an
$H$- and transpose-stable open subset of $\GL_6$.
Direct block multiplication gives
\[
 N((s_1,s_2)x)=a_1N(x)a_1^{-1}.
\]
Under transpose, $N(x^t)=y^tN(x)^t y^{-t}$; hence its trace and determinant
are also unchanged.  There is therefore a disjoint locally closed
decomposition
\begin{equation}\label{repair:jlx4-partition}
 \begin{gathered}
 \mathcal O_1=\mathcal V_2\sqcup\mathcal V_3\sqcup\mathcal V_4,\qquad
 \mathcal V_2=\{\det N\ne0\},\\
 \mathcal V_3=\{\det N=0,\ \operatorname{tr}N\ne0\},\qquad
 \mathcal V_4=\{\det N=\operatorname{tr}N=0,\ N\ne0\}.
 \end{gathered}
\end{equation}
Here $\mathcal V_2$ is open, $\mathcal V_3$ is closed in
$\mathcal O_1\setminus\mathcal V_4$, and $\mathcal V_4$ is closed in
$\mathcal O_1$.  The map $x\mapsto N$ is a submersion: with $y$ and the
other blocks fixed one can vary $x_{22}$ freely.  A nonzero nilpotent
$2\times2$ matrix is a smooth point of the equations
$\operatorname{tr}N=\det N=0$ (their differentials are independent).
Thus $\mathcal V_4$ is smooth of dimension $34$, and $\mathcal V_3$ is
smooth of dimension $35$.

\paragraph{Slices and density conventions.}
The action maps from the following slices are surjective submersions:
\[
 M_2=\{\operatorname{diag}(A,D,I_2):A,D\in\GL_2\},\qquad
 M_3^0=\{\operatorname{diag}(X,\lambda,I_2):X\in\GL_3,
            X_{33}=0,\ \lambda\in F^\times\},
\]
and
\[
 M_4^0=\{z(X)=\operatorname{diag}(X,1,w):X\in\GL_3,
                 X_{33}=0\},\qquad
 w=\begin{pmatrix}0&1\\1&0\end{pmatrix}.
\]
To obtain the respective targets $\mathcal V_2,\mathcal V_3,\mathcal V_4$,
we construct local analytic normalizers at every rational point. First use the $(4,2)$ Gaussian chart to kill the blocks meeting
the last block, leaving $\operatorname{diag}(g_4,y)$, and use a common
left $y^{-1}$ to put $y=I_2$.  Write $g_4=\left(\begin{smallmatrix}A&B\\
C&D\end{smallmatrix}\right)$.  If $D$ is invertible, upper and lower
$U_{2,2}$ eliminate $B,C$ and leave
$\operatorname{diag}(A-BD^{-1}C,D,I_2)$.
If $D$ has rank one and nonzero trace $\lambda$, choose a kernel and an
image vector to conjugate it to $\operatorname{diag}(0,\lambda)$.
Locally these vectors can be chosen as
$(I-D/\lambda)v$ and $(D/\lambda)u$ for fixed vectors $v,u$; their independence is
an open condition.  Upper and lower $U_{2,2}$ then kill the second
column of $B$ and second row of $C$, giving $M_3^0$.
If $D$ is nonzero nilpotent, a cyclic basis $(v,Dv)$ conjugates it to
$J=\left(\begin{smallmatrix}0&0\\1&0\end{smallmatrix}\right)$.
Right multiplication by $D_3(w)$ puts $y=w$ and $D=Jw=E_{22}$.
The same column and row elimination gives $M_4^0$. Because all
inverted quantities are nonzero on the indicated charts, these
choices give local analytic normalizers.  On $M_3^0$ choose $v=e_1,u=e_2$, so the kernel--image
basis is the identity.  On $M_4^0$, the preliminary left $w$ replaces
$y=w$ by $I_2$ and $E_{22}$ by $E_{12}$; choose the cyclic vector $e_2$,
whose cyclic basis is $w$.  The resulting common conjugation by $w$
followed by the right $D_3(w)$ has total left and right multiplier the
identity. The row and column eliminations also fix each slice point,
so we can choose the local normalizer to be the identity on the slice.  Differentiating
$x=h(x)m(x)$ at a slice point shows
$T_x\mathcal V_i=\operatorname{Lie}(H)x+T_xM_i^0$.
This proves submersivity.

We first take slice restriction for generalized functions, the dual
of compactly supported smooth densities, for which pullback by a
surjective submersion is injective.  On $H\times M$, equivariance for
left translation on $H$ factors the pullback as the smooth character on
$H$ times a generalized function on $M$; conjugation on $H$ leaves that
character unchanged.  A subgroup preserving $M$ therefore induces its restricted
character on the slice. We return to distributions by choosing a
density invariant under that subgroup and transpose. In the present
applications these densities can be chosen as follows:
on $\mathcal O_1$ use Haar density of $\GL_6$;
on $\mathcal V_3$ divide that density by the transverse density of
$\det N$; on $\mathcal V_4$ divide by the transverse densities of
$(\operatorname{tr}N,\det N)$.  Both functions are invariant and have
independent differentials where used.  On $M_2$ use product Haar density.
On $M_3^0$ use Lebesgue density on the eight free entries of $X$ times
$d^\times\lambda$.  The subgroup below preserves it: its unipotent
actions are shears, and its torus weights on the free entries multiply
to one. On $M_4^0$, the subsequent support reductions use ordinary
orbital criteria, while the wavefront argument uses the invariant
density on the ambient $\mathcal V_4$.

\paragraph{The ordinary part: $\mathcal V_2$ and $\mathcal V_3$.}
The transporter of $M_2$ consists of simultaneous conjugation by $a$
on $(A,D)$, embedded as
$(D_3(a),D_3(a^{-t}))$; all upper unipotent blocks vanish.
Its character is one since
$\chi(\det a\det a^{-t})=1$.
Transpose acts as $(A,D)\mapsto(A^t,D^t)$.
Lemma~\ref{m8:s:c:two-matrix} therefore kills every odd equivariant
distribution on $\mathcal V_2$.

For $M_3^0$, take the subgroup whose action on $X$ is
\[
 X\longmapsto l n_1 X n_2^t l^{-1},\quad
 l=\operatorname{diag}(t,1,t),\quad
 n_i=I+c_iE_{23}+d_iE_{13}.
\]
It fixes $\lambda$, has character $\psi(d_1+d_2)$, and transpose sends
$(X,\lambda)$ to $(X^t,\lambda)$ and exchanges the two unipotent
factors while inverting $t$. This subgroup is the restriction of
$S\times S$: $l$ extends to $D_3(\operatorname{diag}(t,1))$ and the
other factor extends to its inverse transpose.  The two determinant
characters cancel.
Write $u=X_{13},v=X_{31},p=X_{23},q=X_{32}$.
The successive closed loci $u=v$ and $u=v=0$ are stable.  On $u\ne v$,
\[
 U=I+s(uE_{13}+pE_{23}),\qquad
 V=I+s(vE_{31}+qE_{32})
\]
satisfy $UX=XV$, and the stabilizer $(U,V^{-t})$ has character
$\psi(s(u-v))$.  Choosing a suitable $s$ makes this character nontrivial. Since
the stabilizer is unipotent, its modular character is one, and the
algebraic stabilizer criterion gives vanishing on this locus.
On $u=v\ne0$, the two unipotent factors remove $X_{11},X_{12},X_{21}$,
giving the slice
\[
 \begin{pmatrix}0&0&u\\0&a&p\\u&q&0\end{pmatrix},\qquad
 u,a\ne0,\quad p,q\in F.
\]
Since the elimination is analytic with pivots $u$, this slice is
submersive.  The torus acts by
$(q,p)\mapsto(tq,t^{-1}p)$ and fixes $u,a$; transpose exchanges $p,q$.
The whole opposite-weight plane lemma~\ref{m8:s:c:plane} kills its odd
distributions.  On $u=v=0$, invertibility gives $a=X_{11},p,q\ne0$;
the unipotents remove $X_{12},X_{21},X_{22}$, leaving
\[
 \begin{pmatrix}a&0&0\\0&0&p\\0&q&0\end{pmatrix}.
\]
The same torus and transpose, now on $pq\ne0$, give vanishing (localize
the plane lemma by $pq$, or use multiplicative Haar measure on each
nonzero product fiber).  Ordinary open--closed exactness proves
vanishing on $M_3^0$, and injective slice restriction proves vanishing
on $\mathcal V_3$.

\paragraph{The nilpotent part and its smooth residual support.}
For $z(X)\in M_4^0$, write
\[
 X=\begin{pmatrix}a&b&u\\c&d&p\\v&q&0\end{pmatrix}.
\]
An element transporting $M_4^0$ into itself must have its $(1,3)$ and
$(2,3)$ unipotent blocks and their right counterparts zero, because
$y=w$ is invertible.  If its first reductive block is $a_0$, its right
reductive block satisfies $a_2^t=wa_0^{-1}w$.  The middle block then
requires $a_0J a_0^{-1}=J$, so
\[
 a_0=\begin{pmatrix}\alpha&0\\\beta&\alpha\end{pmatrix},\qquad
 \alpha\ne0,\qquad t=\beta/\alpha.
\]
Its determinant character is one.  Set
$A=\left(\begin{smallmatrix}a&b\\c&d\end{smallmatrix}\right)$,
$B=(u,p)^t$, $C=(v,q)$, and
$L=\left(\begin{smallmatrix}1&0\\t&1\end{smallmatrix}\right)$.
After the unique second-column and second-row eliminations, the exact
remaining action has the form
\begin{equation}\label{repair:jlx4-transporter}
 u'=u,\quad v'=v,\quad p'=p+tu,\quad q'=q-tv,
\qquad
 A'=LAL^{-t}+\ell\,(v,q-tv)+(u,p+tu)^t\rho,
\end{equation}
where $\ell\in F^2$ and $\rho\in(F^2)^*$ are free.  To derive these formulas, multiply the upper block $b_0$ and lower block
$R$ with
$b_{0,*2}=t(u,p+tu)^t$ and
$R_{2*}=-t(v,q-tv)$; absorb the resulting quadratic term into $\rho$.
The scalar $\alpha$ acts trivially on $X$.
In particular $u,v$ are invariants; on $u=v=0$, $a,p,q$ are invariants.
These transporter formulas hold over $F$ throughout the slice.

The sets
\[
 \begin{gathered}
 W_1=\{u=v\},\qquad W_2=\{u=v=0\},\\
 W_{3,1}=\{u=v=0,\ p=q\},\qquad
 W_{3,2}=\{u=v=0,\ a+pq=0\}.
 \end{gathered}
\]
are relatively $H$-stable in $M_4^0$.
Since $H\times M_4^0\to\mathcal V_4$ is an open surjection, their
saturations are closed: the inverse image of each saturation is
exactly $H$ times the displayed closed slice locus.
Transpose exchanges $u,v$ and $p,q$, so the saturations are also
transpose-stable. To check smoothness, use the local analytic
normalizers above: they give a retraction to $M_4^0$ that is the
identity on the slice and hence a submersion.  Near a slice point
$HW_{3,2}$ is the inverse image of
$u=v=0,a+pq=0$ under this retraction.  These are three independent
equations on $M_4^0$; indeed the coefficient of $da$ in the last is one.
Thus
\[
 Y:=HW_{3,2}\quad\hbox{is a closed smooth submanifold of }\mathcal O_1,
 \qquad \dim Y=34-3=31,
\]
and
$T_xY=\operatorname{Lie}(H)x+T_xW_{3,2}$ at slice points.

The stabilizer $UX=XV$ just used kills $\mathcal V_4\setminus HW_1$
with phase $\psi(s(u-v))$.
On $W_1\setminus W_2$ the common value $u=v$ is nonzero.
Take $t=(q-p)/(2u)$ in \eqref{repair:jlx4-transporter}, so $p'=q'$;
the free $\ell,\rho$ remove $A_{11},A_{12},A_{21}$ and leave a symmetric
$X$.  The pivot $u$ is nonzero at every point under consideration, so
each orbit contains a symmetric representative. We can therefore
apply the odd orbital criterion: if $h x=z=z^t$, then
$\tau(x)=\tau(h)^{-1}h x$ and
$\eta(\tau(h)^{-1}h)=1$.  The resulting involution acts by a positive
modulus on the one-dimensional orbital density line, so cannot act
by $-1$.  Lemma~\ref{cf:constructive-restriction} gives vanishing on
$HW_1\setminus HW_2$.

On $W_2$, invertibility of $X$ says
\[
 \det X=-apq\ne0;\qquad a,p,q\in F^\times.
\]
The free $\ell,\rho$ eliminate $b,c,d$, leaving
$X(a,p,q)=\left(\begin{smallmatrix}a&0&0\\0&0&p\\0&q&0\end{smallmatrix}\right)$.
For a stabilizer with parameter $t$, take
$\ell_1=ta/q$ and $\rho_1=-ta/p$; the last diagonal equation is solved
by choosing $\ell_2,\rho_2$.  Multiplication of its two unipotent blocks
gives the exact phase
\[
 \psi\bigl(t(1/q-1/p)(a+pq)\bigr).
\]
It is nontrivial whenever $p\ne q$ and $a+pq\ne0$.  The corresponding
stabilizer is unipotent, including its reductive block
$\left(\begin{smallmatrix}1&0\\t&1\end{smallmatrix}\right)$, so its
modular character is one.  This kills the complement of
$HW_{3,1}\cup HW_{3,2}$ in $HW_2$.
The normal forms with $p=q$ are symmetric, so the same orbital
criterion kills $HW_{3,1}\setminus HW_{3,2}$.
The ordinary open--closed sequences now leave an ambient distribution
supported on the closed smooth $Y$.

\paragraph{The ambient nilpotent-wavefront contradiction.}
At a normal form $x=z(X(-pq,p,q))$, both $p$ and $q$ are nonzero.
Under the trace pairing on $\operatorname{Mat}_6$, set
\[
 C=E_{11}+pE_{23}+qE_{32}-pqE_{44}.
\]
For a tangent vector to $W_{3,2}$,
$\operatorname{tr}(C\,\delta x)=\delta a+q\,\delta p+p\,\delta q=0$.
Moreover
\[
 Cx=xC=pq\operatorname{diag}(-1,1,1,-1,0,0).
\]
For $L\in\operatorname{Lie}(S)$, its three diagonal $2\times2$ blocks
are the same matrix.  The displayed diagonal entries have zero sum
in each repeated coordinate, and all off-diagonal entries vanish;
hence
$\operatorname{tr}(CxL^t)=\operatorname{tr}(xCL)=0$.
Thus $C$ annihilates both infinitesimal orbit directions and lies
in the ambient conormal $CN^{\mathcal O_1}_{Y,x}$.
Its image under either group trivialization of the cotangent fiber is
$Cx$ or $xC$, with characteristic polynomial
\[
 \lambda^2(\lambda-pq)^2(\lambda+pq)^2,
\]
so $Cx$ is nonnilpotent.
At $hx=s_1xs_2^t$, the transported cotangent vector is
$s_2^{-t}Cs_1^{-1}$; its product with $hx$ is conjugate to $Cx$.
It is consequently non-nilpotent at every point of $Y$.
Every orbit in $Y$ meets such a normal form, since elimination of
$b,c,d$ preserves $a,p,q$.

If a nonzero ambient equivariant distribution is supported on $Y$,
choose a support point and use the direct Fourier conormal proof of Proposition~\ref{cf:conormal} there.
With the wavefront convention including the zero cotangent vector at a
support point, the entire conormal fiber is in its wavefront fiber.
The explicitly non-nilpotent
cotangent vector above contradicts the assumed ambient nilpotent wavefront.
Thus every odd equivariant NW distribution on $\mathcal O_1$ is zero.

\paragraph{The second open by an exact symmetry.}
Let $J_6$ be the $6\times6$ antidiagonal permutation matrix and put
$z=\operatorname{diag}(1,1,-1,-1,1,1)$, $K=zJ_6$.
The map
\[
 \Theta(x)=Kx^{-t}K
\]
is an involutive group automorphism, since $K^t=K$ and $K^2=I$,
and commutes with transpose.  It preserves $S$.  In normalized
unipotent coordinates its two adjacent upper blocks are
$w c^t w$ and $w b^t w$; their traces sum to
$\operatorname{tr}b+\operatorname{tr}c$.
Its reductive block is $wa^{-t}w$, of determinant $(\det a)^{-1}$.
Thus it preserves $\epsilon$ and replaces $\chi$ by $\chi^{-1}$.
The complementary-minor nullity identity for an invertible matrix gives
\[
 R(\Theta x)=
 \begin{pmatrix}r_{22}+2&r_{24}\\r_{42}&r_{44}-2\end{pmatrix}
 \quad\text{if}\quad
 R(x)=\begin{pmatrix}r_{44}&r_{42}\\r_{24}&r_{22}\end{pmatrix}.
\]
For completeness, this identity follows by viewing a corner of
$x^{-1}$ as the induced inverse map between coordinate subspaces:
its kernel is the intersection of the image of the complementary
subspace with the opposite coordinate subspace; rank--nullity then
gives $\operatorname{rank}(x^{-1})_{I,J}
=\operatorname{rank}x_{J^c,I^c}+|I|+|J|-6$.
Hence $\Theta$ fixes the rank matrix
$\left(\begin{smallmatrix}4&2\\2&2\end{smallmatrix}\right)$ and
exchanges $\left(\begin{smallmatrix}3&2\\2&2\end{smallmatrix}\right)$
with $\left(\begin{smallmatrix}4&2\\2&1\end{smallmatrix}\right)$.
It maps $G^2_{\mathrm{open},1}$ onto $G^2_{\mathrm{open},2}$.
Because its differential at the identity is $A\mapsto-KA^tK$, it
preserves the nilpotent cone on cotangent fibers.  Pullback by
$\Theta$ transfers the second open problem to the first with
$\chi^{-1}$.  The proof above is uniform in $\chi$, so the second
open-local NW vanishing follows as well.

All constructions use operations in the characteristic-zero field.
The only division by $2$ is $t=(q-p)/(2u)$, which is legitimate also
when the residue characteristic is two.

\subsection{Distributional Ginzburg--Rallis arguments following JLX \S\S5--7}
\label{cf:jlx57-repair}\label{repair:jlx57}
Throughout this subsection $F$ has characteristic zero, with arbitrary residue
characteristic. Write $U=U_{2,2,2}$, $S=D_3(\GL_2)U$, $H=S\times S$,
$N_H=U\times U$, and
\[
 \ell(B,C)=\tr B_{12}+\tr B_{23}+\tr C_{12}+\tr C_{23}.
\]
We use the point action $(s_1,s_2)x=s_1xs_2^t$ and work with
complex-valued distributions.

\paragraph{A closed compatible locus.}
On a parabolic double coset $Z$ define the vector-bundle map
\[
 A_x:\mathfrak u\oplus\mathfrak u\longrightarrow\Mat_6,
 \qquad A_x(B,C)=Bx+xC^t.
\]
Since $N_H$ is normal in $P\times P$, the dimension of $K_x=\ker A_x$
is constant on $Z$, and the kernels form an algebraic vector bundle $K$.
Put
\[
 Y_0=\{x\in Z:\ell|_{K_x}=0\}.
\]
As the zero locus of a section of $K^*$, it is closed in $Z$,
$H$-stable, and transpose-stable. In characteristic zero the stabilizer
in $N_H$ is connected unipotent, exponential identifies its Lie algebra with
$K_x$, and its character is $\psi\circ\ell$. Thus every point outside $Y_0$
has a unipotent stabilizer on which the character is nontrivial; the ordinary
orbital criterion clears $Z\setminus Y_0$.

\paragraph{The codimension-two locus in JLX \S5.}
Here $Z=Z_6$ has rank matrix $\left(\begin{smallmatrix}3&2\\2&1\end{smallmatrix}\right)$.
It has dimension $34$ and $\dim K_x=2$. These dimensions follow directly
from the permutation with nonzero entries in columns $(1,4,5,2,3,6)$:
the $P\cap wP^tw^{-1}$ Lie algebra has dimension $14$, and
$\mathfrak u\cap w\mathfrak u^tw^{-1}$ has the two root spaces
$(2,4),(3,5)$. The normality argument above makes the count valid throughout
this single rational parabolic double coset.

To bring every point to the slice $\mathfrak z_6$ of JLX \S5 by
$H(F)$-operations, first factor the
rank-one lower-right block $R=pq^t$ by taking a nonzero column $p$
and dividing its corresponding row by a nonzero pivot. For any nonzero
vector $v=(v_1,v_2)^t$, the following invertible rational matrix sends
$v$ to $e_2$:
\[
 L(v)=\begin{cases}
 \left(\begin{smallmatrix}-v_2/v_1&1\\1/v_1&0\end{smallmatrix}\right),
       &v_1\ne0,\\
 \left(\begin{smallmatrix}1&-v_1/v_2\\0&1/v_2\end{smallmatrix}\right),
       &v_2\ne0.
 \end{cases}
\]
Use $(D_3(L(p)),D_3(L(q)))$ to obtain $R=E_{22}$.
First subtract the current $x_{6j}$ times column six from columns
$j\le4$, using $I-\sum_{j\le4}x_{6j}E_{j6}$ as the right factor;
then subtract $x_{i6}$ times row six from rows $i\le4$, using
$I-\sum_{i\le4}x_{i6}E_{i6}$ on the left. Both factors belong to $U$.
Write the remaining lower corner as
\(
 \left(\begin{smallmatrix}M&v\\w&0\end{smallmatrix}\right)
\), with $v,w$ nonzero by the two rectangular rank conditions. The remaining
unipotents change $M$ by $\alpha w+v\beta$. The space of these changes is
exactly the kernel of the functional $M\mapsto u^tMz$, where $u^tv=wz=0$.
The lower-corner rank-three condition says $u^tMz=0$, so $M$ is killed.
Explicitly choose $v_iw_j\ne0$, put $\beta=-M_{i,*}/v_i$ and
$N=M+v\beta$. Row $i$ of $N$ is zero, and its other row, indexed
$k\ne i$, is proportional to $w$. Set $\alpha_i=0$ and
$\alpha_k=-N_{k,j}/w_j$. Then $M+v\beta+\alpha w=0$.
The additions use the $U_{23}$ entries $E_{35},E_{45}$ on the two
sides. The four choices of $(i,j)$ cover all points.

On this slice set
\[
 v=(x_{35},x_{45})^t,\quad w=(x_{53},x_{54}),\quad
 u=(x_{45},-x_{35})^t,\quad z=(-x_{54},x_{53})^t,
 \quad \Phi=u^t(X_{12}-X_{21})z.
\]
A frame of $K$ consists of the element in JLX Lemma 5.1 and the element
whose only nonzero upper blocks are
\[
 B_{12}=X_{12}zu^t,\qquad C_{12}=-X_{21}^tuz^t.
\]
The two character values are $x_{35}-x_{53}$ and $\Phi$, and the
frame is independent everywhere: $v,w$ are nonzero, and
$X_{12}z\ne0$, since otherwise $x(0,0,z_1,z_2,0,0)^t=0$.
The differential of $\Phi$ contains the nonzero functional
$u^t(dX_{12})z$; it is independent of $d(x_{35}-x_{53})$.
Consequently the zero section is transverse at every point of the slice.
By $H$-equivariance and its rational coverage, $Y_0$ is smooth of dimension
$32$ throughout $Z_6$.

Following the notation of JLX, the first character value defines
the closed locus $Z_{6,1}$ and its zero sublocus $Z_{6,2}$. For an
intrinsic algebraic description, put $R=X_{33}$ and $J_R=\operatorname{adj}(R)$; on $Z_6$,
$\rank R=1$. Then
\[
 Z_{6,1}=\{J_RX_{23}J_R=J_RX_{32}J_R\},\qquad
 Z_{6,2}=\{J_RX_{23}J_R=J_RX_{32}J_R=0\}.
\]
Indeed, for $R=E_{22}$ these are precisely the displayed scalar equations.
Under $(D_3(g)u,D_3(h)v)$, each of the two matrix expressions transforms
by $M\mapsto(\det g\det h)^2h^{-t}Mg^{-1}$: the extra unipotent terms
vanish because $RJ_R=J_RR=0$. Thus both loci are $H$-stable and
algebraically closed in $Z_6$; transposition exchanges the two equations.
On $Z_{6,1}\setminus Z_{6,2}$ write $a=x_{35}=x_{53}\ne0$.
Lower-triangular diagonal block operations kill $x_{45},x_{54}$; the
$b,b',d,d'$ unipotents, with pivot $a$, give the slice $\mathfrak z$ in JLX
Lemma 5.2. On it
\(\Phi=-a^2(y_2-w_2)\), so the surviving locus is precisely
$Y_1=H\mathfrak z'$ with $y_2=w_2=b$.

This locus has the further rational slice
\[
 z=\begin{pmatrix}
 A&0&0&y&0&0\\
 0&B&0&b&0&0\\
 0&0&0&0&a&0\\
 w&b&0&0&0&0\\
 0&0&a&0&0&0\\
 0&0&0&0&0&1
 \end{pmatrix},\qquad
 a\ne0,\quad D=Ab^2+Byw\ne0.                   \tag{JLX57.1}
\]
For $b\ne0$, adding row four to row one and column four to column one
kills the two off-diagonal entries of the upper-left block. For $b=0$,
invertibility gives $Byw\ne0$, and adding row four to row two and column
four to column two does the same. These are allowed $b$ and $b'$ operations.
Thus $H\times\mathfrak z''\to Y_1$ is onto. It is a submersion at every
point: its differential has rank $32$ on both charts $b\ne0$ and $b=0$.
To verify this rank on both charts, \texttt{jlx57\_allpoints.py}
builds the $38$ tangent columns
$Bz,zB^t,\partial z/\partial(A,B,y,w,a,b)$ from the sixteen elementary
Lie algebra matrices of $S$. Eliminating with $a,b$ on the first
chart and $a,B,y,w$ on the second gives $32$ unit pivots and a zero
remainder over the corresponding localized rings. The identities
therefore specialize at every allowed $F$-point.

\paragraph{Distributional restriction and the density line.}
Let
$\rho:H\times\mathfrak z''\to Y_1$ be the preceding submersion and
$K_0=F^\times\rtimes\{1,\tau\}$, acting on the slice by
\[
 t(y,w)=(ty,t^{-1}w),\qquad \tau(y,w)=(w,y),
 \qquad t,\tau\text{ fix }(A,B,a,b).
\]
It acts on $H\times\mathfrak z''$ by conjugation on $H$ and this slice action.
Choose a positive relative density on the fibres of $\rho$ invariant under
left $H$ and $K_0$. At $(1,z)$, the relative density line is that
of the six-dimensional space $\ker d\rho$. On $b\ne0$, the unit elimination above
identifies that kernel with its six coordinates
\[
 B_{\rm left,diag,11},\quad
 C_{\rm right,diag,22},\quad C_{\rm right,12,22},\quad
 C_{\rm right,23,11},\quad dB,\quad da.
\]
All six have torus weight zero, so the density from these coordinates is
torus-invariant. On $yw\ne0$ the torus action is free and has the rational
section $y=1$; transport a density from this section by the torus. These two
opens cover the slice, since $b=0$ implies $Byw\ne0$. Patch them on the
torus-invariant clopen partition determined by $|b|\ge |yw|$ and its
complement, and average with its transpose to ensure $\tau$-invariance.
Extend the resulting densities at $(1,z)$ to $(h,z)$ by left translation.
The construction preserves both invariances, since $K_0$ normalizes $H$.

Now fibre integration of compactly supported locally constant test
functions defines a surjection
\(\rho_*:\mathcal D(H\times\mathfrak z'')\to\mathcal D(Y_1)\).
Local submersion coordinates prove both local constancy and surjectivity;
compact support makes each fibre integral finite. Its dual is an injective
pullback of distributions. With the chosen density, this pullback is $H$- and $K_0$-equivariant.
A left $(H,\eta)$-equivariant distribution on
$H\times\mathfrak z''$ factors uniquely as twisted Haar measure on $H$
tensored with a distribution on $\mathfrak z''$. The group $H$ is
unimodular: conjugation on each $\Mat_2$ unipotent block has determinant
one. Conjugation by the displayed torus, and exchange of the two factors
by transpose, preserve its Haar measure. Since $\eta$ restricts to one on this torus, an odd distribution on
$Y_1$ therefore injects into \emph{torus-invariant, transpose-odd}
distributions on the open
slice (JLX57.1).

\paragraph{Localization on the actual open slice.}
First consider torus-invariant, swap-odd distributions on $F^2$ for
$t(y,w)=(ty,t^{-1}w)$, and localize on the product $yw$.
For a nonzero product $c$, the fibre is one torus orbit, its invariant
measure is $d^\times y$, and $y\mapsto c/y$ preserves that measure.
On the zero fibre, restriction to the punctured cross is necessarily
$C(d^\times y-d^\times w)$. For
$f=1_{\mathcal O_F^2}$ and a uniformizer $\varpi$,
\[
 (f\circ t-f)|_{w=0}=1_{\varpi^{-1}\mathcal O_F\setminus\mathcal O_F},
 \qquad
 (f\circ t-f)|_{y=0}=-1_{\mathcal O_F^\times},\qquad t=\varpi.
\]
Invariance therefore gives
$2C\operatorname{vol}^\times(\mathcal O_F^\times)=0$.
Thus $C=0$; the distribution remaining at the origin is even.
The factor $2$ is in $\mathbb C$, so this calculation applies in every
residue characteristic.

Localize (JLX57.1) first on the invariant base $(A,B,a,b)$, which transpose
fixes coordinate by coordinate. Each fibre is the open set
\(\{(y,w):Byw\ne-Ab^2\}\).
Localizing that fibre on $yw$ reduces to the fibres just treated, with the
excluded product absent. This proves vanishing on the open slice
and gives the reduction needed for JLX Lemma 5.2. It remains to
treat the closed locus, which we do by the following rational argument.

\begin{lemma}[Intrinsic vanishing on the remaining third-cell locus]
\label{cf:jlx-z62}
With the action and character of Theorem~\ref{cf:gr-nw},
$\Dist_{\widetilde\eta_\chi}(Z_{6,2})=0$.
\end{lemma}
\begin{proof}
Write $E_{ij}$ for the six-by-six elementary matrix and
\[
 J'=\begin{pmatrix}0&1&0\\1&0&0\\0&0&1\end{pmatrix},\qquad
 X=\begin{pmatrix}p&q&a\\r&s&b\\c&d&0\end{pmatrix}\oplus J'.
                                                        \tag{JLX57.4}
\]
Every rational point of $Z_{6,2}$ reaches this form by the following
allowed operations. First use the rational $\mathfrak z_6$ reduction
proved above. Its entries $x_{45}=k$ and $x_{54}=l$ are nonzero:
the two rectangular rank conditions make $(x_{35},x_{45})$ and
$(x_{53},x_{54})$ nonzero, and their first entries vanish here.
Multiply on the left by $D_3(\diag(l^{-1},1))$ and on the right
by $D_3(\diag(k^{-1},1))^t$. Now $x_{45}=x_{54}=1$.
In the following list every coefficient is read from the current matrix:
\begin{enumerate}[label=(\alph*)]
\item subtract $x_{15}$ and $x_{25}$ times row four from rows one
and two, using $I-x_{15}E_{14}-x_{25}E_{24}$ on the left;
\item subtract $x_{51}$ and $x_{52}$ times column four from columns
one and two, using $(I-x_{51}E_{14}-x_{52}E_{24})^t$ on the right;
\item subtract $x_{14}$ and $x_{24}$ times row five from rows one
and two, using $I-x_{14}E_{15}-x_{24}E_{25}$ on the left;
\item subtract $x_{41}$ and $x_{42}$ times column five from columns
one and two, using $(I-x_{41}E_{15}-x_{42}E_{25})^t$ on the right.
\end{enumerate}
The first two operations give unit pivots for the last two. All
belong to $U$, and together they yield (JLX57.4), with $\det X\ne0$,
proving rational coverage.

If $a\ne c$, put
\[
 A=aE_{13}+bE_{23},\qquad B=-cE_{13}-dE_{23}.
\]
Then $A^2=B^2=0$, $AX+XB^t=0$, and
$(I+tA)X(I+tB)^t=X$ for every $t\in F$.
Both factors belong to $U$, and their actual group character is
$\psi(t(a-c))$. Hence it is nontrivial for some $t$.

Suppose next that $a=c\ne0$. Subtract $(p/a)$ and $(r/a)$ times
row three from rows one and two; then subtract the current
$(x_{12}/a)$ times column three from column two. Put $t_0=-b/a$.
Multiply on the left by $D_3(\begin{psmallmatrix}1&0\\t_0&1\end{psmallmatrix})$;
clear the resulting $x_{41},x_{42}$ by subtracting those multiples
of column five from columns one and two, and clear $x_{64}$ by
subtracting its multiple of column six from column four.
Next put $s_0=-d/a$ and multiply on the right by
$D_3(\begin{psmallmatrix}1&0\\s_0&1\end{psmallmatrix})^t$.
Clear $x_{14},x_{24}$ by subtracting their multiples of row five
from rows one and two, and clear $x_{46}$ using row six.
The row/column clearings are respectively the same $U$-operations
as above, or $I-x_{46}E_{46}$ and its right-factor counterpart.
Direct multiplication gives the symmetric representative
\[
 \begin{pmatrix}0&0&a\\0&\beta&0\\a&0&0\end{pmatrix}\oplus J',
 \qquad \beta=s-\frac{rd+qb}{a}+\frac{pdb}{a^2},\qquad
 \det X=a^2\beta\ne0.                              \tag{JLX57.5}
\]
Thus the only nonunit pivot in this case is the explicitly nonzero $a$.

Finally suppose $a=c=0$. Now $\det X=pbd$, so $p,b,d\ne0$.
Subtract $(r/b)$ times column three from column one, then subtract
$(q/d)$ and $(s/d)$ times row three from rows one and two. These
are allowed $U$-operations, and yield
\[
 X(p,b,d)=pE_{11}+bE_{23}+dE_{32}+E_{45}+E_{54}+E_{66}.
                                                        \tag{JLX57.6}
\]
Put $L=E_{21}+E_{43}+E_{65}$. The stabilizer contains the following
two one-parameter algebraic unipotent subgroups:
\[
 \begin{array}{c|c|c}
 A_i&B_i&\eta_\chi(I+tA_i,I+tB_i)\\\hline
 -(p/d)E_{13}-bE_{25}-E_{46}&L&\psi(-t(p+d)/d)\\
 -L&(p/b)E_{13}+dE_{25}+E_{46}&\psi(t(p+b)/b).
 \end{array}                                             \tag{JLX57.7}
\]
Indeed $A_i^2=B_i^2=0$ and $A_iX+XB_i^t=0$, which also give
$A_iXB_i^t=0$ and the exact group identity
$(I+tA_i)X(I+tB_i)^t=X$. The off-block factors are in $U$;
the other factors are $D_3(I_2\pm tE_{21})$.
Their diagonal determinants equal one. Reading the two adjacent
upper-block traces after this factorization gives exactly the
characters in (JLX57.7).
Solving for the full group stabilizer gives $(s_0,v_0^{-t})$, where
\[
 \begin{aligned}
 K&=(p/d)E_{13}+bE_{25}+E_{46},\qquad
 J_0=(p/b)E_{31}+dE_{52}+E_{64},\\
 s_0&=z_0I+z_1L+z_2K+z_3E_{23}+z_4E_{45},\\
 v_0&=z_0I+z_2L^t+z_1J_0+(d/b)z_3E_{32}+z_4E_{54},
 \quad z_0\in F^\times,\quad z_1,z_2,z_3,z_4\in F.
 \end{aligned}
\]
To obtain every solution, write $v_0=X^{-1}s_0X$ and impose its
lower-block triangularity and equality of diagonal blocks. These
linear entry equations give exactly the displayed parametrization;
the certificate checks a rank-$27$ minor $pb^3d^4\ne0$ and all
remaining equations. The diagonal blocks are
$z_0I_2+z_1E_{21}$ and $z_0I_2+z_2E_{12}$, respectively. Factoring
out these blocks gives the full stabilizer character
\[
 \eta_\chi(s_0,v_0^{-t})=
 \psi\left(\frac{z_2}{z_0}\frac{p+d}{d}
           -\frac{z_1}{z_0}\frac{p+b}{b}\right);
\]
the two determinant factors cancel. It is identically trivial
precisely when $b=d=-p$. At these remaining points (JLX57.6) is symmetric.

On each displayed unipotent subgroup, both modular characters
in \eqref{cf:orbital} are trivial: the restriction
of either modular character is a continuous homomorphism
$(F,+)\to\mathbb R_{>0}$, and $F$ is the union of compact additive
subgroups, whose positive real character is trivial. The nontrivial
phases therefore exclude the corresponding ordinary orbital line.
In either symmetric case the orbital line is transpose-even by
the positive density-line argument. More explicitly, if $h x=y=y^t$
is any of the above rational normalizations, let
$\alpha(h_1,h_2)=(h_2,h_1)$. Then
\[
 x^t=\alpha(h)^{-1}h x,\qquad
 \eta_\chi(\alpha(h)^{-1}h)=1.
\]

There are no odd equivariant distributions on any orbit of
$\widetilde H=H\rtimes\langle\tau\rangle$ in $Z_{6,2}$.
This locus is the closed algebraic locus defined above, and
$\widetilde H$ is an algebraic finite extension acting algebraically.
Lemma~\ref{cf:constructive-restriction} verifies the constructive-action
and rank-one $l$-sheaf hypotheses of Bernstein--Zelevinsky
\cite[Theorem 6.9]{BZ76}. It gives the global intrinsic vanishing. Since all operations use
field arithmetic with nonzero stated pivots, and the characters and
distributions take values in $\mathbb C$, the argument applies in
every residue characteristic.
\end{proof}

The open-complement argument preceding Lemma~\ref{cf:jlx-z62} kills
$Z_6\setminus Z_{6,2}$. Since $Z_{6,2}$ is closed in $Z_6$, the
ordinary support sequence and this lemma give intrinsic vanishing
on all of $Z_6$. This is the ordinary third-cell input in
Theorem~\ref{cf:gr-nw}; intrinsic vanishing applies to the
closed-support antecedent. The subsequent NW exclusions are
applied to the original ambient distribution on their indicated opens.

\paragraph{The smooth support in JLX \S6.}
Direct calculation gives the displayed permutation in JLX \S6.1 the rank matrix
$\left(\begin{smallmatrix}3&2\\1&1\end{smallmatrix}\right)$, the stratum
$G^5_{{\rm open},2}\setminus G^4_{\rm open}$.
Let $Z$ be this stratum. Its dimension is $33$ and $\dim K_x=3$:
the permutation columns are $(1,4,3,6,5,2)$, its parabolic stabilizer has
dimension $15$, and its unipotent intersection has root spaces
$(2,3),(2,4),(2,6)$. Every point has the rational normal representative
\[
 g=\diag(x,I_2,y)W\diag(z,I_2,w),\qquad
 W=\begin{pmatrix}E_{11}&E_{22}&0\\0&E_{11}&E_{22}\\E_{22}&0&E_{11}\end{pmatrix}.
\]
Indeed $P=S\{\diag(x,I_2,y)\}$ and the parabolic Bruhat decomposition is
over $F$. In the three-dimensional kernel frame given by the free
parameters $(c_{12},c_{22},d_{22})$ of the stabilizer in JLX \S6.1, the
character-section coordinates are
\[
 \phi_1=-y_{21}/\det y-z_{21},\qquad
 \phi_2=y_{11}/\det y,\qquad \phi_3=-z_{22}.
\]
The row-two, column-one entry $y_{21}$ in $\phi_1$ follows from
the displayed stabilizer calculation.
The zero equations are exactly $y_{11}=z_{22}=0$, $y_{12}z_{21}=1$.
At every such point their differential on the variations
$(dz_{21},dy_{11},dz_{22})$ has determinant $1/\det y\ne0$.
Thus the section is transverse everywhere, and $Y_0$ is a smooth closed
submanifold of $Z$ of dimension $30$. Since $Z$ is closed in
$G^5_{{\rm open},2}$, the same is true in that ambient open.

Every compatible normal representative reaches the following slice.
Write $y=\left(\begin{smallmatrix}0&u\\v&k\end{smallmatrix}\right)$,
$z_{21}=u^{-1}$, $z_{22}=0$. Multiply on the left by
$D_3(\left(\begin{smallmatrix}1&0\\-k/u&1\end{smallmatrix}\right))$,
and on the right by the transpose of $u(0,c',0)$, where
\[
 c'=\begin{pmatrix}(k/u)(-w_{12},w_{11})/\det w\\0\quad0\end{pmatrix}.
\]
Direct multiplication gives
\[
 X=\begin{pmatrix}
 a&b&0&c&0&0\\d&e&0&f&0&0\\0&0&1&0&0&0\\
 0&0&0&0&p&q\\1&0&0&0&0&0\\0&0&0&0&r&s
 \end{pmatrix},\qquad D_1=bf-ce\ne0,\quad D_2=ps-qr\ne0.       \tag{JLX57.2}
\]
Conversely, kill its entries $a,d$ by the left $d$-block
$\left(\begin{smallmatrix}-a&0\\-d&0\end{smallmatrix}\right)$.
The resulting matrix is exactly the compatible normal representative
with $x=\left(\begin{smallmatrix}b&c\\e&f\end{smallmatrix}\right)$,
$y=z=\left(\begin{smallmatrix}0&1\\1&0\end{smallmatrix}\right)$,
and $w=\left(\begin{smallmatrix}r&s\\p&q\end{smallmatrix}\right)$.
This proves $Y_0=H\mathfrak G''_R$ at all rational points.

The tangent sum
\(V=T_X\mathfrak G''_R+\mathfrak sX+X\mathfrak s^t\)
has dimension $30$ at every point. A particularly compact complete
calculation is its annihilator for the entrywise pairing: writing its six
free parameters as $u_1,u_2,u_3,u_4,m,n$, it is exactly
\[
 \begin{pmatrix}
 0&0&0&0&0&0\\0&0&0&0&0&0\\
 0&0&-u_1&-u_2&-(rs m+s^2n)/D_2&(r^2m+rsn)/D_2\\
 0&0&0&0&0&0\\u_1&u_2&u_3&u_4&m&n\\0&0&0&0&0&0
 \end{pmatrix}.                                                   \tag{JLX57.3}
\]
Equivalently, the accompanying all-points elimination builds all $42$
tangent columns from the sixteen Lie algebra matrices and ten slice
coordinates, finds eighteen scalar unit pivots and six two-by-two pivots
with determinant $\pm D_2$, and leaves zero. This proves specialization
at every point, including points where the orbit dimension drops.
Since $V\subseteq T_XY_0$ and both have dimension $30$, they are equal.

Take the matrix $X'$ displayed in JLX \S6.1. Its entrywise linear functional
$(X')^t$ is (JLX57.3) with $u_1=u_2=1$ and all other parameters zero.
It is therefore an ambient conormal direction to $Y_0$, and the
characteristic polynomial of $X'X$ is
$\lambda^4(\lambda-1)(\lambda+1)$, so the matrix is nonnilpotent.
An ambient distribution with nilpotent wavefront supported on $Y_0$
must vanish, by conormal shift invariance (or its local product-coordinate
Fourier proof). This proves the assertion corresponding to Proposition 6.2 in the stated
NW scope; transpose gives the assertion corresponding to Proposition 6.1.

\paragraph{Explicit rational transpose transporters for JLX \S7.1 and \S7.4.}
The compatible representatives contain a common arbitrary block
$x\in\GL_2(F)$, and transpose acts by $x\mapsto x^t$.
Choose $a\in\GL_2(F)$ with $axa^{-1}=x^t$ and use
\[
 h=(D_3(a),D_3(a)^{-t})\in H.
\]
Then $hg=g^t$ and its product character equals one. The rational choice
is explicit for the two compatible block matrices
\[
 g_1=\begin{pmatrix}I_2&0&0\\0&0&x\\0&x&0\end{pmatrix},\qquad
 g_4=\begin{pmatrix}0&x&0\\x&0&0\\0&0&I_2\end{pmatrix}:
 \quad D_3(a)g_iD_3(a)^{-1}=g_i^t.
\]
Indeed both occurrences of $x$ become $x^t$, and
$\chi(\det a)\chi(\det a^{-1})=1$.
The rational choice
of $a$ is as follows: if $x$ is scalar take $a=I$;
otherwise choose cyclic bases $(v,xv)$ and $(w,x^tw)$, in which $x$ and
$x^t$ have the same companion matrix, and take their change of basis.
The orbital transpose-evenness criterion now applies.

For the remaining JLX cases, direct block multiplication gives the
unipotent obstructions and character-one transporters. The
expressions used here are as follows: in \S7.8
the coefficient of $b_{22}$ is $(y^{-1}x)_{22}$, and $\alpha$ is read as
$a^t$; the detailed
representative in \S7.9 has rank matrix
$\left(\begin{smallmatrix}3&2\\1&0\end{smallmatrix}\right)$.
In \S7.10 an invertible upper-triangular $a$ solving
\[
 a_{11}(yw)_{12}+a_{12}(yw)_{22}=a_{22}(yw)_{21}
\]
always exists: if $(yw)_{22}\ne0$, take $a_{11}=a_{22}=1$; otherwise
$(yw)_{12}(yw)_{21}\ne0$, and take
$a_{11}=(yw)_{21}$, $a_{22}=(yw)_{12}$, $a_{12}=0$.
All transporter pairs have product character one. The orbital
constructibility criterion applies to each whole rank stratum.
The finite rank filtration then proves the intrinsic ordinary assertion
for the final complement in the Ginzburg--Rallis support filtration.
Together with the NW-local \S4 argument, these calculations
justify the ordered support argument used in Theorem~\ref{cf:gr-nw}.

\section{Whittaker and Fourier--Jacobi identities}
\label{app:external-repairs}
\subsection{Bi-Whittaker invariance for the split Levi}
\label{repair:shalika-criterion}
We prove directly the bi-Whittaker distribution invariance used in
\cite[Theorem~7.6 and Corollary~7.7]{Yang} for the split connected
reductive group considered here.

Let $L$ be a split connected reductive $F$-group with a pinning,
$B=TU$ its pinned Borel, and $\psi$ the associated nondegenerate
character of $U(F)$.  Fix a nontrivial additive character $\psi_0$;
in pinned simple-root coordinates,
$\psi(x_\alpha(s))=\psi_0(s)$, and $\psi$ is trivial on the
nonsimple positive root groups.  Let $\iota$ be the pinned duality
involution of \cite[Definition~3.1]{Yang} and
$\theta(g)=\iota(g^{-1})$.  Then every bi-Whittaker distribution is
$\theta$-invariant.  More explicitly, use the acting algebraic group
$J=U\times U^{\mathrm{op}}$, with
\[
 (u_1,u_2)g=u_1gu_2,
 \qquad \eta(u_1,u_2)=\psi(u_1)\psi(u_2).
\]
The inverse on the right when one writes this with ordinary
$U\times U$ instead of $U\times U^{\mathrm{op}}$ is just the usual
left/right translation convention.  We prove that there is no
$(J,\eta)$-equivariant distribution odd under $\theta$.

\paragraph{Normalizer representatives and compatibility.}
Write $n_w$ for the pinned rational Weyl representatives.  Rational
Bruhat factorization gives a representative $n=tn_w$, $t\in T(F)$,
in every $J(F)$-orbit.  One may use the split Bruhat multiplication
isomorphism \cite[Theorem~21.80]{MilneAG}.  The stabilizer in $J(F)$ is the unipotent group
\[
 J_n=\{(nv^{-1}n^{-1},v):
              v\in U(F)\cap n^{-1}U(F)n\}.
\]
Its character is trivial precisely when
\begin{equation}\label{repair:shalika-compatible}
 \psi(nvn^{-1})=\psi(v)
 \qquad(v\in U(F)\cap n^{-1}U(F)n).
\end{equation}

We check that this condition implies $\theta(n)=n$, including all
torus components.  Put
$I=\{\alpha\in\Delta:w\alpha>0\}$.
For $\alpha\in I$, the root group $U_\alpha$ occurs in the
intersection in \eqref{repair:shalika-compatible}.
If $w\alpha$ were nonsimple, the left character would be trivial
while the right character would be nontrivial.  Hence $wI\subset\Delta$.
This also proves directly the needed Weyl-combinatorial fact.
Set $I'=wI$.  If a positive root $\gamma$ lies outside the span of $I$,
its expansion has some simple-root coefficient outside $I$.
The roots $w\alpha$ for $\alpha\notin I$ are negative and none lies
in the span of $I'$; their coefficients outside $I'$ are nonpositive,
with at least one strictly negative.  The roots $w\alpha$ for
$\alpha\in I$ are simple roots in $I'$.  Thus $w\gamma$ has a
strictly negative coefficient outside $I'$ and cannot be positive.
Conversely, $w$ maps every positive root in the subsystem generated
by $I$ to a positive root.  If $w_I$ is the longest element of this
subsystem, it follows that $ww_I$ maps every positive root to a
negative root: $w_I$ reverses the subsystem and preserves the set of
positive roots outside it.  Consequently
\[
 ww_I=w_0,\qquad w=w_0w_I.
\]

The pinned representative identity
\cite[Lemma~2.1, equation~(2.15)]{Yang} gives
$n_wx_\alpha(s)n_w^{-1}=x_{w\alpha}(s)$ for $\alpha\in I$.
This identity of pinned split root homomorphisms holds over $F$.  Substituting it in
\eqref{repair:shalika-compatible} gives
\[
 \psi_0((w\alpha)(t)s)=\psi_0(s)\quad(s\in F),
 \quad\hbox{hence}\quad (w\alpha)(t)=1
 \quad(\alpha\in I).
\]
The last implication uses nontriviality of $\psi_0$: multiplication
by a nonzero scalar is onto $F$.  Therefore
$s_0=w^{-1}(t)$ lies in the full center of the standard Levi $L_I$.
In particular $w_I(s_0)=s_0$, also for disconnected central kernels.
Since $w^{-1}=w_Iw_0$, this implies $w_0(t)=s_0$.

For completeness, the other pinned identities used here are
\[
 \iota(t)=w_0(t^{-1}),\qquad
 \iota(n_w)=(n_{w_0w^{-1}w_0})^{-1}.
\]
The first follows from the torus definition, while the second is
\cite[Theorem~3.4]{Yang}, proved by applying $\iota$ to each pinned
simple-reflection representative and reversing a reduced expression.
For $w=w_0w_I$ one has $w_0w^{-1}w_0=w$, and hence
\[
 \theta(tn_w)=\iota(n_w)^{-1}\iota(t^{-1})
       =n_ww_0(t)=n_ws_0=tn_w.
\]
This proves the normalizer compatibility condition.

\paragraph{The distribution argument and modular factors.}
The involution satisfies $\theta(U)=U$ and
$\psi(\theta(u))=\psi(u)$, since
$\psi(\iota(u))=\psi(u)^{-1}$.
It acts on $J$ by
\[
 \theta_J(u_1,u_2)=(\theta(u_2),\theta(u_1)),
 \qquad \eta\circ\theta_J=\eta.
\]
Thus $\widetilde J=J\rtimes\{1,\theta\}$ is an algebraic group
acting on $L$, and $\eta$ extends to a smooth character
$\widetilde\eta$ with $\widetilde\eta(\theta)=-1$.
Since $J(F)$ is unipotent, it is unimodular. Its finite extension
$\widetilde J(F)$ is also unimodular: the positive Haar modulus of
an involutive automorphism has square one and is one.

At a representative $n$ violating
\eqref{repair:shalika-compatible}, the restriction of
$\widetilde\eta$ to its unipotent stabilizer $J_n(F)$ is nontrivial.
The ambient and stabilizer modular characters are trivial on those
unipotent elements, so the orbital distribution criterion gives zero.
At a compatible representative $n$, the preceding proof gives
$\theta(n)=n$.  The extended stabilizer is
$J_n(F)\rtimes\{1,\theta\}$ and is unimodular, by the same
finite-extension argument.  In particular the modular ratio at its
order-two element $\theta$ is one, while
$\widetilde\eta(\theta)=-1$.  The orbital criterion again gives zero.
The character restriction and this mismatch are unchanged when a
representative is moved within its $J(F)$-orbit, since
$\widetilde\eta$ is a group character.

Finally the algebraic action of $\widetilde J$ on $L$ is constructive
on $F$-points.  Applying the constructibility theorem with the
rank-one sheaf twisted by $\widetilde\eta$ gives zero globally:
the exact inputs are the orbital criterion and constructibility
theorem stated in \cite[Theorems~2.3 and~2.5, Remark~2.4]{JLX}.
Equivalently these are the Bernstein--Zelevinsky results used
throughout the ordinary support arguments of this paper.
Thus the odd part of a bi-Whittaker distribution vanishes, and every
such distribution is $\theta$-invariant.

This proves the distribution input of \cite[Theorem~7.6]{Yang} for
the split Levi, with arbitrary central form and arbitrary
residue characteristic in characteristic zero. 

For its representation-theoretic application here let \(\sigma\) be
an irreducible admissible \emph{unitarizable} \(\psi\)-generic representation;
in particular this applies to the tempered inducing data. Indeed,
\cite[Lemma~III.1.3, Proposition~III.4.1]{Waldspurger} realizes it
as an irreducible subquotient of induction from square-integrable,
hence unitarizable, data. Normalized induction is unitary and an
admissible unitary representation is semisimple
\cite[Propositions~3.1.4 and~2.1.14]{Casselman}, so this subquotient
is unitarizable. Choose an invariant positive Hermitian form linear
in its first variable.
The map \(\overline V\to V^\vee\),
\(\overline v\mapsto (w\mapsto\langle w,v\rangle)\), is complex-linear
and equivariant. A subgroup fixing \(v\) also fixes the functional, so the map
takes values in the smooth dual. To prove surjectivity, observe that a smooth functional
fixed by a compact open \(K\) factors through compact averaging onto
finite-dimensional \(V^K\), where the Hermitian form represents every
linear functional; averaging is self-adjoint. If \(\ell\ne0\) is a
\(\psi\)-Whittaker functional, then
\[
 \overline\ell(\overline v)=\overline{\ell(v)},\qquad
 \overline\ell(\overline\sigma(u)\overline v)
   =\psi(u)^{-1}\overline\ell(\overline v),
\]
since \(\psi\) is unitary. Thus \(\sigma^\vee\) is
\(\psi^{-1}\)-generic. Both hypotheses on the two distinguished
representations in \cite[Proposition~7.5]{Yang} now hold, so its
relative-character argument gives
\(\sigma^\iota\simeq\sigma^\vee\).
Only this tempered, hence unitarizable, instance of
\cite[Corollary~7.7]{Yang} is used in Part~V.

\subsection{The canonical morphism to the opposite standard module}
\begin{proposition}[The image of the opposite-module morphism]
\label{app:konno-image-independent}
Let \(P=MN\) be a parabolic subgroup of a connected reductive group over
a non-Archimedean local field of characteristic zero, let \(\sigma\) be irreducible tempered,
and let \(\lambda\) be real and strictly \(P\)-positive. Put
\(\tau=\sigma e^\lambda\) and \(I=i_P^G(\tau)\). There is a nonzero
\(G\)-morphism
\[
 J':I\longrightarrow i_{\bar P}^G(\tau)
\]
whose image is irreducible and is the unique irreducible quotient of \(I\).
It is defined below from the normalized Jacquet projection, using only
a compactly supported open-cell integral.
\end{proposition}
\begin{proof}
We construct the morphism and identify its image using the
asymptotics of a sequence of central elements.
All Jacquet modules in this proof are normalized.

Choose \(a_0\) in the split center \(A_M(F)\) with
\(|\alpha(a_0)|>1\) for every \(N\)-root and set
\(a_n=a_0^n\), \(\zeta=\omega_\tau(a_0)\). Such an element exists by
taking an integral multiple of an interior rational cocharacter and
evaluating on an inverse uniformizer.
Then \(a_n\bar n a_n^{-1}\to1\) uniformly on compact subsets of
\(\bar N\), and \(|\zeta|=e^\lambda(a_0)\).

Let \(I^\circ\) be the sections supported compactly on the
\(\bar N\)-coordinate of the open cell \(P\bar N\).
It is a \(\bar P\)-submodule, and
\[
 r_{\bar P}(I^\circ)\simeq\tau,\qquad
 f\longmapsto\int_{\bar N}f(u)\,du.
\]
The kernel and surjectivity follow from the compact-induction
coinvariant calculation \cite[Proposition 6.2.1]{Casselman}.
For the normalization, right translation by \(m\) multiplies the
integral by \(\delta_P(m)^{-1/2}\tau(m)\), since conjugation on
\(\bar N\) has Jacobian \(\delta_P(m)^{-1}\);
the normalized \(\bar P\)-Jacquet factor
\(\delta_{\bar P}(m)^{-1/2}\) cancels it.
Exactness embeds this module as \(S\simeq\tau\) in
\(D=r_{\bar P}(I)\).

We isolate \(S\) spectrally. Choose a minimal parabolic inside \(P\).
The normalized geometric lemma
\cite[\S I.3, pp.~243--244]{Waldspurger}, whose integral and
determinant calculation are proved in
\cite[Propositions 6.3.1--6.3.3]{Casselman}, filters \(D\) by
\(W_M\backslash W/W_M\). The identity double coset is \(S\).
A different representative \(w\) gives the grade
\[
 i_{M\cap wPw^{-1}}^M
 \bigl(w\,r_{M\cap w^{-1}\bar Pw}^{M}(\tau)\bigr).
\]
For an exponent \(\chi\) of
\(r_{M\cap w^{-1}\bar Pw}^{M}(\sigma)\), the element
\(b=w^{-1}a_0w\) lies in its Levi's split center and contracts its
radical: each radical root \(\beta\) satisfies
\(w\beta\in\Sigma(\bar N)\).
Applying the normalized tempered criterion
\cite[Proposition III.2.2, pp.~265--266]{Waldspurger}
gives \(|\chi(b)|\le1\).
Induction inside \(M\) has modulus one on central \(a_0\).
Thus every exponent of this grade on \(a_0\) has modulus at most
\(e^{w\lambda}(a_0)<e^\lambda(a_0)\).

To prove the strict inequality, work in the relative Weyl root system. Let \(J\) be the simple roots of \(M\).
The padded \(\lambda\) has zero coroot pairings on \(J\) and positive
ones off \(J\). For a reduced word \(w=s_{i_1}\cdots s_{i_k}\),
\[
 \lambda-w\lambda
 =\sum_{j=1}^k\langle\lambda,\alpha_{i_j}^\vee\rangle
 s_{i_1}\cdots s_{i_{j-1}}\alpha_{i_j}.
\]
All displayed roots are positive. For \(w\notin W_M\) the sum is
nonzero, and it cannot lie in the span of \(J\): \(\lambda\) is
orthogonal to that span and \(w\) preserves its norm.
It therefore has a positive simple-root coefficient outside \(J\).
Roots in \(J\) vanish on \(H(a_0)\), while all others are positive.
This proves strictness; a global-central part of \(\lambda\) cancels.

For each compact open \(K_M\subset M\), exactness of invariants gives
\((D/S)^{K_M}=D^{K_M}/S^{K_M}\). Admissibility makes these spaces
finite dimensional, and the central operator \(a_0\) preserves them.
The quotient spectrum has modulus below \(|\zeta|\), so the generalized
\(\zeta\)-eigenspace of \(D^{K_M}\) is exactly \(S^{K_M}\), with scalar action. Inclusions between fixed spaces
intertwine \(a_0\), so their primary-space complements agree on
intersections. Taking their union gives \(D=S\oplus E\); because
\(M\) commutes with \(a_0\), \(E\) is \(M\)-stable. For each vector in \(E\), its finite-dimensional
\(a_0\)-orbit has eigenvalues of modulus below \(|\zeta|\);
hence \(\zeta^{-n}r_{\bar P}(a_n)\) tends to zero on that orbit.
Let \(p:D\to S\simeq\tau\) be this \(M\)-equivariant projection,
\(q:I\to D\) the Jacquet quotient, and \(j'=p\circ q\).
The normalization is explicit: for \(m\in M\), \(\bar n\in\bar N\),
\[
 j'(I(m\bar n)f)=\delta_{\bar P}(m)^{1/2}\tau(m)j'(f).
\]
Thus
\[
 (J'f)(g)=j'(I(g)f)
\]
is a section of \(i_{\bar P}^G(\tau)\), is right-invariant by every
compact open subgroup fixing \(f\), and defines a \(G\)-morphism.
This is precisely normalized Frobenius reciprocity
\(\Hom_G(I,i_{\bar P}^G\tau)\simeq
\Hom_M(r_{\bar P}I,\tau)\), obtained from
\cite[Theorem~3.2.4, p.~34]{Casselman}; the same \(\bar P\) occurs
in the target induction and source Jacquet module.
Since \(p|_S=\mathrm{id}\) and \(q\) is onto, \(J'\ne0\).
The only integral in this construction is on \(I^\circ\), where
compact support proves its convergence.

By \cite[Proposition~3.1.2]{Casselman},
\(I^\vee=i_P^G(\tau^\vee)\). Put \(D'=r_P(I^\vee)\). The canonical pairing theorem
\cite[Theorem I.4.1, p.~245]{Waldspurger}, proved in
\cite[Proposition 4.2.3, Theorems 4.2.4 and 4.3.3]{Casselman},
gives, for every \(f\in I\), \(\phi\in I^\vee\), and all sufficiently
large \(n\),
\[
 \langle r_{\bar P}(a_n)\bar f,\bar\phi\rangle_{\bar P}
 =\delta_P(a_n)^{1/2}\langle I(a_n)f,\phi\rangle .
\]
We identify its \(S\)-coefficient explicitly. For \(f_0\in I^\circ\),
the open-cell pairing and the substitution
\(\bar n=a_nua_n^{-1}\) give
\[
 \begin{aligned}
 \delta_P(a_n)^{1/2}\zeta^{-n}\langle I(a_n)f_0,\phi\rangle
 &=\gamma(G/M)^{-1}\int_{\bar N}
       \langle f_0(u),\phi(a_nua_n^{-1})\rangle\,du
 \\
 &=\gamma(G/M)^{-1}
       \left\langle\int_{\bar N}f_0(u)\,du,\phi(1)\right\rangle
 \end{aligned}
\]
for all sufficiently large \(n\).
Here the Jacobian is \(\delta_P(a_n)^{-1}\), while the inducing
factor is \(\delta_P(a_n)^{1/2}\zeta^n\).
The second equality is exact because \(f_0\) has compact support and
\(\phi\) is constant near \(1\).
Here \(\gamma(G/M)\) is the measure constant of
\cite[\S 2.1, pp.~385--386]{Konno}.
Every \(s\in S\) has such a lift \(f_0\).
Consequently, for every smooth dual vector, the canonical pairing of
\(s\) with \(\bar\phi\) is
\(\gamma(G/M)^{-1}\langle s,\phi(1)\rangle\).
This identifies the pairing on the leading summand. Since the
\(E\)-term decays spectrally, for arbitrary \(f,\phi\) we obtain
\[
 \lim_{n\to\infty}\delta_P(a_n)^{1/2}\zeta^{-n}
     \langle I(a_n)f,\phi\rangle
 =\gamma(G/M)^{-1}\langle j'(f),\phi(1)\rangle. \tag{*}
\]

Let \(U\subset I\) be a \(G\)-submodule with \(J'(U)\ne0\).
If \(U\) were proper, its admissible quotient would have a nonzero
smooth dual, so choose \(0\ne\phi\in U^\perp\).
Indeed a nonzero quotient vector is fixed by a compact open subgroup;
a functional on that fixed space, composed with compact averaging,
is a nonzero smooth functional. Compact averaging also proves
admissibility of the quotient.
For \(f\in U\), (*) gives \(\langle j'(f),\phi(1)\rangle=0\).
The image \(j'(U)\) is a nonzero \(M\)-submodule of irreducible
\(\tau\): translate a section with a nonzero value of \(J'\) to \(1\).
Thus \(j'(U)=\tau\), and \(\phi(1)=0\).
The same reasoning for every \(G\)-translate of \(\phi\), which
still belongs to \(U^\perp\), gives \(\phi(g)=0\) for all \(g\).
This contradiction proves \(U=I\).
Therefore every proper submodule lies in the proper kernel of \(J'\).
The kernel is the unique maximal submodule and its quotient
\(\operatorname{im}J'\) is irreducible, as claimed.

Whenever the full unnormalized opposite integral is absolutely
convergent and defines an intertwiner, its value at \(1\) induces an
\(M\)-map \(D\to\tau\). It is identity on \(S\) by the open-cell
calculation and zero on \(E\), whose \(a_0\)-spectrum has smaller
modulus. It therefore equals \(p\), and that intertwiner equals \(J'\).

The proof uses complex-coefficient compact averaging, relative roots,
canonical Jacquet pairing, and the sequence limit (*); all cited
results apply in residue characteristic two.
\end{proof}

\bibliographystyle{amsalpha}
\bibliography{references}
\end{document}